\documentclass[preprint,12pt]{elsarticle}
 
\usepackage{amsmath,amssymb,amsthm}
\usepackage{graphicx}
\usepackage[margin=1in]{geometry} 
\usepackage{subfigure}
 \usepackage{url}
\usepackage{placeins}
\usepackage{multirow}
 \usepackage{lineno}
\numberwithin{equation}{section}
\usepackage{color}

\begin{document}
\journal{J Computational and Applied Mathematics}
\begin{frontmatter}

\title{Non-negatively constrained least squares and parameter choice by the residual periodogram for the inversion of dielectric relaxation spectra: Supplementary Materials}

\author[label1]{Jakob Hansen}
\ead{jkhanse2@asu.edu}
\author[label1]{Jarom Hogue}
\ead{jdhogue@asu.edu}
\author[label1]{Grant Sander}
\ead{gksander@asu.edu}
\author[label1]{Rosemary Renaut\corref{cor1}}
\ead{renaut@asu.edu}
\ead[url]{math.asu.edu/~rosie}
\cortext[cor1]{Corresponding Author: Rosemary Renaut, 480 965 3795 }
\address[label1]{School of Mathematical and Statistical Sciences, Arizona State University, Tempe, AZ 85287-1804, USA}

\begin{abstract}
This  document contains supplementary derivations and discussions not provided in the submitted paper. Additional results for the NCP and L-Curve comparisons with higher noise levels are given. 
\end{abstract}


\end{frontmatter}

\section{Tables for the NLS Fitting}
To carry out the NLS fitting data were chosen to provide aligned DRTs in $s$-space.  To obtain this we note that the lognormal DRT is given by 
$$g_{\mathrm{LN}}(t | \mu, \sigma)= \frac{1}{t \sigma \sqrt{2 \pi}} \exp(- \frac{(\ln (t) -\mu)^2}{2 \sigma^2}).$$
It is centered at $t_0=\exp(\mu-\sigma^2)$ and can be written  in terms of $t_0$. We have 
\begin{align}
\ln(t_0)&= \mu -\sigma^2 \\
\mu&=\ln(t_0)+\sigma^2 \\
 \ln(t)-\mu&= \ln(t)-\ln(t_0) -\sigma^2 =\ln(t/t_0)-\sigma^2 \\
  \frac{(\ln (t) -\mu)^2}{2 \sigma^2}& = \frac{(\ln(t/t_0)-\sigma^2)^2}{2\sigma^2} \\
  &= (\frac{\ln(t/t_0)-\sigma^2}{\sqrt{2}\sigma})^2
  \end{align}
 Let $s=\ln(t/t_0)$ then
\begin{align}g_{\mathrm{LN}}(t | \mu, \sigma)&= \frac{1}{t \sigma \sqrt{2 \pi}} \exp (-(\frac{s-\sigma^2}{\sqrt{2}\sigma})^2)\\
&= \frac{1}{t_0 \sigma \sqrt{2 \pi}} \exp (-(\frac{s-\sigma^2}{\sqrt{2}\sigma})^B-s) \\
f_{\mathrm{LN}}(s)&=  \frac{1}{t_0 \sigma \sqrt{2 \pi}}  \exp (-\frac{(\sigma^4+s^2)}{2 \sigma^2}) \end{align}
At $t=t_0$ $s=0$ and 
$$g_{\mathrm{LN}}(t | \mu, \sigma)= f_{\mathrm{LN}}(0)=\frac{1}{t_0 \sigma \sqrt{2 \pi}}  \exp (-\frac{\sigma^2}{2})$$
We consider  the Cole-Cole
\begin{align}
g_{\mathrm{RQ}}(t|t_0,\beta) &= \frac{1}{2\pi t}\frac{\sin \beta \pi}{\cosh\left(\beta \ln \left|\frac{t}{t_0}\right|\right) + \cos \beta \pi}  \\
f_{\mathrm{RQ}}(s)& =  \frac{1}{2\pi t_0}\frac{\sin \beta \pi \exp(-|s|)}{\cosh\left(\beta  s\right) + \cos \beta \pi}  \label{cc}
\end{align}
Suppose that the center points $t_0$ are the same in each case. In \eqref{cc} when $s=0$
\begin{align}
f_{\mathrm{RQ}}(0) & = \frac{1}{2\pi t_0}\frac{\sin \beta \pi }{1+ \cos \beta \pi}  \\
f_{\mathrm{RQ}}(0)=f_{\mathrm{LN}}(0)&= \frac{1}{2\pi t_0}\frac{\sin \beta \pi }{1+ \cos \beta \pi}  =\frac{1}{t_0 \sigma \sqrt{2 \pi}}  \exp (-\frac{\sigma^2}{2})\\
 \frac{\sin \beta \pi }{1+ \cos \beta \pi}  & = \frac{\sqrt{2 \pi}}{  \sigma}  \exp (-\frac{\sigma^2}{2})\\
 \tan( \frac{ \beta \pi }{2})&= \frac{\sqrt{2 \pi}}{  \sigma}  \exp (-\frac{\sigma^2}{2})\\
\frac{ \beta \pi }{2}&=  \arctan\left( \frac{\sqrt{2 \pi}}{  \sigma}  \exp (-\frac{\sigma^2}{2})\right) \\
\beta & =\frac{2}{\pi}  \arctan\left( \frac{\sqrt{2 \pi}}{  \sigma}  \exp (-\frac{\sigma^2}{2})\right) 
\end{align}

The parameters for the fitting were chosen to create matching DRTs as given above. The results here expand on the paper in that more noise levels are given. 
The data are initialized for the LN fitting with $t_0$ from $1/\omega_0(peak)$, see Section~\ref{peaks}, and  $\sigma_0=.69$, with $scale=1$.   The bounds prescribed are
$0<t_0<100$, $.1<\sigma<1$ and $0<scale<1.1$.
For the RQ fitting the equivalent information is 
$t_0=1/\omega_0(peak)$, $\beta_0=.8$, and $scale=1$, with bounds 
$0<t_0<100$, $.1<\beta<1$ and $0<scale<1.1$.
\clearpage
\begin{table}[h!]\begin{center}\begin{tabular}{|*{6}{c|}}\hline
&$-6.00$&$-5.12$&$-4.25$&$-3.38$&$-2.50$\\
\hline$0.72$&$0.86$($4e-07$)&$0.86$($3e-06$)&$0.86$($2e-05$)&$0.86$($2e-04$)&$0.86$($1e-03$)\\
\hline$0.10$&$0.20$($2e-07$)&$0.20$($1e-06$)&$0.20$($1e-05$)&$0.20$($7e-05$)&$0.20$($6e-04$)\\
\hline$1.00$&$1.01$($3e-07$)&$1.01$($2e-06$)&$1.01$($1e-05$)&$1.01$($1e-04$)&$1.01$($8e-04$)\\
\hline$0.83$&$0.83$($2e-06$)&$0.83$($1e-05$)&$0.83$($1e-04$)&$0.83$($7e-04$)&$0.83$($5e-03$)\\
\hline$0.10$&$0.10$($3e-07$)&$0.10$($2e-06$)&$0.10$($2e-05$)&$0.10$($1e-04$)&$0.10$($9e-04$)\\
\hline$1.00$&$1.00$($2e-07$)&$1.00$($2e-06$)&$1.00$($1e-05$)&$1.00$($1e-04$)&$1.00$($8e-04$)\\
\hline
\end{tabular}\caption{LN.tex\label{LN.tex}}\end{center}\end{table}

\begin{table}[h!]\begin{center}\begin{tabular}{|*{6}{c|}}\hline
&$-6.00$&$-5.12$&$-4.25$&$-3.38$&$-2.50$\\
\hline$0.72$&$0.72$($5e-07$)&$0.72$($4e-06$)&$0.72$($3e-05$)&$0.72$($2e-04$)&$0.72$($1e-03$)\\
\hline$0.10$&$0.10$($1e-07$)&$0.10$($8e-07$)&$0.10$($6e-06$)&$0.10$($5e-05$)&$0.10$($4e-04$)\\
\hline$1.00$&$1.00$($3e-07$)&$1.00$($2e-06$)&$1.00$($1e-05$)&$1.00$($1e-04$)&$1.00$($8e-04$)\\
\hline$0.83$&$1.00$($8e-16$)&$1.00$($4e-15$)&$1.00$($3e-14$)&$1.00$($2e-11$)&$1.00$($5e-09$)\\
\hline$0.10$&$0.03$($4e-08$)&$0.03$($3e-07$)&$0.03$($2e-06$)&$0.03$($2e-05$)&$0.03$($1e-04$)\\
\hline$1.00$&$0.97$($2e-07$)&$0.97$($2e-06$)&$0.97$($1e-05$)&$0.97$($9e-05$)&$0.97$($7e-04$)\\
\hline
\end{tabular}\caption{RQ.tex\label{RQ.tex}}\end{center}\end{table}
\section{Peaks in $Z_2$}\label{peaks}
Here we use $g(t)$ to refer  to  the $t$-space function and use $g_1(t)=tg(t)$, to refer to the  $s$-space function. Using the Log-normal model, 3 simulations are used. Also, 3 simulations were used for the RQ model. The simulation parameters are shown in Table~\ref{table:modelParams}.

\begin{table}[h]
\centering
\caption{Simulation Parameters}
\begin{tabular}{|c|ccc|}
\hline
RQ simulations & $\beta$ & $t_0$ & $scale$\\ \hline 
A-RQ & $[0.8]$ & $[e^{-1.5}]$ &$[1]$ \\ \hline
B-RQ & $[0.7, 0.5]$ & $[e^{-4}, e^0]$ & $[0.5, 0.5]$\\ \hline
C-RQ & $[0.8, 0.6]$&$[e^{-1.5}, e^{-.5}]$ &$[0.5, 0.5]$ \\ \hline \hline
 Log-normal simulations& $\mu$ & $s$ & $scale$\\ \hline
A-LN &$[-3.5]$ &$[0.8]$ &$[1]$ \\ \hline
B-LN & $[-7, 1]$&$[\log(1.7),\log(1.5)]$ &$[0.7,0.3]$ \\ \hline 
C-LN & $[-5,-3.25]$ &$[\log(1.7),\log(1.5)]$ &$[0.7,0.3]$ \\ \hline
\end{tabular}
\label{table:modelParams}
\end{table}

\subsection{Correlations}
When there is one process, or multiple processes spread far enough in time, there is a correlation between the imaginary part of impedance and the time of the process. The $\omega$ values corresponding to peaks in $Z_2$, the imaginary part of impedance, are the reciprocals of the times of the processes. That is $t_i^\ast=1/\omega_i$ where $t_i^\ast$ is the time of the $i^{th}$ process and $\omega_i$ is the $i^{th}$ frequency value corresponding to a peak in $Z_2$. 

\begin{figure}[h!]
   \centering
  \subfigure[Set A-LN]{\label{fig:1LN} \includegraphics[width=3in]{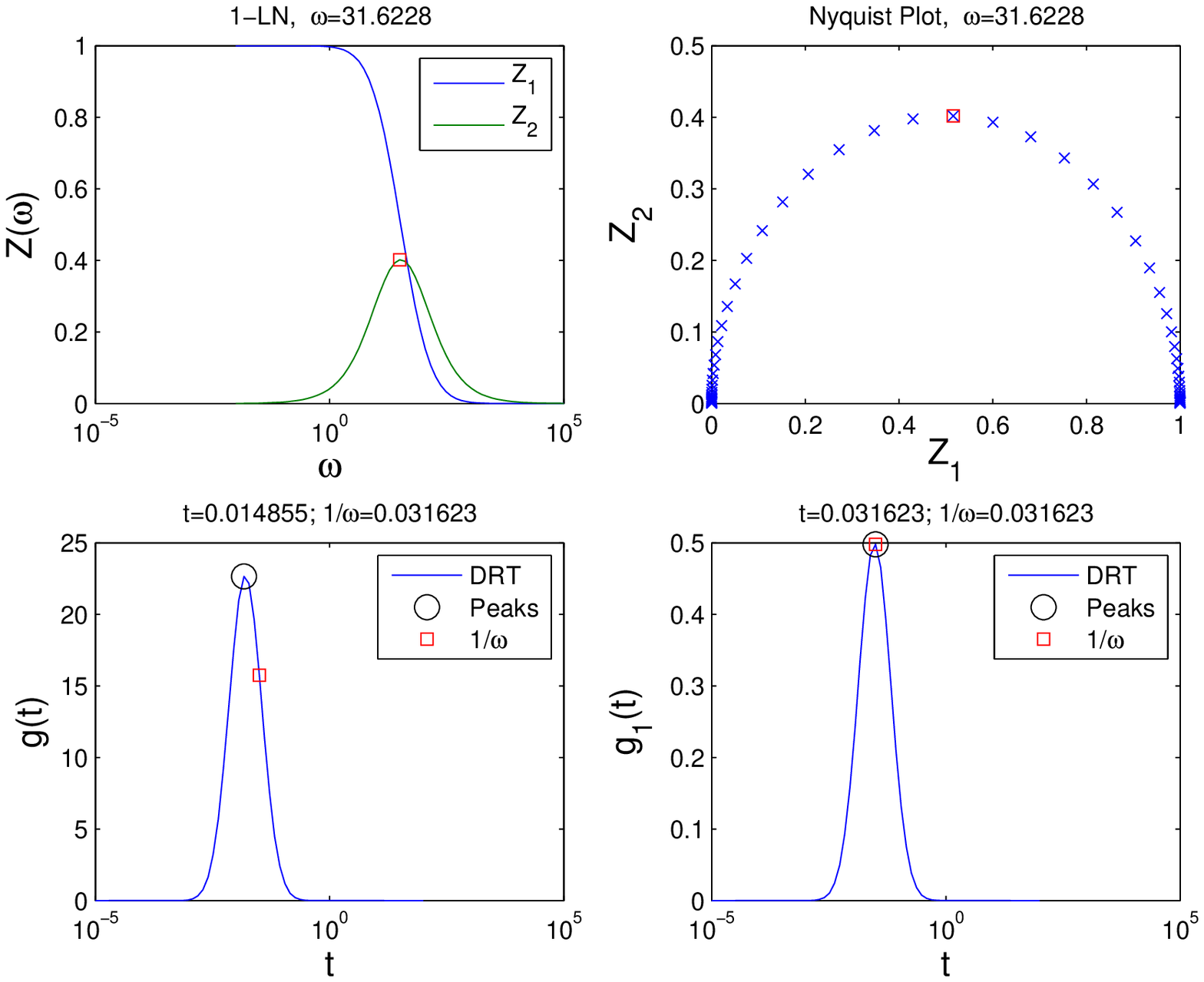}} 
    \subfigure[Set A-RQ]{\label{fig:1RQ} \includegraphics[width=3in]{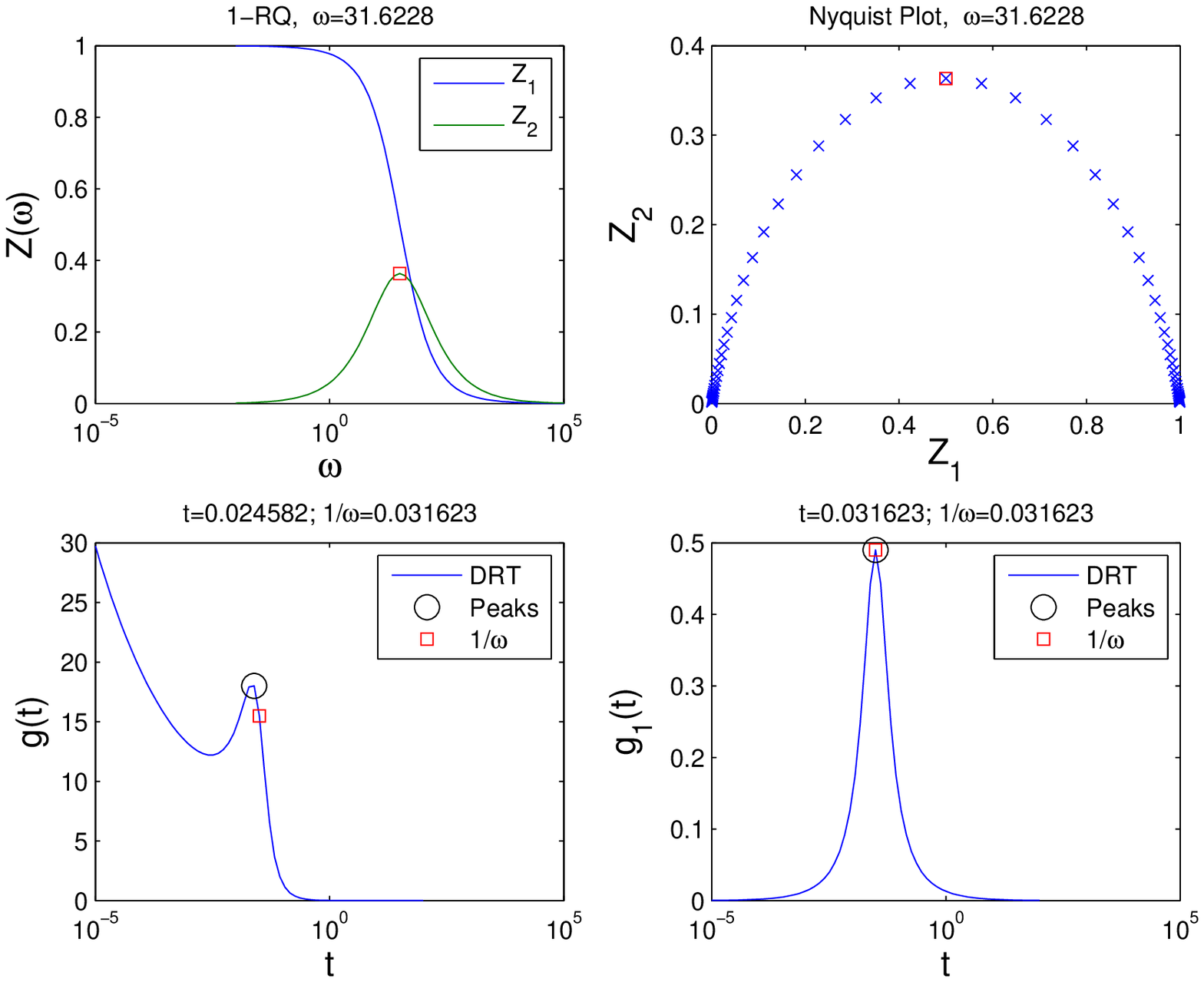}} 
  \subfigure[Set B-LN]{\label{fig:2LN} \includegraphics[width=3in]{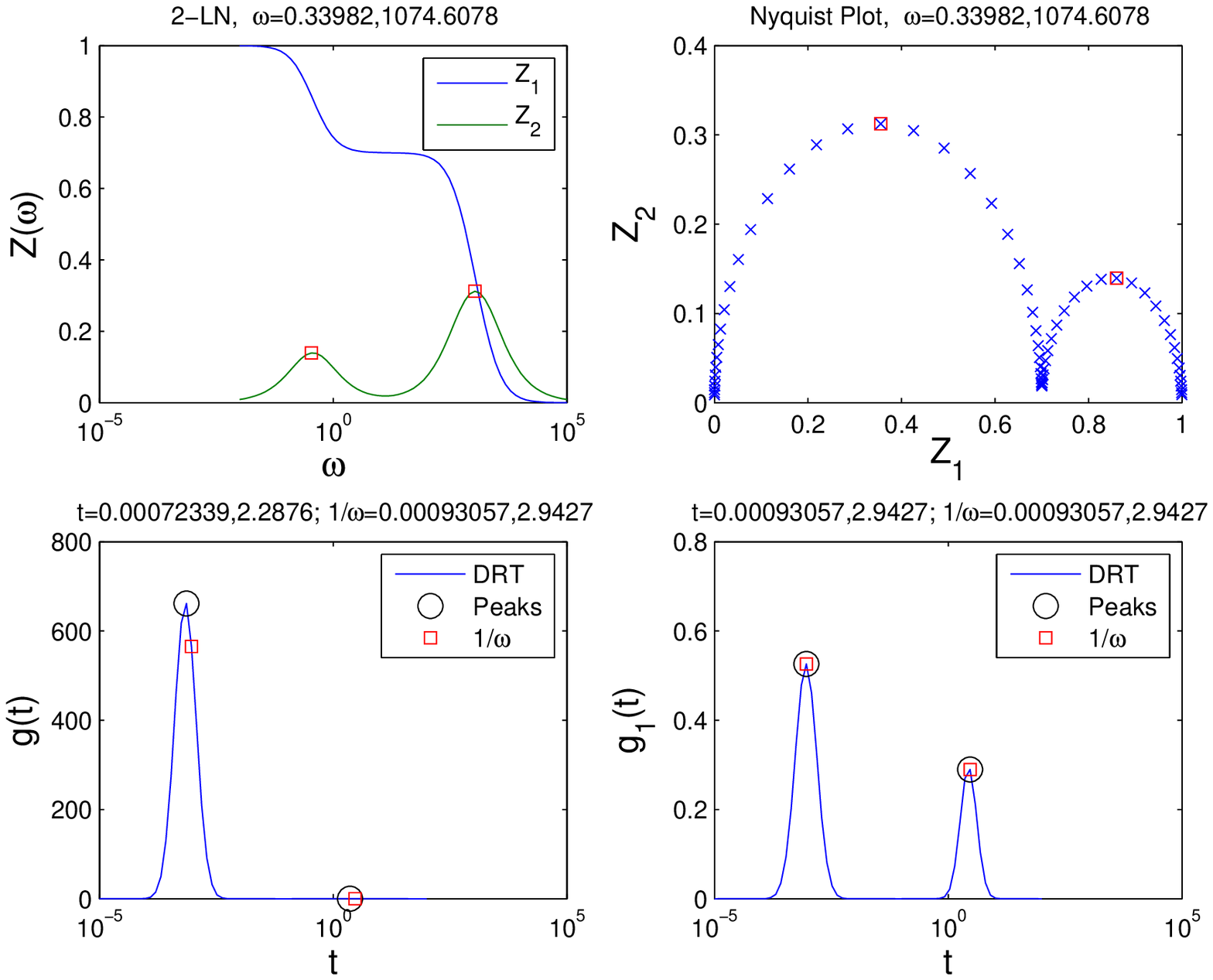}} 
      \subfigure[Set B-RQ]{\label{fig:2RQ} \includegraphics[width=3in]{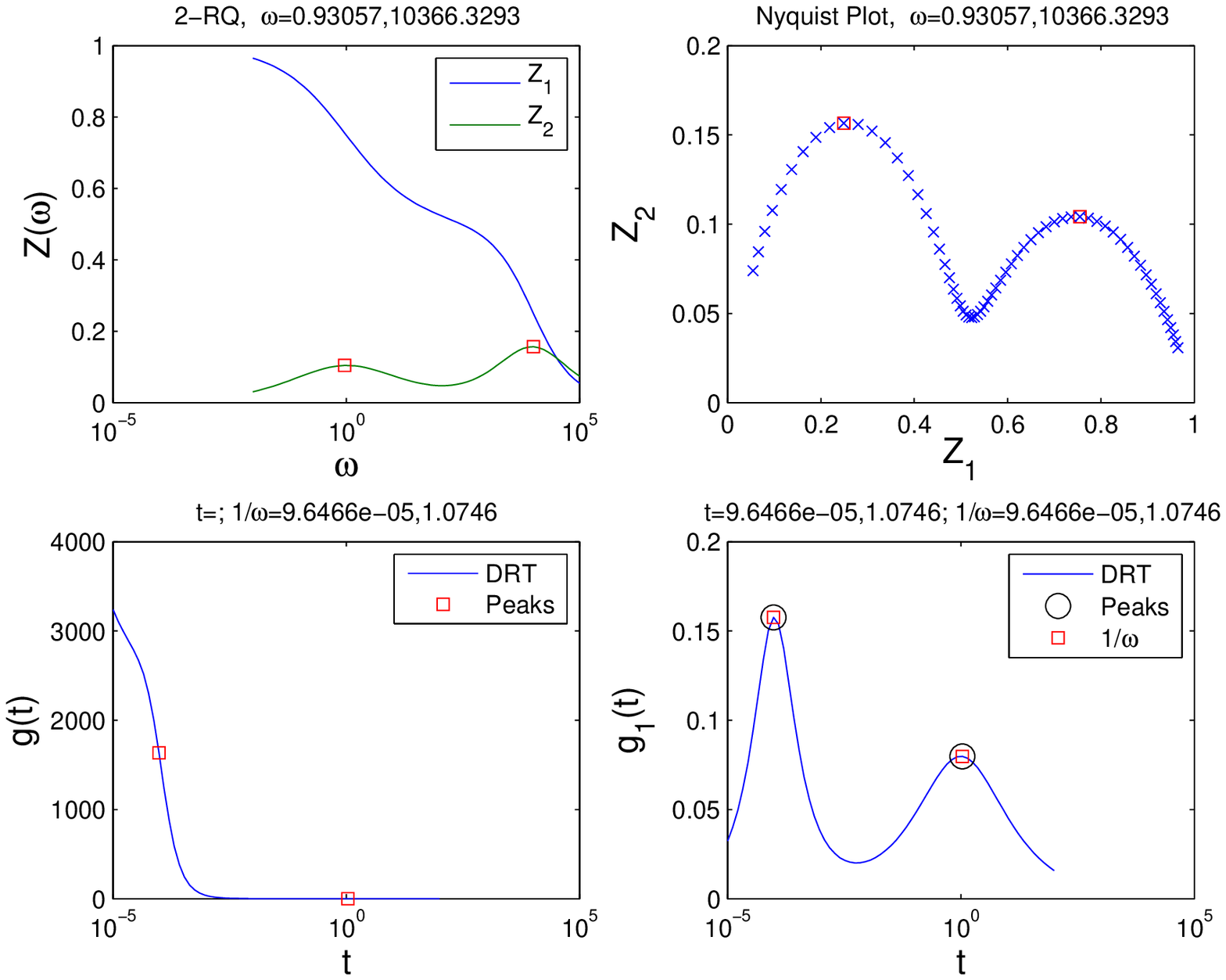}} 
    \subfigure[Set C-LN]{\label{fig:3LN} \includegraphics[width=3in]{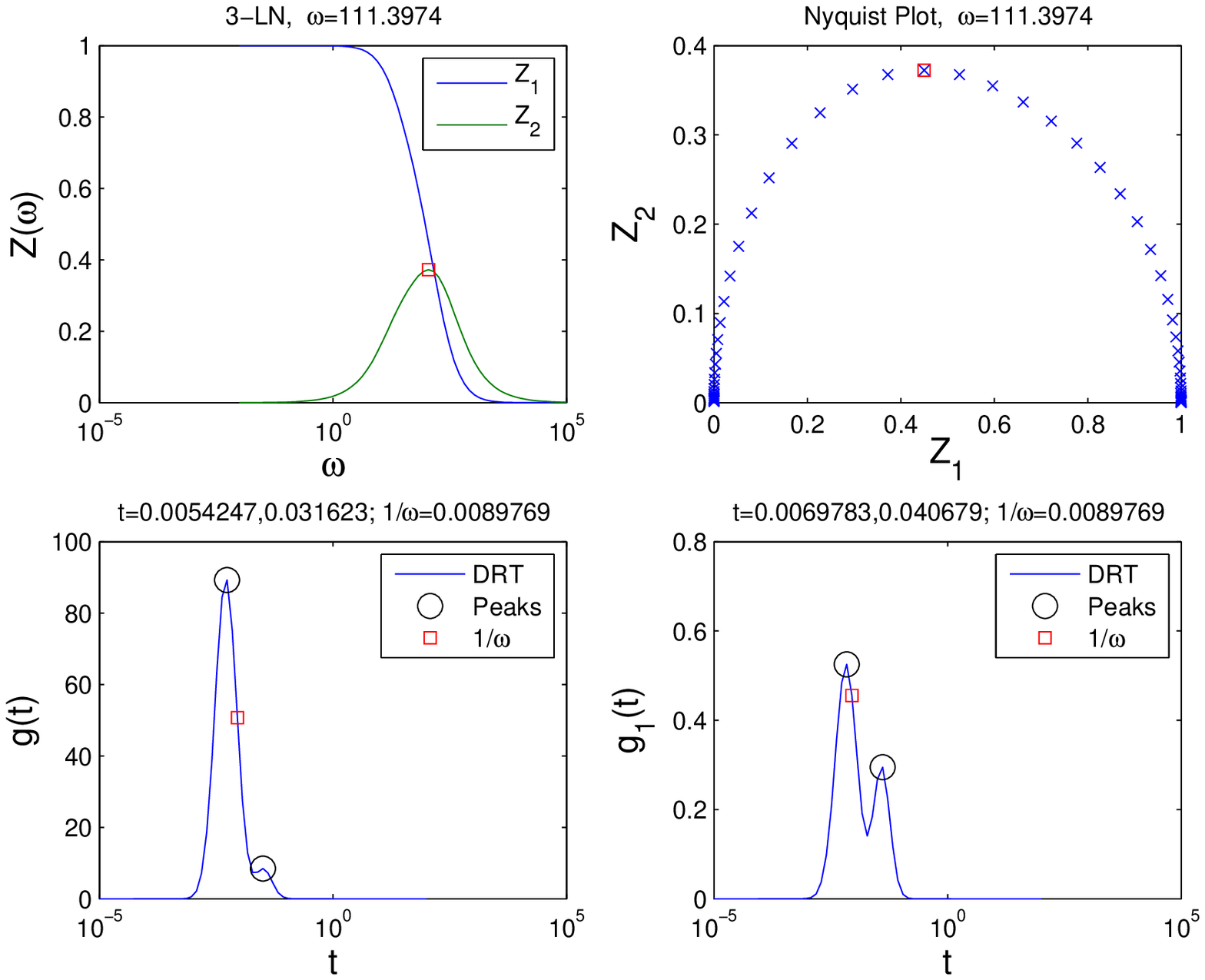}} 
  \subfigure[Set B-RQ]{\label{fig:3RQ} \includegraphics[width=3in]{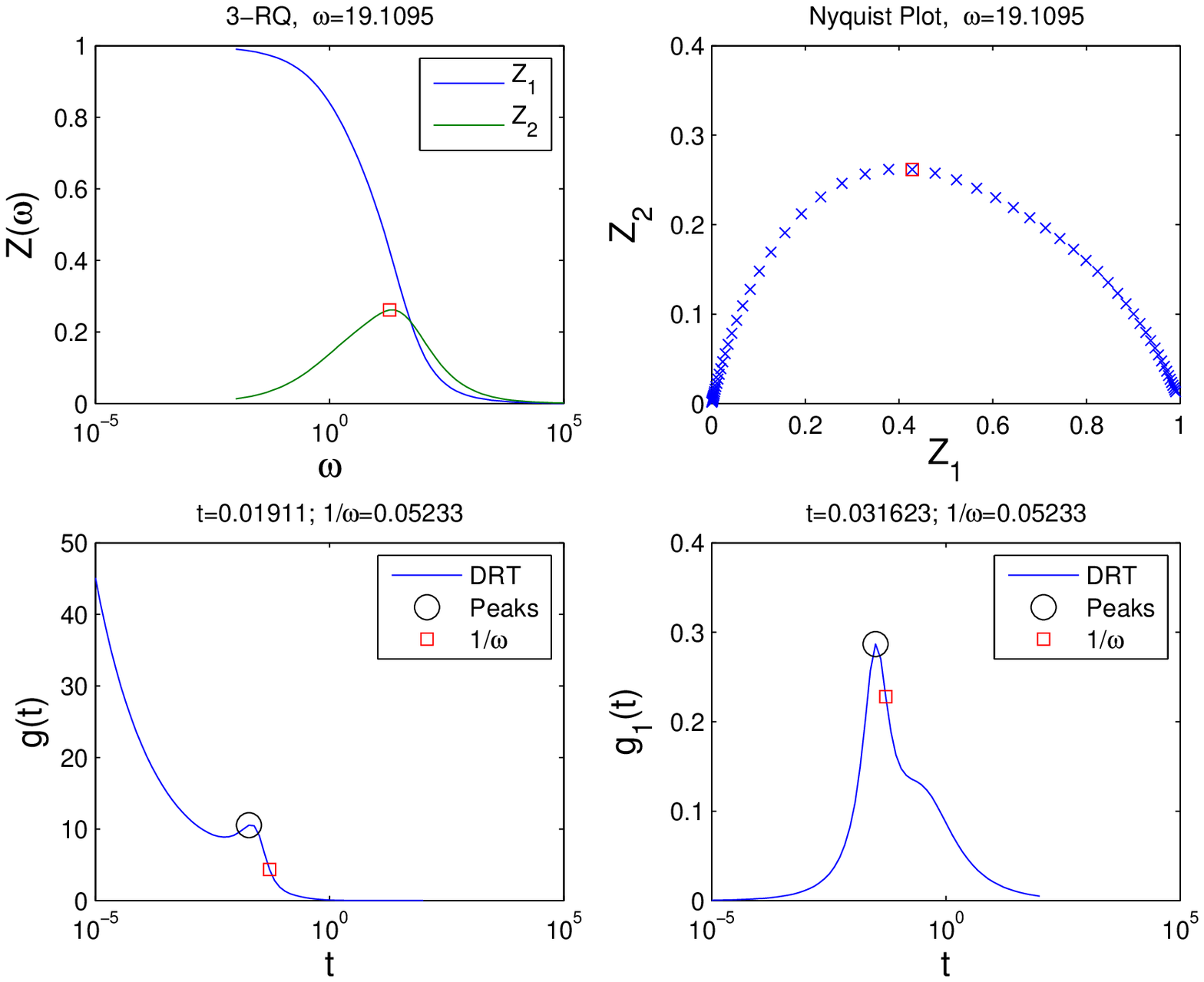}} 
\end{figure}

Analytically, $Z_1$ is monotonically decreasing under the assumption that $g(t)$, and thus $g_1(t)$, is nonnegative and $\omega$ is increasing. We have
\begin{displaymath}
Z_1(\omega )=\int_0^\infty \frac{g_1(t)}{1+\omega^2 t^2} dt
\end{displaymath}
For each $n\in \{2,3,\ldots ,N\}$ it follows that $Z_1(\omega_n)\leq Z_1(\omega_{n-1})$ since for $t\geq 0$ we have 
\begin{displaymath}
0\leq \frac{g_1(t)}{1+\omega_{n-1}^2 t^2} \leq \frac{g_1(t)}{1+\omega_n^2 t^2}
\end{displaymath}
and thus
\begin{displaymath}
Z_1(\omega_n) = \int_0^\infty \frac{g_1(t)}{1+\omega_n^2 t^2} dt \leq \int_0^\infty \frac{g_1(t)}{1+\omega_{n-1}^2 t^2} dt = Z_1(\omega_{n-1})
\end{displaymath} 
Therefore, since $Z_1$ is monotonically decreasing and peaks in $Z_2$ correspond to peaks in $g_1(t)$ if the processes are spread apart far enough, it follows that if the processes are spread far enough apart the peaks in $g_1(t)$ will be equal to the reciprocal of the $\omega$ values corresponding to peaks in the Nyquist Plot. Simply stated, $\omega$ values for peaks in $Z_2$ are the same values for peaks in the Nyquist plot, and the reciprocal of these $\omega$ values are the time points in which the processes of $g_1(t)$ peak. This is shown in Figures~\ref{fig:1LN},\ref{fig:2LN},\ref{fig:1RQ},\ref{fig:2RQ}.

 In the case where there are multiple processes but that are not spread far enough apart, $Z_2$ has one peak rather than two. Due to this, the Nyquist Plot will also only have one peak and the predicted $t$ value for process peak will lie somewhere between the two processes. This is shown in Figure~\ref{fig:3LN} and Figure~\ref{fig:3RQ}.  
\FloatBarrier
\section{Right Preconditioning}
Consider the quadrature rule for the function $g(t)=g_1(t)/t$ 
\[\int_0^\infty h(\omega,t) \frac{g_1(t)}{t}\,dt \approx \sum_{i=1}^{N-1} \Delta t_i \frac{h(\omega,t_i) \frac{g_1(t_i)}{t_i}+h(\omega,t_{i+1})\frac{g_1(t_{i+1})}{t_{i+1}}}{2} = \sum_{i=1}^{N} w_i h(\omega,t_i)\frac{g_1(t_i)}{t_i}\]
where 
\[w_i = \begin{cases}
		\Delta t_1 / 2 &\quad i = 1 \\
		(\Delta t_{i-1} + \Delta t_i)/2 &\quad 2\leq i \leq N-1 \\
		\Delta t_{N-1} / 2 &\quad i = N
		\end{cases}.
		\]
Dividing each $w_i$ by $t_i$ to move the factor of $1/t$ away from $g_1(t)$ gives 
\[w_i/t_i = \begin{cases} 
		\frac{t_2-t_1}{2t_1} &\quad i = 1\\
		\frac{t_{i+1}-t_{i-1}}{2t_i} &\quad 2 \leq i \leq N-1 \\
		\frac{t_N-t_{N-1}}{2t_{N}} &\quad i = N
		\end{cases}.
		\]
Now for the logarithmic spacing (base 10) for the $t$ we have $t_i=t_{i-1}10^{\Delta t}$ hence
\[w_i/t_i = \begin{cases} 
		\frac{10^{\Delta t}-1}{2} = \frac{\sinh(\ln (10) \Delta t )}{1+10^{-\Delta t}}&\quad i = 1\\
		\frac{10^{\Delta t} - 10^{-\Delta t}}{2} =\sinh(\ln (10) \Delta t) &\quad 2 \leq i \leq N-1 \\
		\frac{1-10^{-\Delta t}}{2} = \frac{\sinh(\ln (10) \Delta t )}{1+10^{\Delta t}}&\quad  i = N
		\end{cases}.
		\]
Now we want $\exp(s_i)=t_i$ so that the sampling matches from $s$ to $t$. Then $s_{i+1}-s_i = \ln t_{i+1} - \ln t_{i} = \ln (10)( \log(t_{i+1})-\log(t_i))= \ln (10)\Delta t $, or $\Delta t = \Delta s /\ln(10)$. Thus
\[w_i/t_i =\sinh(\Delta s)  \begin{cases} 
		\frac{1}{1+e^{-\Delta s}} =a_1&\quad i = 1\\
		 1 =a_i&\quad 2 \leq i \leq N-1 \\
		\frac{1}{1+e^{\Delta s}} =a_N &\quad i = N
		\end{cases},
		\]
and we have the quadrature formula 
\begin{equation} 	 \sum_{i=1}^{N} w_i h(\omega,t_i)\frac{g_1(t_i)}{t_i}	= \sinh(\Delta s)   \sum_{i=1}^{N} a_i h(\omega, \exp(s_i)) g_1(\exp(s_i)). \end{equation}
Meanwhile, if we first perform the change of variables $s = \ln(t)$ and then do the same trapezoidal quadrature with equally spaced intervals, we have with $f(s)=g_1(\exp(s))$
\[\int_{-\infty}^\infty h_1(\omega,s) f(s) \,ds \approx  \sum_{i=1}^N v_i h_1(\omega,s_i) f(s_i), \]
where 
\[v_i = \begin{cases}
		\Delta s / 2 &\quad i = 1 \\
		\Delta s & \quad 2 \leq i \leq N-1\\
		\Delta s /2 & \quad i = N
		\end{cases}. \]
\FloatBarrier
\section{Numerical Results for NNLS Fitting}
In the tables the numbers are given as triples mean(standard deviation) and number of samples out of $100$ ($\mu (\sigma), n$) used in calculating the mean and variance. When the third number is missing, all cases in the table generated results with relative errors less than $100\%$. We briefly list the key observations.
\begin{description}
\item[LC-NCP A4] Tables~\ref{tab:NNLSderivLCA4HN} and \ref{tab:NNLSderivNCPA4HN} compare the results using matrix $A_4$ with increasing noise levels, for the two parameter choice criteria, L-Curve and NCP. As anticipated the solutions are more robust for lower noise (fewer missed samples in calculations of mean) for both NCP and LC. On the other hand, whereas the decrease in reliability is significant for the LC, the NCP starts out worse for low noise but is far more robust to increasing noise, indeed better than the LC for higher noise for both the $I$ and $L_1$ operators. For $L_2$ the NCP also drops off in robustness. Keep in mind when comparing the means and standard deviations, that when taken over a smaller set, it does mean a better result, but the reduction of samples is significant in estimating how often the method fails. Hence comparable values ($\mu (\sigma)$) with larger $n$ suggest the case with larger $n$ is more robust.  For low noise $0.1\%$ the L-Curve results are best. A clear case for one operator over another cannot be made. It is clear that the approach is more reliable for the LN fitting than the RQ fitting, probably due to the more significant truncation of the RQ processes than the LN processes as $t\rightarrow 0$.
\item[LC-NCP A3] Tables~\ref{tab:NNLSderivLCA3HN} and \ref{tab:NNLSderivA3NCPHN} compare the results using matrix $A_3$ with increasing noise levels, for the two parameter choice criteria, L-Curve and NCP. The conclusions are similar for the $A_4$ case but overall the results for higher noise are less robust for the LC but comparable for NCP. It is interesting that even though $A_3$ has significantly better conditioning than $A_4$ the resolution for $A_4$ may be beneficial.  We deduce that if $A_4$ wins completely for just one case, and other results are comparable,  that is sufficient to indicate that one should use $A_4$.
\item[LC-NCP A3] Tables~\ref{tab:NNLSderivNCPA3LN}-\ref{tab:NNLSderivA3NCPLN}  provide the results equivalent to those given for matrix $A_4$ in the paper, namely almost the same as Tables~\ref{tab:NNLSderivLCA3HN} and \ref{tab:NNLSderivA3NCPHN} but including the lower noise level $.03\%$ in place of the $5\%$ results. Note that the times are the total times for the runs over all entries in a given table. Generally the NCP is slightly cheaper to run and $A_4$ is considerably more expensive than $A_3$. However, given the problem size and small numbers of experiments the timings are not significant for the given application. From these tables the results for $1\%$ noise make it clear that the NCP is more robust than the LC. This is borne out also for the $A_4$ matrices.
\item[LC $A_3$ and $A_4$] Tables~\ref{tab:LSderivA3HN}-\ref{tab:LSderivA4HN} compare the LS solutions for matrices $A_3$ and $A_4$ for the $s-$quadrature matrices. These results demonstrate the slightly greater stability of the $A_4$ matrices. We see that the increased resolution provides results with an often reduced variance. 
\item[LC $A_4$ for LS and NNLS] Tables~\ref{tab:NNLSderivLCA4HN}-\ref{tab:LSderivA4HN}  For the comparison of NNLS and LS results it is apparent that for stable solutions at higher noise levels the overall mean errors are reduced. On the other hand the LS is more stable in generating solutions with relative errors consistently less than $100\%$. The LS algorithm is so fast that one might use an LC algorithm and if the solution appears to be unstable, the solution should be then found with NNLS.  These results complement the similar paper in the table, but include instead $5\%$ noise over the $.03\%$ noise.
\end{description}
These observations concerning the comparisons of the two matrix sizes confirms the results in the original paper, that the extra resolution of $A_4$ can be helpful.
\vspace{-.5cm}
\begin{table}[!ht]\begin{small}
\begin{center}
\caption{NNLS: Matrix $A_4$ using the L-Curve Criterion. Higher noise.}
\label{tab:NNLSderivLCA4HN}
\begin{tabular}{|c|l||c|c|c|}\hline
Simulation &Method &$0.1\%$ &$1\%$ &$5\%$\\ 
\hline
(1,RQ) &NNLS ($L=I$) &$19 ( 2.3) 99 $&$20 (3.7)65 $&$35 (3.7)66 $\\
 &NNLS ($L=L_1$) &$12 ( 0.8)    $&$23 (1.8) 77 $&$38 (3.2)68 $\\
 &NNLS ($L=L_2$) &$13 (0.5)    $&$25 ( 2.7) 91 $&$29 ( 3.8 57 $\\
\hline
(1,LN) &NNLS ($L=I$) &$8 ( 1.9)    $&$14  (3.0) 66 $&$33 (5.5)69 $\\
 &NNLS ($L=L_1$) &$4 (1.2)    $&$12 ( 2.5) 74 $&$41 ( 3.8) 67 $\\
 &NNLS ($L=L_2$) &$5 ( 0.4)    $&$12  (3.0) 93 $&$27 (4.3) 71 $\\
\hline
\hline
(2,RQ) &NNLS ($L=I$) &$10( 2.2 ) 99 $&$16 (3.5) 76 $&$27 (3.4) 66 $\\
 &NNLS ($L=L_1$) &$8 (0.8)    $&$22 ( 0.8) 86 $&$32  (2.0) 68 $\\
 &NNLS ($L=L_2$) &$9 ( 0.8)    $&$23 (0.5) 94 $&$26 (2.9) 58 $\\
\hline
(2,LN) &NNLS ($L=I$) &$7 (1.8 ) 99 $&$13 (3.4) 77 $&$27 (3.6) 68 $\\
 &NNLS ($L=L_1$) &$4  (1.0)    $&$22 ( 1.6) 85 $&$34  (2.0) 62 $\\
 &NNLS ($L=L_2$) &$5 ( 0.7)    $&$26 (0.5) 95 $&$27 (2.3 ) 58 $\\
\hline
\hline
(3,RQ) &NNLS ($L=I$) &$14 (3.1 ) 99 $&$20 (3.4)  75 $&$34 (3.8)  68 $\\
 &NNLS ($L=L_1$) &$9 (0.9)    $&$23 ( 1.5) 74 $&$37 ( 3.2) 68 $\\
 &NNLS ($L=L_2$) &$12 ( 0.7)    $&$25 (1.1) ( 89 $&$30 (4.3 ) 69 $\\
\hline
(3,LN) &NNLS ($L=I$) &$9 (1.9 ) 99 $&$15 (2.6) ( 74 $&$33 ( 4.3) 77 $\\
 &NNLS ($L=L_1$) &$5 (1.1)    $&$14 (2.1) 82 $&$42 (3.7) 75 $\\
 &NNLS ($L=L_2$) &$5 (0.4)    $&$14 (2.6) 90 $&$29 (3.1) 72 $\\ \hline
\end{tabular}
\end{center}\end{small}
\end{table}

\vspace{-.5cm}
\begin{table}[!ht]\begin{small}
\begin{center}
\caption{NNLS: Matrix $A_4$ using the NCP Criterion. Higher noise.}
\label{tab:NNLSderivNCPA4HN}
\begin{tabular}{|c|l||c|c|c|}\hline
Simulation &Method &$0.1\%$ &$1\%$ &$5\%$\\
\hline
(1,RQ) &NNLS ($L=I$) &$14  (8.3) 91 $&$23 (4.0)83 $&$39  (7.9) 80 $\\
 &NNLS ($L=L_1$) &$14 (9.1) 92 $&$24 (4.0)86 $&$42  (8.1) 80 $\\
 &NNLS ($L=L_2$) &$14 (8.5) 92 $&$25(3.9) 88 $&$31  (8.2) 60 $\\
\hline
(1,LN) &NNLS ($L=I$) &$6 (1.3) 82 $&$18 (9.8) 81 $&$36  (7.4) 78 $\\
 &NNLS ($L=L_1$) &$6 (1.1) 88 $&$16 (10.0) 83 $&$37  (8.7) 81 $\\
 &NNLS ($L=L_2$) &$6 (0.9) 88 $&$15 (10.0) 85 $&$29 (5.0)71 $\\
\hline
\hline
(2,RQ) &NNLS ($L=I$) &$9 (4.9) 88 $&$19(2.7) 85 $&$30(6.7) 79 $\\
 &NNLS ($L=L_1$) &$10 (5.2) 88 $&$22(2.7) 83 $&$33(6.5) 78 $\\
 &NNLS ($L=L_2$) &$10 (4.5) 91 $&$22(2.6) 82 $&$26(6.4) 55 $\\
\hline
(2,LN) &NNLS ($L=I$) &$6 (4.5) 84 $&$18(3.4 81 $&$31 (5.0)75 $\\
 &NNLS ($L=L_1$) &$7 (5.1) 86 $&$19(4.4 82 $&$34(5.3) 78 $\\
 &NNLS ($L=L_2$) &$7 (4.9) 87 $&$19 (5.0)83 $&$28(4.3) 63 $\\
\hline
\hline
(3,RQ) &NNLS ($L=I$) &$11  (8.1) 87 $&$22(3.5) 84 $&$38(7.5) 81 $\\
 &NNLS ($L=L_1$) &$11 (7.5) 90 $&$24(3.5) 86 $&$40(7.5) 79 $\\
 &NNLS ($L=L_2$) &$11  (7.9) 92 $&$25(3.4) 85 $&$31  (8.7) 53 $\\
\hline
(3,LN) &NNLS ($L=I$) &$7 (1.8) 86 $&$18(4.6) 84 $&$37(5.6) 78 $\\
 &NNLS ($L=L_1$) &$7 (1.4) 85 $&$18(4.7) 86 $&$40(6.5) 80 $\\
 &NNLS ($L=L_2$) &$6 (1.4) 87 )$ &$17(4.8) 85 )$ &$30(3.1) 68 $\\ \hline
\end{tabular}
\end{center}\end{small}
\end{table}

\begin{table}[h!]
\centering
\caption{Percentage relative errors for NNLS with matrix $A_4$. NCP.  Lower noise Time $ 25677$s }
\label{tab:NNLSderiv}
\begin{tabular}{|c|c||c|c|c|}
\hline
Simulation &Method &$0.1\%$ &$0.3\%$ &$1\%$\\
\hline
(1,RQ) &NNLS ($L=I$)$\,\,\,$ &$14   \,(8.3 )  \,91 $ &$17   \,(2.8 )  \,88 $ &$23   \,(4 )  \,83 $\\
 &NNLS ($L=L_1$) &$14   \,(9.1 )  \,92 $ &$17   \,(2.5 )  \,90 $ &$24   \,(4 )  \,86 $\\
 &NNLS ($L=L_2$) &$14   \,(8.5 )  \,92 $ &$17   \,(2.8 )  \,90 $ &$25   \,(3.9 )  \,88 $\\
\hline
(1,LN) &NNLS ($L=I$)$\,\,\,$ &$6   \,(1.3 )  \,82 $ &$11   \,(4.8 )  \,80 $ &$18   \,(9.8 )  \,81 $\\
 &NNLS ($L=L_1$) &$6   \,(1.1 )  \,88 $ &$10   \,(5 )  \,85 $ &$16   \,(10 )  \,83 $\\
 &NNLS ($L=L_2$) &$6   \,(0.9 )  \,88 $ &$9   \,(5.5 )  \,86 $ &$15   \,(10 )  \,85 $\\
\hline
\hline
(2,RQ) &NNLS ($L=I$)$\,\,\,$ &$9   \,(4.9 )  \,88 $ &$14   \,(8.8 )  \,84 $ &$19   \,(2.7 )  \,85 $\\
 &NNLS ($L=L_1$) &$10   \,(5.2 )  \,88 $ &$16   \,(8.1 )  \,90 $ &$22   \,(2.7 )  \,83 $\\
 &NNLS ($L=L_2$) &$10   \,(4.5 )  \,91 $ &$16   \,(7.3 )  \,90 $ &$22   \,(2.6 )  \,82 $\\
\hline
(2,LN) &NNLS ($L=I$)$\,\,\,$ &$6   \,(4.5 )  \,84 $ &$11   \,(3.1 )  \,83 $ &$18   \,(3.4 )  \,81 $\\
 &NNLS ($L=L_1$) &$7   \,(5.1 )  \,86 $ &$11   \,(3.5 )  \,84 $ &$19   \,(4.4 )  \,82 $\\
 &NNLS ($L=L_2$) &$7   \,(4.9 )  \,87 $ &$11   \,(3.3 )  \,85 $ &$19   \,(5 )  \,83 $\\
\hline
\hline
(3,RQ) &NNLS ($L=I$)$\,\,\,$ &$11   \,(8.1 )  \,87 $ &$15   \,(2 )  \,85 $ &$22   \,(3.5 )  \,84 $\\
 &NNLS ($L=L_1$) &$11   \,(7.5 )  \,90 $ &$16   \,(2.2 )  \,89 $ &$24   \,(3.5 )  \,86 $\\
 &NNLS ($L=L_2$) &$11   \,(7.9 )  \,92 $ &$17   \,(2.1 )  \,90 $ &$25   \,(3.4 )  \,85 $\\
\hline
(3,LN) &NNLS ($L=I$)$\,\,\,$ &$7   \,(1.8 )  \,86 $ &$12   \,(7.2 )  \,86 $ &$18   \,(4.6 )  \,84 $\\
 &NNLS ($L=L_1$) &$7   \,(1.4 )  \,85 $ &$11   \,(8 )  \,87 $ &$18   \,(4.7 )  \,86 $\\
 &NNLS ($L=L_2$) &$6   \,(1.4 )  \,87 $ &$10   \,(8.8 )  \,87 $ &$17   \,(4.8 )  \,85 $\\
\hline
\end{tabular}
\end{table}

\begin{table}[!ht]\begin{small}
\begin{center}
\caption{NNLS:Matrix $A_3$ using the L-Curve Criterion. Higher noise.}
\label{tab:NNLSderivLCA3HN}
\begin{tabular}{|c|l||c|c|c|}\hline
Simulation &Method &$0.1\%$ &$1\%$ &$5\%$\\
\hline
(1,RQ) &NNLS ($L=I$) &$17 (1.7) 76 $ &$20 (3.3) 71 $&$35 (4.3) 70 $\\
 &NNLS ($L=L_1$) &$13 (0.8) 87 $ &$23 (1.8) 63 $ &$38 (3.2) 64 $\\
 &NNLS ($L=L_2$) &$13 (0.6) 99 $ &$26 (3.5) 65 $&$34 (5.6) 56 $\\
\hline
(1,LN) &NNLS ($L=I$) &$9  (8.2) 81 $ &$14 (3.1) 72 $ &$34 (4.9) 62 $\\
 &NNLS ($L=L_1$) &$4 (1.2) 92 $ &$12 (2.4) 68 $ &$41 (3.4) 55 $\\
 &NNLS ($L=L_2$) &$5 (0.4) 99 $ &$12 (3.1) 67 $ &$35 (4.4) 51 $\\
\hline
\hline
(2,RQ) &NNLS ($L=I$) &$10 (2.3) 71 $ &$16 (3.1) 58  $ &$27(3.3) 62  $\\
 &NNLS ($L=L_1$) &$8 (0.8) 81  $ &$22 (1.0) 62 $&$32( 2.3) 68 $\\
 &NNLS ($L=L_2$) &$9 ( 0.7) 96   $ &$25 (0.8) 66 $&$28 (1.3) 58 $\\
\hline
(2,LN) &NNLS ($L=I$) &$7 (2.0) 76  $ &$13 (3.3) 56 $&$27 (3.7) 59 $\\
 &NNLS ($L=L_1$) &$3 (1.1) 88  $ &$22 (1.5) 61  $ &$34 (1.9) 61  $\\
 &NNLS ($L=L_2$) &$5 (0.8) 96  $ &$28 (1.2) 67 $ &$32 (0.9) 56 $\\
\hline
\hline
(3,RQ) &NNLS ($L=I$) &$14 (3.1) 79  $ &$20 (3.4) 69  $ &$34 (3.6) 64  $\\
 &NNLS ($L=L_1$) &$9 (0.9) 88  $ &$23 (1.5) 66  $ &$37 (2.8) 58  $\\
 &NNLS ($L=L_2$) &$13 (0.7) 99   $ &$27 (3.0) 67 $ &$34 (2.4) 58  $\\
\hline
(3,LN) &NNLS ($L=I$) &$9 (2.0) 72  $ &$15 (2.6) 56 $ &$34(4.3)  61  $\\
 &NNLS ($L=L_1$) &$5 (1.1) 91  $ &$14 (2.0) 58  $ &$41 (4.0) 56  $\\
 &NNLS ($L=L_2$) &$5 (0.5)     $ &$15 (2.7) 62 $ &$40 (2.4) 55  $\\ \hline
\end{tabular}
\end{center}\end{small}
\end{table}

\begin{table}[!ht]\begin{small}
\begin{center}
\caption{NNLS:Matrix $A_3$ using the NCP  Criterion. Higher noise.}
\label{tab:NNLSderivA3NCPHN}
\begin{tabular}{|c|l||c|c|c|}\hline
Simulation &Method &$0.1\%$ &$1\%$ &$5\%$\\
\hline
(1,RQ) &NNLS ($L=I$) &$16 (7.3) 88 $&$23 (4.0)83 $&$39 (8.0) 81 $\\
 &NNLS ($L=L_1$) &$15  (7.4) 90 $&$24 (4.0) 85 $&$42  (8.2) 80 $\\
 &NNLS ($L=L_2$) &$15 ( 6.9)  90 $&$25 (4.2) 86 $&$38 (7.0) 74 $\\
\hline
(1,LN) &NNLS ($L=I$) &$6 (1.2) 85 $&$18 (9.6) 81 $&$35 (7.6) 78 $\\
 &NNLS ($L=L_1$) &$6 (1.0)  88 $&$17 (10.0) 84 $&$37 (8.6) 80 $\\
 &NNLS ($L=L_2$) &$6 (0.9) 88 $&$15 (10.2) 82 $&$36 (9.0) 77 $\\
\hline
\hline
(2,RQ) &NNLS ($L=I$) &$9 (4.9) 87 $&$19 (2.7) 85 $&$30 (6.8) 81 $\\
 &NNLS ($L=L_1$) &$10 (5.3) 89 $&$22 (2.8) 84 $&$33 (6.6) 78 $\\
 &NNLS ($L=L_2$) &$10(4.7) 91 $&$23 (3.0) 88 $&$29 (5.5) 67 $\\
\hline
(2,LN) &NNLS ($L=I$) &$6 (5.5) 83 $&$18 (3.5) 80 $&$31(4.9) 74 $\\
 &NNLS ($L=L_1$) &$7 (6.7) 85 $&$19 (4.4) 82 $&$34 (5.3) 75 $\\
 &NNLS ($L=L_2$) &$7 (6.1) 86 $&$20 (4.9) 83 $&$32 (4.8) 68 $\\
\hline
\hline
(3,RQ) &NNLS ($L=I$) &$11(3.2) 87 $&$22 (3.6) 84 $&$38 (7.7) 79 $\\
 &NNLS ($L=L_1$) &$11(2.4) 90 $&$24 (3.6) 87 $&$41 (7.7) 79 $\\
 &NNLS ($L=L_2$) &$10 (1.6) 91 $&$25 (3.6) 87 $&$36 (6.8) 70 $\\
\hline
(3,LN) &NNLS ($L=I$) &$7 (1.5) 85 $&$18 (4.6) 84 $&$37 (5.6) 78 $\\
 &NNLS ($L=L_1$) &$7 (1.2) 85 $&$17 (4.7)  85 $&$40 (6.5) 80 $\\
 &NNLS ($L=L_2$) &$6 (1.0) 87 $&$17 (4.6) 85 $&$38 (4.8) 74 $\\ \hline 
\end{tabular}
\end{center}\end{small}
\end{table}

\begin{table}[!ht]\begin{small}
\begin{center}
\caption{NNLS:Matrix $A_3$ using the L-Curve Criterion. Total time $8790$s for 100 runs. }
\label{tab:NNLSderivNCPA3LN}
\begin{tabular}{|c|c||c|c|c|}
\hline
Simulation &Method &$0.1\%$ &$0.3\%$ &$1\%$\\
\hline
(1,RQ) &NNLS ($L=I$)$\,\,\,$ &$17   \,(1.7 )  \,76 $ &$17   \,(3 )  \,64 $ &$20   \,(3.3 )  \,71 $\\
 &NNLS ($L=L_1$) &$13   \,(0.8 )  \,87 $ &$15   \,(1.3 )  \,69 $ &$23   \,(1.8 )  \,63 $\\
 &NNLS ($L=L_2$) &$13   \,(0.6 )  \,99 $ &$18   \,(2.8 )  \,76 $ &$26   \,(3.5 )  \,65 $\\
\hline
(1,LN) &NNLS ($L=I$)$\,\,\,$ &$9   \,(8.2 )  \,81 $ &$10   \,(2.9 )  \,71 $ &$14   \,(3.1 )  \,72 $\\
 &NNLS ($L=L_1$) &$4   \,(1.2 )  \,92 $ &$9   \,(10 )  \,79 $ &$12   \,(2.4 )  \,68 $\\
 &NNLS ($L=L_2$) &$5   \,(0.4 )  \,99 $ &$7   \,(0.8 )  \,85 $ &$12   \,(3.1 )  \,67 $\\
\hline
\hline
(2,RQ) &NNLS ($L=I$)$\,\,\,$ &$10   \,(2.3 )  \,71 $ &$11   \,(2.8 )  \,65 $ &$16   \,(3.1 )  \,58 $\\
 &NNLS ($L=L_1$) &$8   \,(0.8 )  \,81 $ &$14   \,(1.1 )  \,64 $ &$22   \,(1 )  \,62 $\\
 &NNLS ($L=L_2$) &$9   \,(0.7 )  \,96 $ &$18   \,(2.1 )  \,76 $ &$25   \,(0.8 )  \,66 $\\
\hline
(2,LN) &NNLS ($L=I$)$\,\,\,$ &$7   \,(2 )  \,76 $ &$8   \,(2.5 )  \,77 $ &$13   \,(3.3 )  \,56 $\\
 &NNLS ($L=L_1$) &$3   \,(1.1 )  \,88 $ &$9   \,(1.2 )  \,64 $ &$22   \,(1.5 )  \,61 $\\
 &NNLS ($L=L_2$) &$5   \,(0.8 )  \,96 $ &$9   \,(1.5 )  \,72 $ &$28   \,(1.2 )  \,67 $\\
\hline
\hline
(3,RQ) &NNLS ($L=I$)$\,\,\,$ &$14   \,(3.1 )  \,79 $ &$16   \,(2.7 )  \,59 $ &$20   \,(3.4 )  \,69 $\\
 &NNLS ($L=L_1$) &$9   \,(0.9 )  \,88 $ &$15   \,(1.1 )  \,65 $ &$23   \,(1.5 )  \,66 $\\
 &NNLS ($L=L_2$) &$13   \,(0.7 )  \,99 $ &$17   \,(1.3 )  \,75 $ &$27   \,(3 )  \,67 $\\
\hline
(3,LN) &NNLS ($L=I$)$\,\,\,$ &$9   \,(2 )  \,72 $ &$12   \,(3.1 )  \,62 $ &$15   \,(2.6 )  \,56 $\\
 &NNLS ($L=L_1$) &$5   \,(1.1 )  \,91 $ &$8   \,(0.9 )  \,73 $ &$14   \,(2 )  \,58 $\\
 &NNLS ($L=L_2$) &$5   \,(0.5 )$ &$8   \,(0.9 )  \,82 $ &$15   \,(2.7 )  \,62 $\\
\hline
\end{tabular}
\end{center}\end{small}
\end{table}

\begin{table}[!ht]\begin{small}
\begin{center}
\caption{NNLS: Matrix $A_3$ using the NCP Criterion.  Total time $6767$s for 100 runs.}
\label{tab:NNLSderivA3NCPLN}
\begin{tabular}{|c|c||c|c|c|}
\hline
Simulation &Method &$0.1\%$ &$0.3\%$ &$1\%$\\
\hline
(1,RQ) &NNLS ($L=I$)$\,\,\,$ &$16   \,(7.3 )  \,88 $ &$17   \,(2.9 )  \,87 $ &$23   \,(4 )  \,83 $\\
 &NNLS ($L=L_1$) &$15   \,(7.4 )  \,90 $ &$17   \,(3 )  \,90 $ &$24   \,(4 )  \,85 $\\
 &NNLS ($L=L_2$) &$15   \,(6.9 )  \,90 $ &$17   \,(3 )  \,90 $ &$25   \,(4.2 )  \,86 $\\
\hline
(1,LN) &NNLS ($L=I$)$\,\,\,$ &$6   \,(1.2 )  \,85 $ &$12   \,(10.8 )  \,83 $ &$18   \,(9.6 )  \,81 $\\
 &NNLS ($L=L_1$) &$6   \,(1 )  \,88 $ &$10   \,(5.4 )  \,85 $ &$17   \,(10 )  \,84 $\\
 &NNLS ($L=L_2$) &$6   \,(0.9 )  \,88 $ &$10   \,(10.7 )  \,86 $ &$15   \,(10.2 )  \,82 $\\
\hline
\hline
(2,RQ) &NNLS ($L=I$)$\,\,\,$ &$9   \,(4.9 )  \,87 $ &$14   \,(8.7 )  \,84 $ &$19   \,(2.7 )  \,85 $\\
 &NNLS ($L=L_1$) &$10   \,(5.3 )  \,89 $ &$16   \,(8.3 )  \,90 $ &$22   \,(2.8 )  \,84 $\\
 &NNLS ($L=L_2$) &$10   \,(4.7 )  \,91 $ &$16   \,(2.3 )  \,87 $ &$23   \,(3 )  \,88 $\\
\hline
(2,LN) &NNLS ($L=I$)$\,\,\,$ &$6   \,(5.5 )  \,83 $ &$12   \,(8.9 )  \,83 $ &$18   \,(3.5 )  \,80 $\\
 &NNLS ($L=L_1$) &$7   \,(6.7 )  \,85 $ &$12   \,(8.9 )  \,85 $ &$19   \,(4.4 )  \,82 $\\
 &NNLS ($L=L_2$) &$7   \,(6.1 )  \,86 $ &$12   \,(9.9 )  \,85 $ &$20   \,(4.9 )  \,83 $\\
\hline
\hline
(3,RQ) &NNLS ($L=I$)$\,\,\,$ &$11   \,(3.2 )  \,87 $ &$15   \,(2.4 )  \,84 $ &$22   \,(3.6 )  \,84 $\\
 &NNLS ($L=L_1$) &$11   \,(2.4 )  \,90 $ &$16   \,(2.3 )  \,88 $ &$24   \,(3.6 )  \,87 $\\
 &NNLS ($L=L_2$) &$10   \,(1.6 )  \,91 $ &$17   \,(2.3 )  \,88 $ &$25   \,(3.6 )  \,87 $\\
\hline
(3,LN) &NNLS ($L=I$)$\,\,\,$ &$7   \,(1.5 )  \,85 $ &$12   \,(5.5 )  \,85 $ &$18   \,(4.6 )  \,84 $\\
 &NNLS ($L=L_1$) &$7   \,(1.2 )  \,85 $ &$11   \,(5.5 )  \,86 $ &$17   \,(4.7 )  \,85 $\\
 &NNLS ($L=L_2$) &$6   \,(1 )  \,87 $ &$10   \,(5.5 )  \,86 $ &$17   \,(4.6 )  \,85 $\\
\hline
\end{tabular}
 \end{center}\end{small}
\end{table}

\begin{table}[h!]
\centering
\caption{Percentage relative errors for NNLS with matrix $A_4$.  L-Curve. Time $41078$s}
\label{tab:NNLSderiv}
\begin{tabular}{|c|c||c|c|c|}
\hline
Simulation &Method &$0.1\%$ &$0.3\%$ &$1\%$\\
\hline
(1,RQ) &NNLS ($L=I$)$\,\,\,$ &$19   \,(2.3 )  \,99 $ &$17   \,(2.8 )  \,90 $ &$20   \,(3.7 )  \,65 $\\
 &NNLS ($L=L_1$) &$12   \,(0.8 )$ &$15   \,(1.3 )  \,98 $ &$23   \,(1.8 )  \,77 $\\
 &NNLS ($L=L_2$) &$13   \,(0.5 )$ &$17   \,(3 )$ &$25   \,(2.7 )  \,91 $\\
\hline
(1,LN) &NNLS ($L=I$)$\,\,\,$ &$8   \,(1.9 )$ &$10   \,(2.8 )  \,89 $ &$14   \,(3 )  \,66 $\\
 &NNLS ($L=L_1$) &$4   \,(1.2 )$ &$7   \,(1 )  \,98 $ &$12   \,(2.5 )  \,74 $\\
 &NNLS ($L=L_2$) &$5   \,(0.4 )$ &$6   \,(0.8 )$ &$12   \,(3 )  \,93 $\\
\hline
\hline
(2,RQ) &NNLS ($L=I$)$\,\,\,$ &$10   \,(2.2 )  \,99 $ &$11   \,(2.7 )  \,88 $ &$16   \,(3.5 )  \,76 $\\
 &NNLS ($L=L_1$) &$8   \,(0.8 )$ &$14   \,(1 )  \,97 $ &$22   \,(0.8 )  \,86 $\\
 &NNLS ($L=L_2$) &$9   \,(0.8 )$ &$18   \,(2.1 )$ &$23   \,(0.5 )  \,94 $\\
\hline
(2,LN) &NNLS ($L=I$)$\,\,\,$ &$7   \,(1.8 )  \,99 $ &$8   \,(2.4 )  \,91 $ &$13   \,(3.4 )  \,77 $\\
 &NNLS ($L=L_1$) &$4   \,(1 )$ &$9   \,(1.4 )  \,97 $ &$22   \,(1.6 )  \,85 $\\
 &NNLS ($L=L_2$) &$5   \,(0.7 )$ &$9   \,(1.4 )$ &$26   \,(0.5 )  \,95 $\\
\hline
\hline
(3,RQ) &NNLS ($L=I$)$\,\,\,$ &$14   \,(3.1 )  \,99 $ &$15   \,(2.8 )  \,81 $ &$20   \,(3.4 )  \,75 $\\
 &NNLS ($L=L_1$) &$9   \,(0.9 )$ &$15   \,(1.1 )  \,95 $ &$23   \,(1.5 )  \,74 $\\
 &NNLS ($L=L_2$) &$12   \,(0.7 )$ &$17   \,(1 )$ &$25   \,(1.1 )  \,89 $\\
\hline
(3,LN) &NNLS ($L=I$)$\,\,\,$ &$9   \,(1.9 )  \,99 $ &$11   \,(2.7 )  \,85 $ &$15   \,(2.6 )  \,74 $\\
 &NNLS ($L=L_1$) &$5   \,(1.1 )$ &$8   \,(1 )  \,95 $ &$14   \,(2.1 )  \,82 $\\
 &NNLS ($L=L_2$) &$5   \,(0.4 )$ &$7   \,(0.9 )$ &$14   \,(2.6 )  \,90 $\\
\hline
\end{tabular}
\end{table}

\begin{table}[!ht]\begin{small}
\begin{center}
\caption{Percentage relative errors for LS with matrix $A_3$. L-Curve.}
\label{tab:LSderivA3HN}
\begin{tabular}{|c|c||c|c|c|}
\hline
Simulation &Method &$0.1\%$ &$1\%$ &$5\%$\\
\hline
(1,RQ) &LS ($L=I$)$\,\,\,$ &$27   \,(6.8 )$ &$24   \,(4.7 )$ &$35   \,(5.2 )$\\
 &LS ($L=L_1$) &$13   \,(1.1 )$ &$21   \,(1.8 )$ &$38   \,(3.2 )  \,99 $\\
 &LS ($L=L_2$) &$40   \,(5.8 )$ &$40   \,(5.2 )$ &$45   \,(1.2 )  \,99 $\\
\hline
(1,LN) &LS ($L=I$)$\,\,\,$ &$18   \,(5.5 )$ &$20   \,(4.4 )$ &$35   \,(5.8 )$\\
 &LS ($L=L_1$) &$6   \,(1.4 )$ &$18   \,(2.7 )$ &$41   \,(3.7 )  \,99 $\\
 &LS ($L=L_2$) &$29   \,(12.9 )$ &$29   \,(12.7 )$ &$52   \,(0.8 )$\\
\hline
\hline
(2,RQ) &LS ($L=I$)$\,\,\,$ &$18   \,(5.3 )$ &$20   \,(4.1 )$ &$30   \,(4.7 )$\\
 &LS ($L=L_1$) &$7   \,(1.1 )$ &$20   \,(1.3 )$ &$31   \,(2.3 )$\\
 &LS ($L=L_2$) &$24   \,(2.4 )$ &$25   \,(1.7 )$ &$39   \,(6 )  \,96 $\\
\hline
(2,LN) &LS ($L=I$)$\,\,\,$ &$18   \,(5.5 )$ &$19   \,(4 )$ &$32   \,(5.2 )$\\
 &LS ($L=L_1$) &$5   \,(1.2 )$ &$21   \,(2 )$ &$34   \,(2.5 )$\\
 &LS ($L=L_2$) &$27   \,(4 )$ &$28   \,(3.2 )$ &$39   \,(7.2 )  \,98 $\\
\hline
\hline
(3,RQ) &LS ($L=I$)$\,\,\,$ &$20   \,(4.9 )$ &$22   \,(4.2 )$ &$34   \,(5.2 )$\\
 &LS ($L=L_1$) &$10   \,(1.4 )$ &$21   \,(1.9 )$ &$37   \,(2.8 )  \,98 $\\
 &LS ($L=L_2$) &$25   \,(2.4 )$ &$25   \,(1.7 )$ &$43   \,(5 )  \,91 $\\
\hline
(3,LN) &LS ($L=I$)$\,\,\,$ &$18   \,(5.3 )$ &$21   \,(4.2 )$ &$35   \,(6 )$\\
 &LS ($L=L_1$) &$6   \,(1.4 )$ &$19   \,(2.7 )$ &$40   \,(3.8 )  \,99 $\\
 &LS ($L=L_2$) &$25   \,(5.6 )$ &$26   \,(5.2 )$ &$46   \,(7.6 )  \,93 $\\
\hline
\end{tabular}
\end{center}\end{small}
\end{table}

\begin{table}[!ht]\begin{small}
\begin{center}
\caption{Percentage relative errors for LS with matrix $A_4$. L-Curve.}
\label{tab:LSderivA4HN}
\begin{tabular}{|c|c||c|c|c|}
\hline
Simulation &Method &$0.1\%$ &$1\%$ &$5\%$\\
\hline
(1,RQ) &LS ($L=I$)$\,\,\,$ &$18   \,(3.2 )$ &$21   \,(3.7 )$ &$36   \,(4.2 )$\\
 &LS ($L=L_1$) &$12   \,(0.9 )$ &$24   \,(1.9 )$ &$39   \,(3 )$\\
 &LS ($L=L_2$) &$42   \,(2.8 )$ &$42   \,(0.2 )$ &$45   \,(1.3 )$\\
\hline
(1,LN) &LS ($L=I$)$\,\,\,$ &$12   \,(2.6 )$ &$18   \,(3.9 )$ &$37   \,(4.6 )$\\
 &LS ($L=L_1$) &$5   \,(1.2 )$ &$20   \,(2.2 )$ &$43   \,(2.9 )$\\
 &LS ($L=L_2$) &$44   \,(11.5 )$ &$44   \,(10.9 )$ &$52   \,(0.9 )$\\
\hline
\hline
(2,RQ) &LS ($L=I$)$\,\,\,$ &$12   \,(2.6 )$ &$18   \,(3.1 )$ &$29   \,(3.3 )$\\
 &LS ($L=L_1$) &$8   \,(0.8 )$ &$22   \,(0.9 )$ &$33   \,(1.7 )$\\
 &LS ($L=L_2$) &$26   \,(1.7 )$ &$26   \,(0.4 )$ &$41   \,(4.2 )$\\
\hline
(2,LN) &LS ($L=I$)$\,\,\,$ &$12   \,(2.7 )$ &$17   \,(3.3 )$ &$31   \,(3.5 )$\\
 &LS ($L=L_1$) &$5   \,(1 )$ &$24   \,(1.3 )$ &$36   \,(2.5 )$\\
 &LS ($L=L_2$) &$30   \,(2.6 )$ &$30   \,(0.8 )$ &$43   \,(7.8 )$\\
\hline
\hline
(3,RQ) &LS ($L=I$)$\,\,\,$ &$14   \,(2.4 )$ &$20   \,(3.4 )$ &$34   \,(3.7 )$\\
 &LS ($L=L_1$) &$9   \,(0.9 )$ &$23   \,(1.3 )$ &$38   \,(2.4 )$\\
 &LS ($L=L_2$) &$26   \,(1.5 )$ &$27   \,(0.8 )$ &$45   \,(3.2 )  \,99 $\\
\hline
(3,LN) &LS ($L=I$)$\,\,\,$ &$12   \,(2.5 )$ &$18   \,(3.6 )$ &$36   \,(4.2 )$\\
 &LS ($L=L_1$) &$5   \,(1.2 )$ &$23   \,(2.1 )$ &$42   \,(2.6 )$\\
 &LS ($L=L_2$) &$30   \,(2.5 )$ &$30   \,(0.8 )$ &$50   \,(5.8 )$\\
\hline
\end{tabular}
\end{center}\end{small}
\end{table}

\begin{table}[h!]
\centering
\caption{Percentage relative errors for LS with matrix $A_4$. L-Curve.}
\label{tab:LSderiv}
\begin{tabular}{|c|c||c|c|c|}
\hline
Simulation &Method &$0.1\%$ &$0.3\%$ &$1\%$\\
\hline
(1,RQ) &LS ($L=I$)$\,\,\,$ &$18   \,(3.2 )$ &$17   \,(2.8 )$ &$21   \,(3.7 )$\\
 &LS ($L=L_1$) &$12   \,(0.9 )$ &$15   \,(1.3 )$ &$24   \,(1.9 )$\\
 &LS ($L=L_2$) &$27   \,(12 )$ &$27   \,(11.9 )$ &$42   \,(0.2 )$\\
\hline
(1,LN) &LS ($L=I$)$\,\,\,$ &$12   \,(2.6 )$ &$13   \,(2.5 )$ &$18   \,(3.9 )$\\
 &LS ($L=L_1$) &$5   \,(1.2 )$ &$10   \,(1.6 )$ &$20   \,(2.2 )$\\
 &LS ($L=L_2$) &$12   \,(1.7 )$ &$12   \,(1.6 )$ &$44   \,(10.9 )$\\
\hline
\hline
(2,RQ) &LS ($L=I$)$\,\,\,$ &$12   \,(2.6 )$ &$13   \,(2.5 )$ &$18   \,(3.1 )$\\
 &LS ($L=L_1$) &$8   \,(0.8 )$ &$14   \,(1 )$ &$22   \,(0.9 )$\\
 &LS ($L=L_2$) &$18   \,(1.9 )$ &$18   \,(1.6 )$ &$26   \,(0.4 )$\\
\hline
(2,LN) &LS ($L=I$)$\,\,\,$ &$12   \,(2.7 )$ &$13   \,(2.5 )$ &$17   \,(3.3 )$\\
 &LS ($L=L_1$) &$5   \,(1 )$ &$11   \,(1.5 )$ &$24   \,(1.3 )$\\
 &LS ($L=L_2$) &$14   \,(3 )$ &$14   \,(2.9 )$ &$30   \,(0.8 )$\\
\hline
\hline
(3,RQ) &LS ($L=I$)$\,\,\,$ &$14   \,(2.4 )$ &$15   \,(2.9 )$ &$20   \,(3.4 )$\\
 &LS ($L=L_1$) &$9   \,(0.9 )$ &$15   \,(1.1 )$ &$23   \,(1.3 )$\\
 &LS ($L=L_2$) &$16   \,(1 )$ &$17   \,(1 )$ &$27   \,(0.8 )$\\
\hline
(3,LN) &LS ($L=I$)$\,\,\,$ &$12   \,(2.5 )$ &$13   \,(2.7 )$ &$18   \,(3.6 )$\\
 &LS ($L=L_1$) &$5   \,(1.2 )$ &$11   \,(1.5 )$ &$23   \,(2.1 )$\\
 &LS ($L=L_2$) &$13   \,(1.4 )$ &$13   \,(1.2 )$ &$30   \,(0.8 )$\\
\hline
\end{tabular}
\end{table}

\clearpage
\section{L-curve and NCP Parameter Choice Comparisons}
For each noise realization the following information was recorded: the optimal solution obtained by the NCP and L-curve parameter choice methods, with the optimally found $\lambda_{\mathrm{NCP}}$ and $\lambda_{\mathrm{LC}}$, and the optimal solution over all $50$ choices for $\lambda$, with the respective $\lambda_{\mathrm{opt}}$, as measured with respect to the absolute error in the $s$ space.  The geometric means  of  $\lambda_{\mathrm{NCP}}$ and $\lambda_{\mathrm{LC}}$  were  calculated over all $50$ noise realizations. The absolute error for each choice of $\lambda$ was also recorded for each noise realization, and  the mean of these absolute errors taken to give an average error for a given $\lambda$ which can be visualized against $\lambda$. In the plots we thus show the average error against $\lambda$ indicated by the $\circ$ plot. On the same plot we indicate by the vertical lines the minimum $\lambda_{\mathrm{opt}}$, and the geometric means for  $\lambda_{\mathrm{NCP}}$ and $\lambda_{\mathrm{LC}}$, as the solid (red), dashed (green) and dot-dashed $\circ$ (blue) vertical lines, respectively.  For each simulation set the same procedure was performed for all smoothing norms $L$. To demonstrate the dependence of the obtained solution on the optimal parameter, an example noise realization was chosen in each case and the solutions found using the chosen optimal parameters and compared with the exact solution. These are indicated by the solid line (black), $\diamond$ (red) , $\times$ (green)  and $\circ$ (blue), for the exact, $\lambda_{\mathrm{opt}}$, $\lambda_{\mathrm{NCP}}$, and $\lambda_{\mathrm{LC}}$ solutions, respectively. 

\begin{figure}[!ht]
\centering
\subfigure[$L=I$]{\includegraphics[width=1.7in]{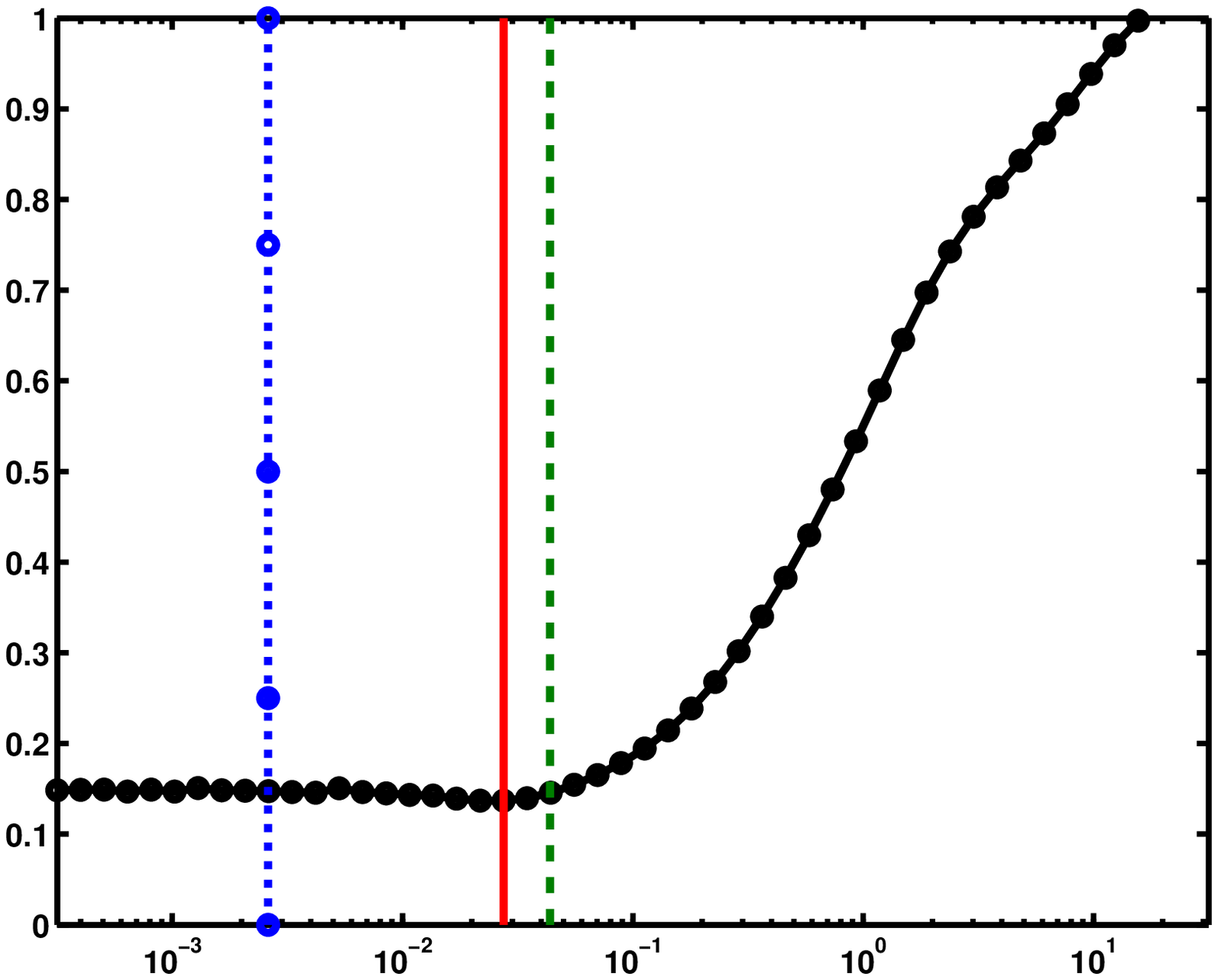}}
\subfigure[$L=L_1$]{\includegraphics[width=1.7in]{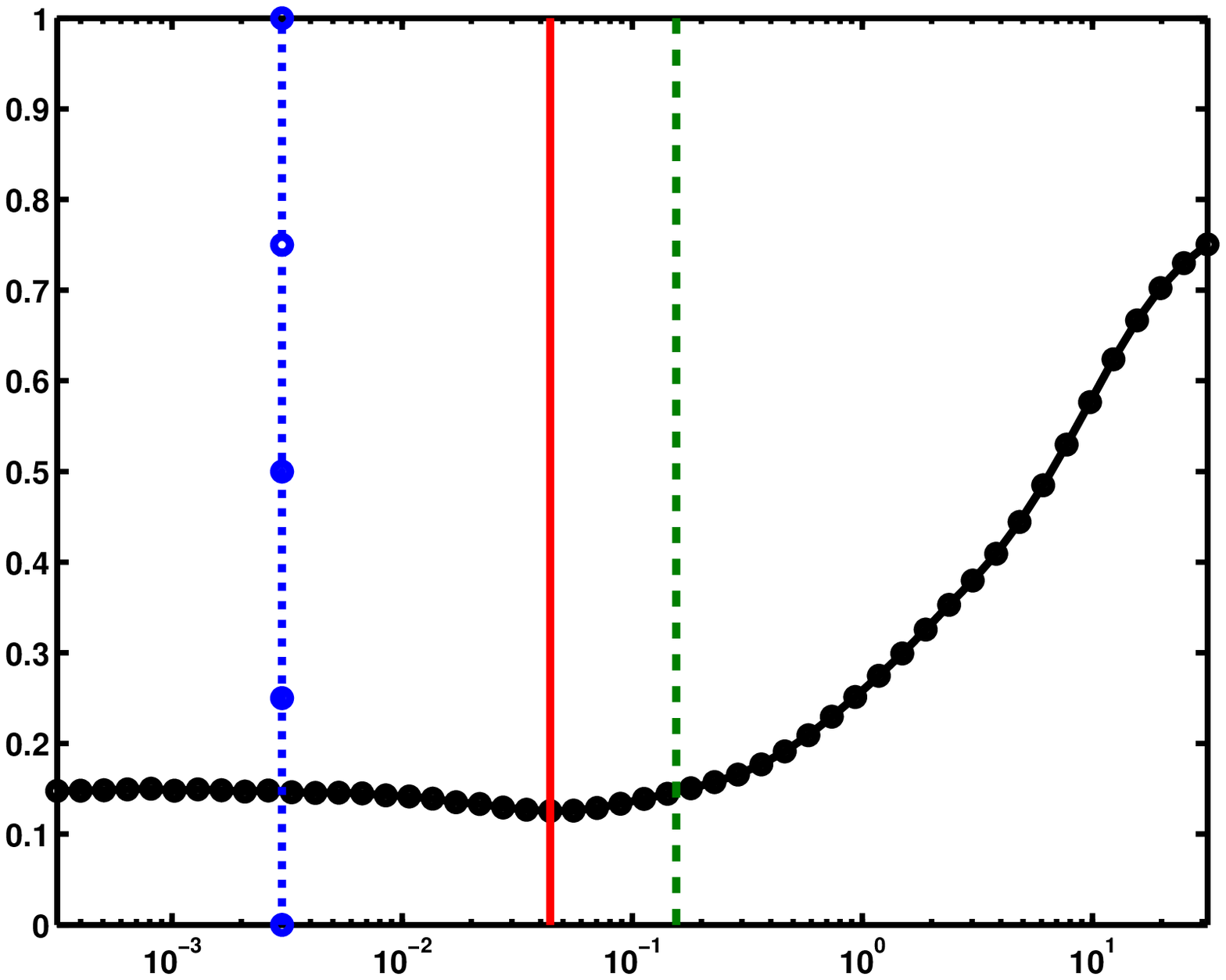}}
\subfigure[$L=L_2$]{\includegraphics[width=1.7in]{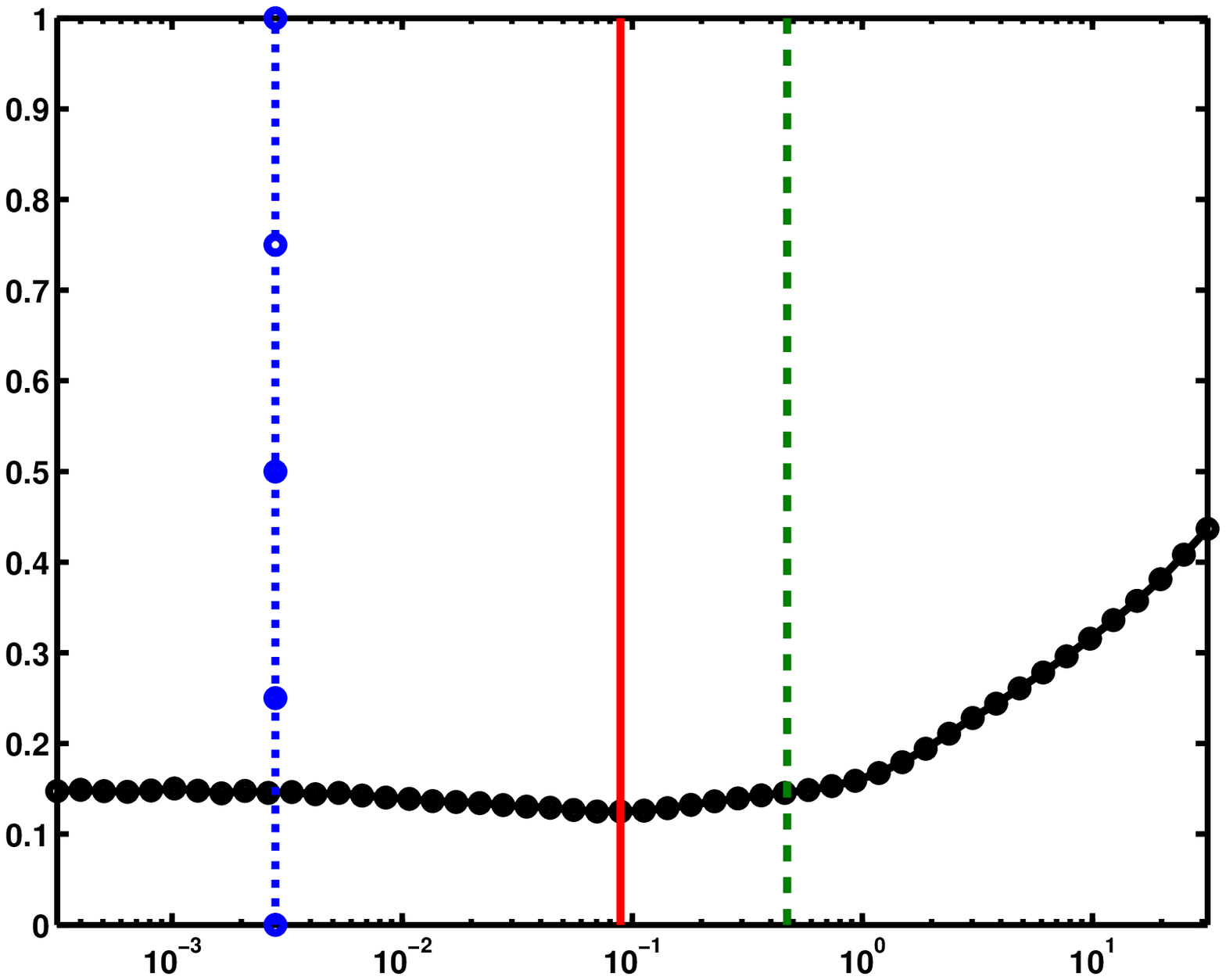}}
\subfigure[$L=I$]{\includegraphics[width=1.7in]{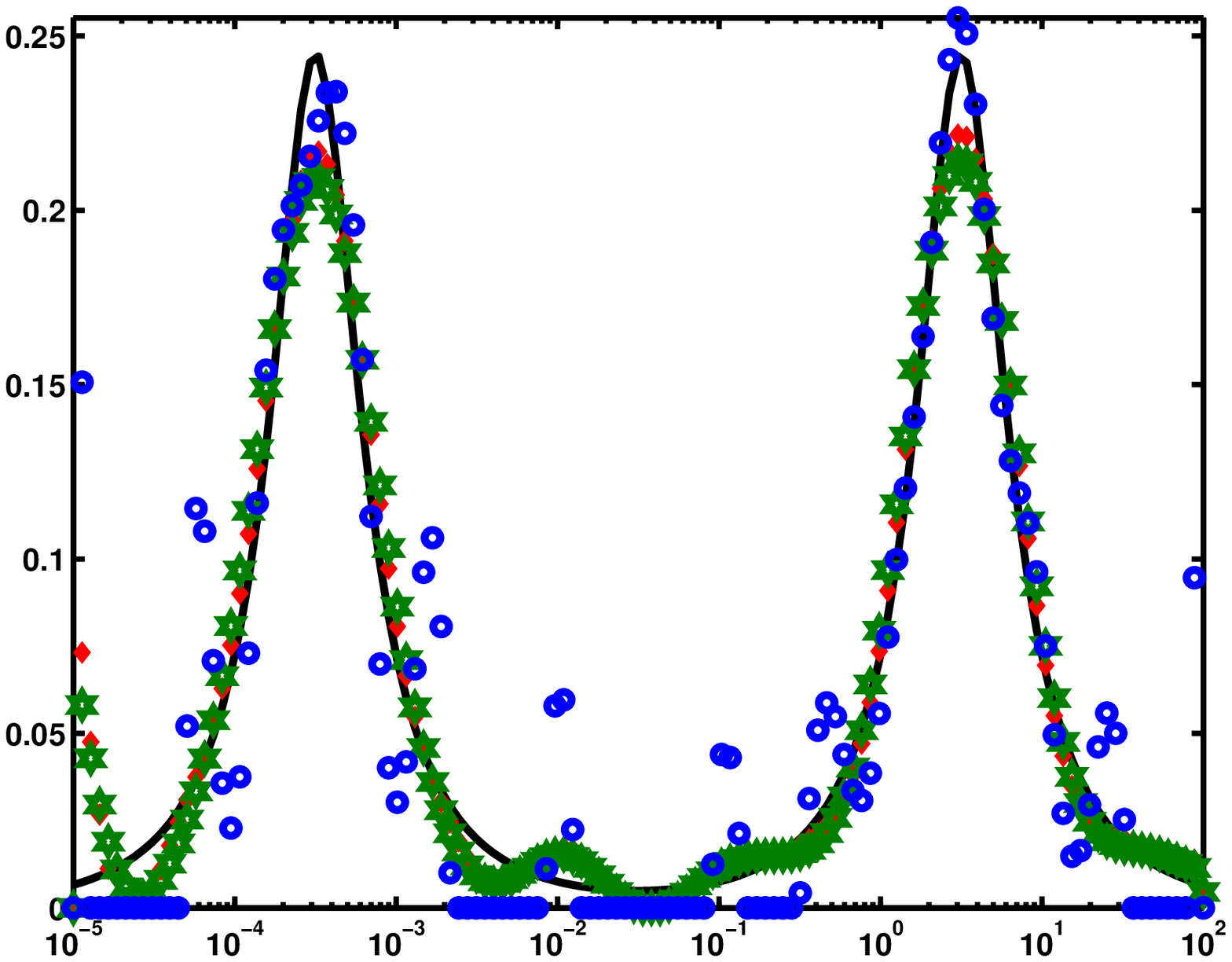}}
\subfigure[$L=L_1$]{\includegraphics[width=1.7in]{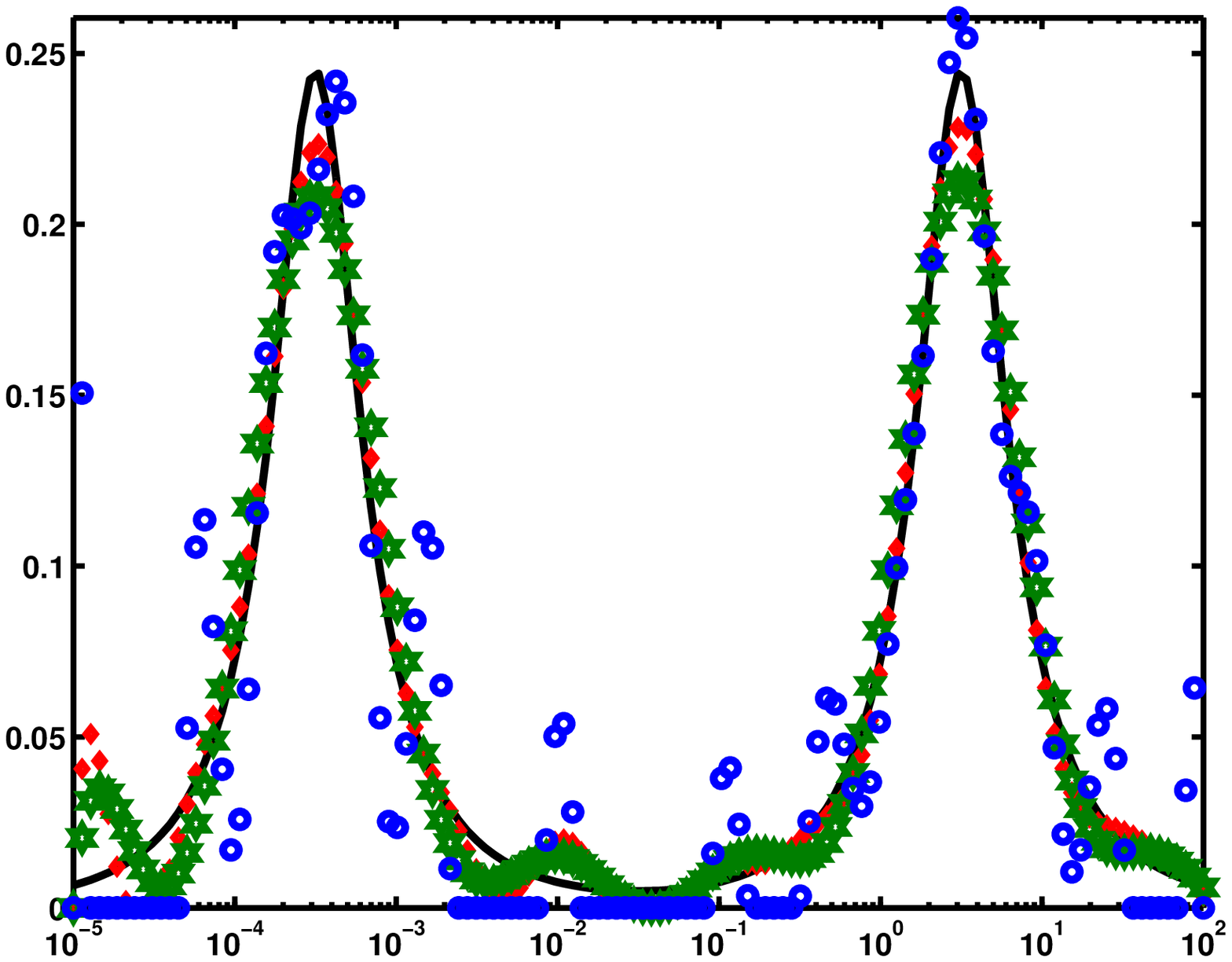}}
\subfigure[$L=L_2$]{\includegraphics[width=1.7in]{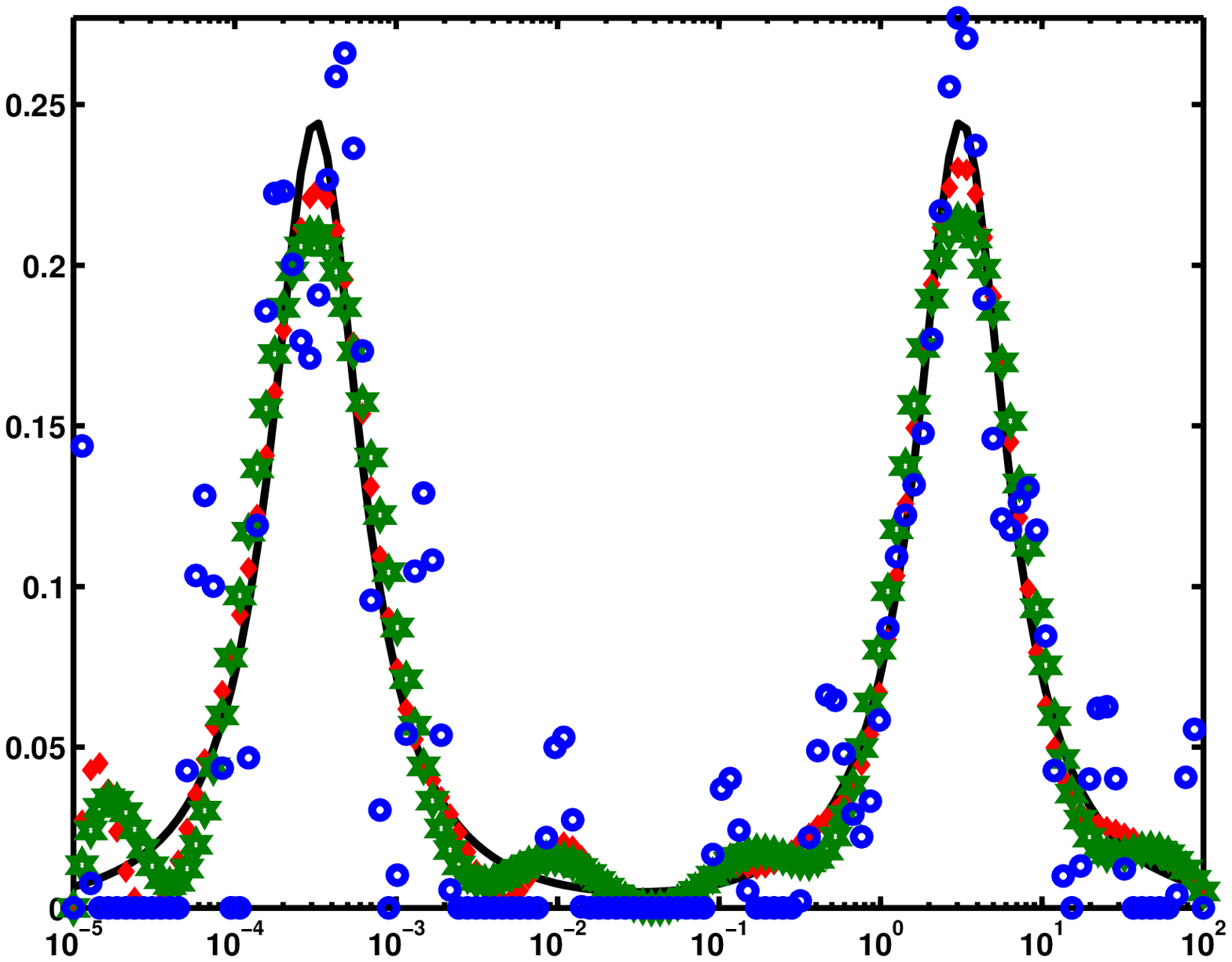}}
\caption{NNLS solutions of RQ-A matrix $A_4$. Noise level $1\%$ using the SBB algorithm}
\label{lnfig-lambdachoiceRQ1A4LNSBB}
\end{figure}

\begin{figure}[!ht]
\centering
\subfigure[$L=I$]{\includegraphics[width=1.7in]{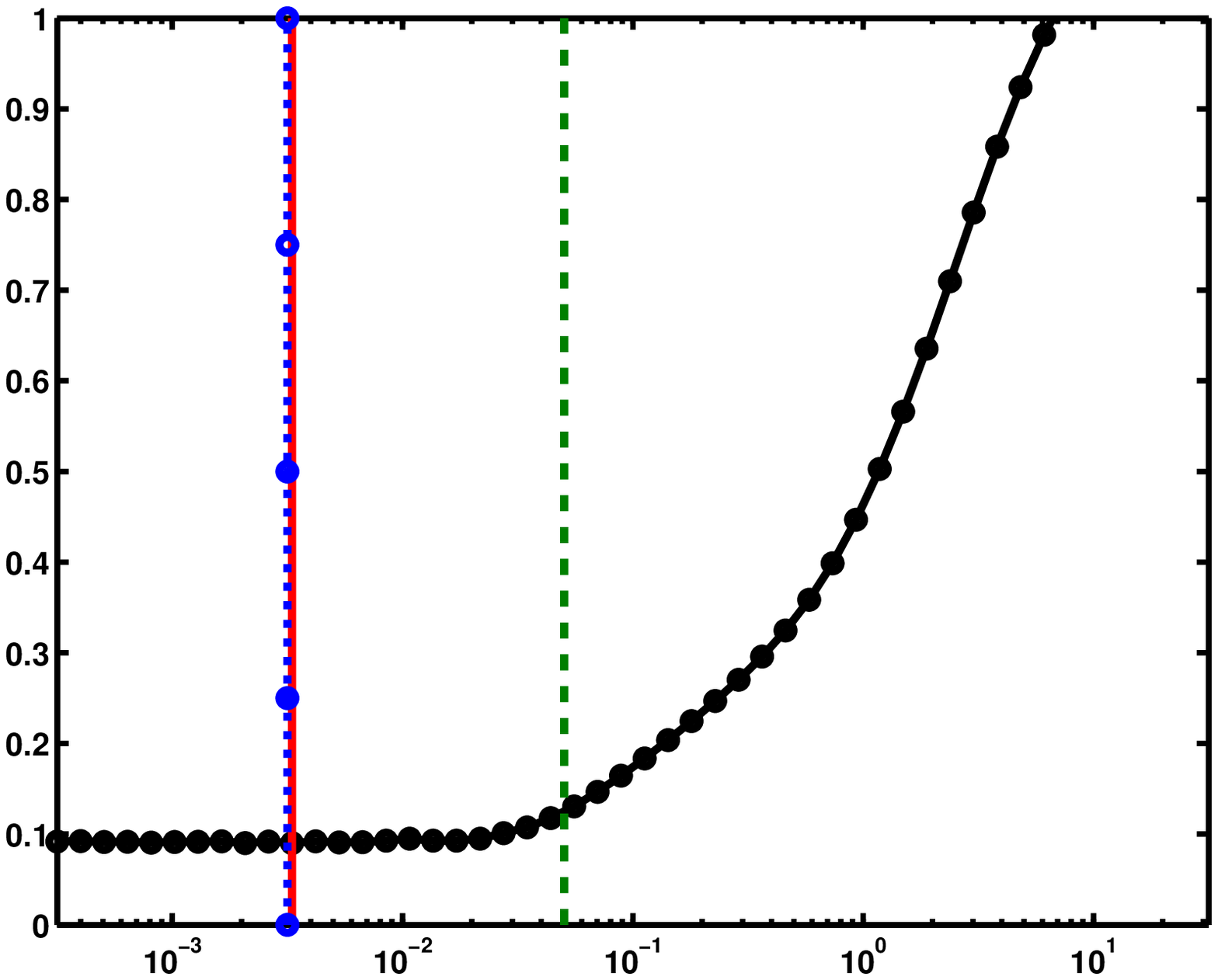}}
\subfigure[$L=L_1$]{\includegraphics[width=1.7in]{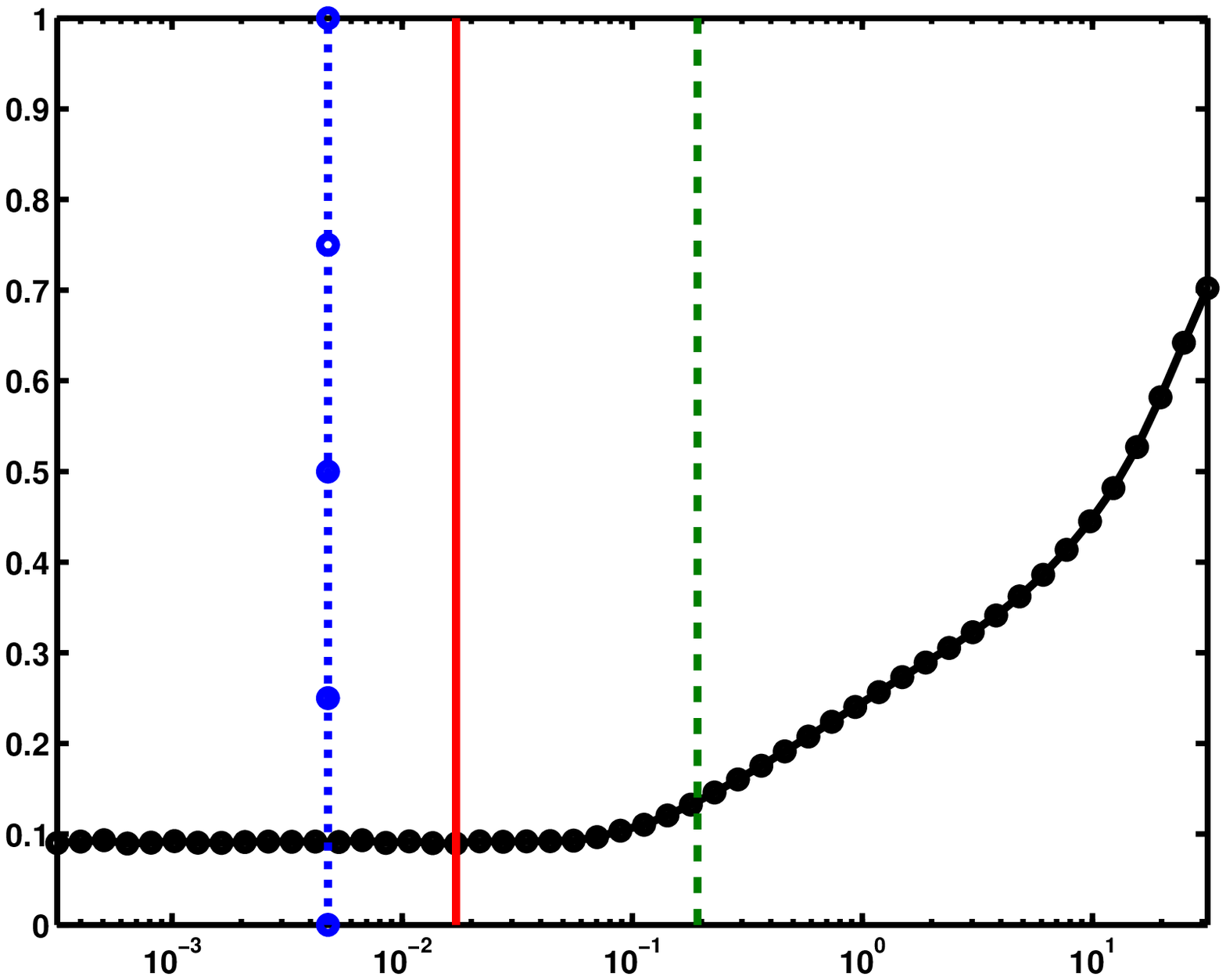}}
\subfigure[$L=L_2$]{\includegraphics[width=1.7in]{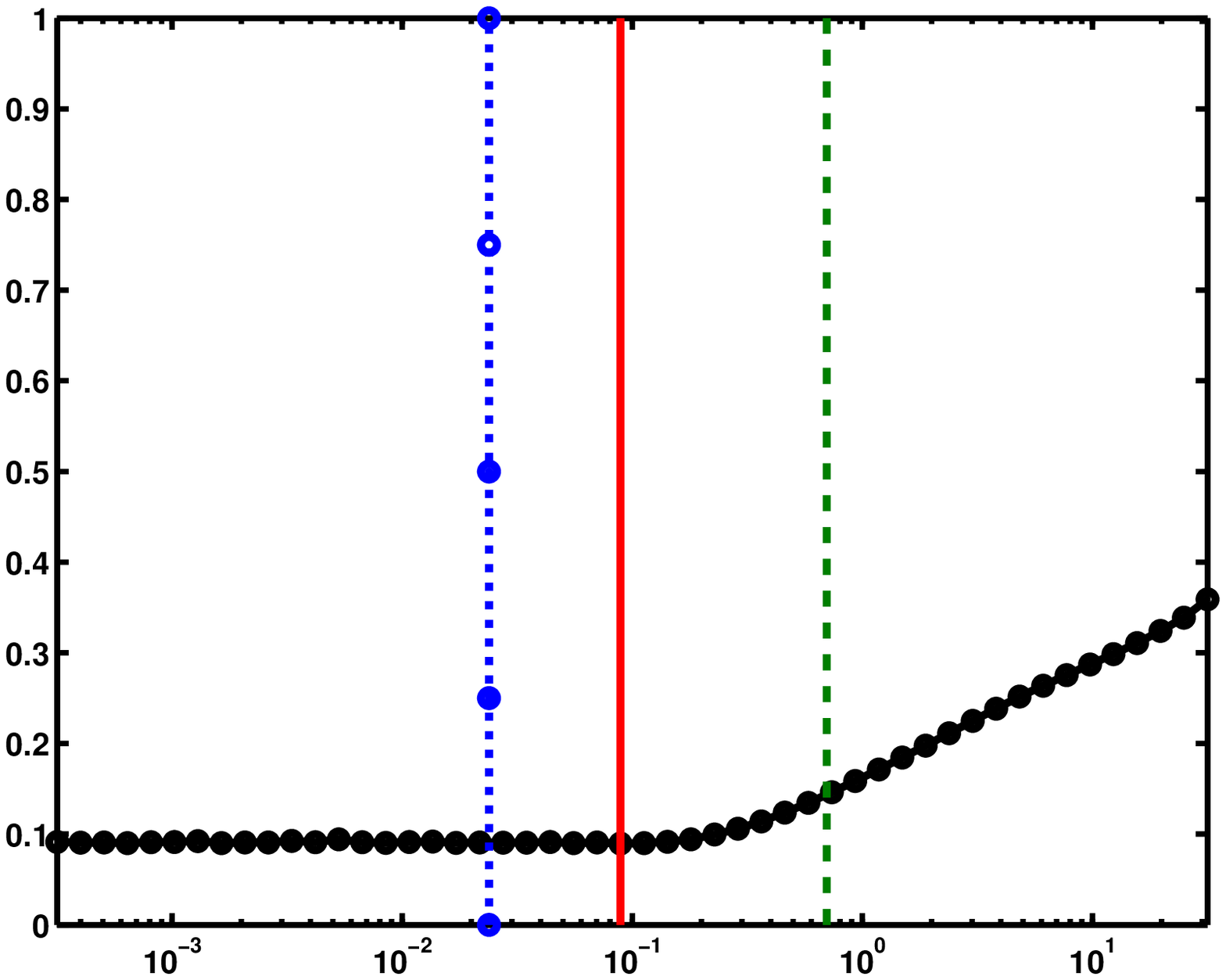}}
\subfigure[$L=I$]{\includegraphics[width=1.7in]{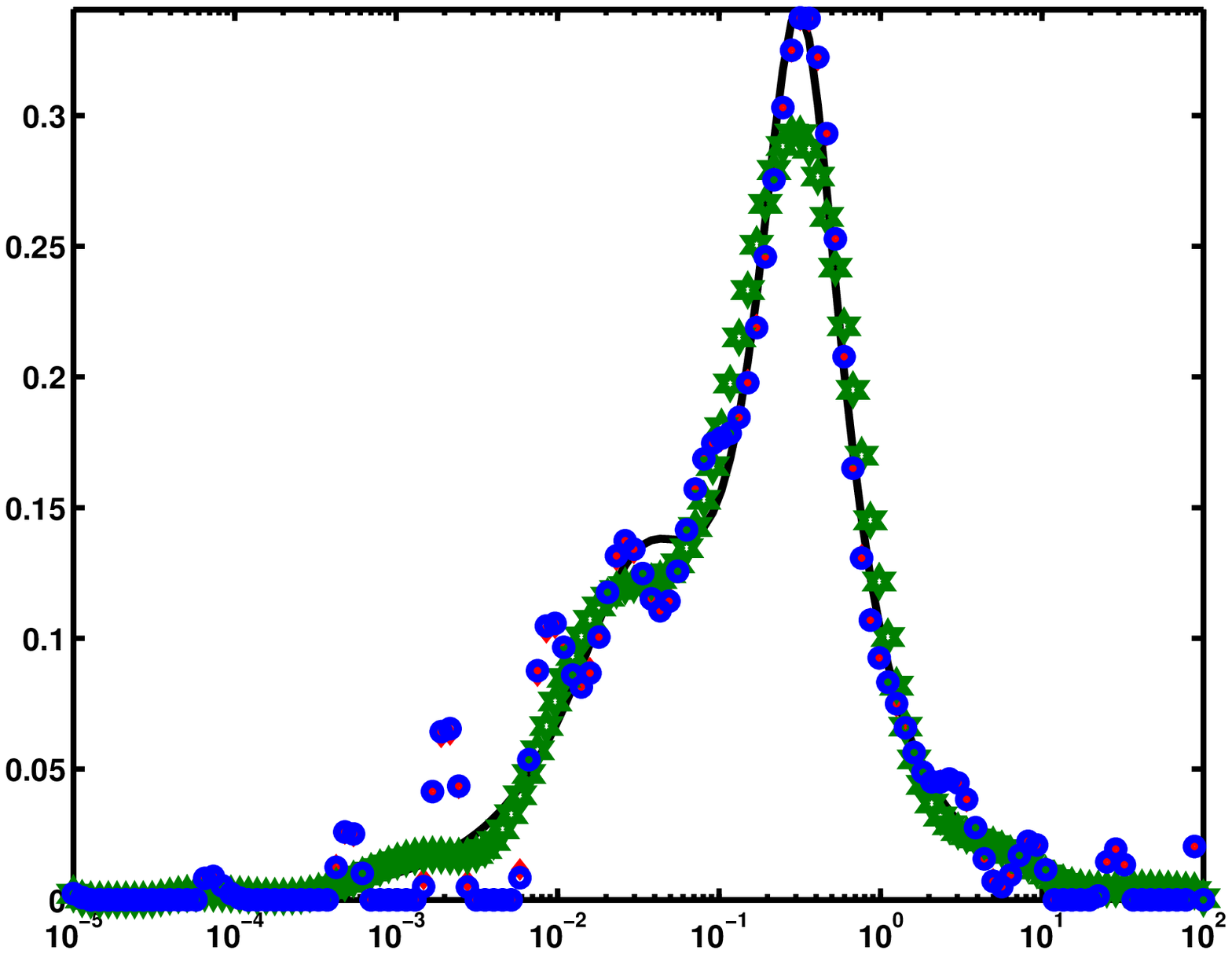}}
\subfigure[$L=L_1$]{\includegraphics[width=1.7in]{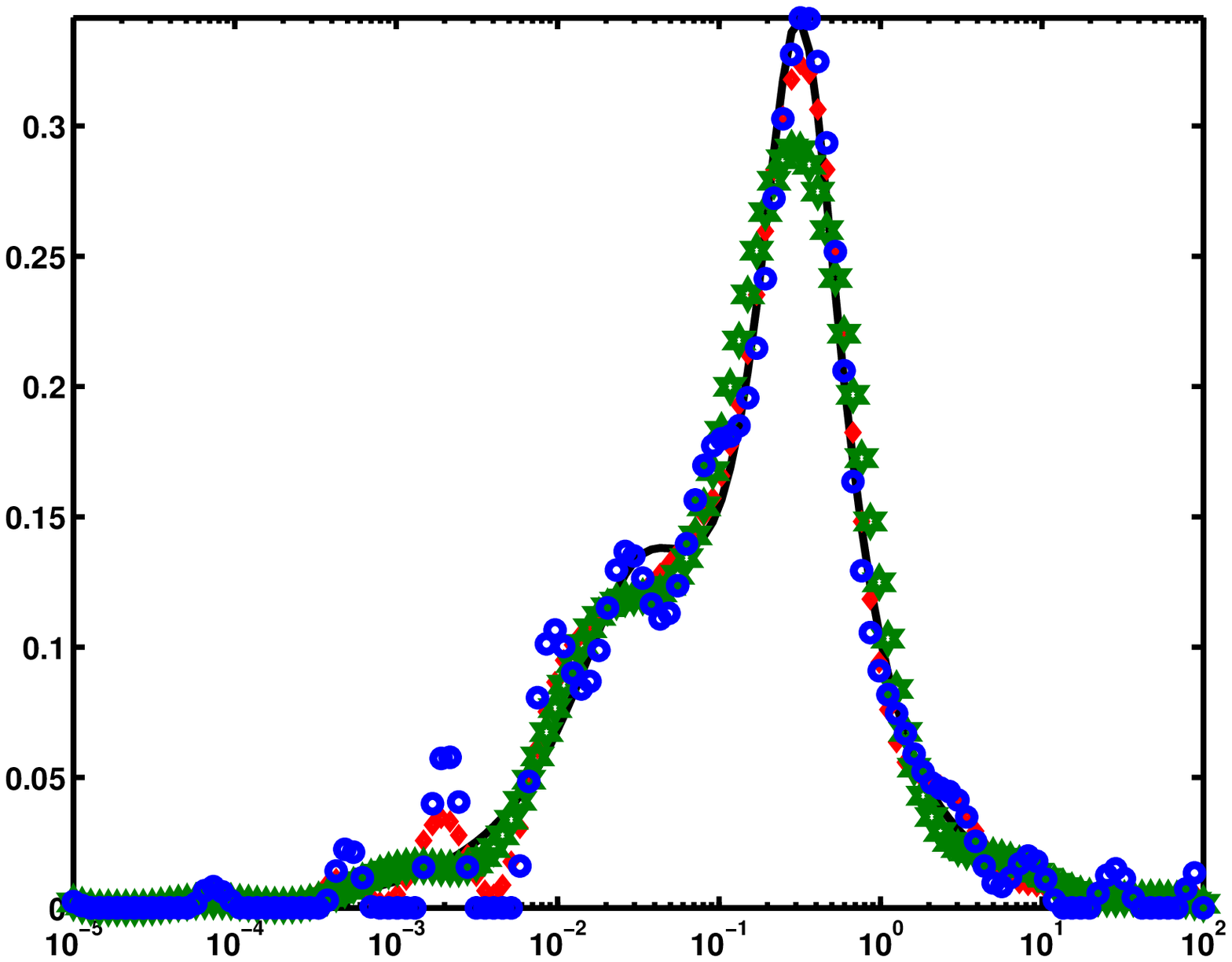}}
\subfigure[$L=L_2$]{\includegraphics[width=1.7in]{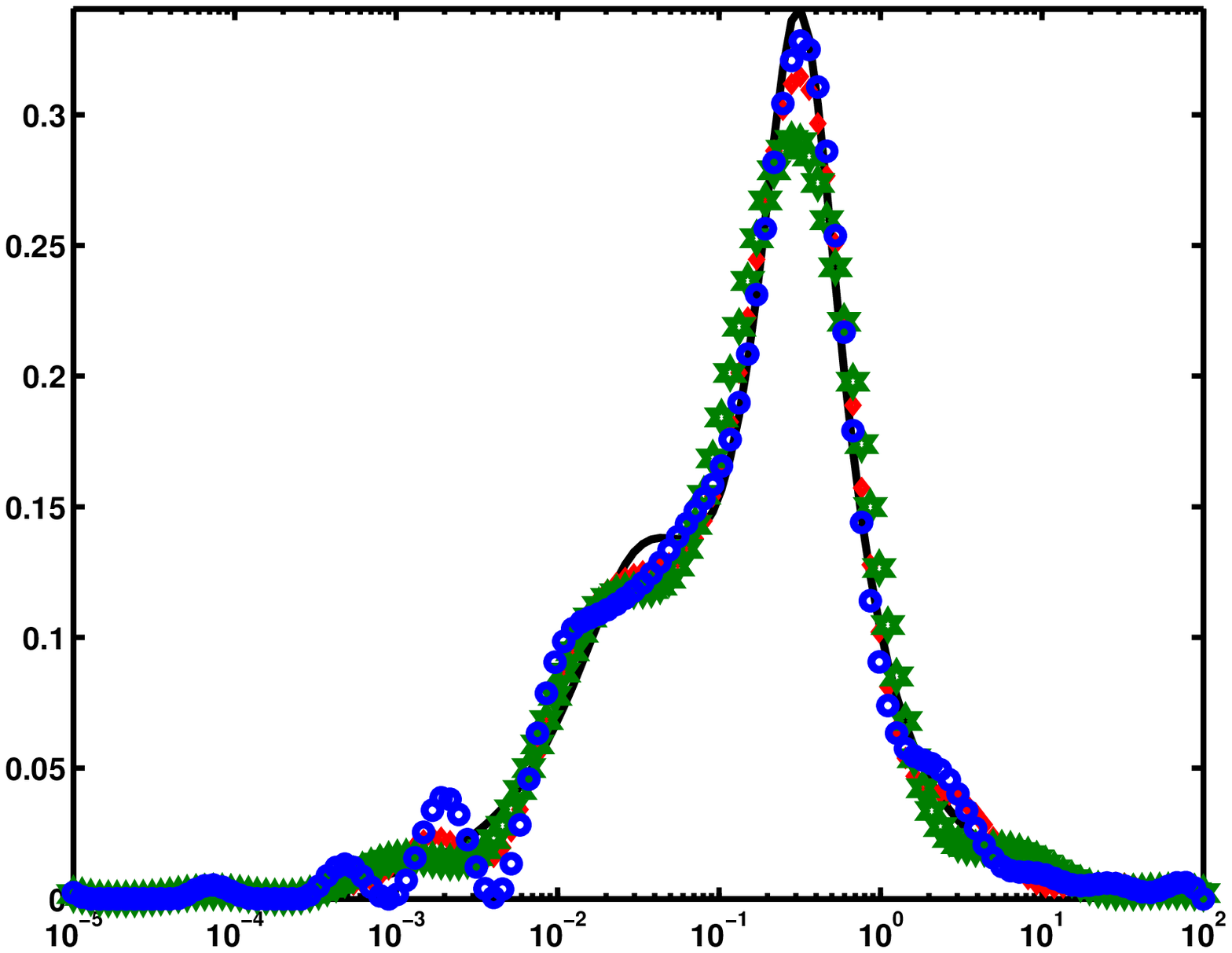}}
\caption{NNLS solutions of RQ-B matrix $A_4$. Noise level $.1\%$ using the SBB algorithm.}
\label{lnfig-lambdachoiceRQ5A4LNSBB}
\end{figure}

\begin{figure}[!ht]
\centering
\subfigure[$L=I$]{\includegraphics[width=1.7in]{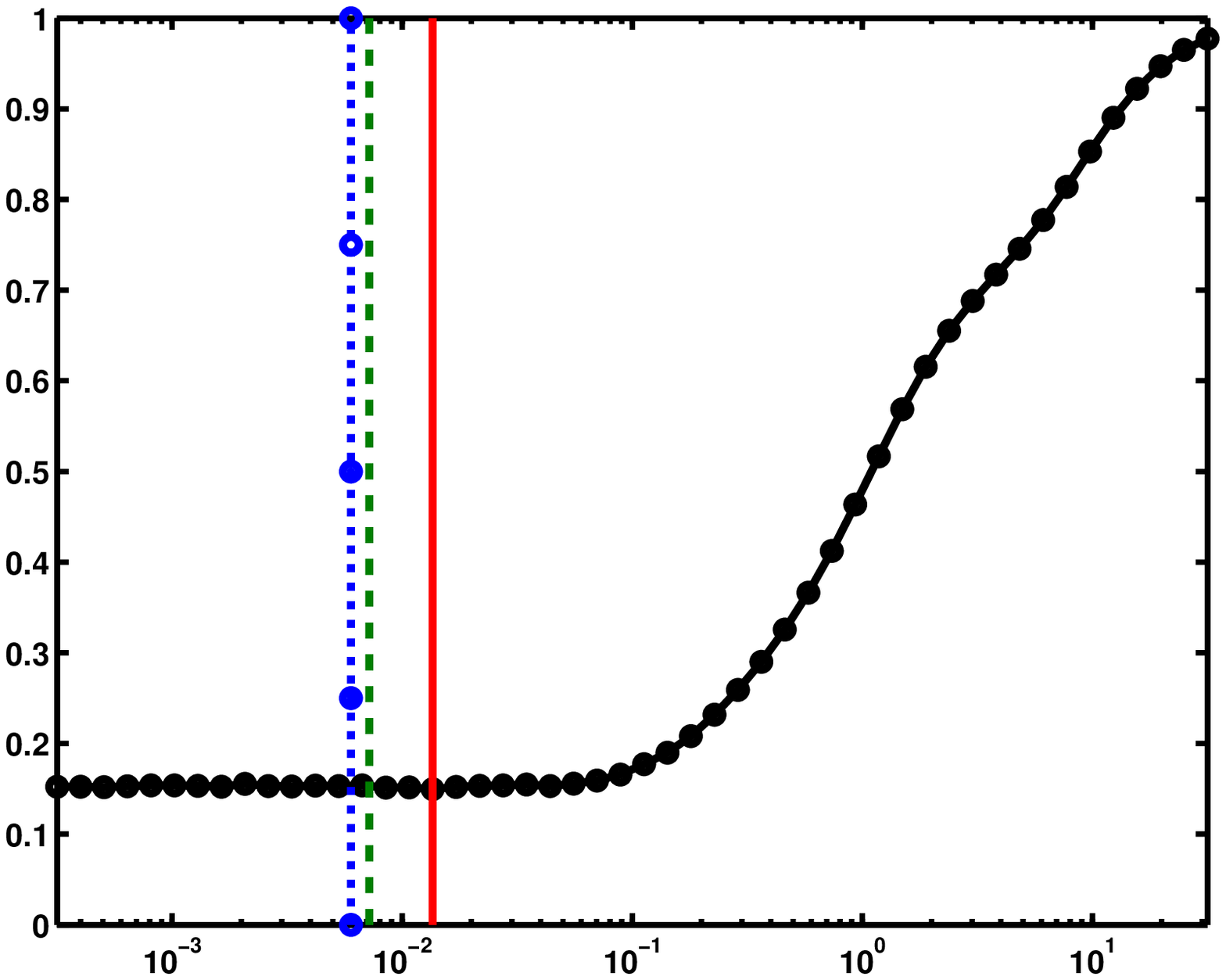}}
\subfigure[$L=L_1$]{\includegraphics[width=1.7in]{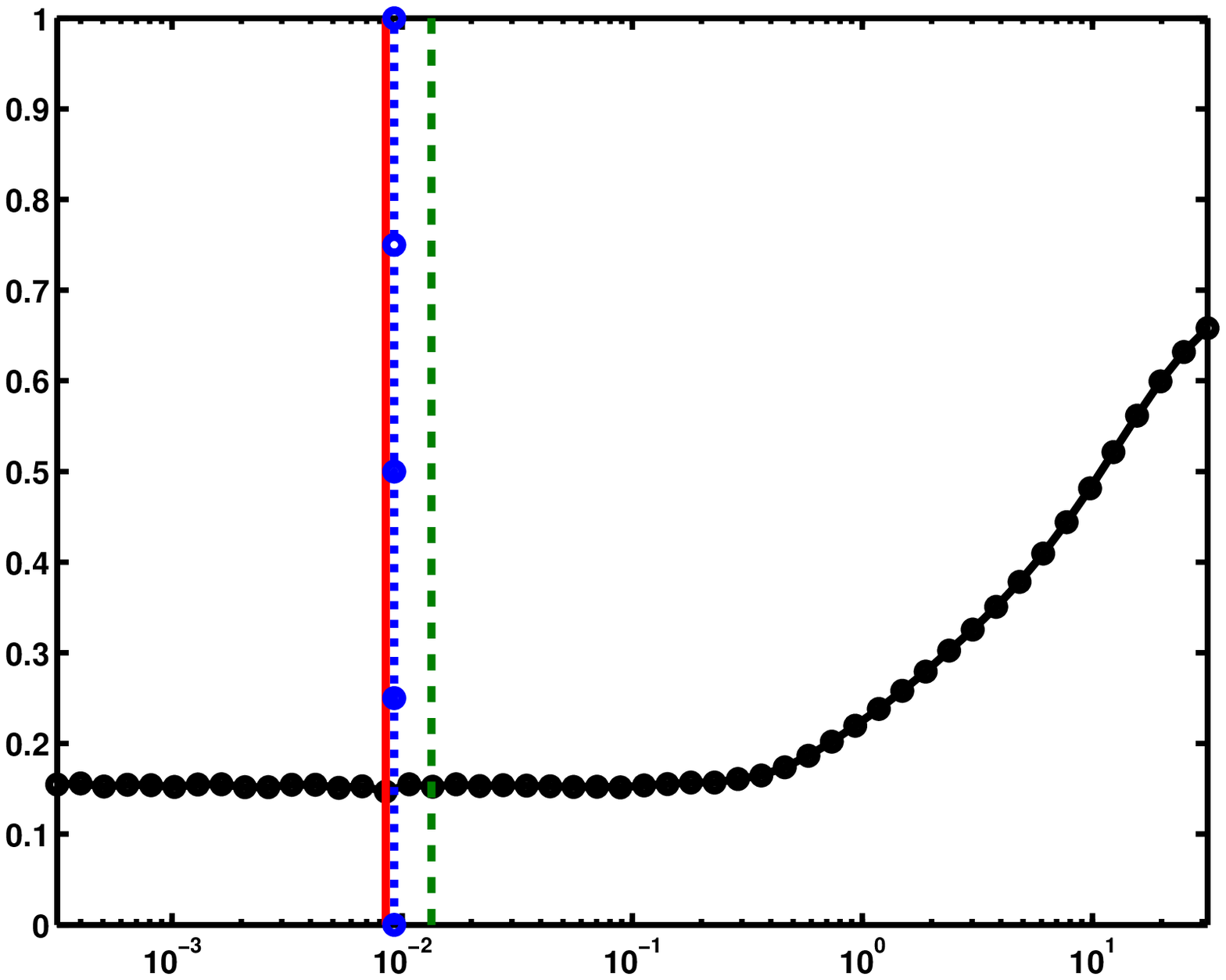}}
\subfigure[$L=L_2$]{\includegraphics[width=1.7in]{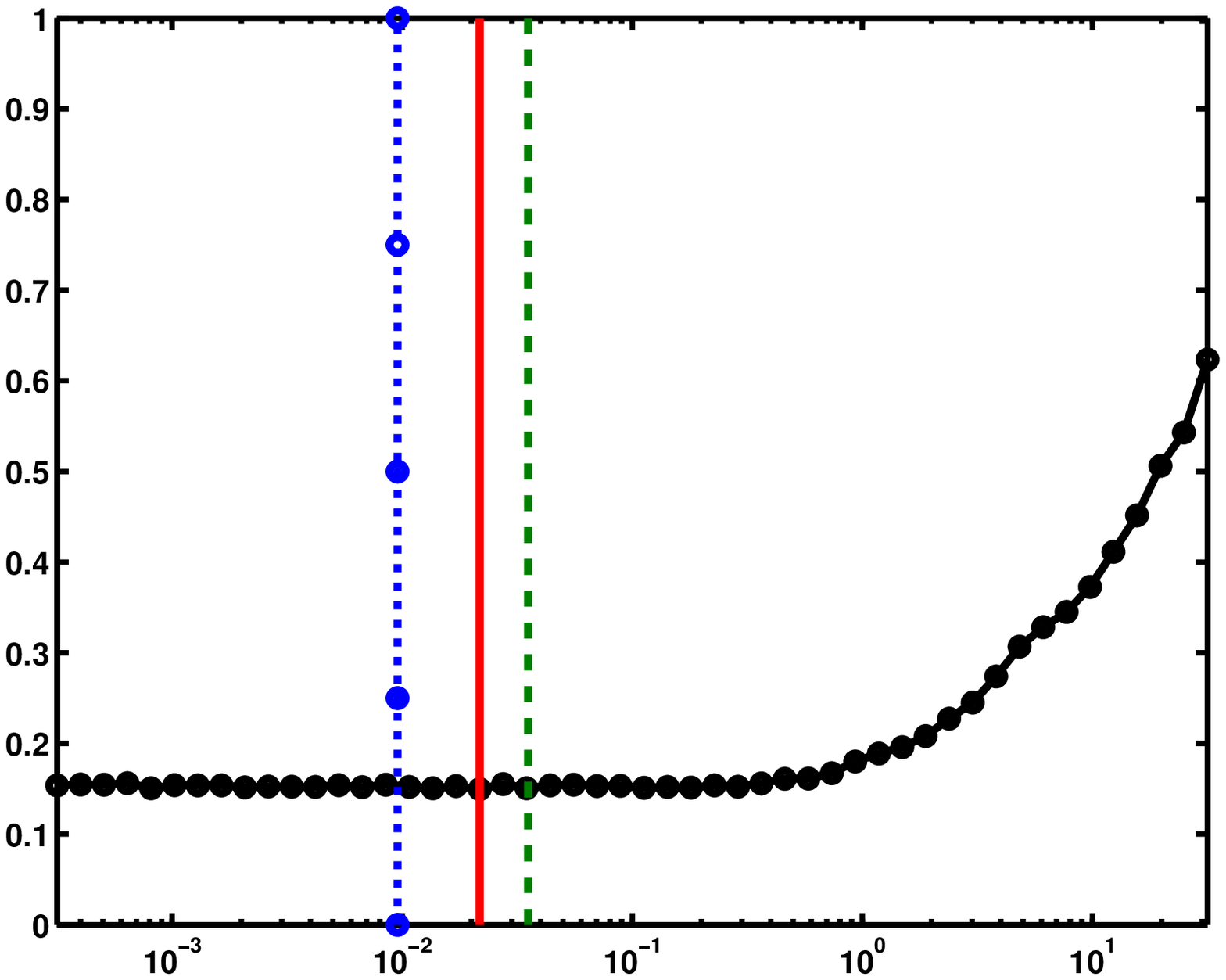}}
\subfigure[$L=I$]{\includegraphics[width=1.7in]{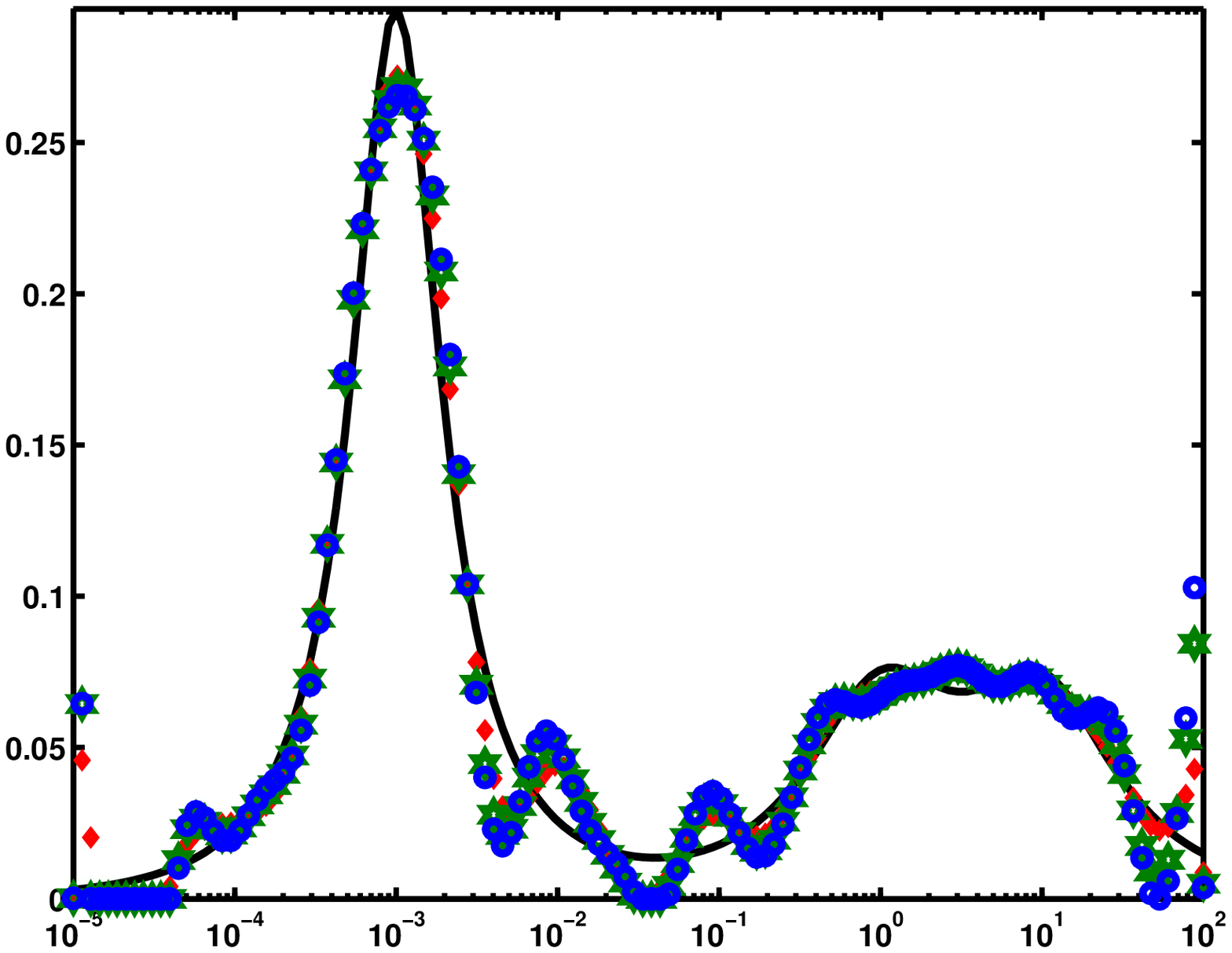}}
\subfigure[$L=L_1$]{\includegraphics[width=1.7in]{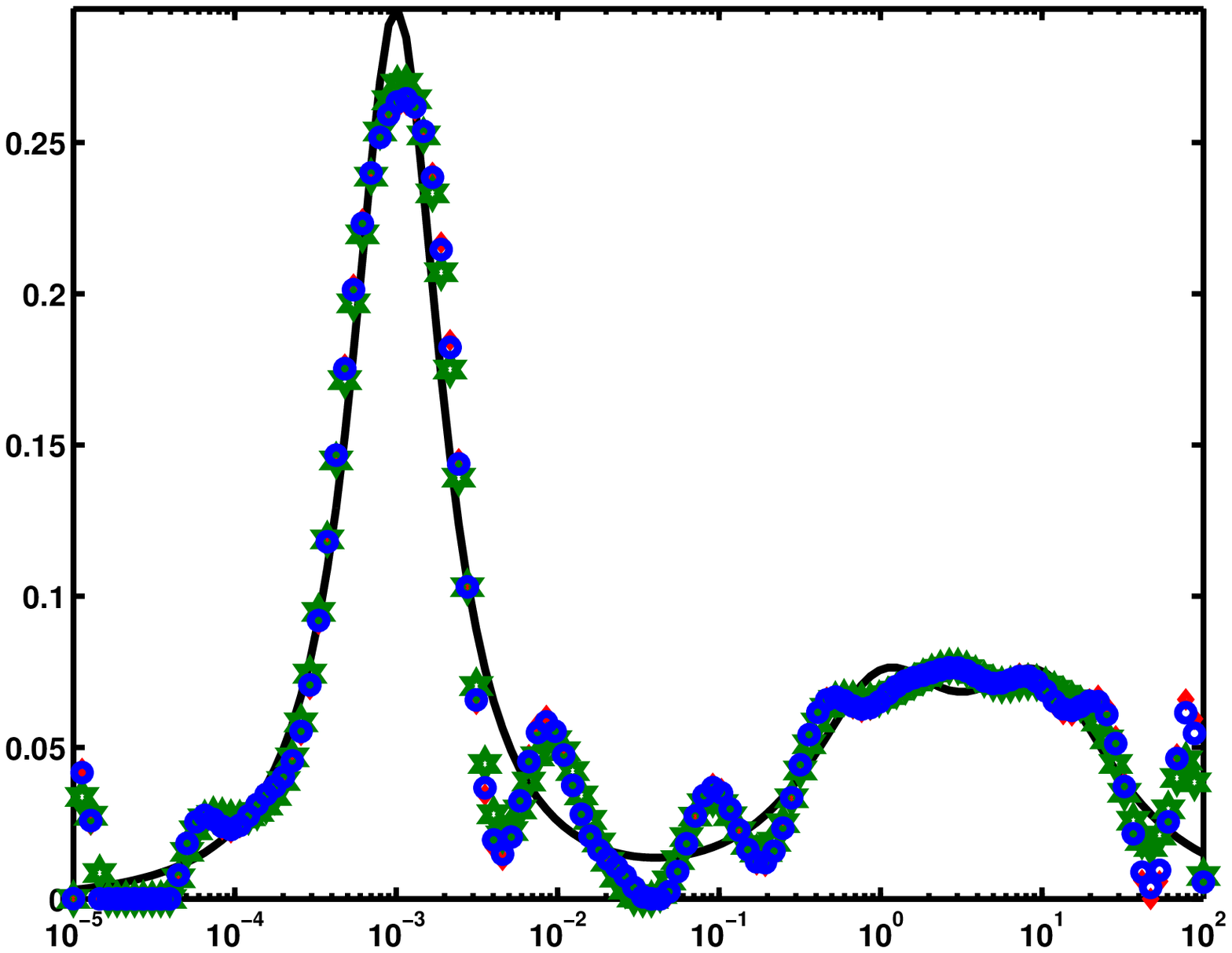}}
\subfigure[$L=L_2$]{\includegraphics[width=1.7in]{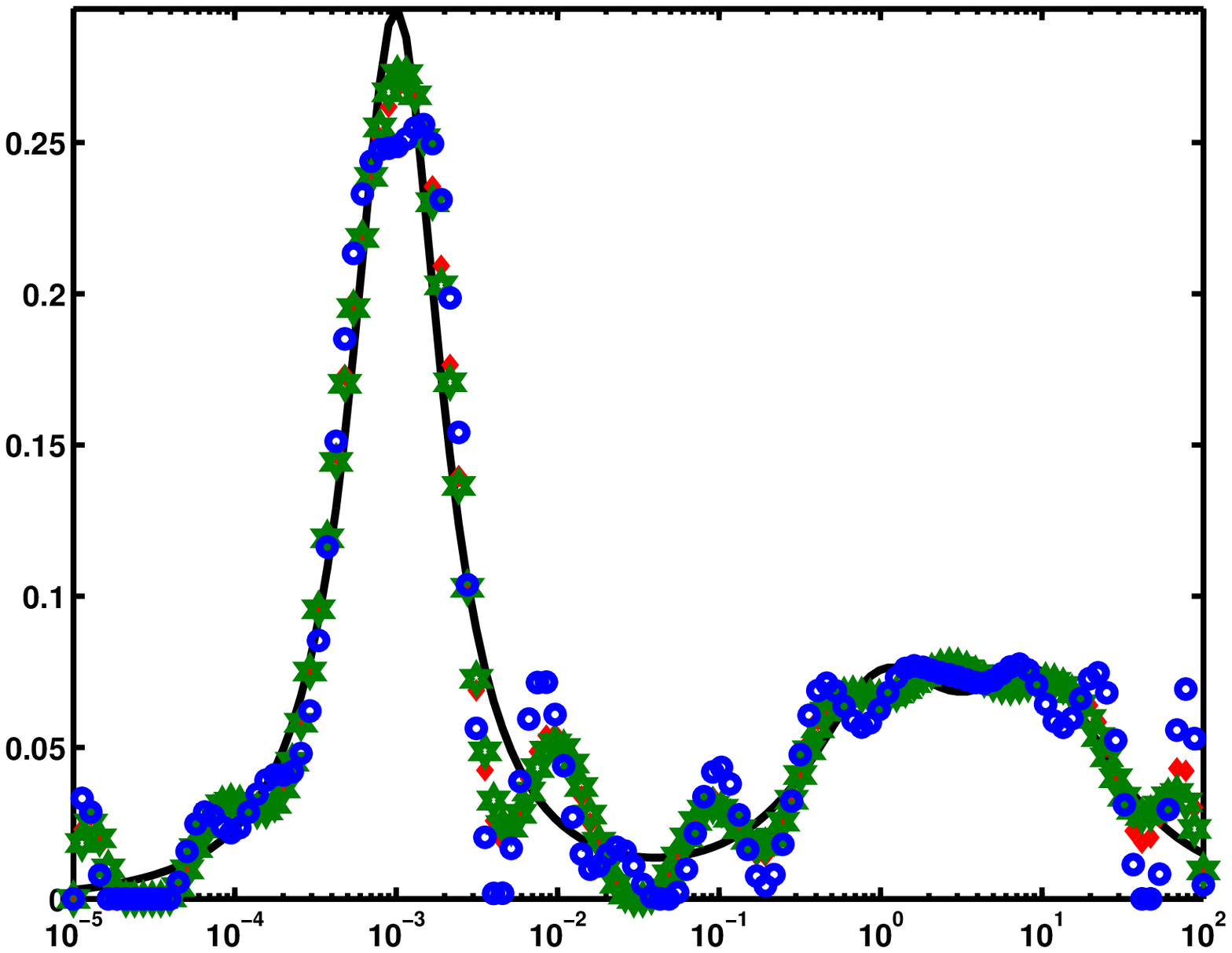}}
\caption{NNLS solutions of RQ-C matrix $A_4$. Noise level $.1\%$ using the SBB algorithm.}
\label{lnfig-lambdachoiceRQ6A4LNSBB}
\end{figure}

\begin{figure}[!ht]
\centering
\subfigure[$L=I$]{\includegraphics[width=1.7in]{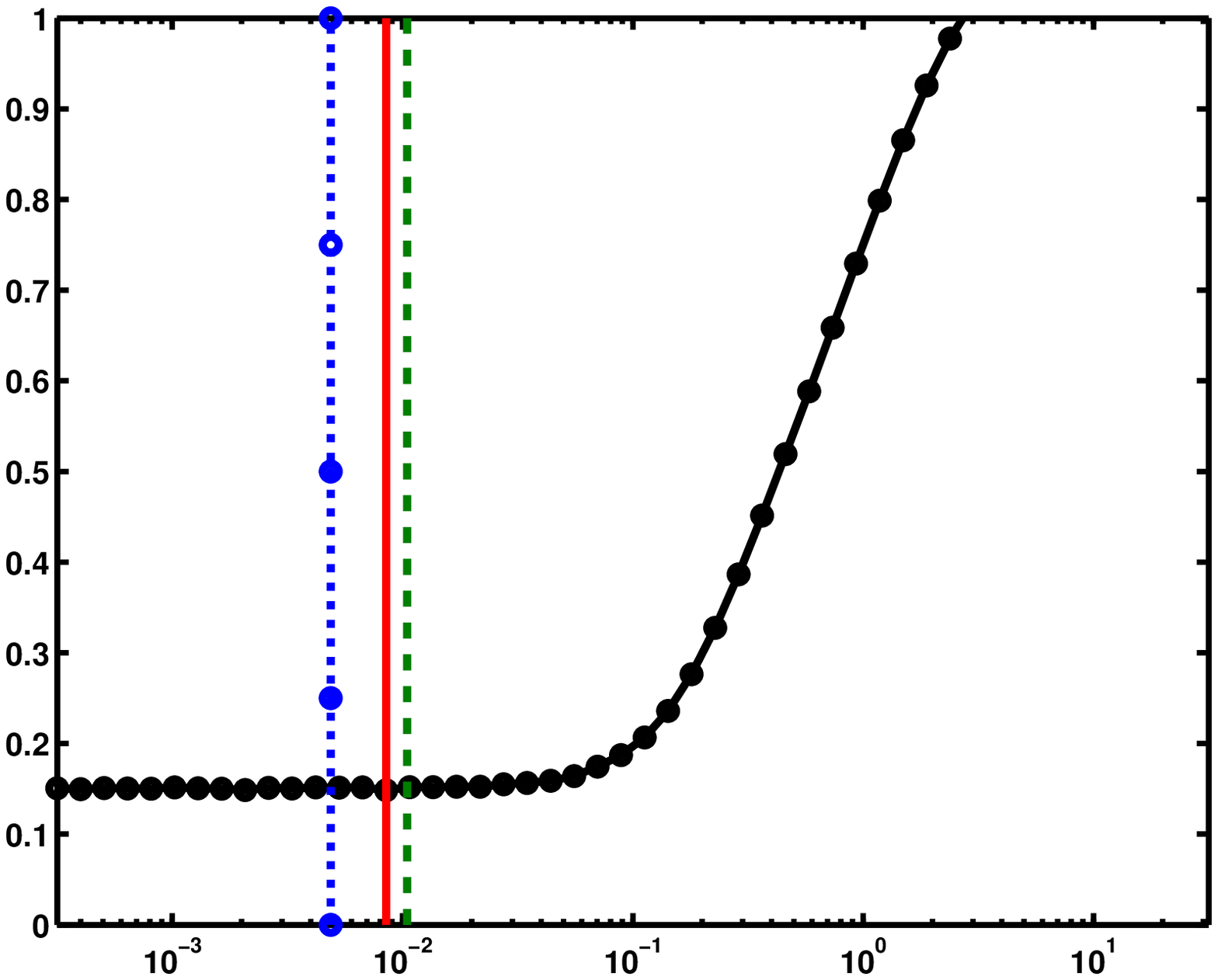}}
\subfigure[$L=L_1$]{\includegraphics[width=1.7in]{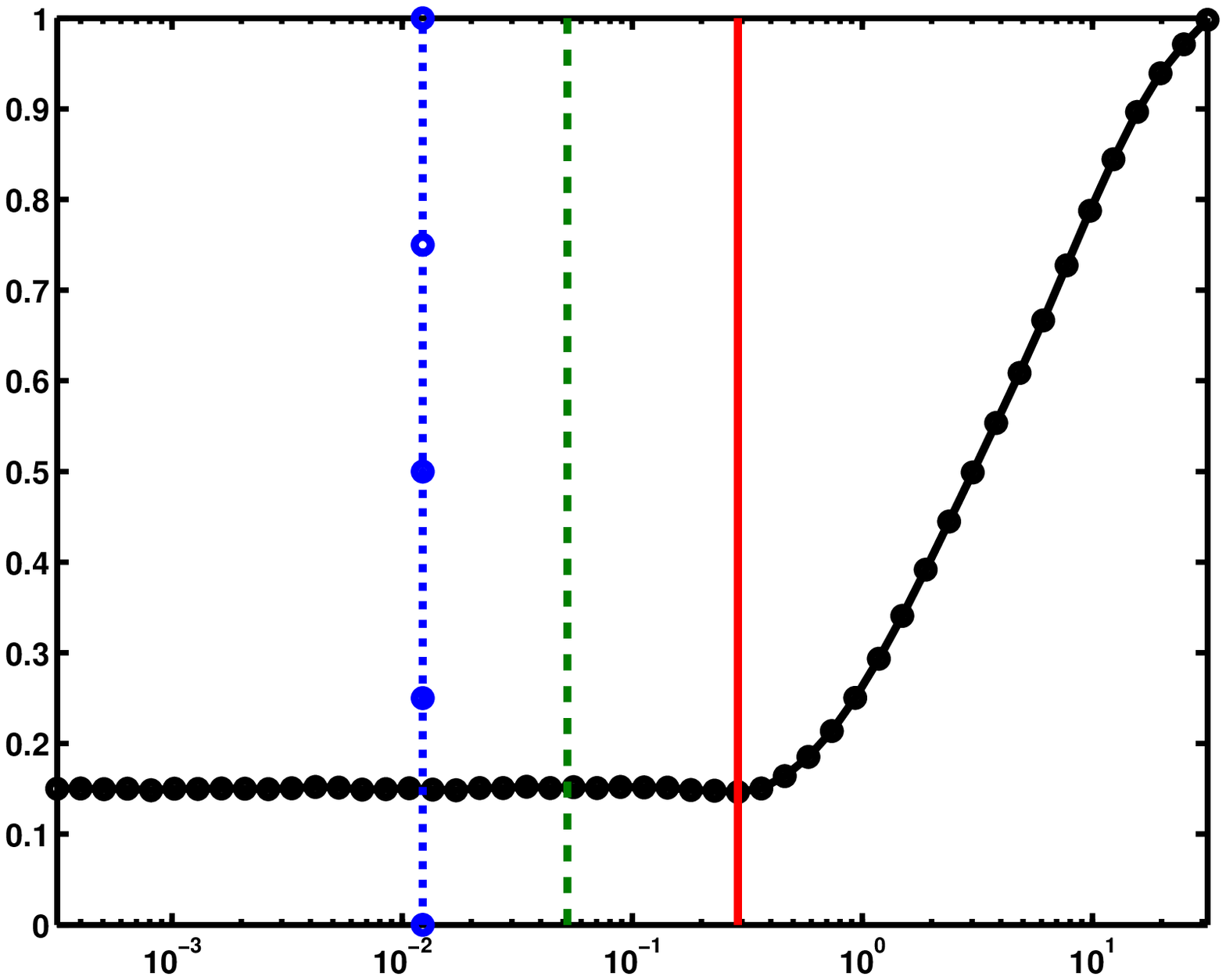}}
\subfigure[$L=L_2$]{\includegraphics[width=1.7in]{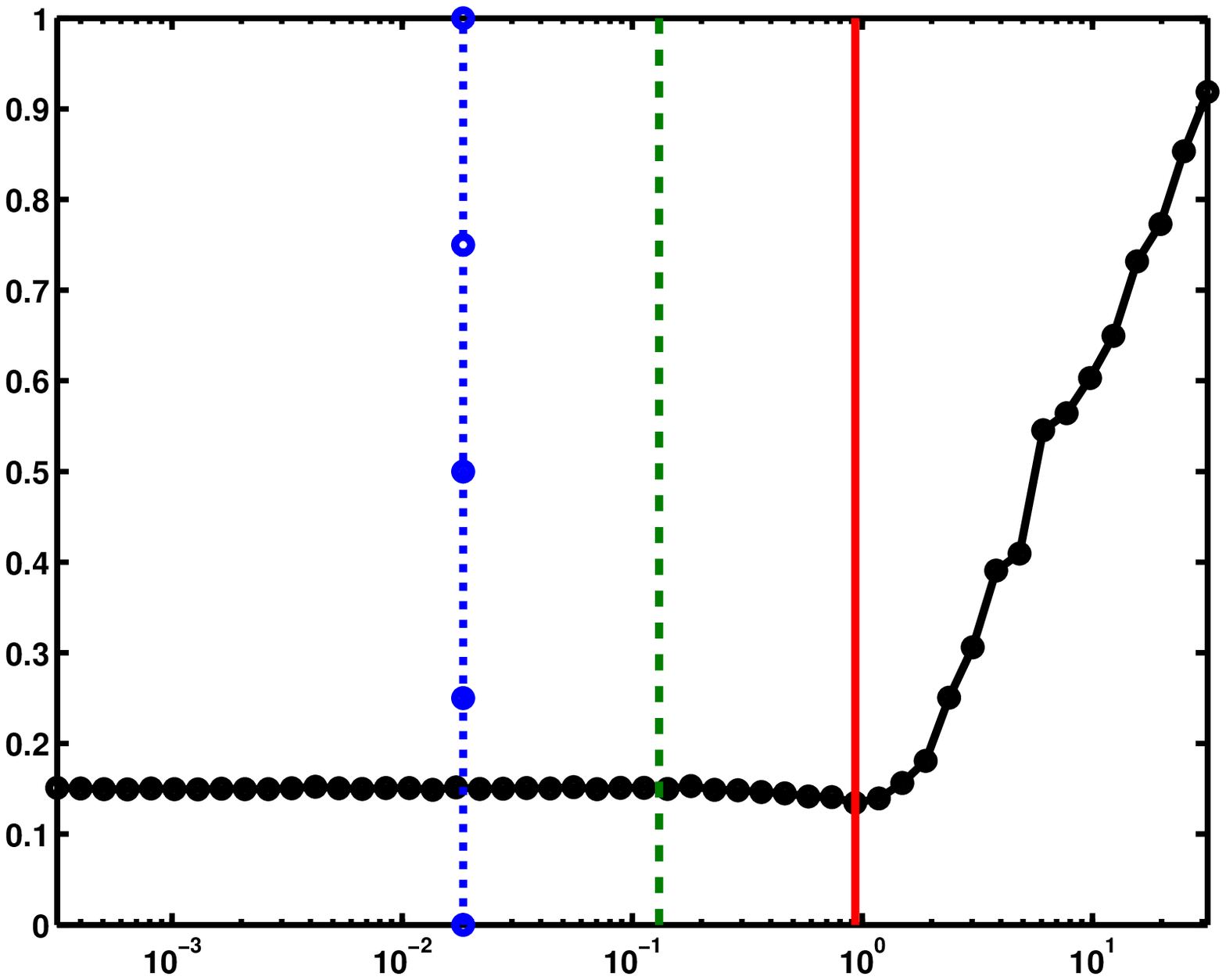}}
\subfigure[$L=I$]{\includegraphics[width=1.7in]{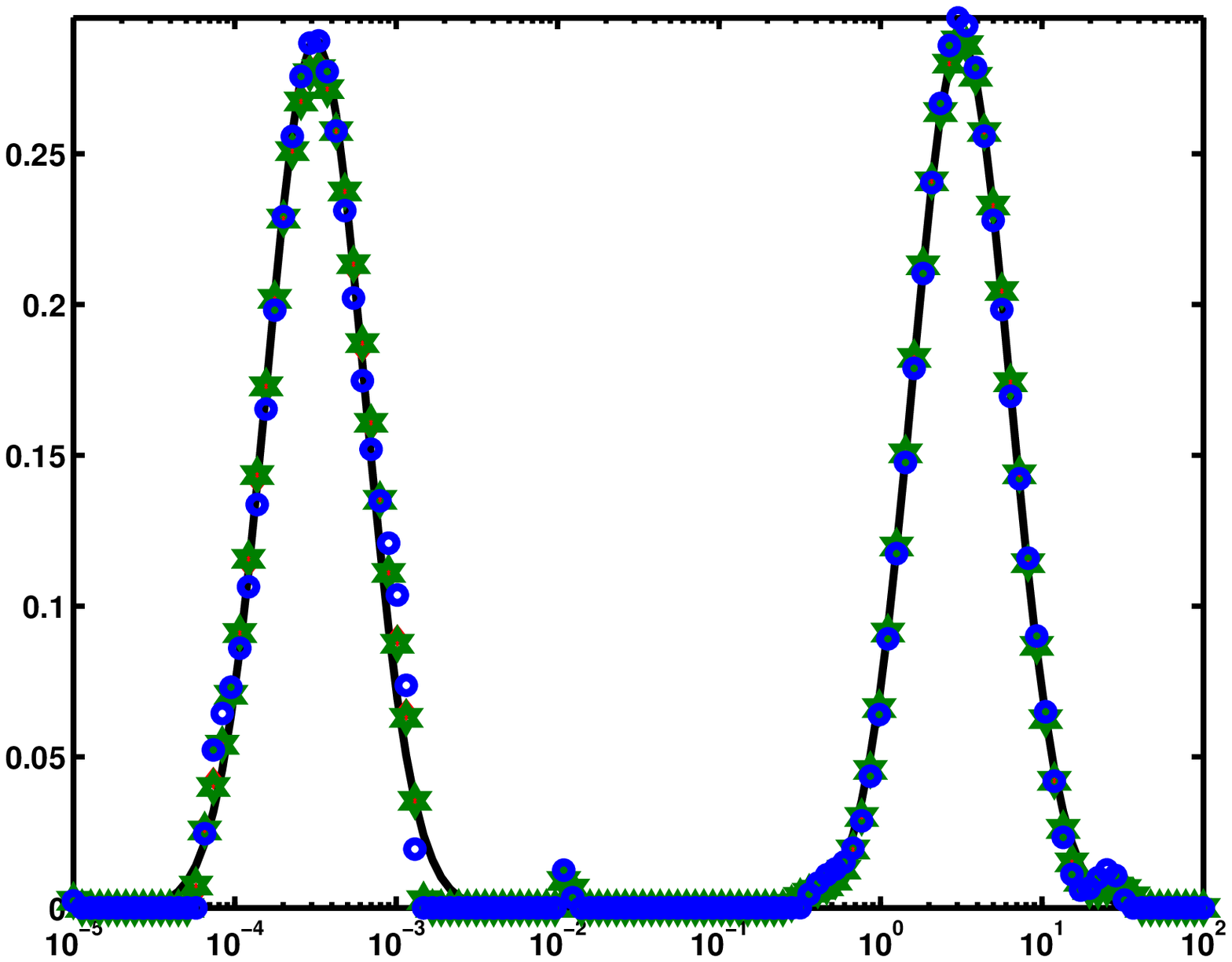}}
\subfigure[$L=L_1$]{\includegraphics[width=1.7in]{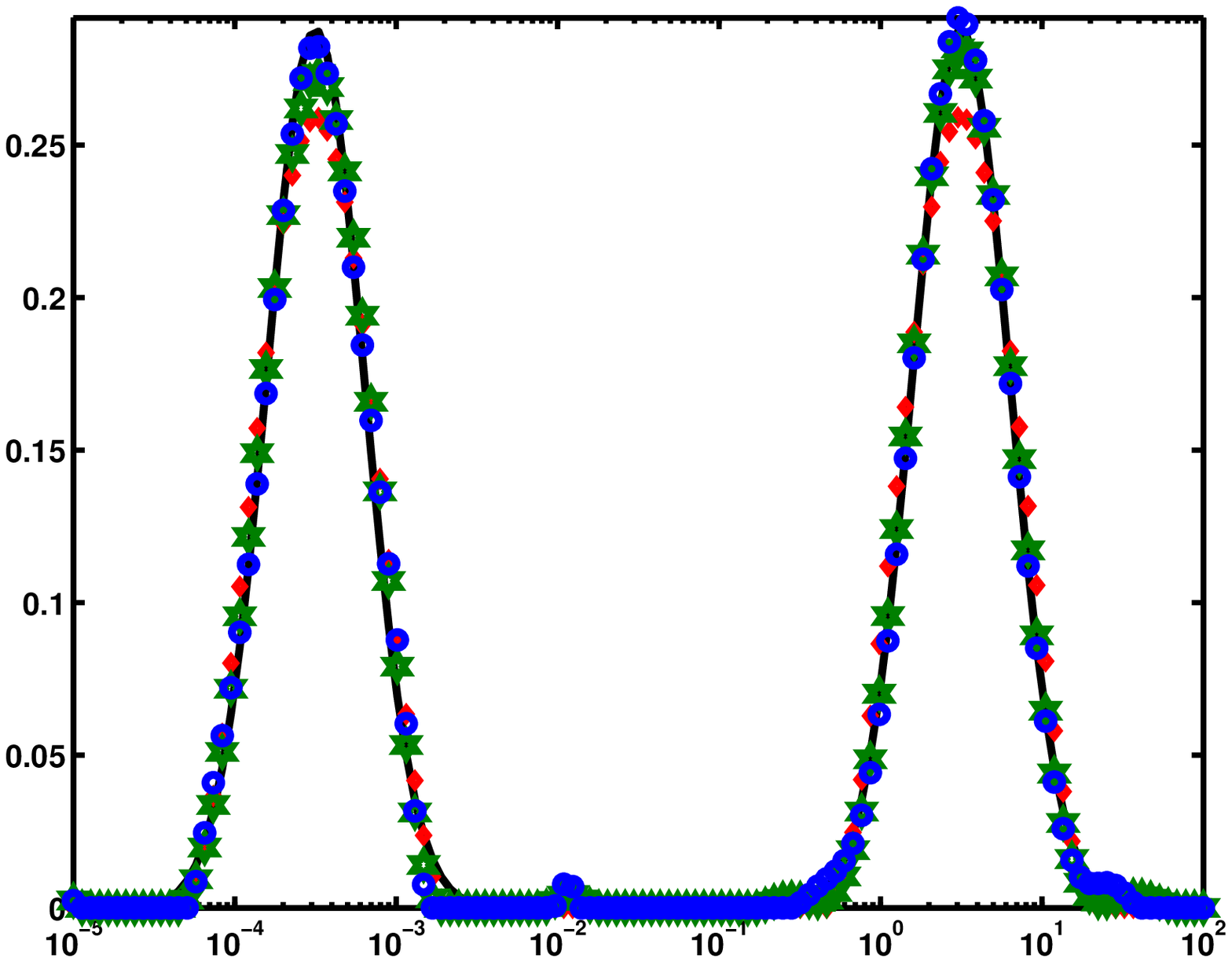}}
\subfigure[$L=L_2$]{\includegraphics[width=1.7in]{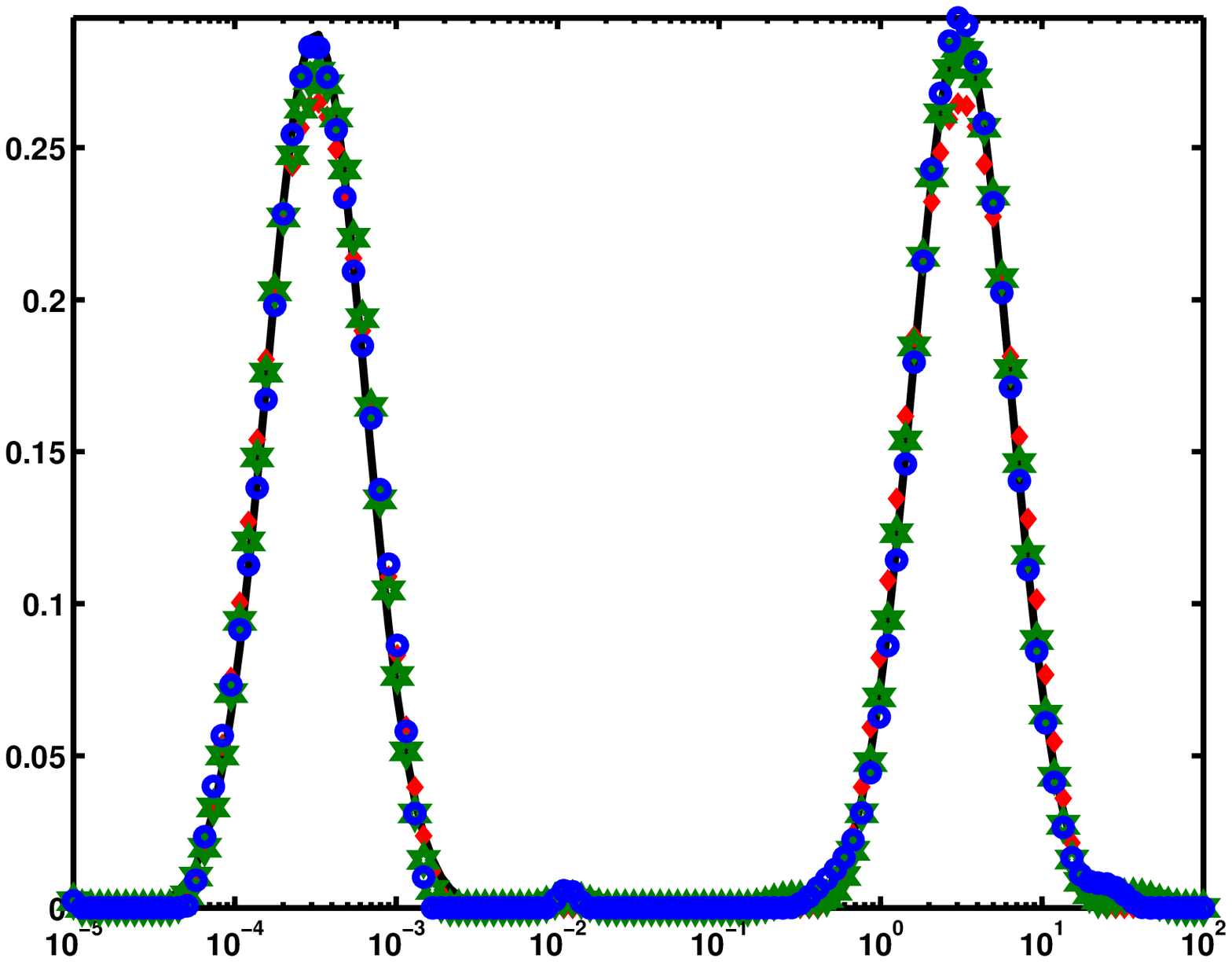}}
\caption{NNLS solutions of LN-A matrix $A_4$. Noise level $.1\%$ using the SBB algorithm.}
\label{lnfig-lambdachoiceLN2A4LNSBB}
\end{figure}

\begin{figure}[!ht]
\centering
\subfigure[$L=I$]{\includegraphics[width=1.7in]{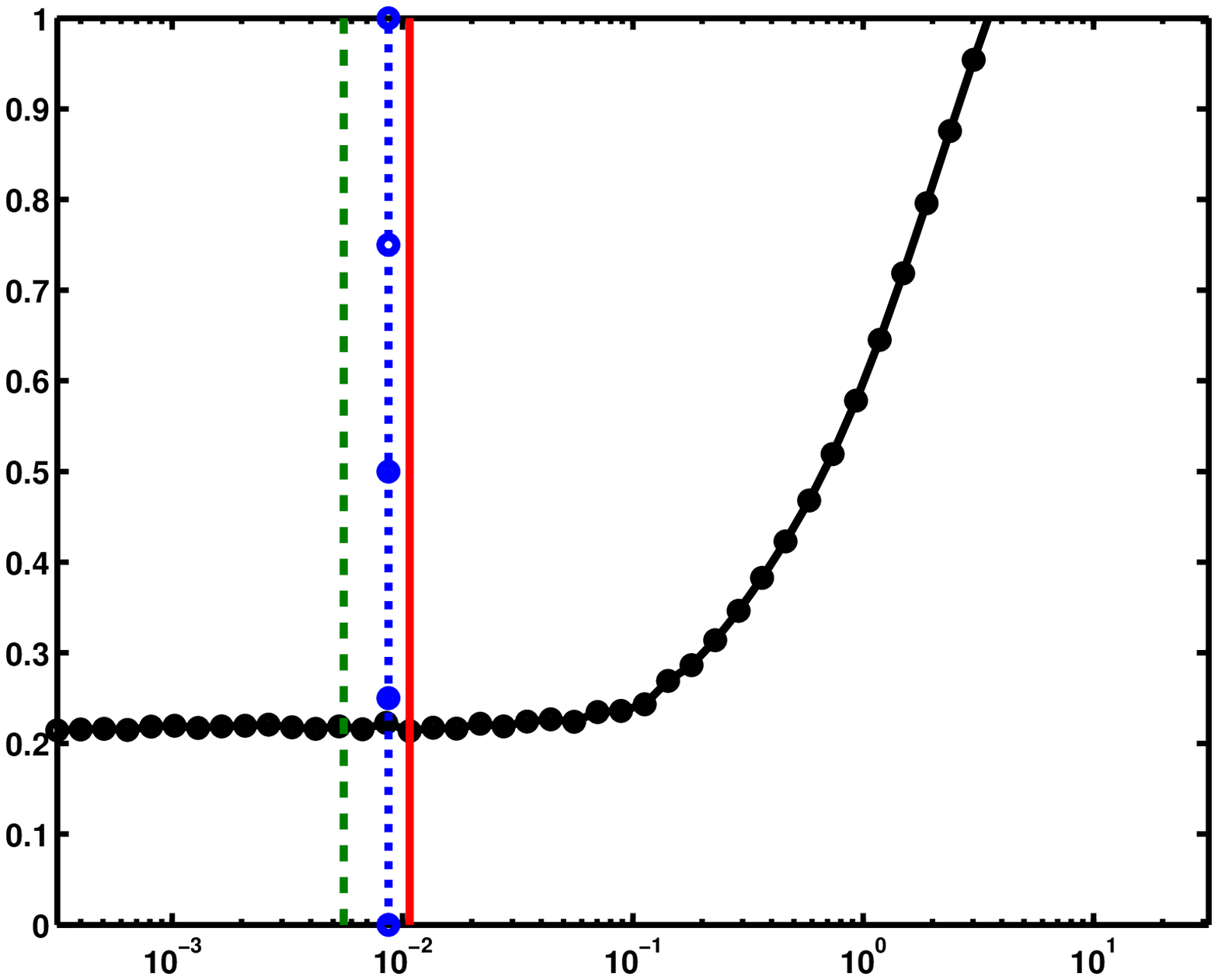}}
\subfigure[$L=L_1$]{\includegraphics[width=1.7in]{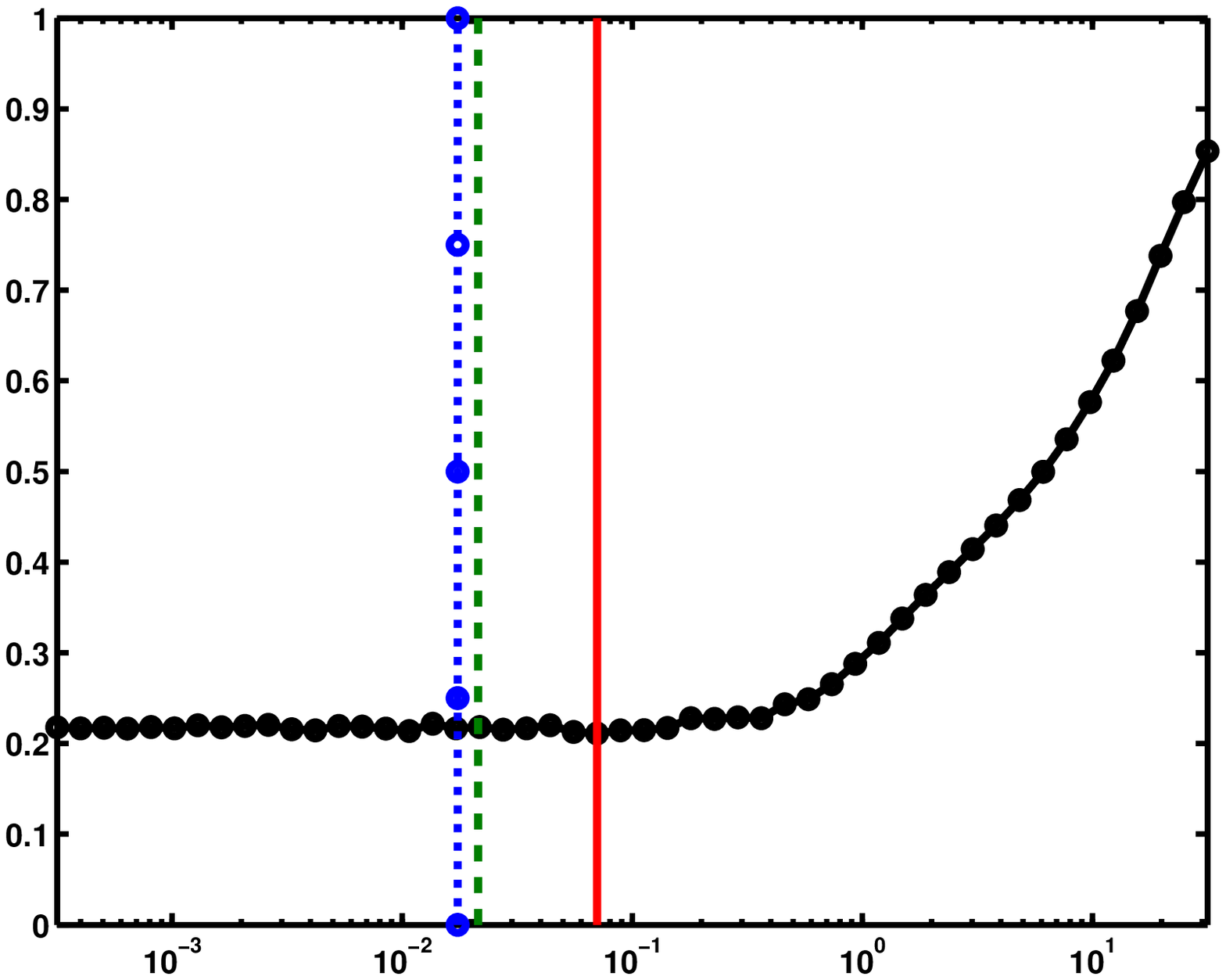}}
\subfigure[$L=L_2$]{\includegraphics[width=1.7in]{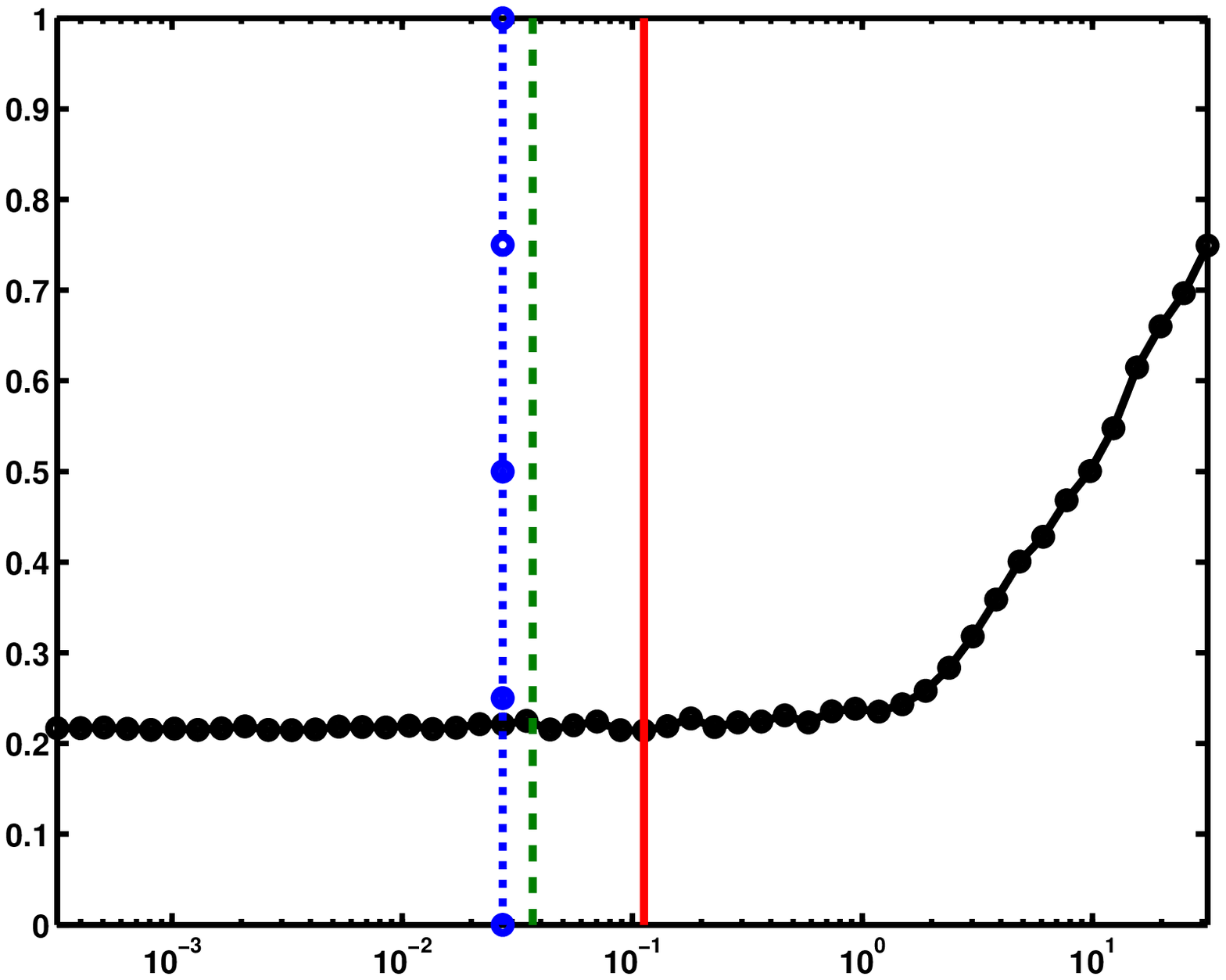}}
\subfigure[$L=I$]{\includegraphics[width=1.7in]{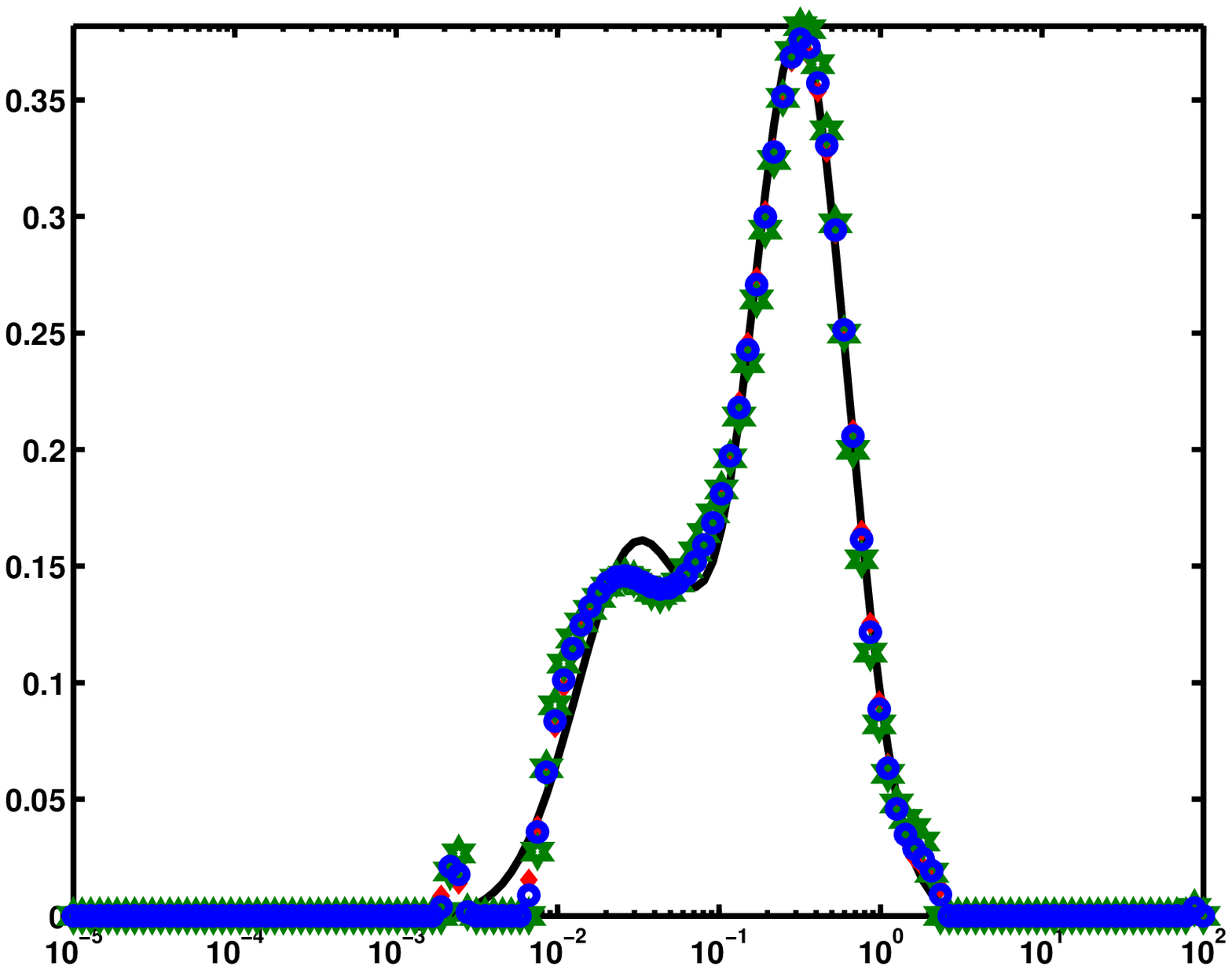}}
\subfigure[$L=L_1$]{\includegraphics[width=1.7in]{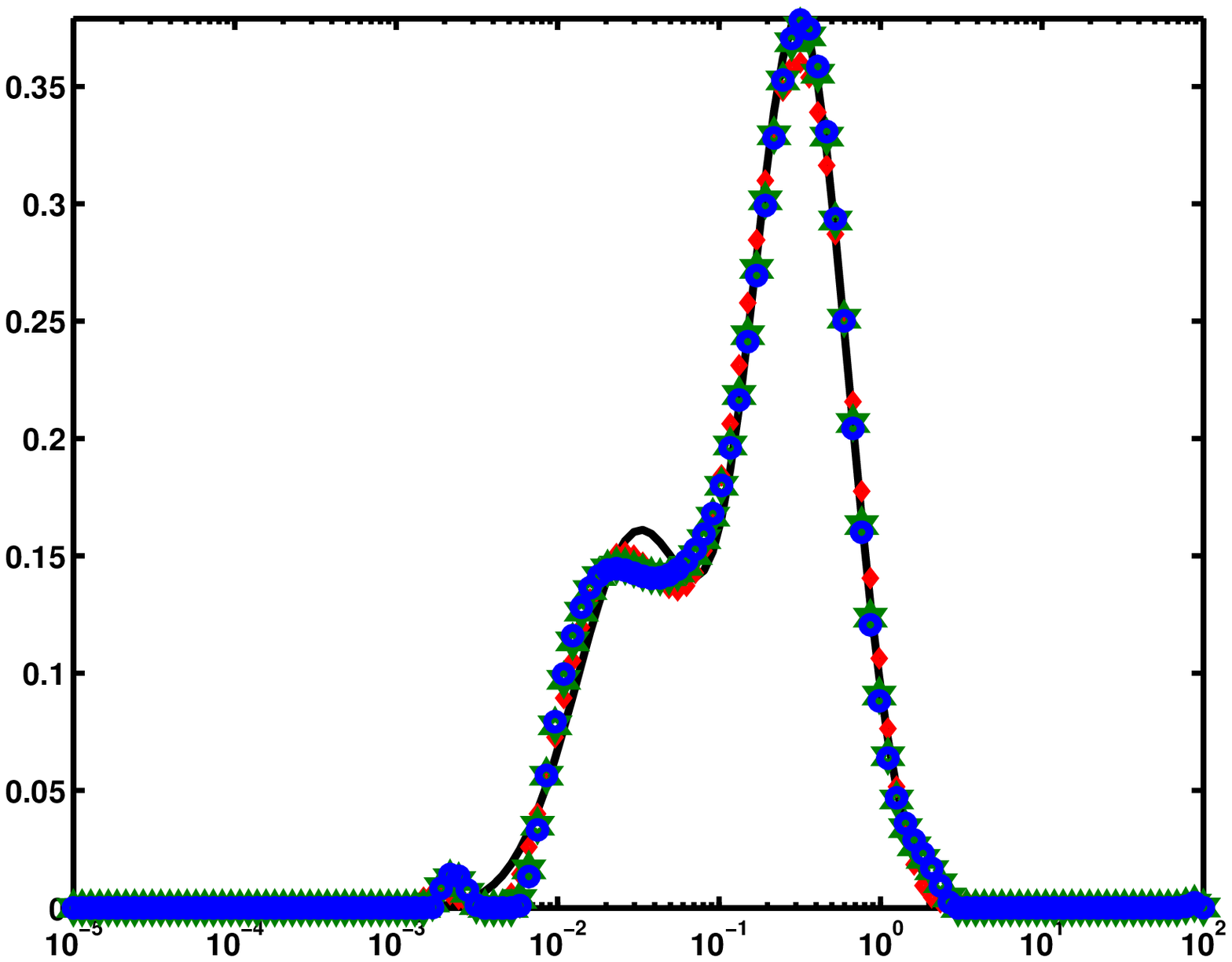}}
\subfigure[$L=L_2$]{\includegraphics[width=1.7in]{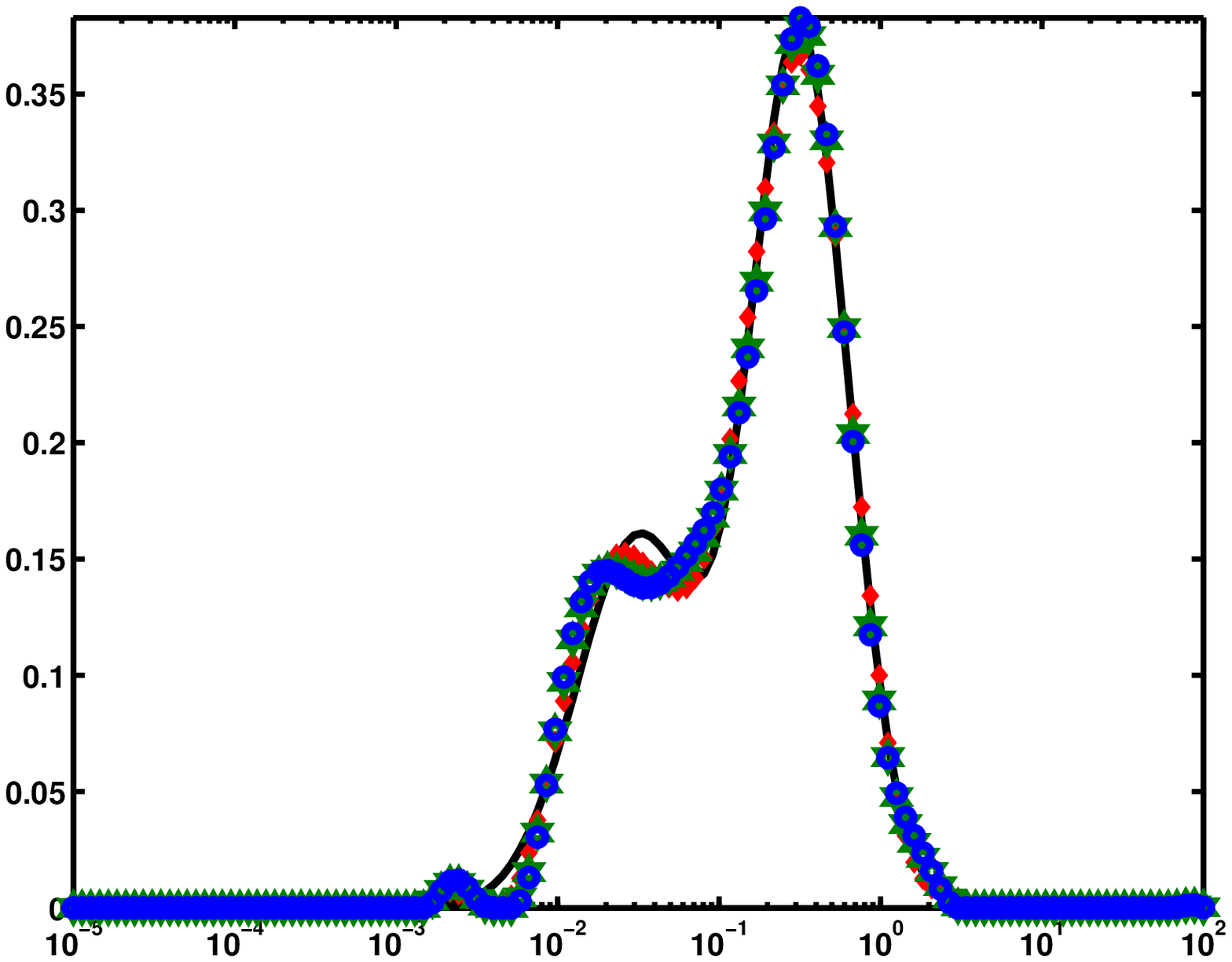}}
\caption{NNLS solutions of LN-B matrix $A_4$. Noise level $.1\%$ using the SBB algorithm.}
\label{lnfig-lambdachoiceLN5A4LNSBB}
\end{figure}

\begin{figure}[!ht]
\centering
\subfigure[$L=I$]{\includegraphics[width=1.7in]{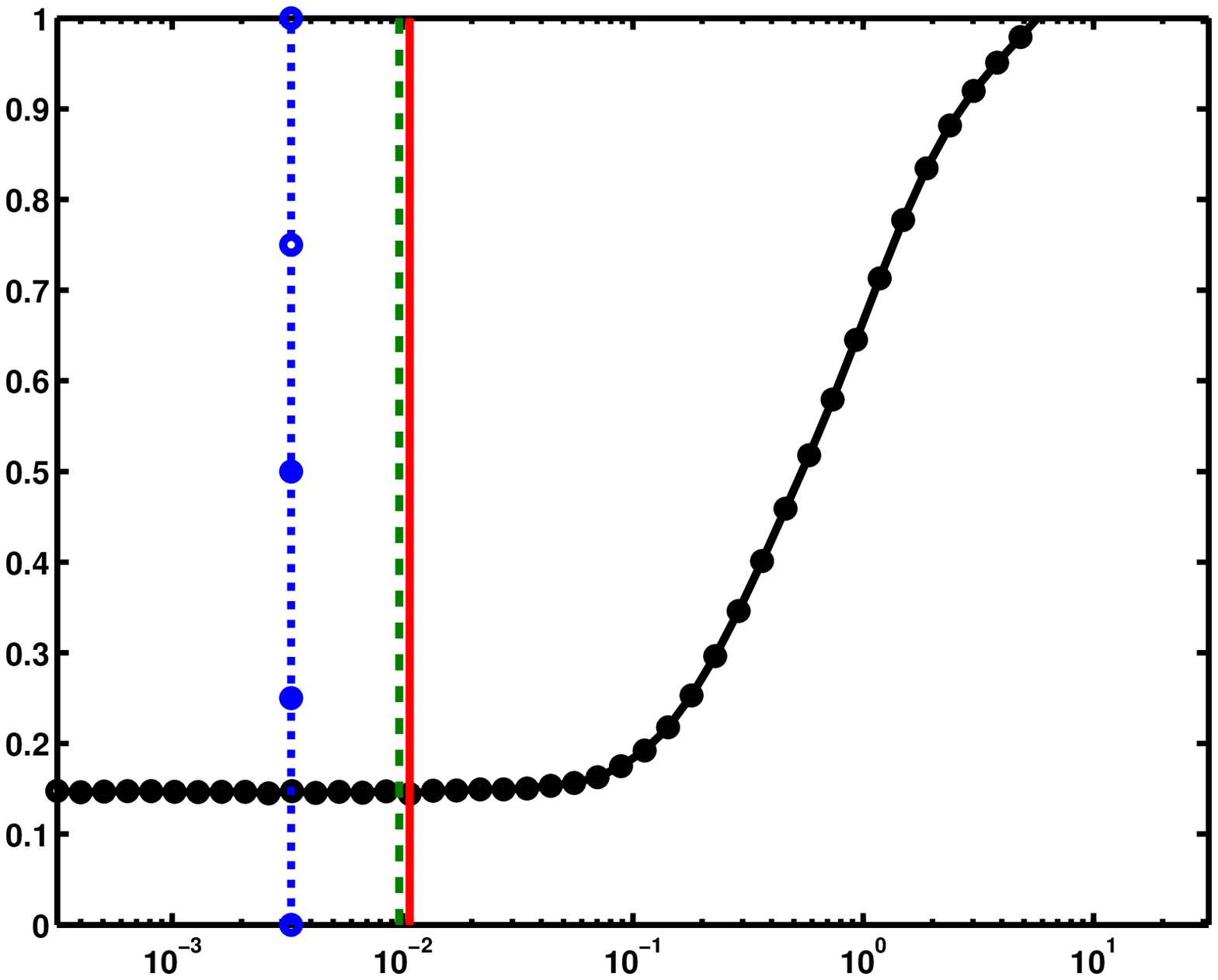}}
\subfigure[$L=L_1$]{\includegraphics[width=1.7in]{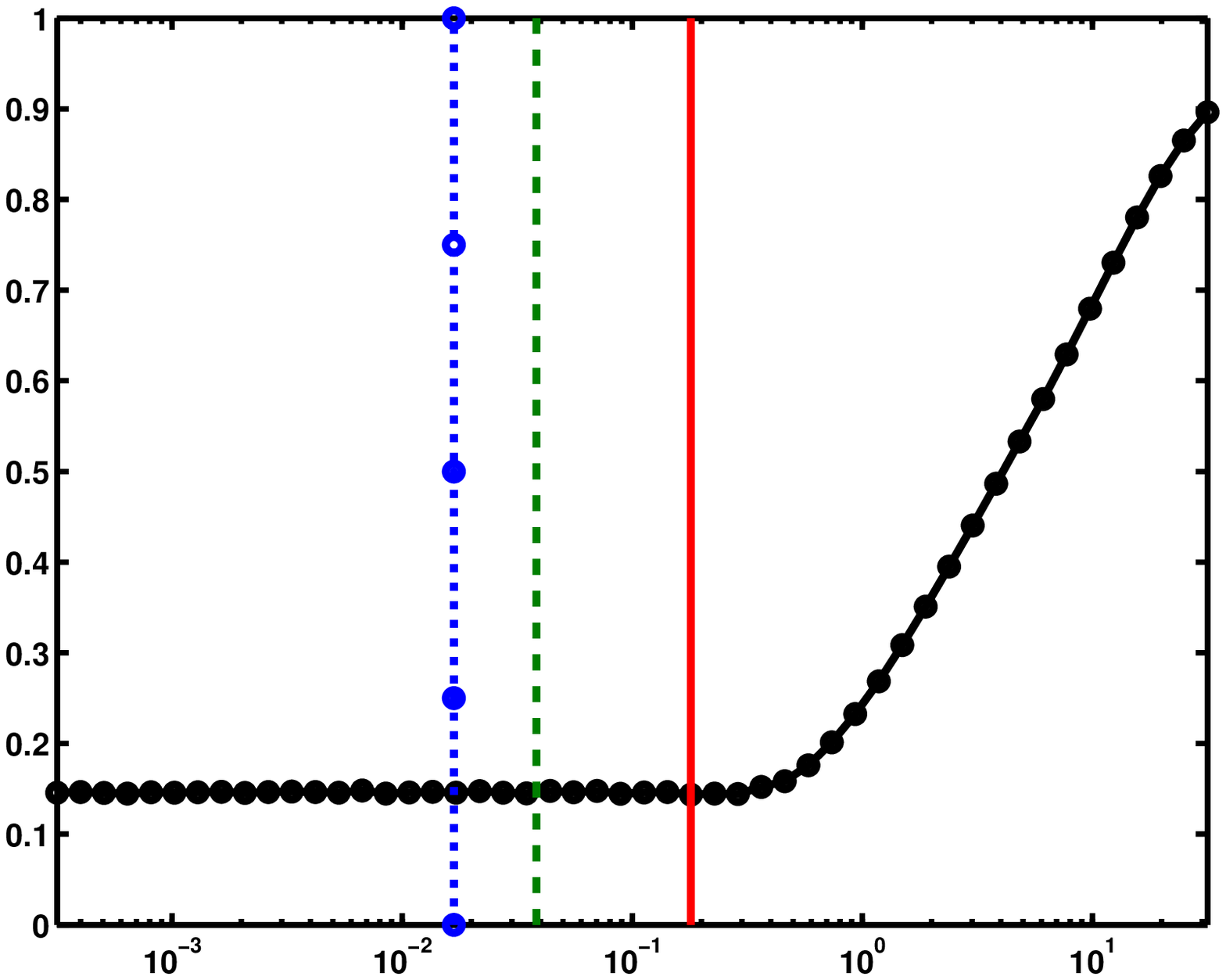}}
\subfigure[$L=L_2$]{\includegraphics[width=1.7in]{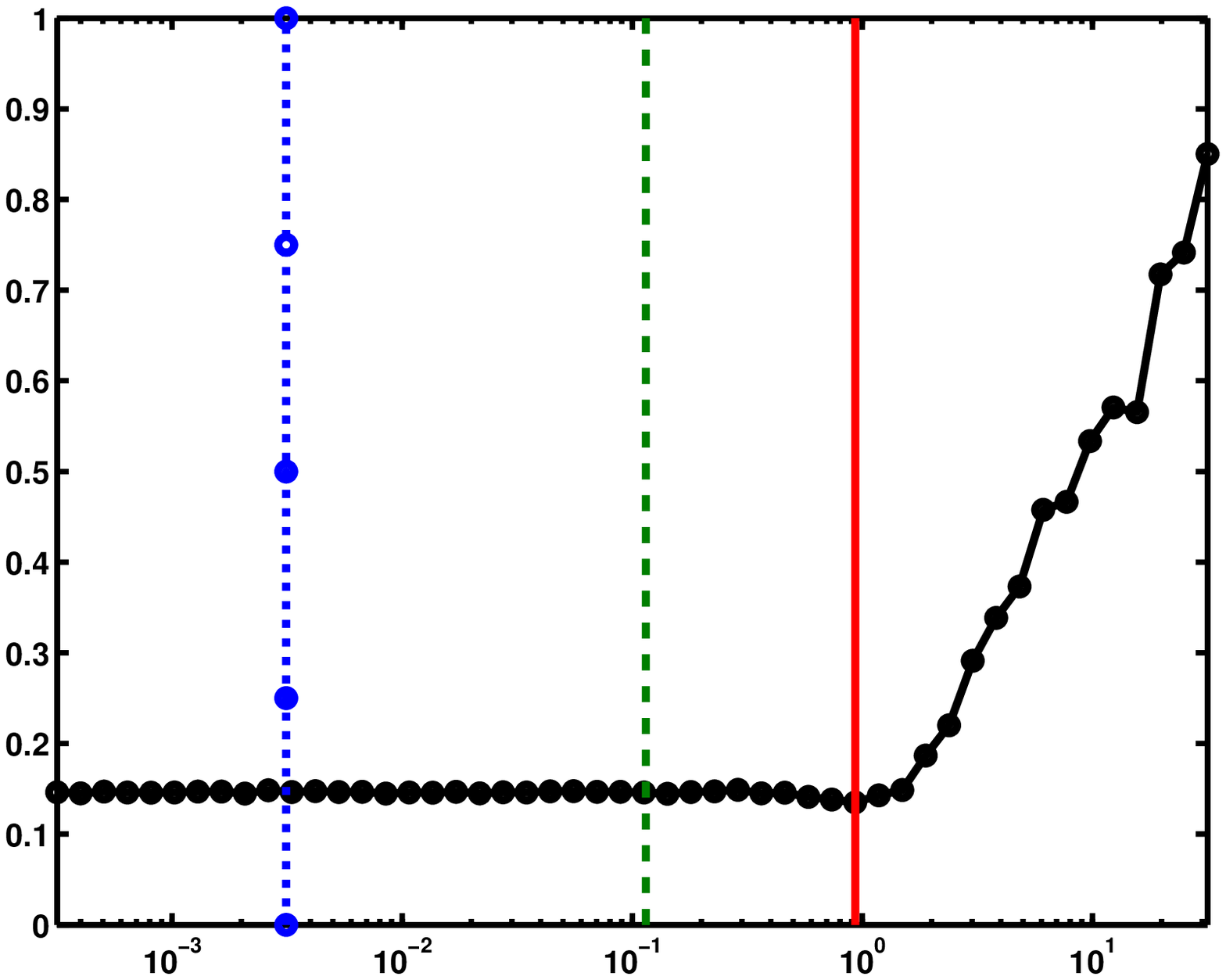}}
\subfigure[$L=I$]{\includegraphics[width=1.7in]{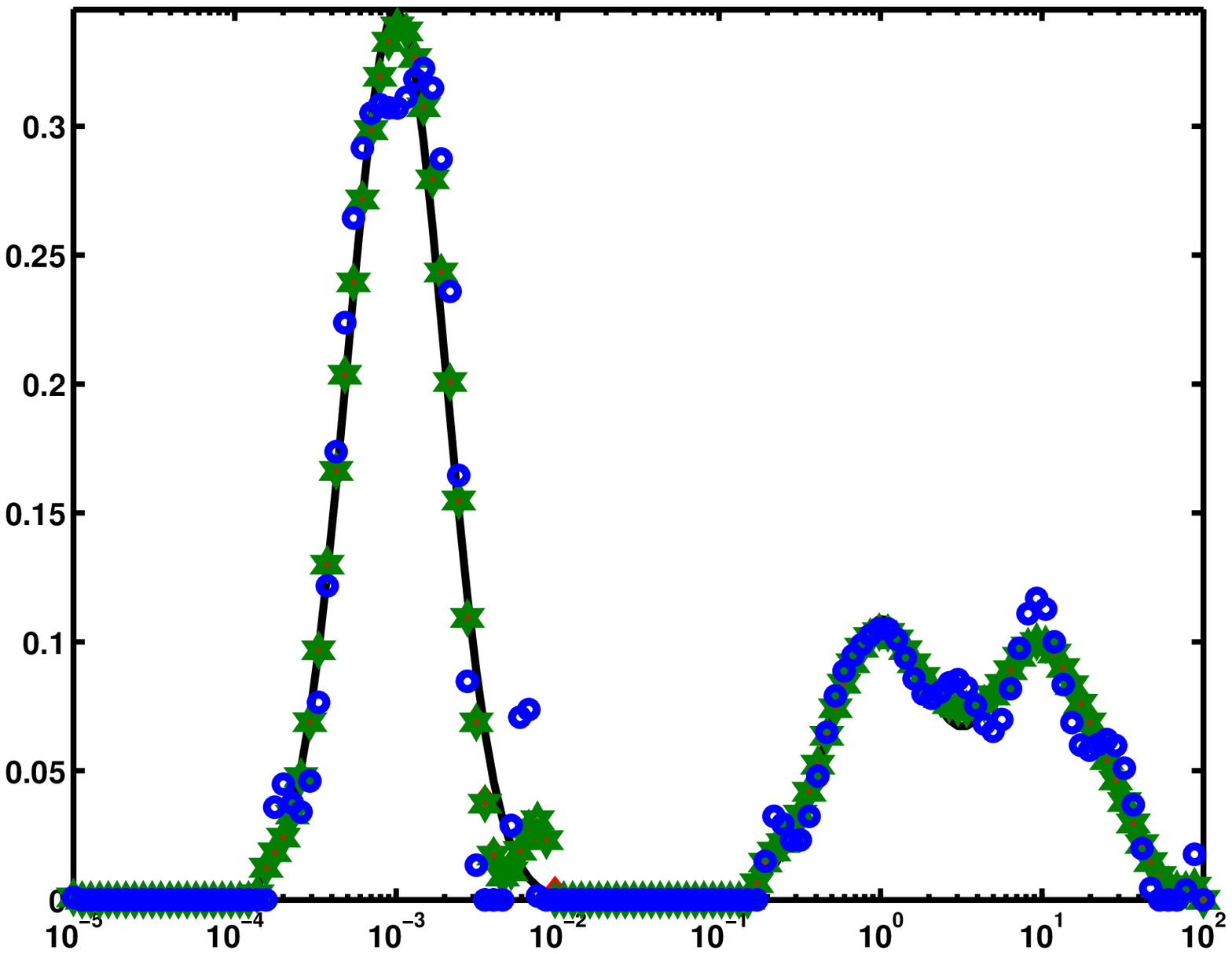}}
\subfigure[$L=L_1$]{\includegraphics[width=1.7in]{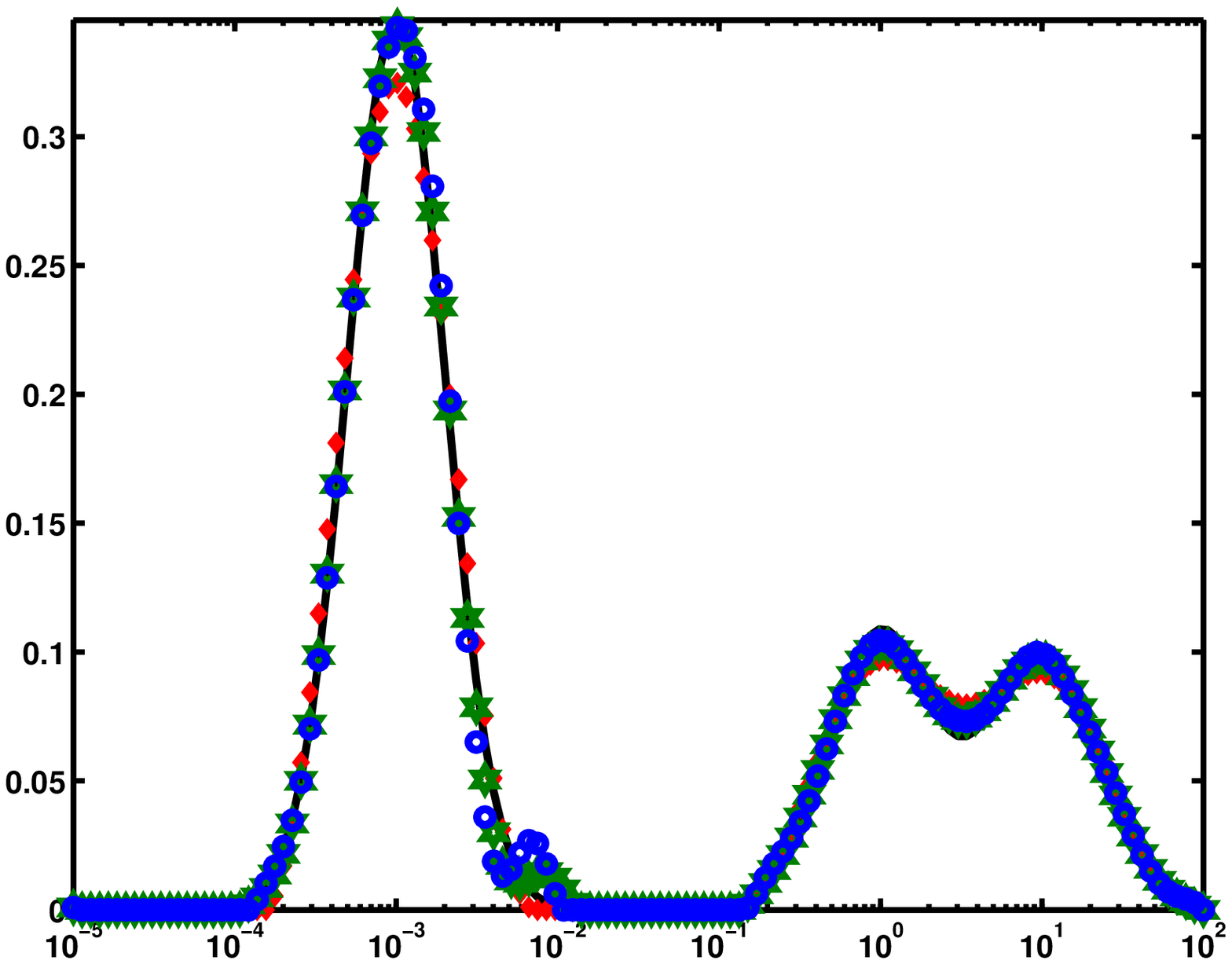}}
\subfigure[$L=L_2$]{\includegraphics[width=1.7in]{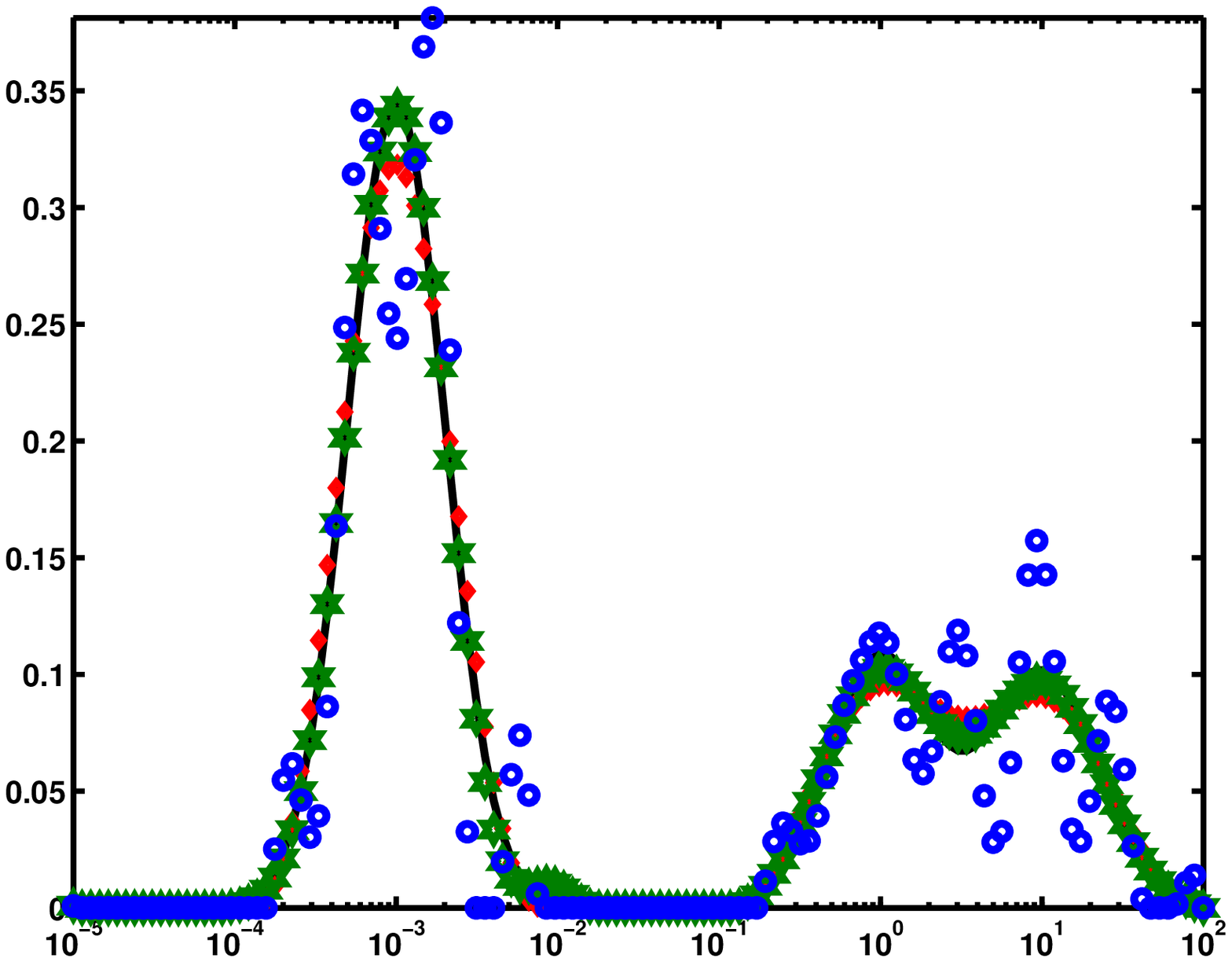}}
\caption{NNLS solutions of LN-C matrix $A_4$. Noise level $.1\%$ using the SBB algorithm.}
\label{lnfig-lambdachoiceLN6A4LNSBB}
\end{figure}

\begin{figure}[!ht]
\centering
\subfigure[$L=I$]{\includegraphics[width=1.7in]{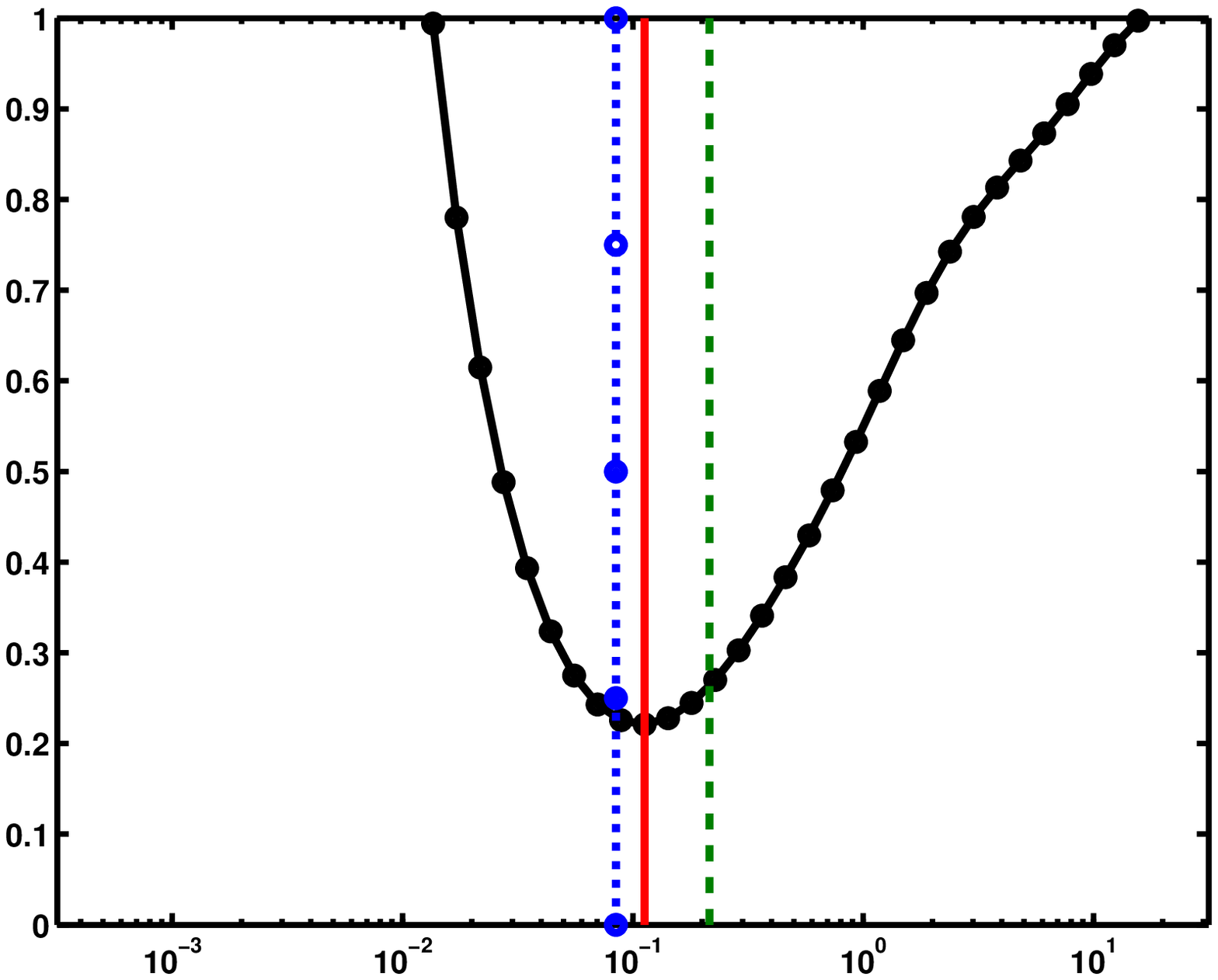}}
\subfigure[$L=L_1$]{\includegraphics[width=1.7in]{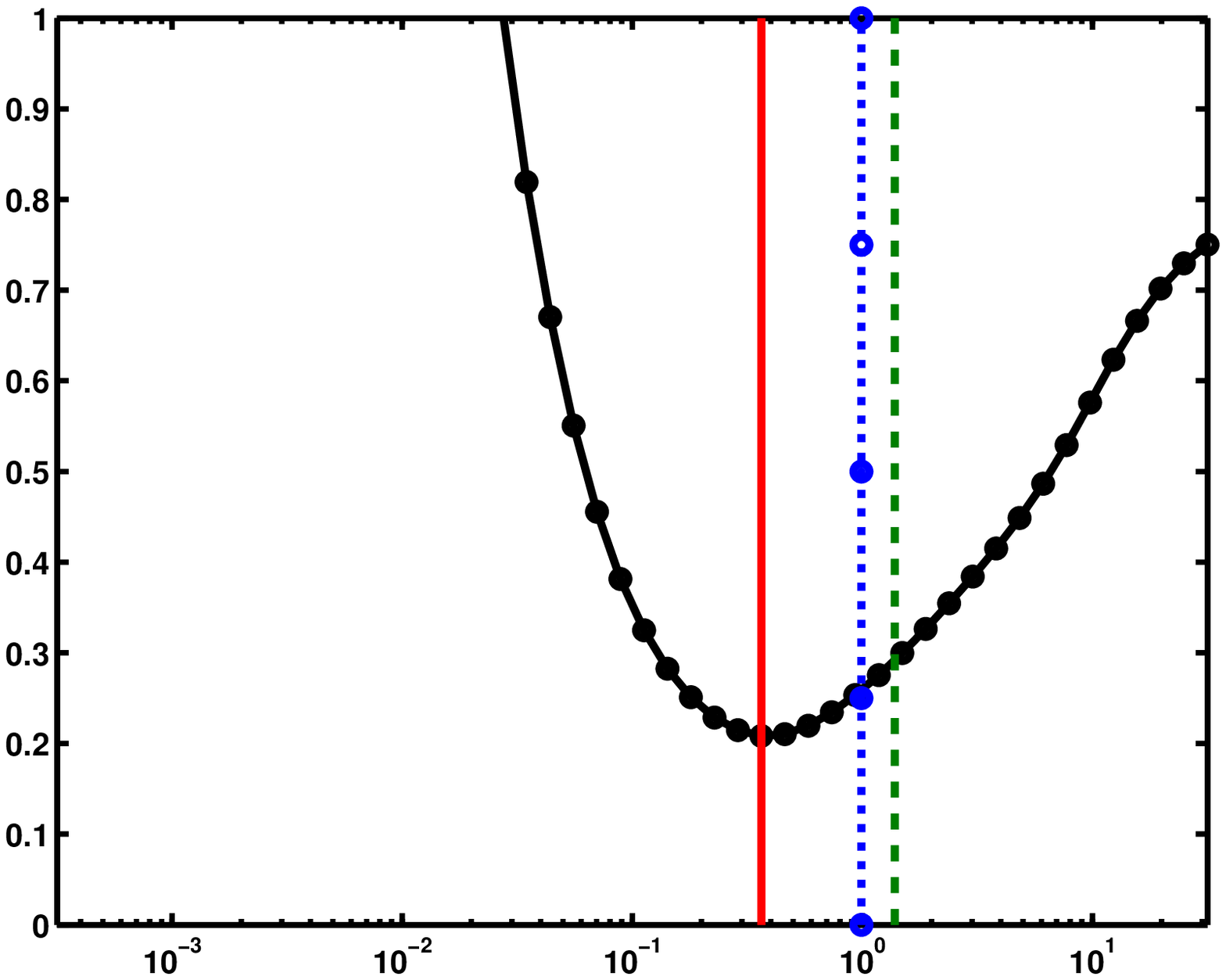}}
\subfigure[$L=L_2$]{\includegraphics[width=1.7in]{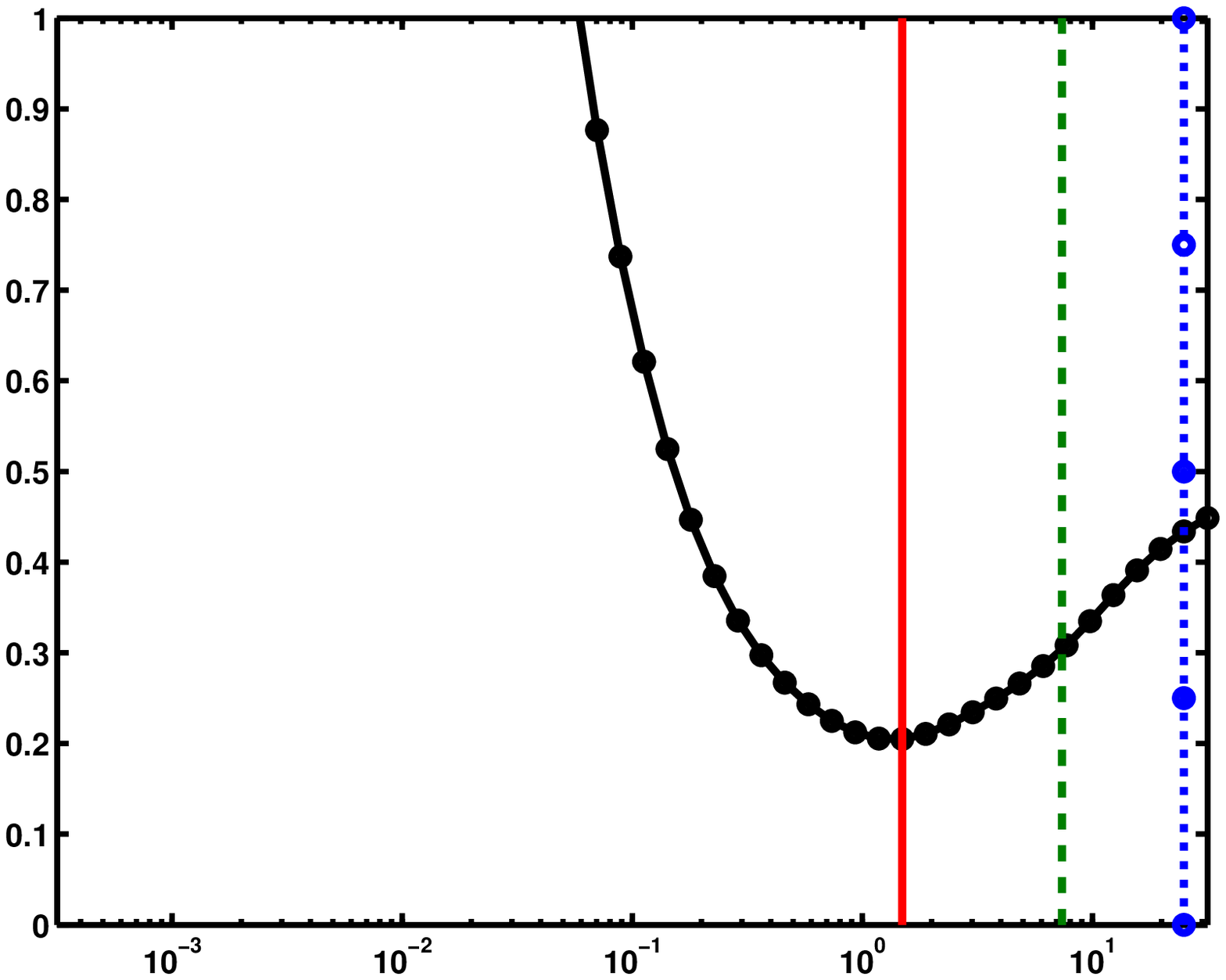}}
\subfigure[$L=I$]{\includegraphics[width=1.7in]{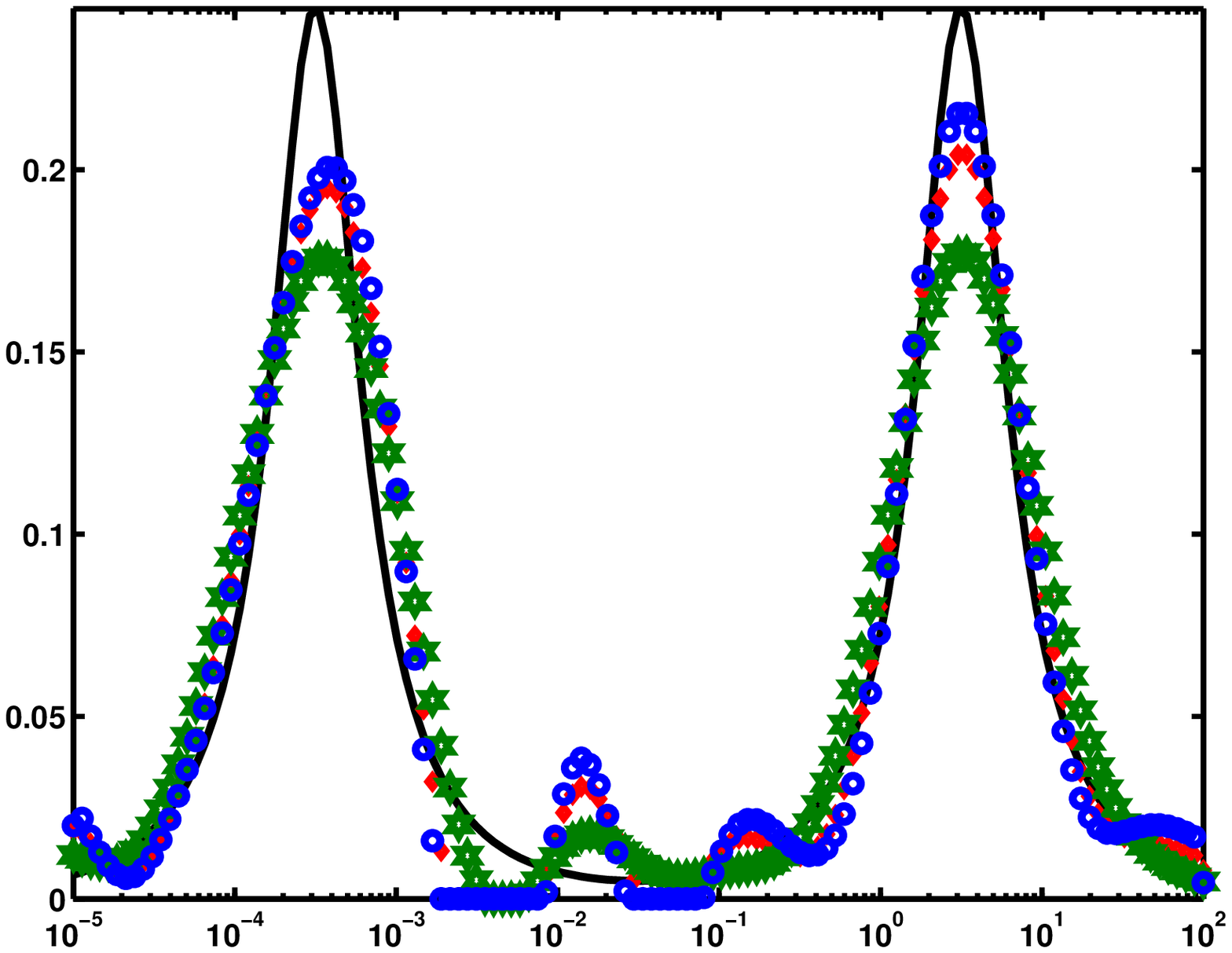}}
\subfigure[$L=L_1$]{\includegraphics[width=1.7in]{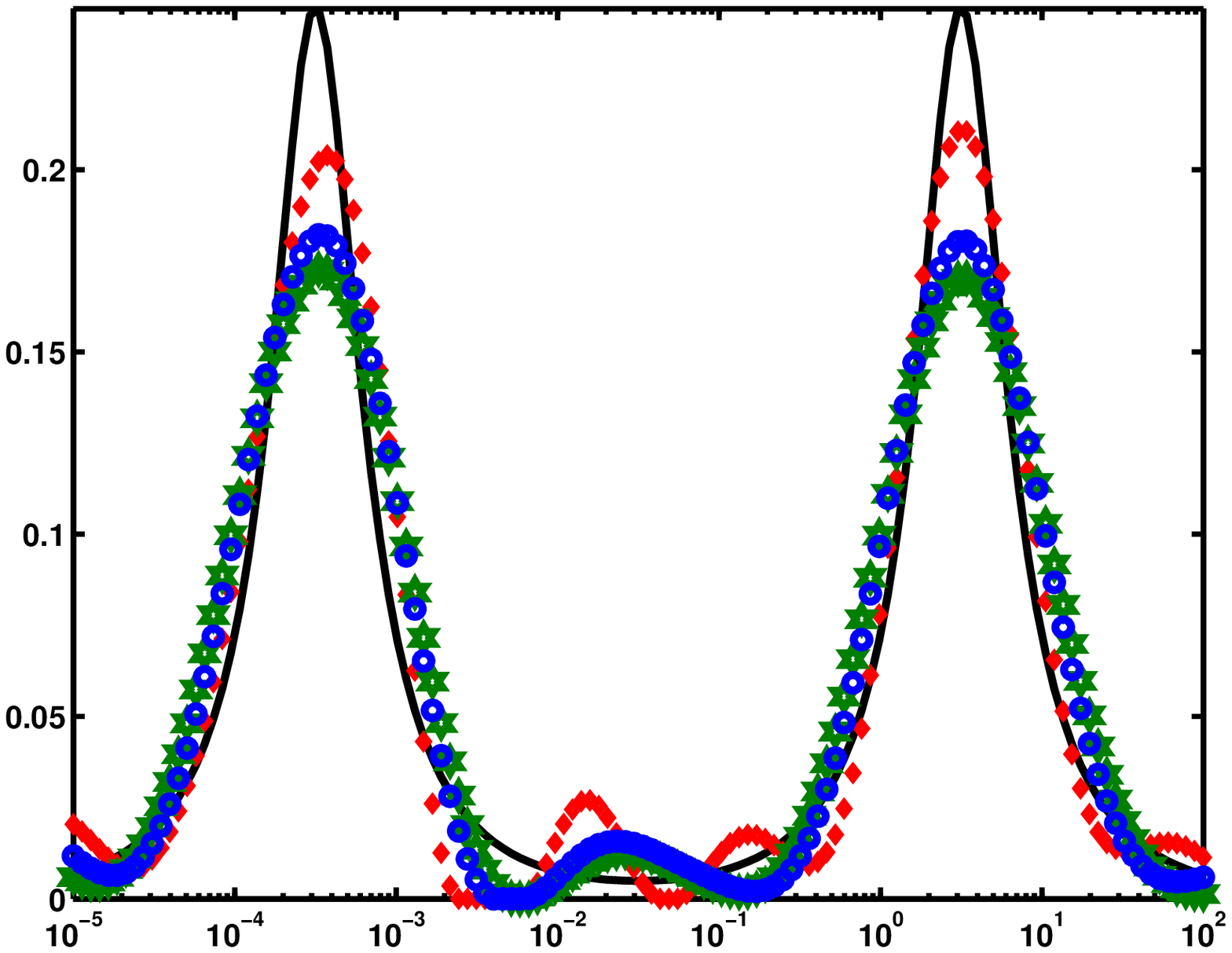}}
\subfigure[$L=L_2$]{\includegraphics[width=1.7in]{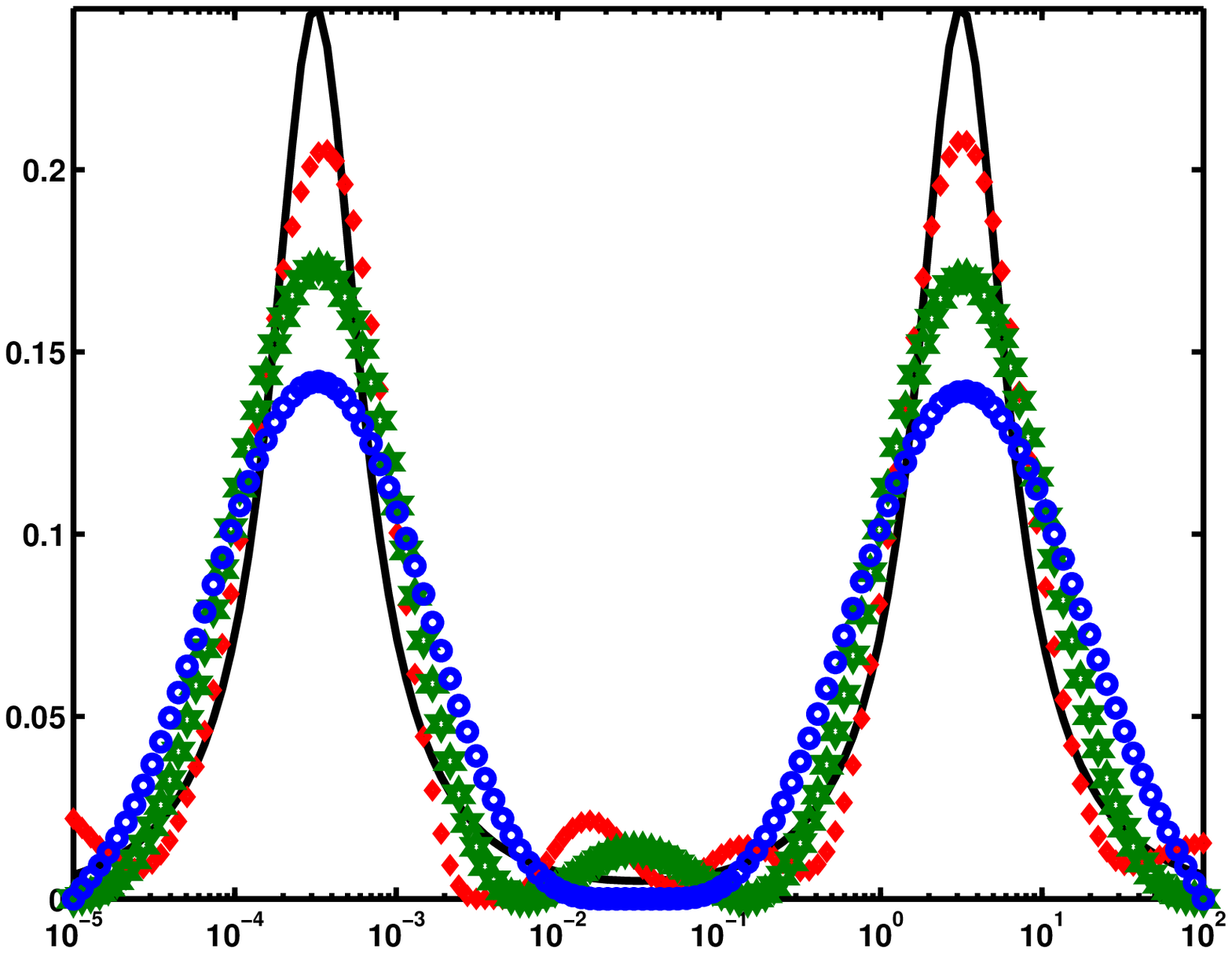}}
\caption{NNLS solutions of RQ-A matrix $A_4$. Noise level $1\%$. Method CVX}
\label{hnfig-lambdachoiceRQ1A4HNCVX}
\end{figure}

\begin{figure}[!ht]
\centering
\subfigure[$L=I$]{\includegraphics[width=1.7in]{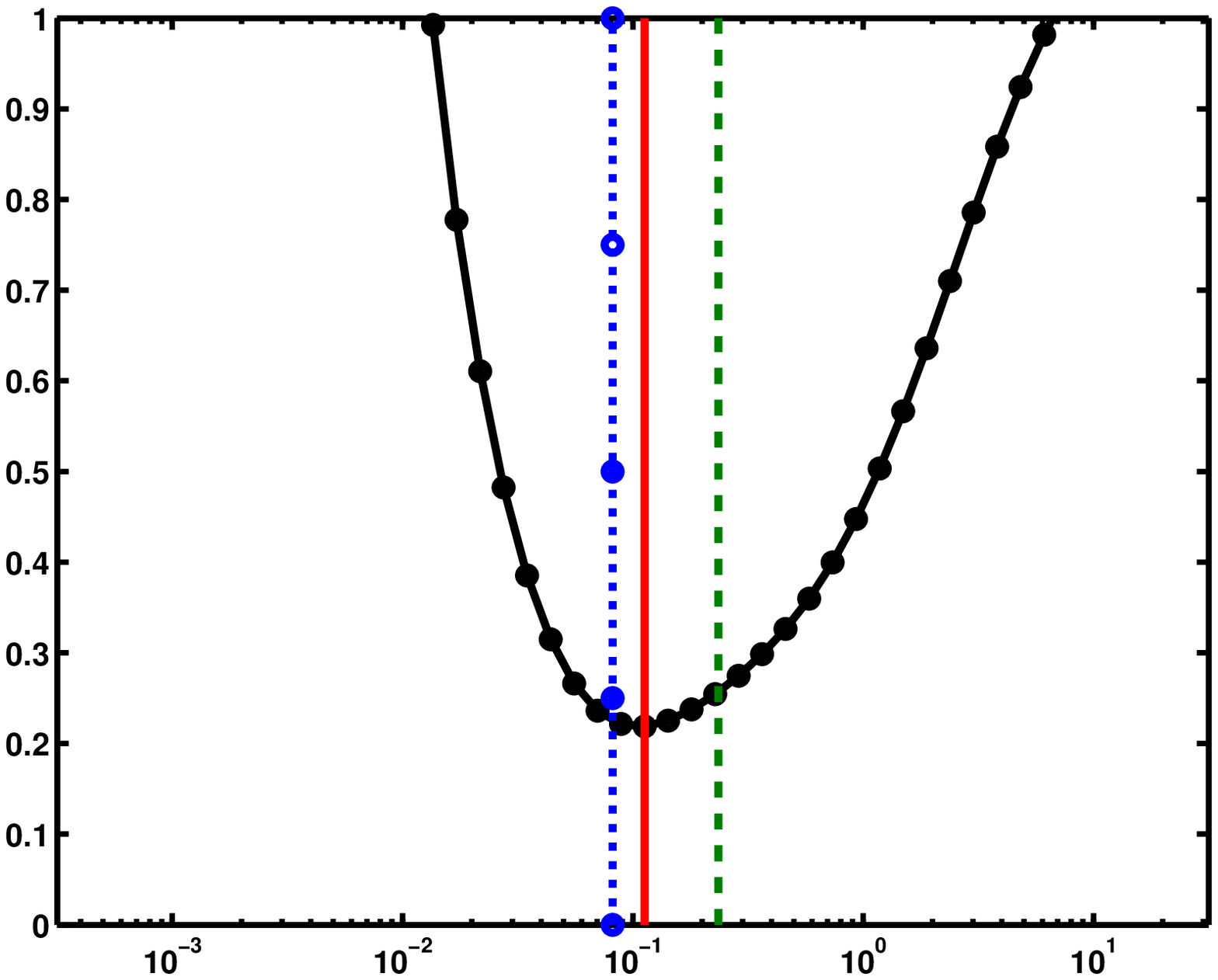}}
\subfigure[$L=L_1$]{\includegraphics[width=1.7in]{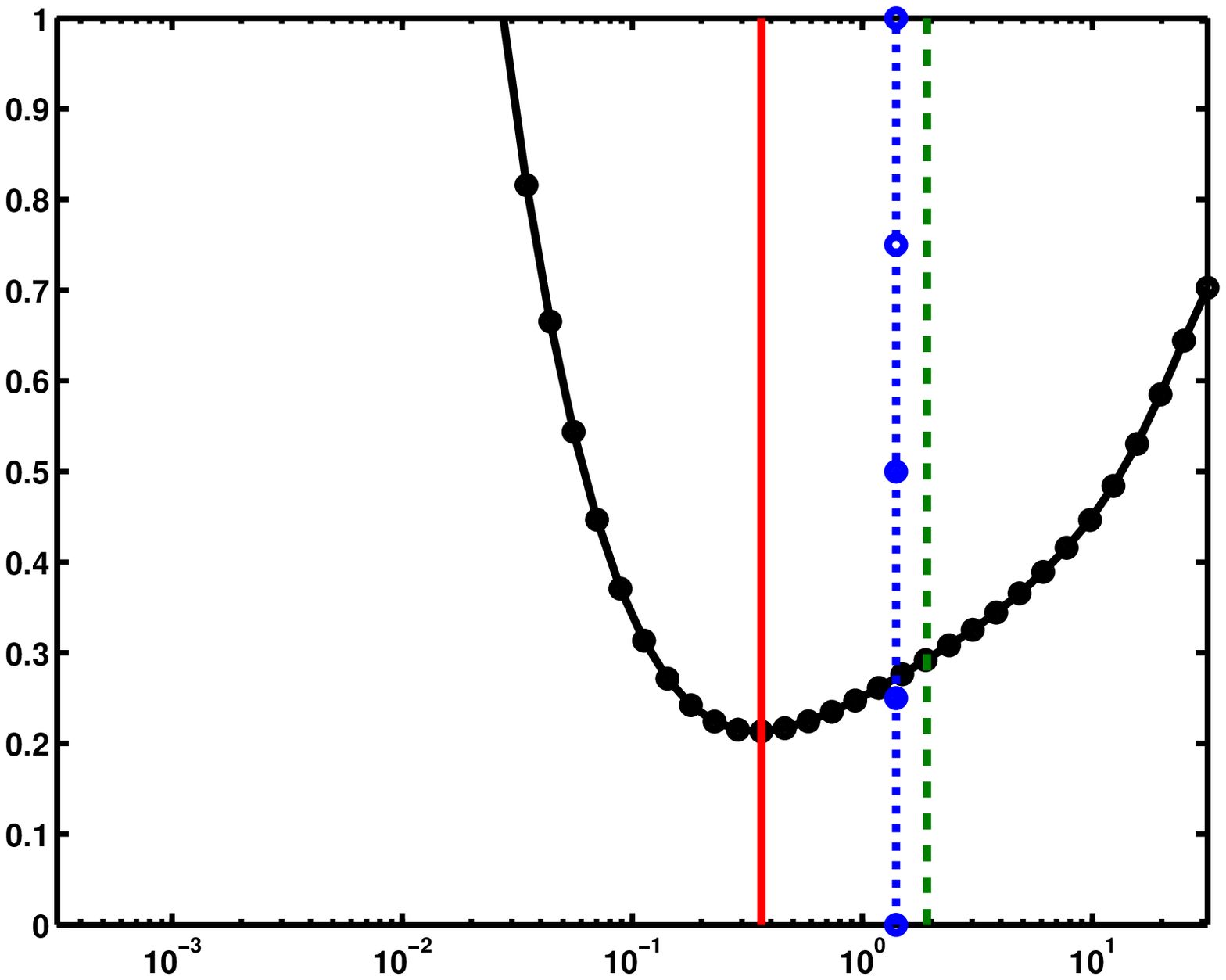}}
\subfigure[$L=L_2$]{\includegraphics[width=1.7in]{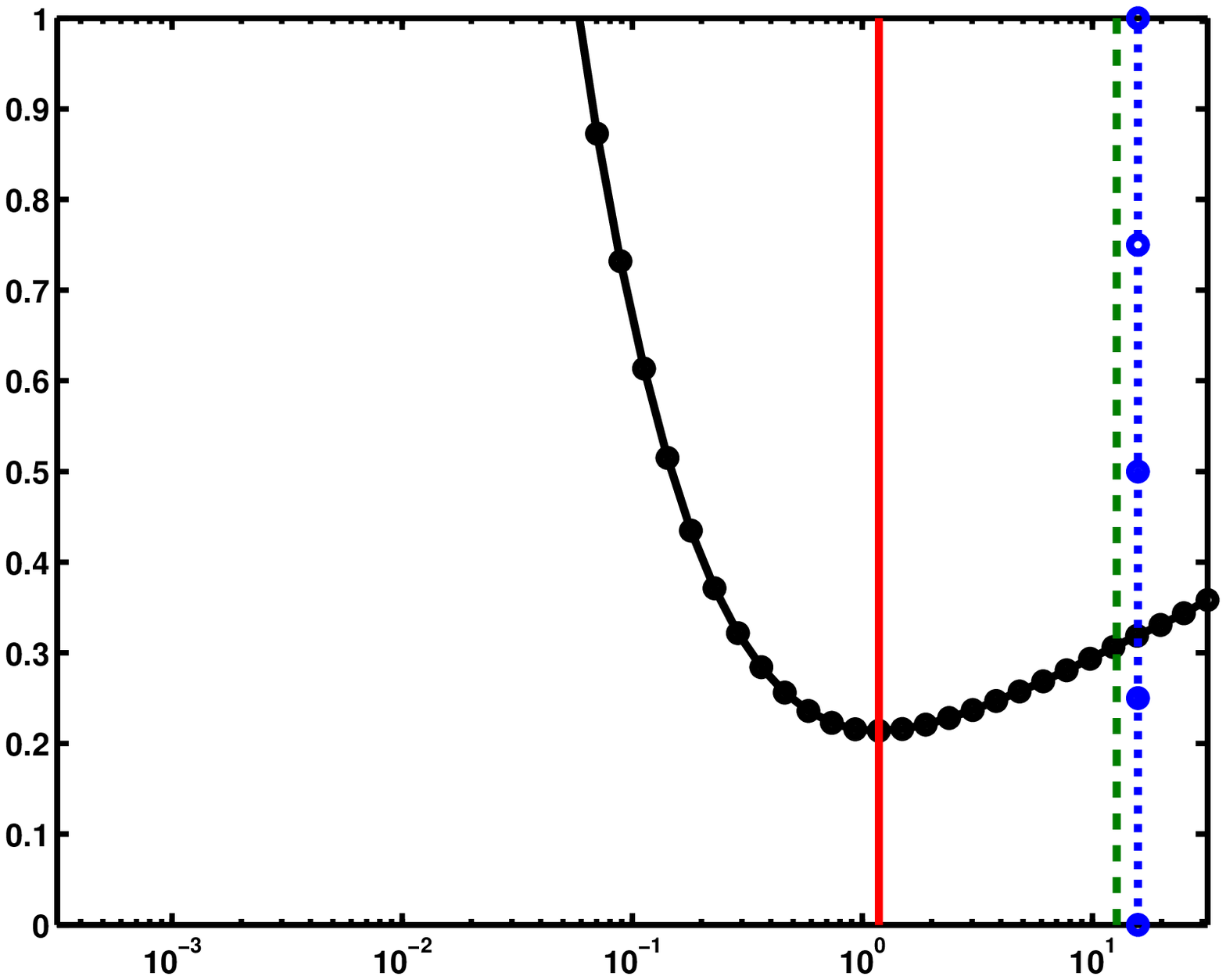}}
\subfigure[$L=I$]{\includegraphics[width=1.7in]{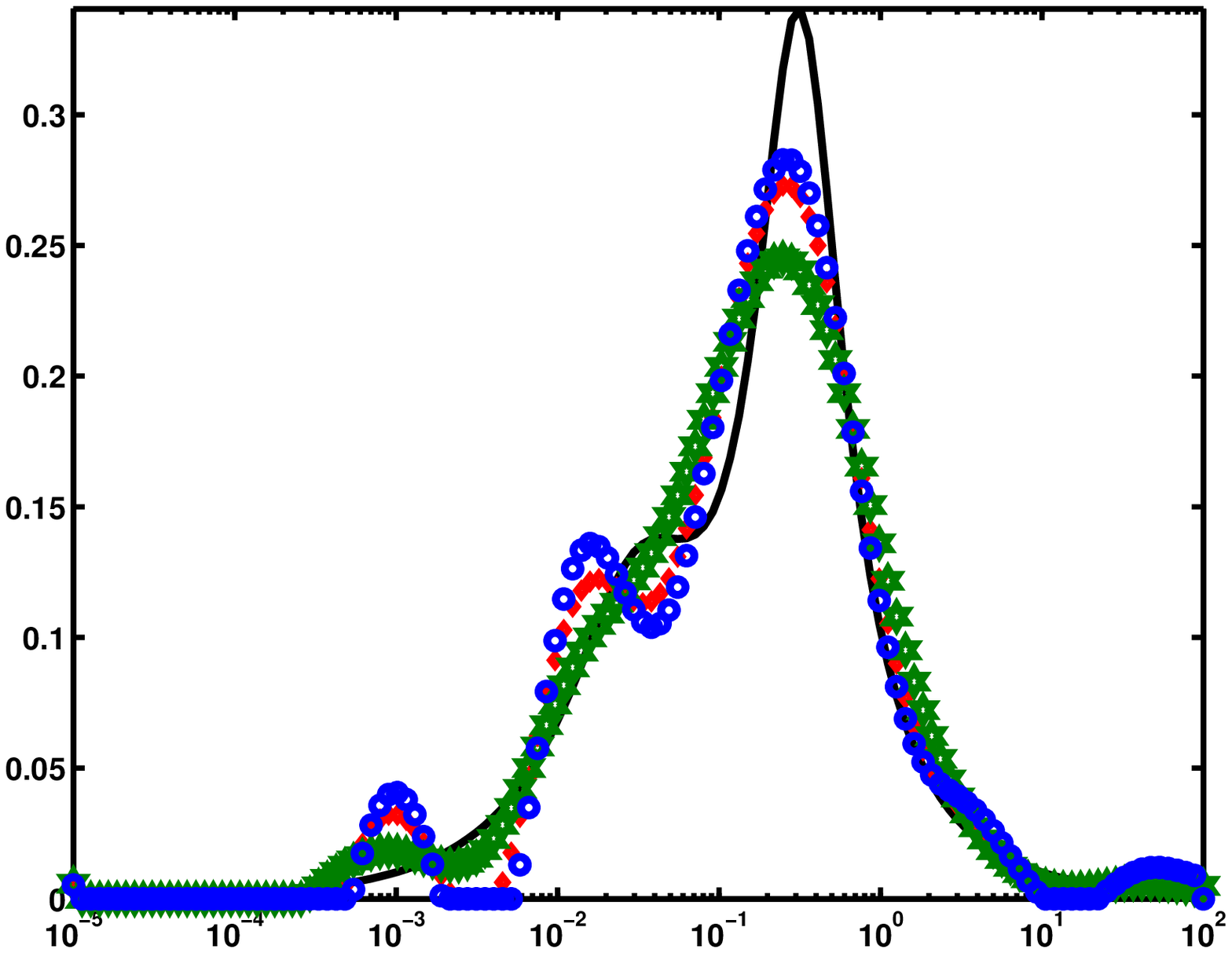}}
\subfigure[$L=L_1$]{\includegraphics[width=1.7in]{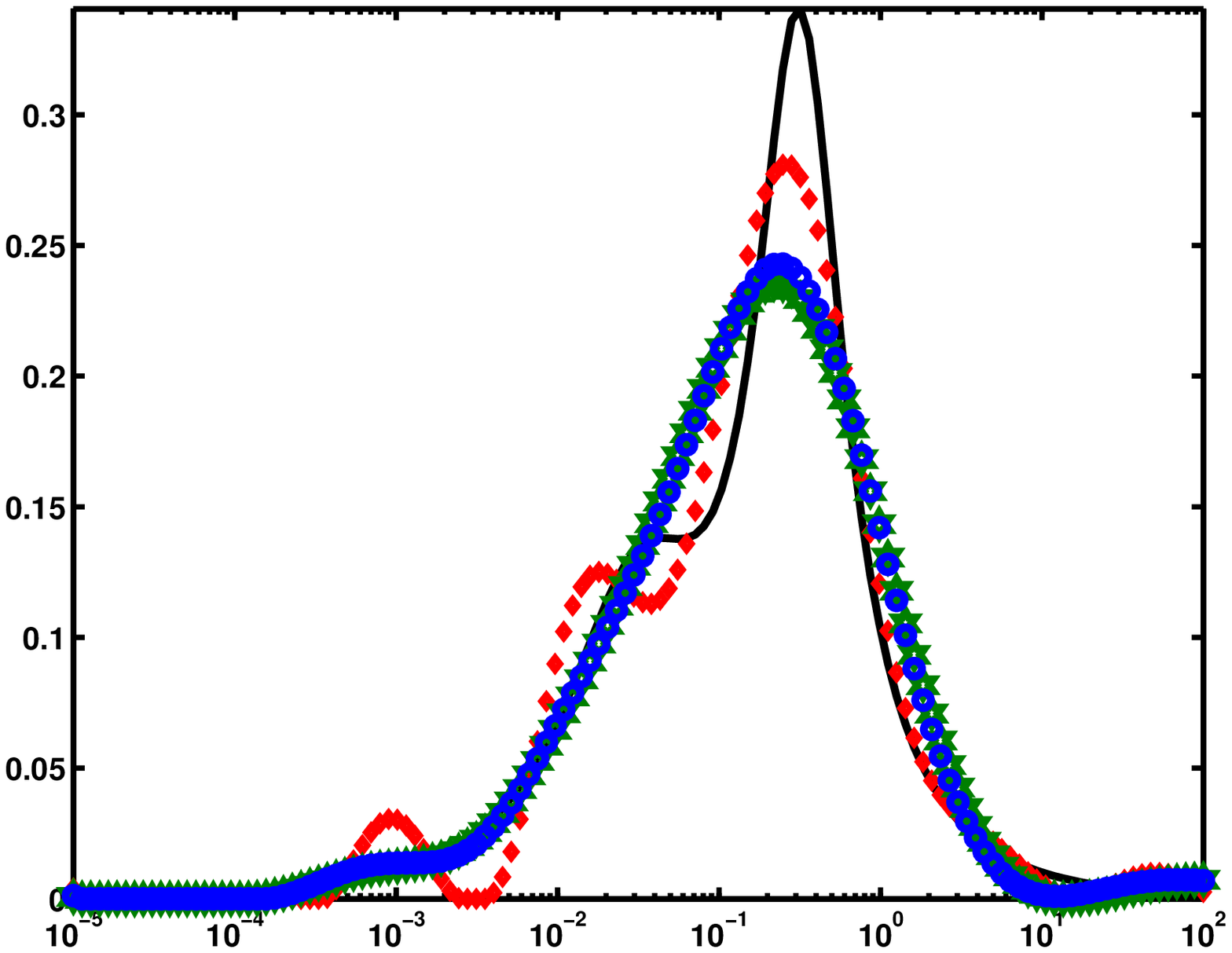}}
\subfigure[$L=L_2$]{\includegraphics[width=1.7in]{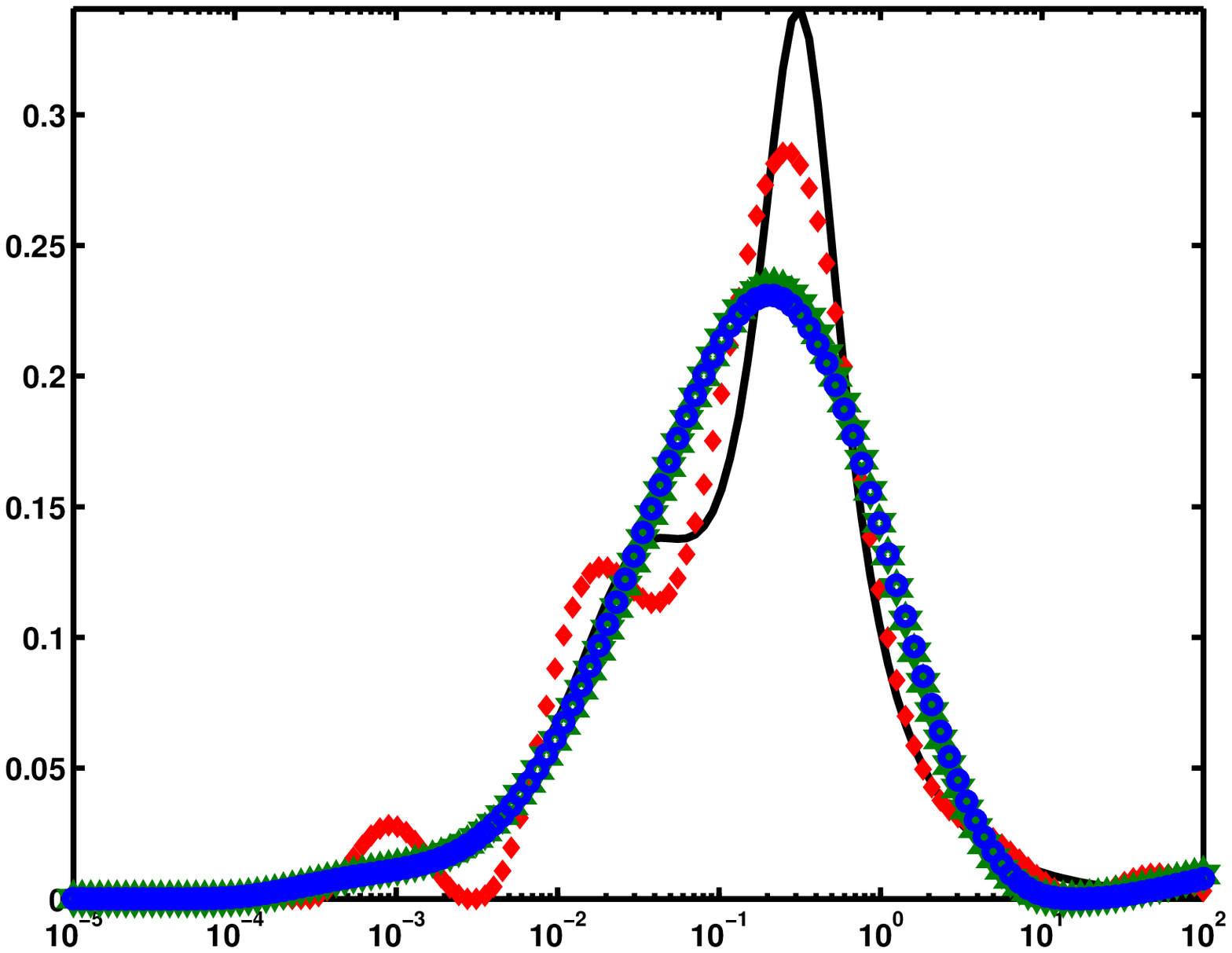}}
\caption{NNLS solutions of RQ-B matrix $A_4$. Noise level $1\%$.Method CVX}
\label{hnfig-lambdachoiceRQ5A4HNCVX}
\end{figure}

\begin{figure}[!ht]
\centering
\subfigure[$L=I$]{\includegraphics[width=1.7in]{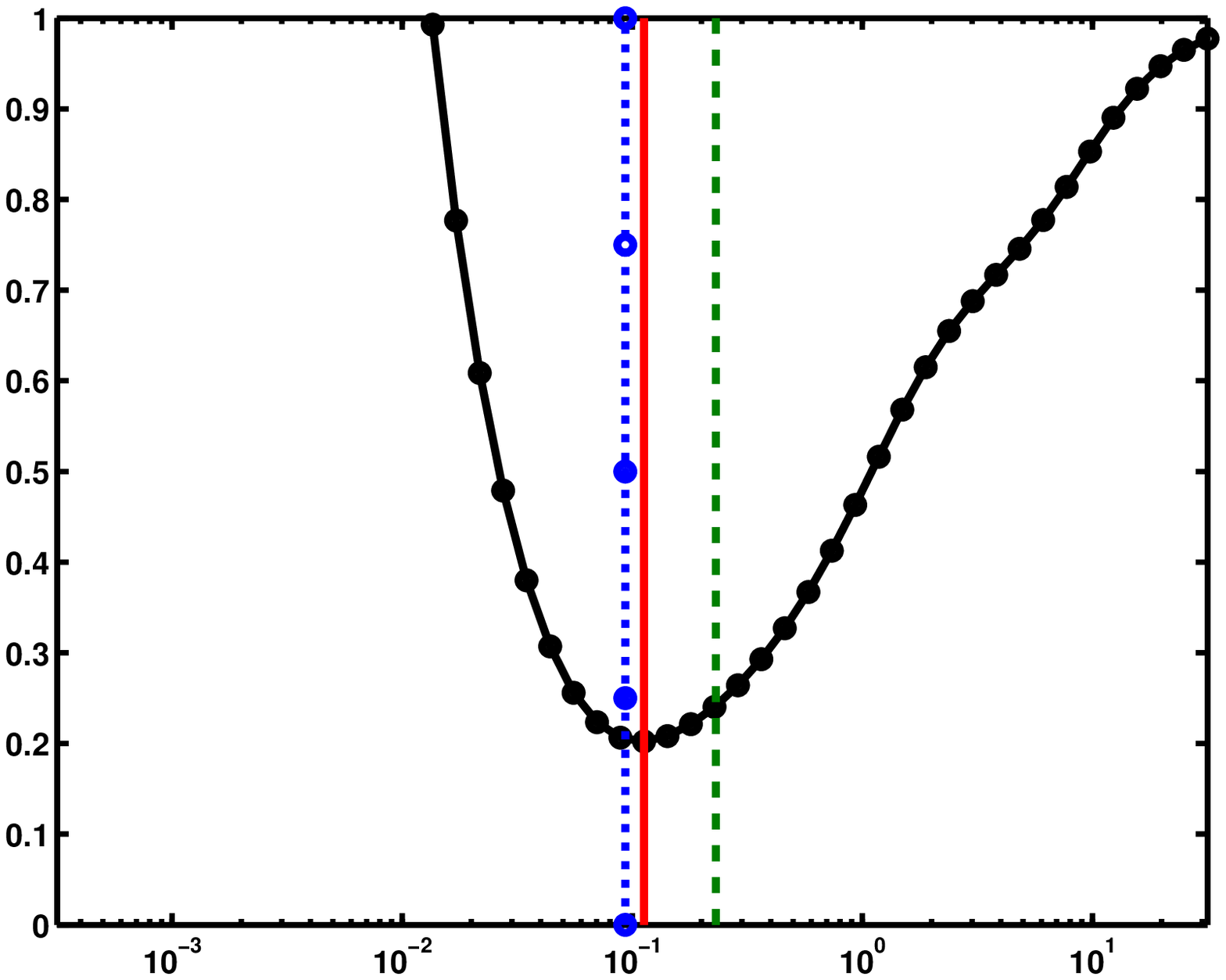}}
\subfigure[$L=L_1$]{\includegraphics[width=1.7in]{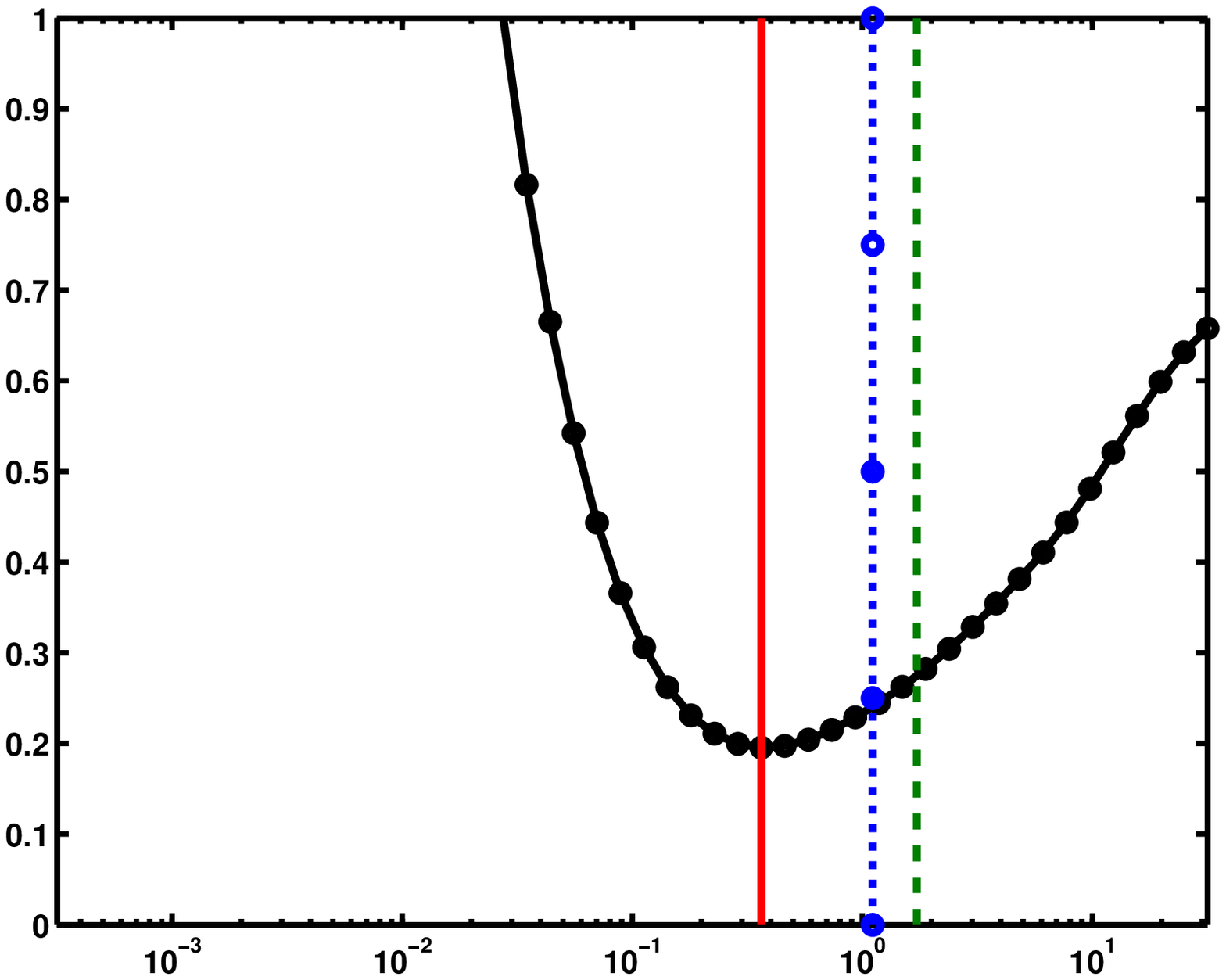}}
\subfigure[$L=L_2$]{\includegraphics[width=1.7in]{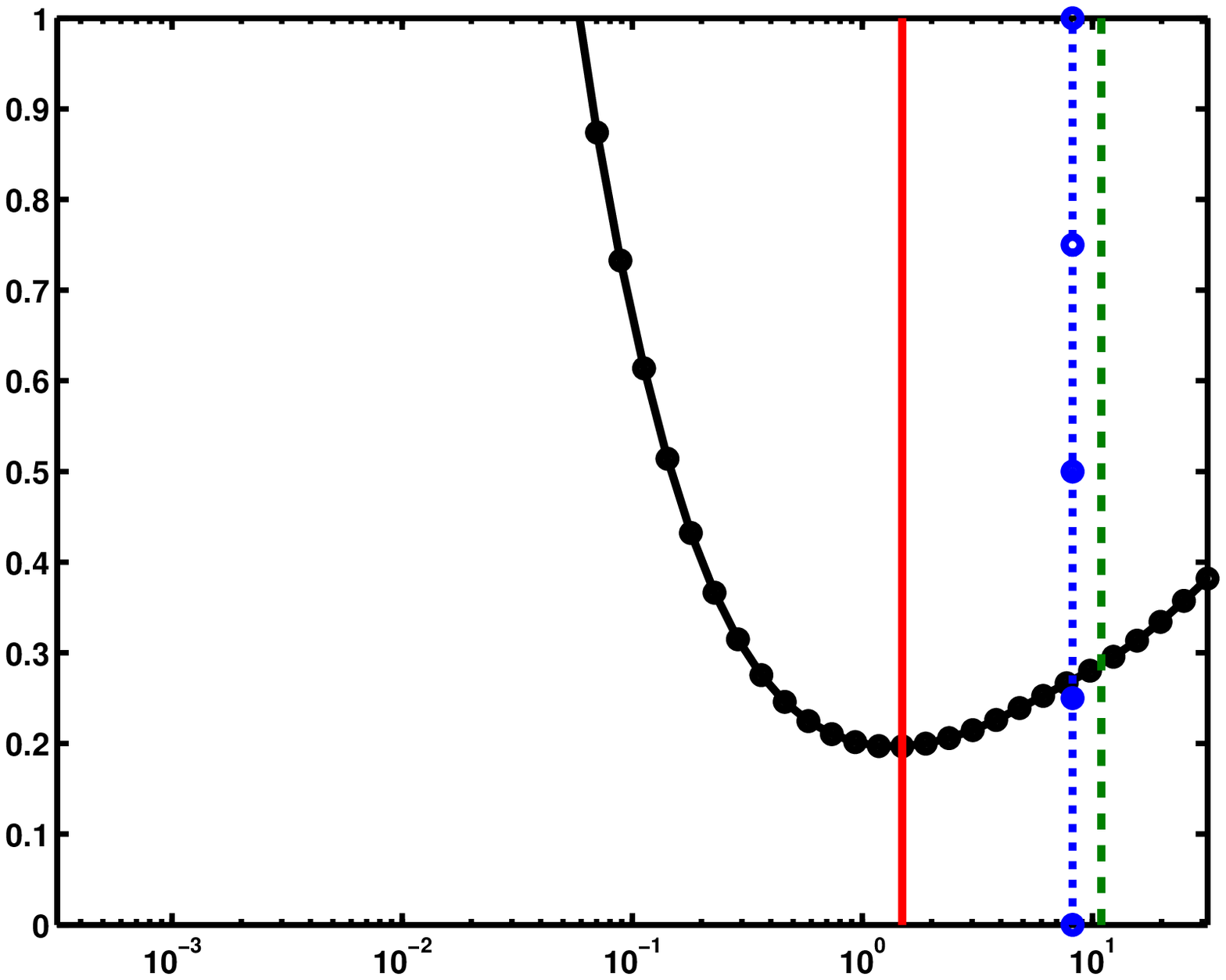}}
\subfigure[$L=I$]{\includegraphics[width=1.7in]{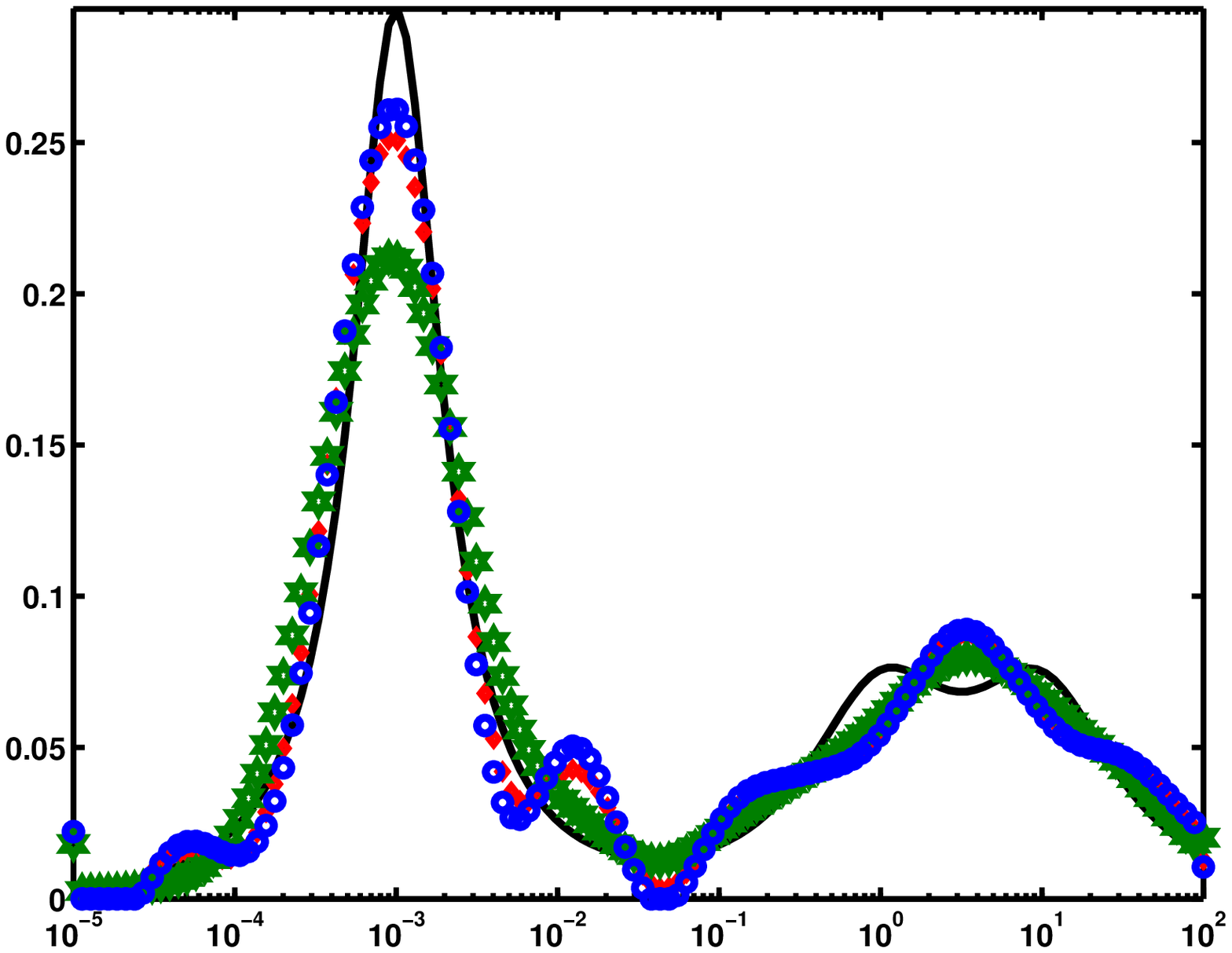}}
\subfigure[$L=L_1$]{\includegraphics[width=1.7in]{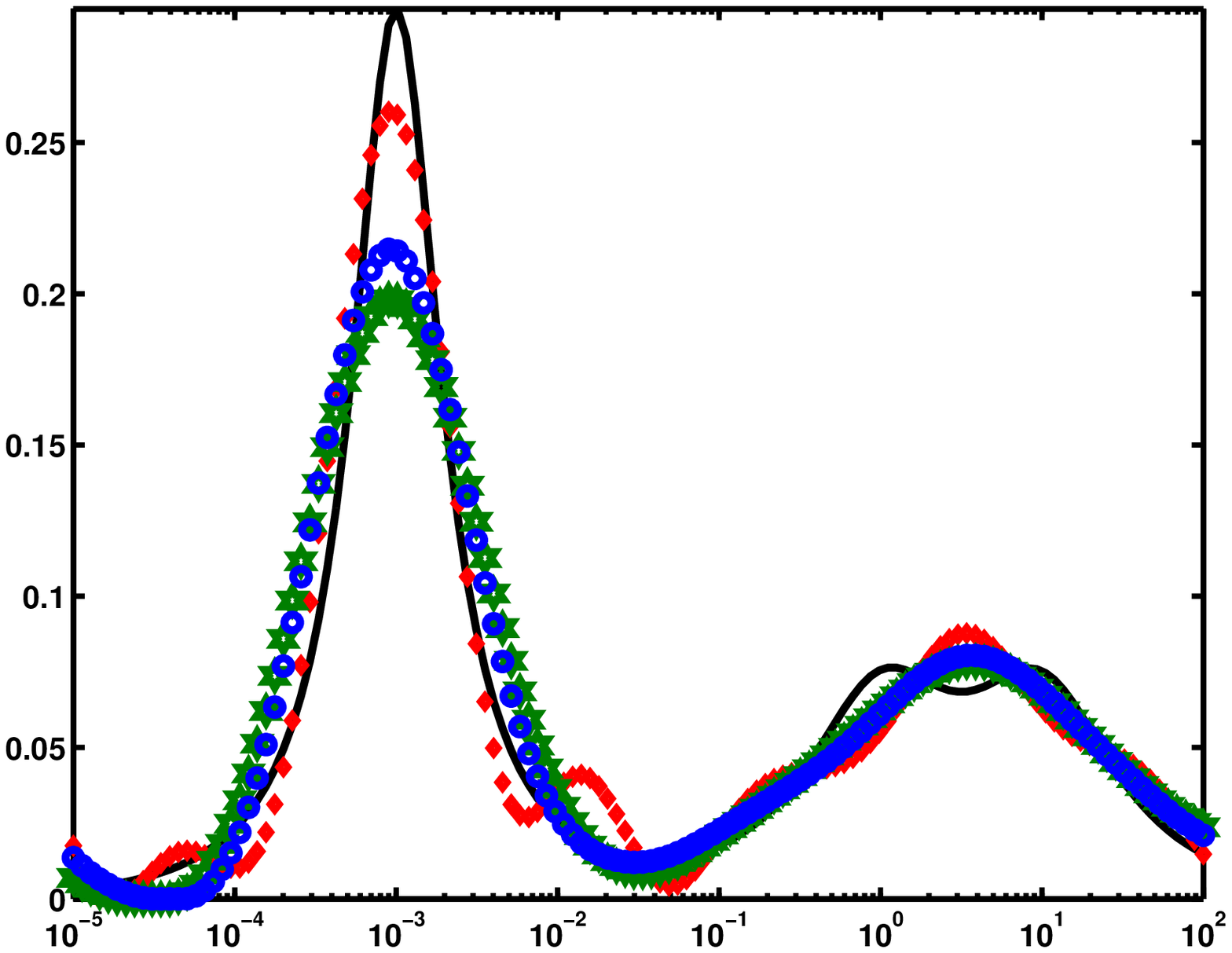}}
\subfigure[$L=L_2$]{\includegraphics[width=1.7in]{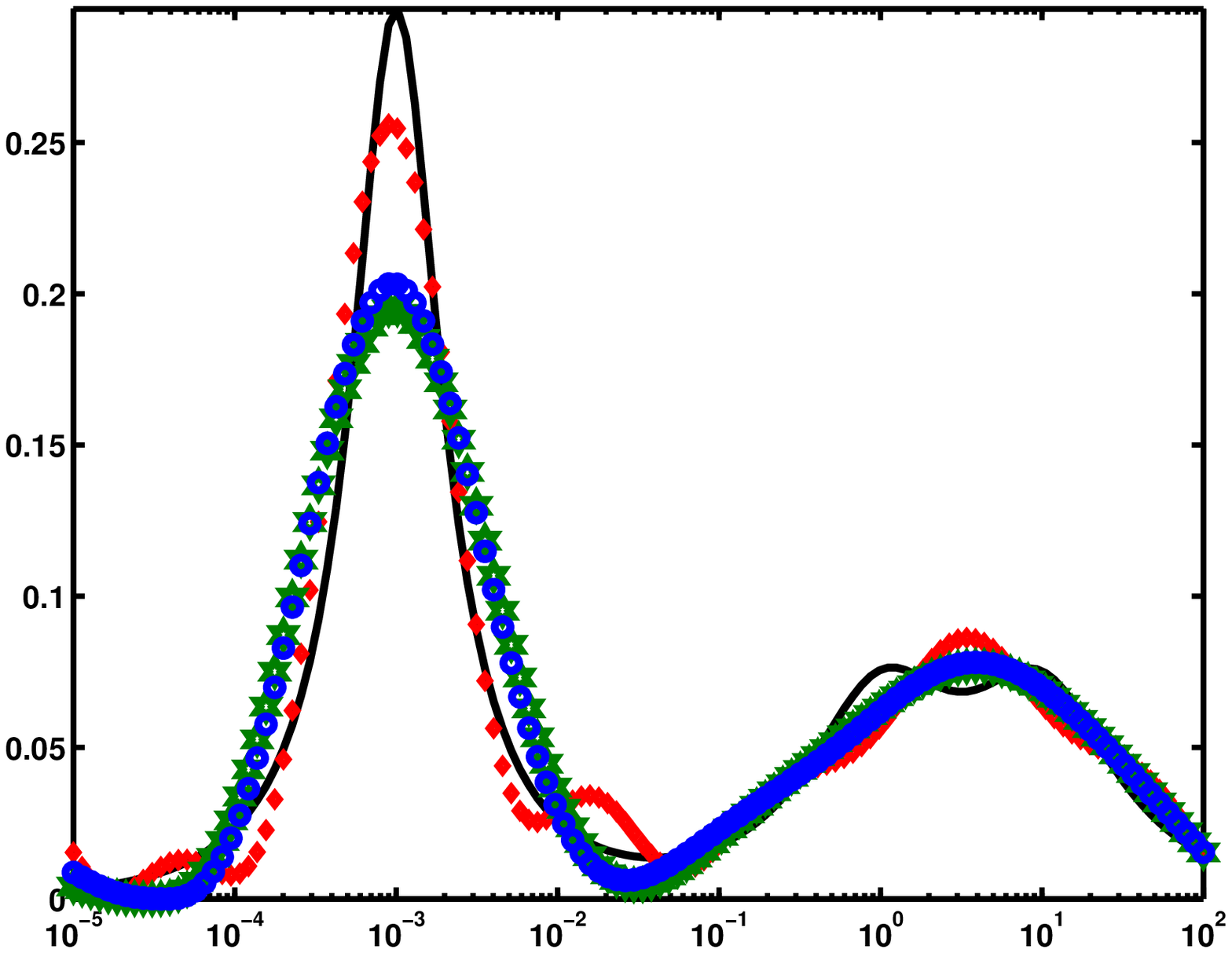}}
\caption{NNLS solutions of RQ-C matrix $A_4$. Noise level $1\%$.Method CVX}
\label{hnfig-lambdachoiceRQ6A4HNCVX}
\end{figure}

\begin{figure}[!ht]
\centering
\subfigure[$L=I$]{\includegraphics[width=1.7in]{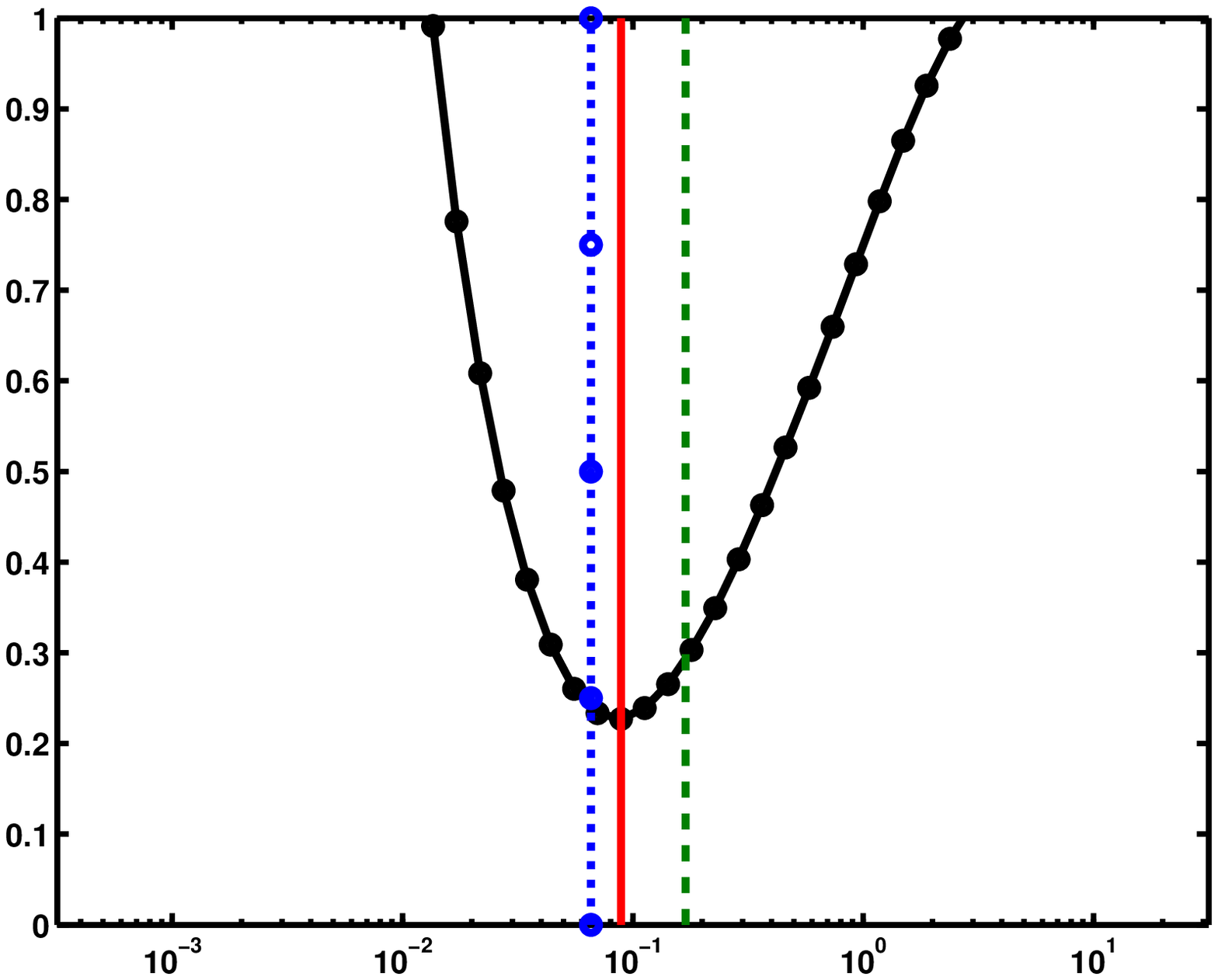}}
\subfigure[$L=L_1$]{\includegraphics[width=1.7in]{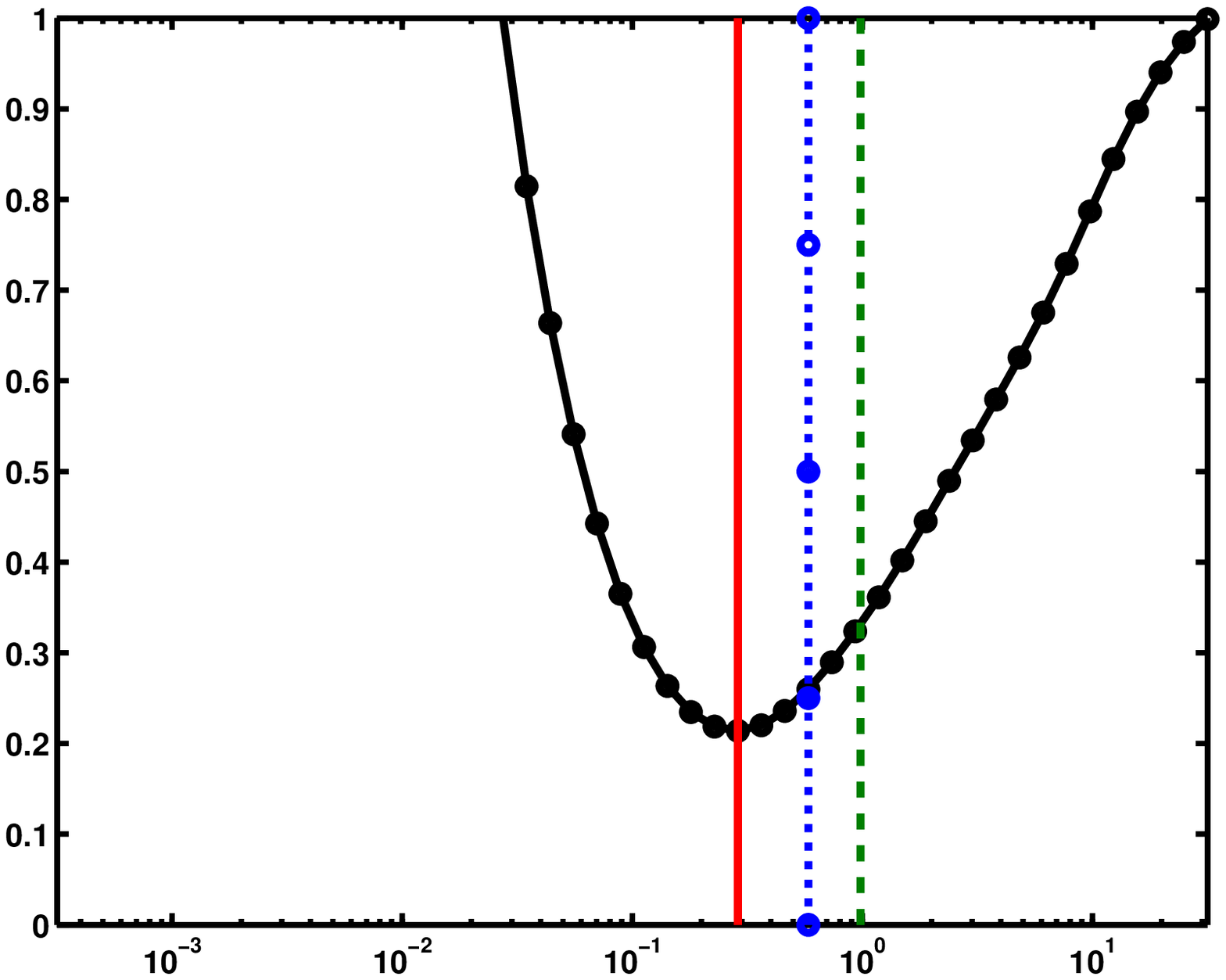}}
\subfigure[$L=L_2$]{\includegraphics[width=1.7in]{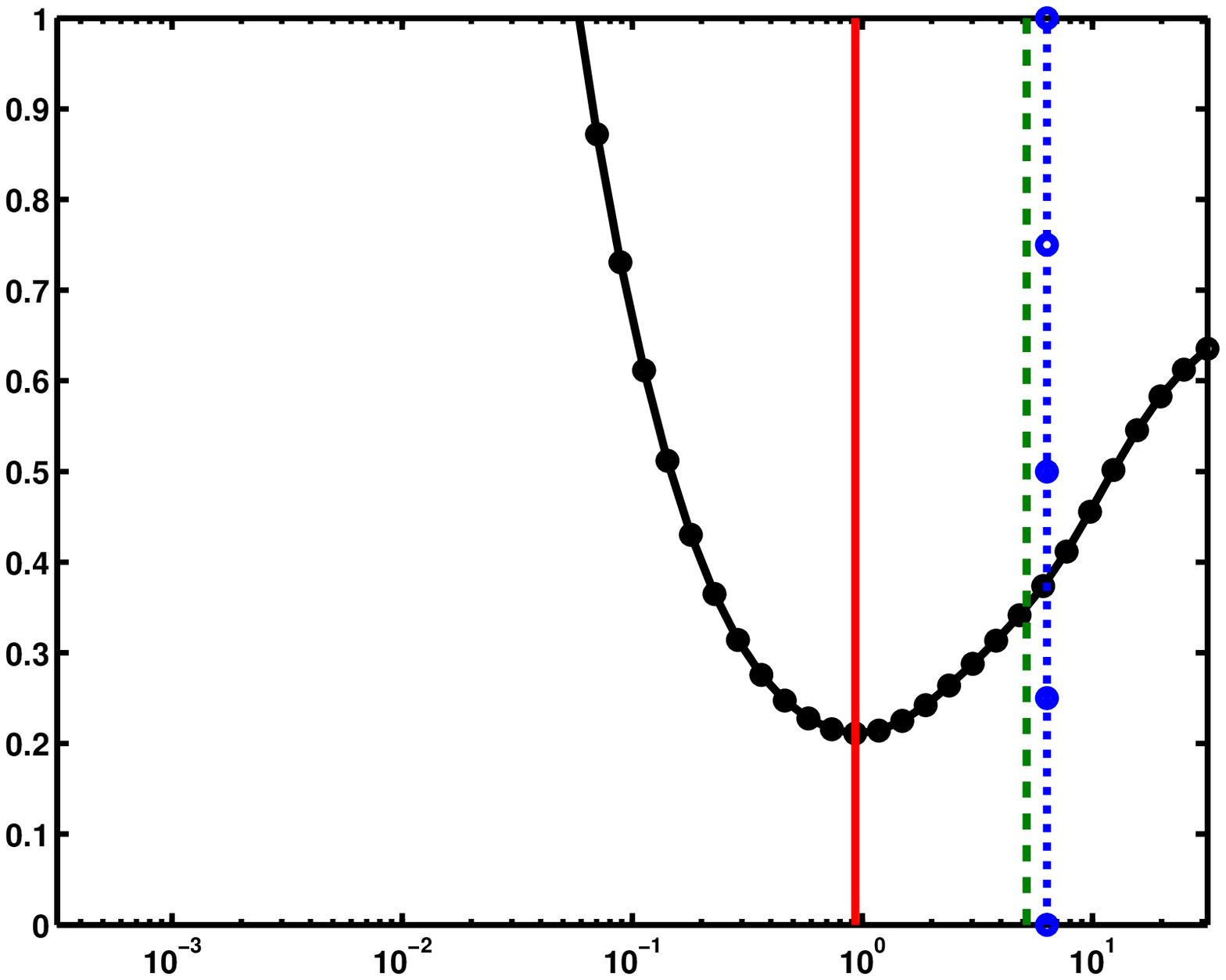}}
\subfigure[$L=I$]{\includegraphics[width=1.7in]{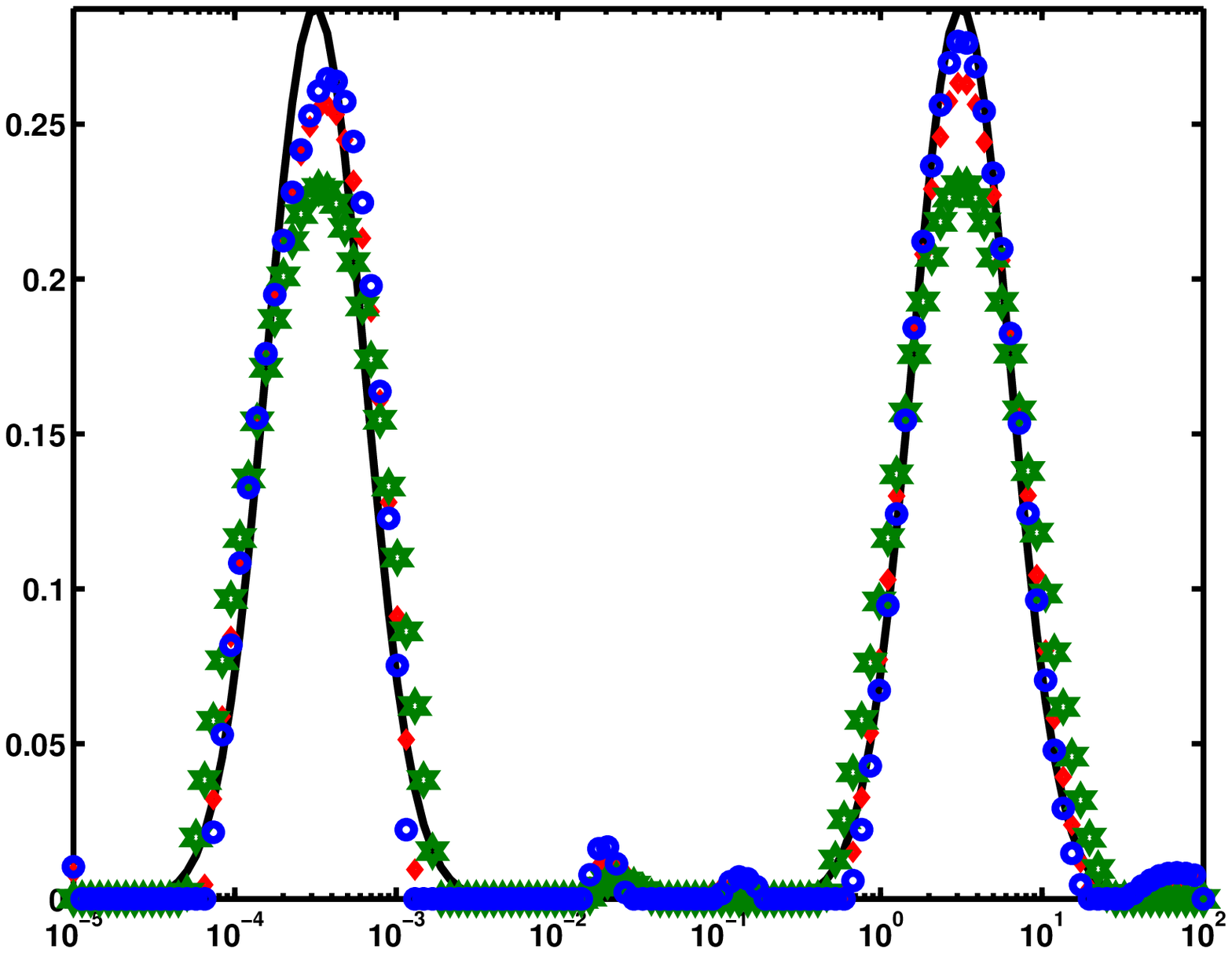}}
\subfigure[$L=L_1$]{\includegraphics[width=1.7in]{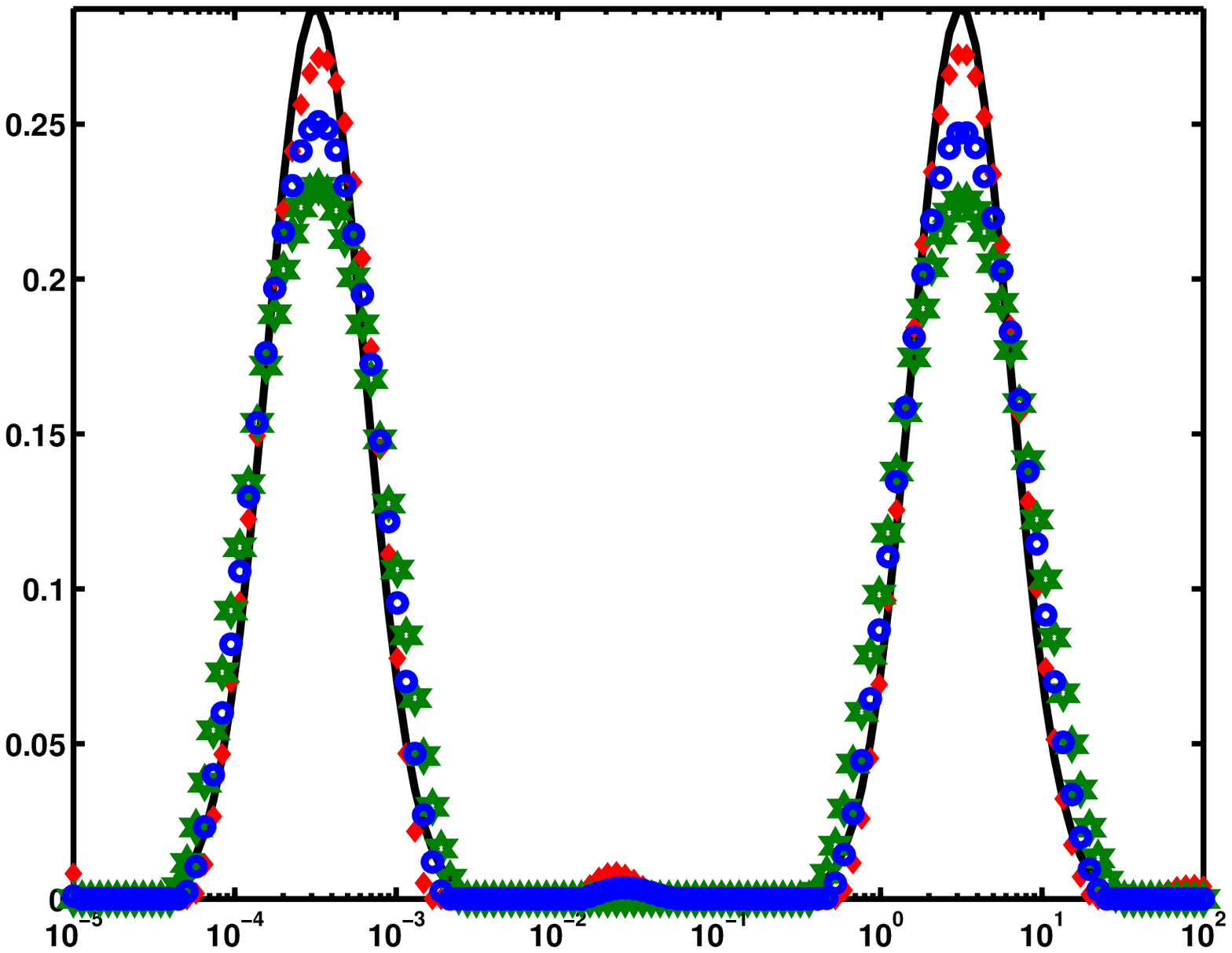}}
\subfigure[$L=L_2$]{\includegraphics[width=1.7in]{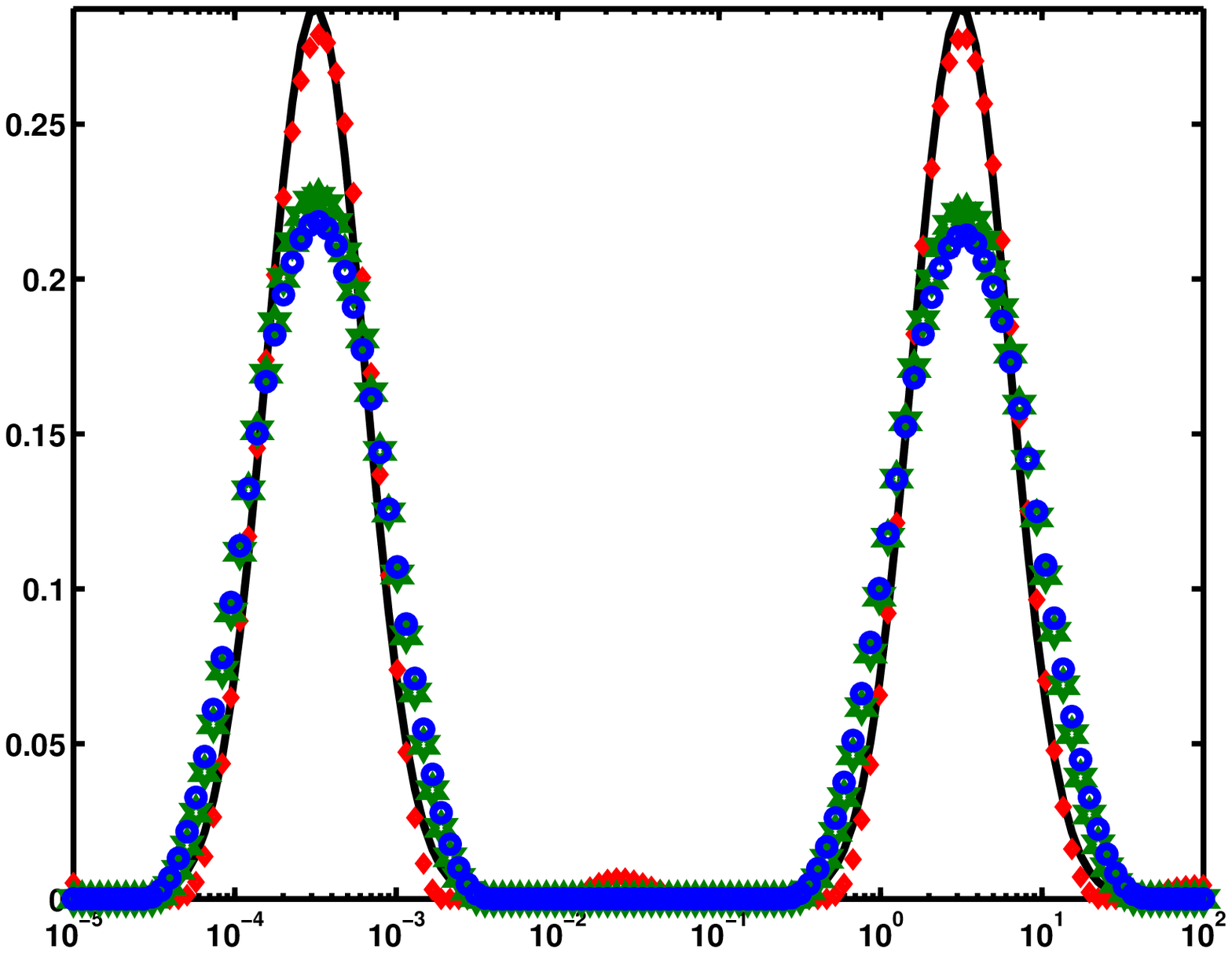}}
\caption{NNLS solutions of LN-A matrix $A_4$. Noise level $1\%$.Method CVX}
\label{hnfig-lambdachoiceLN2A4HNCVX}
\end{figure}

\begin{figure}[!ht]
\centering
\subfigure[$L=I$]{\includegraphics[width=1.7in]{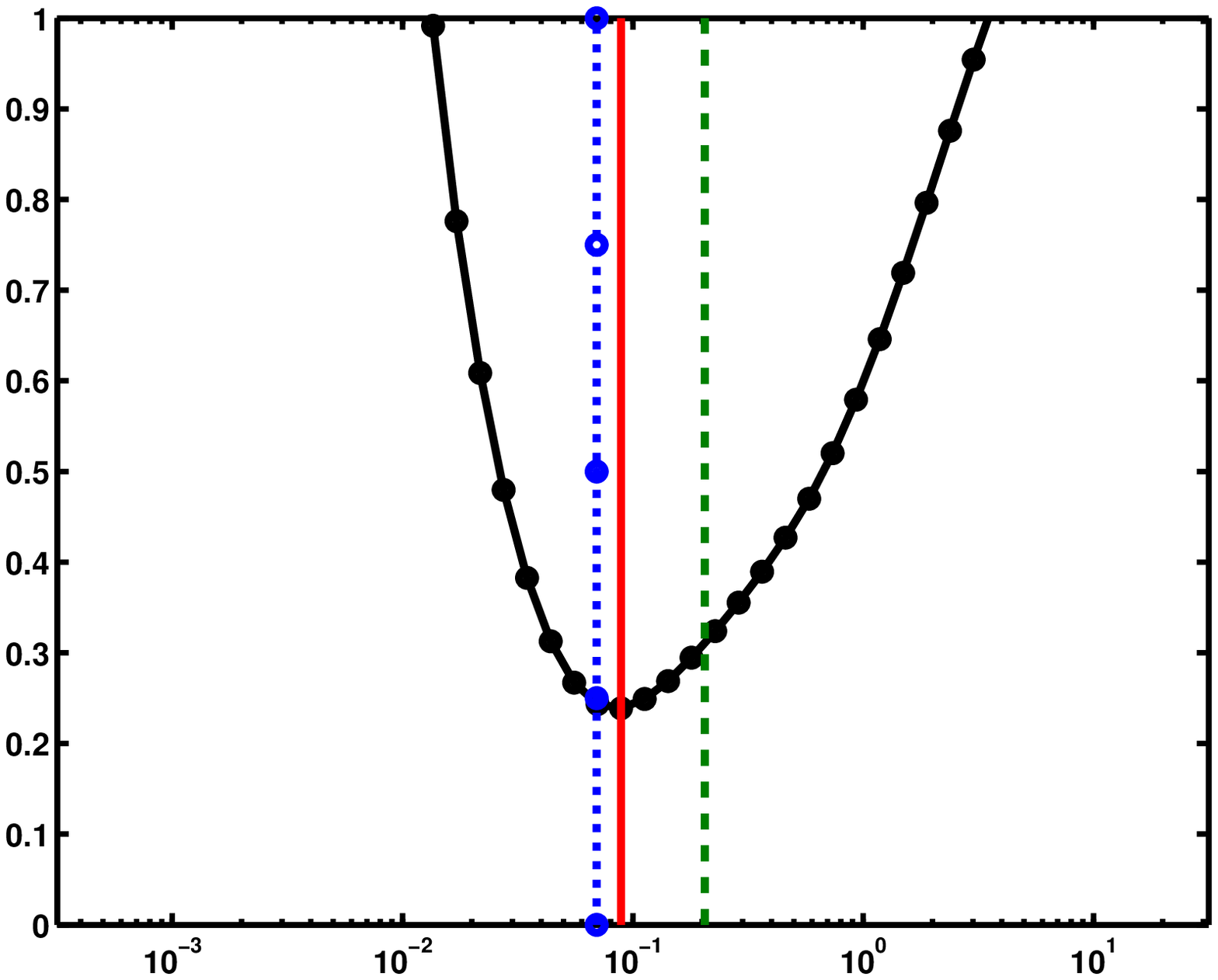}}
\subfigure[$L=L_1$]{\includegraphics[width=1.7in]{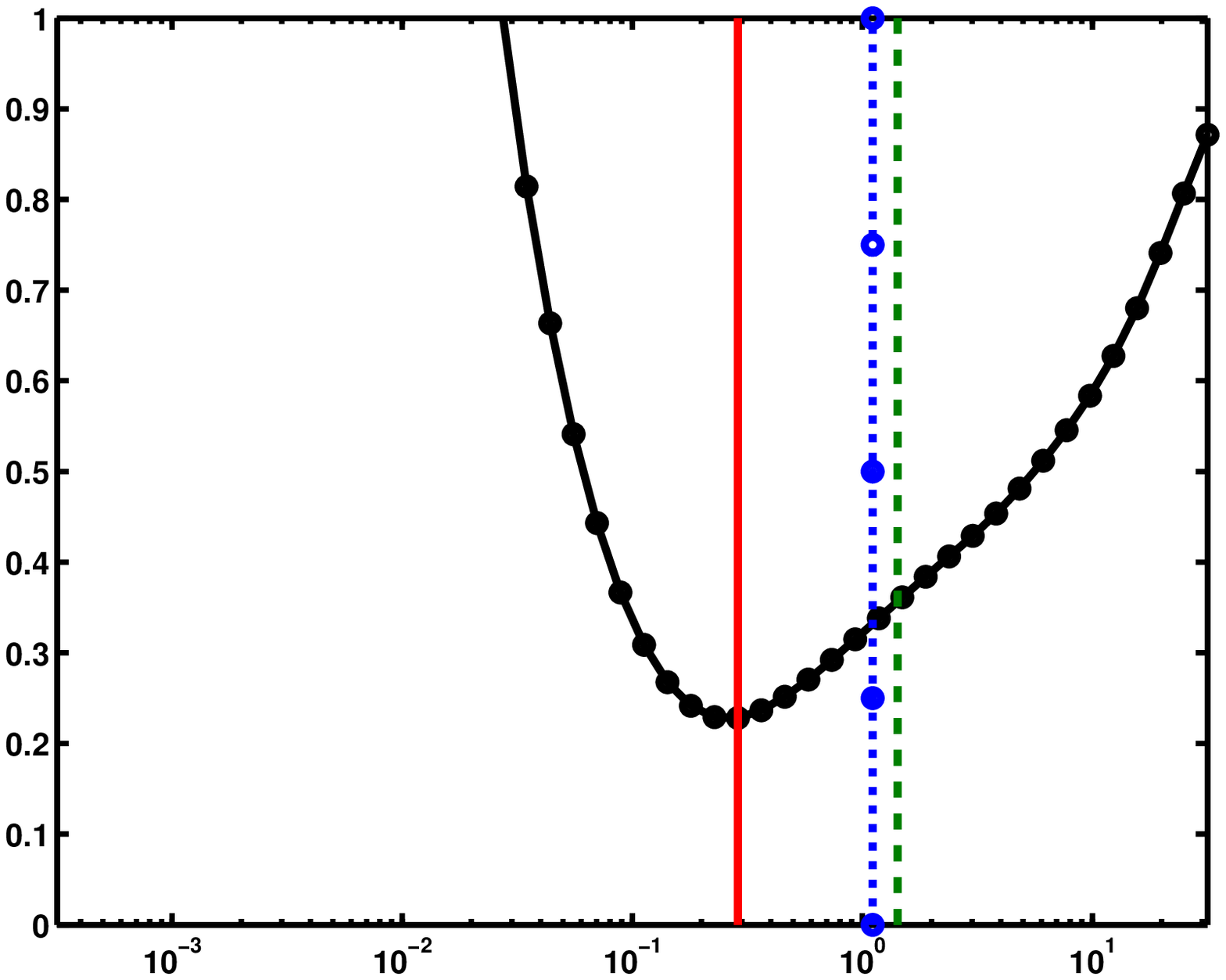}}
\subfigure[$L=L_2$]{\includegraphics[width=1.7in]{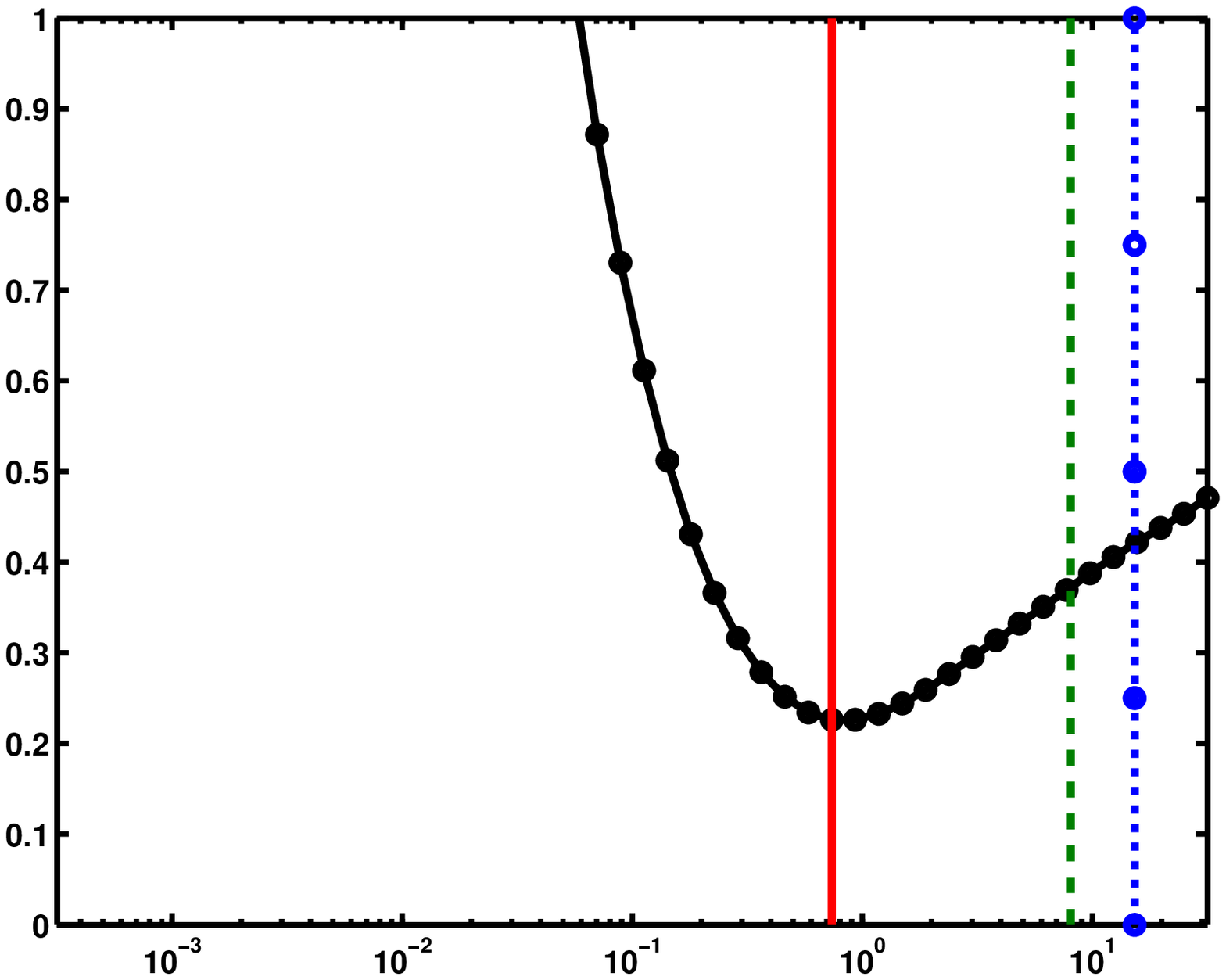}}
\subfigure[$L=I$]{\includegraphics[width=1.7in]{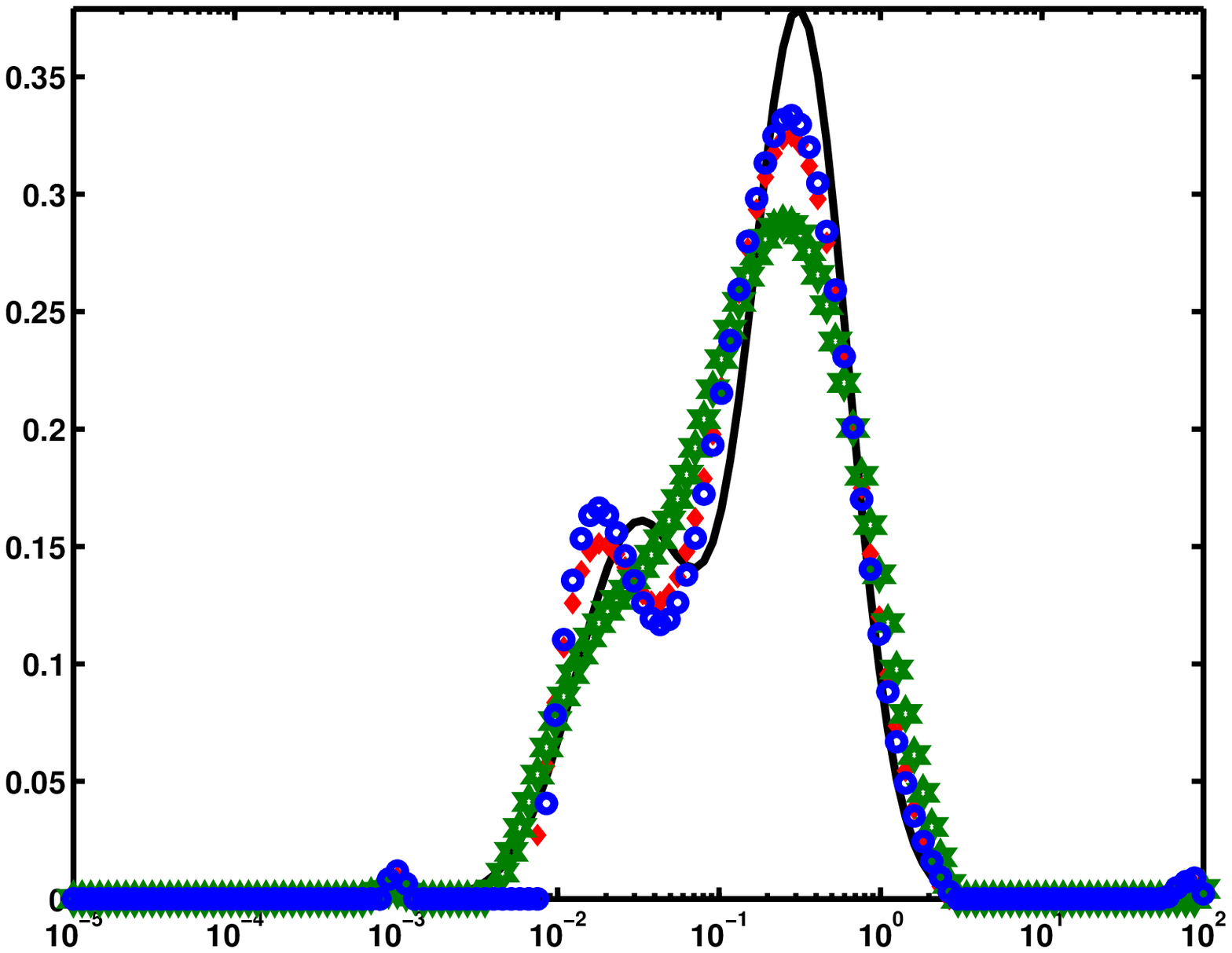}}
\subfigure[$L=L_1$]{\includegraphics[width=1.7in]{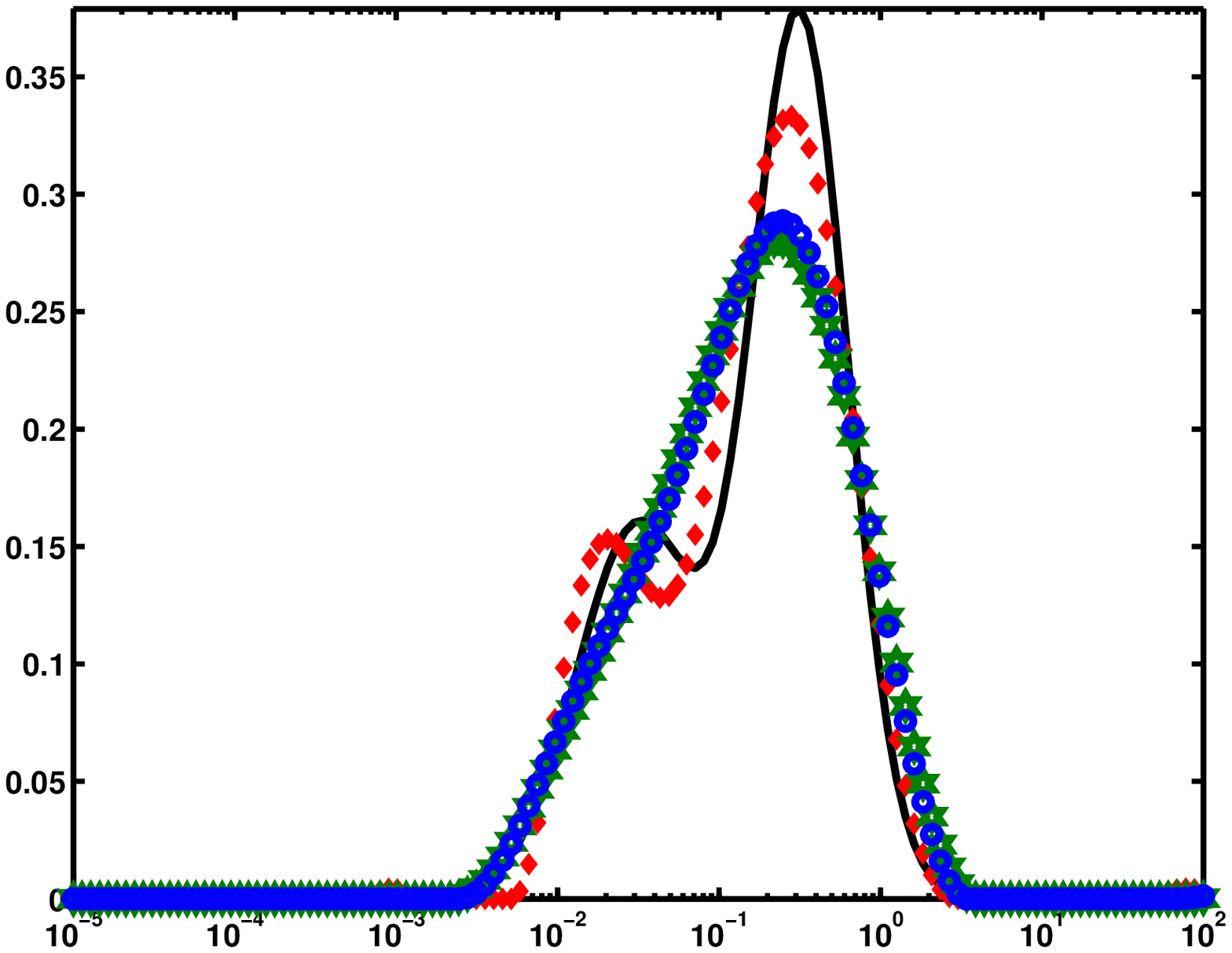}}
\subfigure[$L=L_2$]{\includegraphics[width=1.7in]{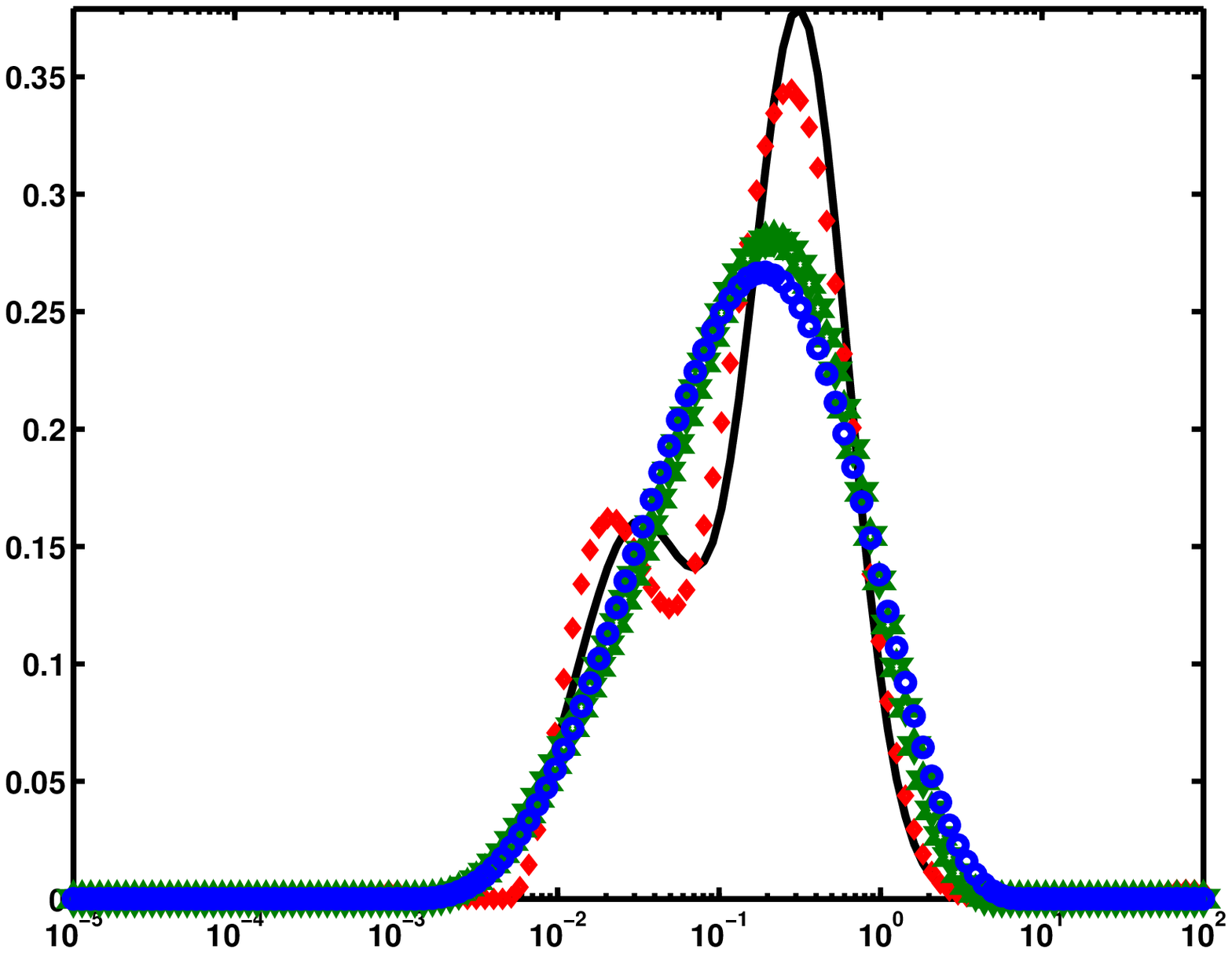}}
\caption{NNLS solutions of LN-B matrix $A_4$. Noise level $1\%$.Method CVX}
\label{hnfig-lambdachoiceLN5A4HNCVX}
\end{figure}

\begin{figure}[!ht]
\centering
\subfigure[$L=I$]{\includegraphics[width=1.7in]{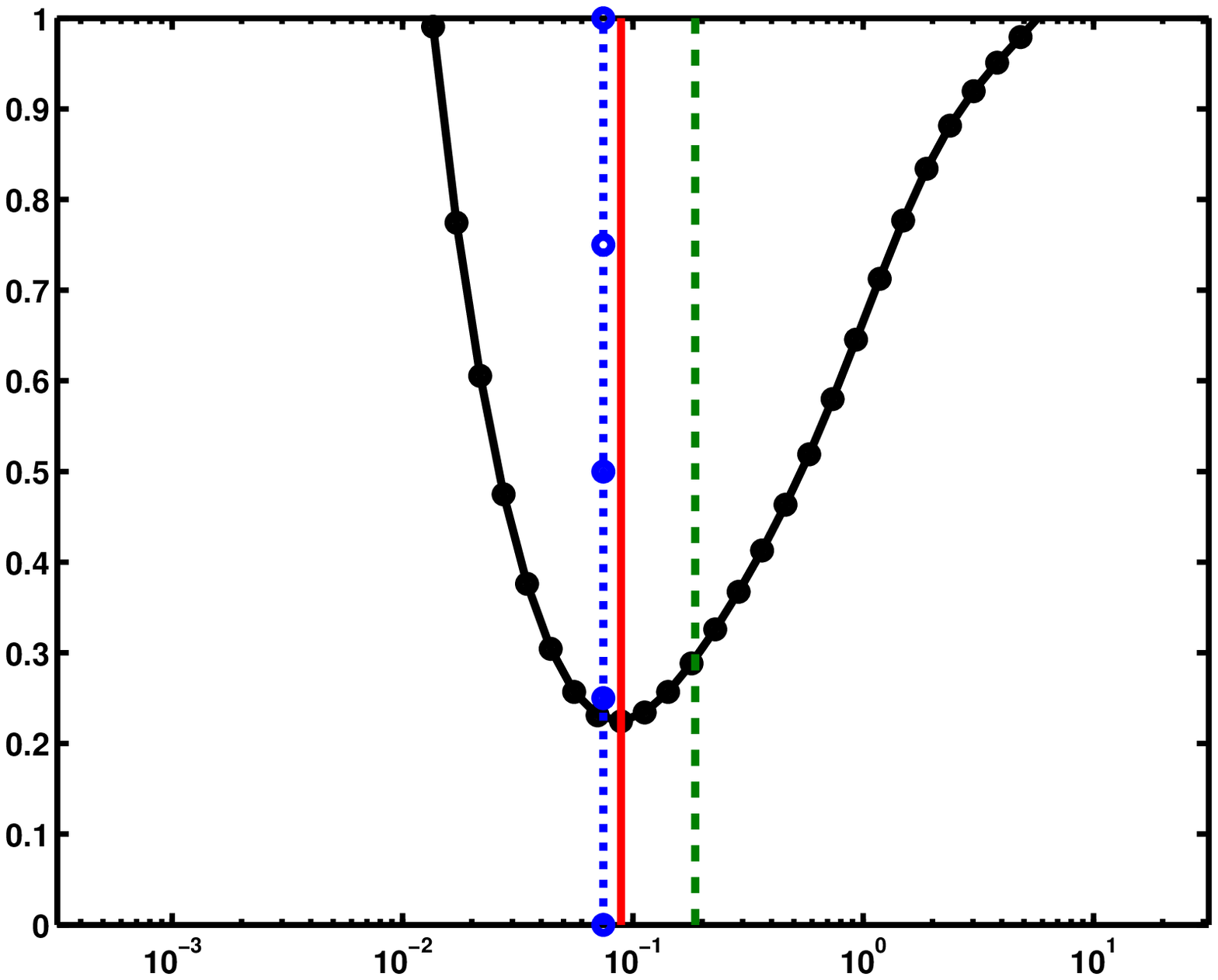}}
\subfigure[$L=L_1$]{\includegraphics[width=1.7in]{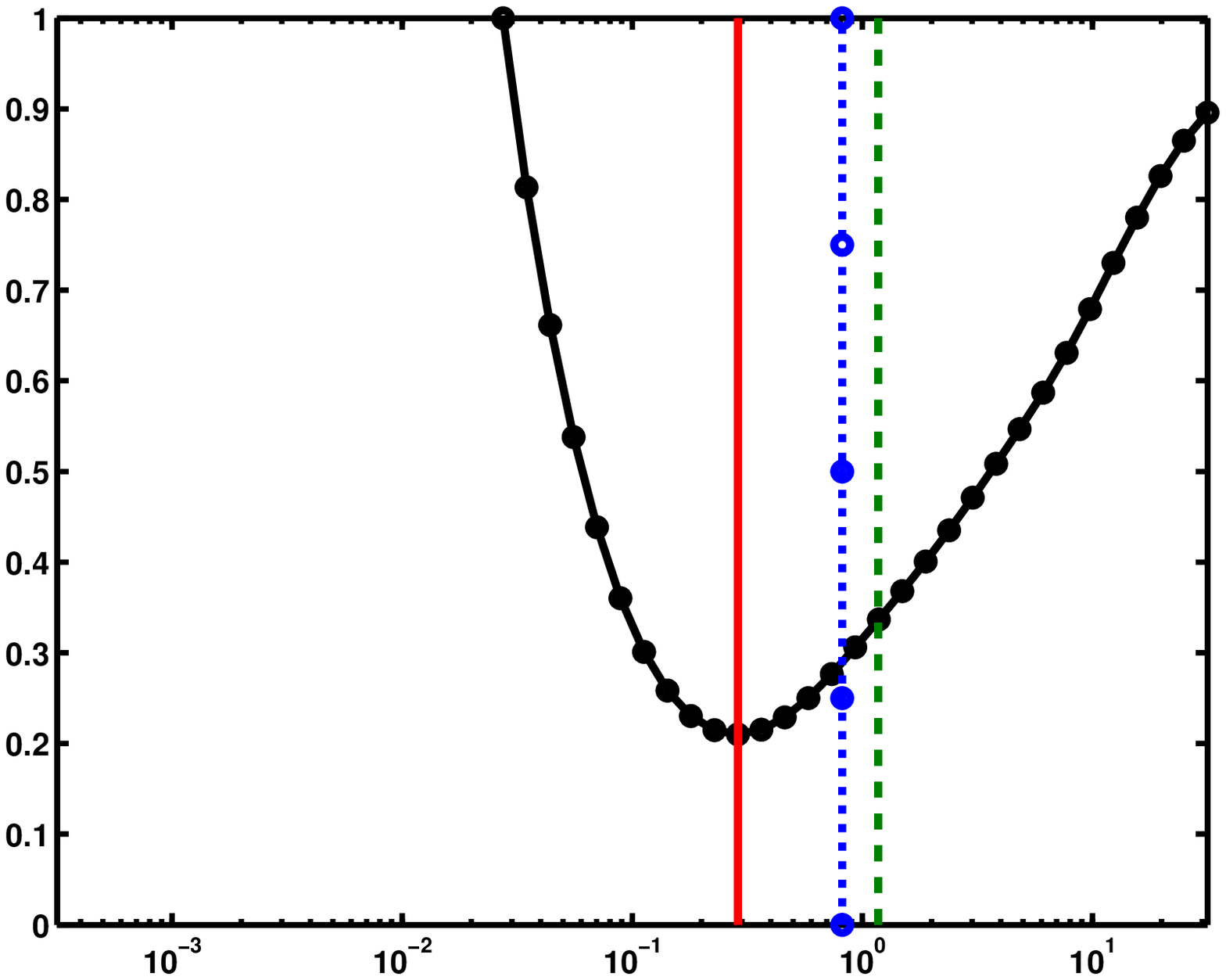}}
\subfigure[$L=L_2$]{\includegraphics[width=1.7in]{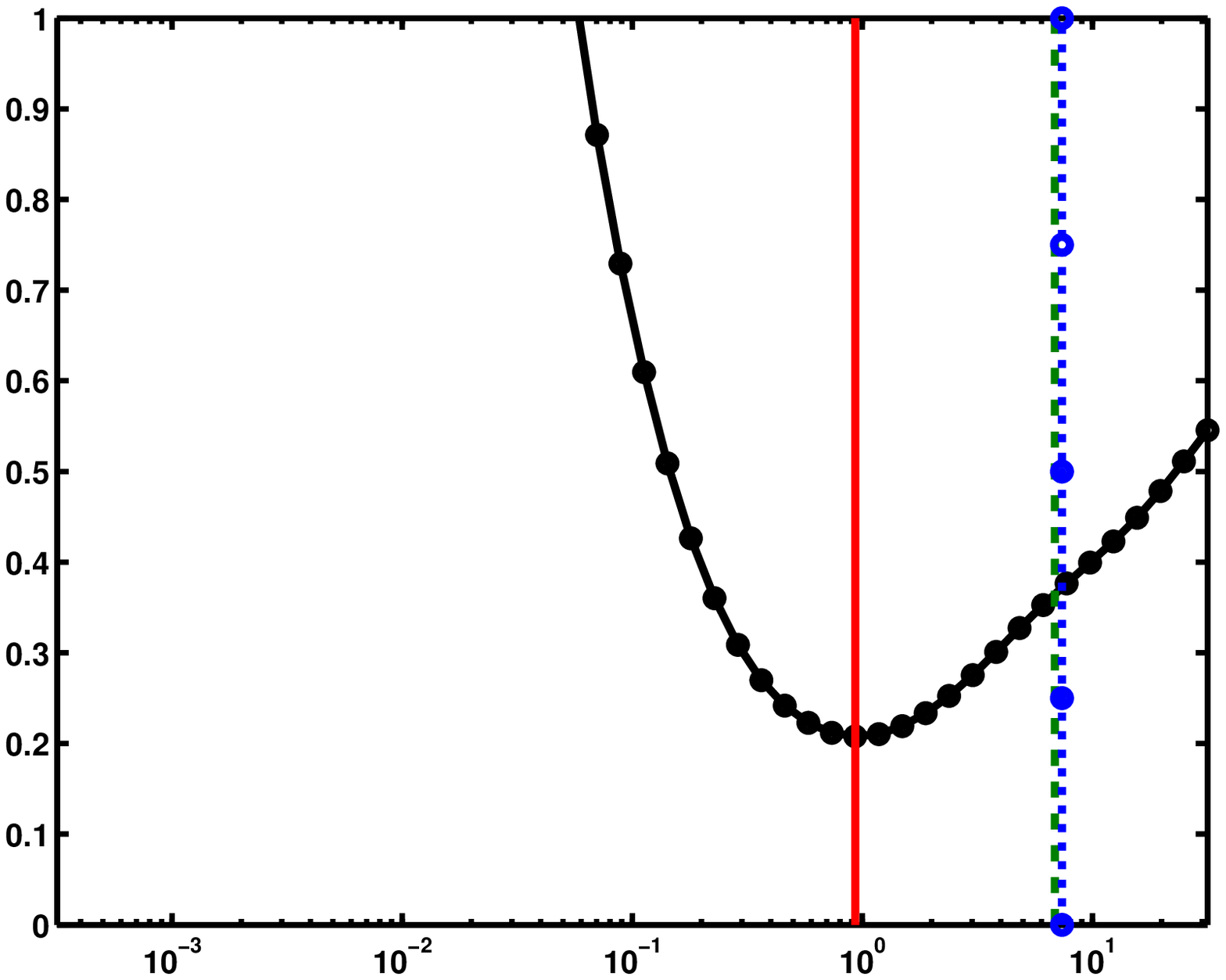}}
\subfigure[$L=I$]{\includegraphics[width=1.7in]{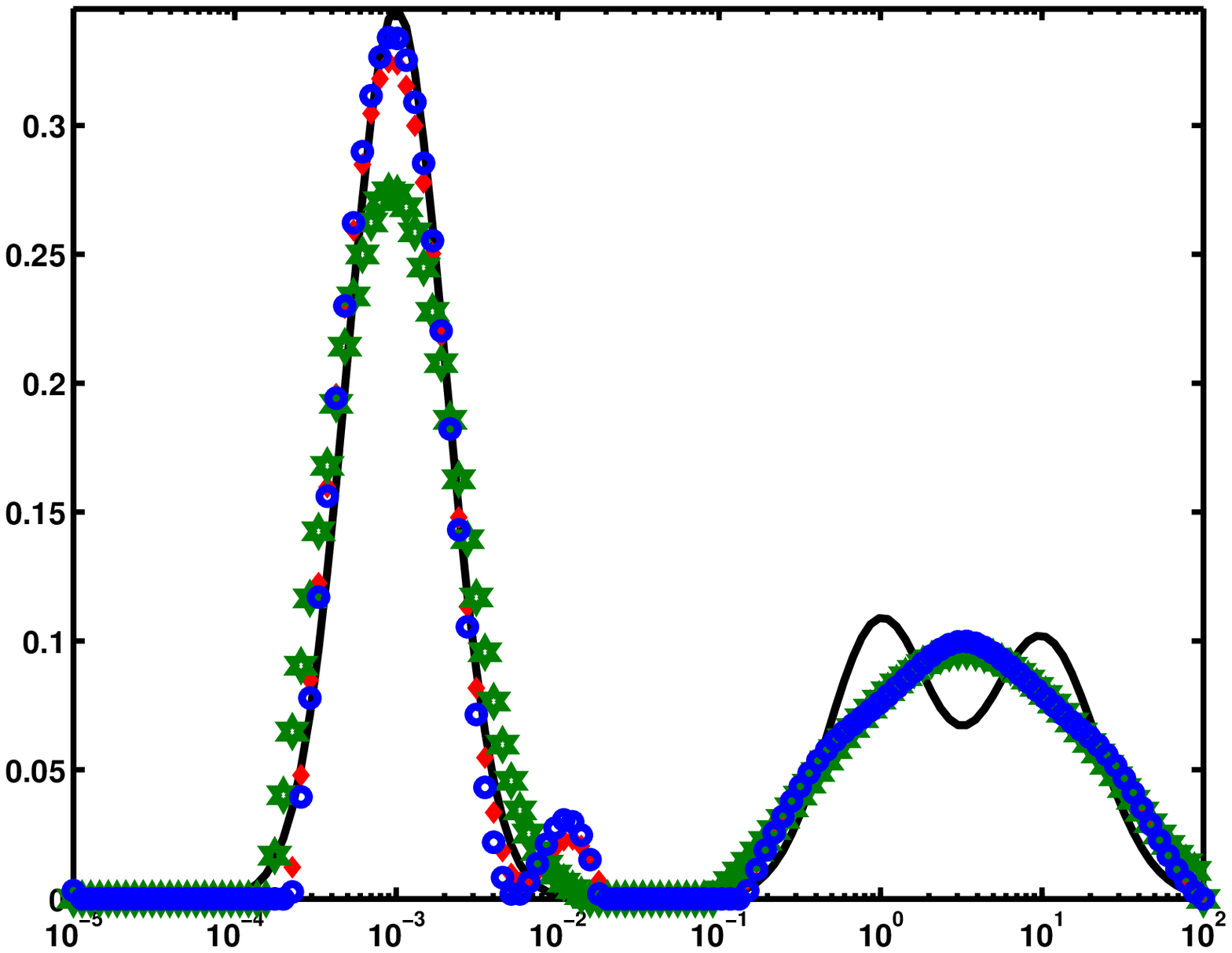}}
\subfigure[$L=L_1$]{\includegraphics[width=1.7in]{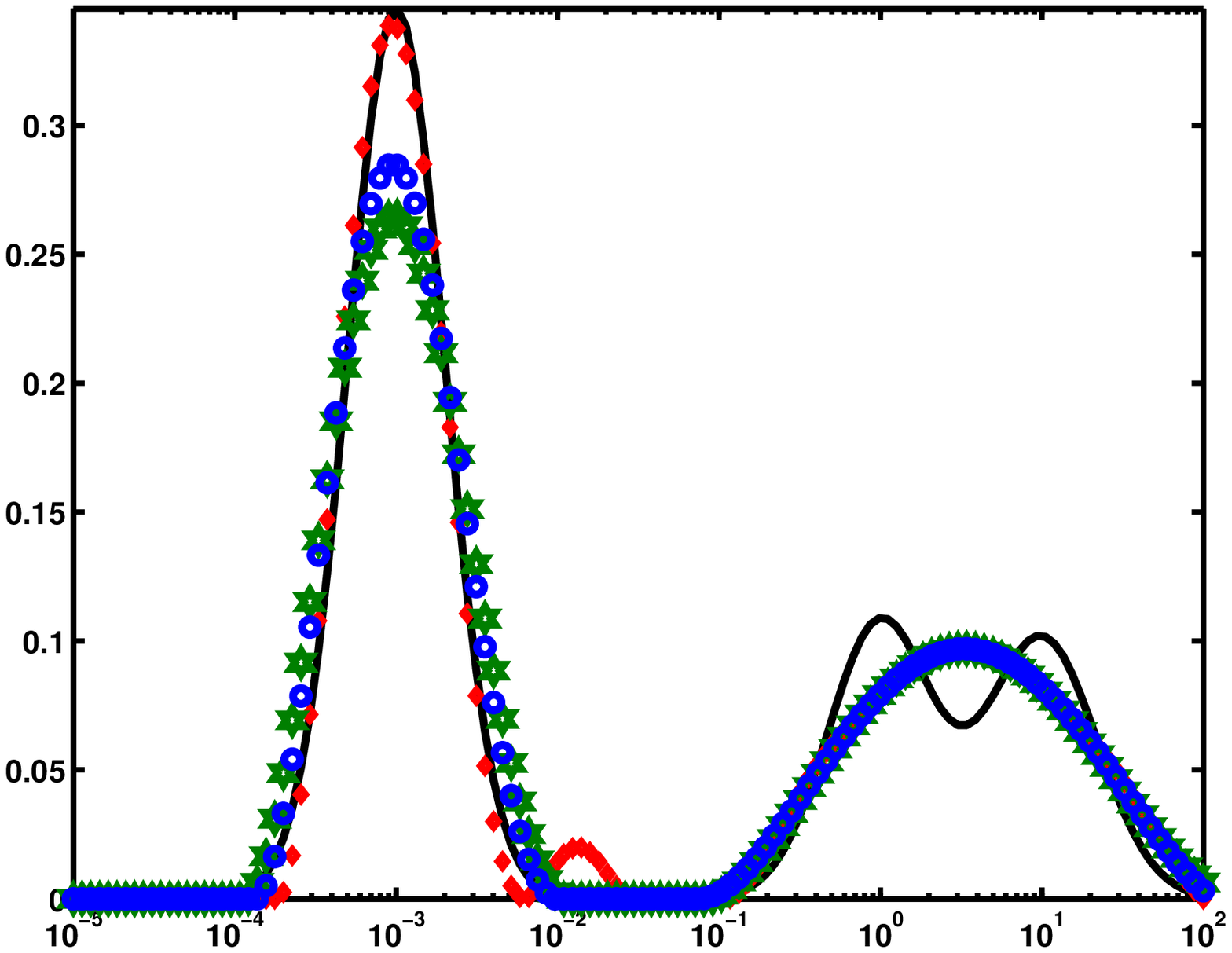}}
\subfigure[$L=L_2$]{\includegraphics[width=1.7in]{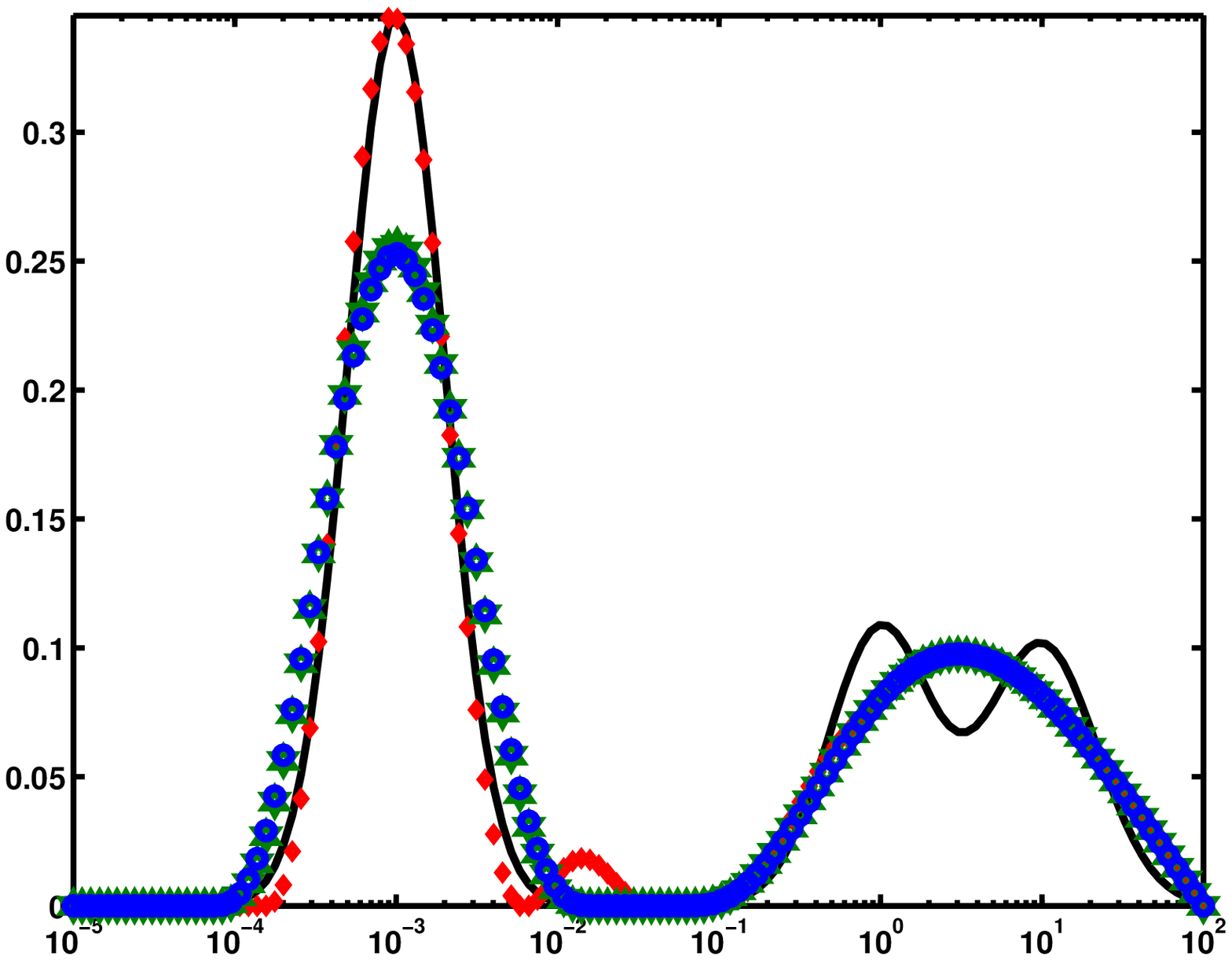}}
\caption{NNLS solutions of LN-C matrix $A_4$. Noise level $1\%$.Method CVX}
\label{hnfig-lambdachoiceLN6A4HNCVX}
\end{figure}

In the figures we compare the parameter choice methods. For each set of results the first row, Figures (a)-(c) in each case indicate the mean error results for the different smoothing norms, and (d)-(f) demonstrate the sensitivity, or lack thereof, of the solution to the choice of $\lambda$ near the optimum. We briefly list the key observations.
\begin{description}
\item[$L=I$] In this case the NCP results are more often close to the optimum
\item[$L=L_1$] Except for high noise the LC results can be close to the optimum
\item[$L=L_2$] It is hard to distinguish between the LC and NCP results in terms of overall best match to the optimal solution
\item[High noise] It is clear that the results with $5\%$ noise are not as good. Moreover the LC results are consistently under smoothed for all operator norms.
\item[$A_3$ or $A_4$] The parameter choice methods perform quite similarly for both matrices. The lack of resolution of $A_3$ is now more apparent. 
\item[NNLS Algorithm] Two different algorithms for NNLS are investigated in Sections~\ref{sec:SBB} and \ref{sec:cvx}, the SBB algorithm in \cite{kimsra} and the CVX algorithm in \cite{cvx,gb08}. 
\item[NNLS SBB] It is immediate from Figures~\ref{lnfig-lambdachoiceRQ1A4LNSBB}-\ref{lnfig-lambdachoiceLN6A4LNSBB} that the constrained Barzilai-Borwein algorithm, which was obtained  with \cite{sbbcode} creates additional difficulties for finding the optimum choice of $\lambda$ when the range of $\lambda$ includes small values. The obtained solutions as $\lambda \to 0$ tend to solutions with constant error, but when regarded in the solution space, the Tikhonov regularization is not sufficiently applied so that solutions have theoretical low error, but are insufficiently smoothed. This is a feature of the fact that the two norm of the error does not always provide a good mechanism for finding a good solution. Indeed, we know that as $\lambda$ decreases less smoothing is applied to the solution, and hence on the average the solution may have less error, (the two norm error), but provides a noisier estimate of the actual solution. No additional results are provided for this algorithm.
\item[NNLS CVX] Results for the CVX implementation are provided in the original article, \cite{JCAMHansen}. They demonstrate for low noise results that are comparable to the use of the \texttt{lsqnonneg} algorithm. Here further results in Figures-\ref{hnfig-lambdachoiceRQ1A4HNCVX}-\ref{hnfig-lambdachoiceLN6A4HNCVX} show that robustness holds for increasing noise levels. 
\end{description}
\clearpage
\subsection{Examples: $.1\%$ noise matrix $A_3$ NNLS }

\begin{figure}[!h]
 \centering
\subfigure[$L=I$]{\includegraphics[width=1.7in]{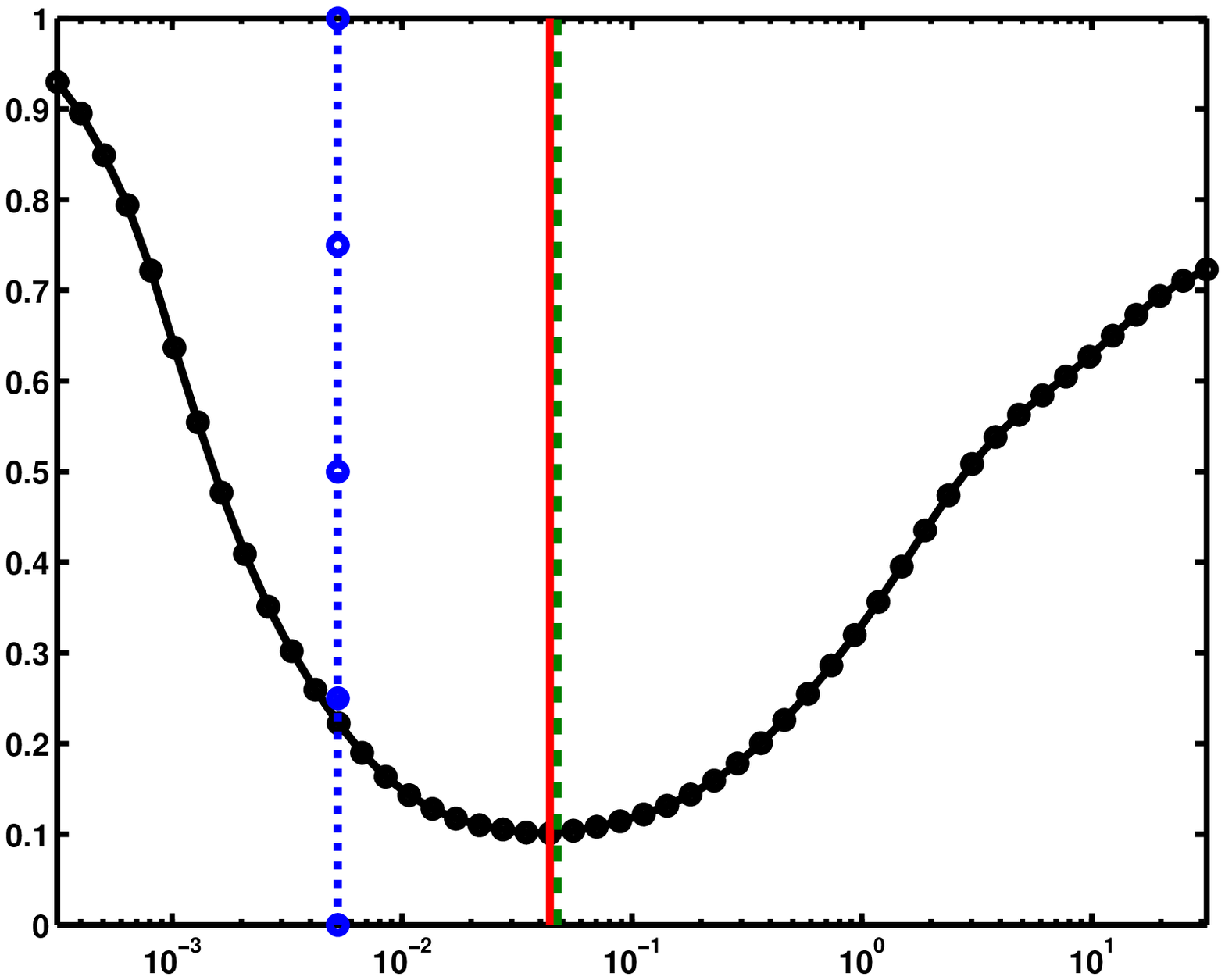}}
\subfigure[$L=L_1$]{\includegraphics[width=1.7in]{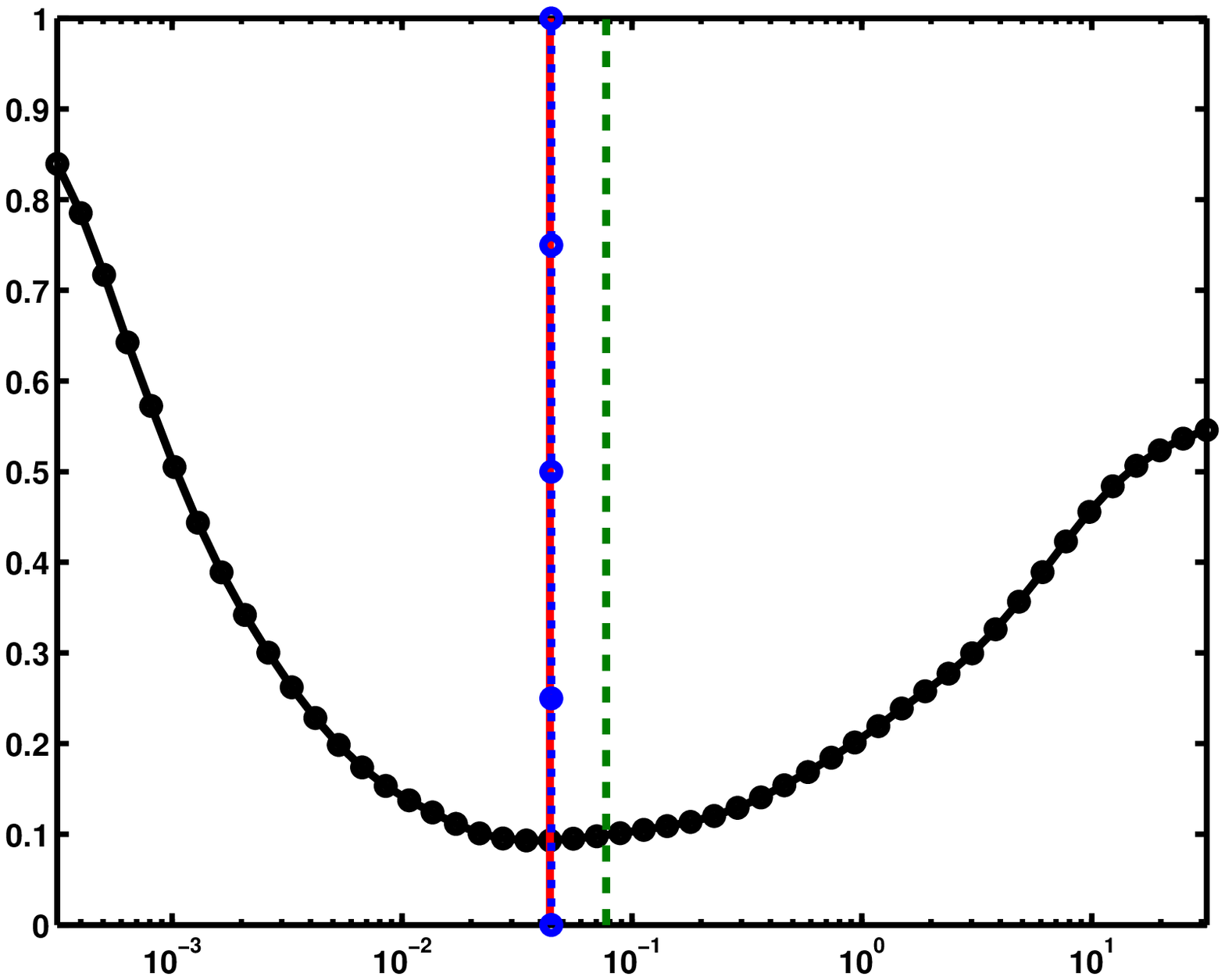}}
\subfigure[$L=L_2$]{\includegraphics[width=1.7in]{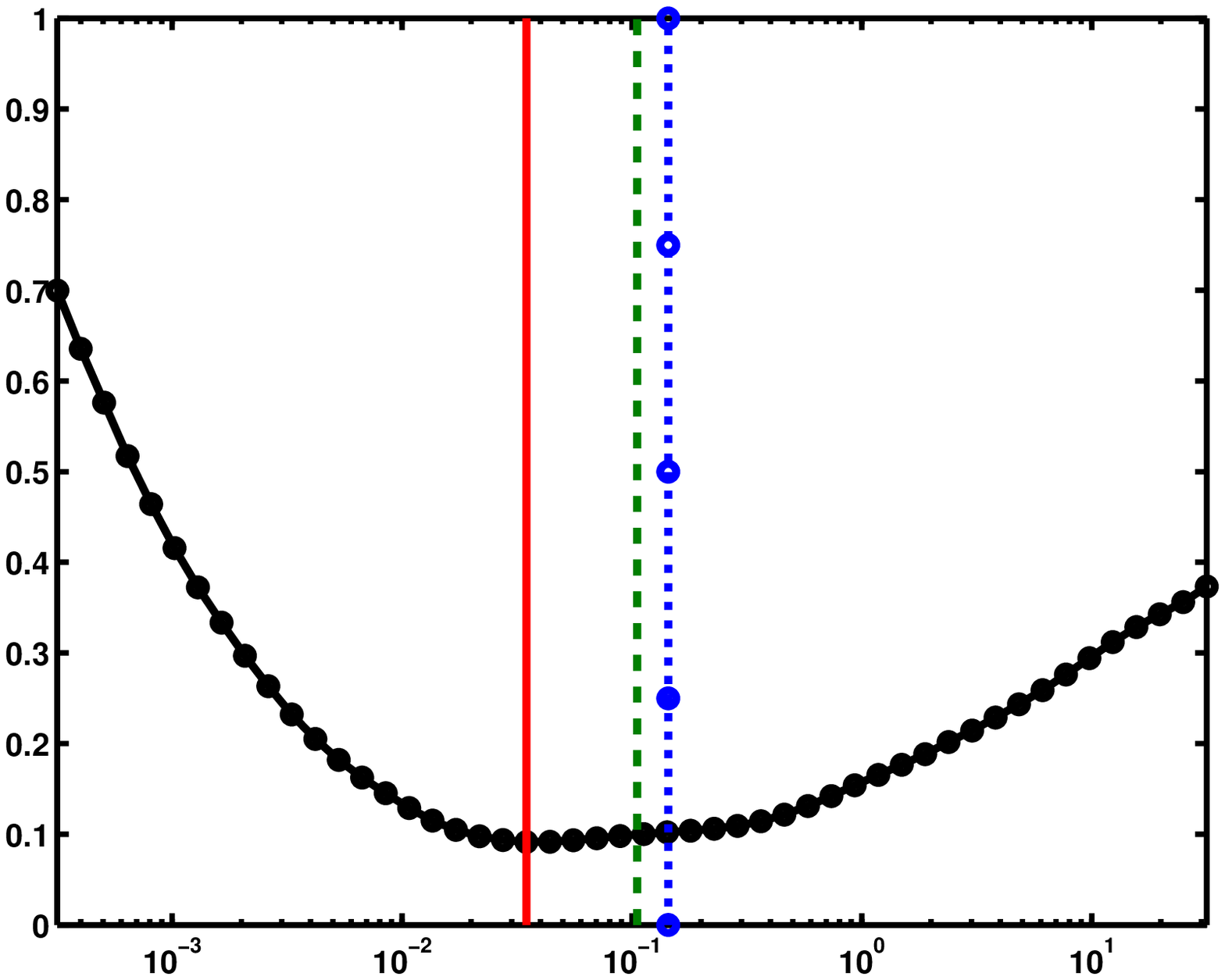}}
\subfigure[$L=I$]{\includegraphics[width=1.7in]{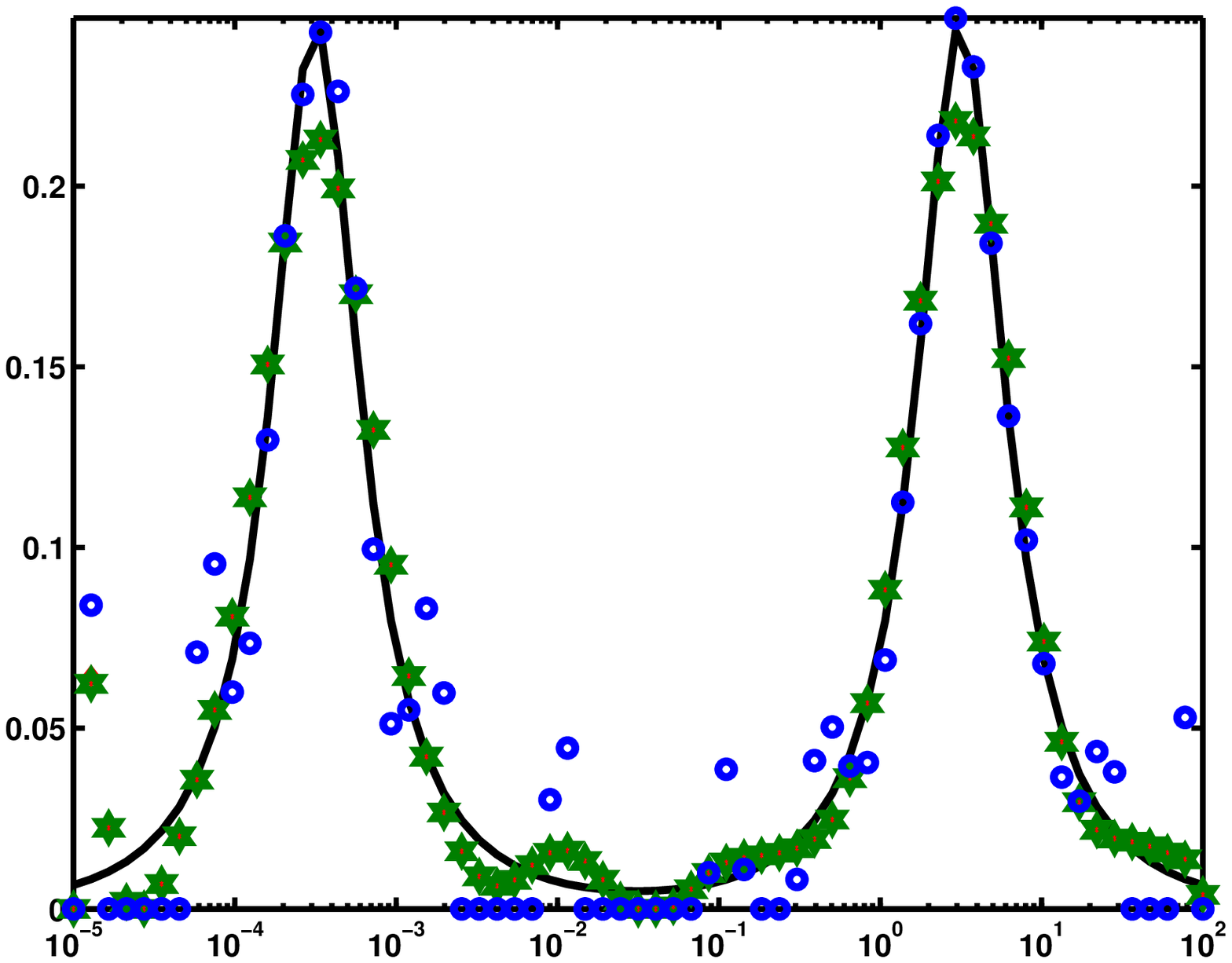}}
\subfigure[$L=L_1$]{\includegraphics[width=1.7in]{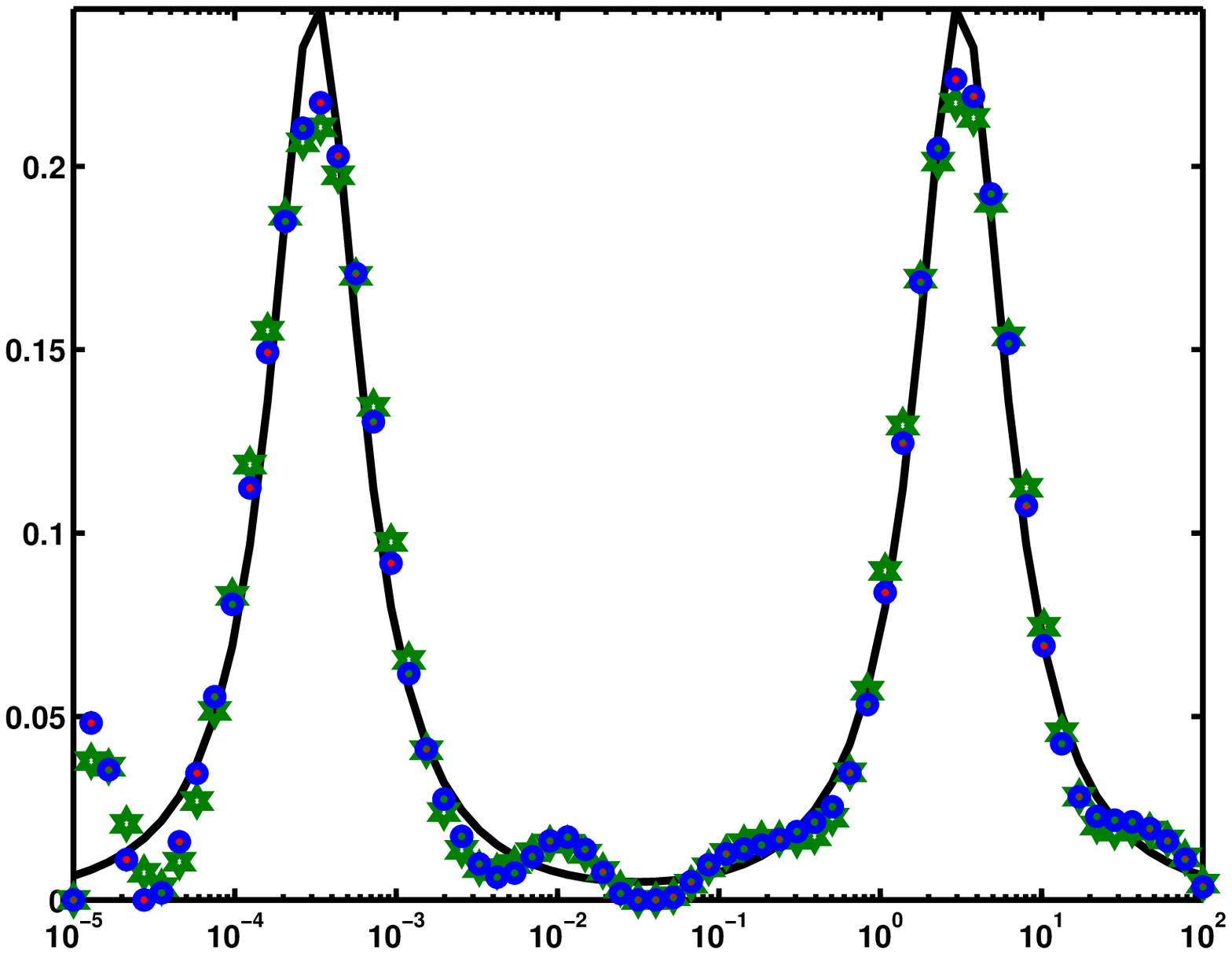}}
\subfigure[$L=L_2$]{\includegraphics[width=1.7in]{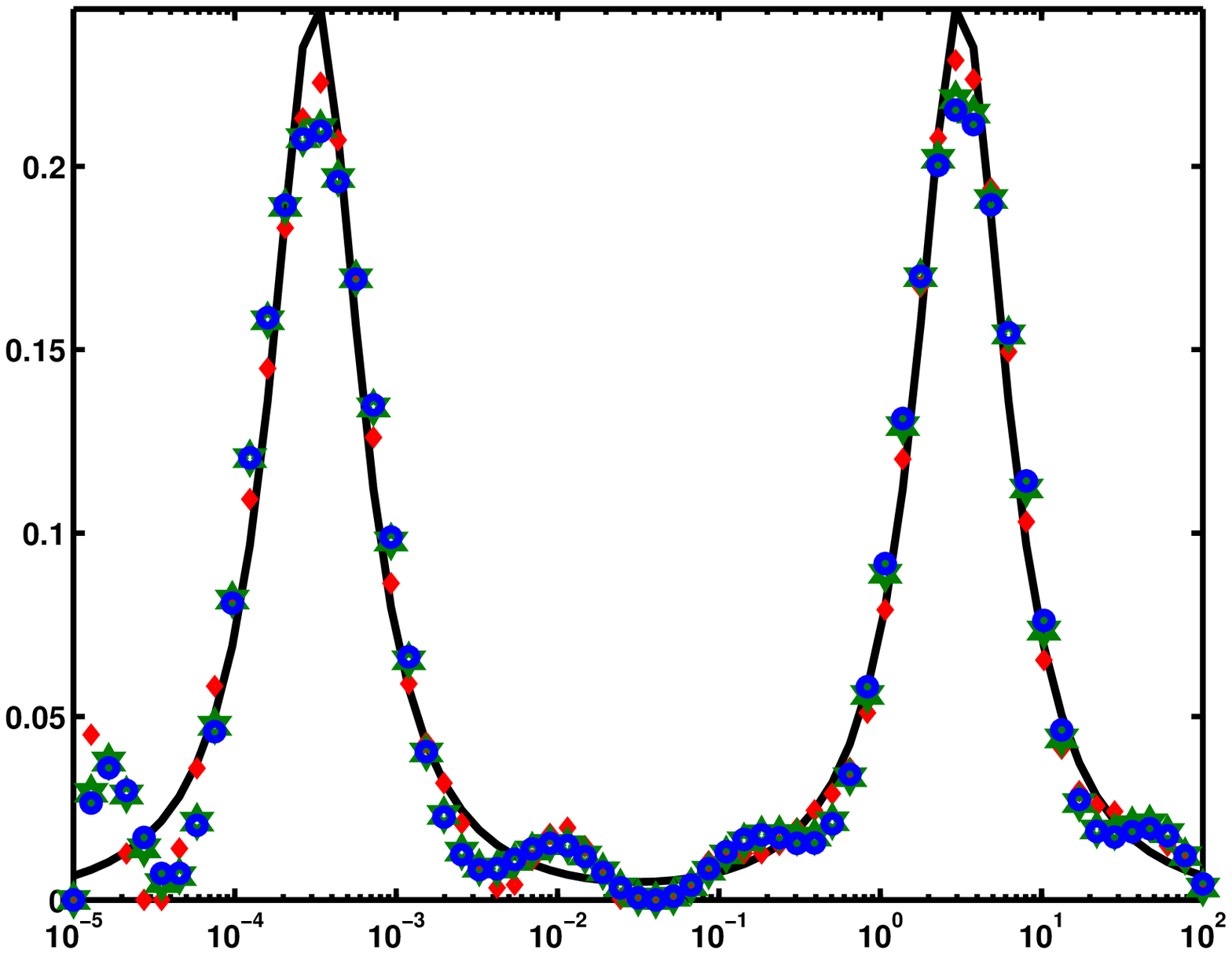}}
\caption{Mean error and example NNLS solutions.  $.1\%$ noise. RQ-A data set matrix $A_3$}
\label{fig-lambdachoiceRQ1A3LN}
\end{figure}

\begin{figure}[!h]
 \centering
\subfigure[$L=I$]{\includegraphics[width=1.7in]{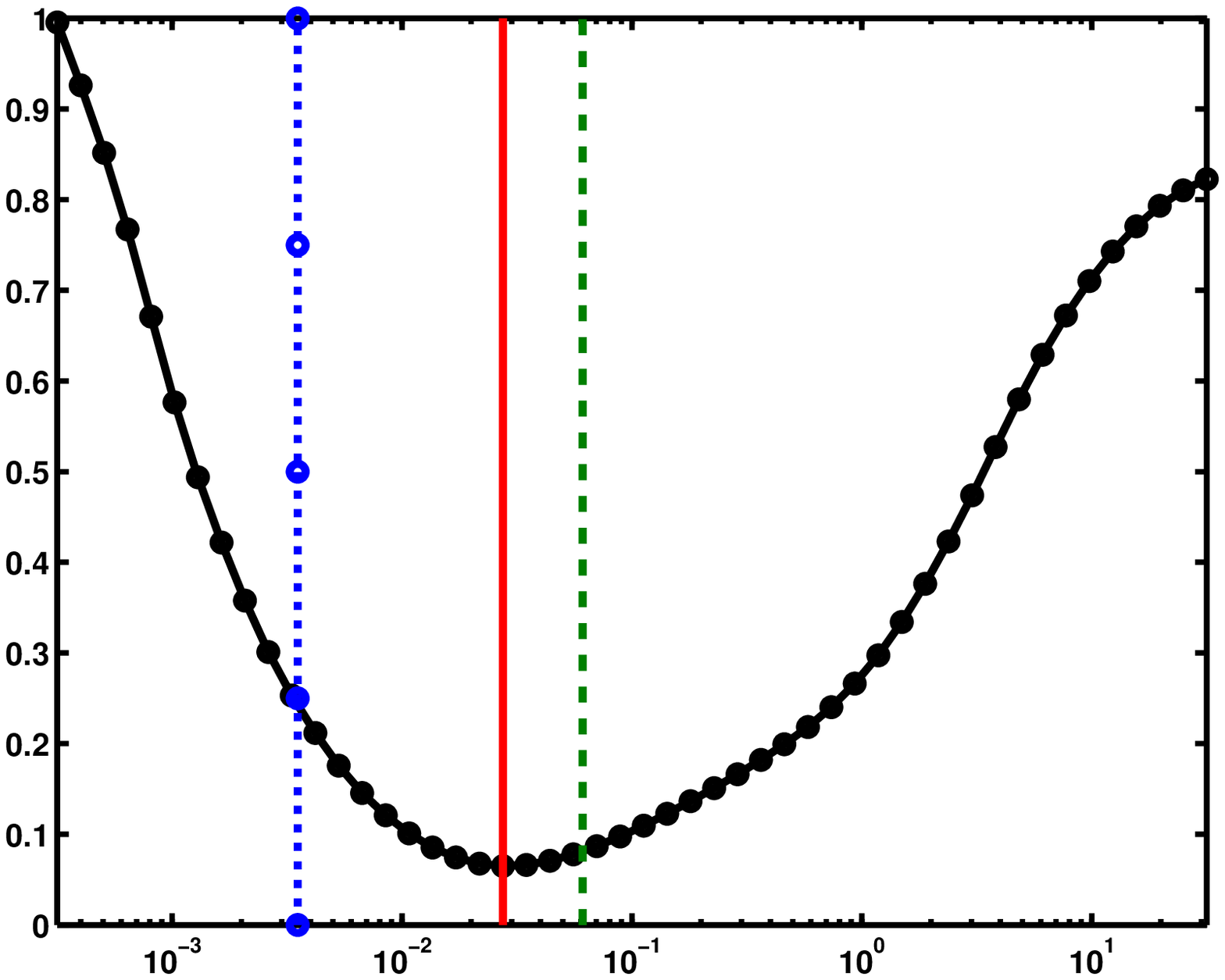}}
\subfigure[$L=L_1$]{\includegraphics[width=1.7in]{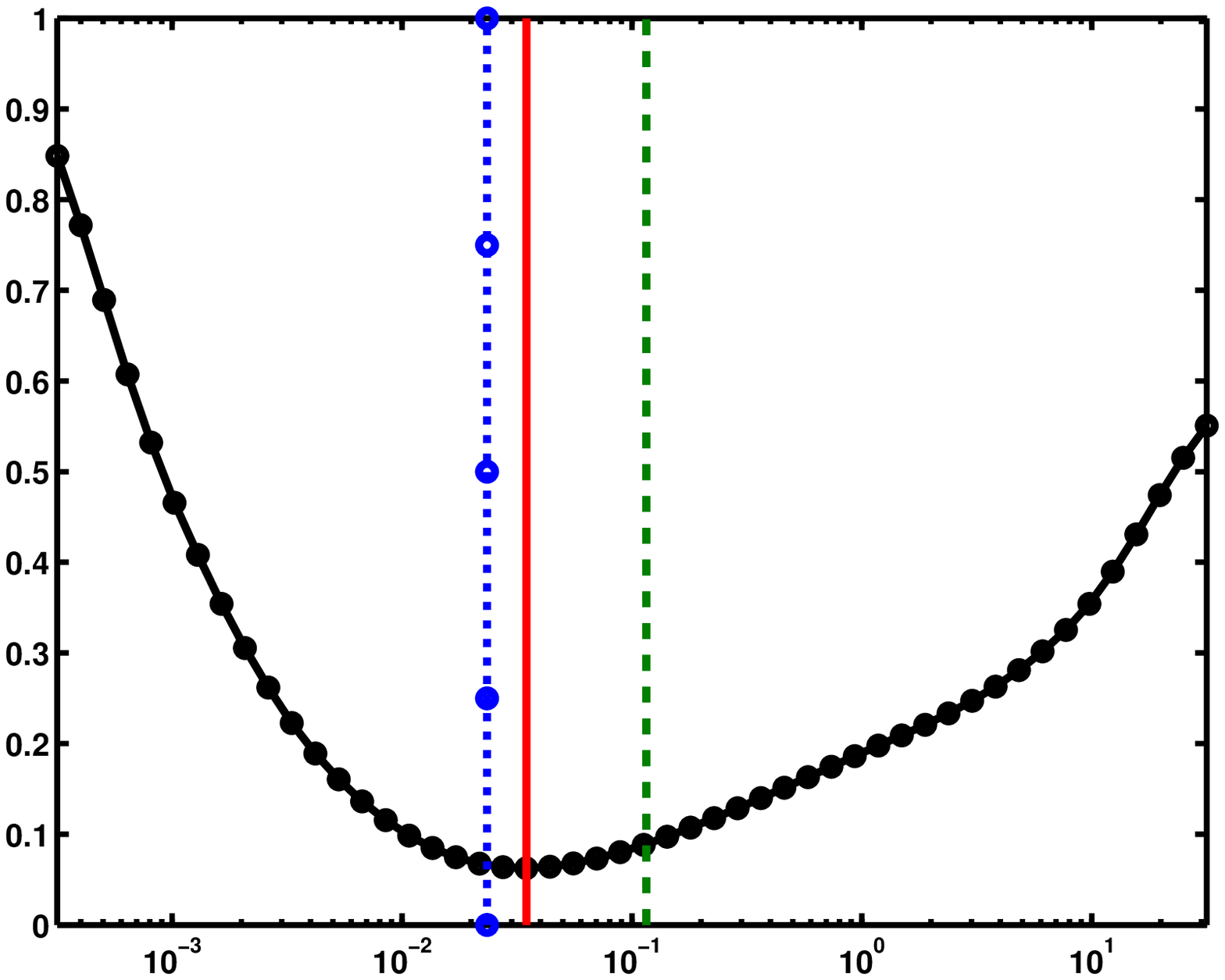}}
\subfigure[$L=L_2$]{\includegraphics[width=1.7in]{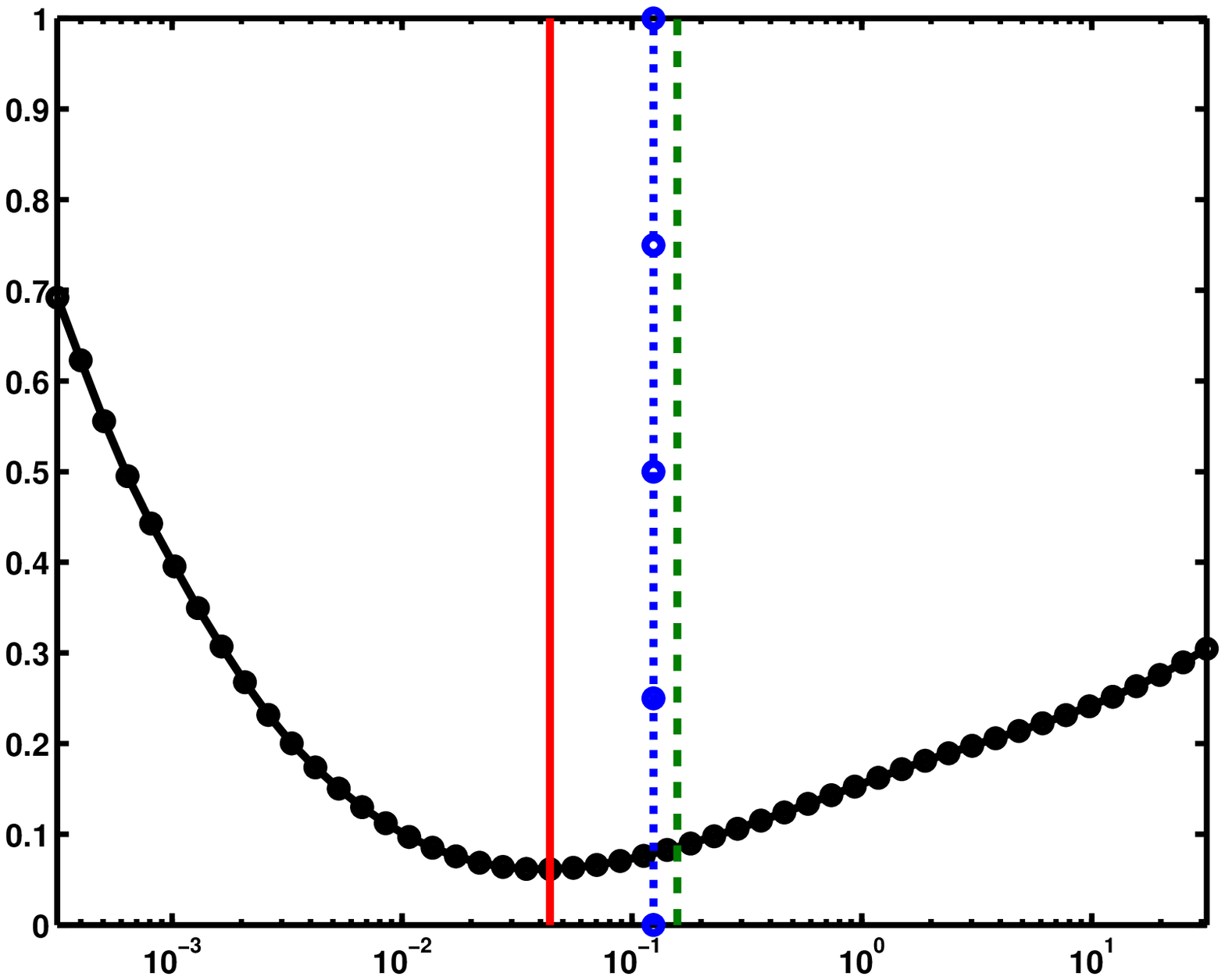}}
\subfigure[$L=I$]{\includegraphics[width=1.7in]{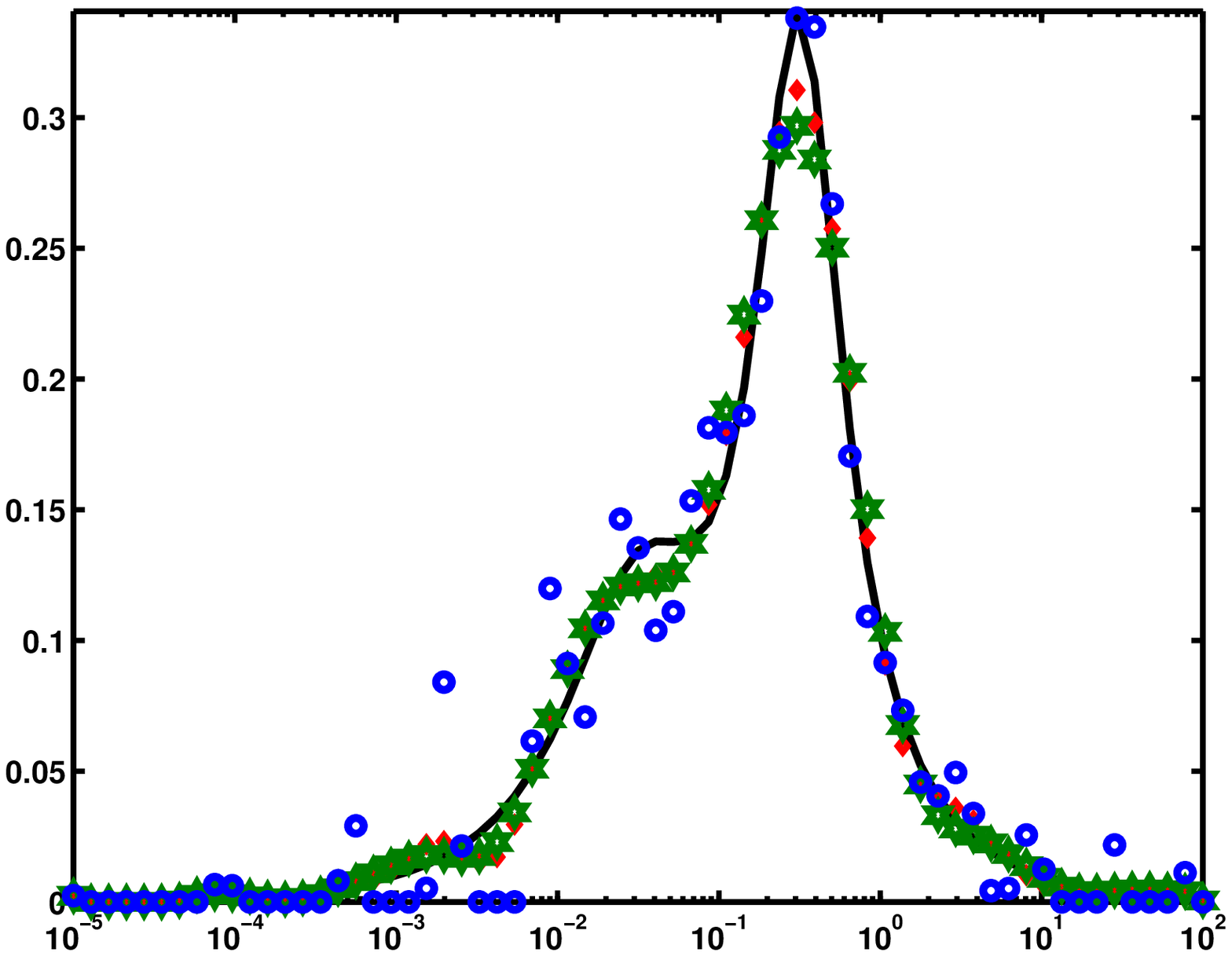}}
\subfigure[$L=L_1$]{\includegraphics[width=1.7in]{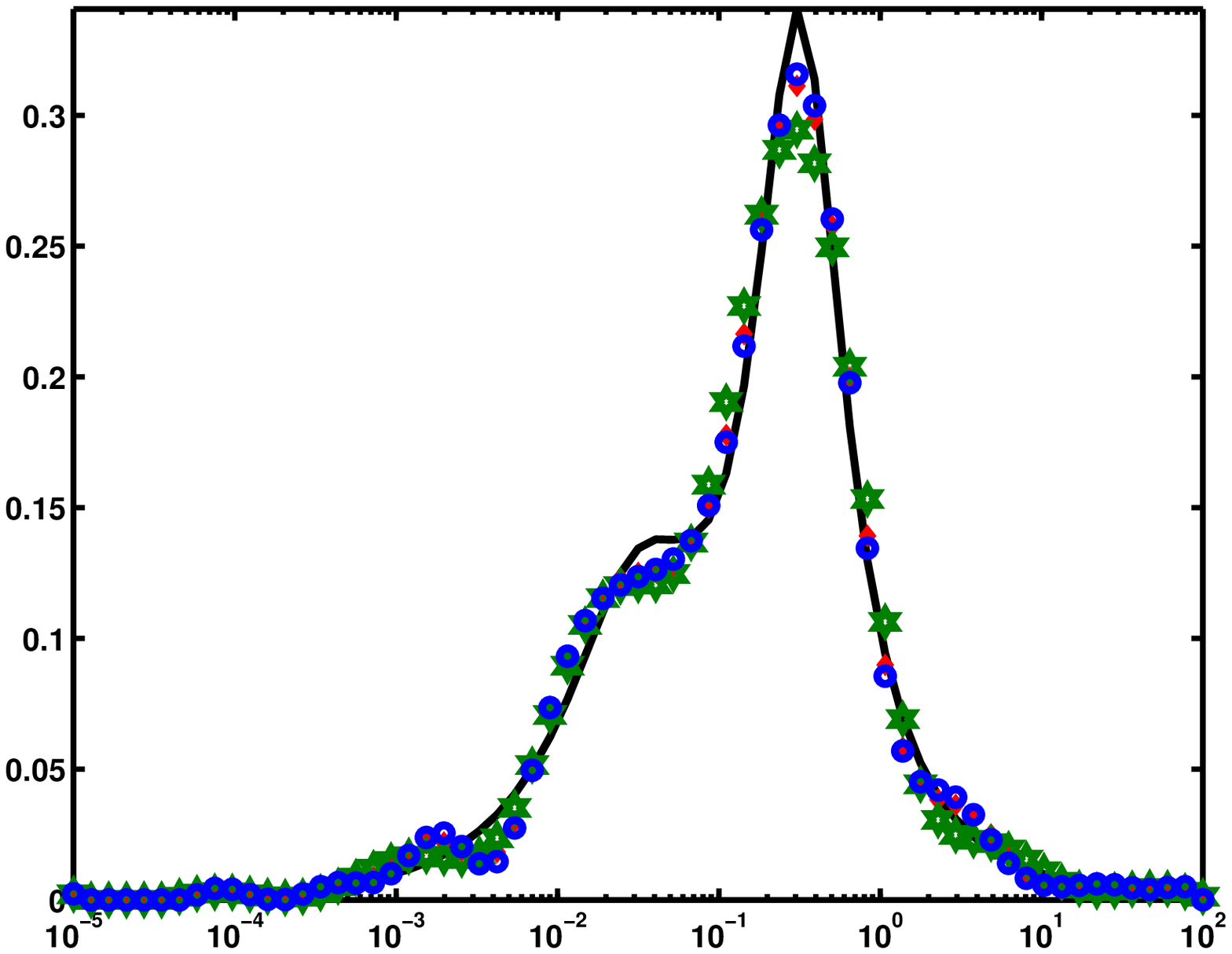}}
\subfigure[$L=L_2$]{\includegraphics[width=1.7in]{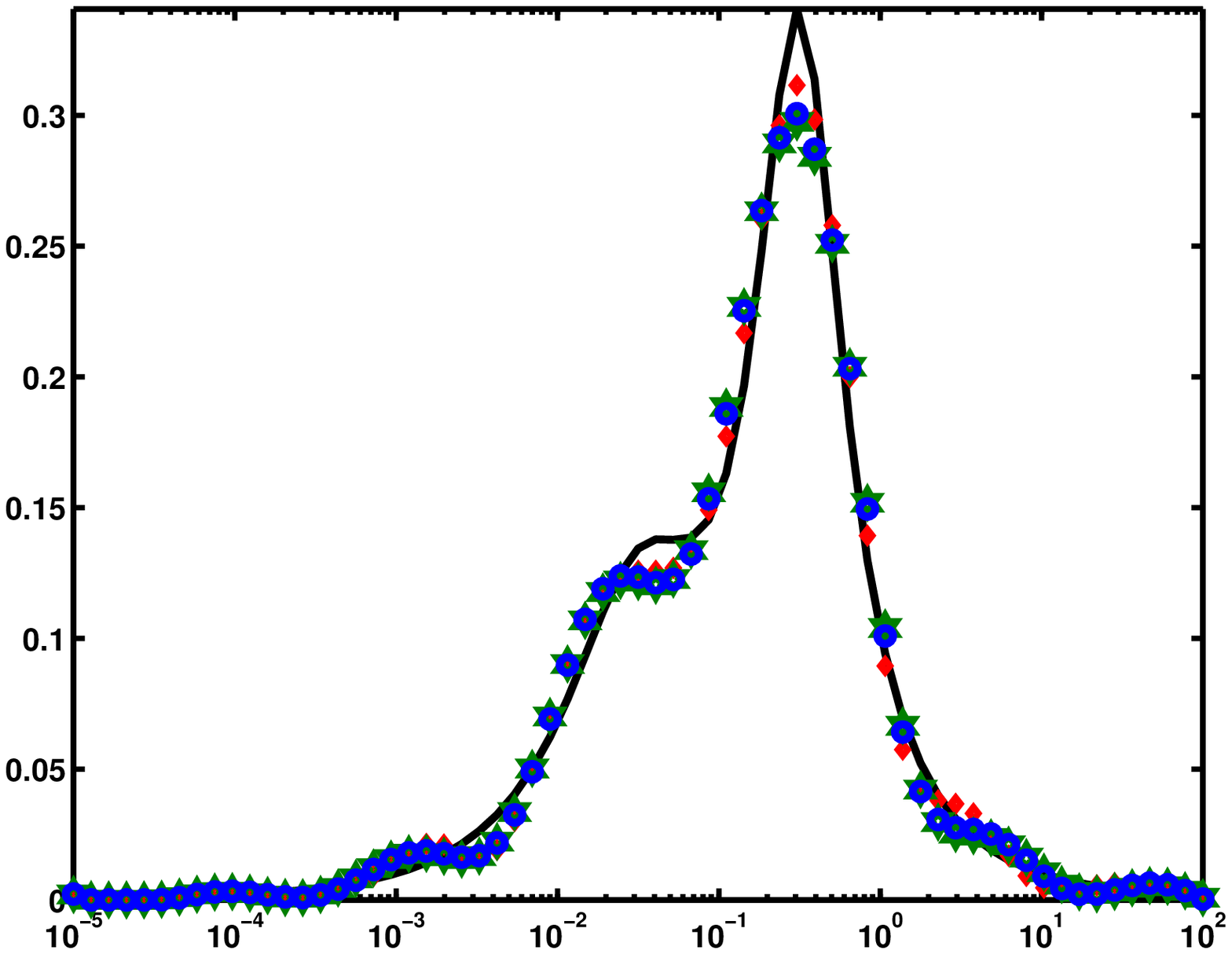}}
\caption{Mean error and example NNLS solutions.  $.1\%$ noise. RQ-B data set matrix $A_3$}
\label{fig-lambdachoiceRQ5A3LN}
\end{figure}

\begin{figure}[!h]
 \centering
\subfigure[$L=I$]{\includegraphics[width=1.7in]{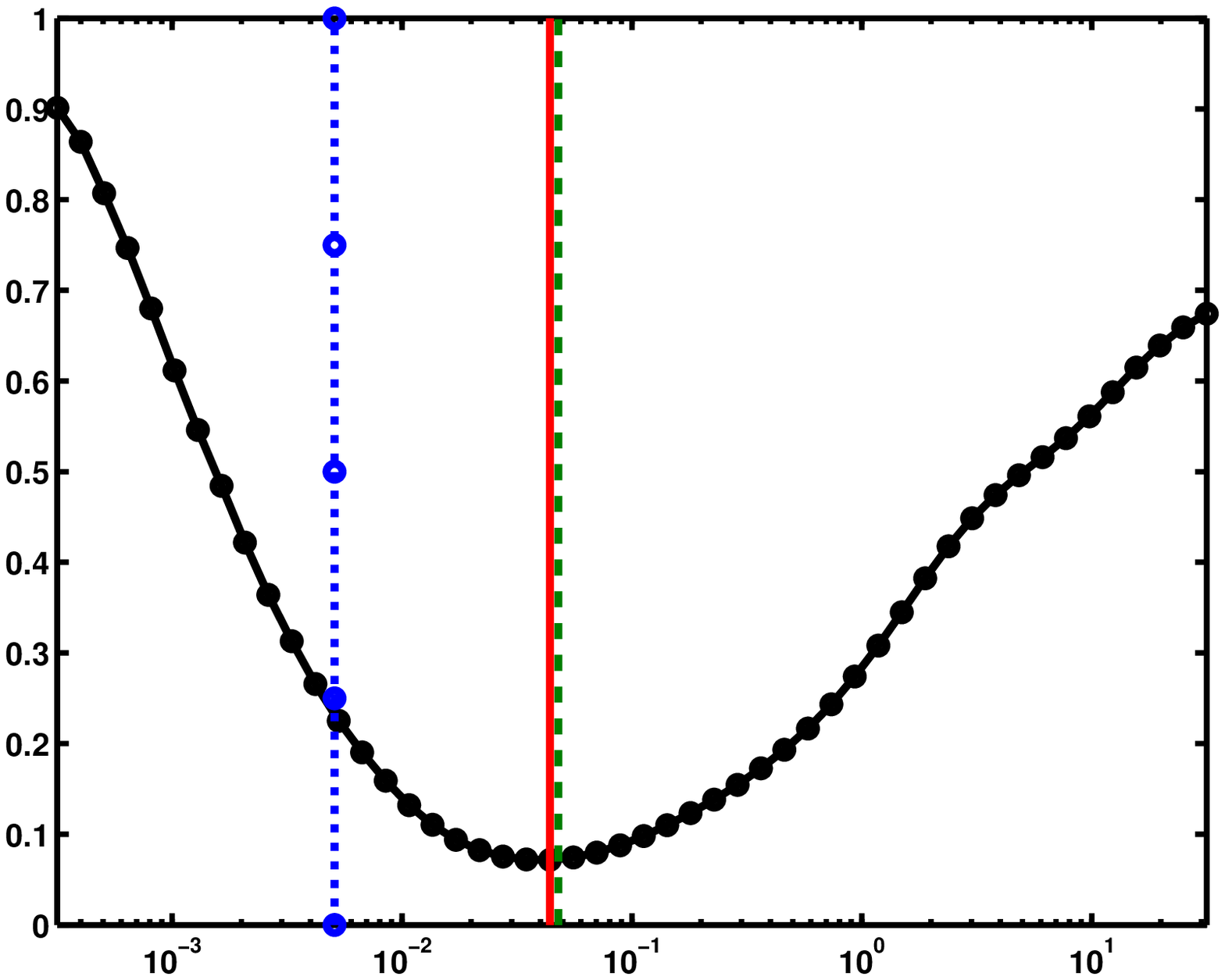}}
\subfigure[$L=L_1$]{\includegraphics[width=1.7in]{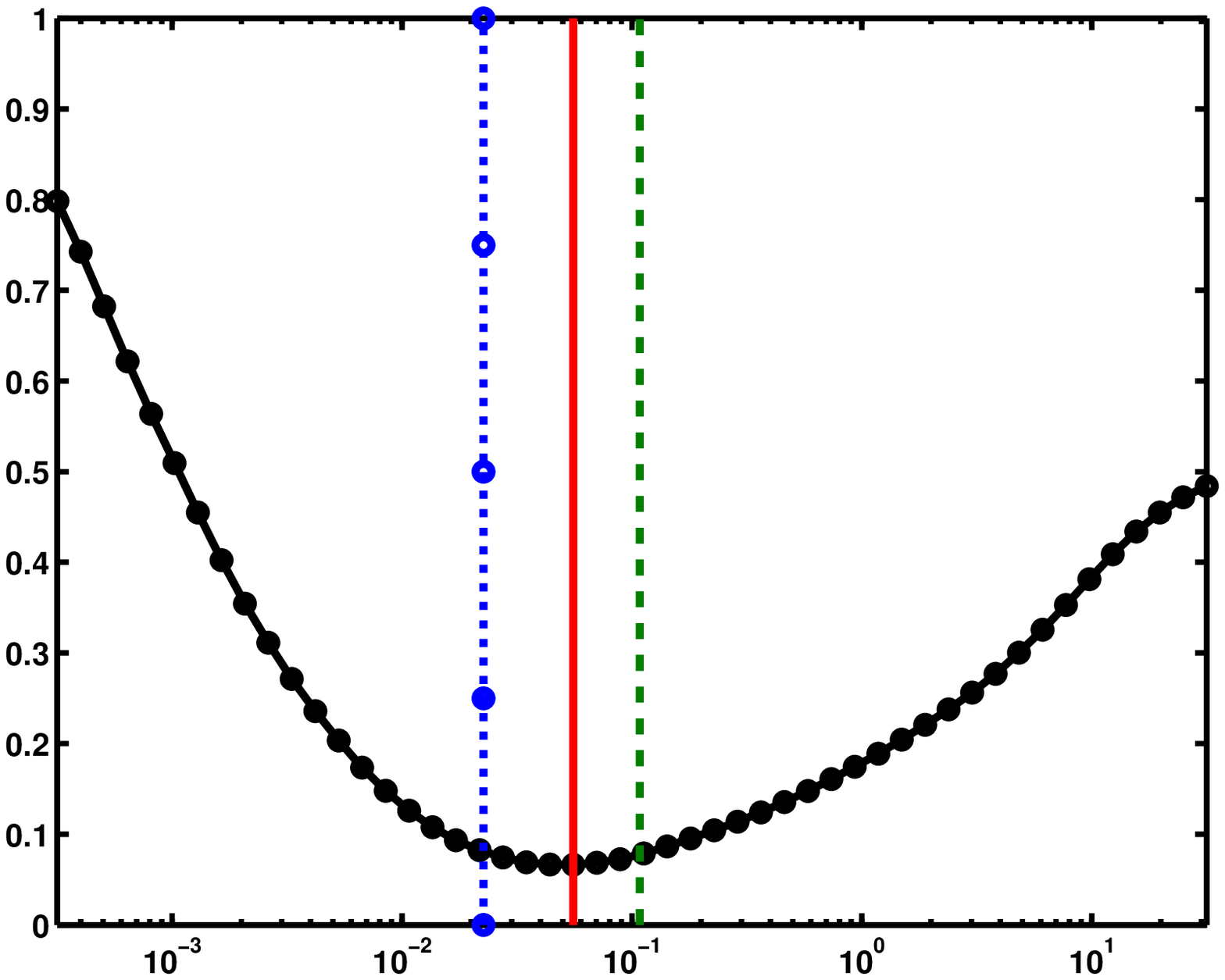}}
\subfigure[$L=L_2$]{\includegraphics[width=1.7in]{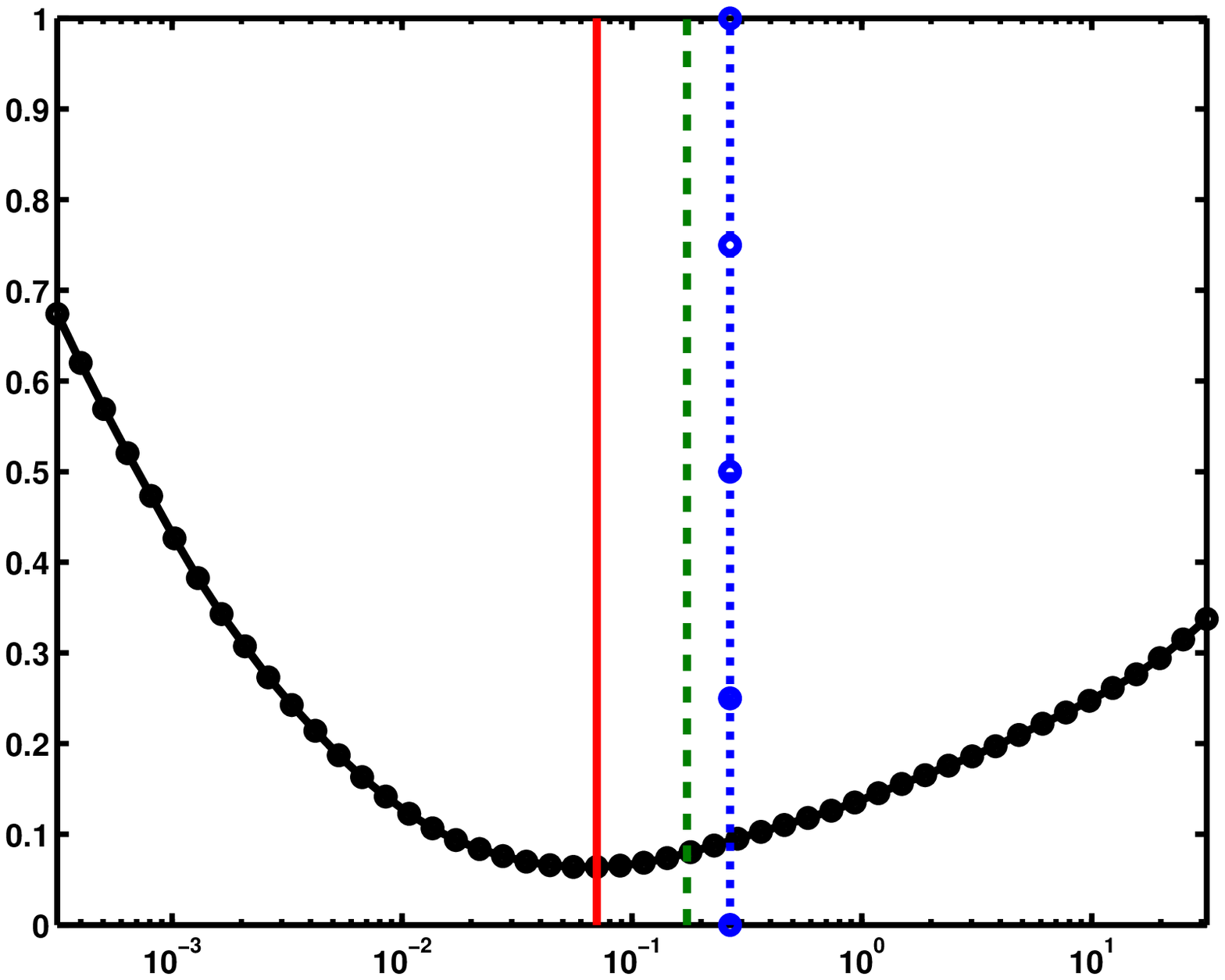}}
\subfigure[$L=I$]{\includegraphics[width=1.7in]{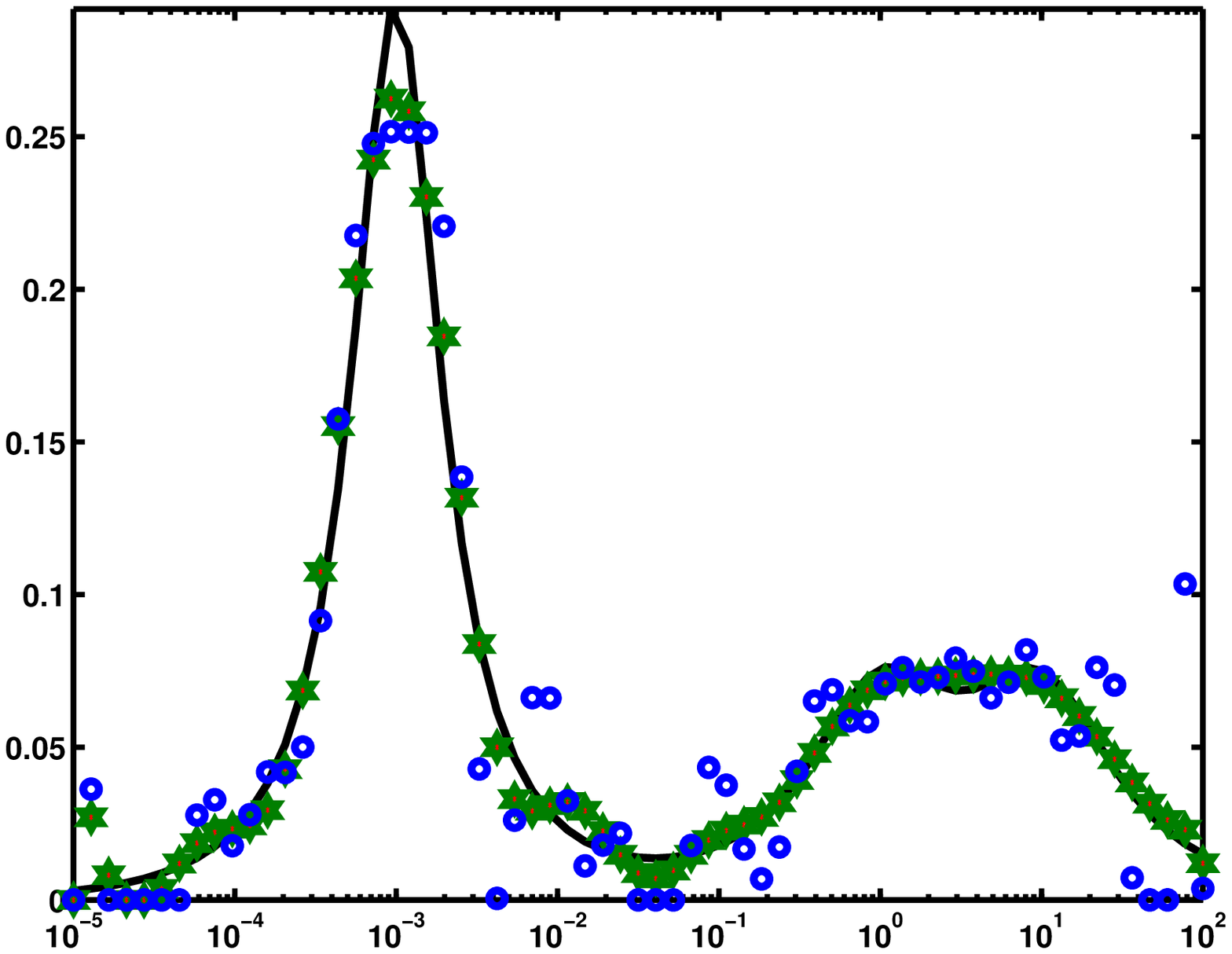}}
\subfigure[$L=L_1$]{\includegraphics[width=1.7in]{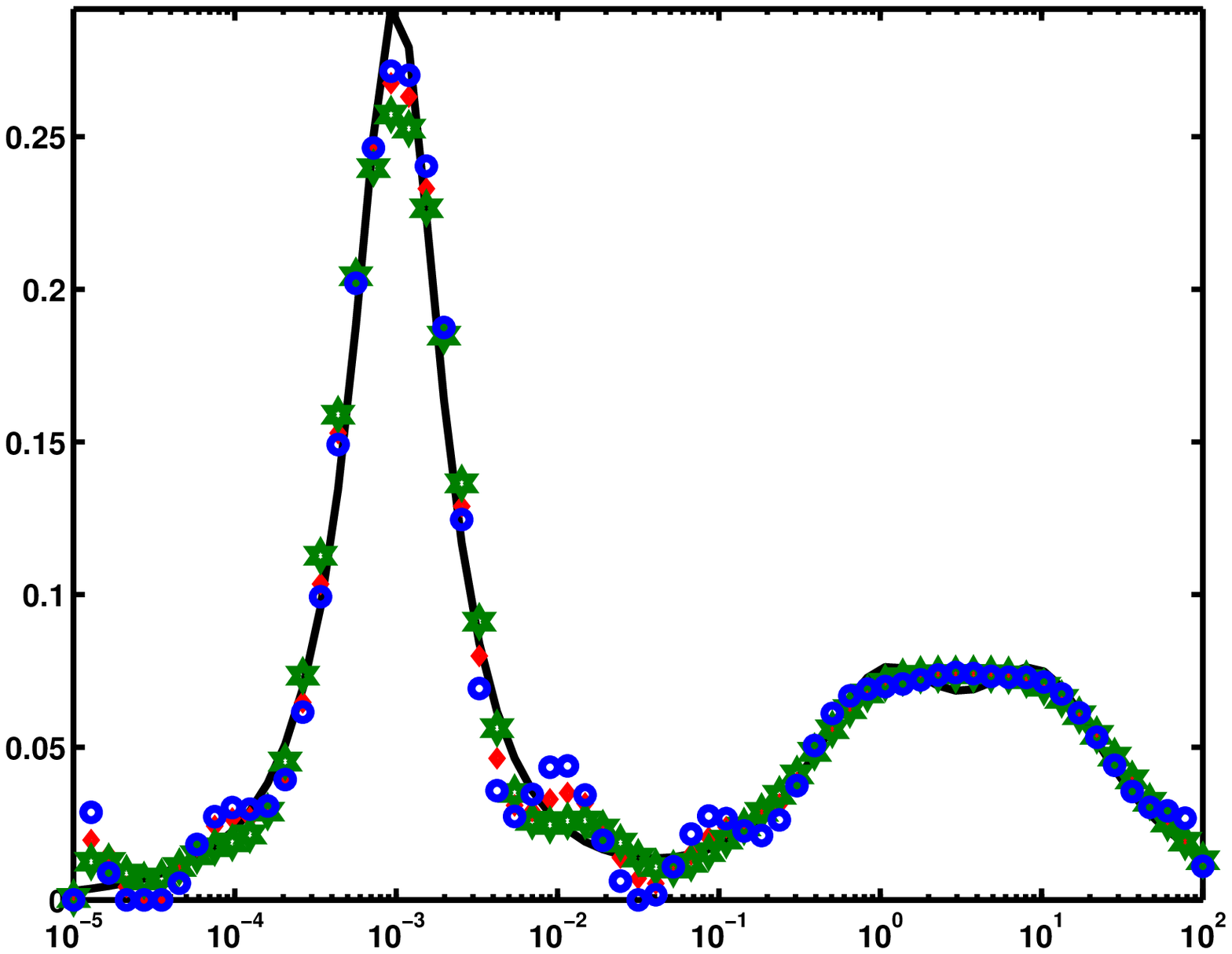}}
\subfigure[$L=L_2$]{\includegraphics[width=1.7in]{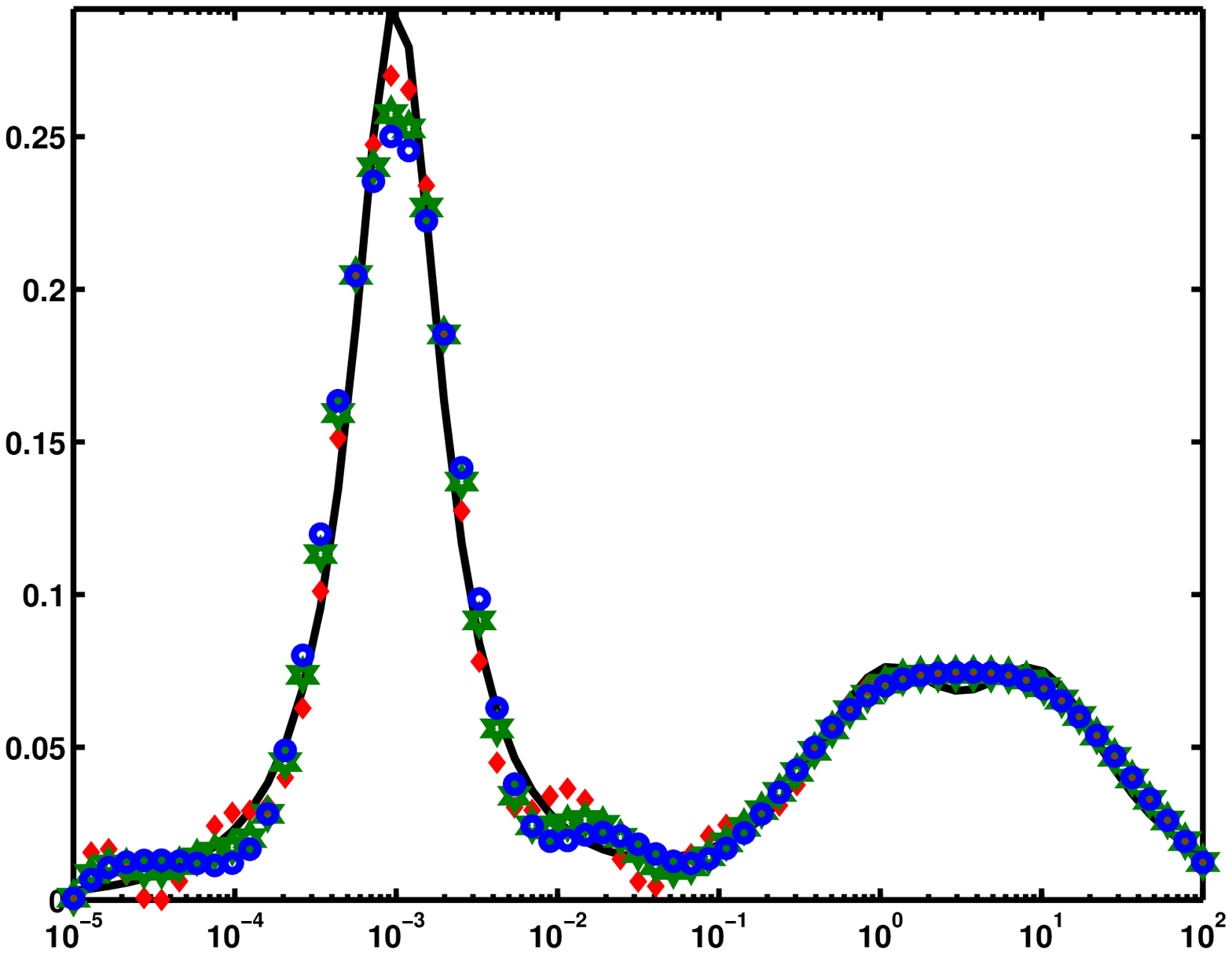}}
\caption{Mean error and example NNLS solutions.  $.1\%$ noise. RQ-C data set matrix $A_3$}
\label{fig-lambdachoiceRQ6A3LN}
\end{figure}

\begin{figure}[!h]
 \centering
\subfigure[$L=I$]{\includegraphics[width=1.7in]{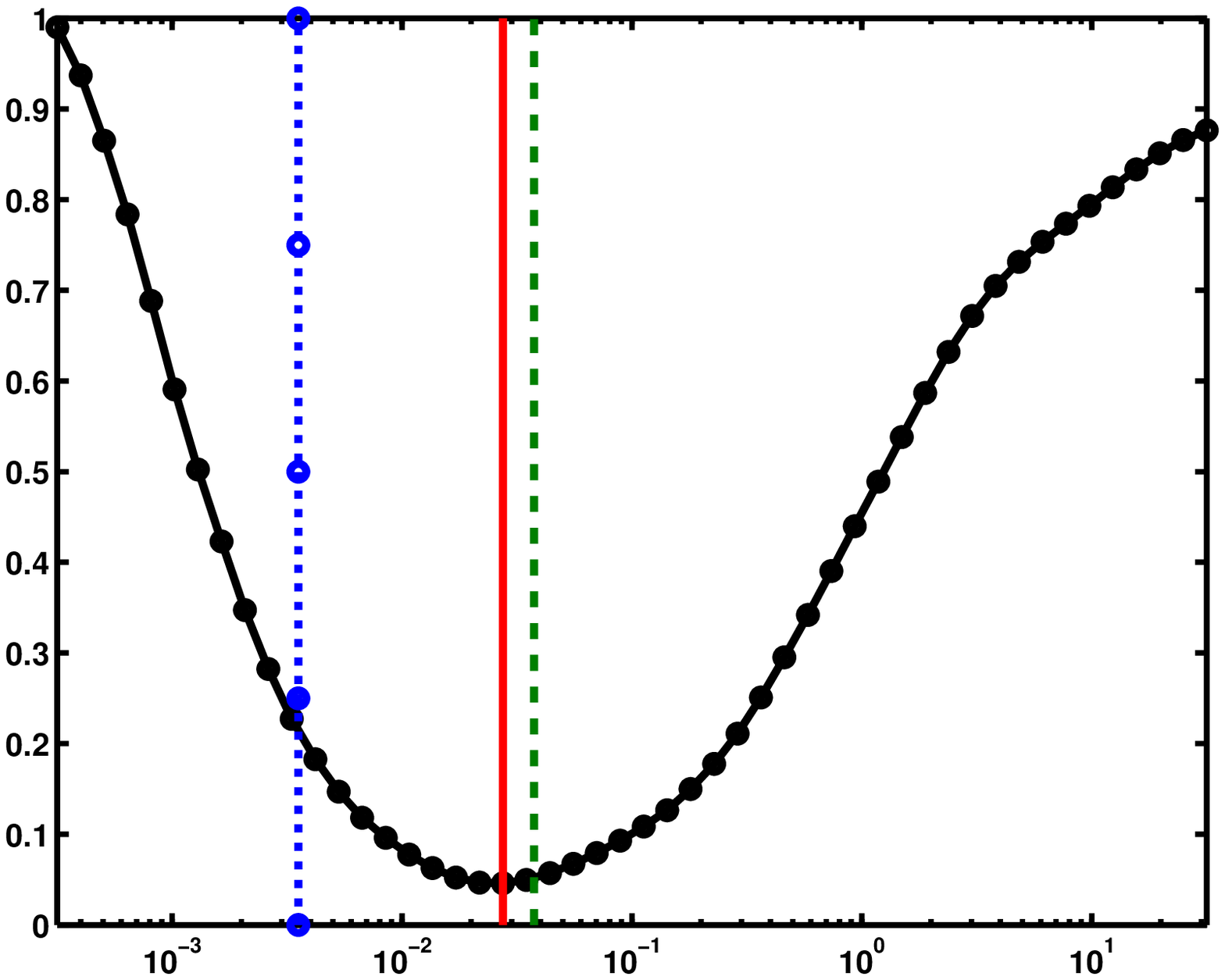}}
\subfigure[$L=L_1$]{\includegraphics[width=1.7in]{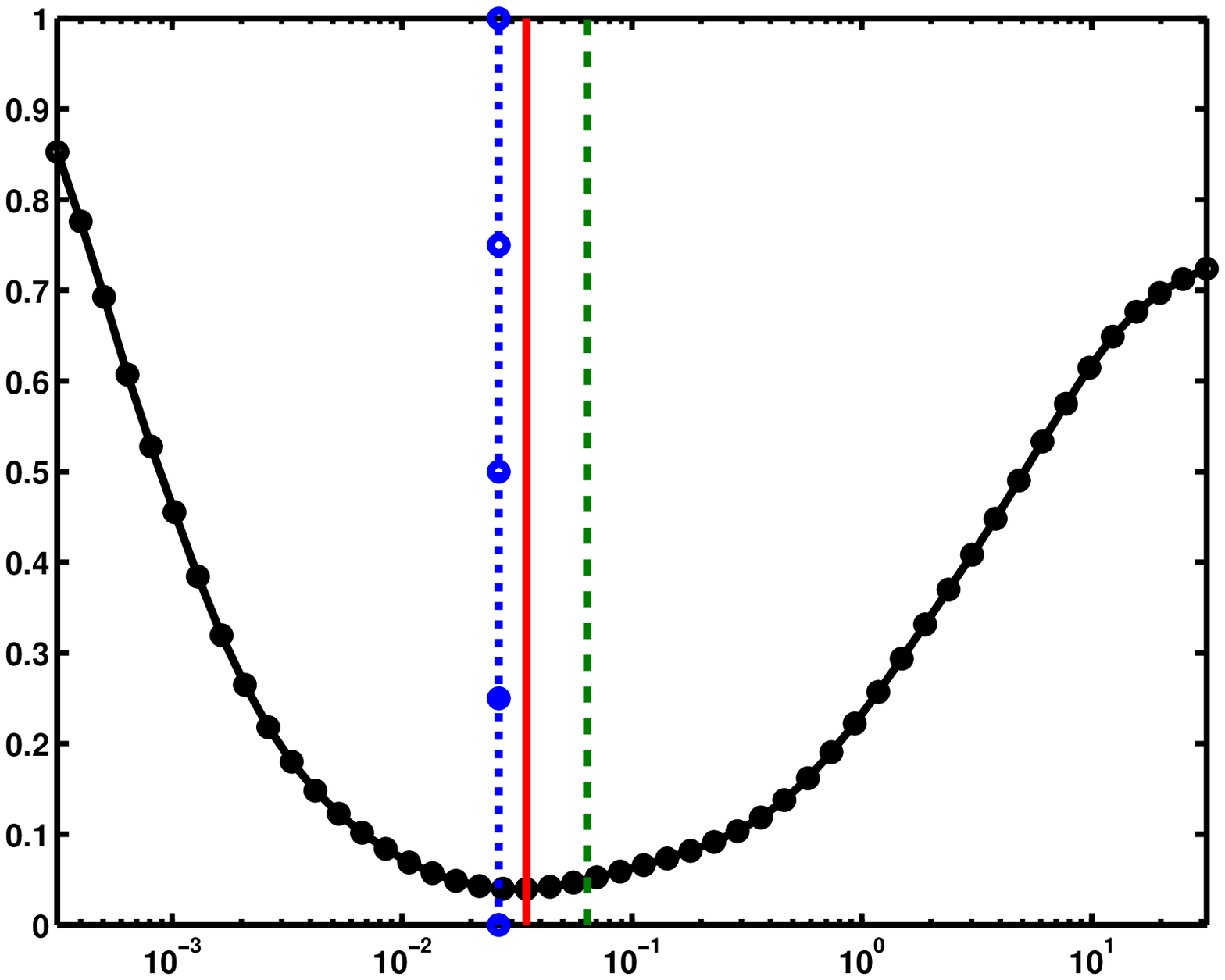}}
\subfigure[$L=L_2$]{\includegraphics[width=1.7in]{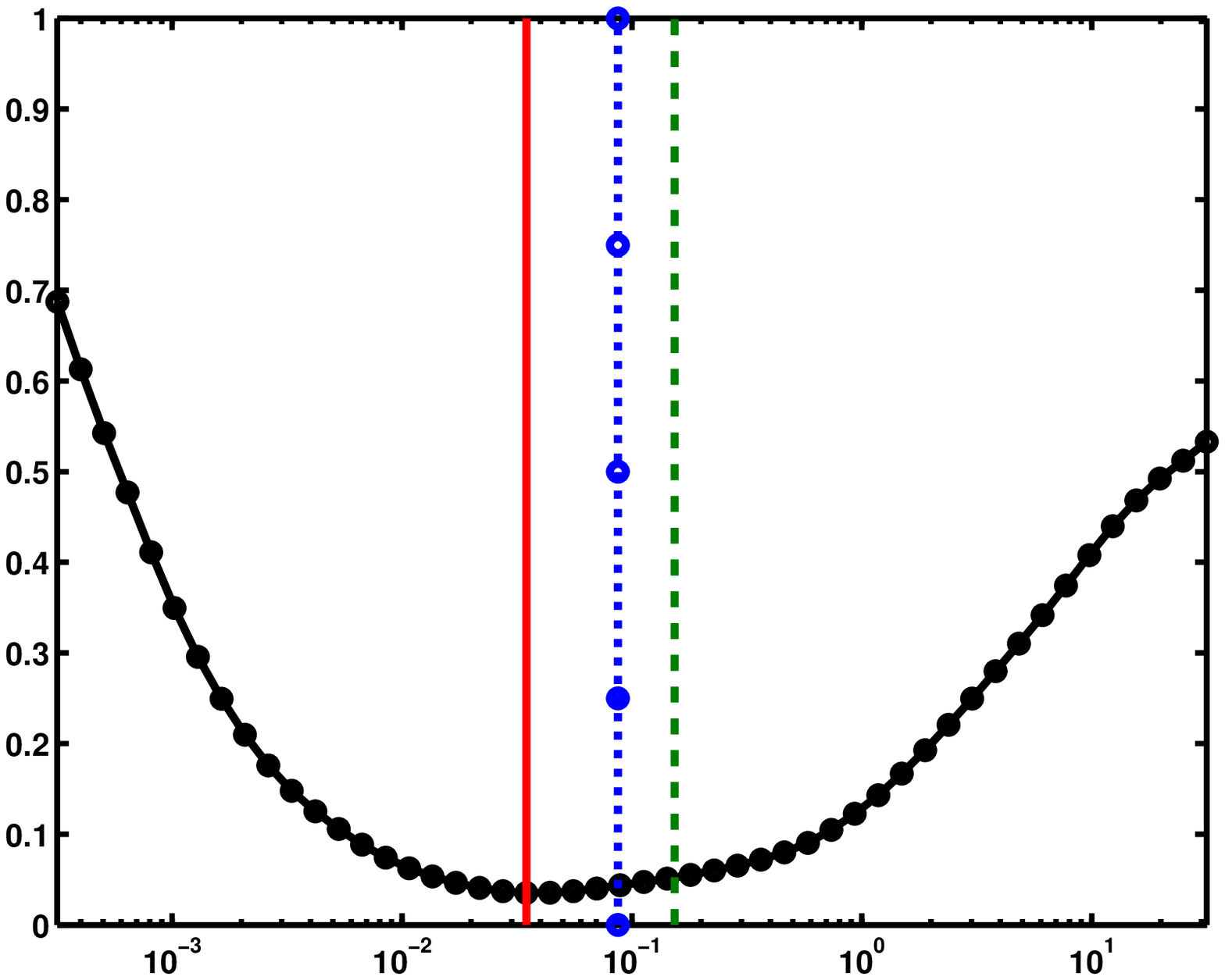}}
\subfigure[$L=I$]{\includegraphics[width=1.7in]{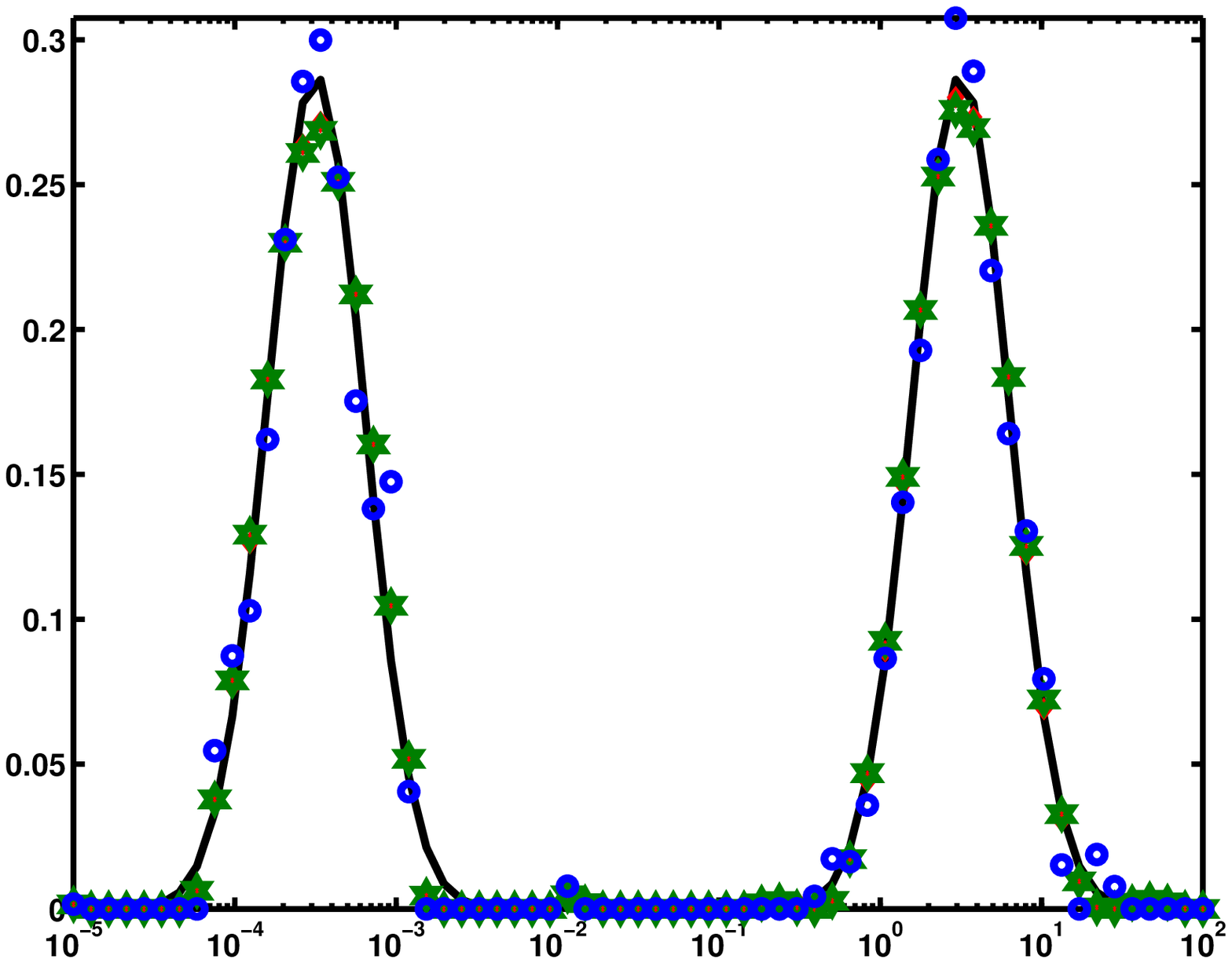}}
\subfigure[$L=L_1$]{\includegraphics[width=1.7in]{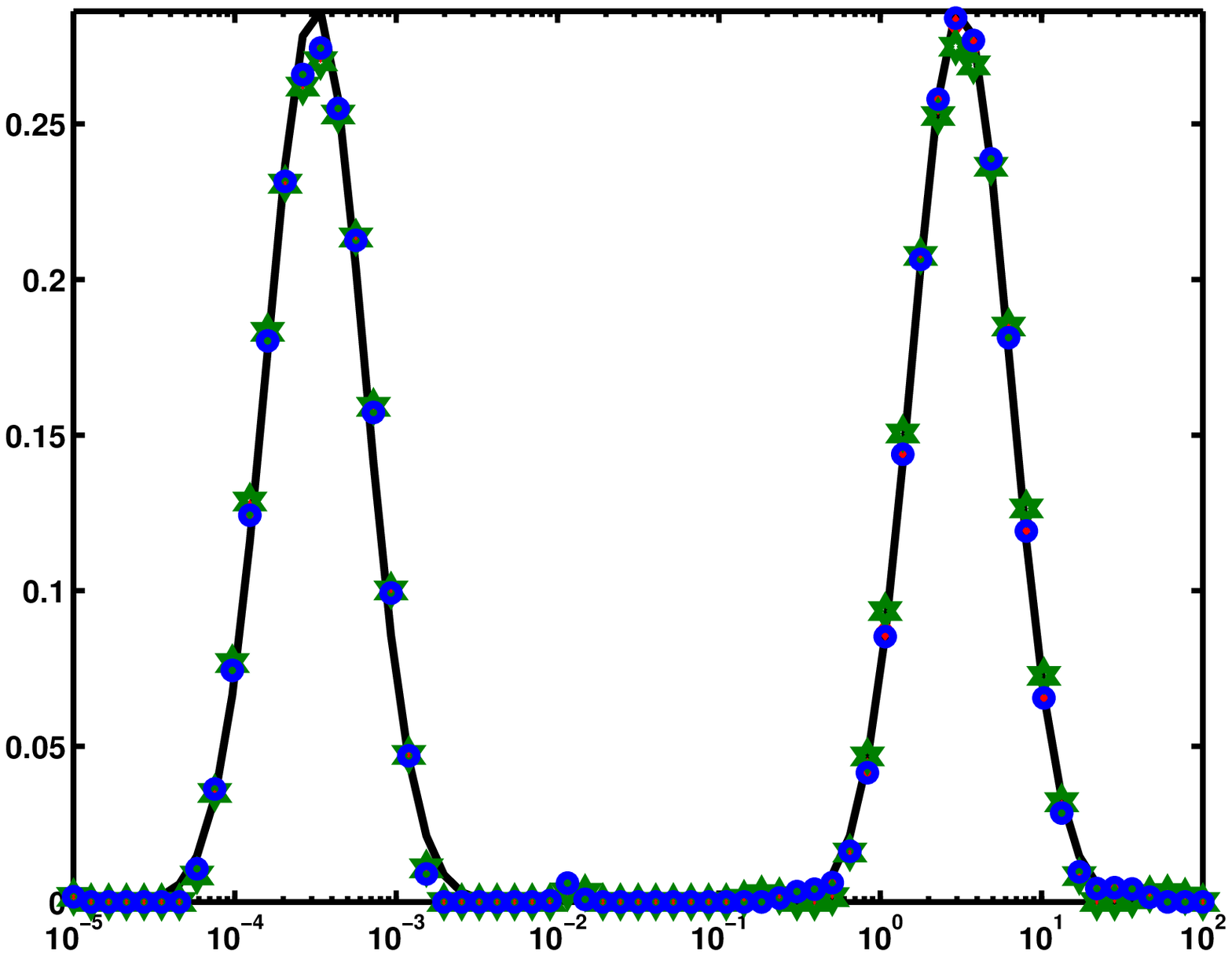}}
\subfigure[$L=L_2$]{\includegraphics[width=1.7in]{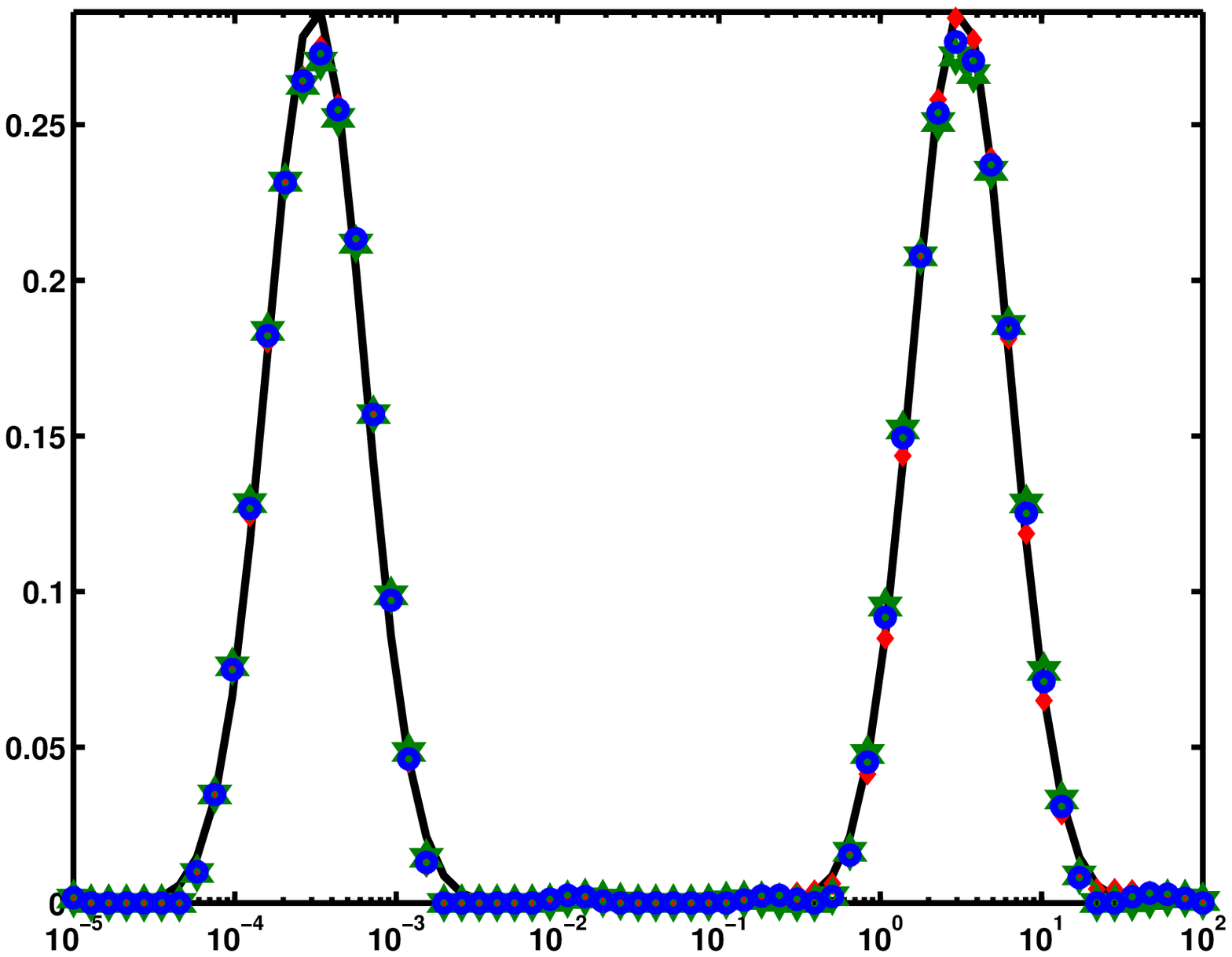}}
\caption{Mean error and example NNLS solutions.  $.1\%$ noise. LN-A data set matrix $A_3$}
\label{fig-lambdachoiceLN2A3LN}
\end{figure}

\begin{figure}[!h]
 \centering
\subfigure[$L=I$]{\includegraphics[width=1.7in]{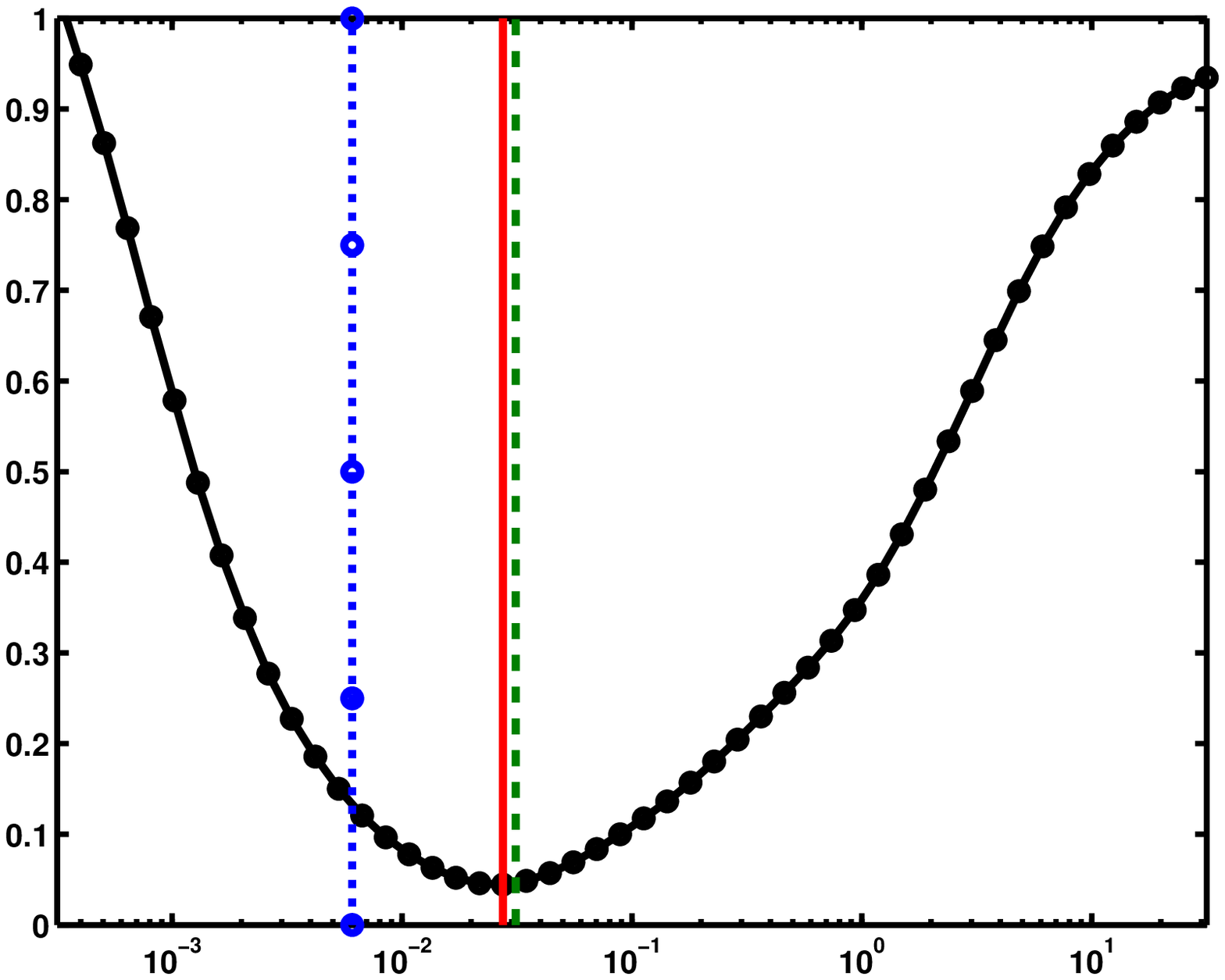}}
\subfigure[$L=L_1$]{\includegraphics[width=1.7in]{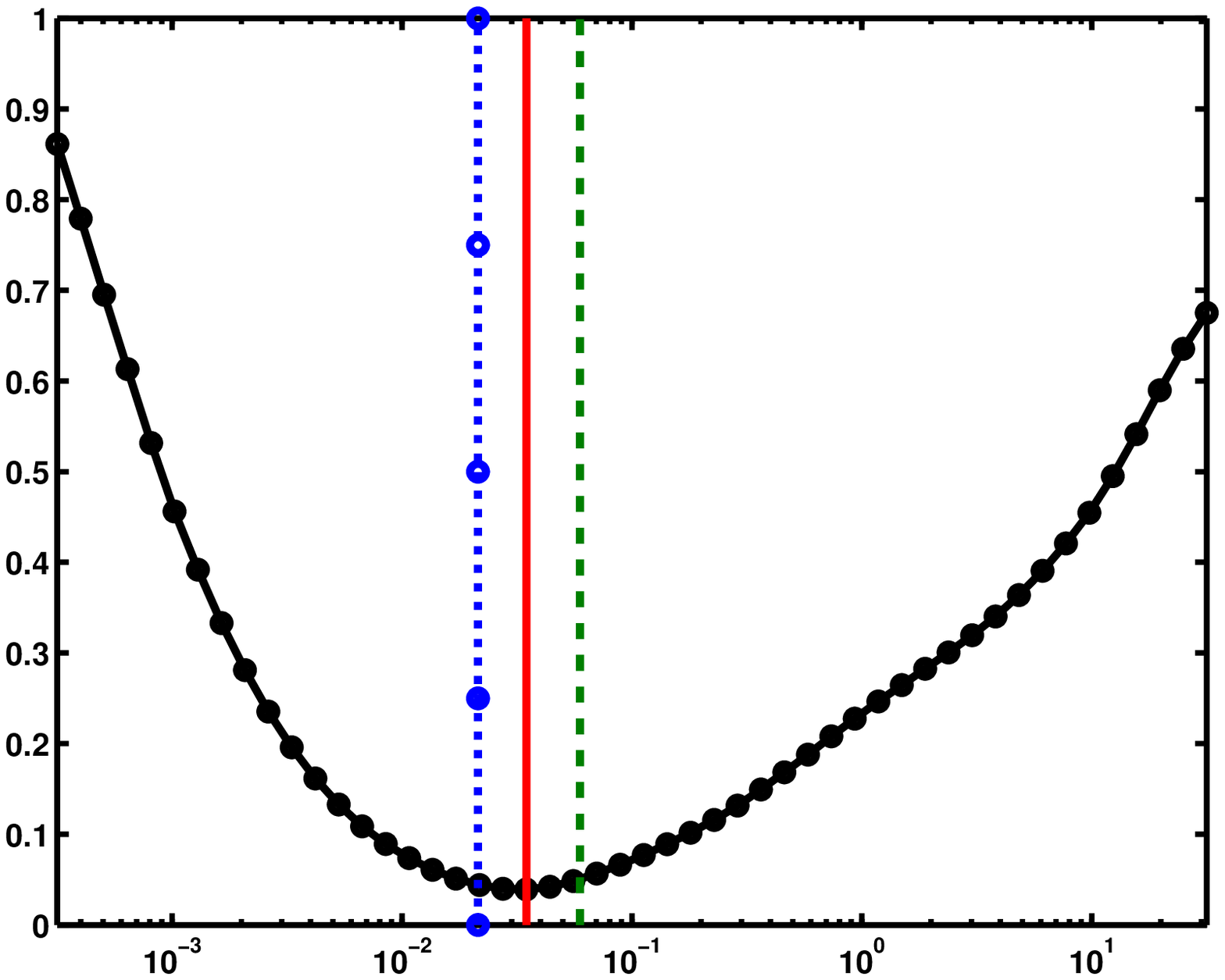}}
\subfigure[$L=L_2$]{\includegraphics[width=1.7in]{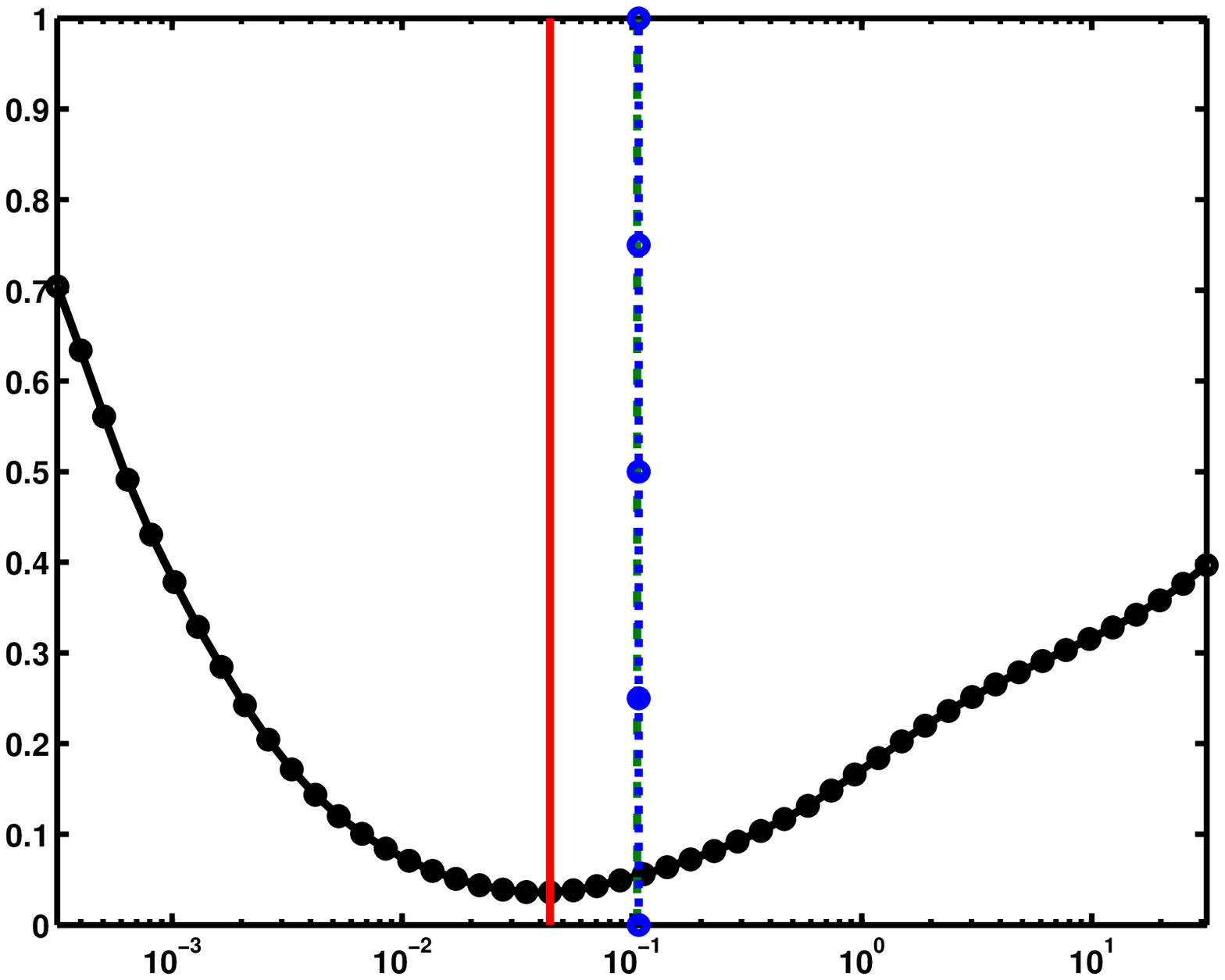}}
\subfigure[$L=I$]{\includegraphics[width=1.7in]{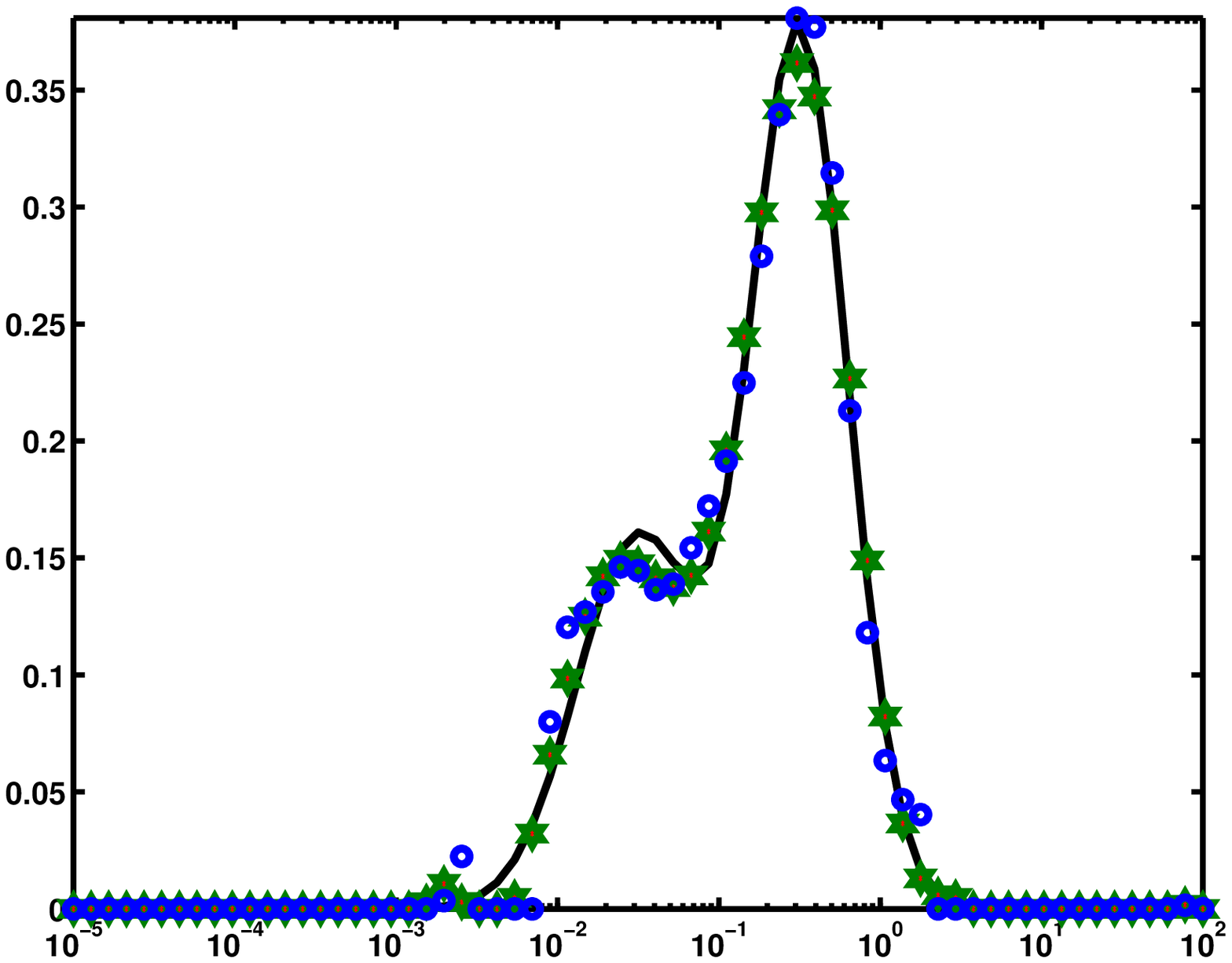}}
\subfigure[$L=L_1$]{\includegraphics[width=1.7in]{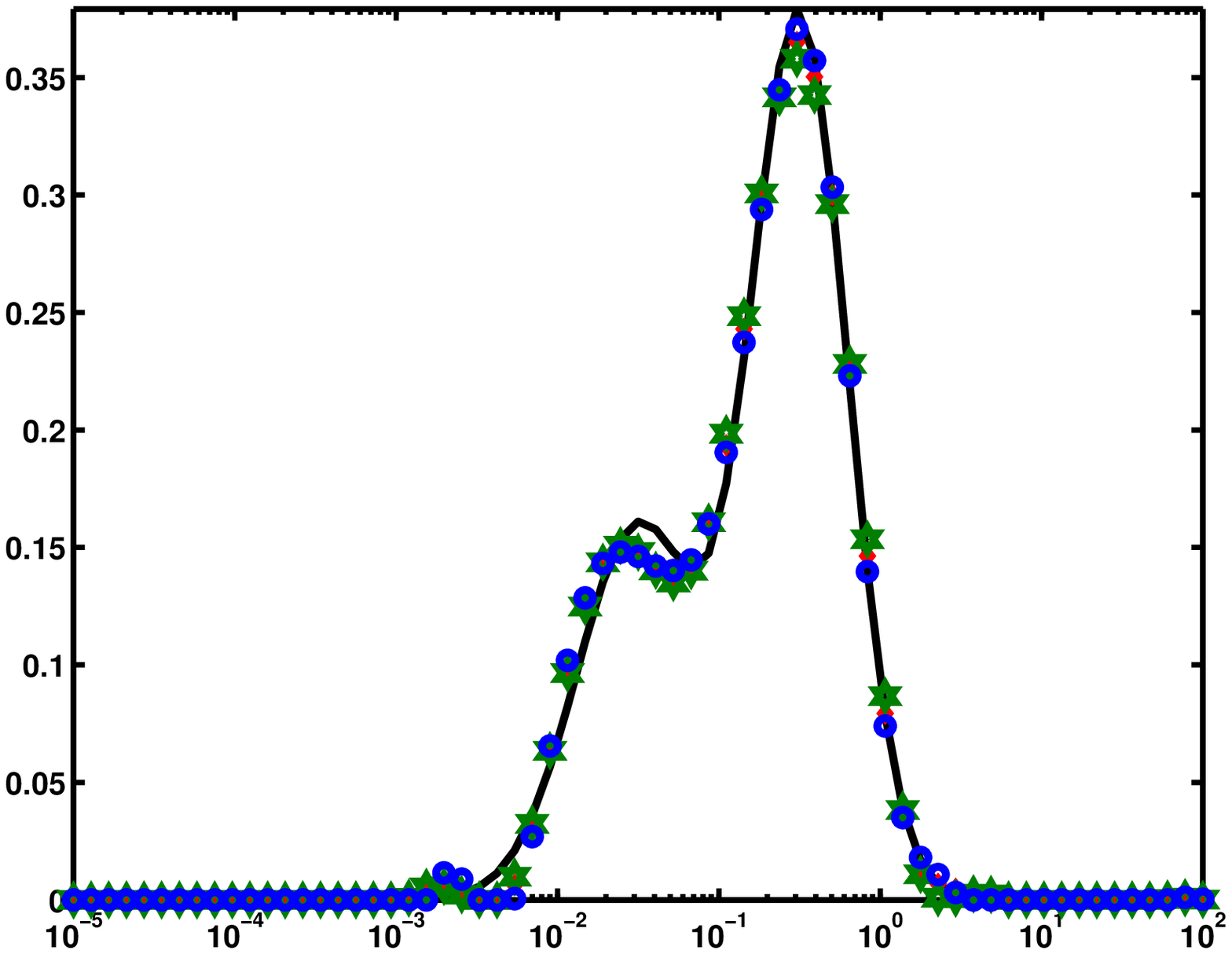}}
\subfigure[$L=L_2$]{\includegraphics[width=1.7in]{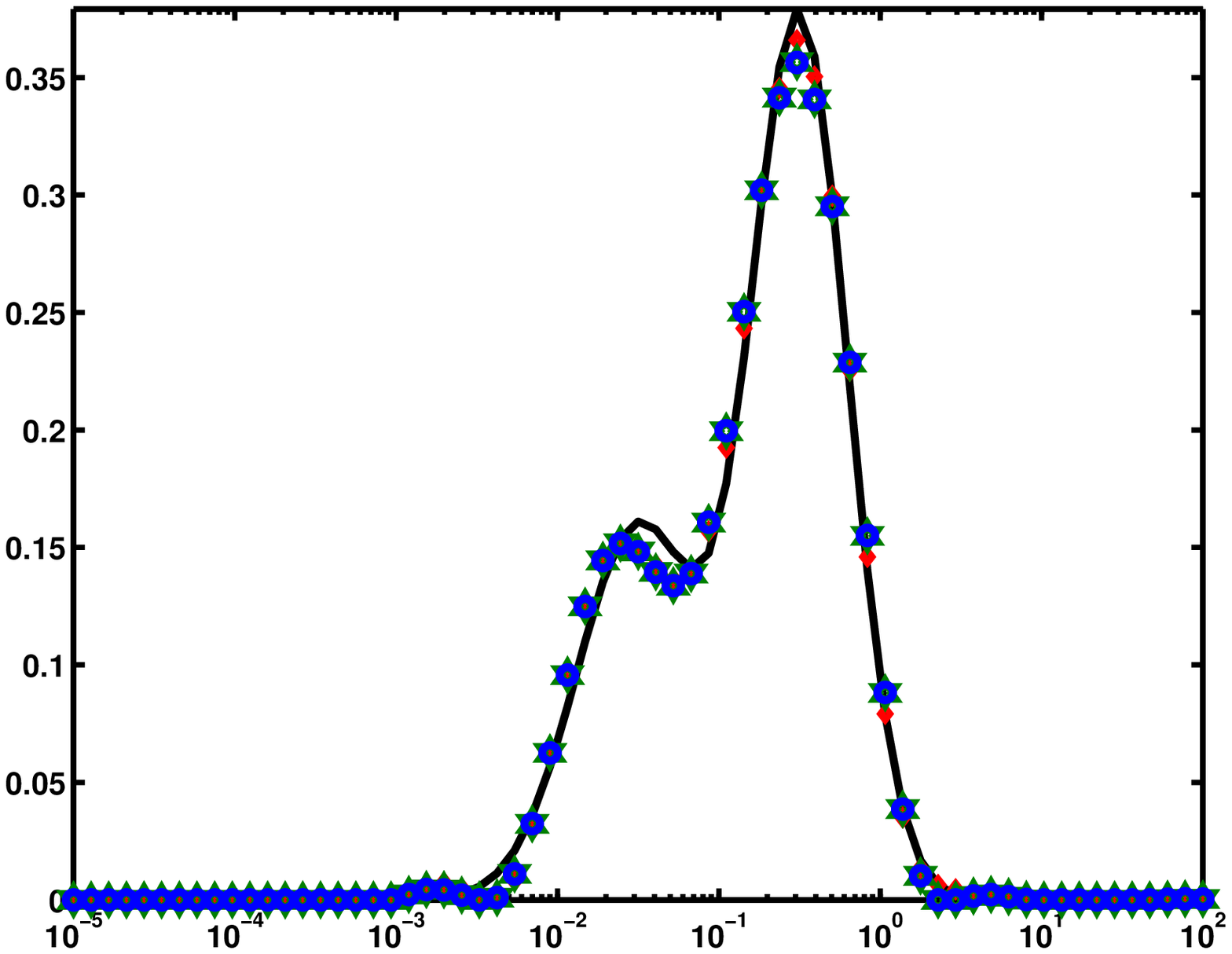}}
\caption{Mean error and example NNLS solutions.  $.1\%$ noise. LN-B data set matrix $A_3$}
\label{fig-lambdachoiceLN5A3LN}
\end{figure}

\begin{figure}[!h]
 \centering
\subfigure[$L=I$]{\includegraphics[width=1.7in]{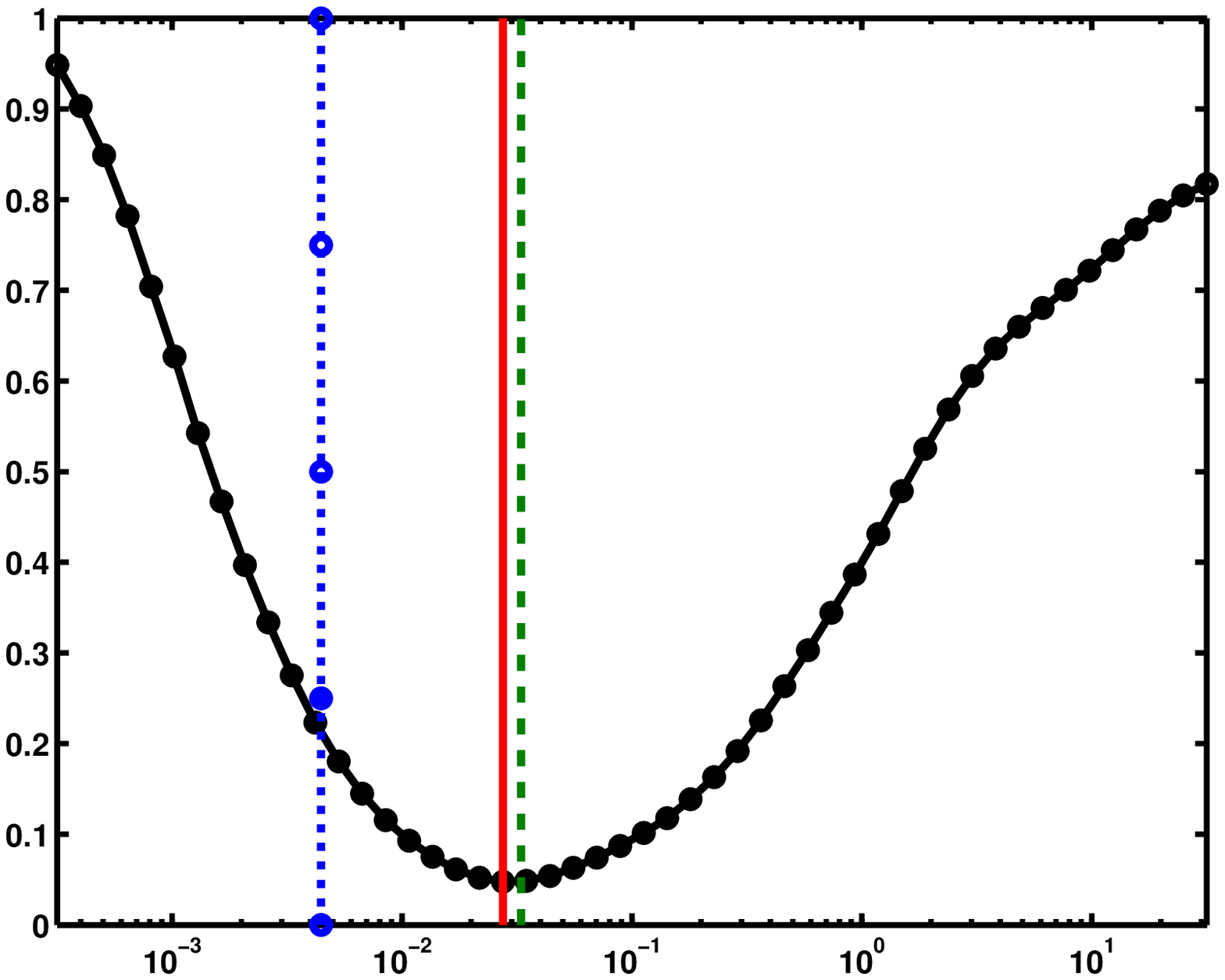}}
\subfigure[$L=L_1$]{\includegraphics[width=1.7in]{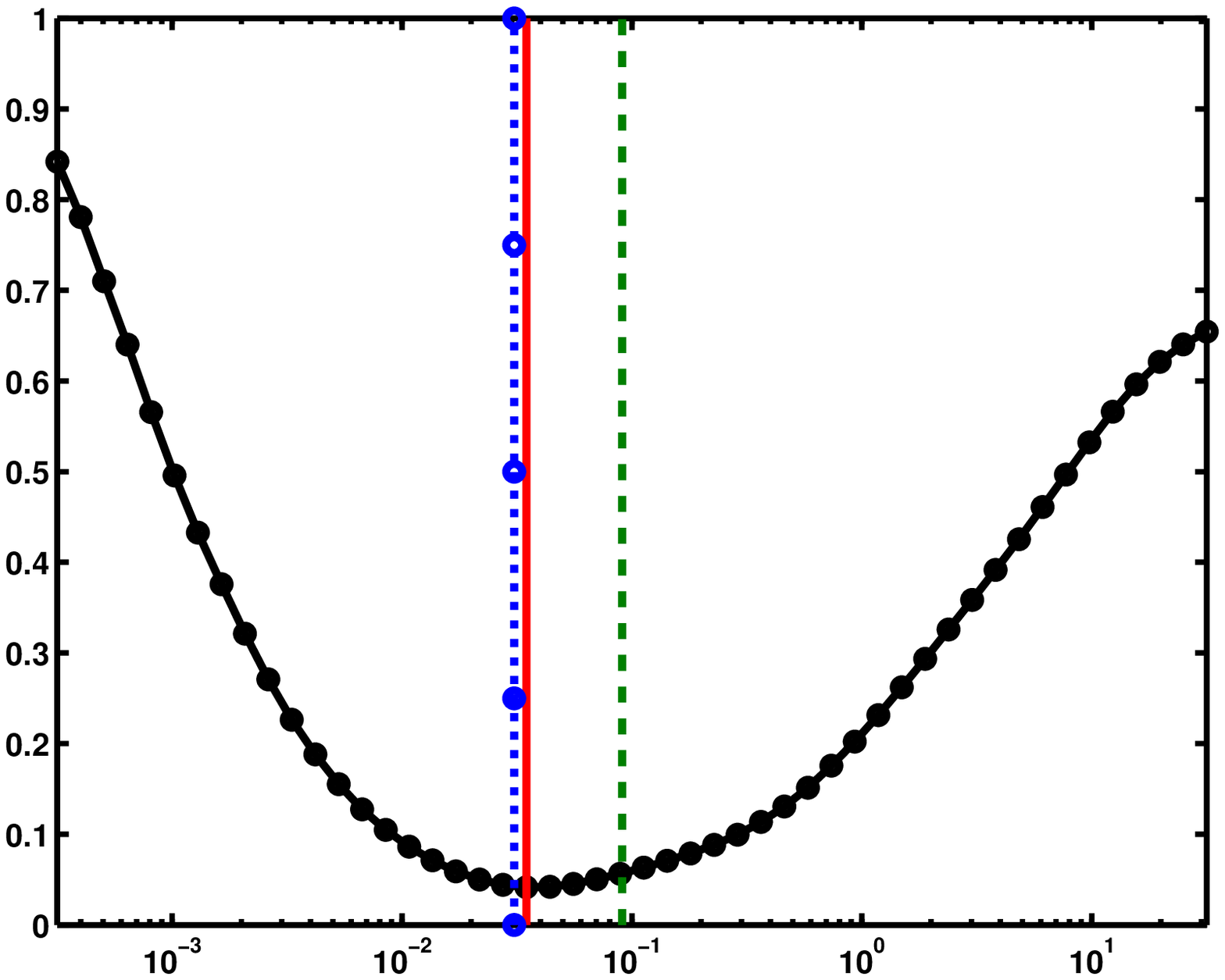}}
\subfigure[$L=L_2$]{\includegraphics[width=1.7in]{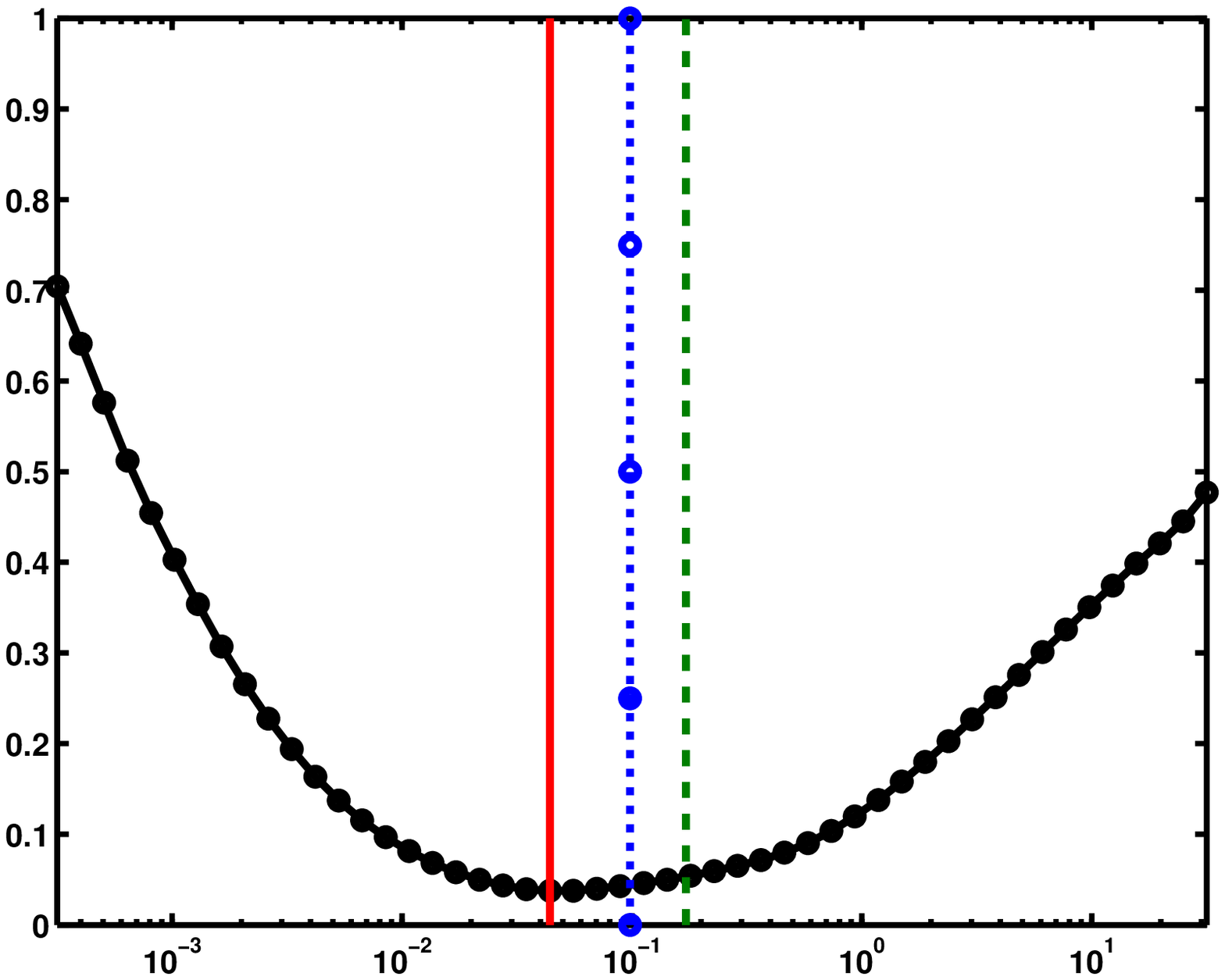}}
\subfigure[$L=I$]{\includegraphics[width=1.7in]{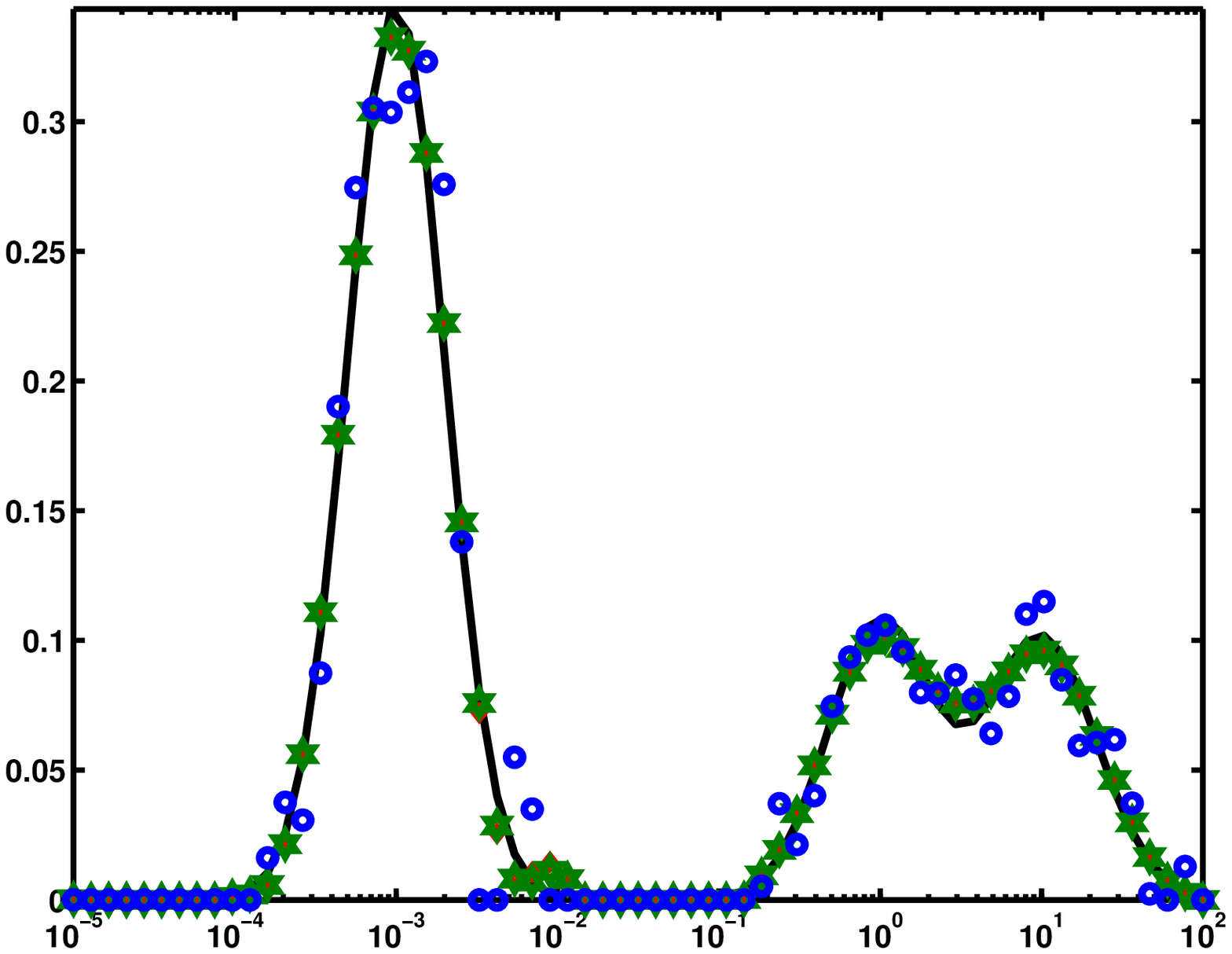}}
\subfigure[$L=L_1$]{\includegraphics[width=1.7in]{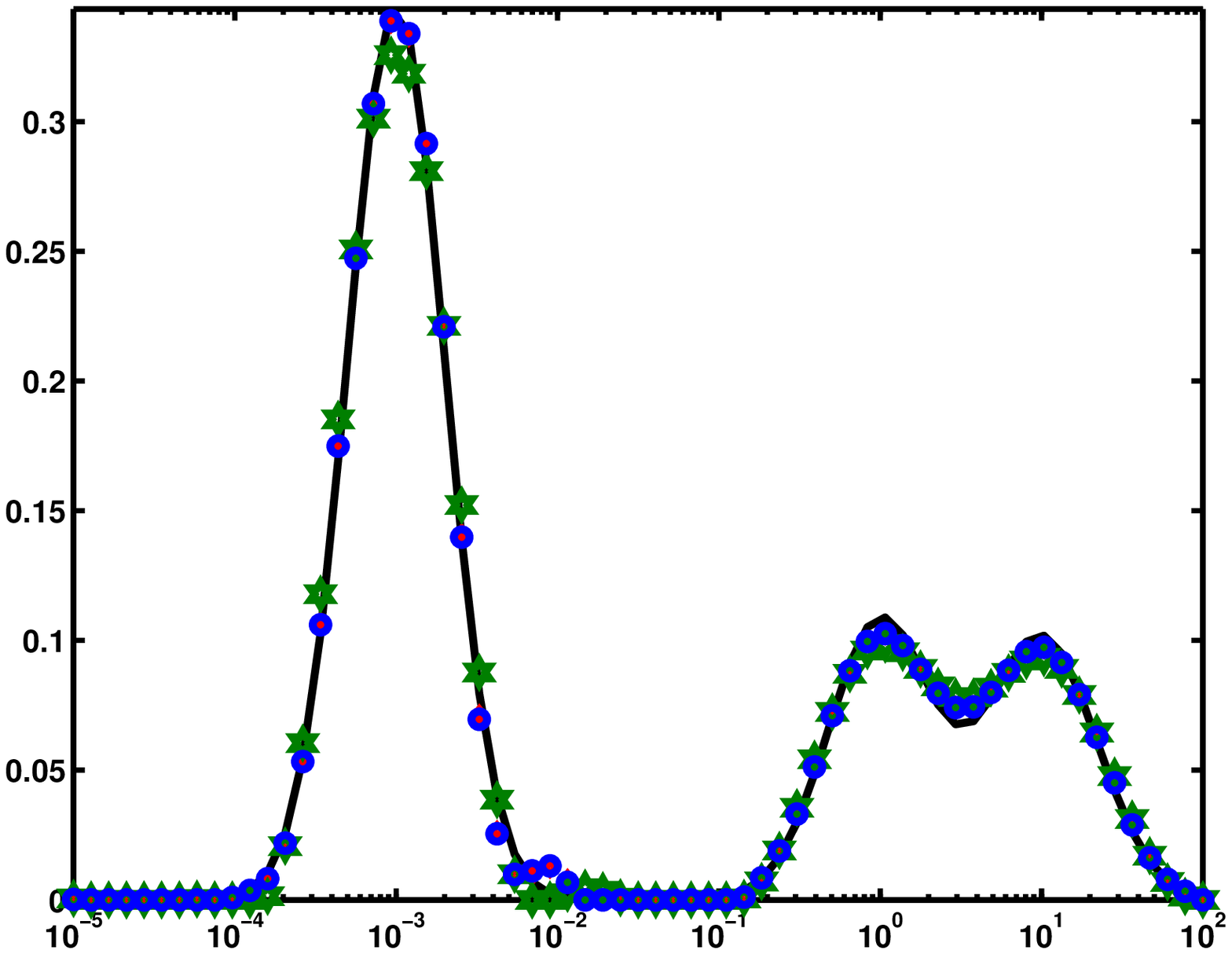}}
\subfigure[$L=L_2$]{\includegraphics[width=1.7in]{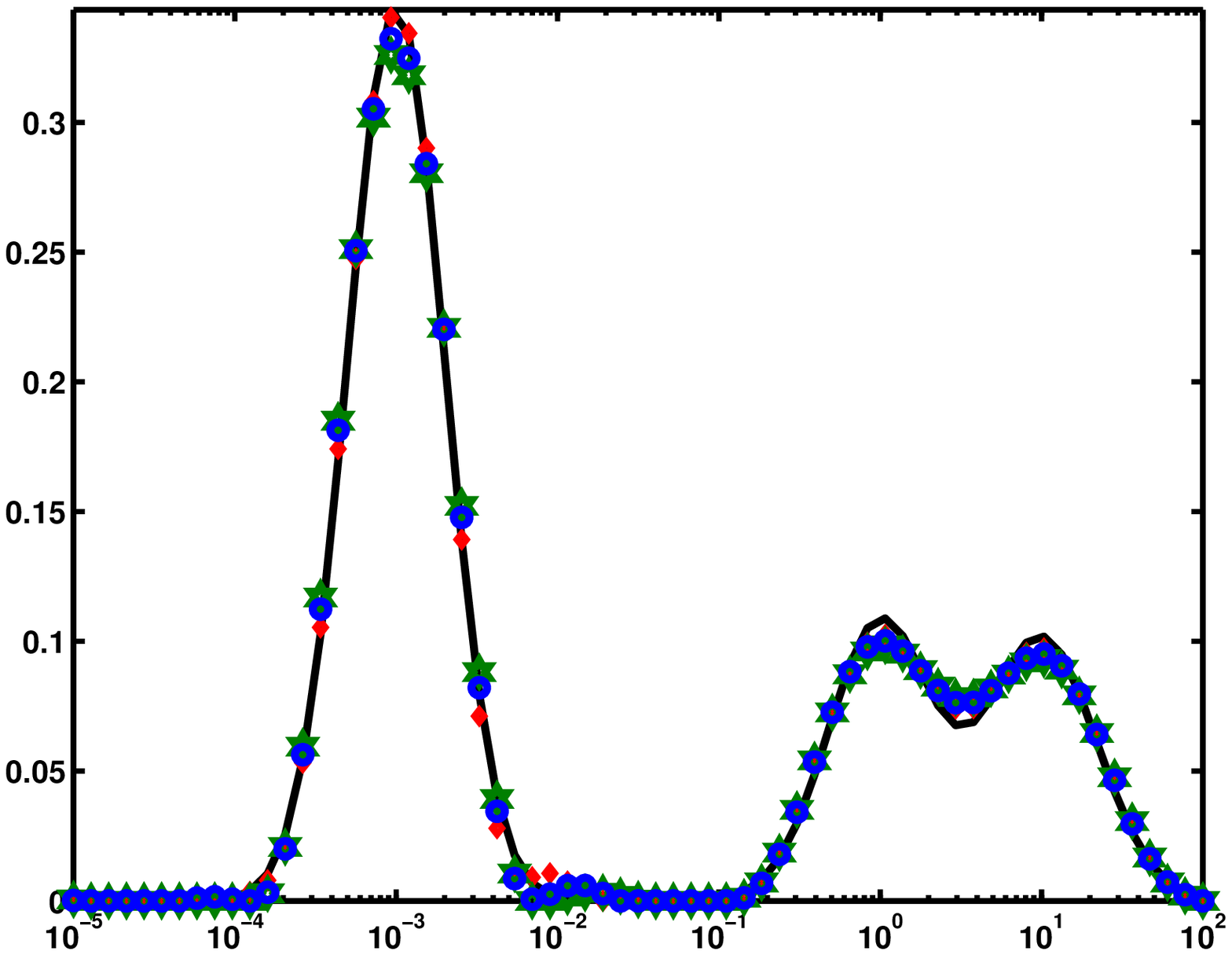}}
\caption{Mean error and example NNLS solutions.  $.1\%$ noise. LN-C data set matrix $A_3$}
\label{fig-lambdachoiceLN6A3LN}
\end{figure}

\clearpage

\subsection{Examples: Noise level $.1\%$ matrix  $A_4$ NNLS with SBB Algorithm}\label{sec:SBB}


\clearpage

\subsection{Examples: Noise level $1\%$ matrix  $A_4$ NNLS}
 \begin{figure}[!ht]
\centering
\subfigure[$L=I$]{\includegraphics[width=1.7in]{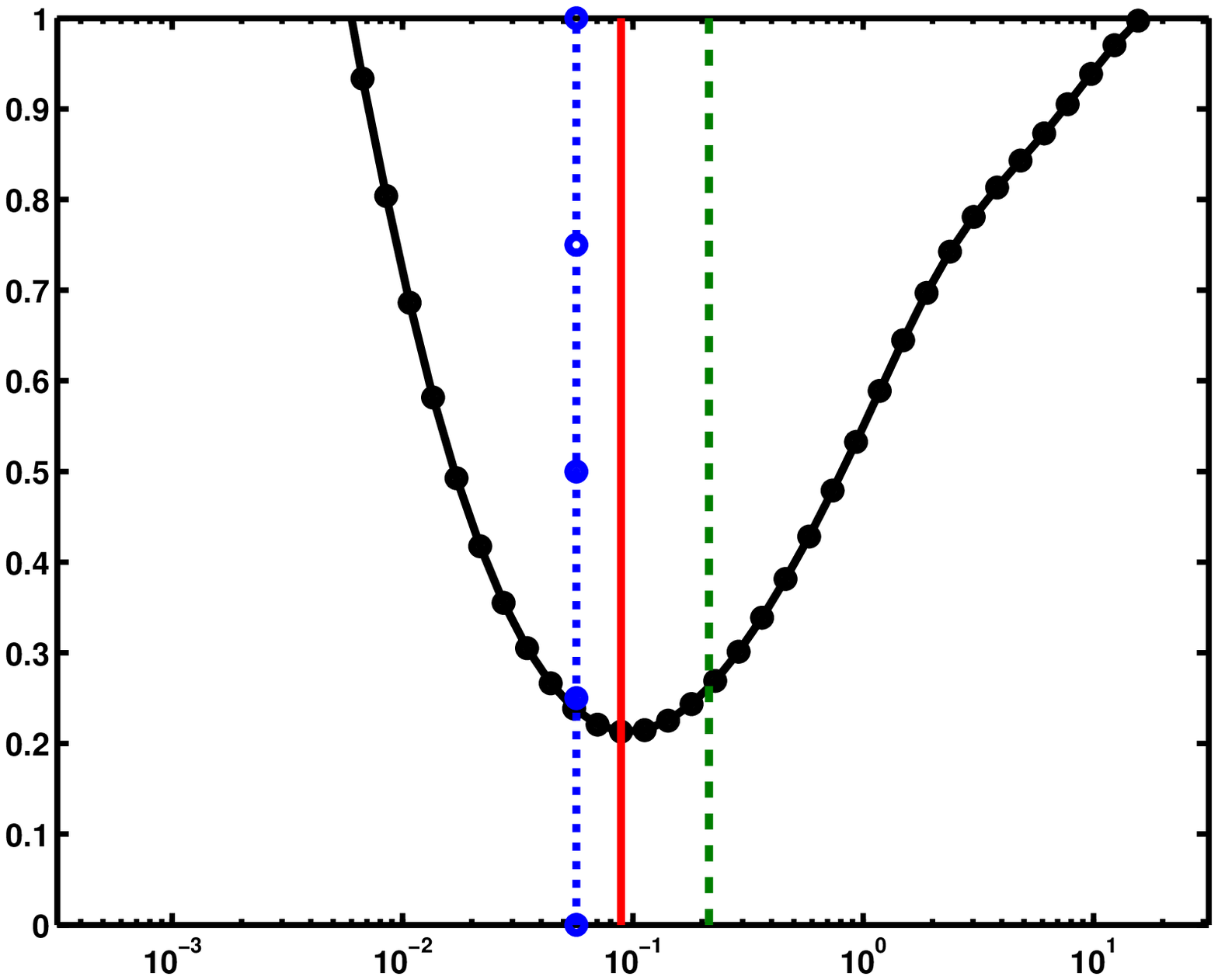}}
\subfigure[$L=L_1$]{\includegraphics[width=1.7in]{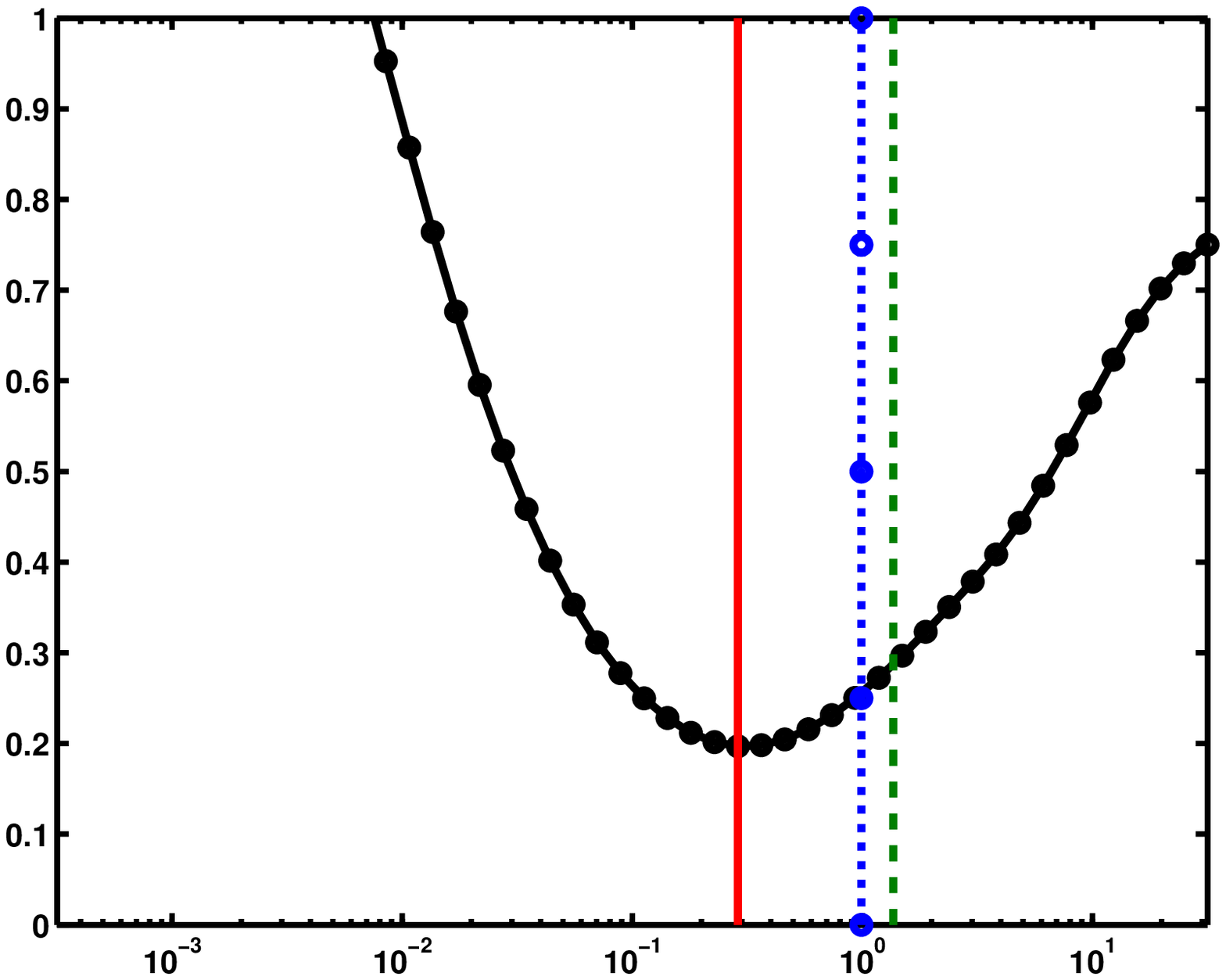}}
\subfigure[$L=L_2$]{\includegraphics[width=1.7in]{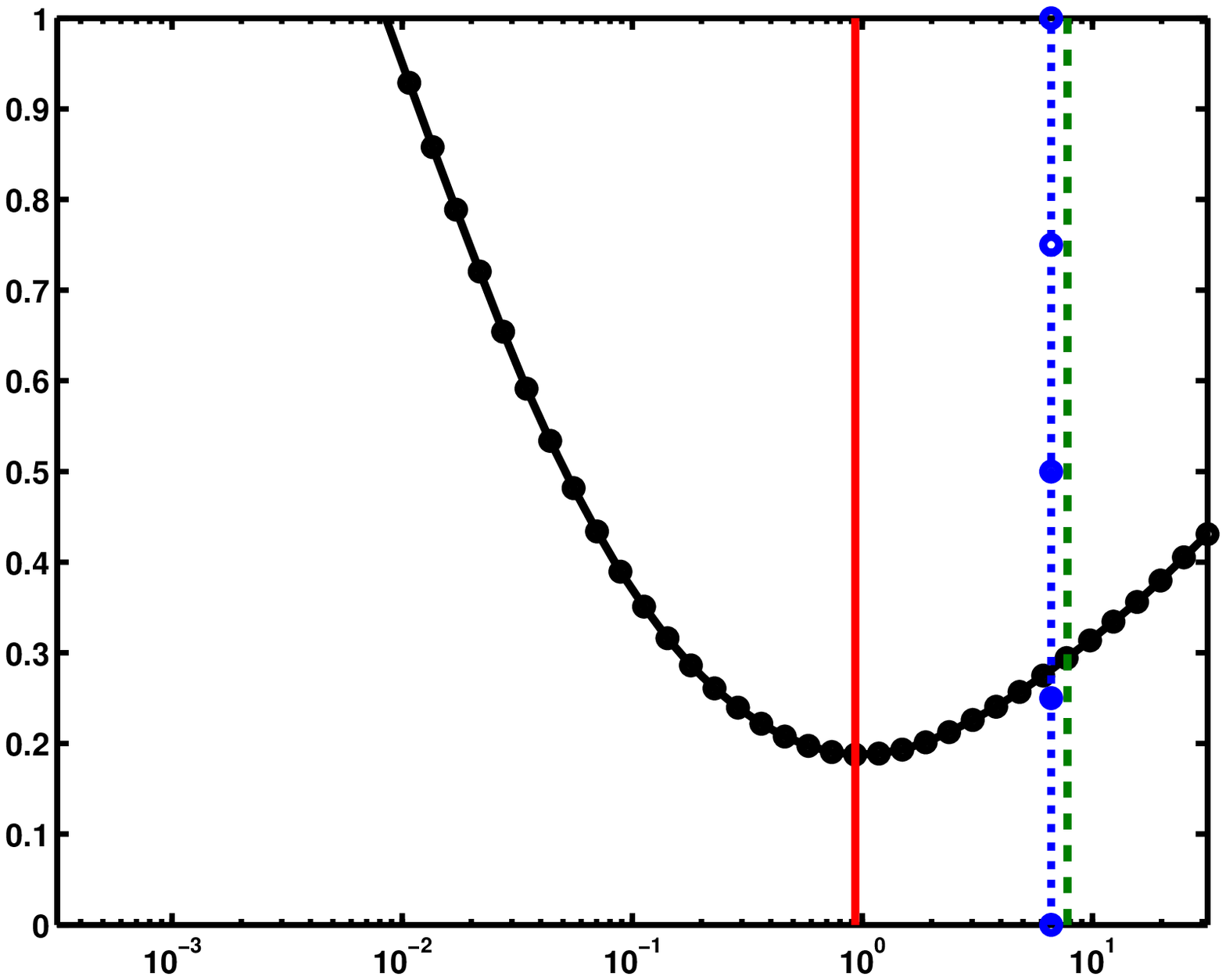}}
\subfigure[$L=I$]{\includegraphics[width=1.7in]{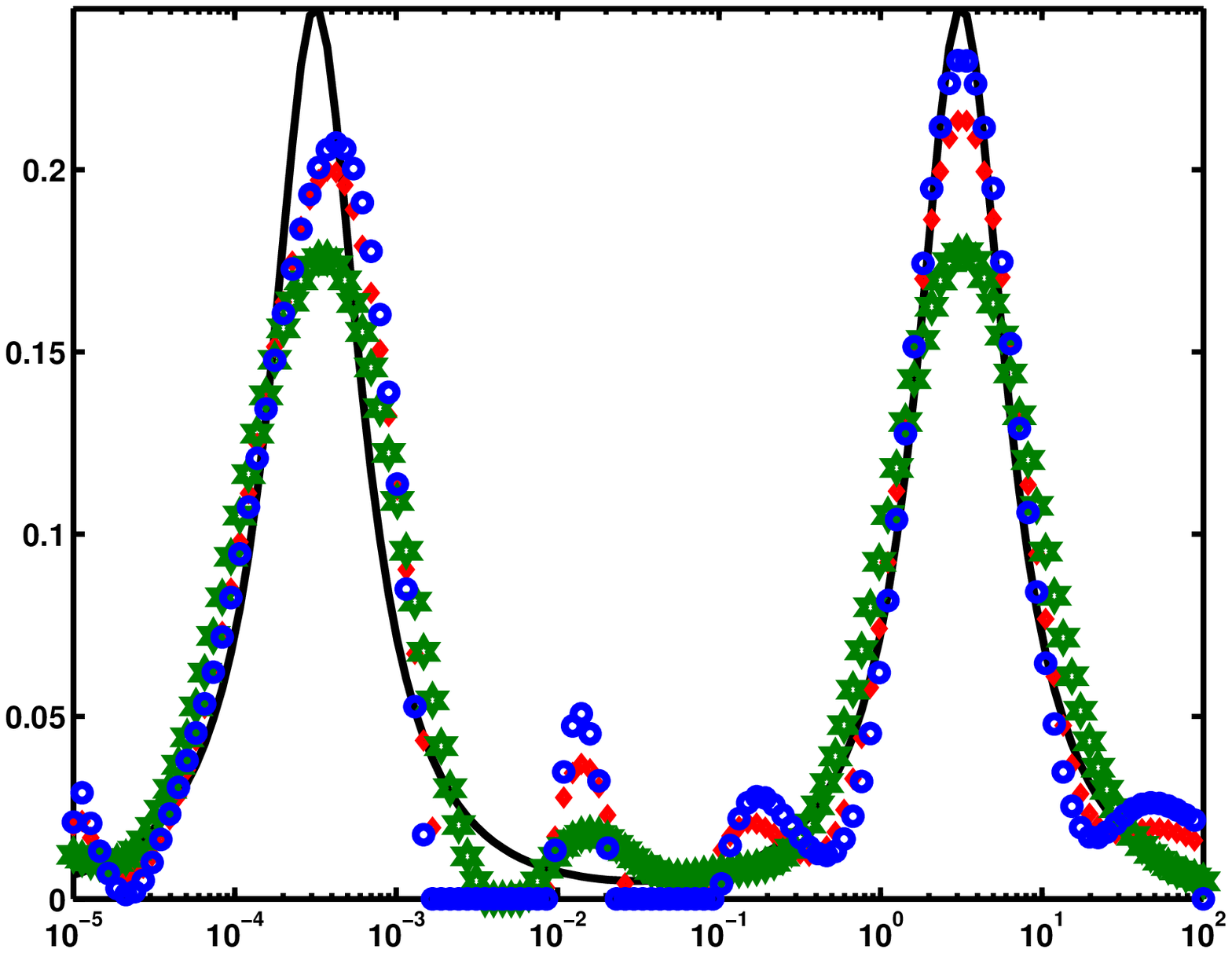}}
\subfigure[$L=L_1$]{\includegraphics[width=1.7in]{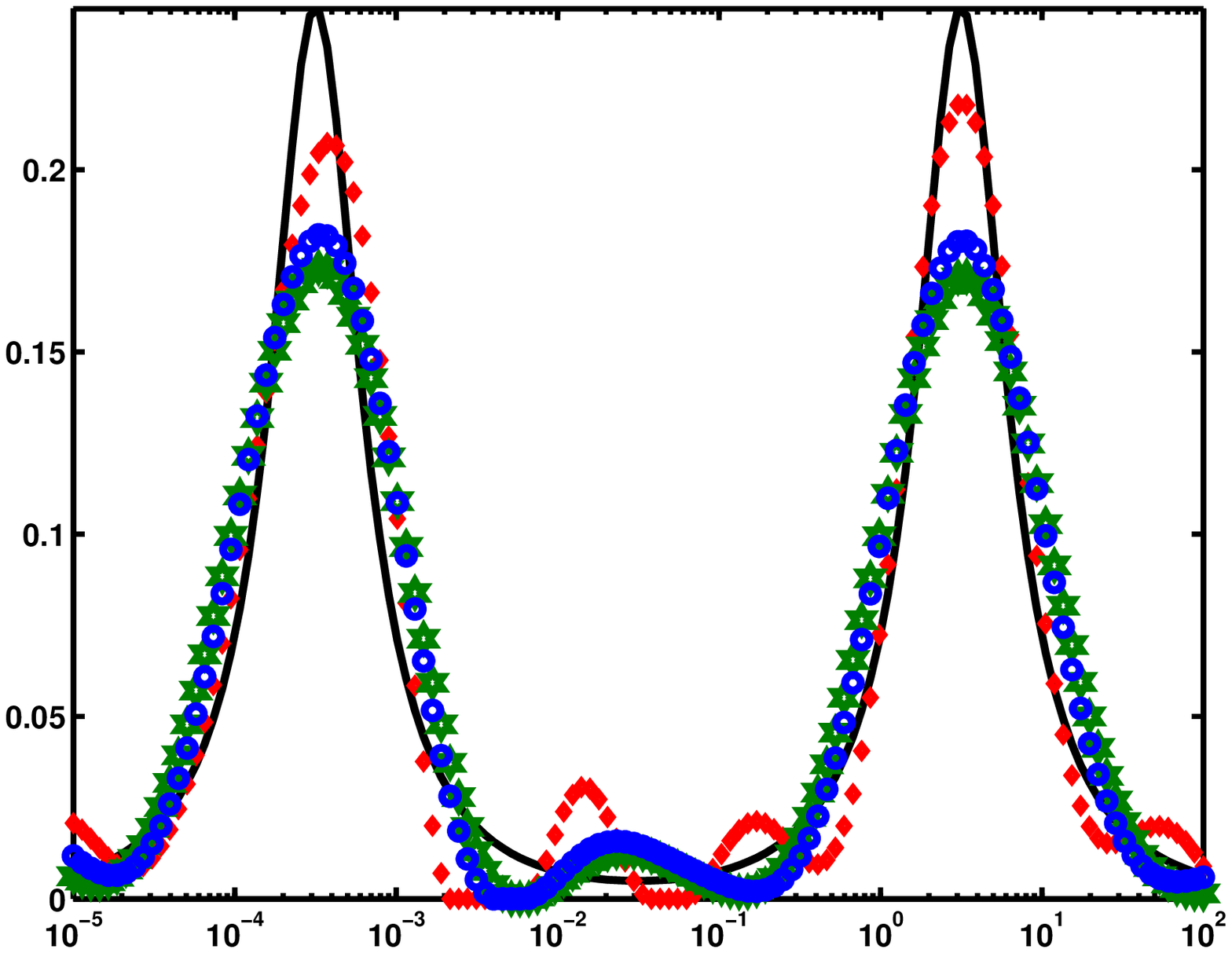}}
\subfigure[$L=L_2$]{\includegraphics[width=1.7in]{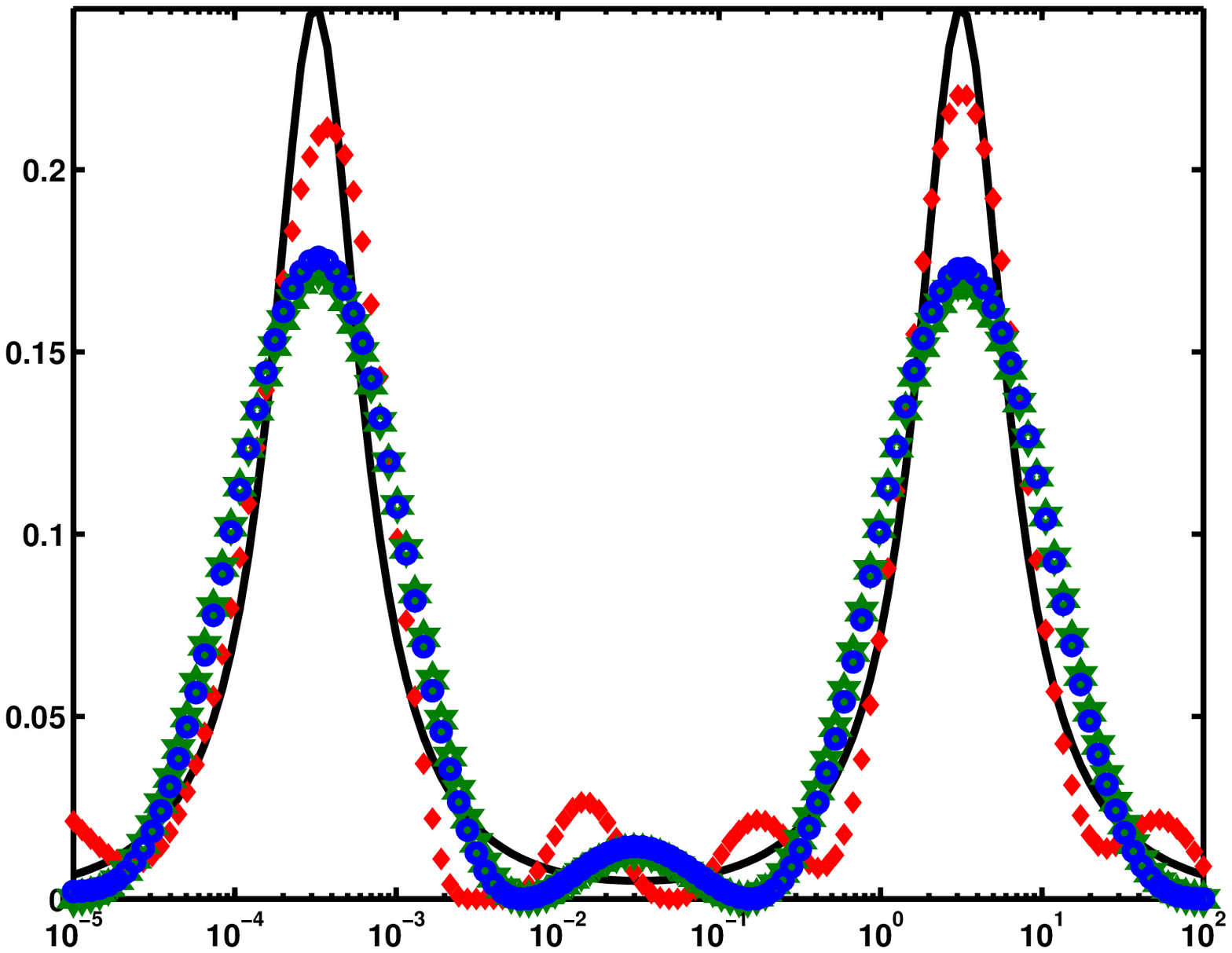}}
\caption{NNLS solutions of RQ-A matrix $A_4$. Noise level $1\%$.}
\label{hnfig-lambdachoiceRQ1A4HN}
\end{figure}

 \begin{figure}[!ht]
\centering
\subfigure[$L=I$]{\includegraphics[width=1.7in]{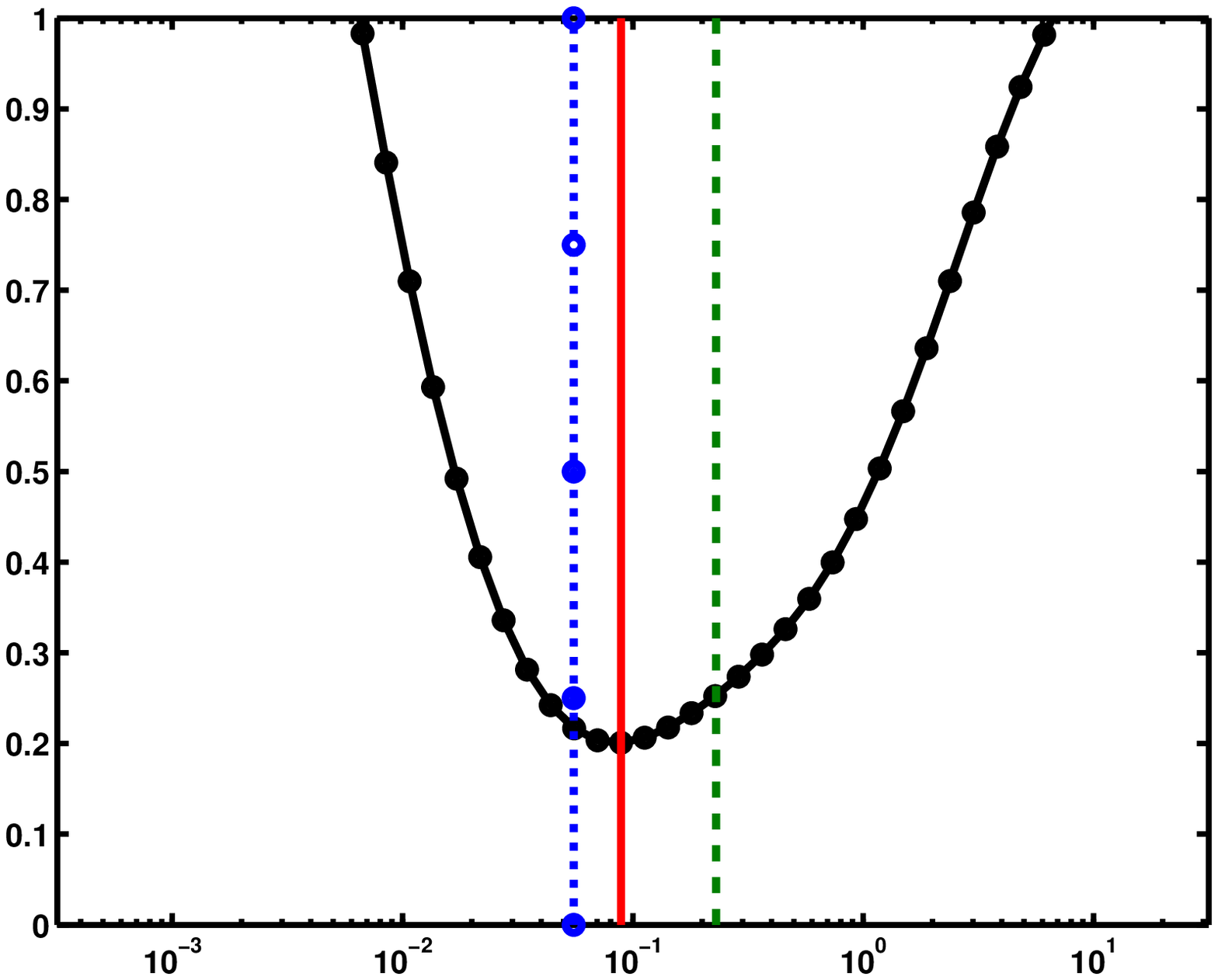}}
\subfigure[$L=L_1$]{\includegraphics[width=1.7in]{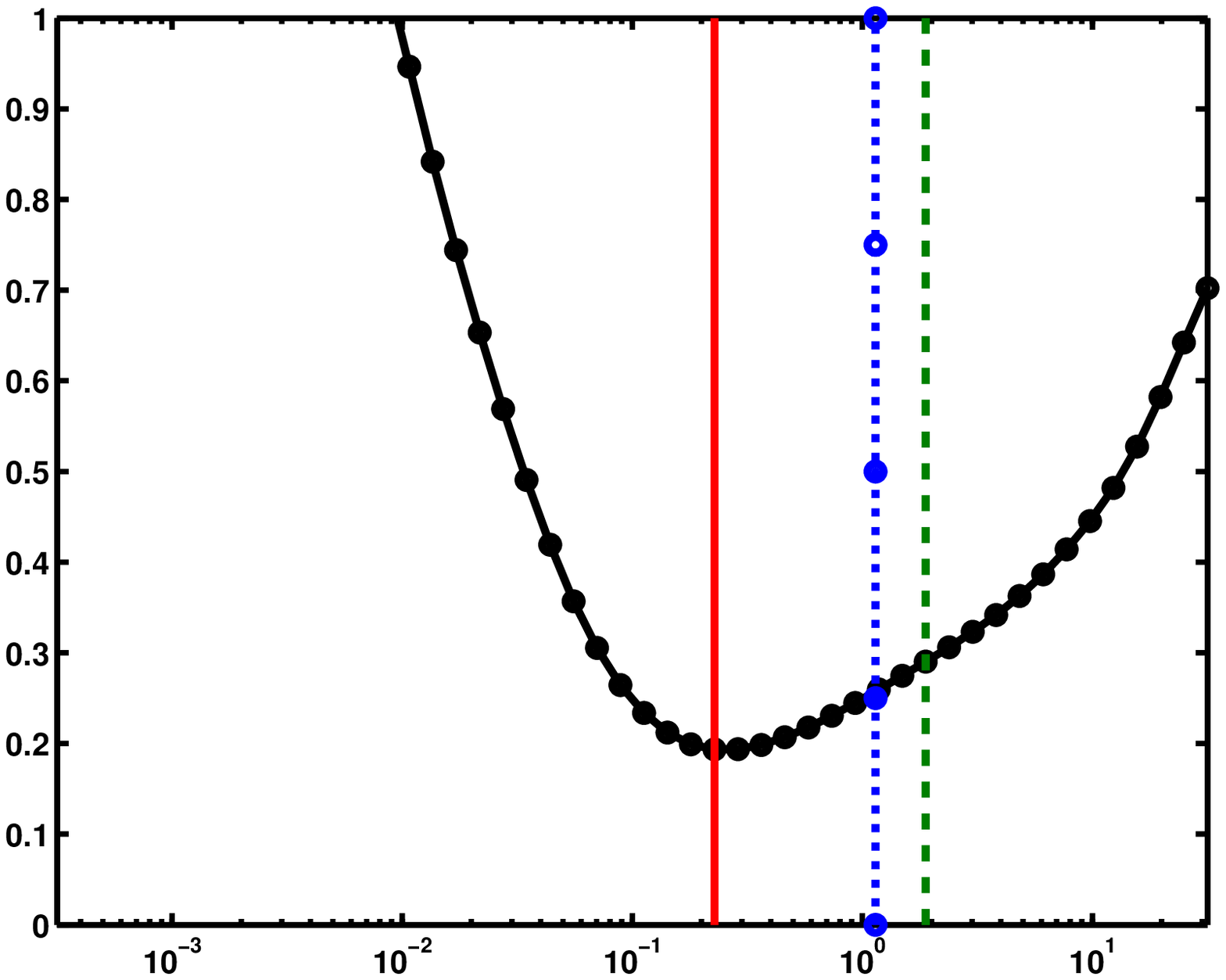}}
\subfigure[$L=L_2$]{\includegraphics[width=1.7in]{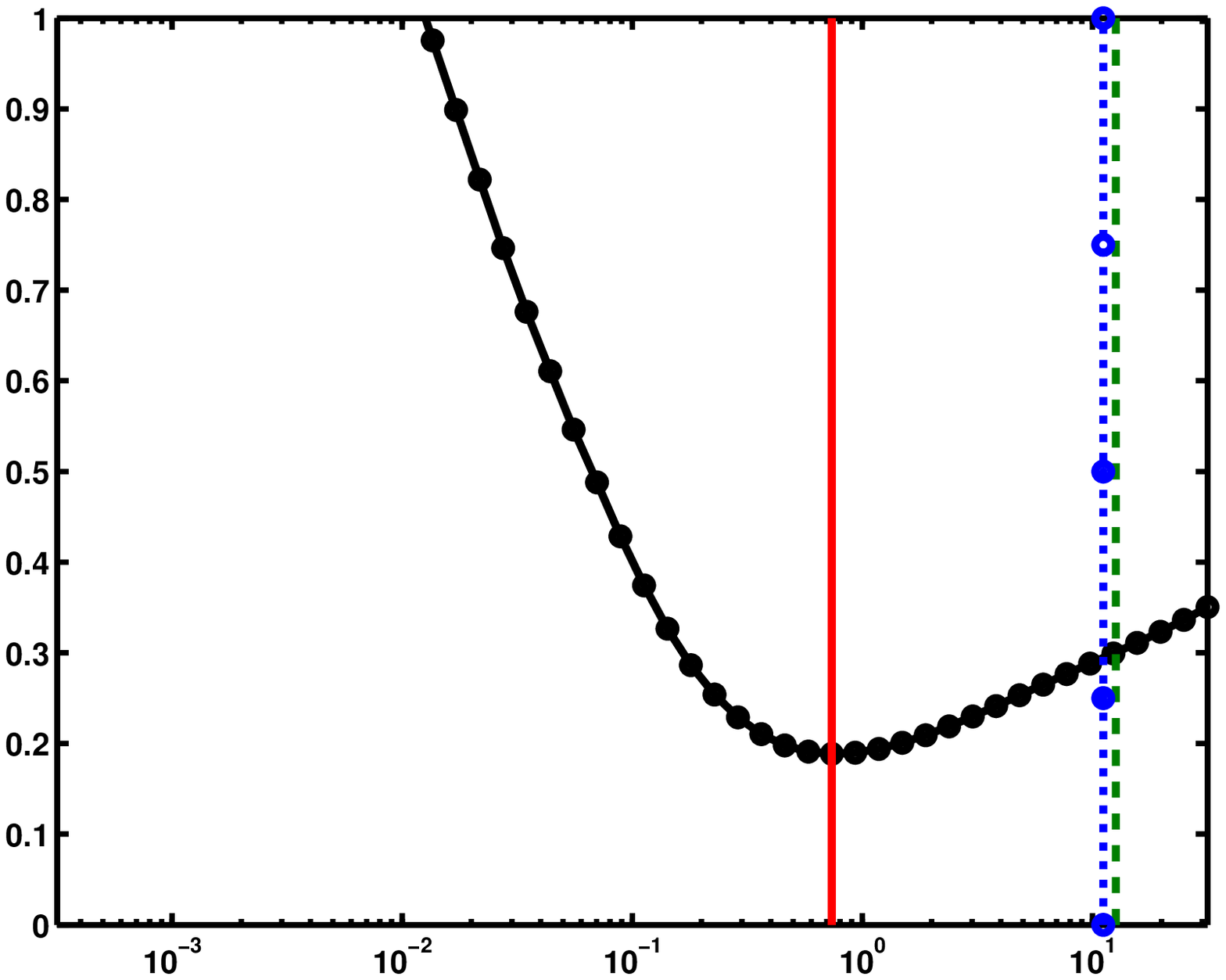}}
\subfigure[$L=I$]{\includegraphics[width=1.7in]{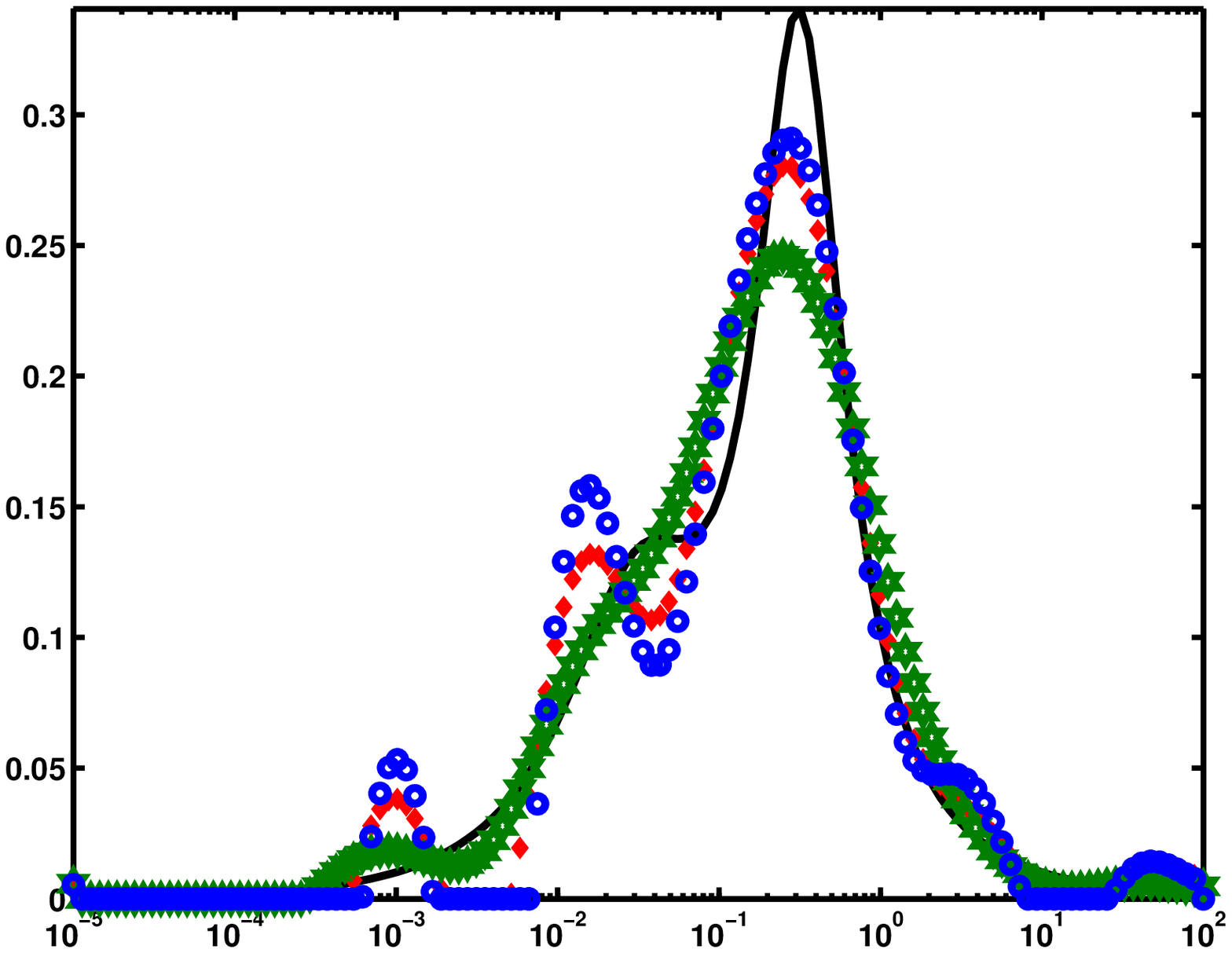}}
\subfigure[$L=L_1$]{\includegraphics[width=1.7in]{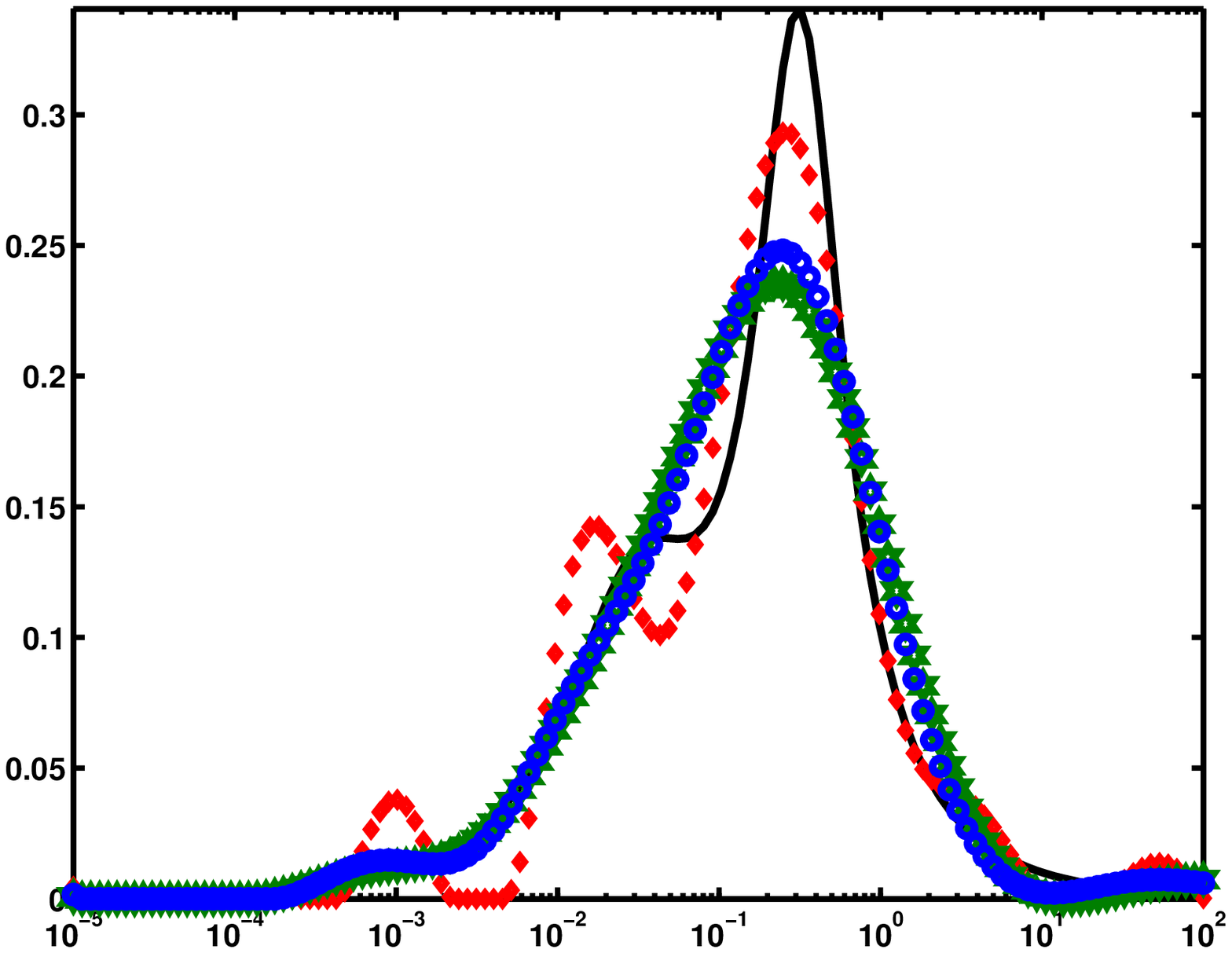}}
\subfigure[$L=L_2$]{\includegraphics[width=1.7in]{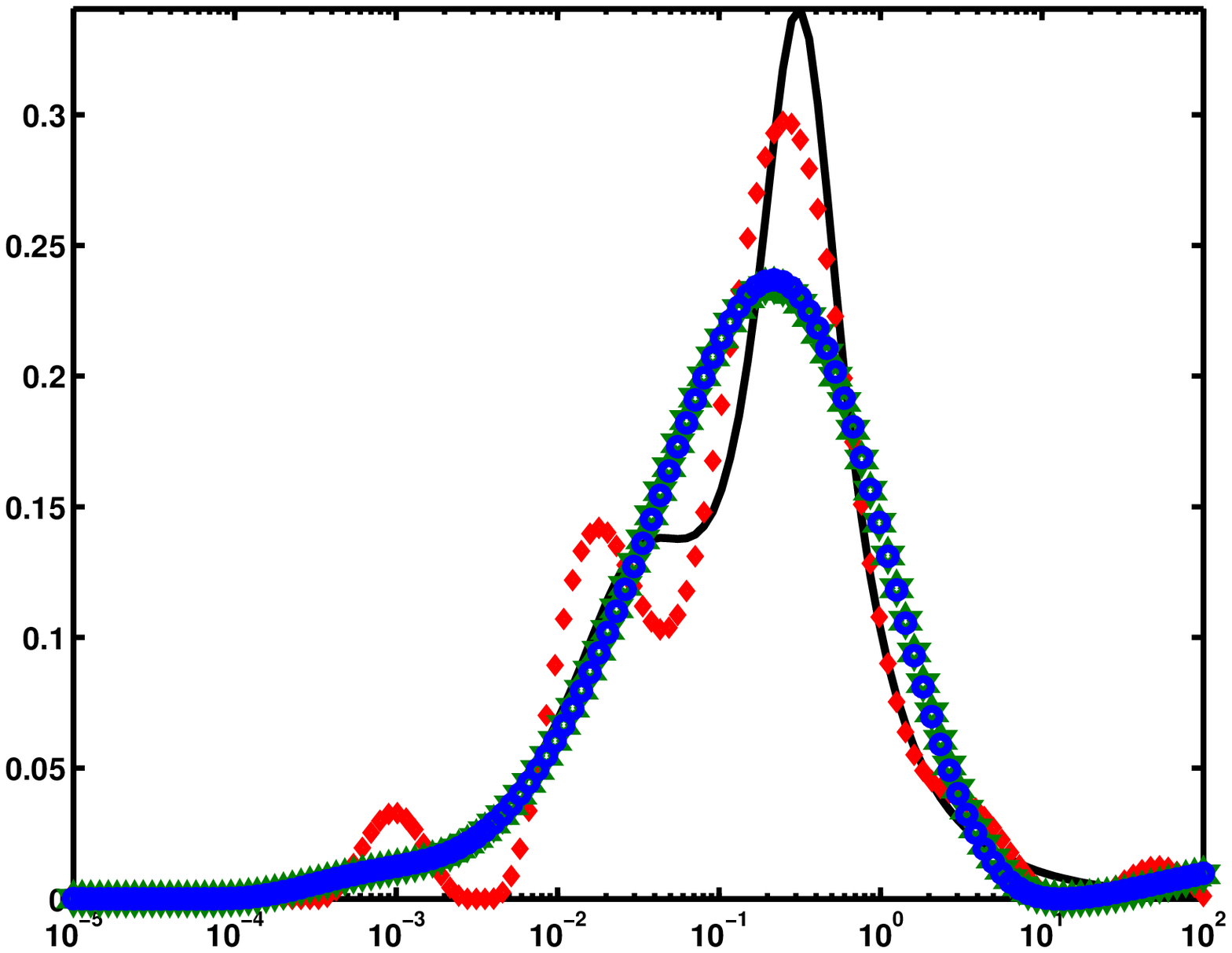}}
\caption{NNLS solutions of RQ-B matrix $A_4$. Noise level $1\%$.}
\label{hnfig-lambdachoiceRQ5A4HN}
\end{figure}
 \begin{figure}[!ht]
\centering
\subfigure[$L=I$]{\includegraphics[width=1.7in]{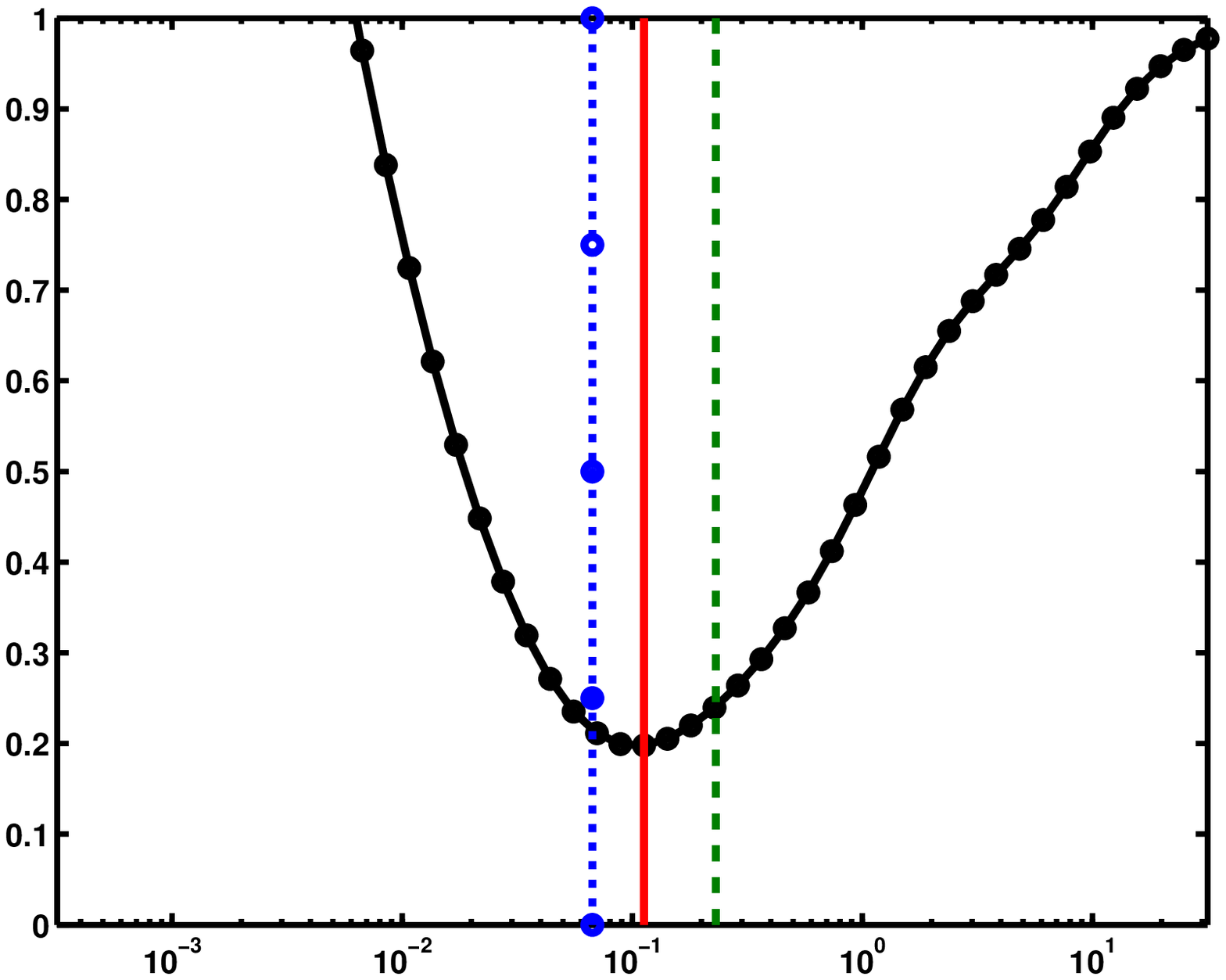}}
\subfigure[$L=L_1$]{\includegraphics[width=1.7in]{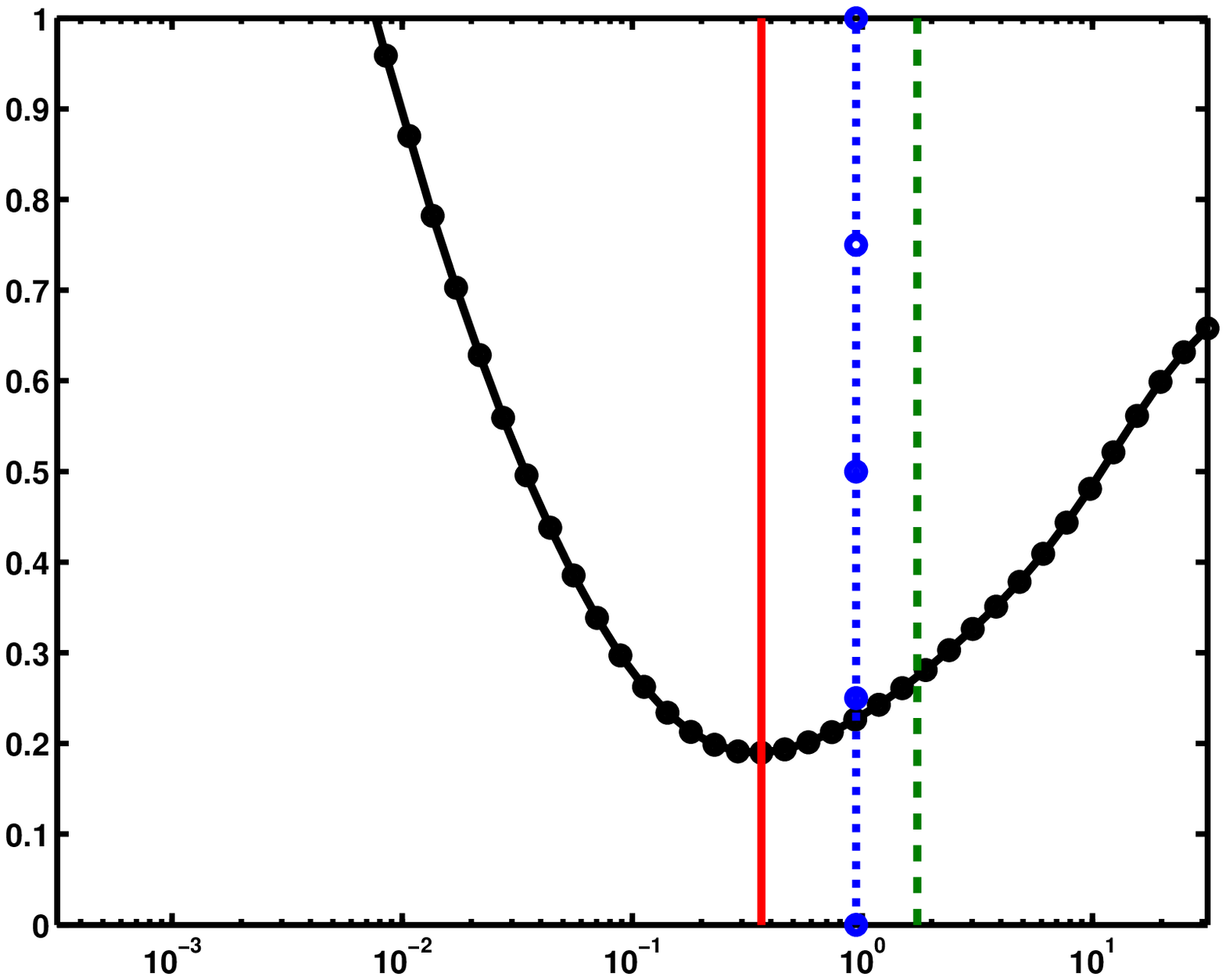}}
\subfigure[$L=L_2$]{\includegraphics[width=1.7in]{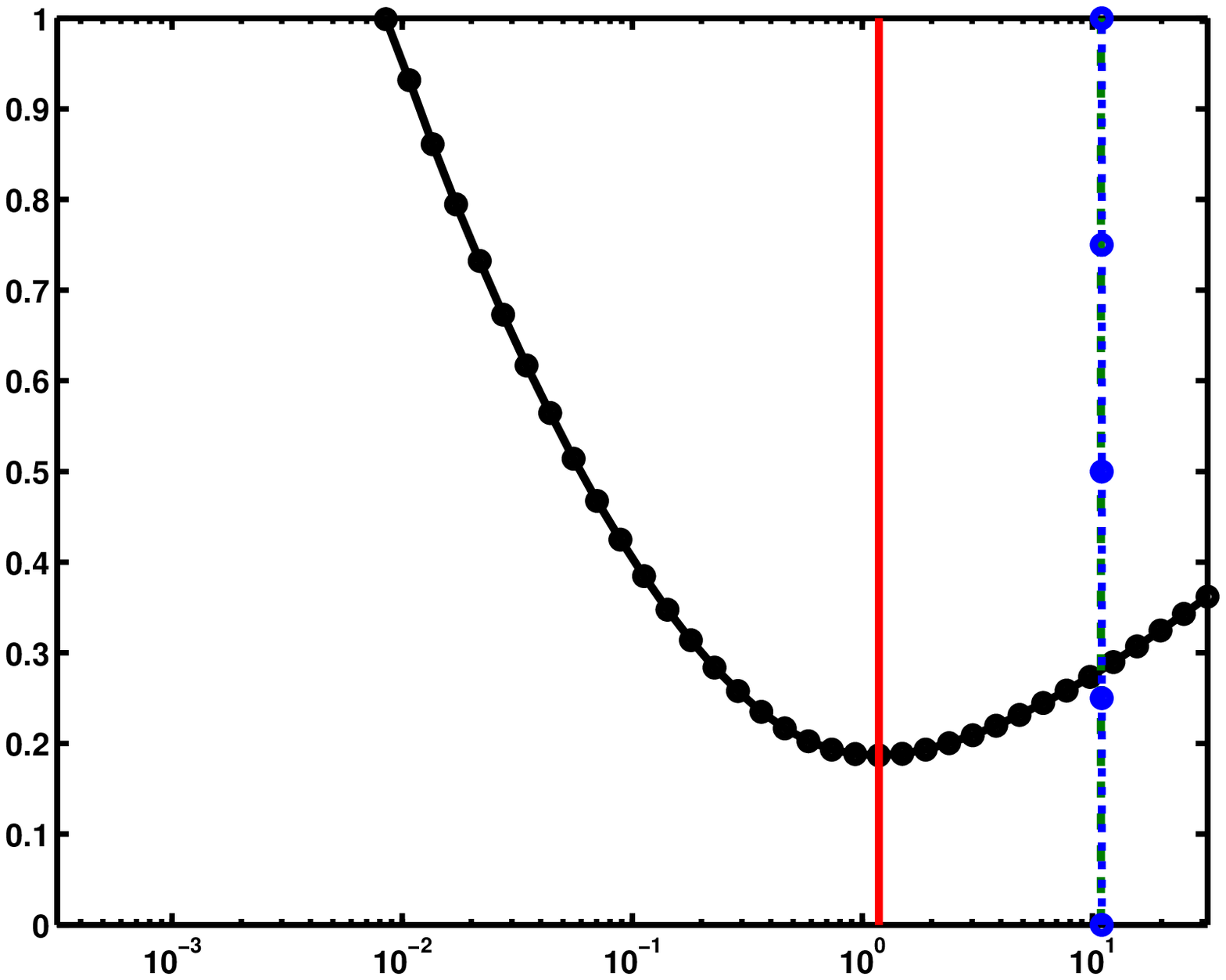}}
\subfigure[$L=I$]{\includegraphics[width=1.7in]{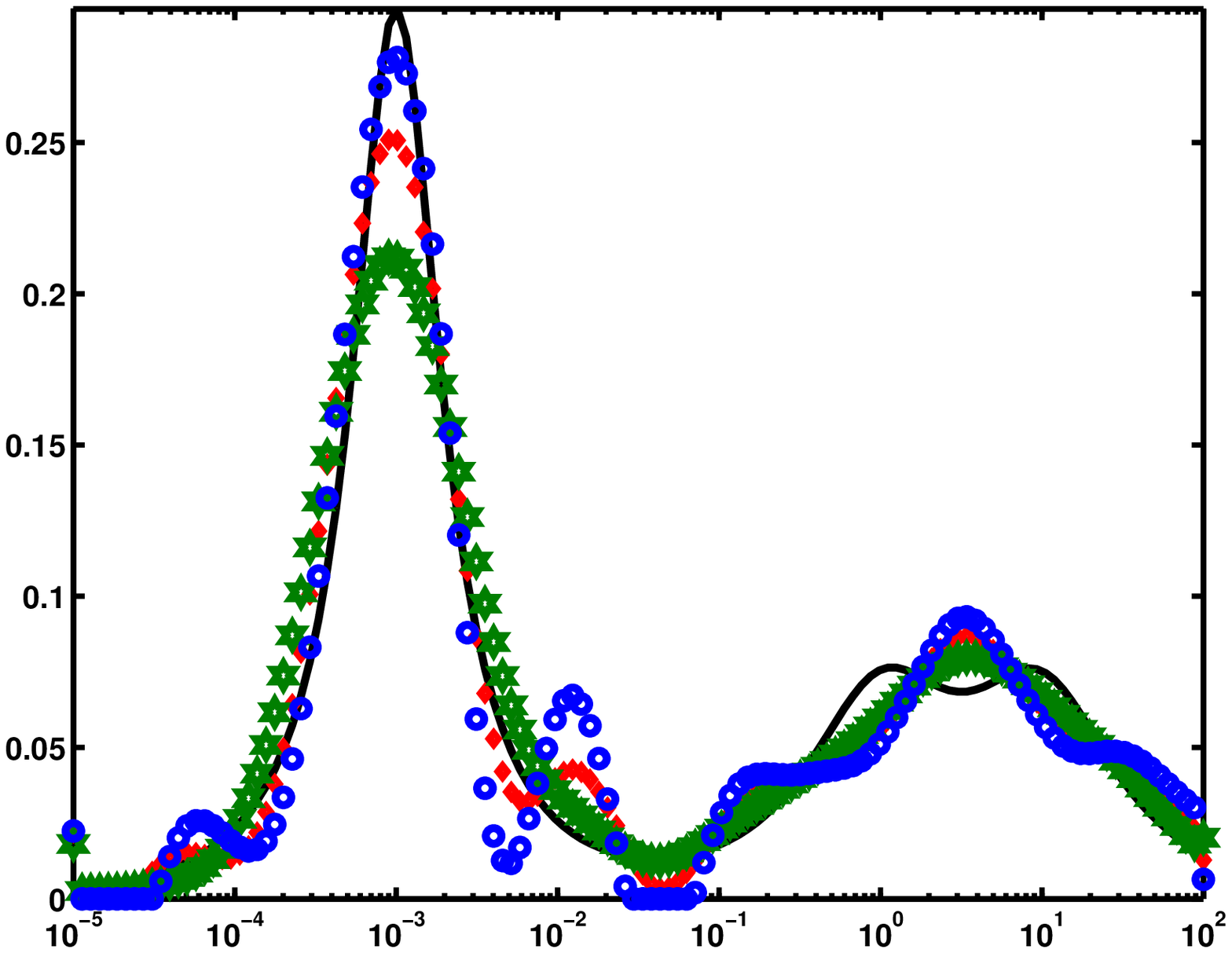}}
\subfigure[$L=L_1$]{\includegraphics[width=1.7in]{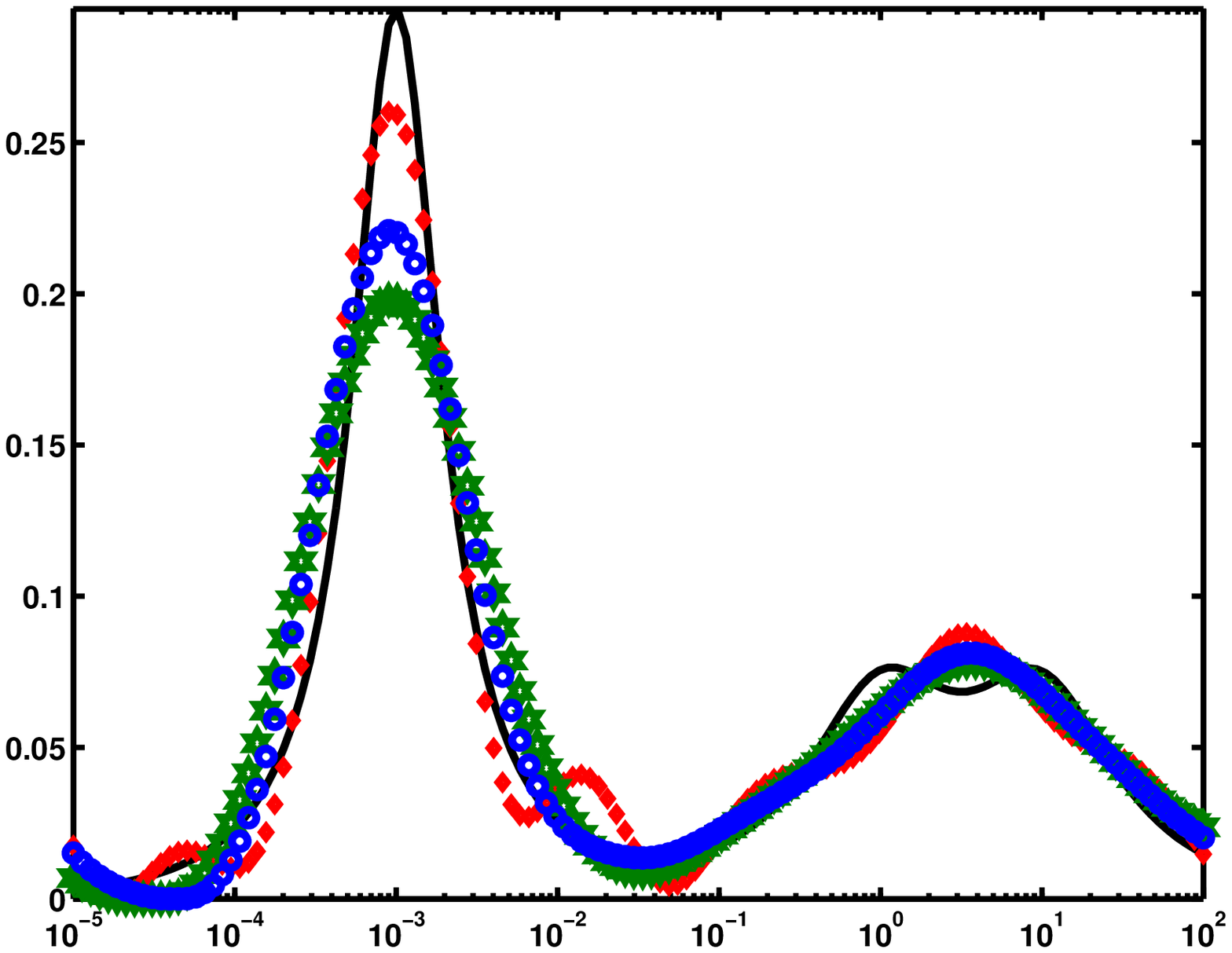}}
\subfigure[$L=L_2$]{\includegraphics[width=1.7in]{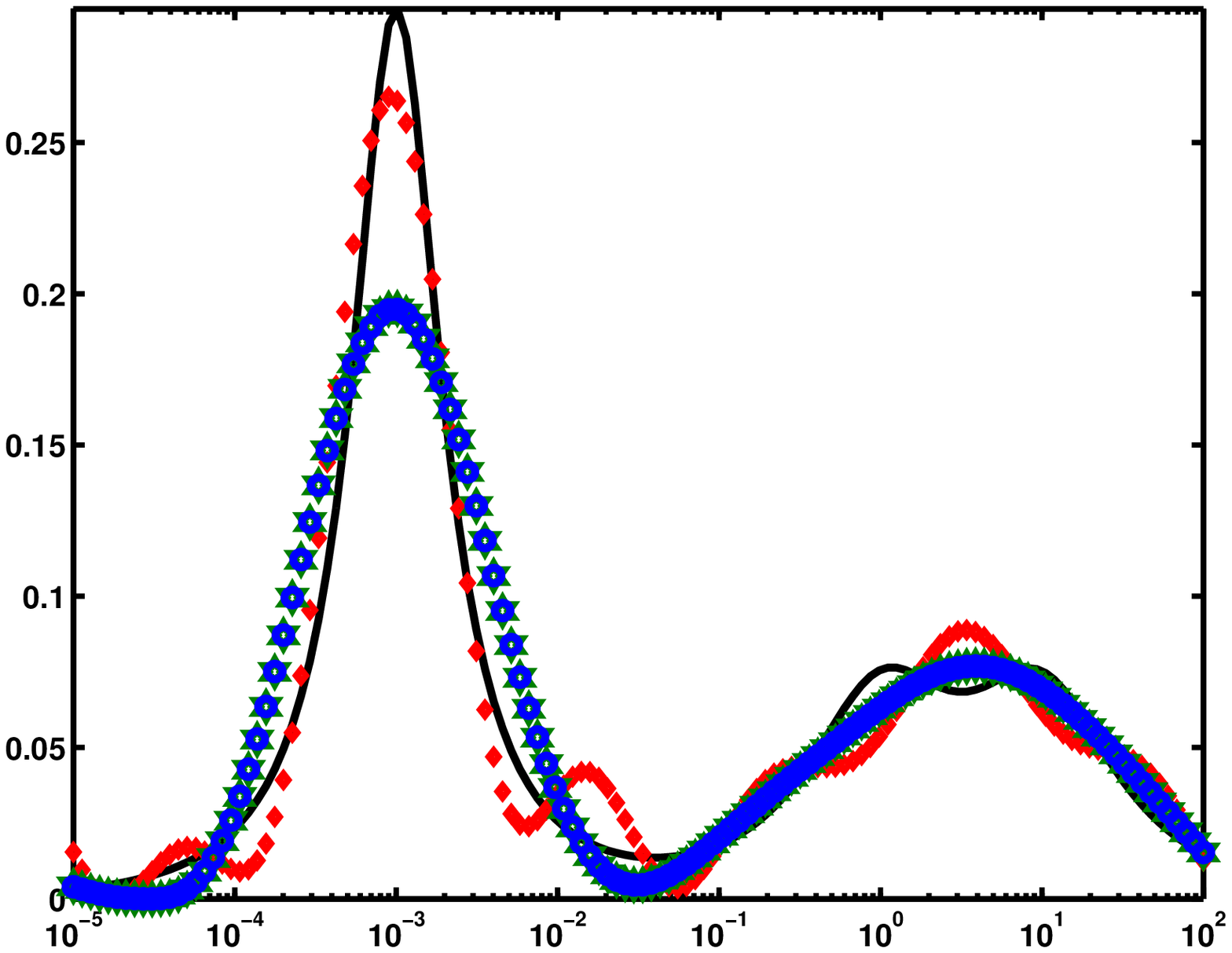}}
\caption{NNLS solutions of RQ-C matrix $A_4$. Noise level $1\%$.}
\label{hnfig-lambdachoiceRQ6A4HN}
\end{figure}
 \begin{figure}[!ht]
\centering
\subfigure[$L=I$]{\includegraphics[width=1.7in]{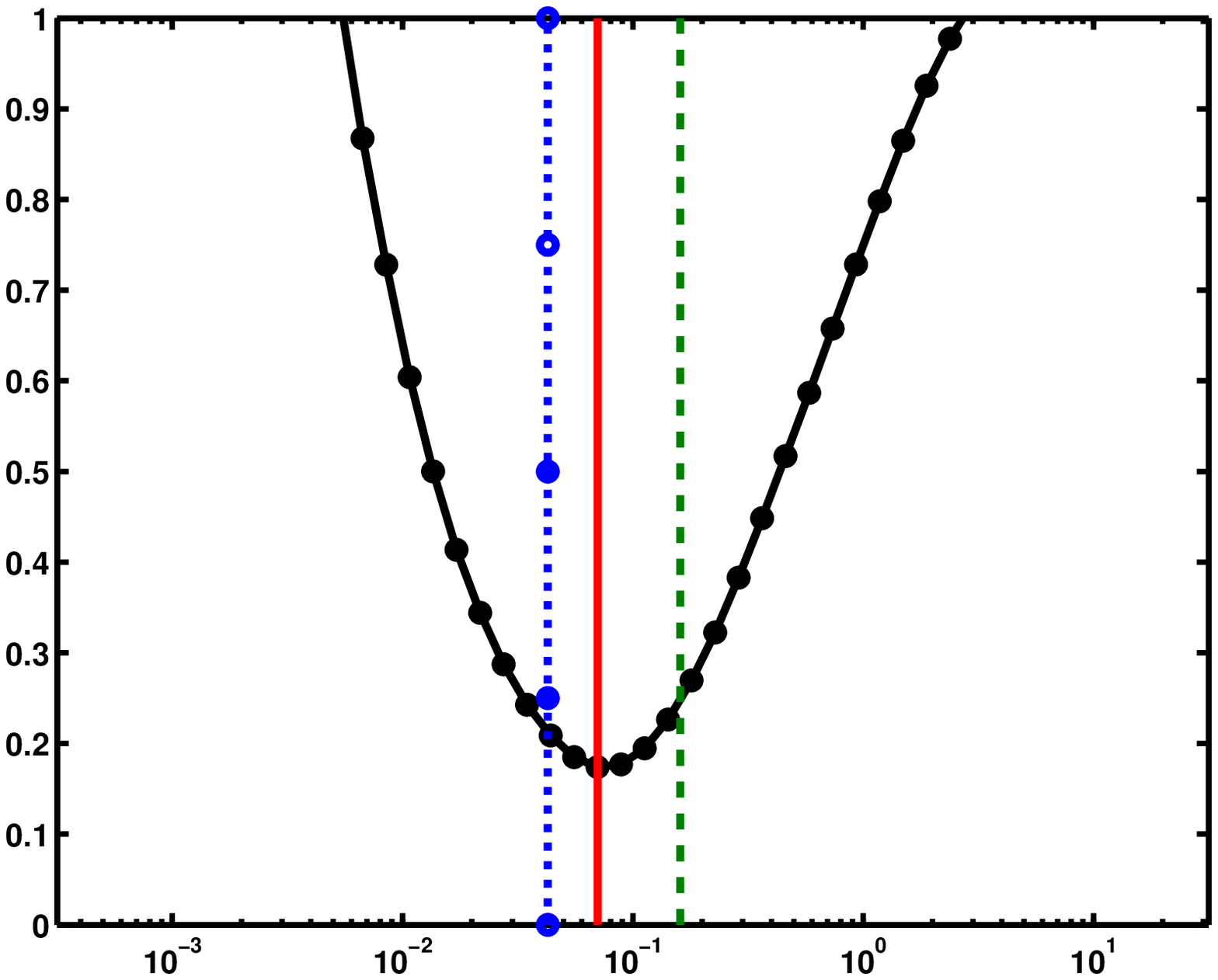}}
\subfigure[$L=L_1$]{\includegraphics[width=1.7in]{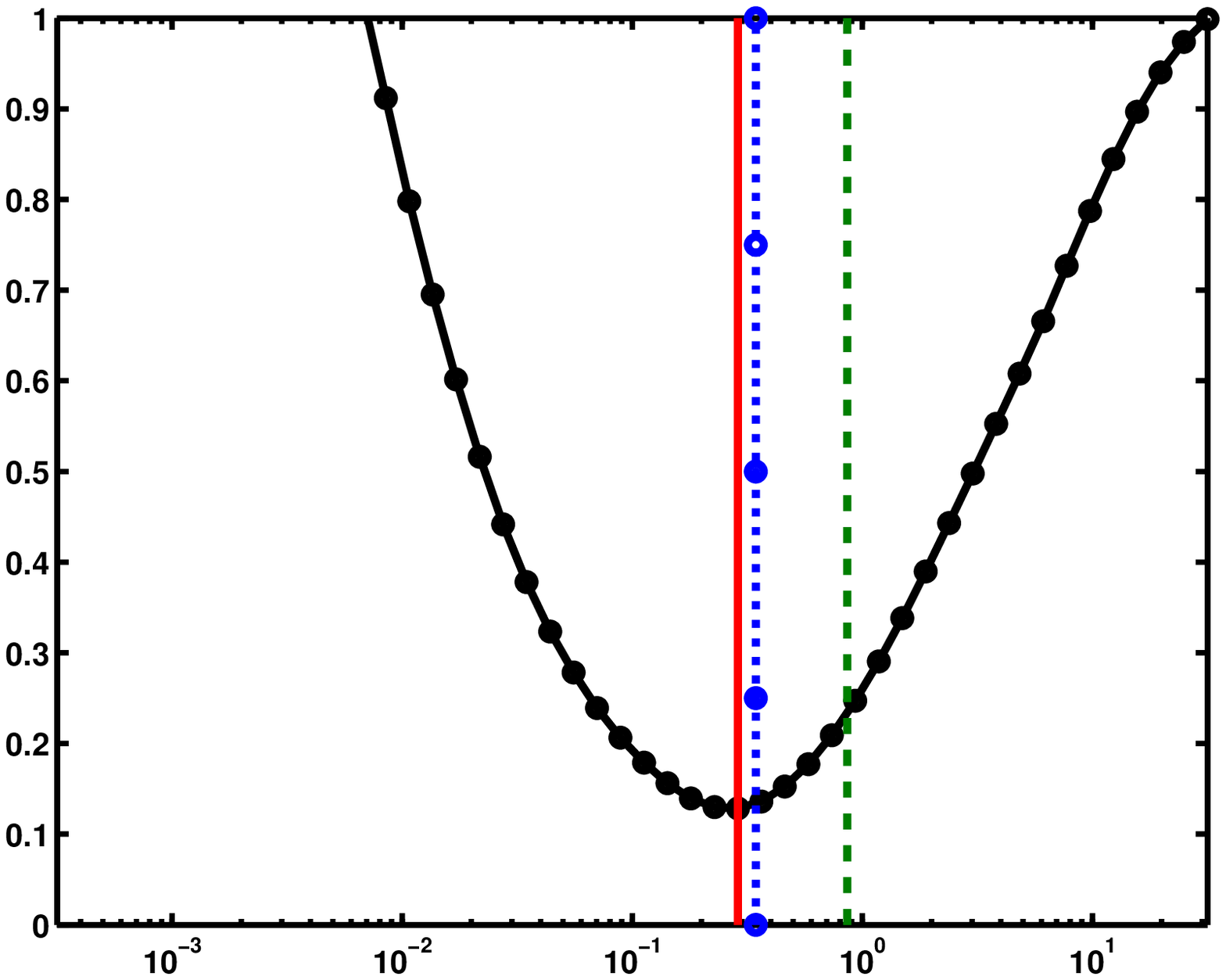}}
\subfigure[$L=L_2$]{\includegraphics[width=1.7in]{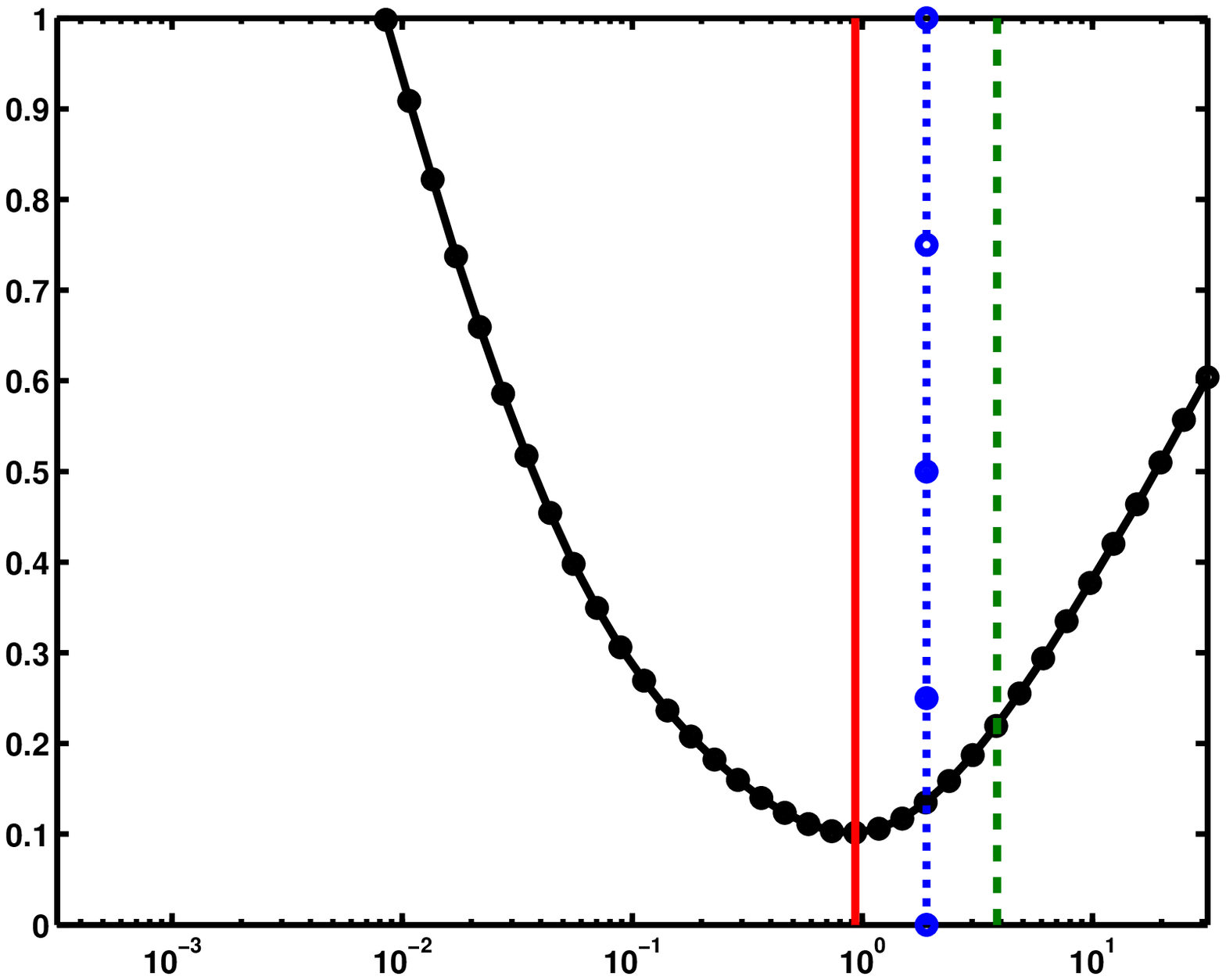}}
\subfigure[$L=I$]{\includegraphics[width=1.7in]{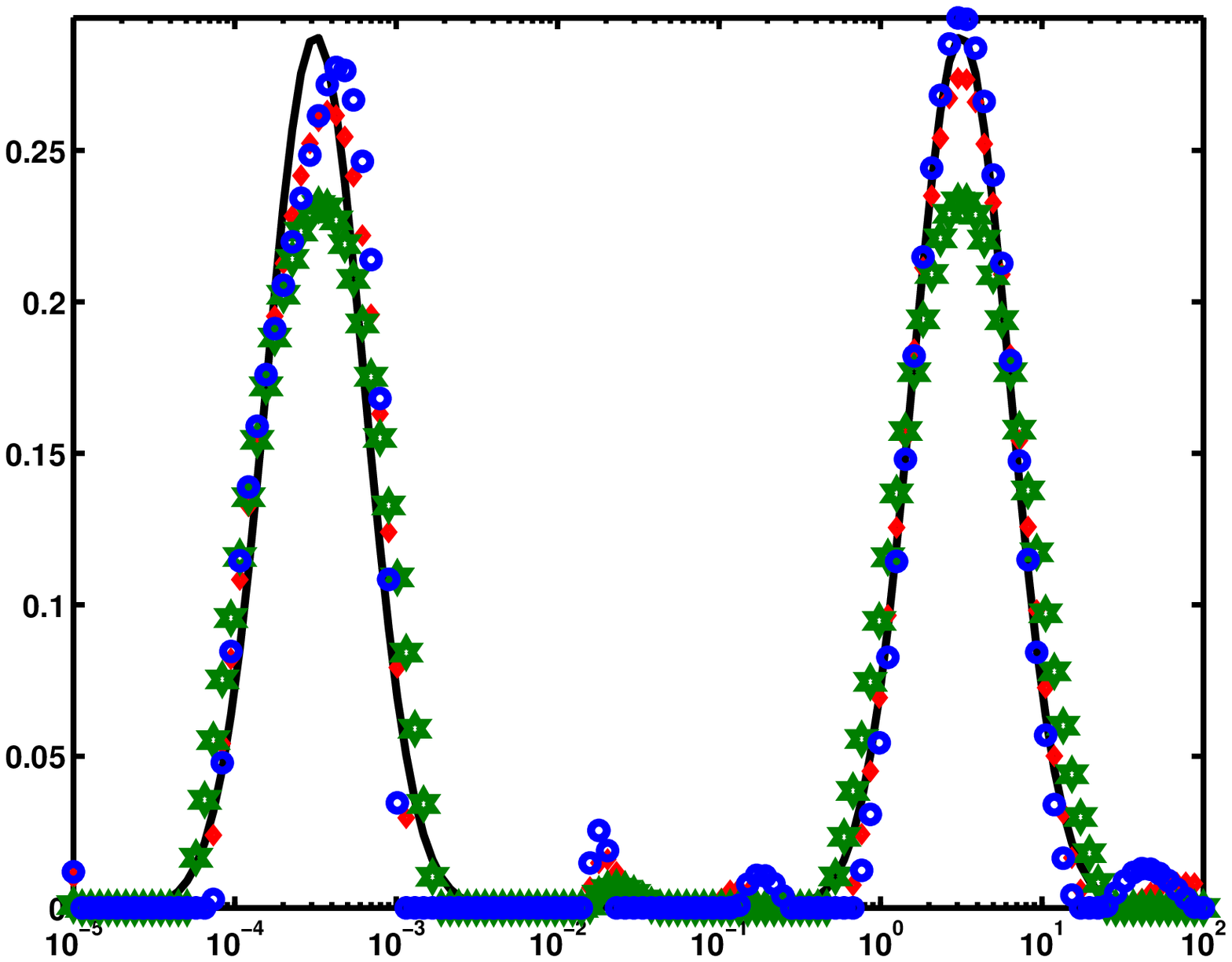}}
\subfigure[$L=L_1$]{\includegraphics[width=1.7in]{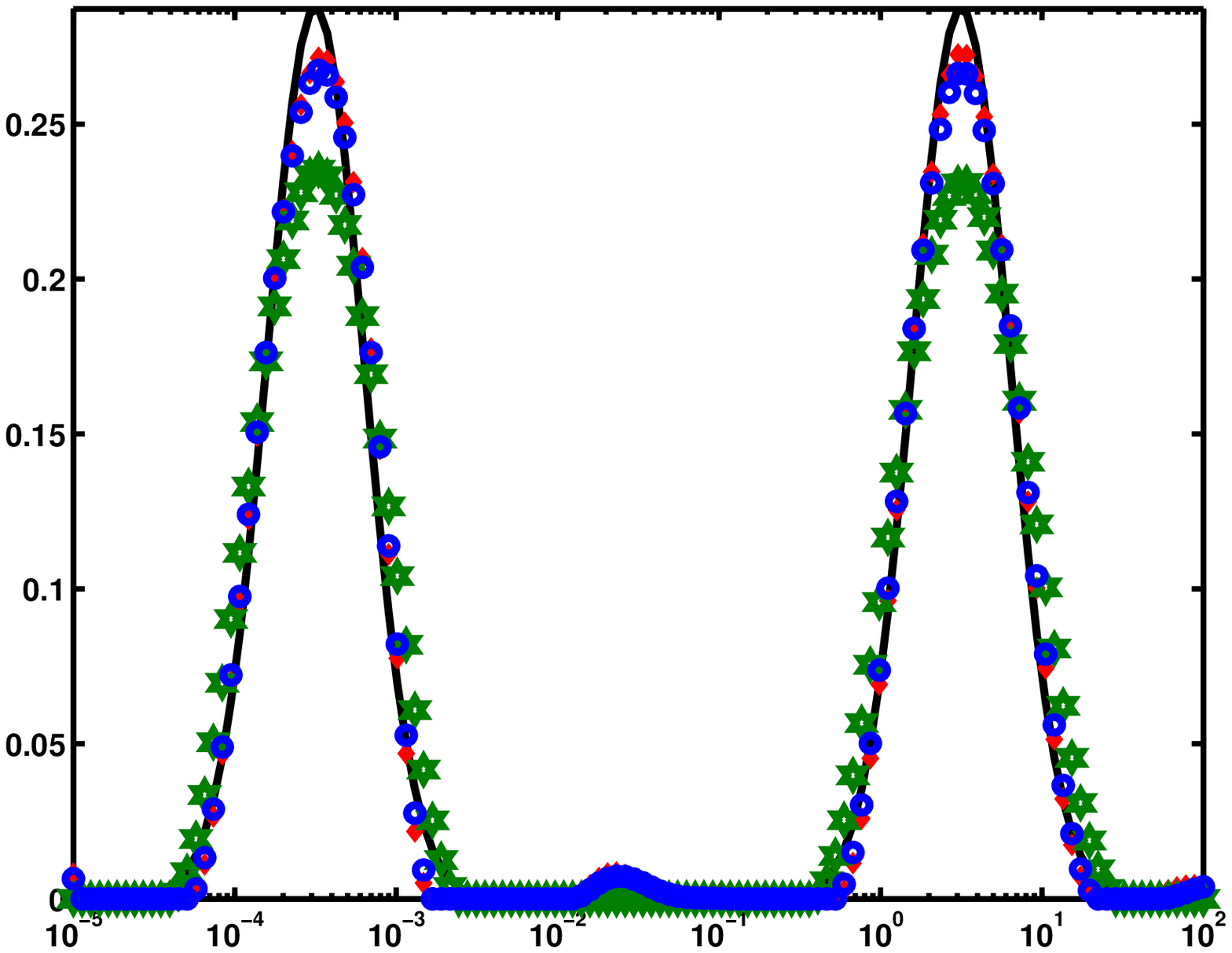}}
\subfigure[$L=L_2$]{\includegraphics[width=1.7in]{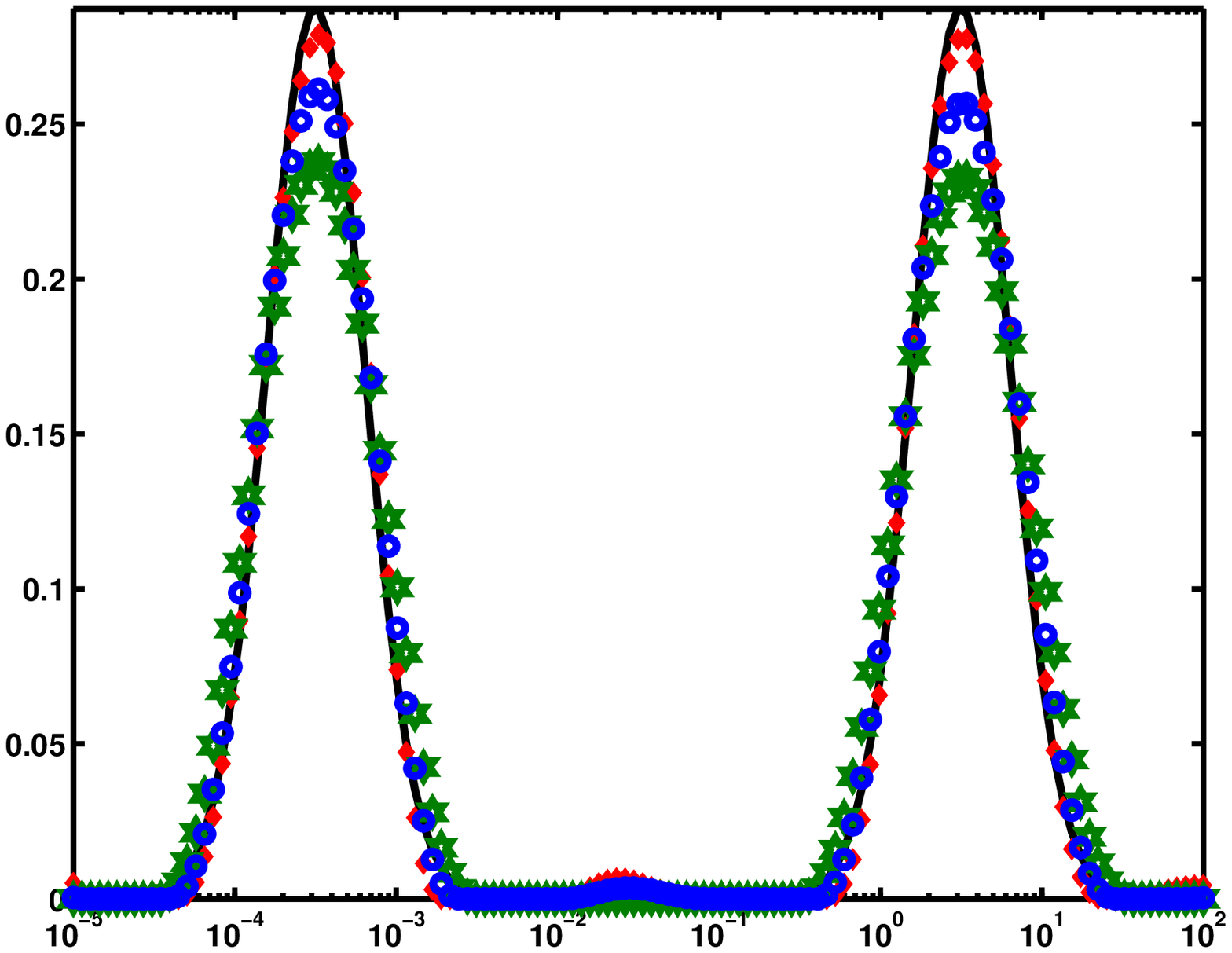}}
\caption{NNLS solutions of LN-A matrix $A_4$. Noise level $1\%$.}
\label{hnfig-lambdachoiceLN2A4HN}
\end{figure}
 \begin{figure}[!ht]
\centering
\subfigure[$L=I$]{\includegraphics[width=1.7in]{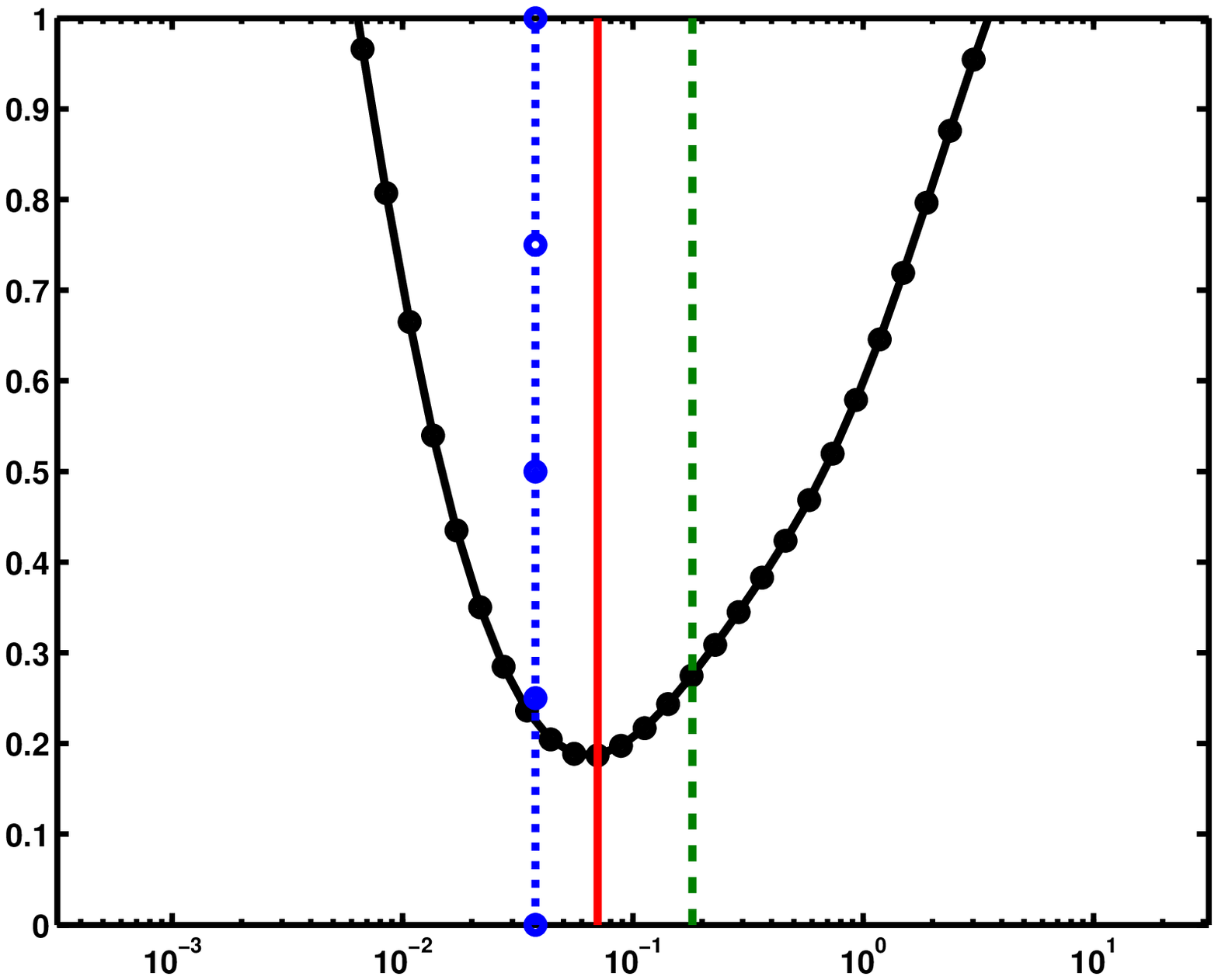}}
\subfigure[$L=L_1$]{\includegraphics[width=1.7in]{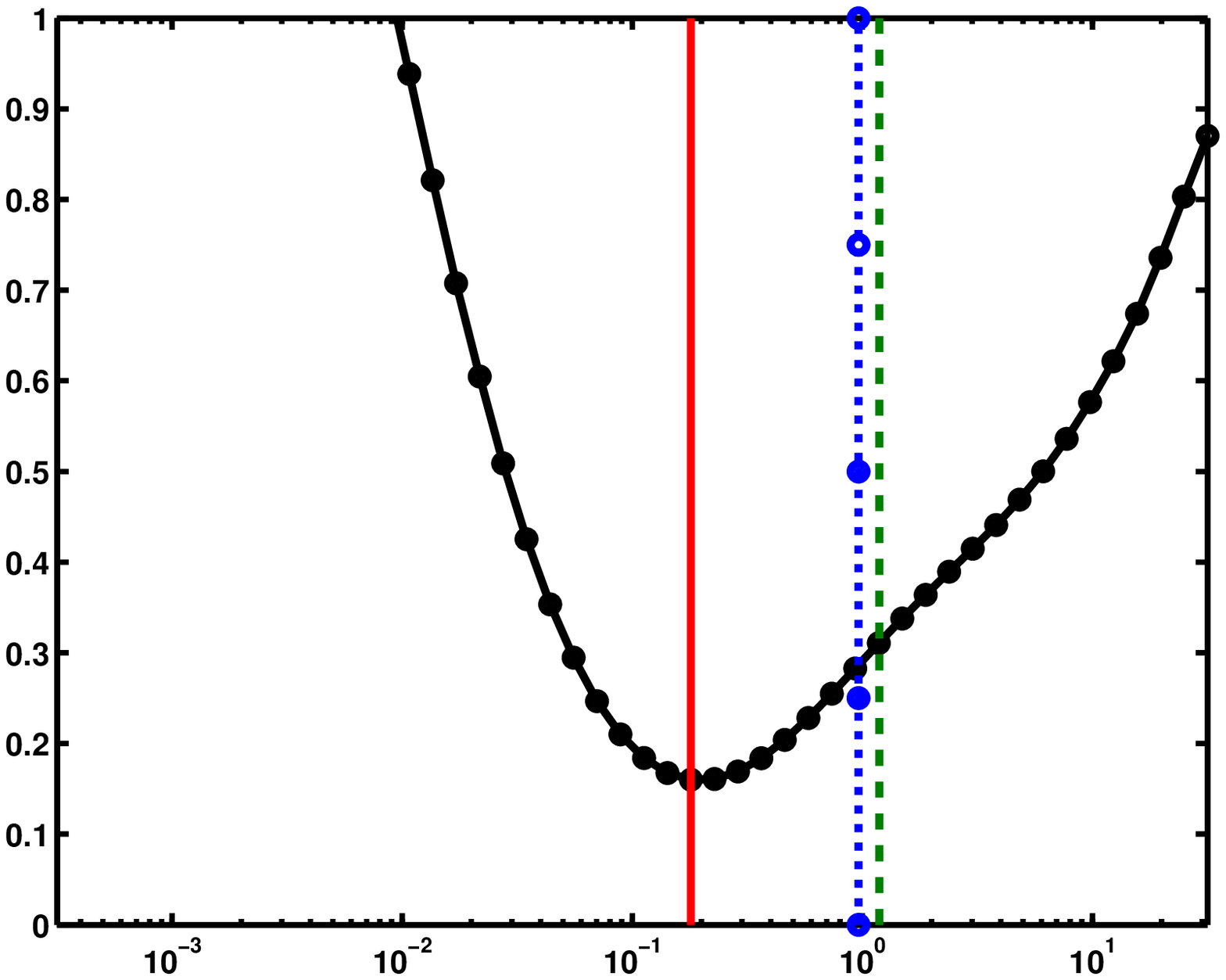}}
\subfigure[$L=L_2$]{\includegraphics[width=1.7in]{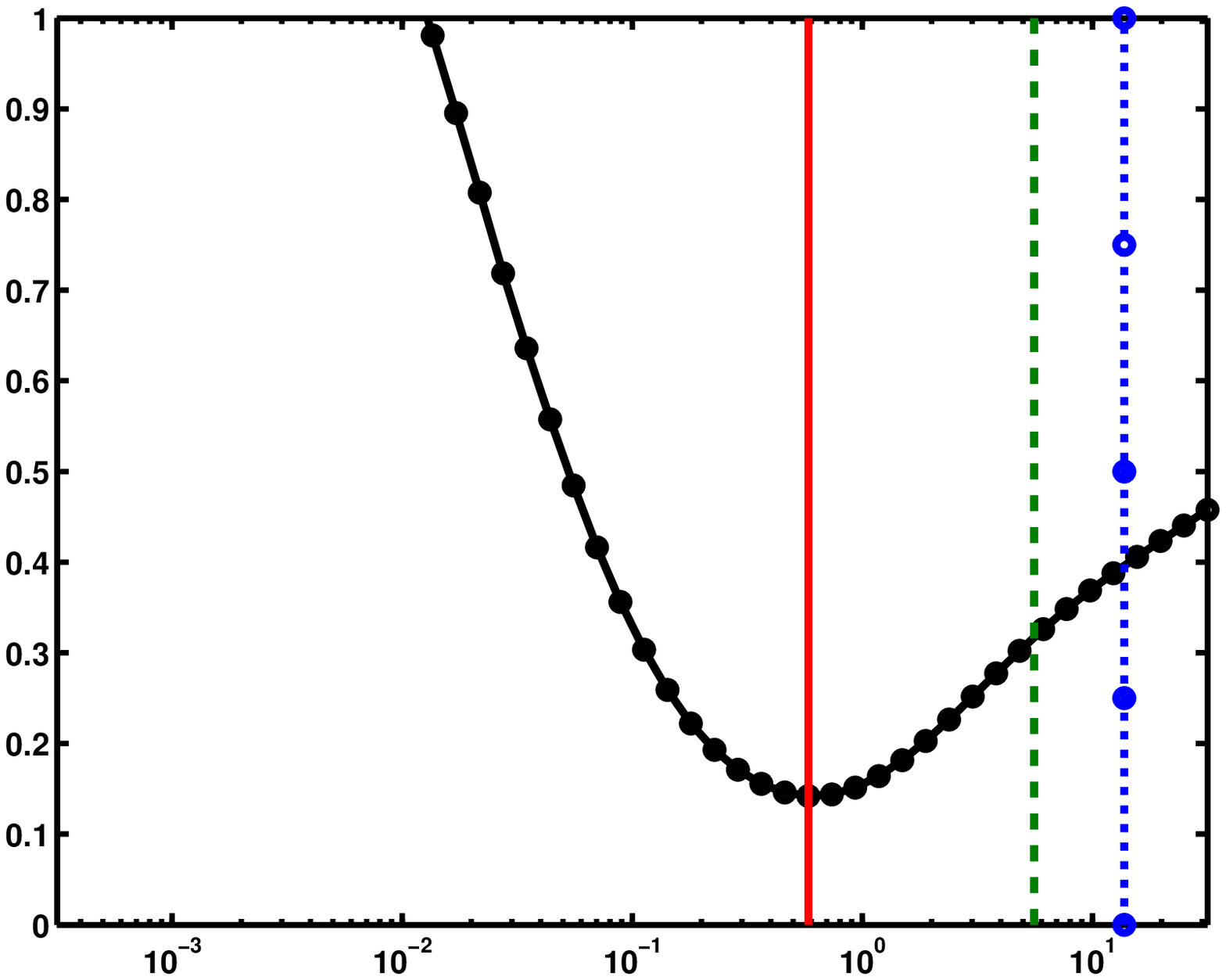}}
\subfigure[$L=I$]{\includegraphics[width=1.7in]{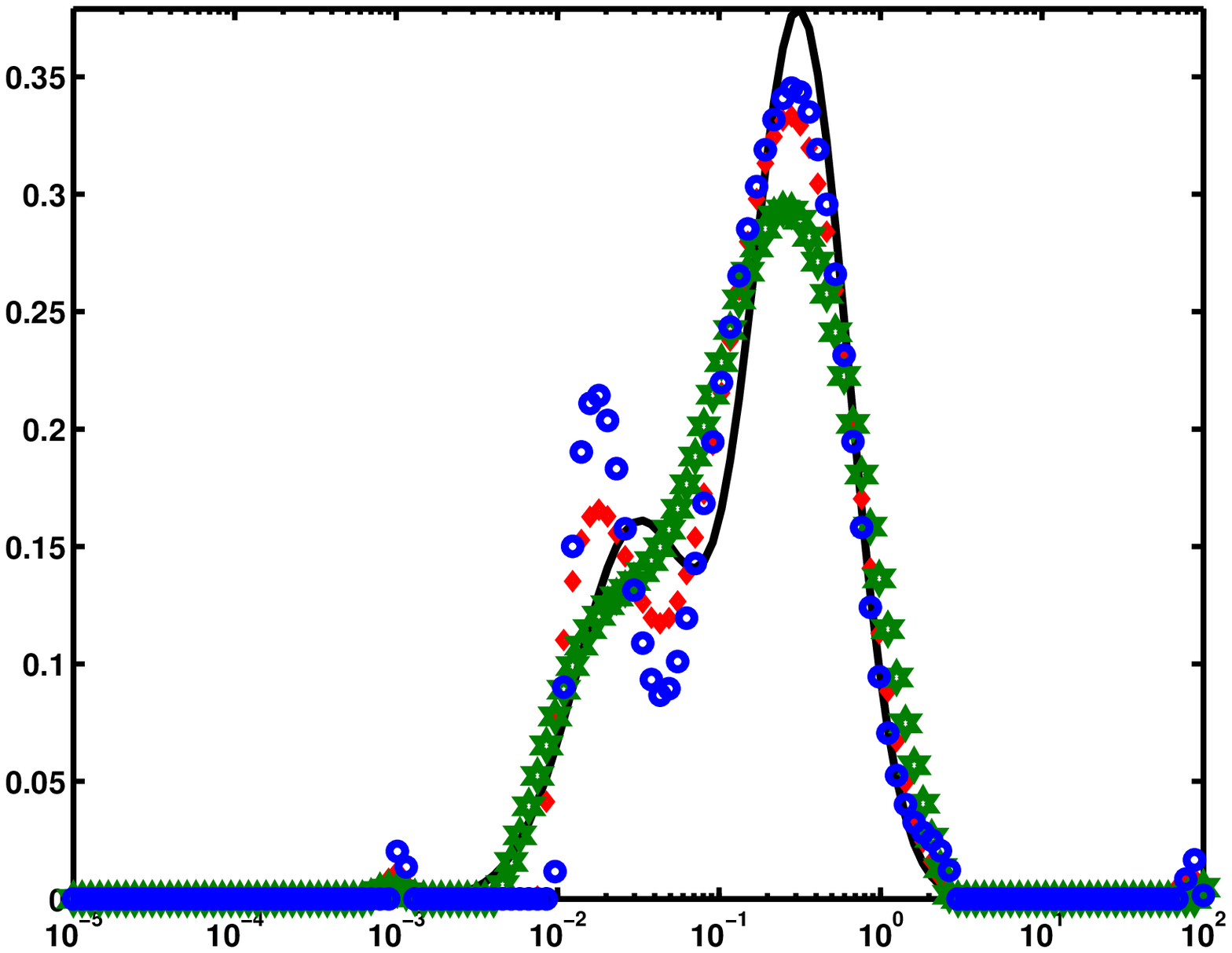}}
\subfigure[$L=L_1$]{\includegraphics[width=1.7in]{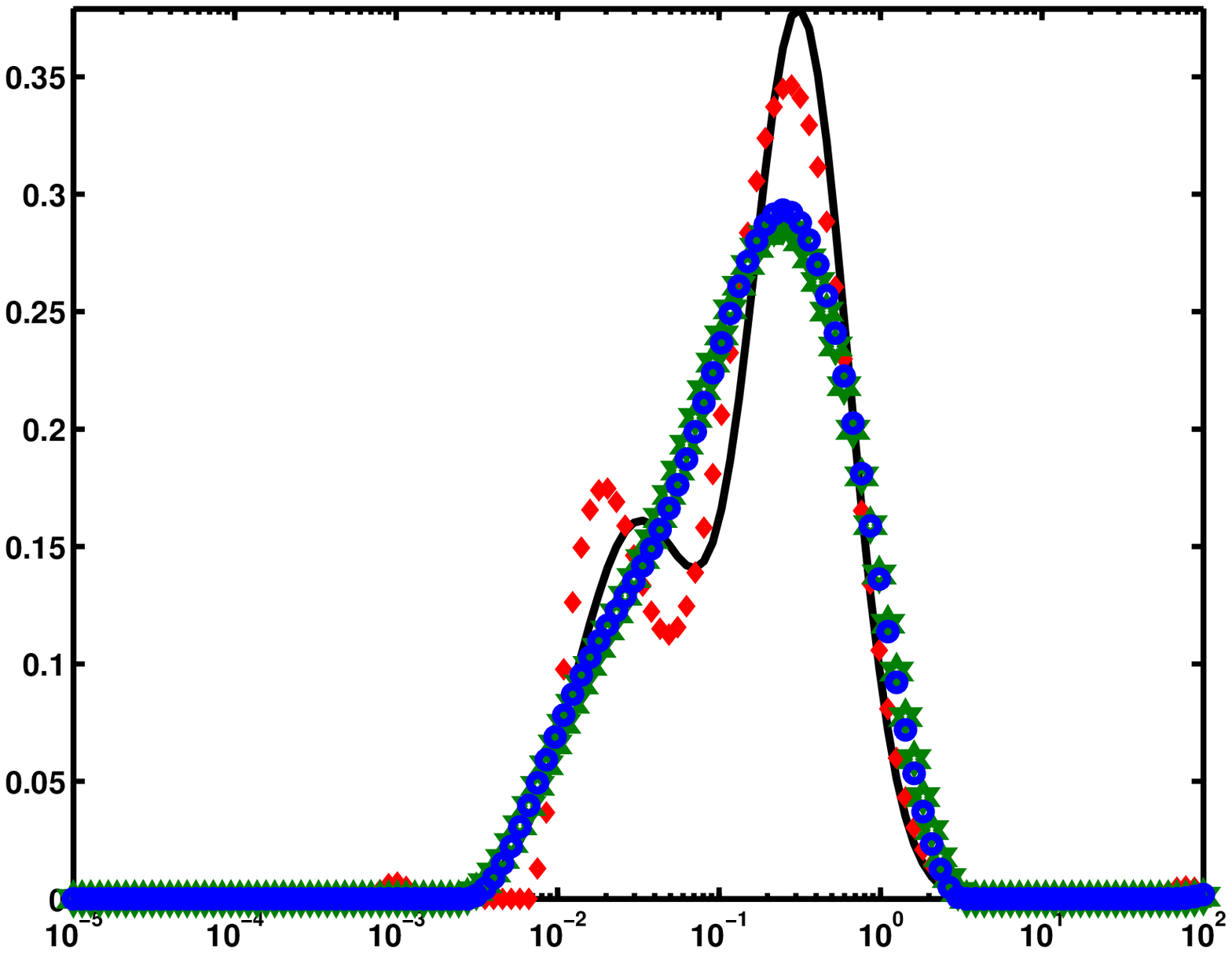}}
\subfigure[$L=L_2$]{\includegraphics[width=1.7in]{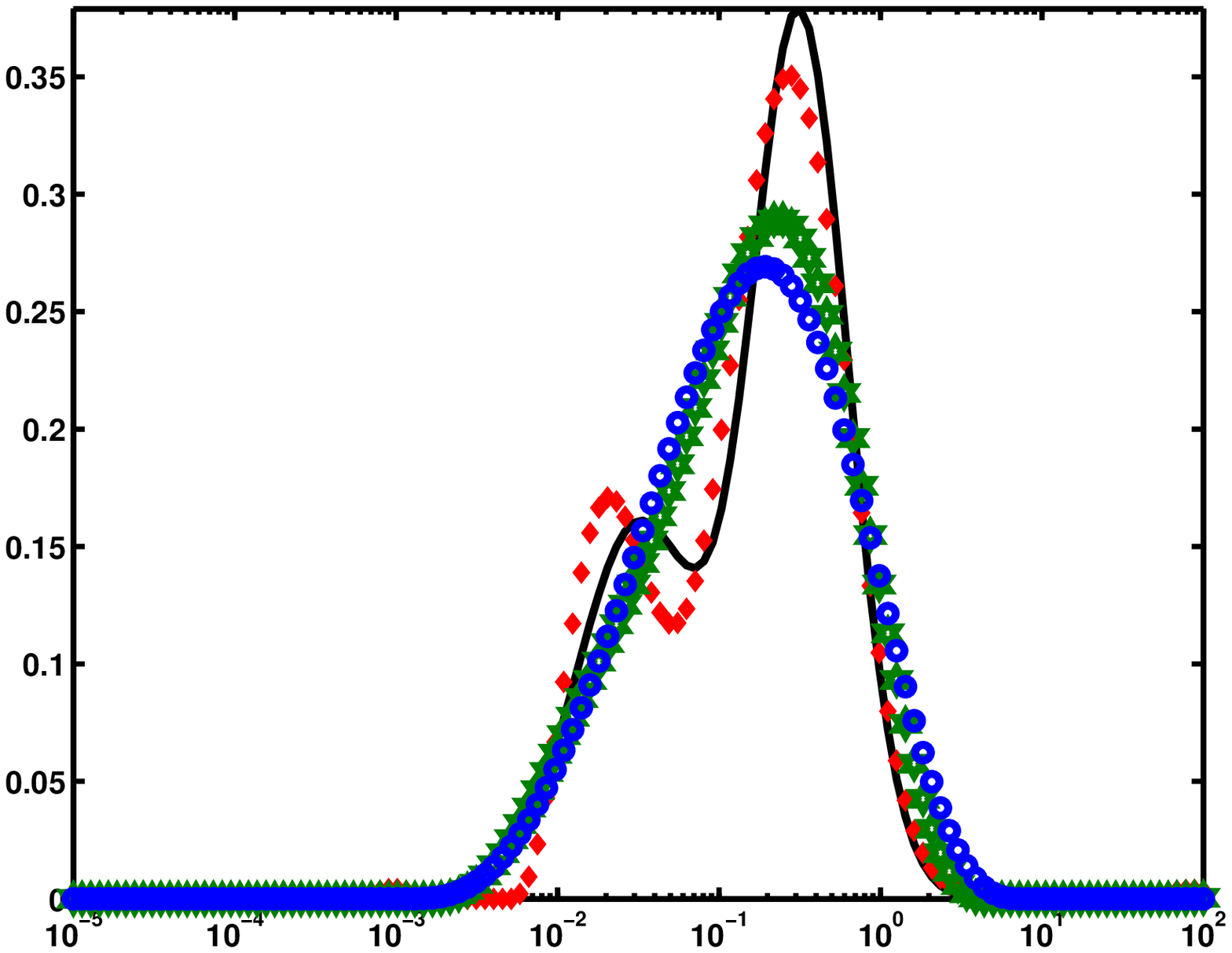}}
\caption{NNLS solutions of LN-B matrix $A_4$. Noise level $1\%$.}
\label{hnfig-lambdachoiceLN5A4HN}
\end{figure}
 \begin{figure}[!ht]
\centering
\subfigure[$L=I$]{\includegraphics[width=1.7in]{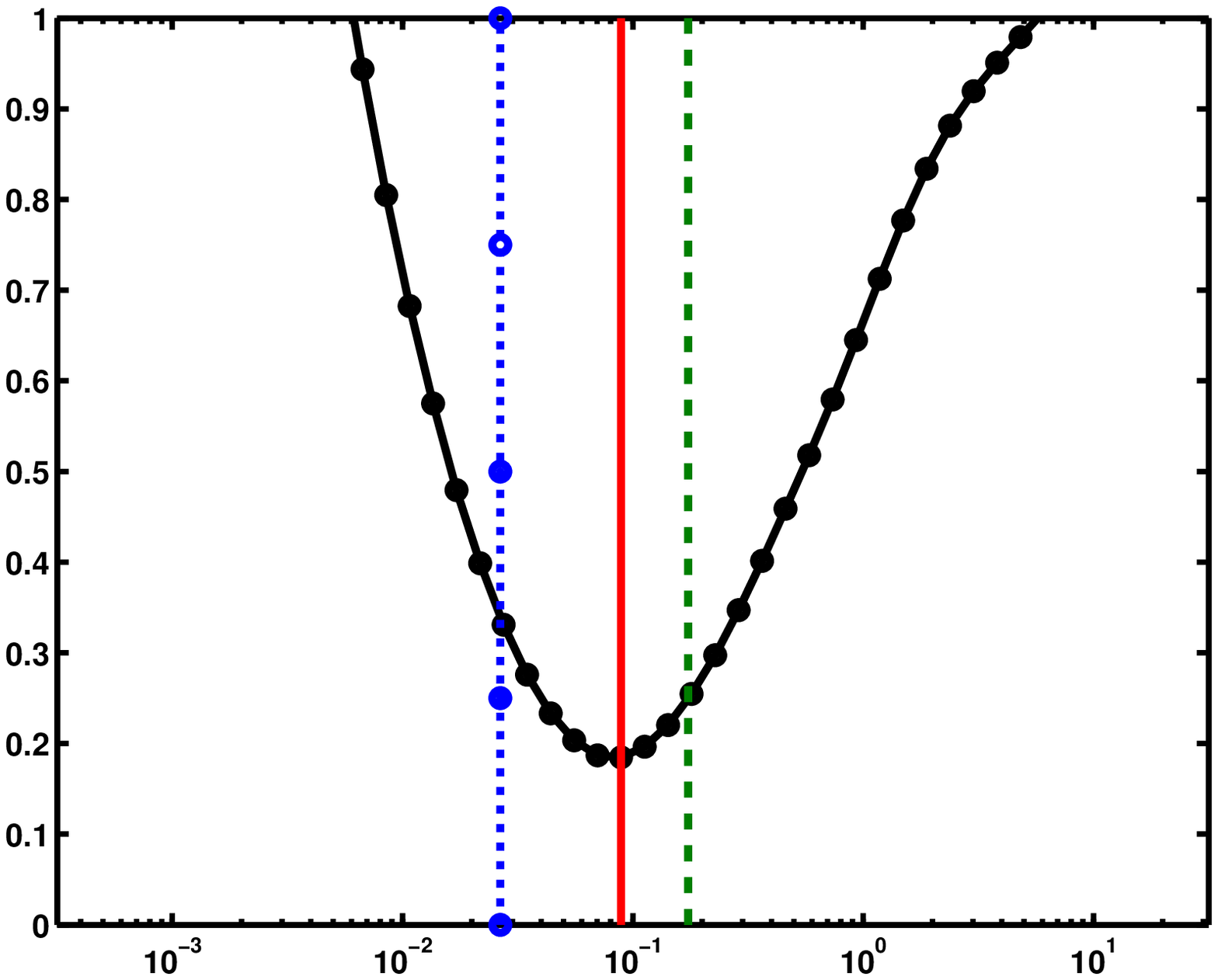}}
\subfigure[$L=L_1$]{\includegraphics[width=1.7in]{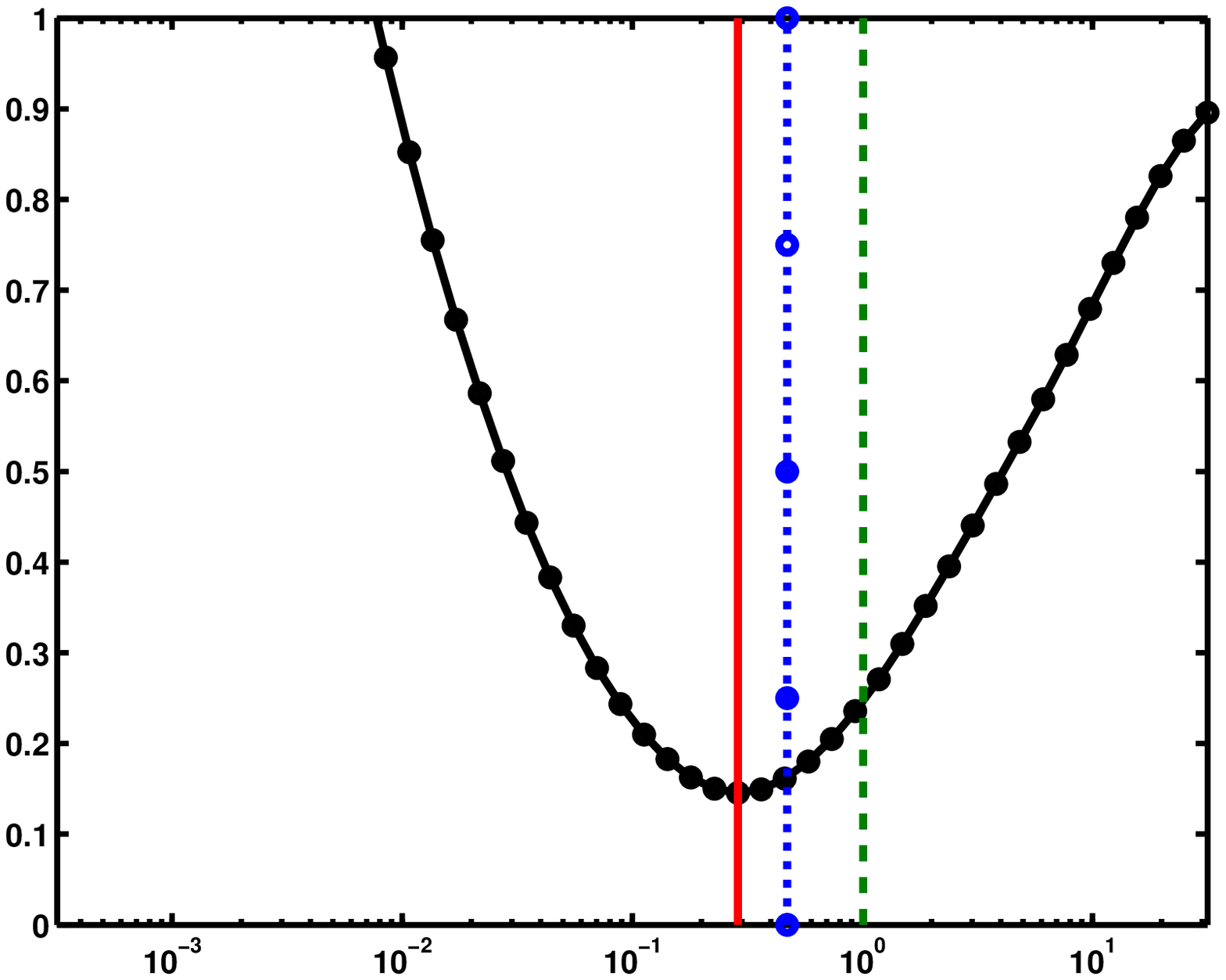}}
\subfigure[$L=L_2$]{\includegraphics[width=1.7in]{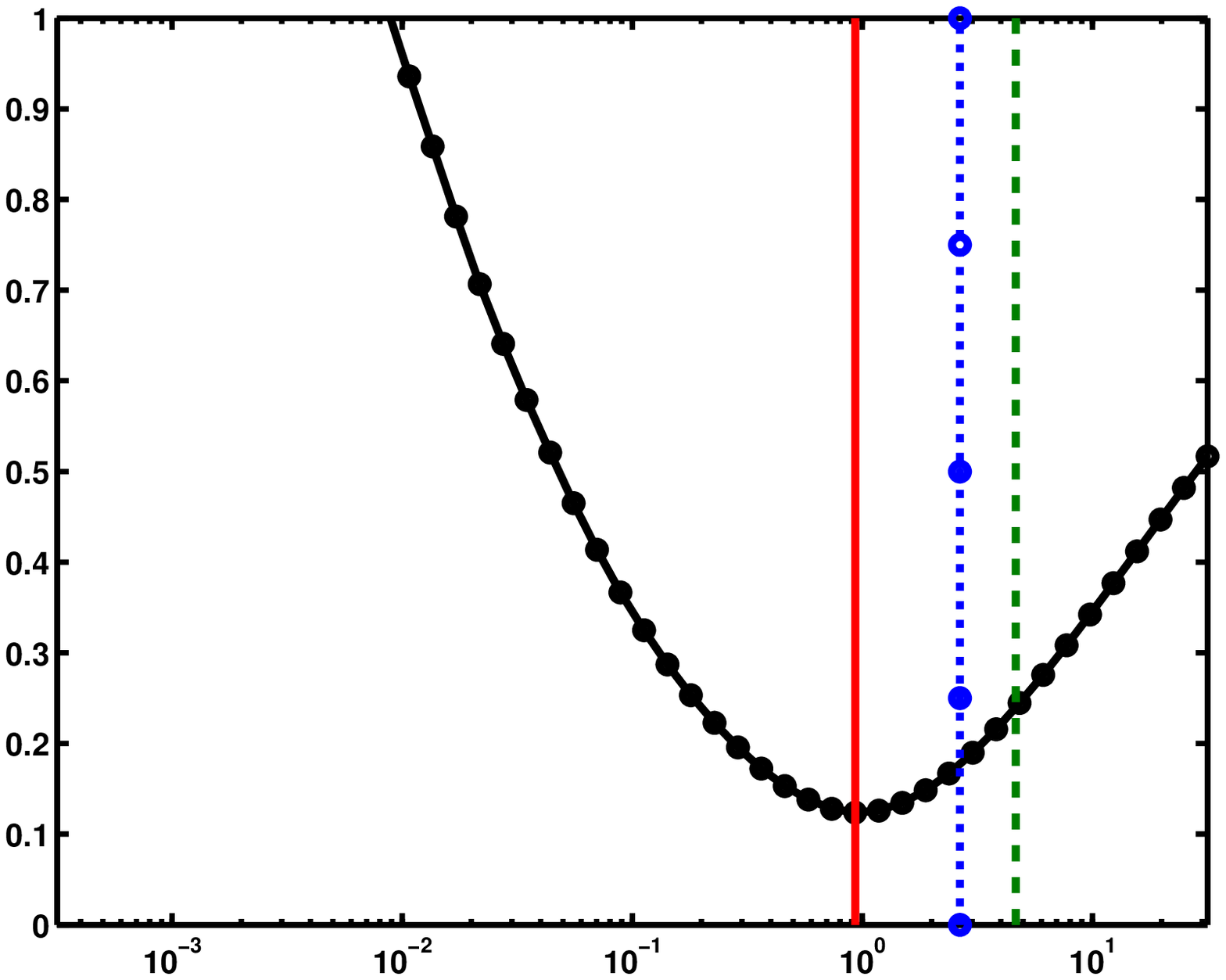}}
\subfigure[$L=I$]{\includegraphics[width=1.7in]{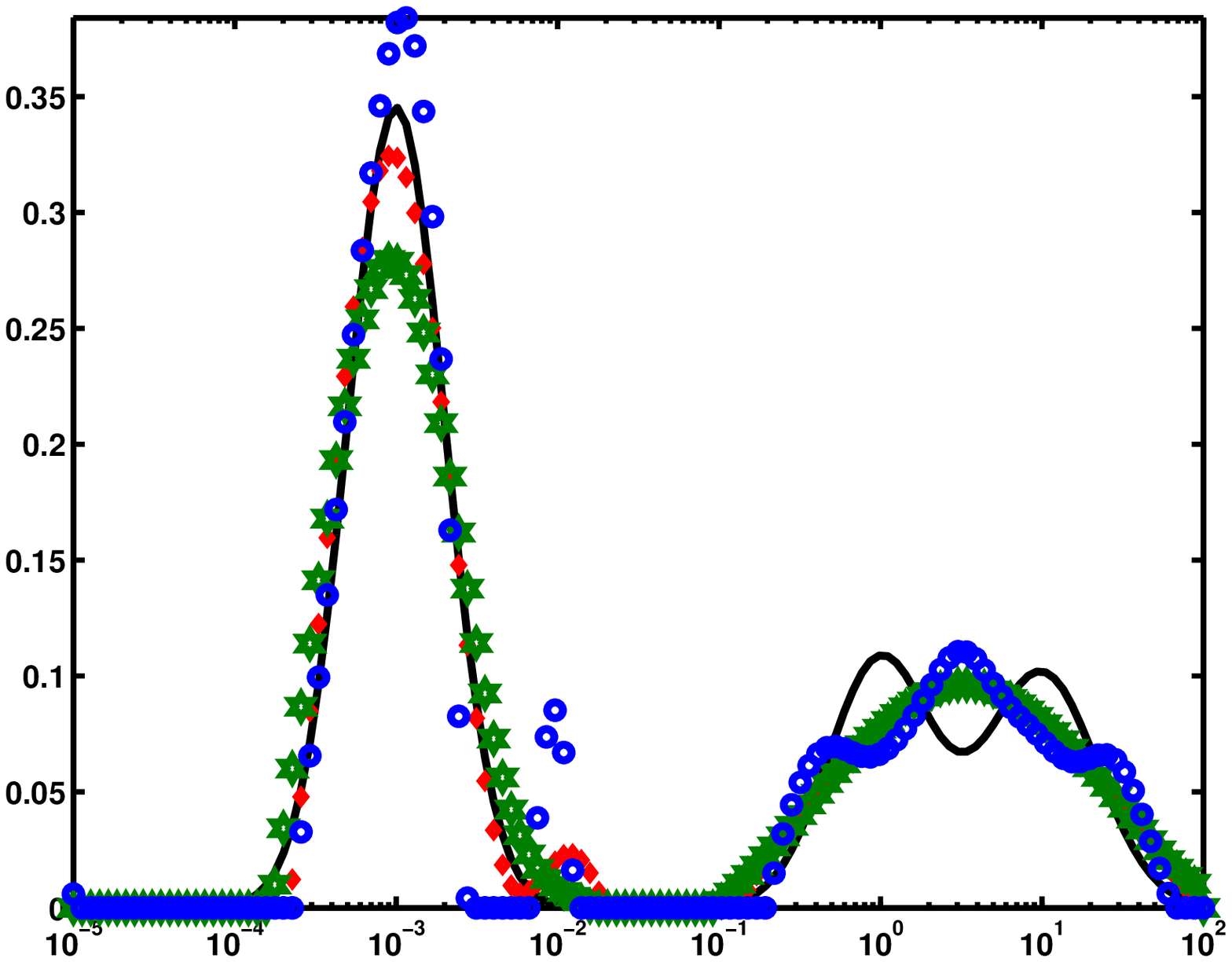}}
\subfigure[$L=L_1$]{\includegraphics[width=1.7in]{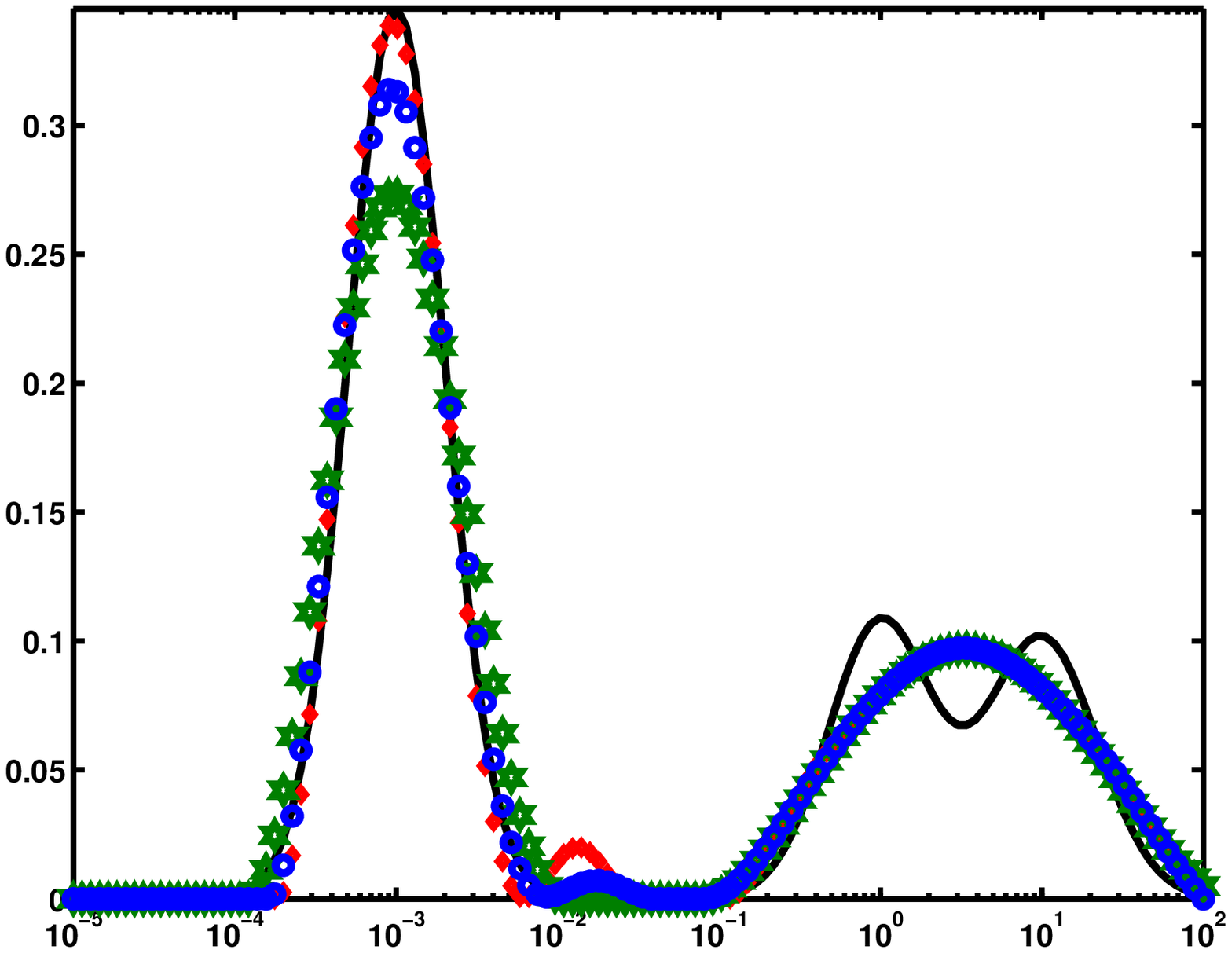}}
\subfigure[$L=L_2$]{\includegraphics[width=1.7in]{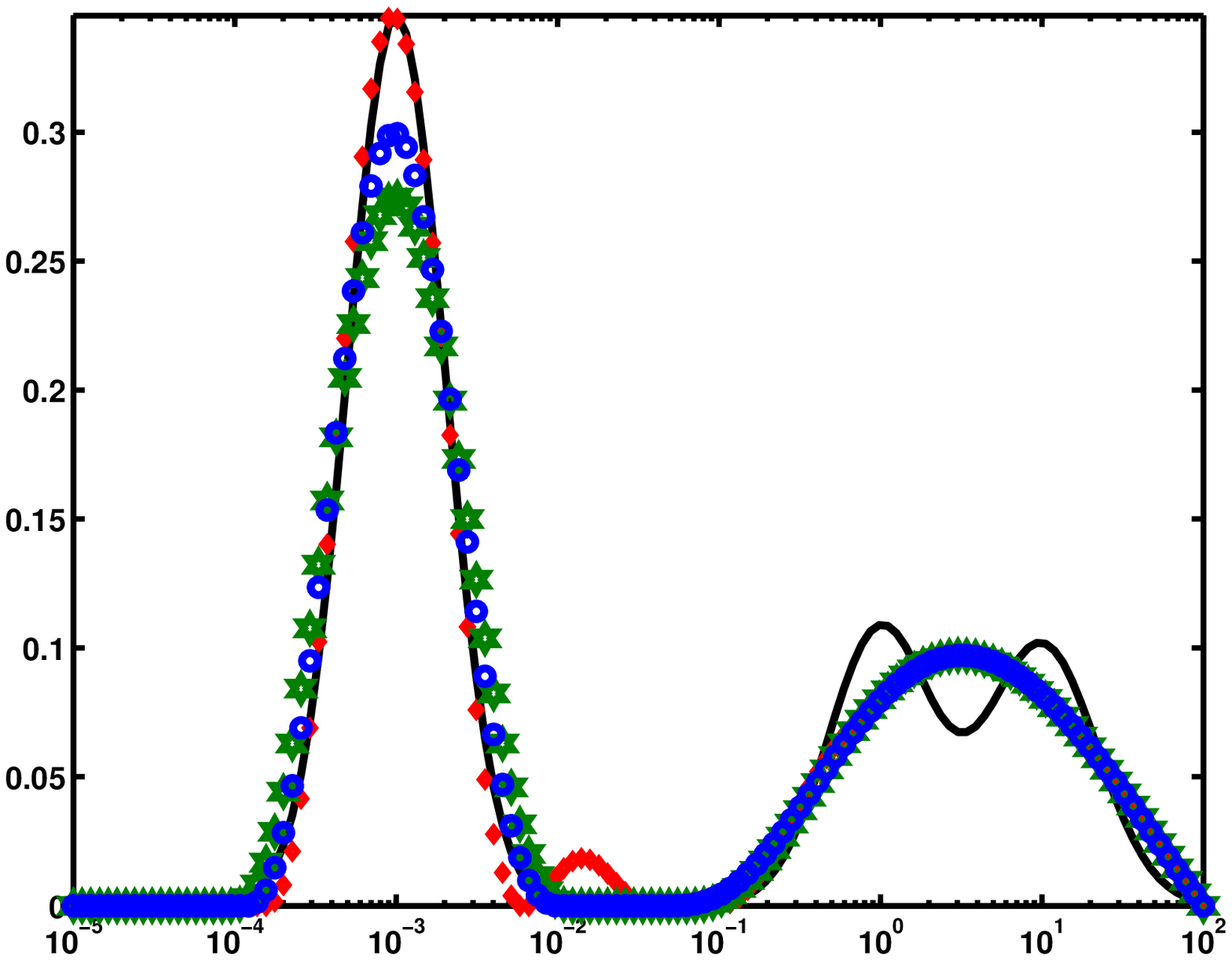}}
\caption{NNLS solutions of LN-C matrix $A_4$. Noise level $1\%$.}
\label{hnfig-lambdachoiceLN6A4HN}
\end{figure}
\clearpage

\subsection{Examples: Noise level $1\%$ matrix  $A_4$ NNLS with CVX}\label{sec:cvx}


\clearpage
\subsection{Examples: Noise level $1\%$ matrix  $A_3$ NNLS}
\begin{figure}[!ht]
\centering
\subfigure[$L=I$]{\includegraphics[width=1.7in]{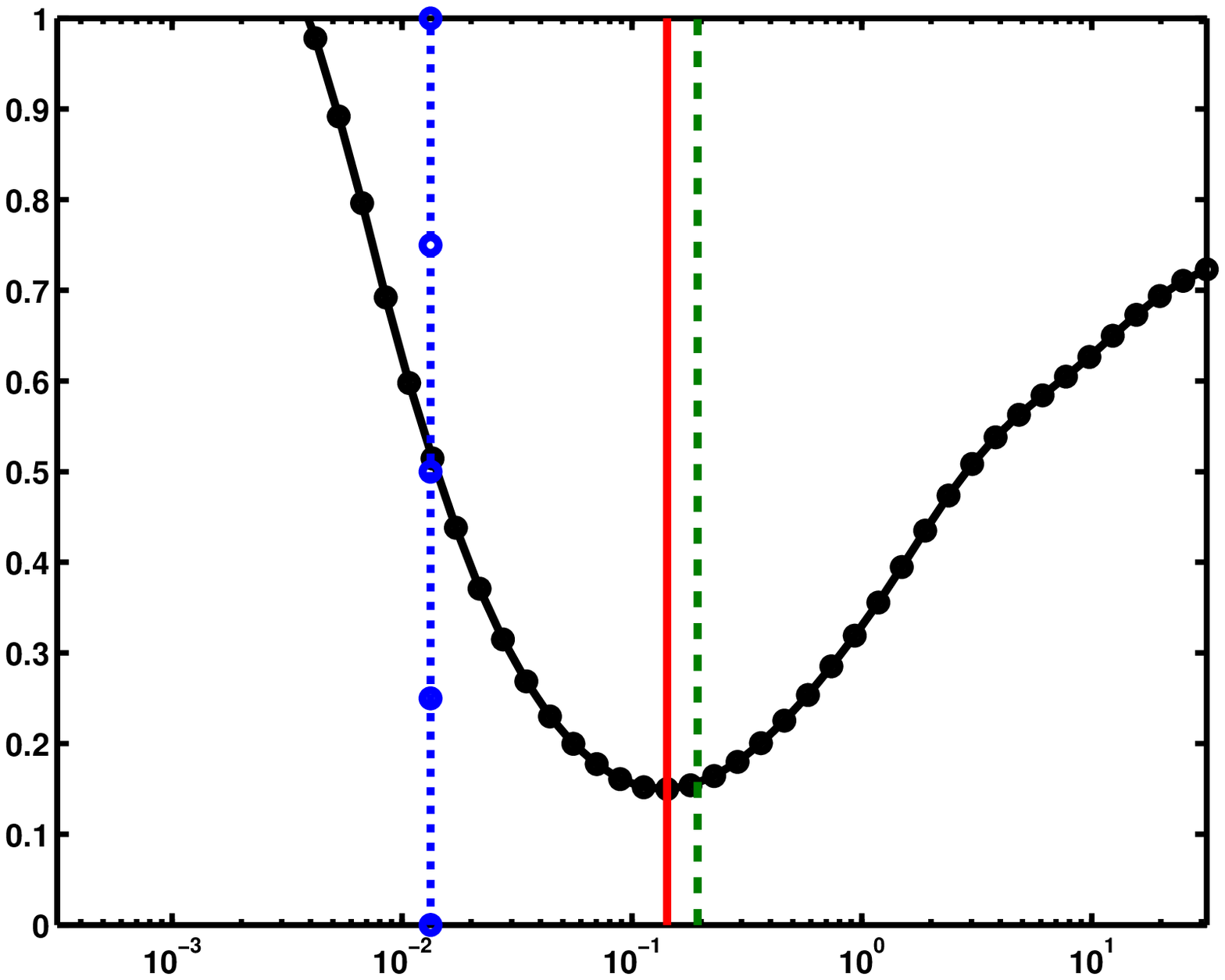}}
\subfigure[$L=L_1$]{\includegraphics[width=1.7in]{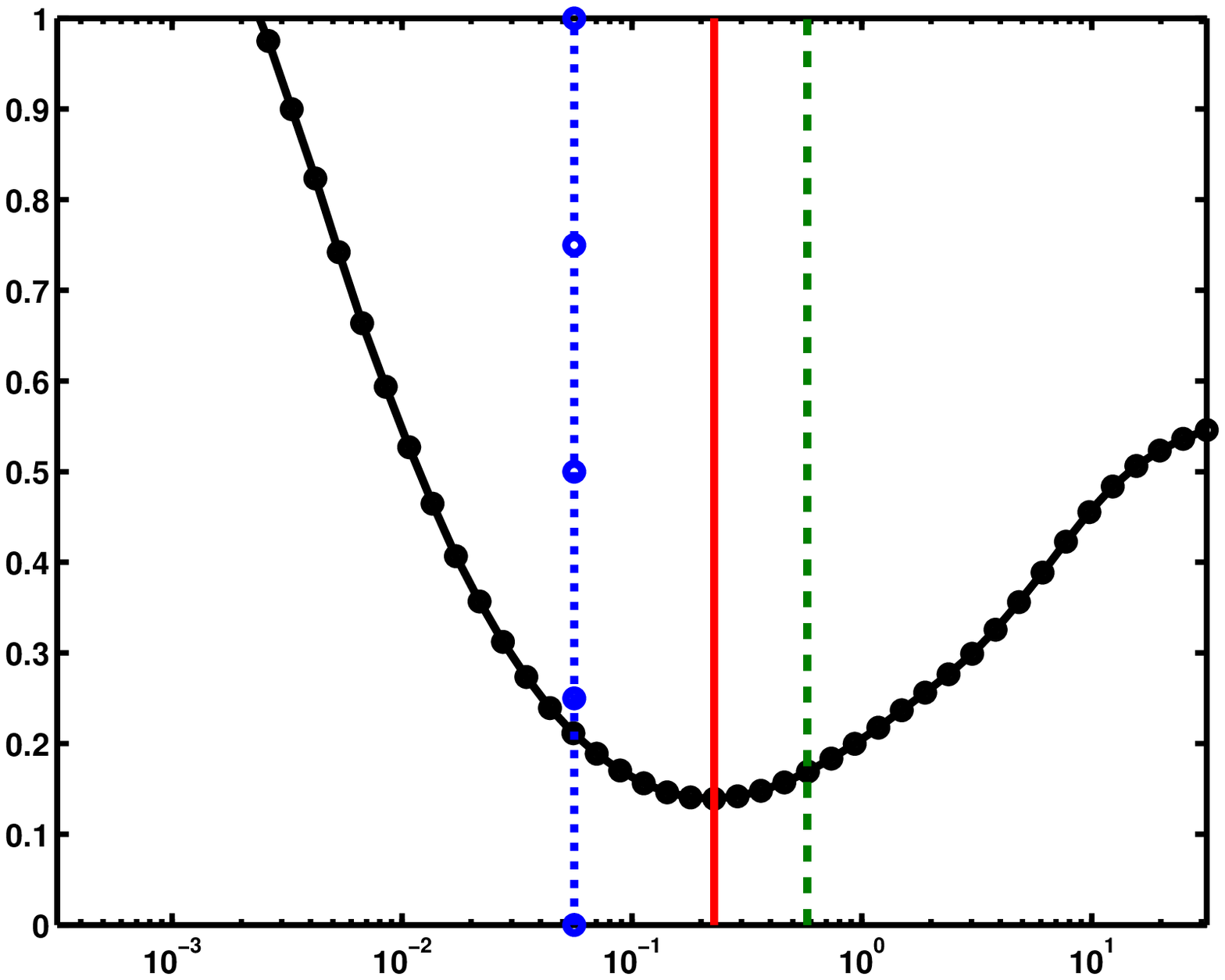}}
\subfigure[$L=L_2$]{\includegraphics[width=1.7in]{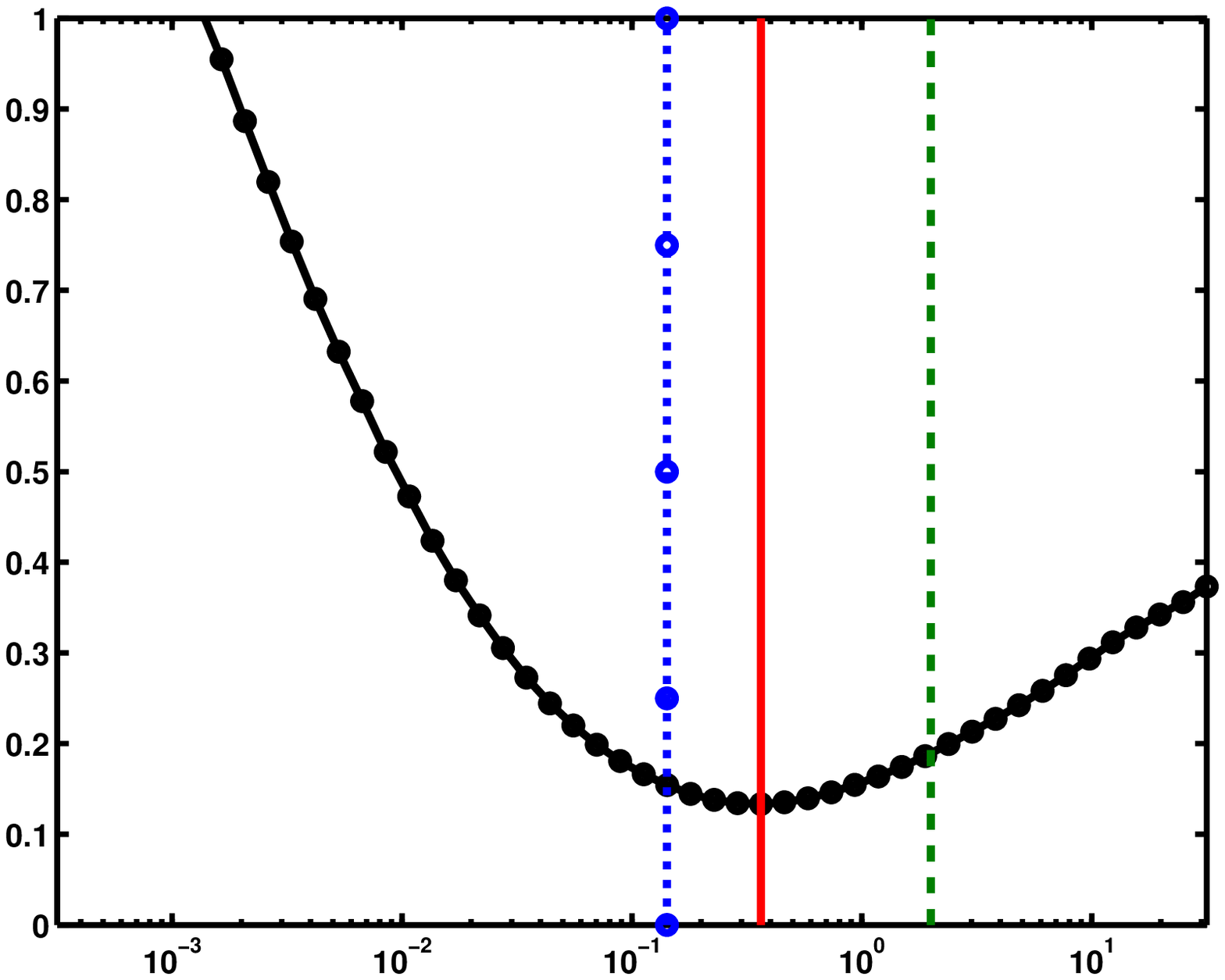}}
\subfigure[$L=I$]{\includegraphics[width=1.7in]{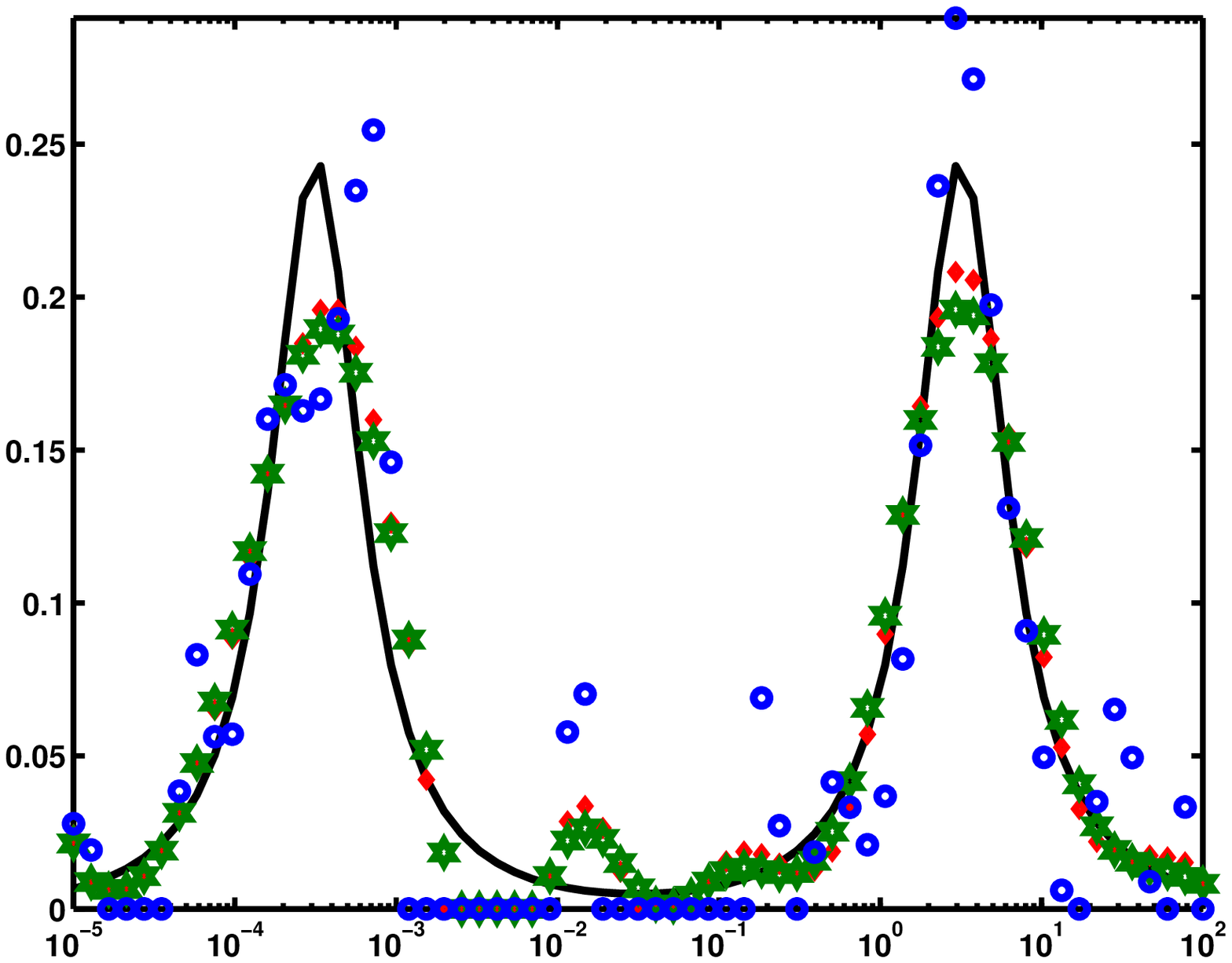}}
\subfigure[$L=L_1$]{\includegraphics[width=1.7in]{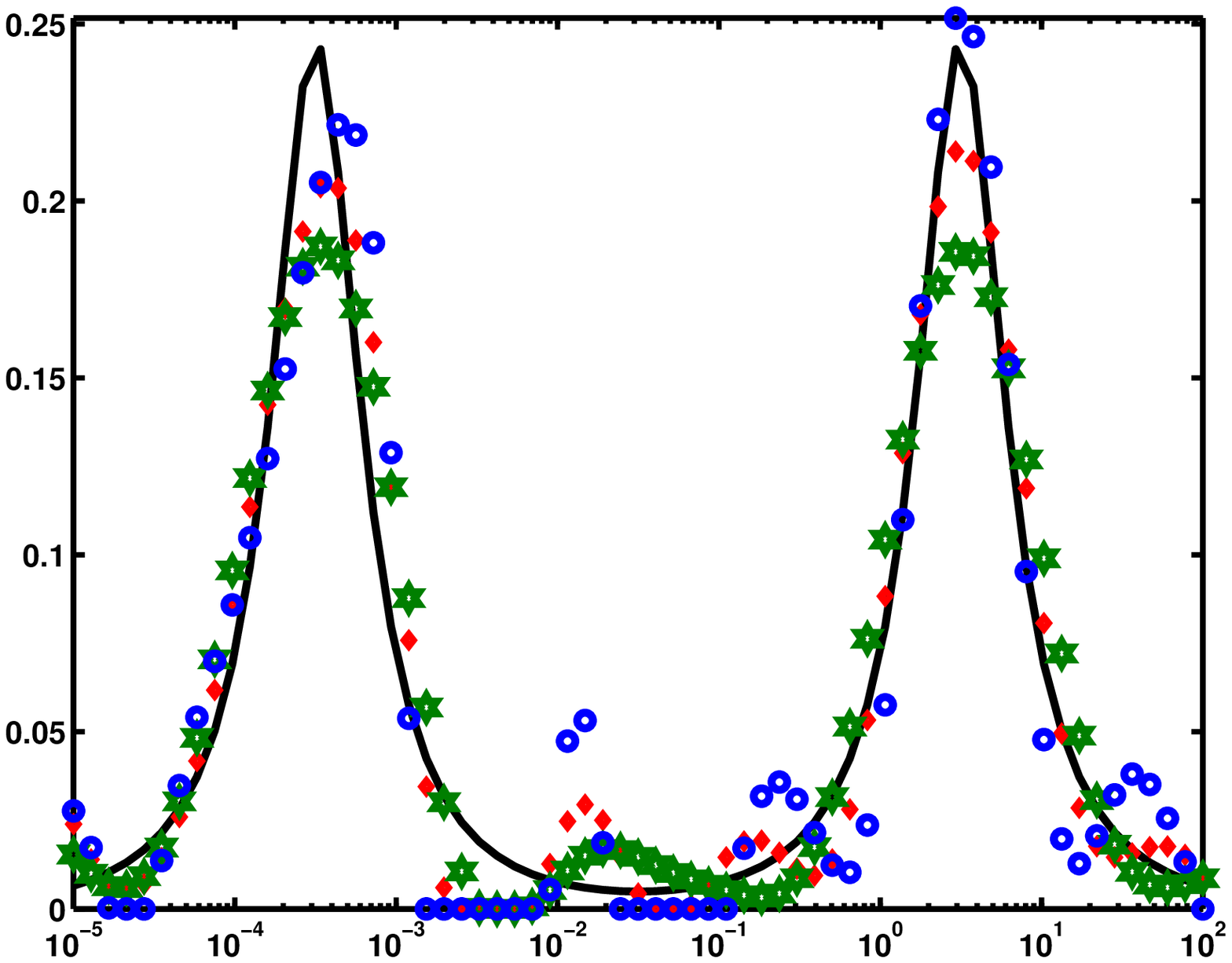}}
\subfigure[$L=L_2$]{\includegraphics[width=1.7in]{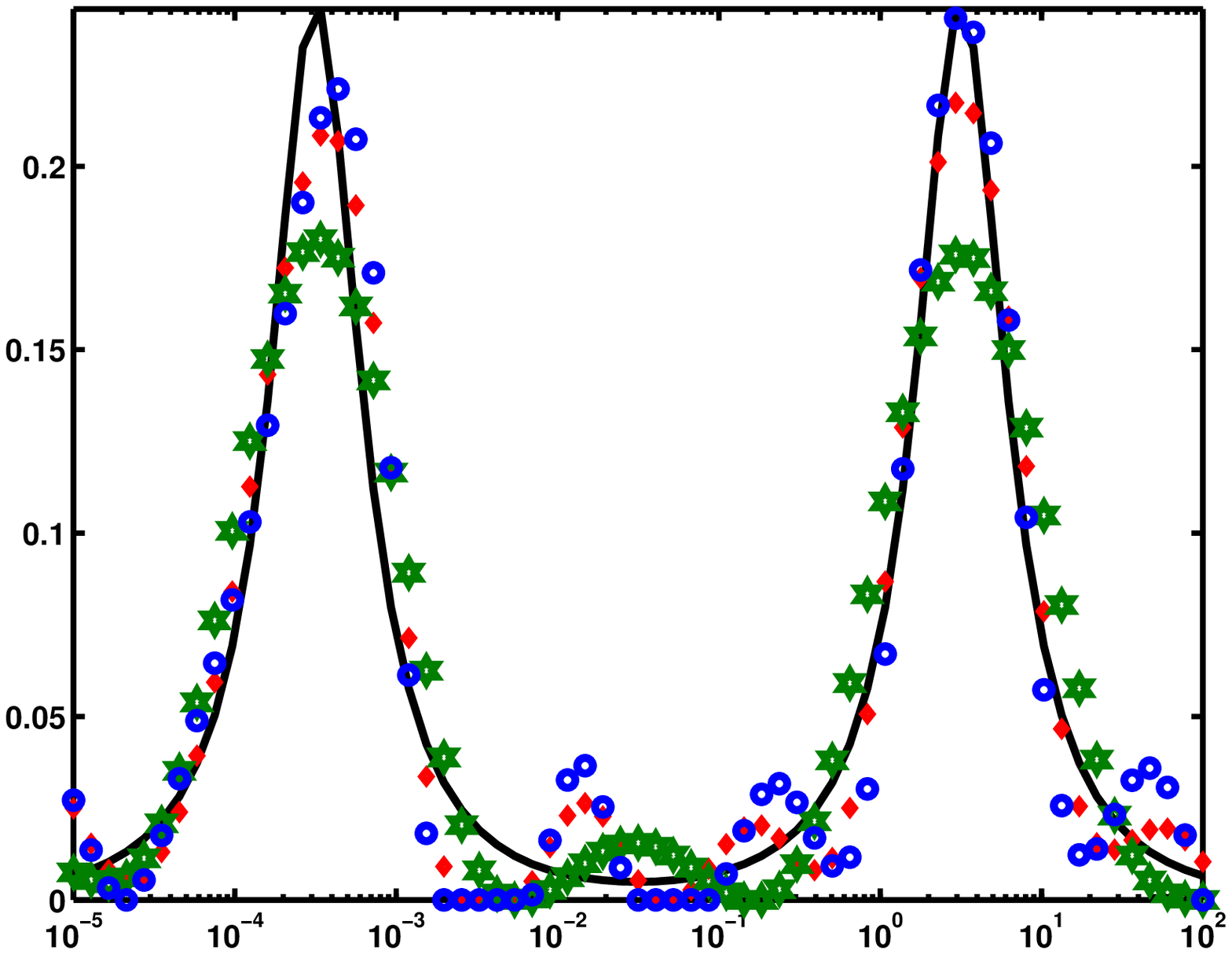}}
\caption{NNLS solutions of RQ-A matrix $A_3$. Noise level $1\%$.}
\label{hnfig-lambdachoiceRQ1A3HN}
\end{figure}
\begin{figure}[!ht]
\centering
\subfigure[$L=I$]{\includegraphics[width=1.7in]{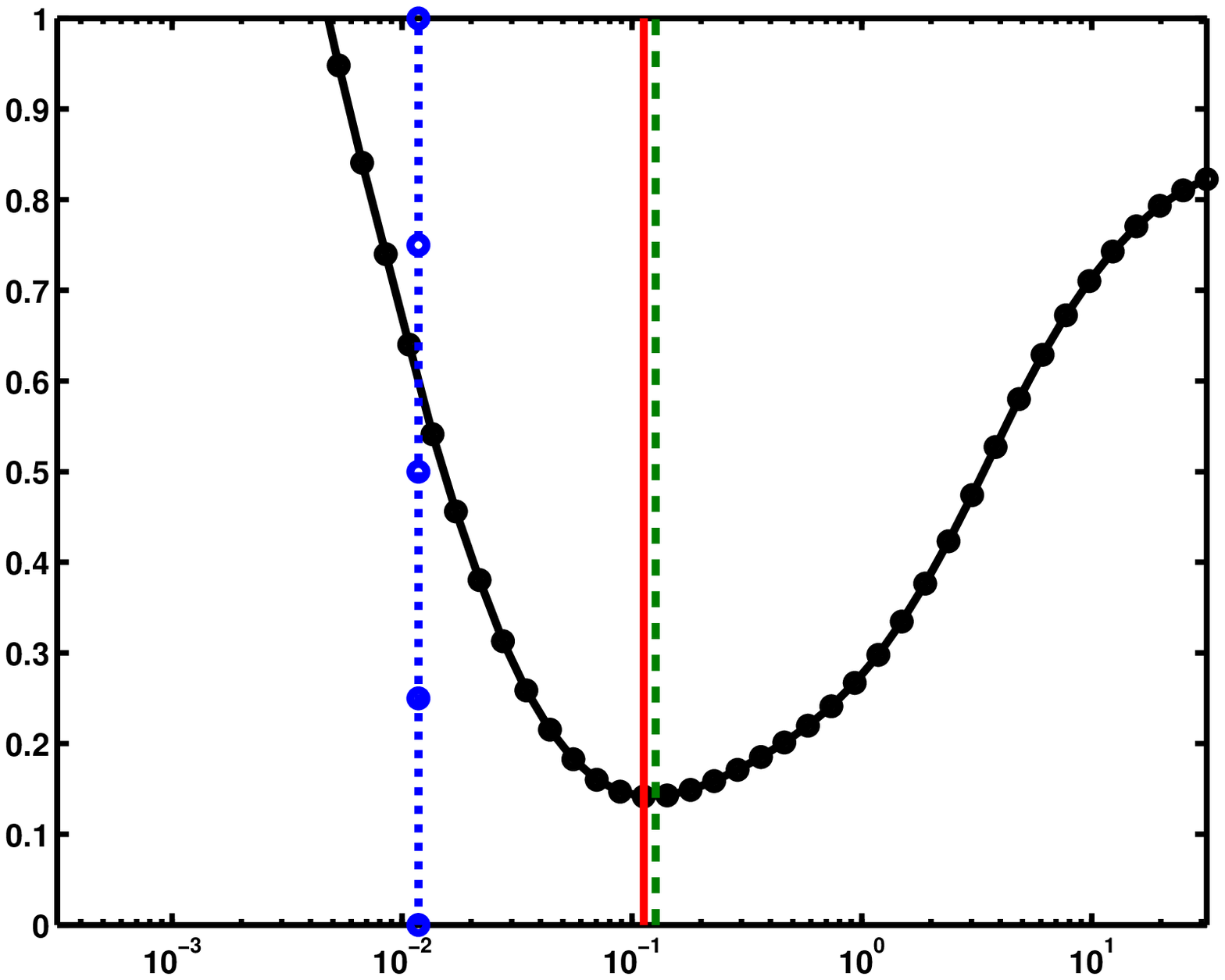}}
\subfigure[$L=L_1$]{\includegraphics[width=1.7in]{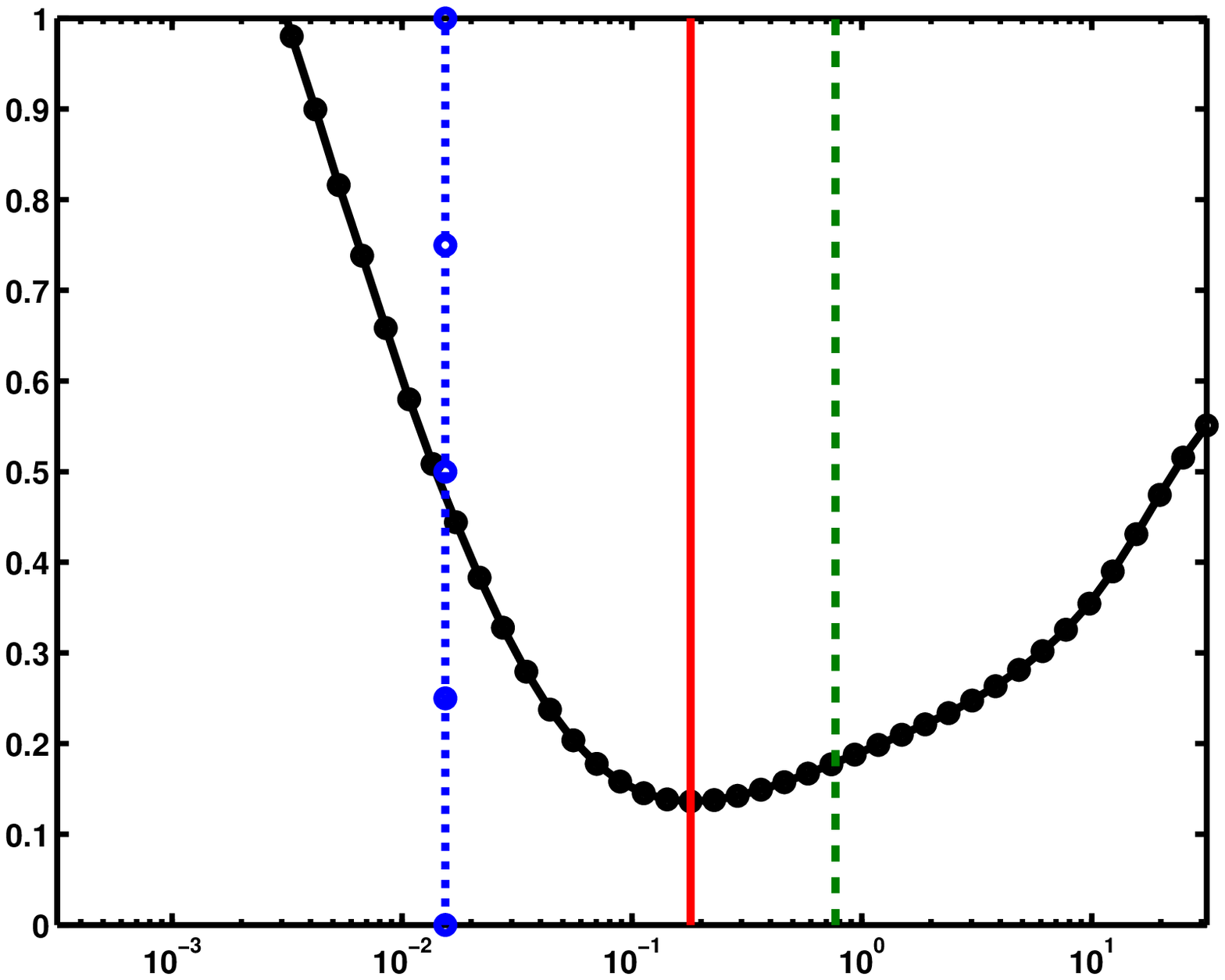}}
\subfigure[$L=L_2$]{\includegraphics[width=1.7in]{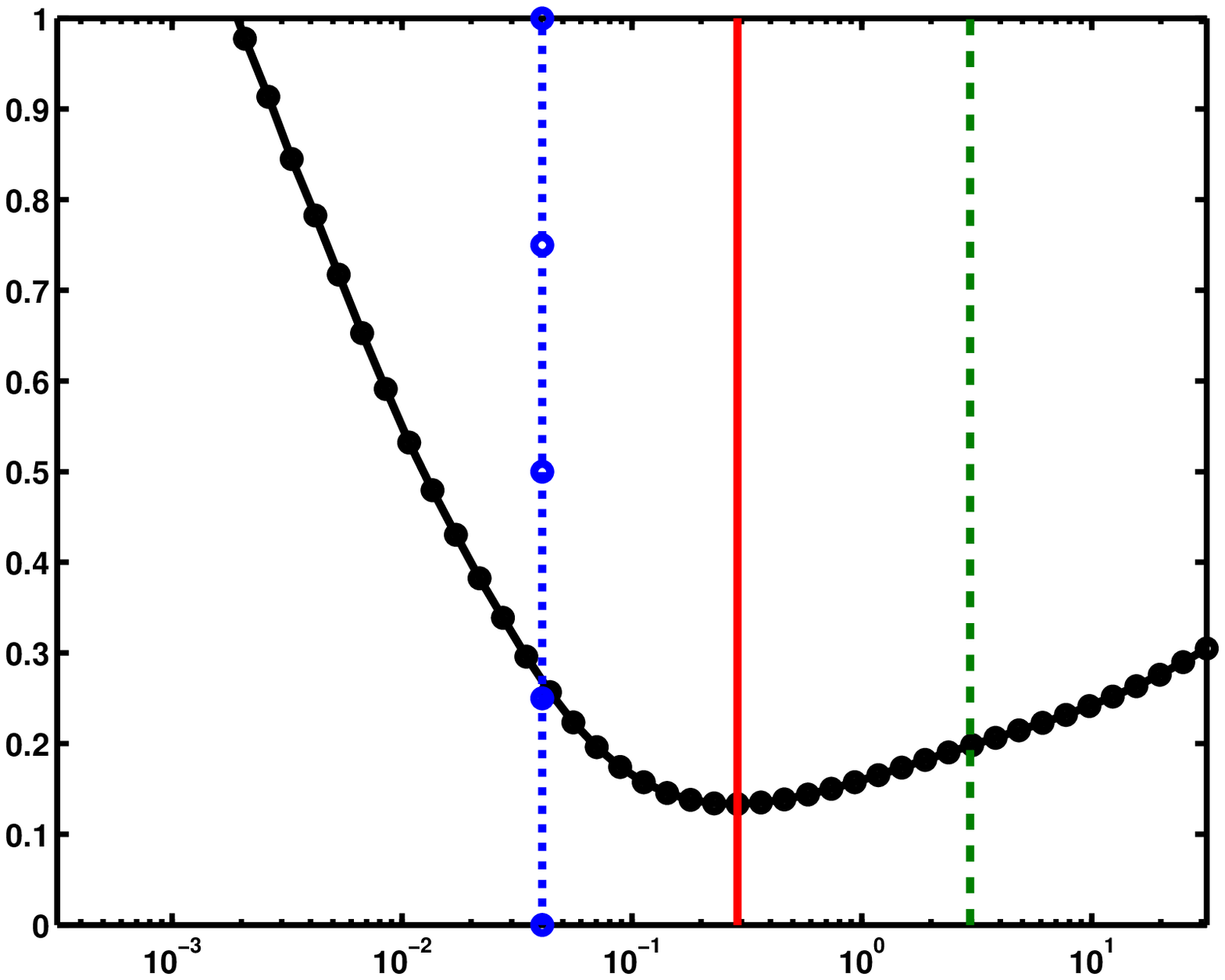}}
\subfigure[$L=I$]{\includegraphics[width=1.7in]{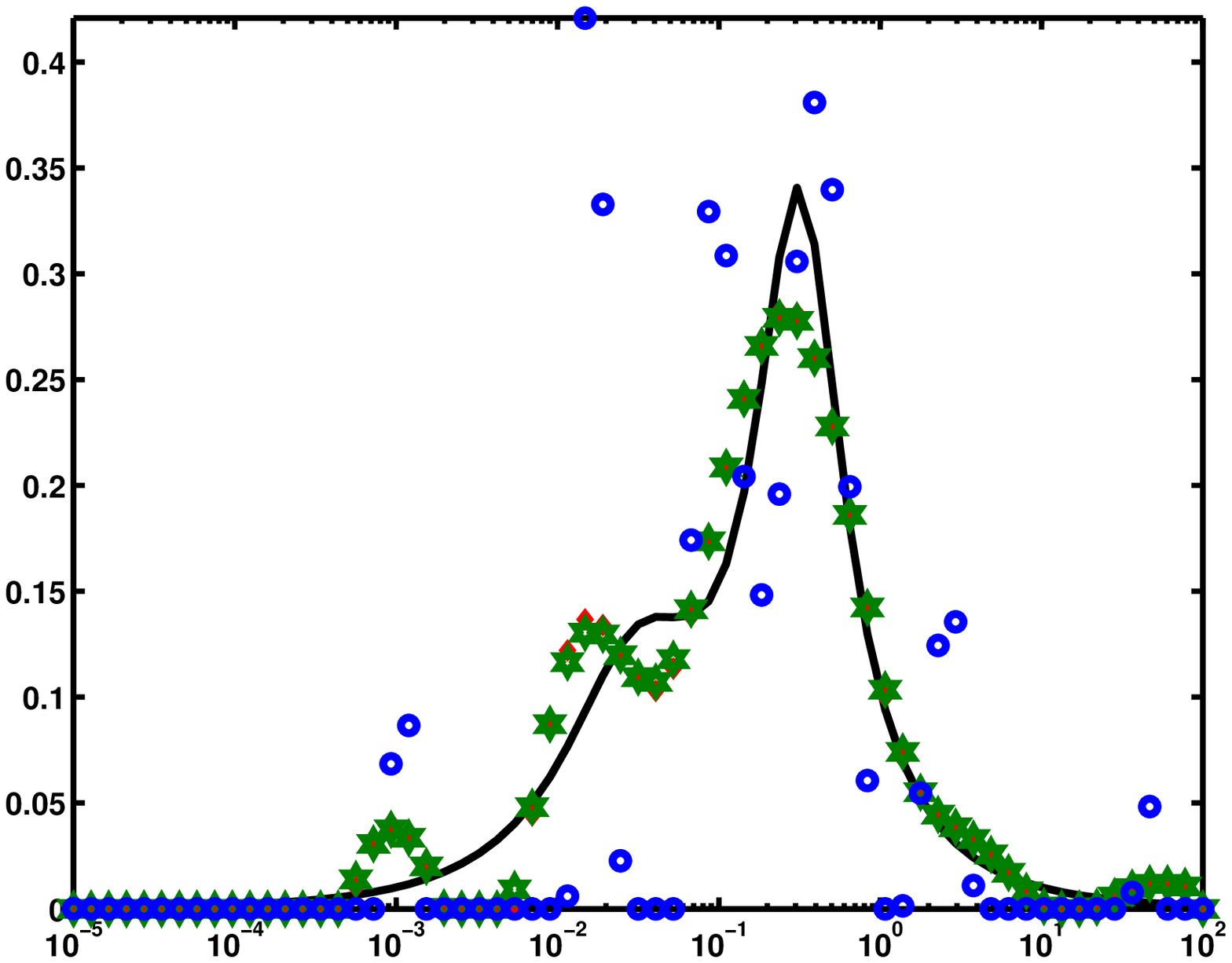}}
\subfigure[$L=L_1$]{\includegraphics[width=1.7in]{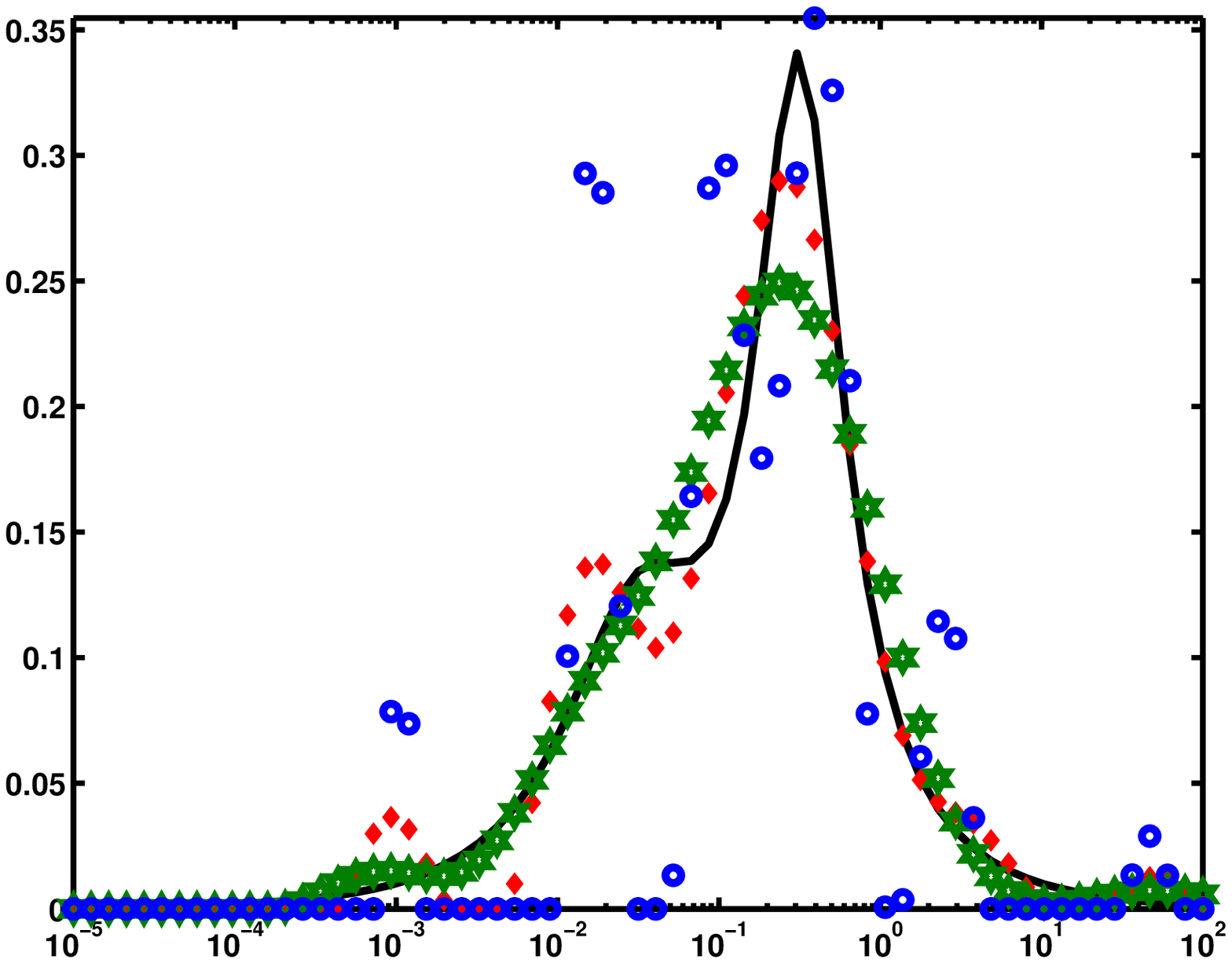}}
\subfigure[$L=L_2$]{\includegraphics[width=1.7in]{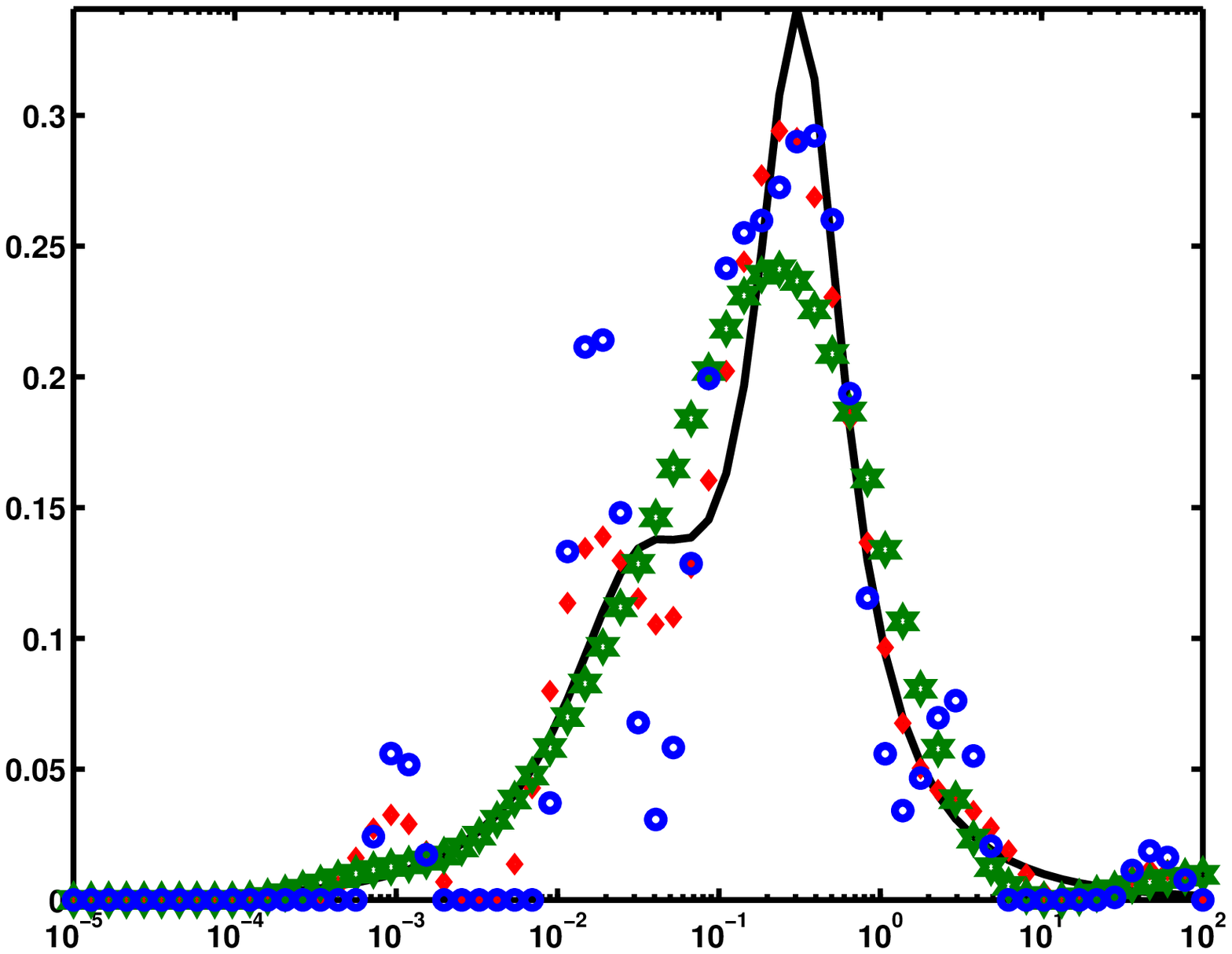}}
\caption{NNLS solutions of RQ-B matrix $A_3$. Noise level $1\%$.}
\label{hnfig-lambdachoiceRQ5A3HN}
\end{figure}

\begin{figure}[!ht]
\centering
\subfigure[$L=I$]{\includegraphics[width=1.7in]{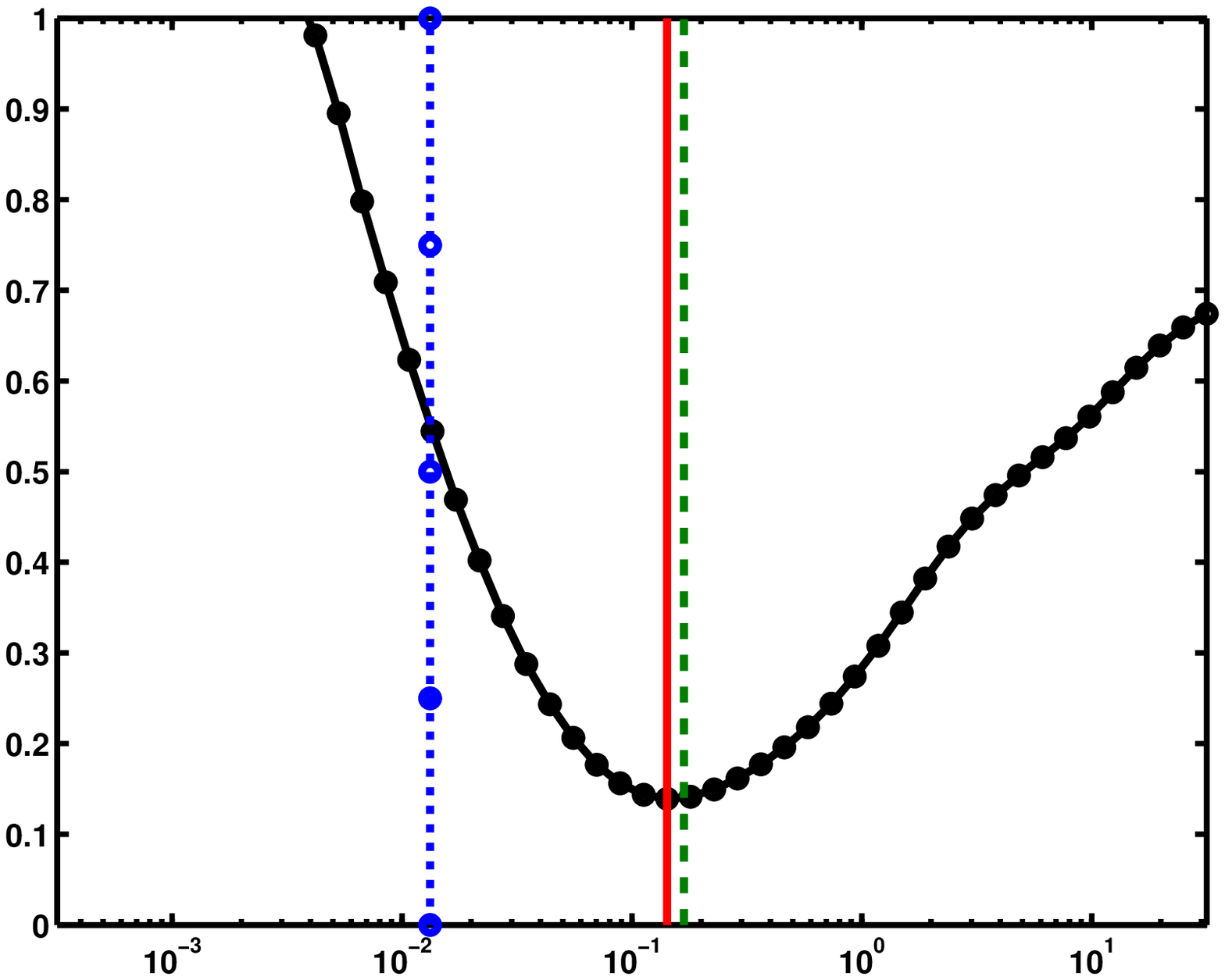}}
\subfigure[$L=L_1$]{\includegraphics[width=1.7in]{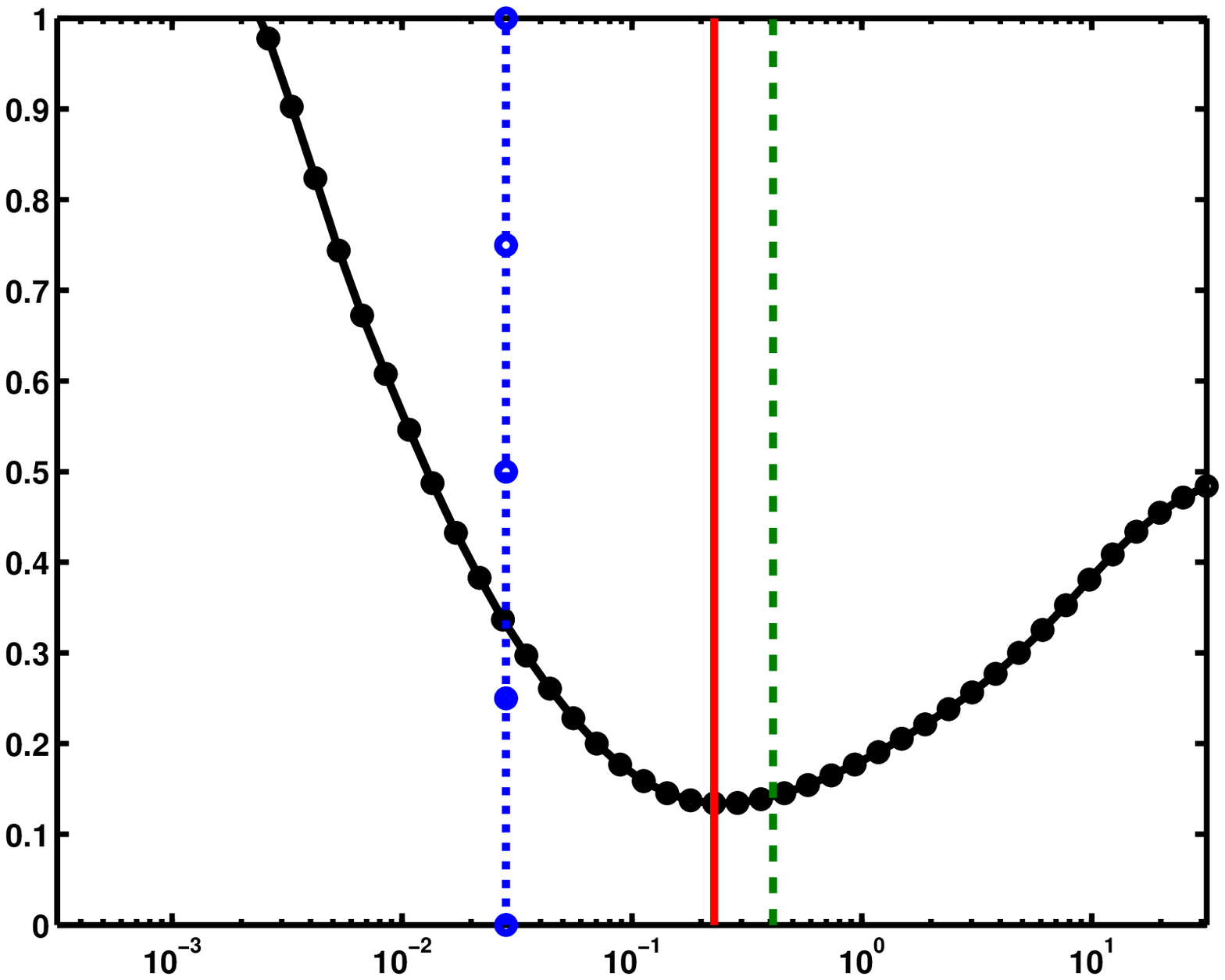}}
\subfigure[$L=L_2$]{\includegraphics[width=1.7in]{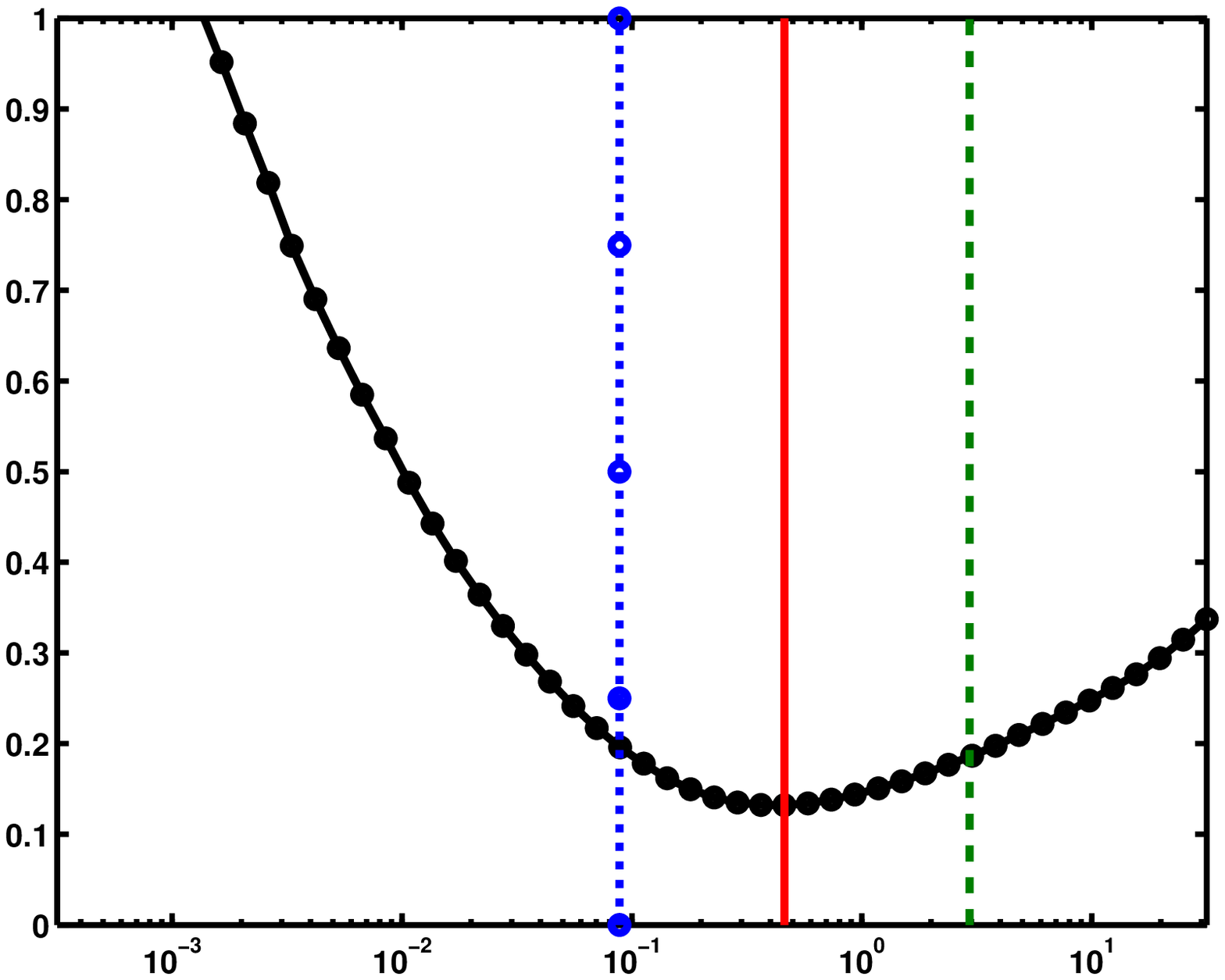}}
\subfigure[$L=I$]{\includegraphics[width=1.7in]{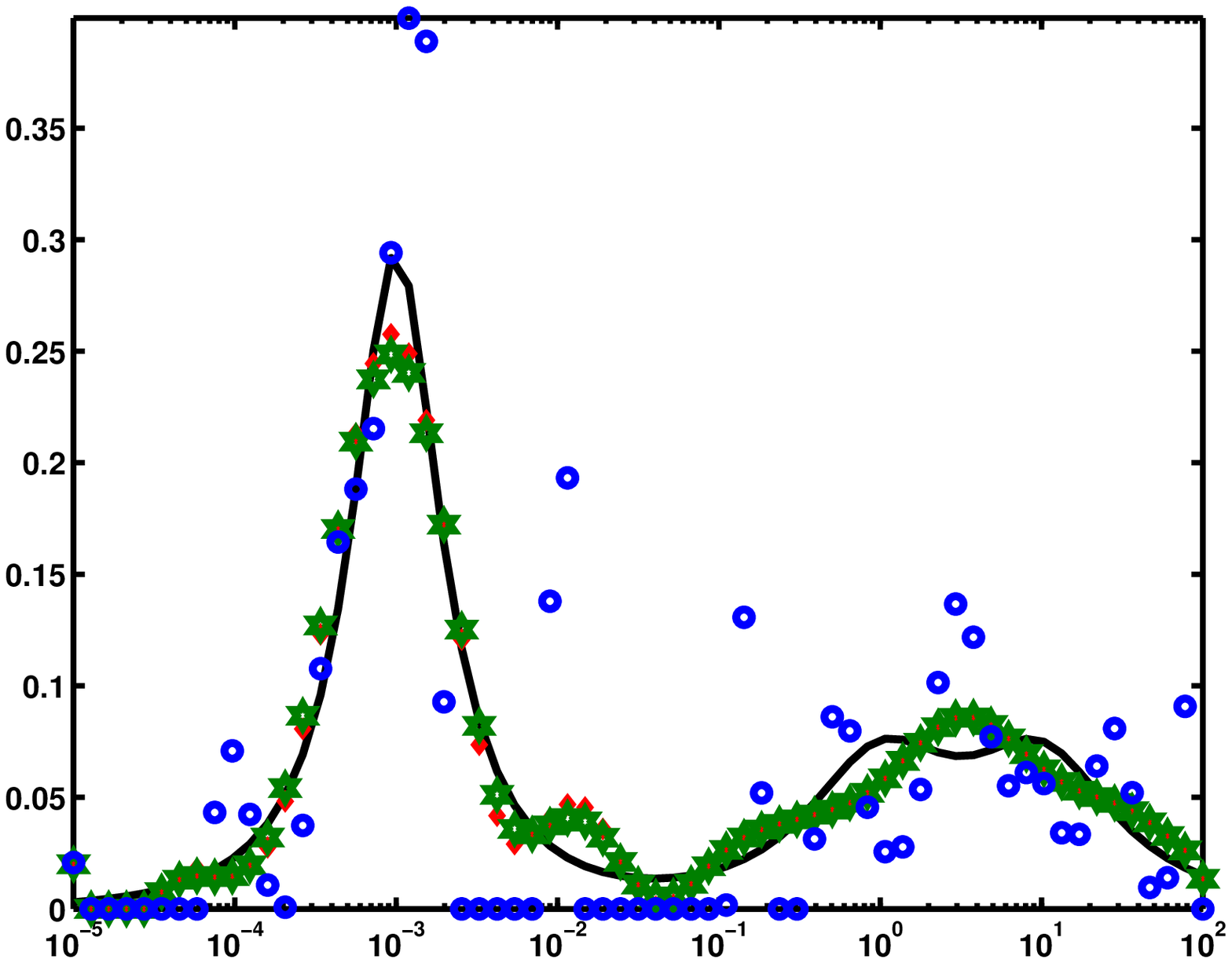}}
\subfigure[$L=L_1$]{\includegraphics[width=1.7in]{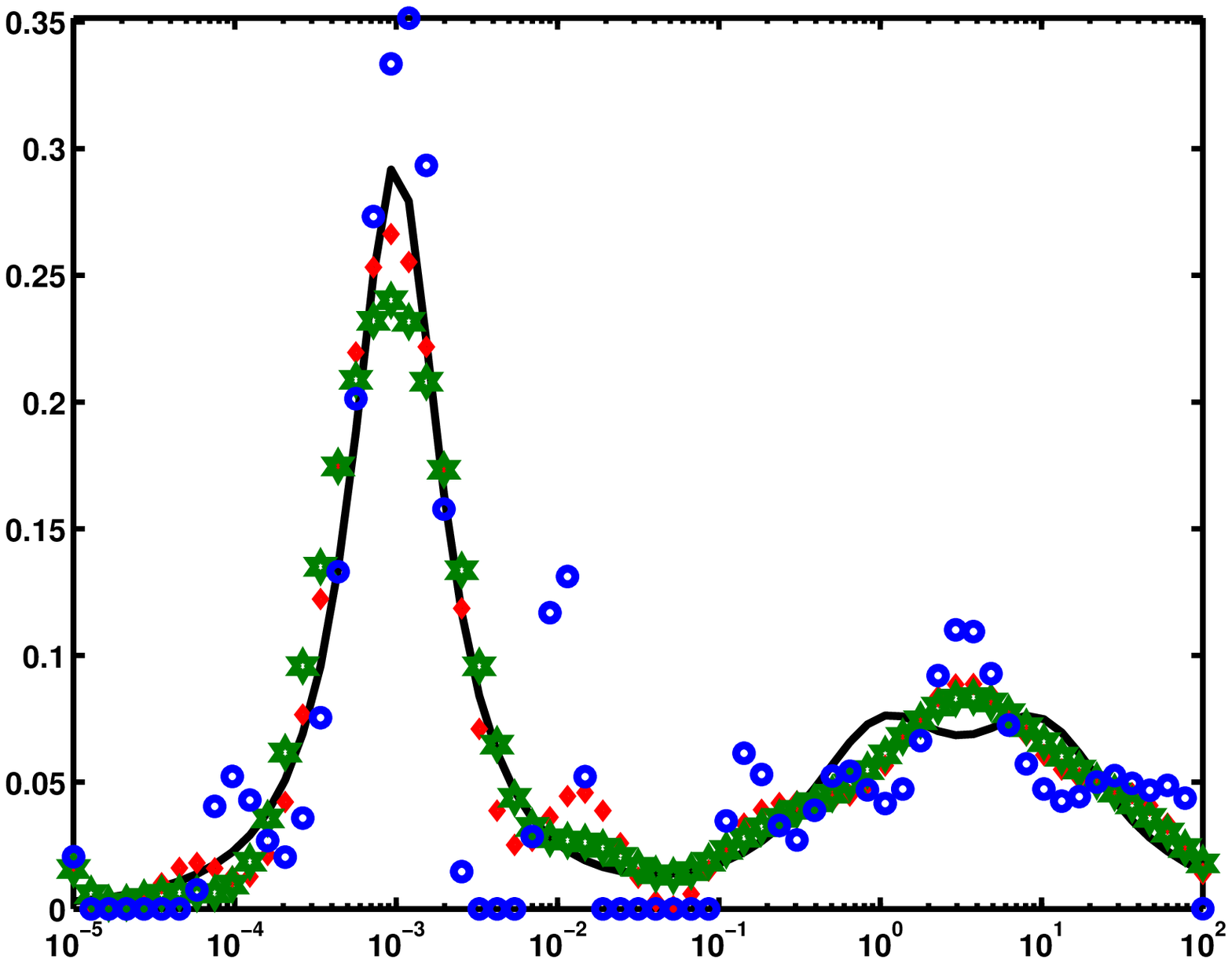}}
\subfigure[$L=L_2$]{\includegraphics[width=1.7in]{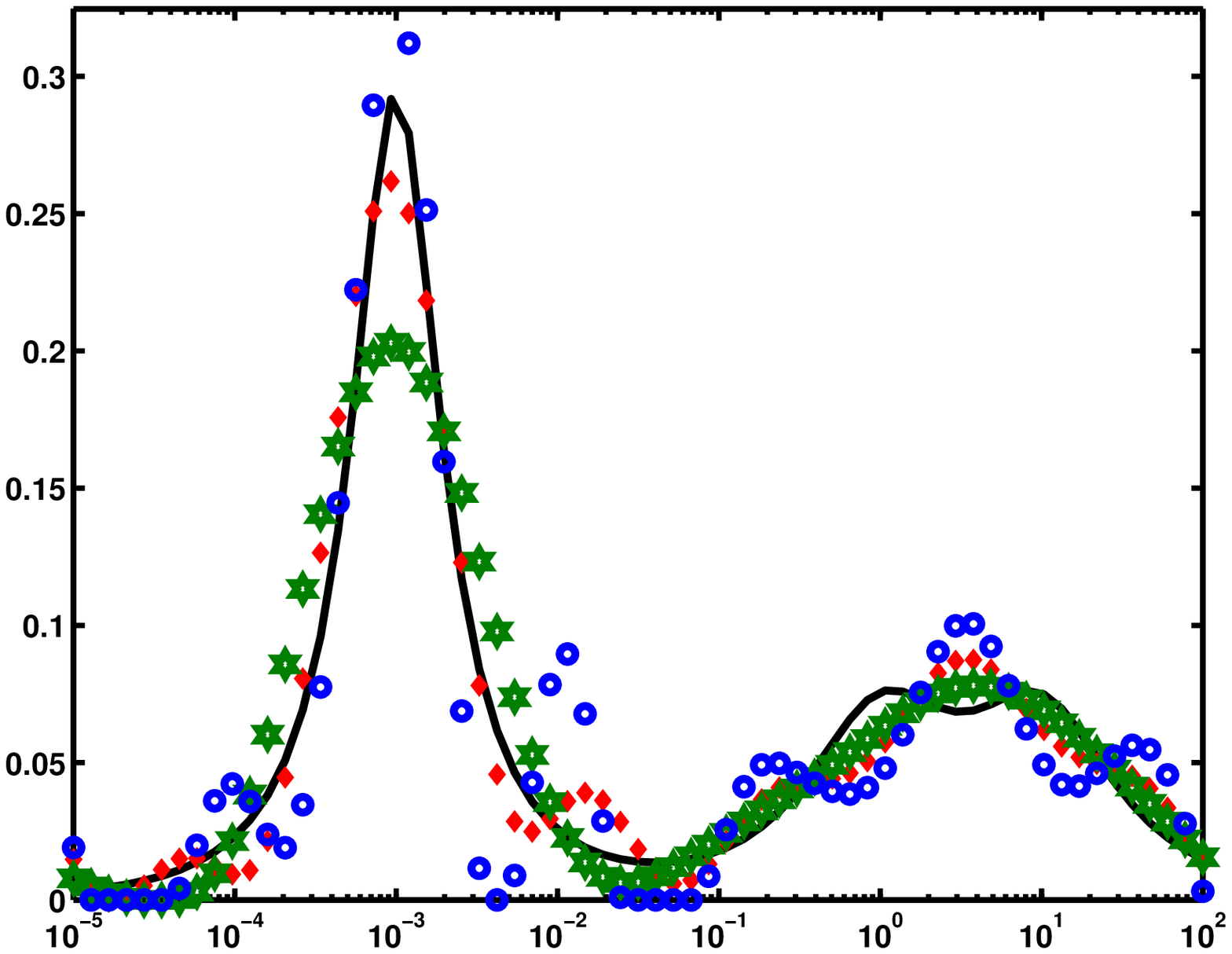}}
\caption{NNLS solutions of RQ-C matrix $A_3$. Noise level $1\%$.}
\label{hnfig-lambdachoiceRQ6A3HN}
\end{figure}

\begin{figure}[!ht]
\centering
\subfigure[$L=I$]{\includegraphics[width=1.7in]{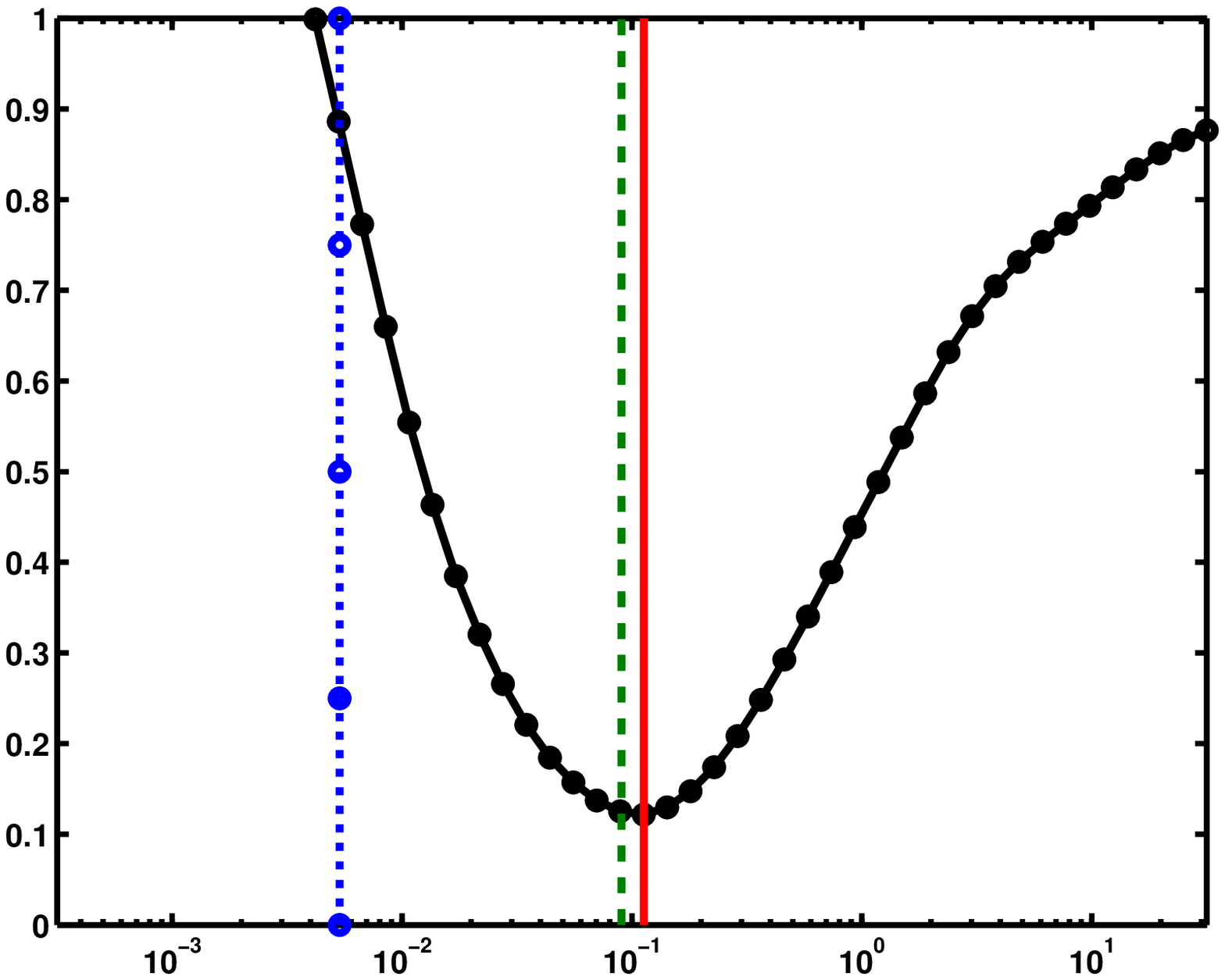}}
\subfigure[$L=L_1$]{\includegraphics[width=1.7in]{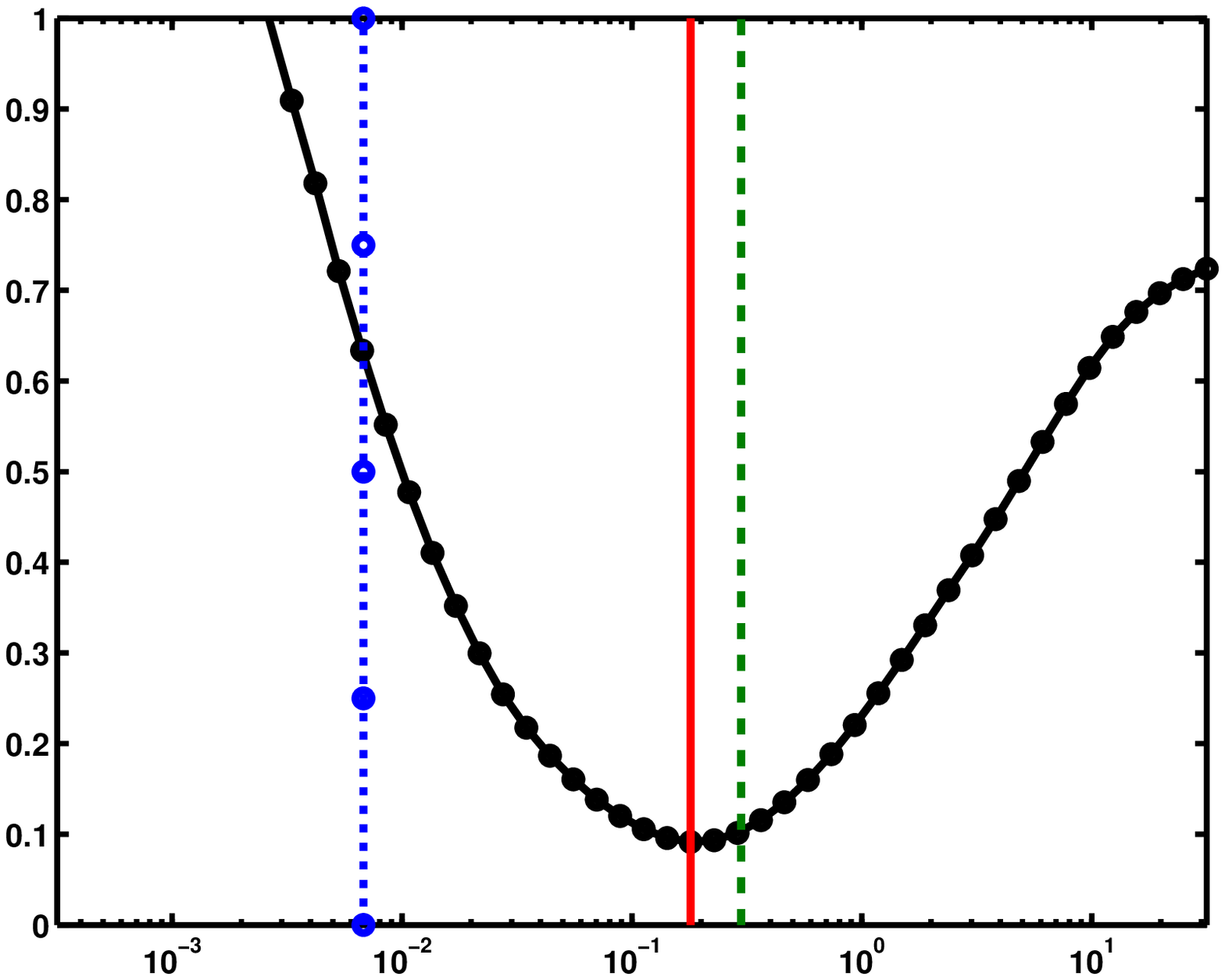}}
\subfigure[$L=L_2$]{\includegraphics[width=1.7in]{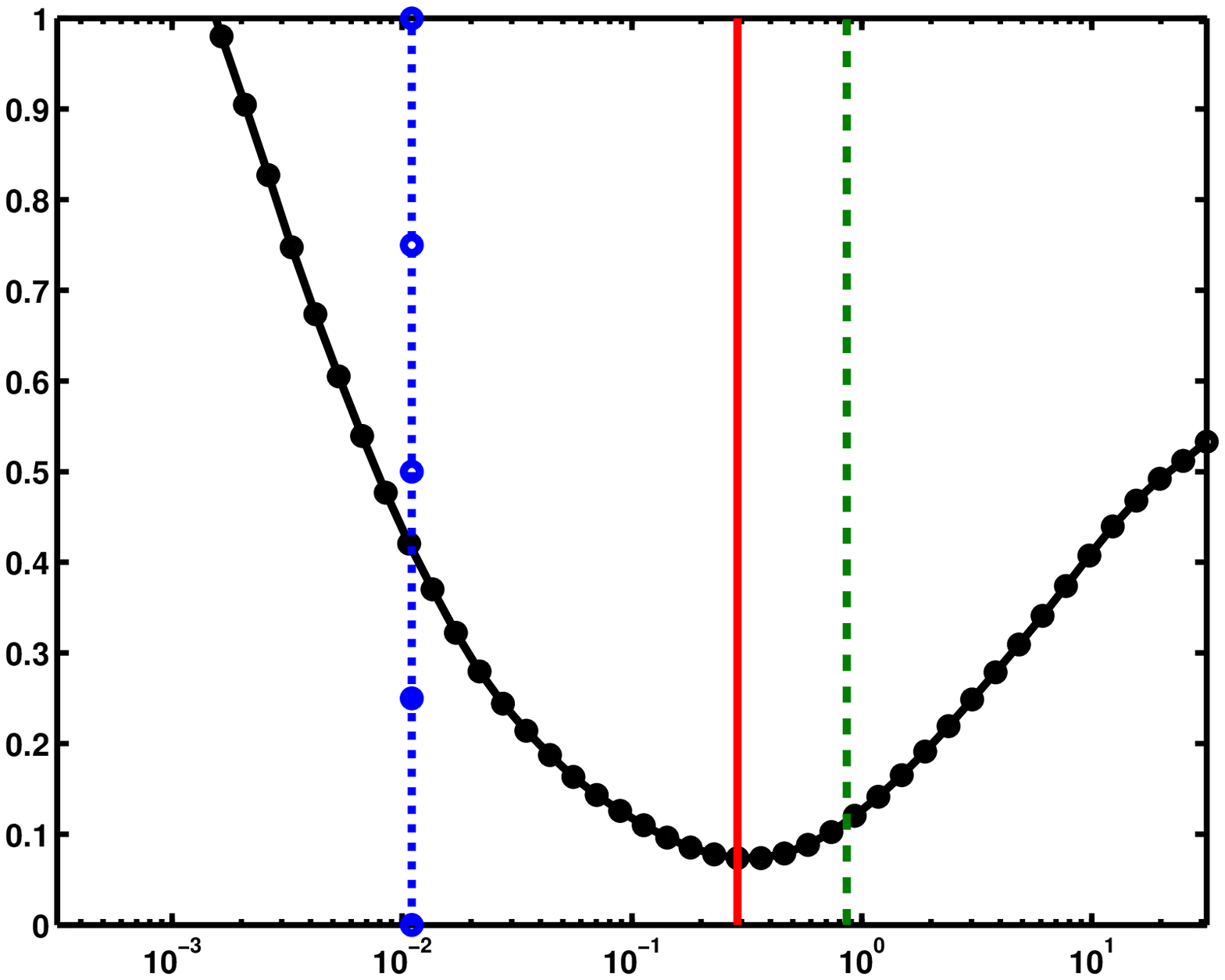}}
\subfigure[$L=I$]{\includegraphics[width=1.7in]{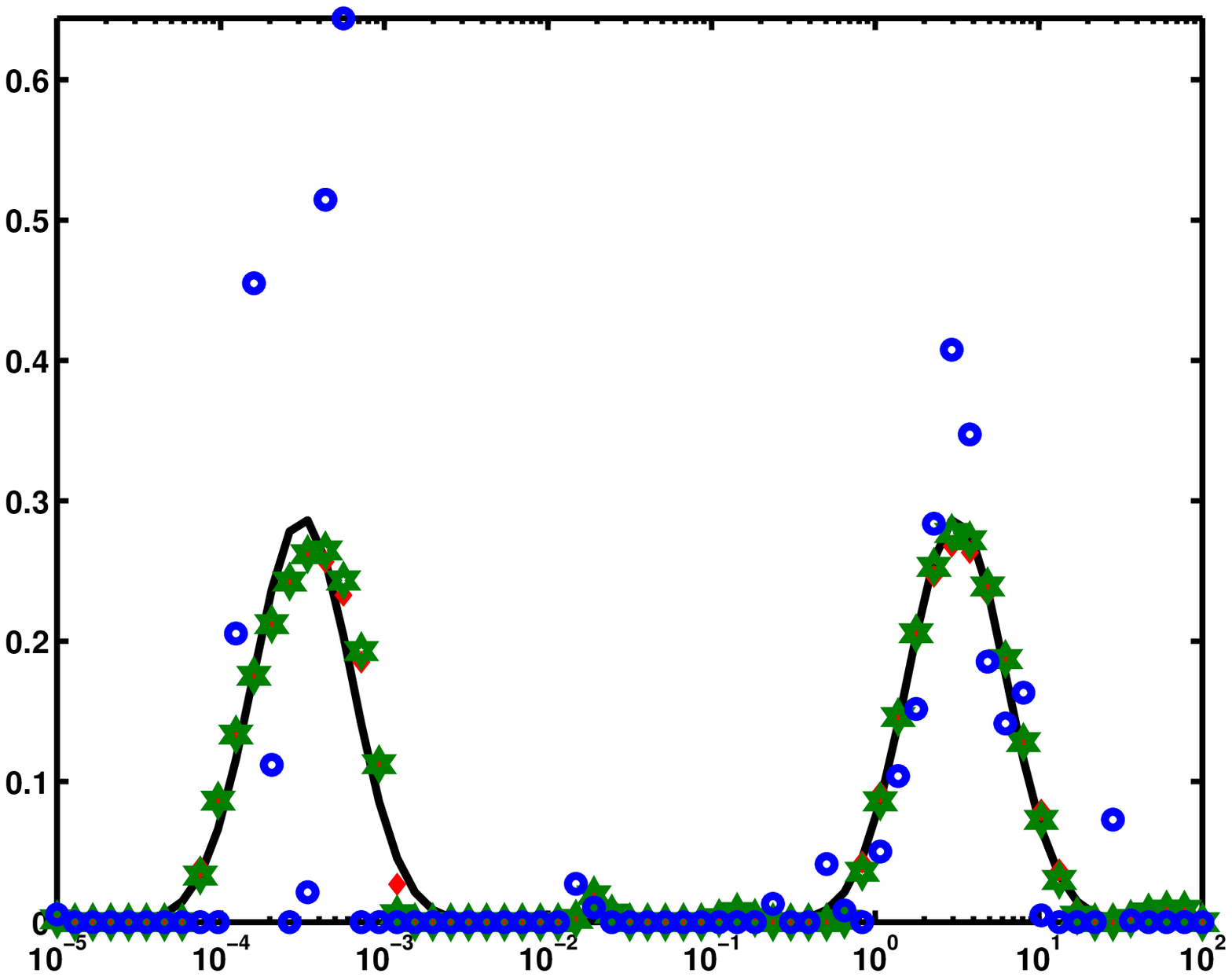}}
\subfigure[$L=L_1$]{\includegraphics[width=1.7in]{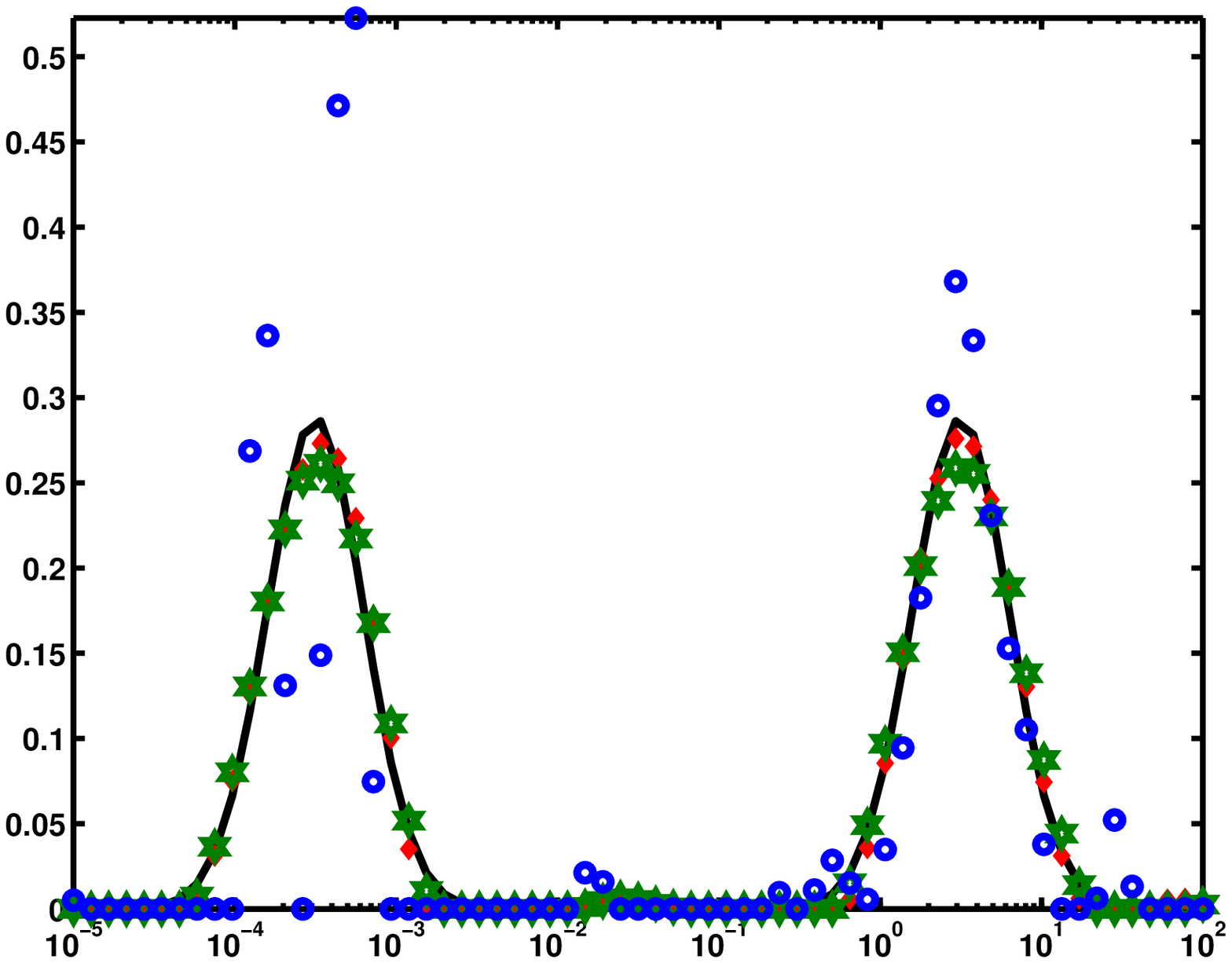}}
\subfigure[$L=L_2$]{\includegraphics[width=1.7in]{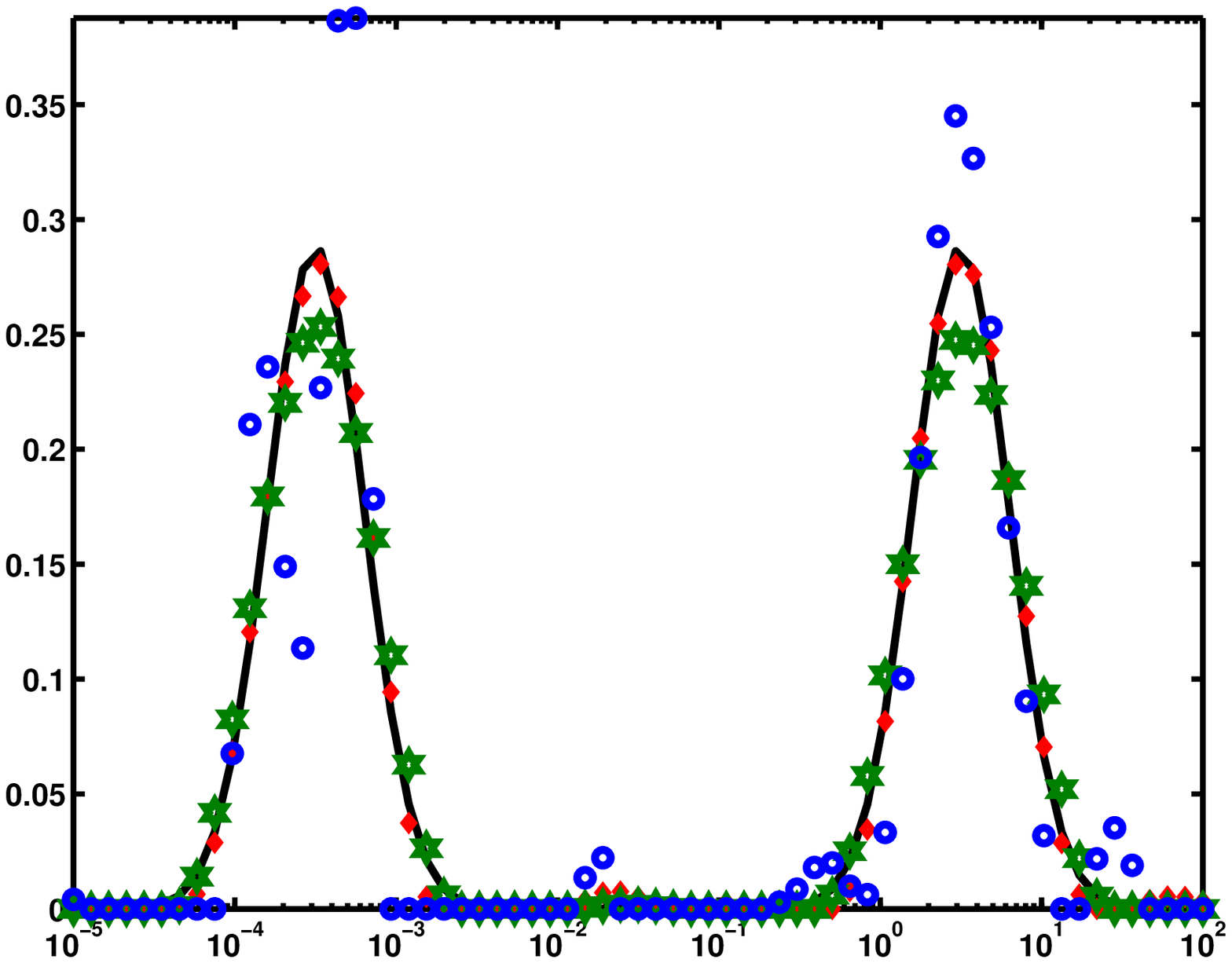}}
\caption{NNLS solutions of LN-A matrix $A_3$. Noise level $1\%$.}
\label{hnfig-lambdachoiceLN2A3HN}
\end{figure}

\begin{figure}[!ht]
\centering
\subfigure[$L=I$]{\includegraphics[width=1.7in]{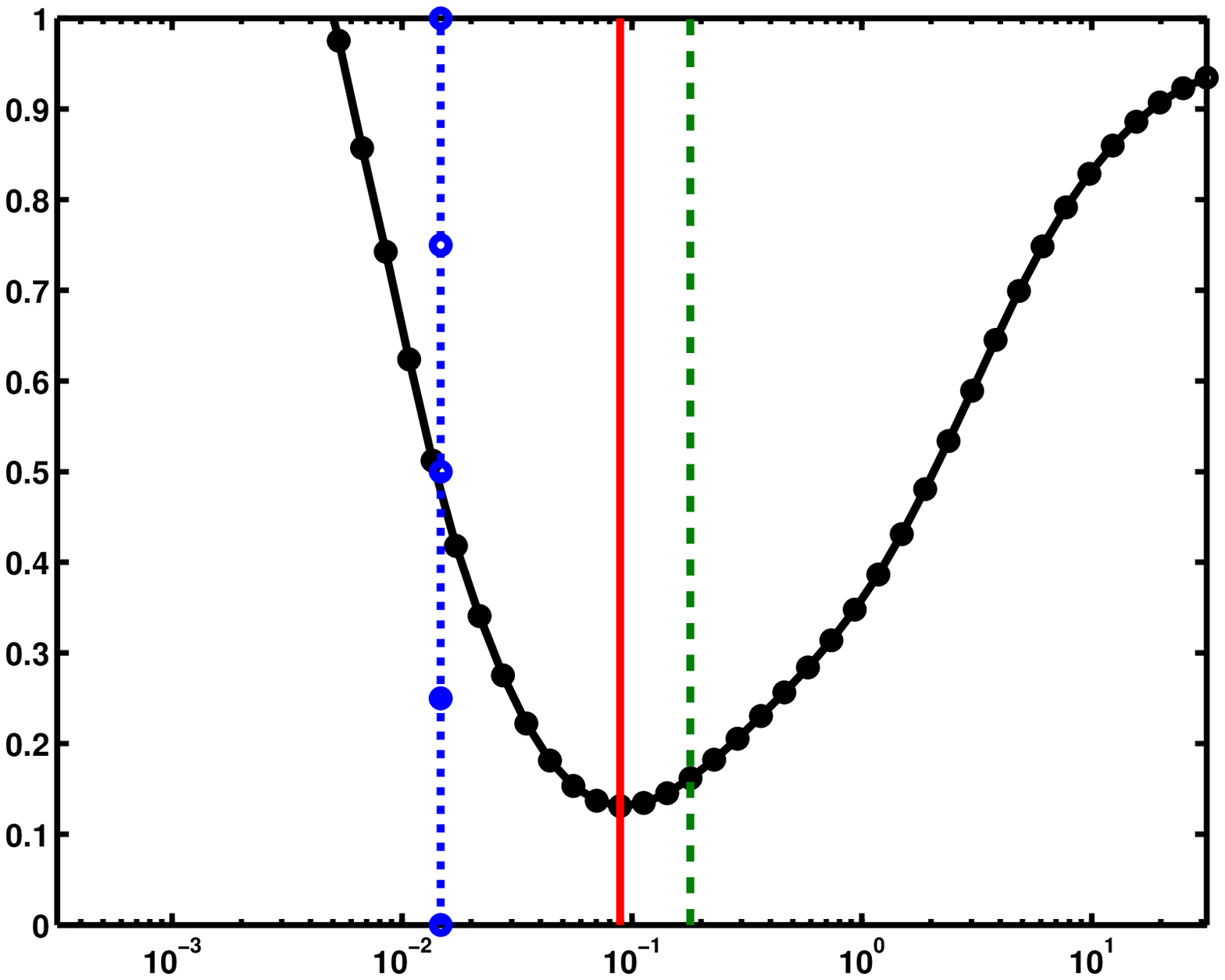}}
\subfigure[$L=L_1$]{\includegraphics[width=1.7in]{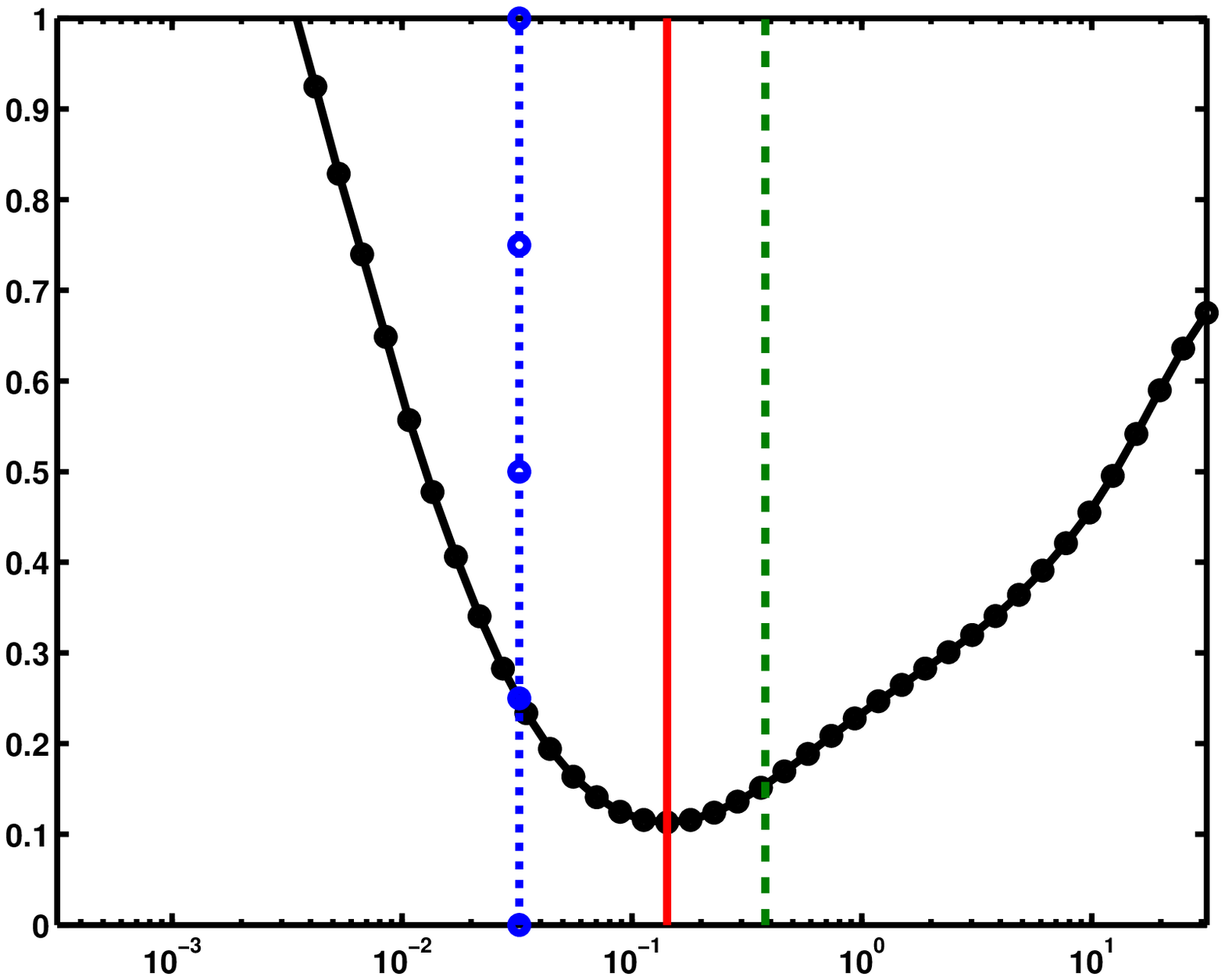}}
\subfigure[$L=L_2$]{\includegraphics[width=1.7in]{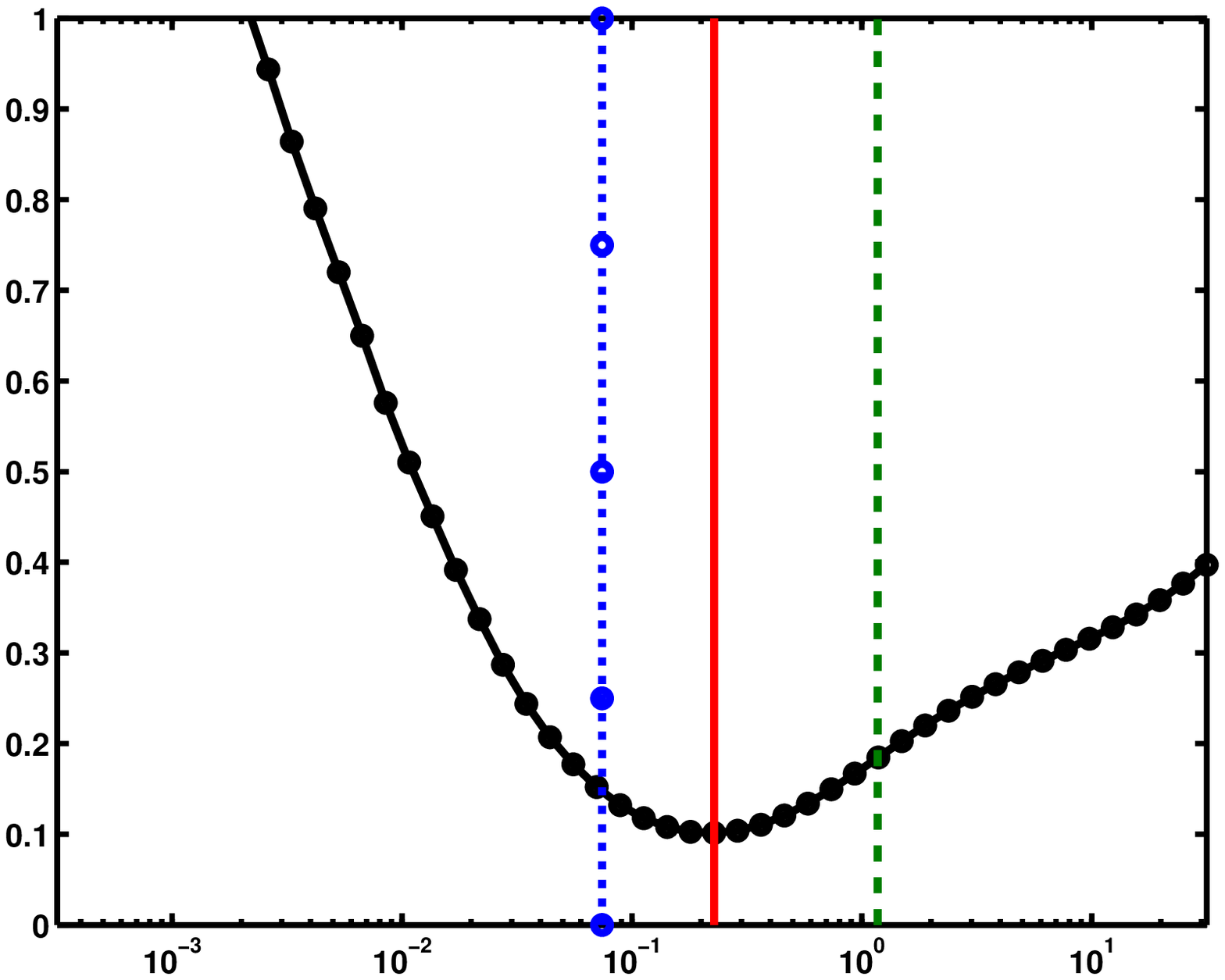}}
\subfigure[$L=I$]{\includegraphics[width=1.7in]{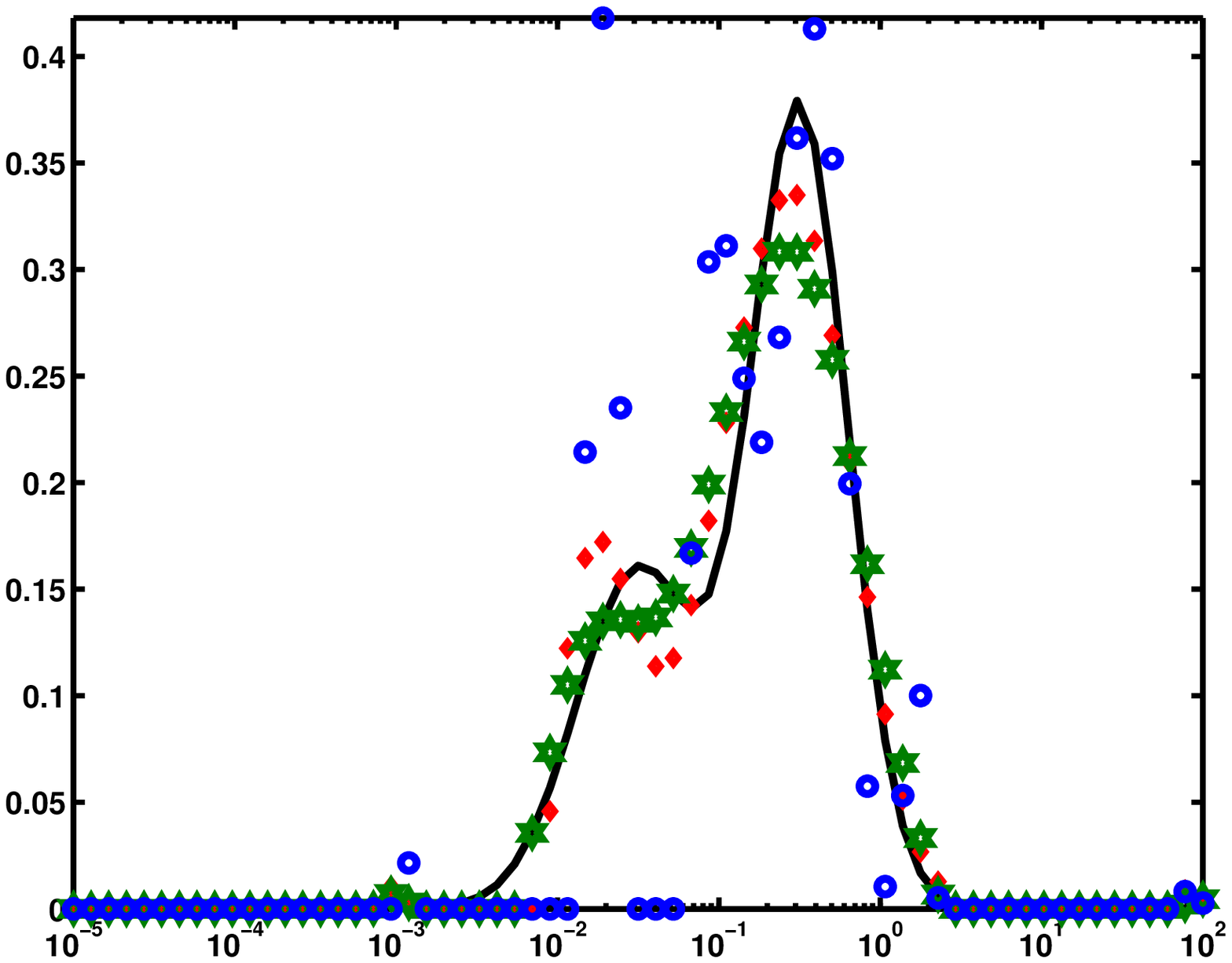}}
\subfigure[$L=L_1$]{\includegraphics[width=1.7in]{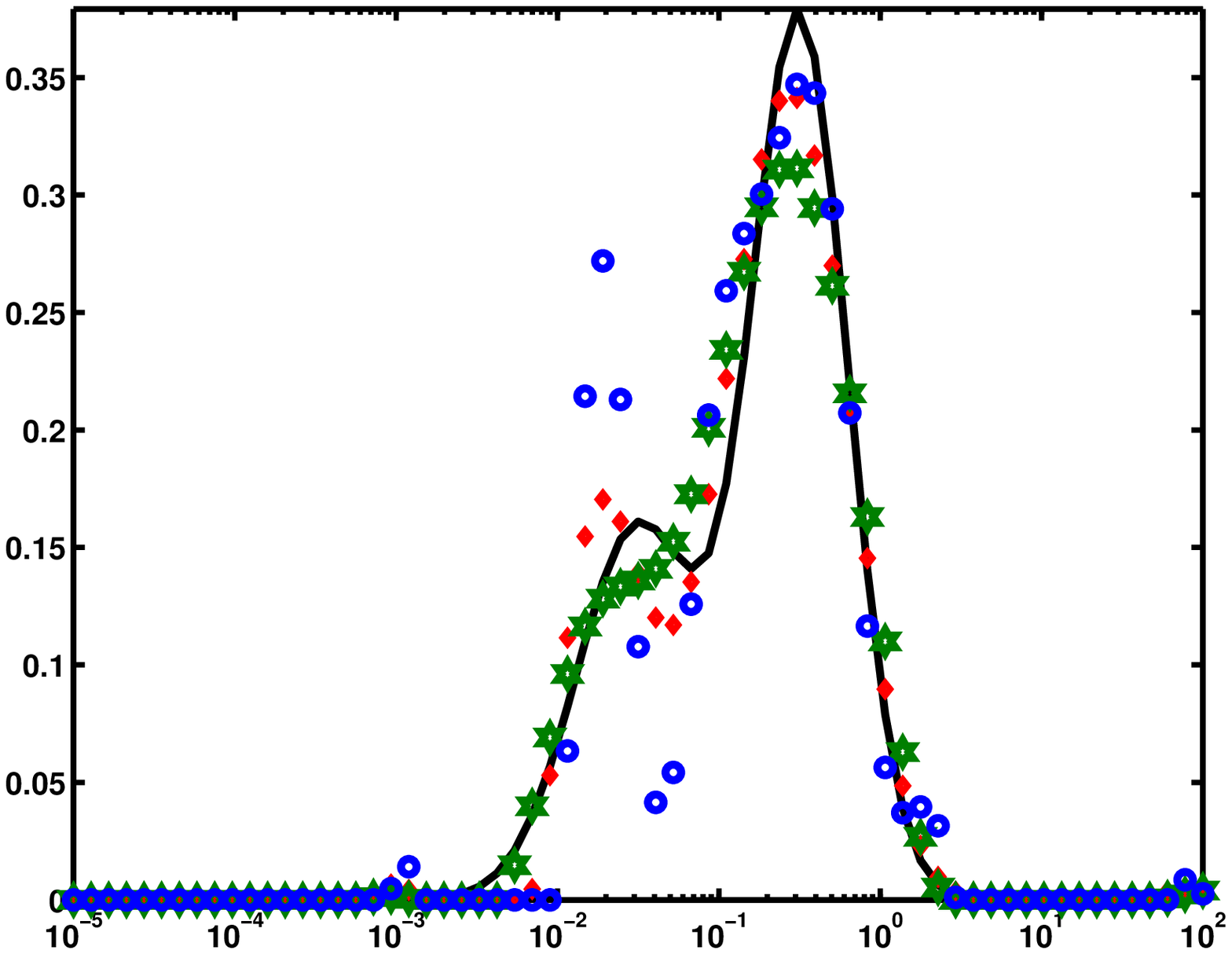}}
\subfigure[$L=L_2$]{\includegraphics[width=1.7in]{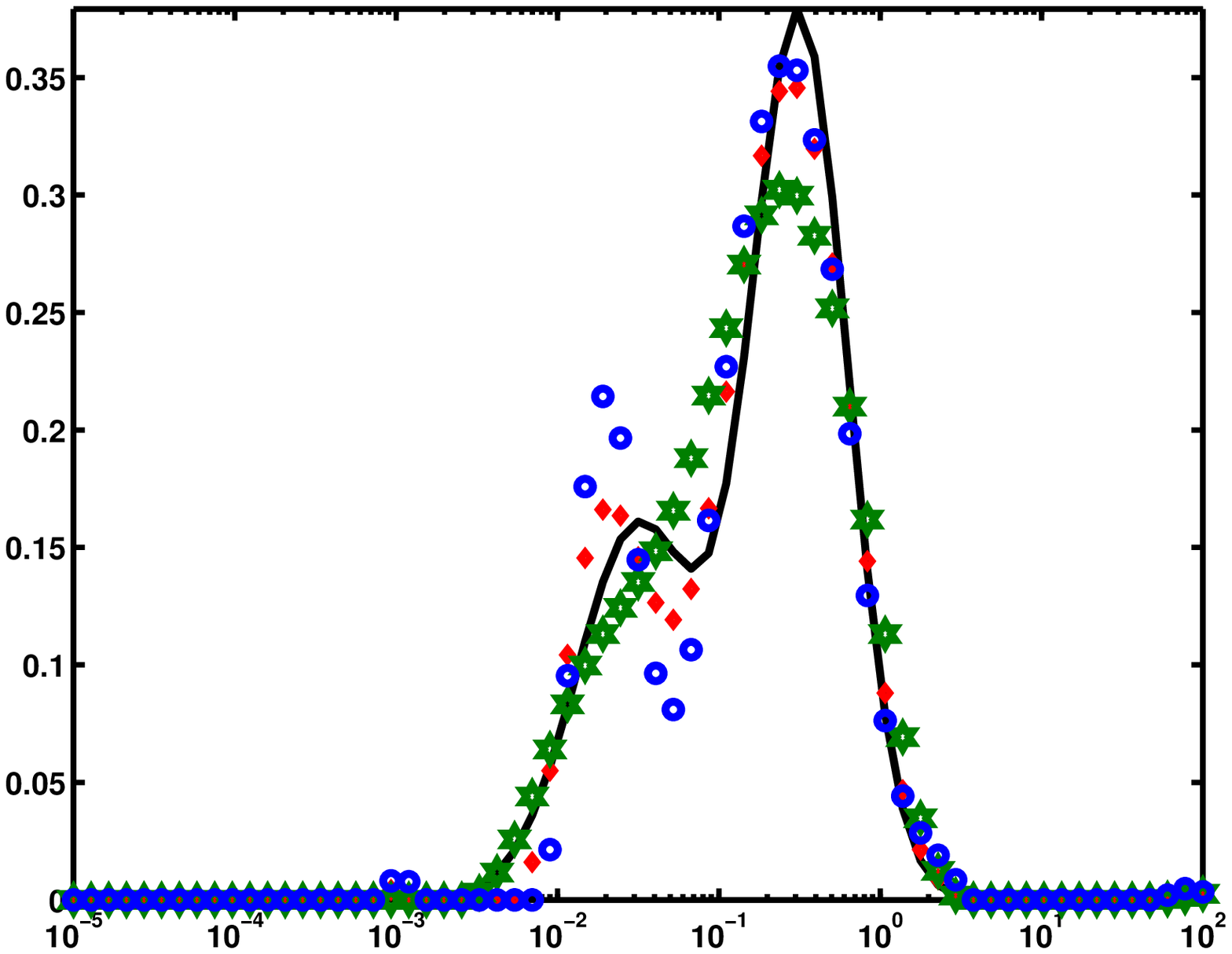}}
\caption{NNLS solutions of LN-B matrix $A_3$. Noise level $1\%$.}
\label{hnfig-lambdachoiceLN5A3HN}
\end{figure}

\begin{figure}[!ht]
\centering
\subfigure[$L=I$]{\includegraphics[width=1.7in]{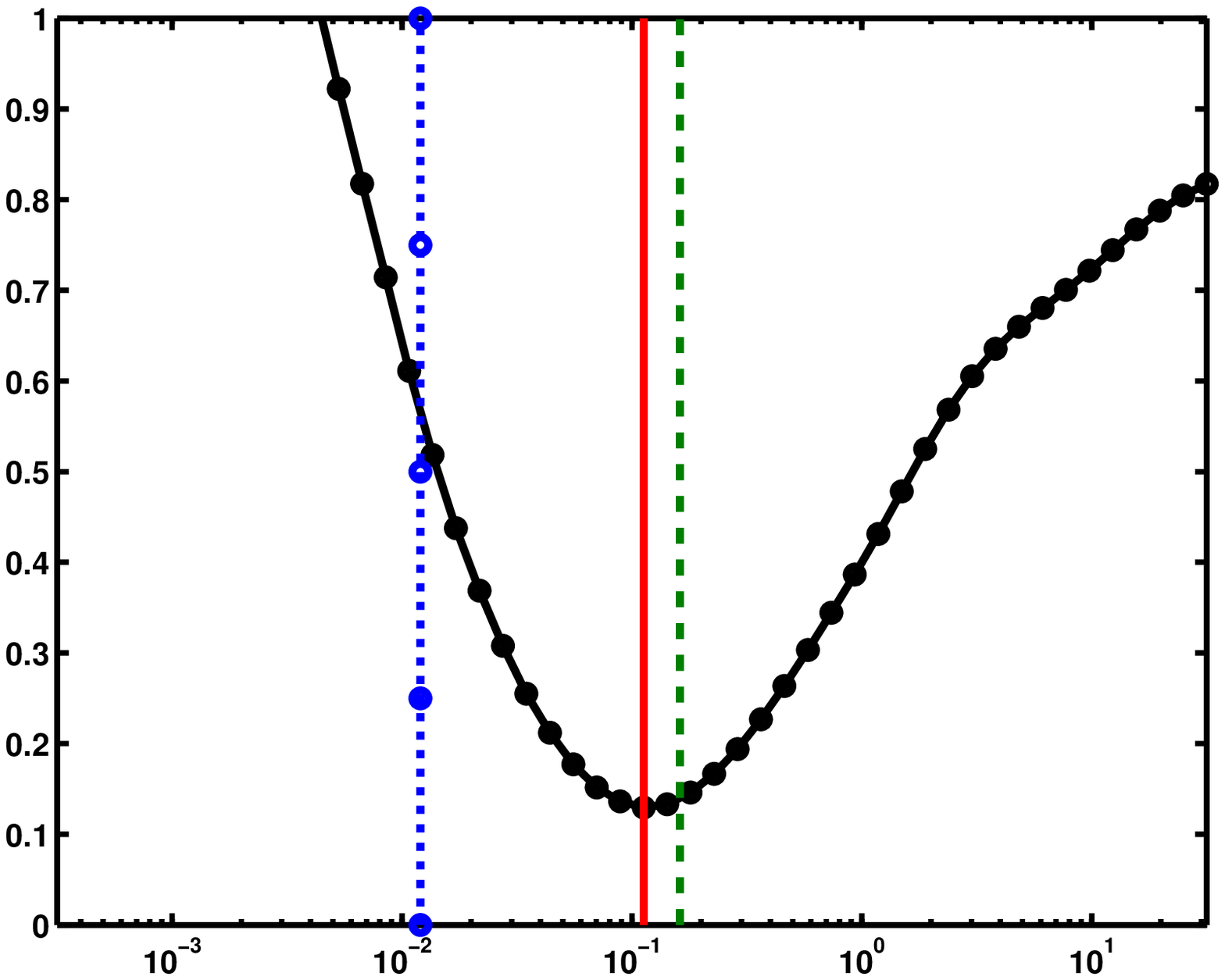}}
\subfigure[$L=L_1$]{\includegraphics[width=1.7in]{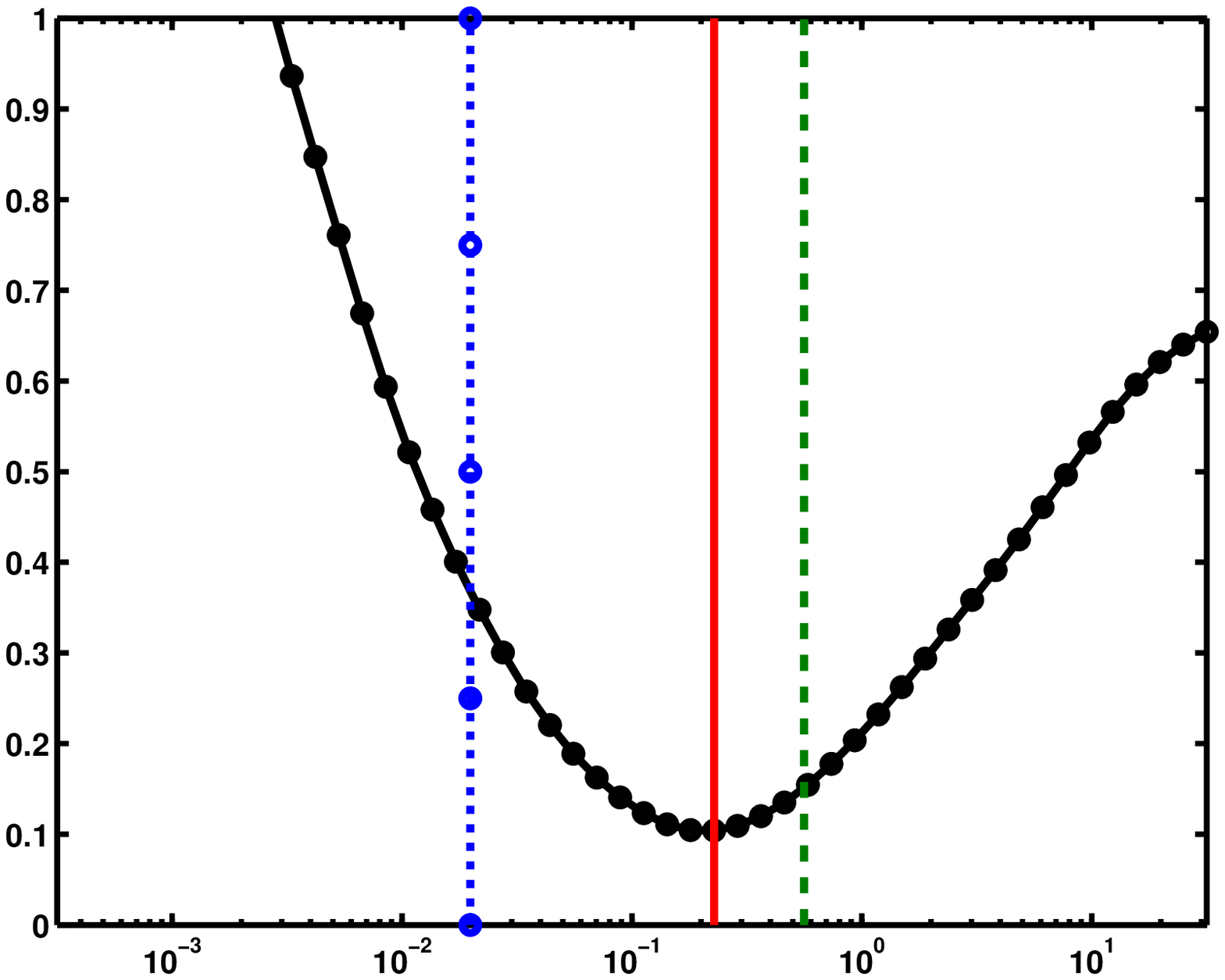}}
\subfigure[$L=L_2$]{\includegraphics[width=1.7in]{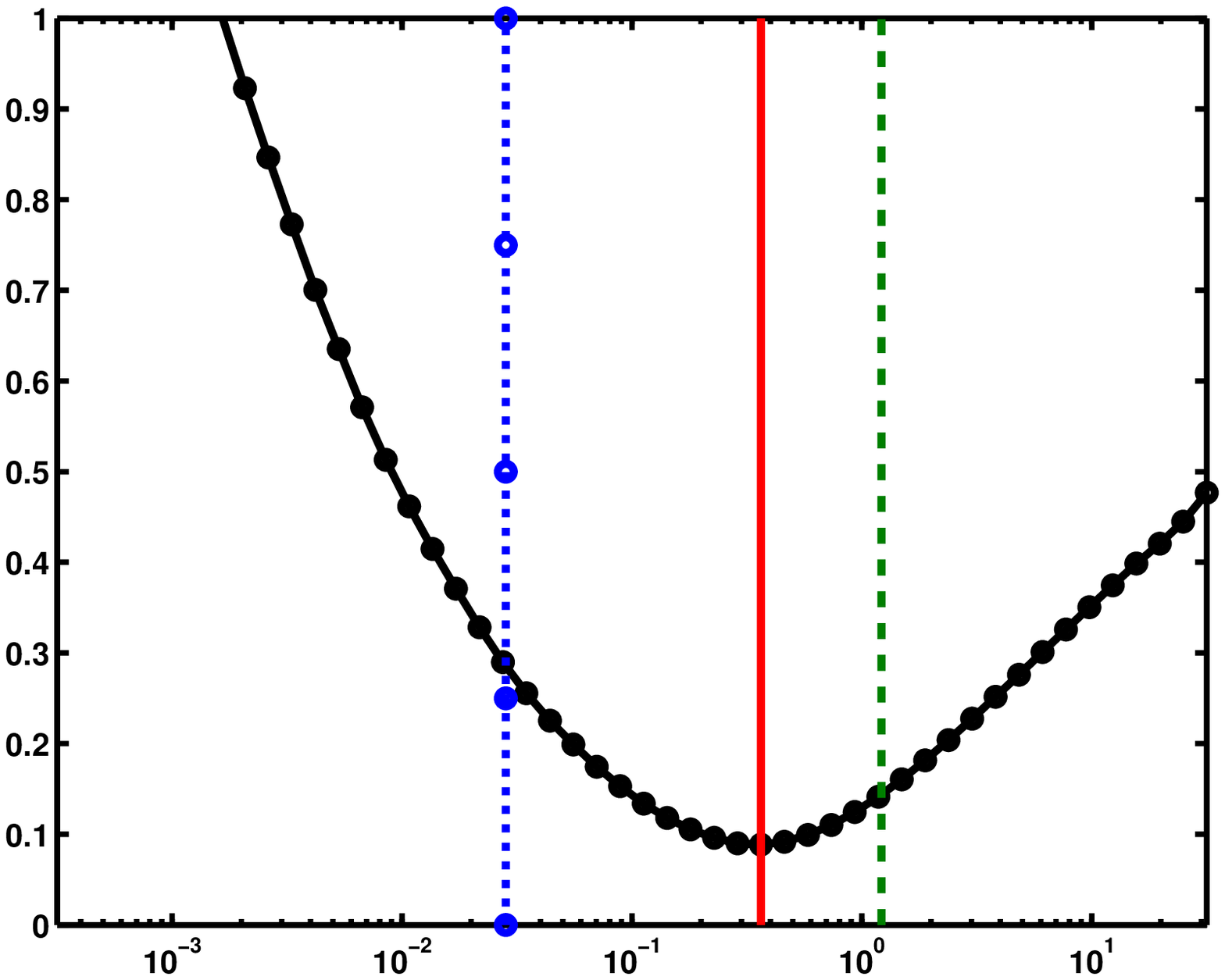}}
\subfigure[$L=I$]{\includegraphics[width=1.7in]{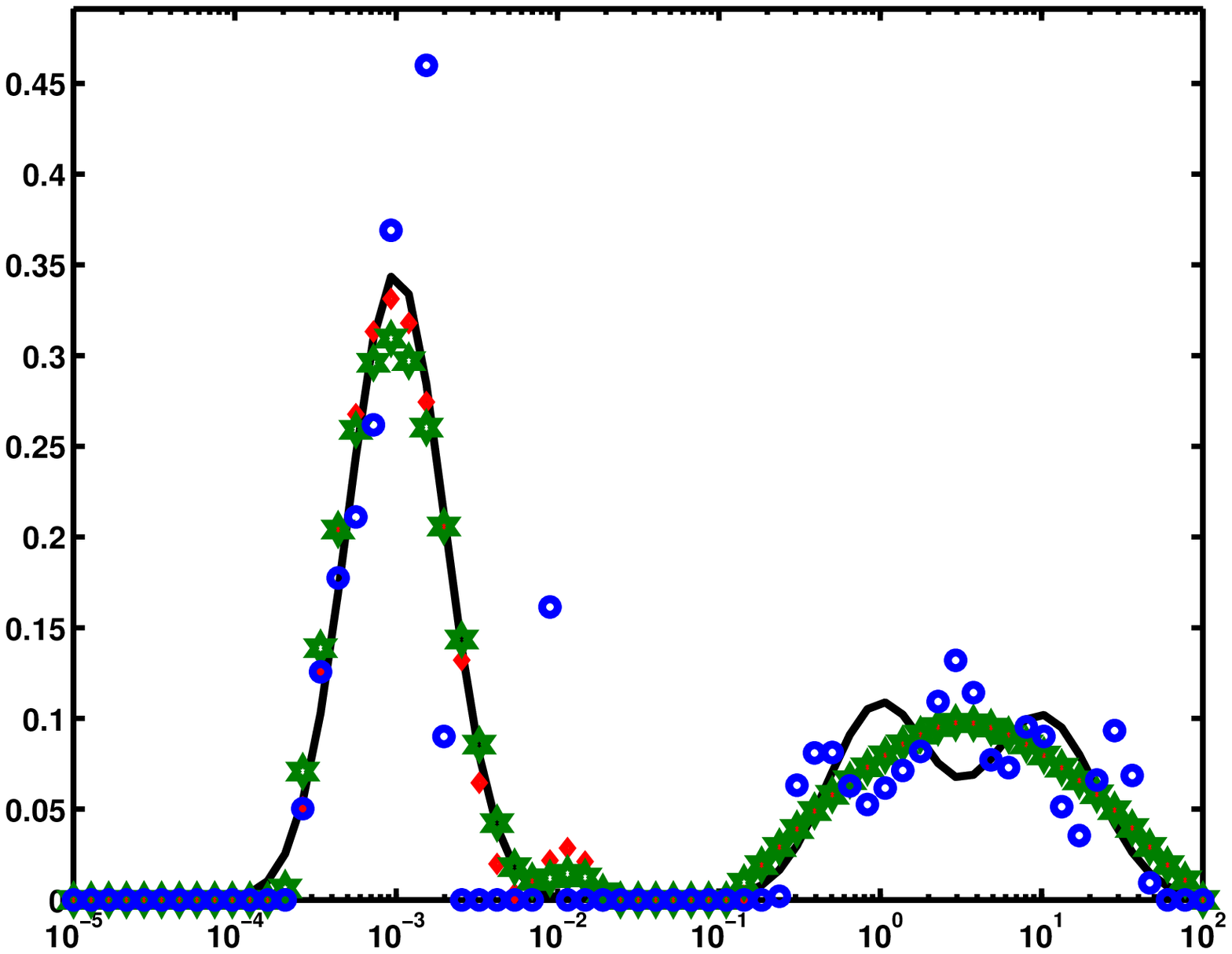}}
\subfigure[$L=L_1$]{\includegraphics[width=1.7in]{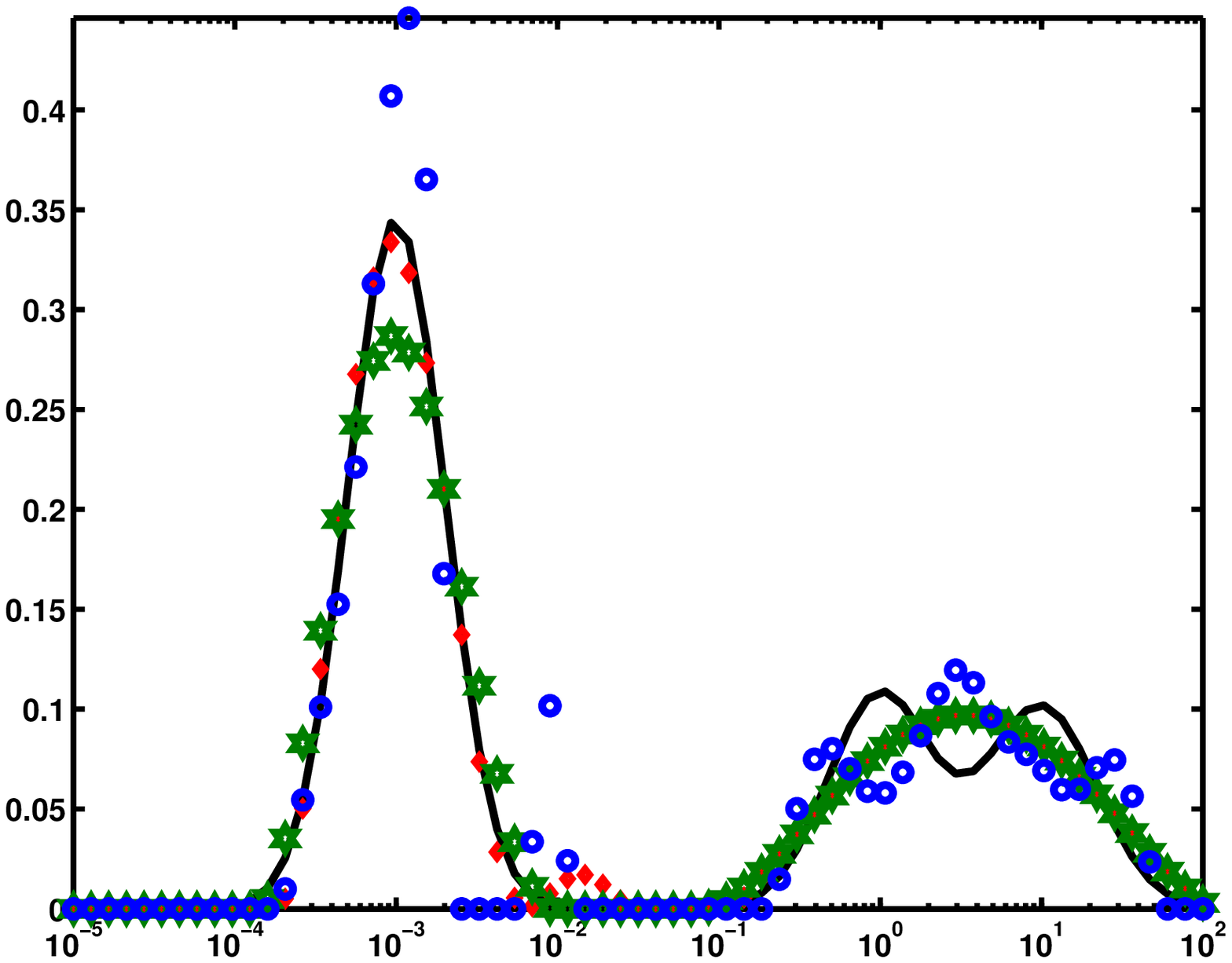}}
\subfigure[$L=L_2$]{\includegraphics[width=1.7in]{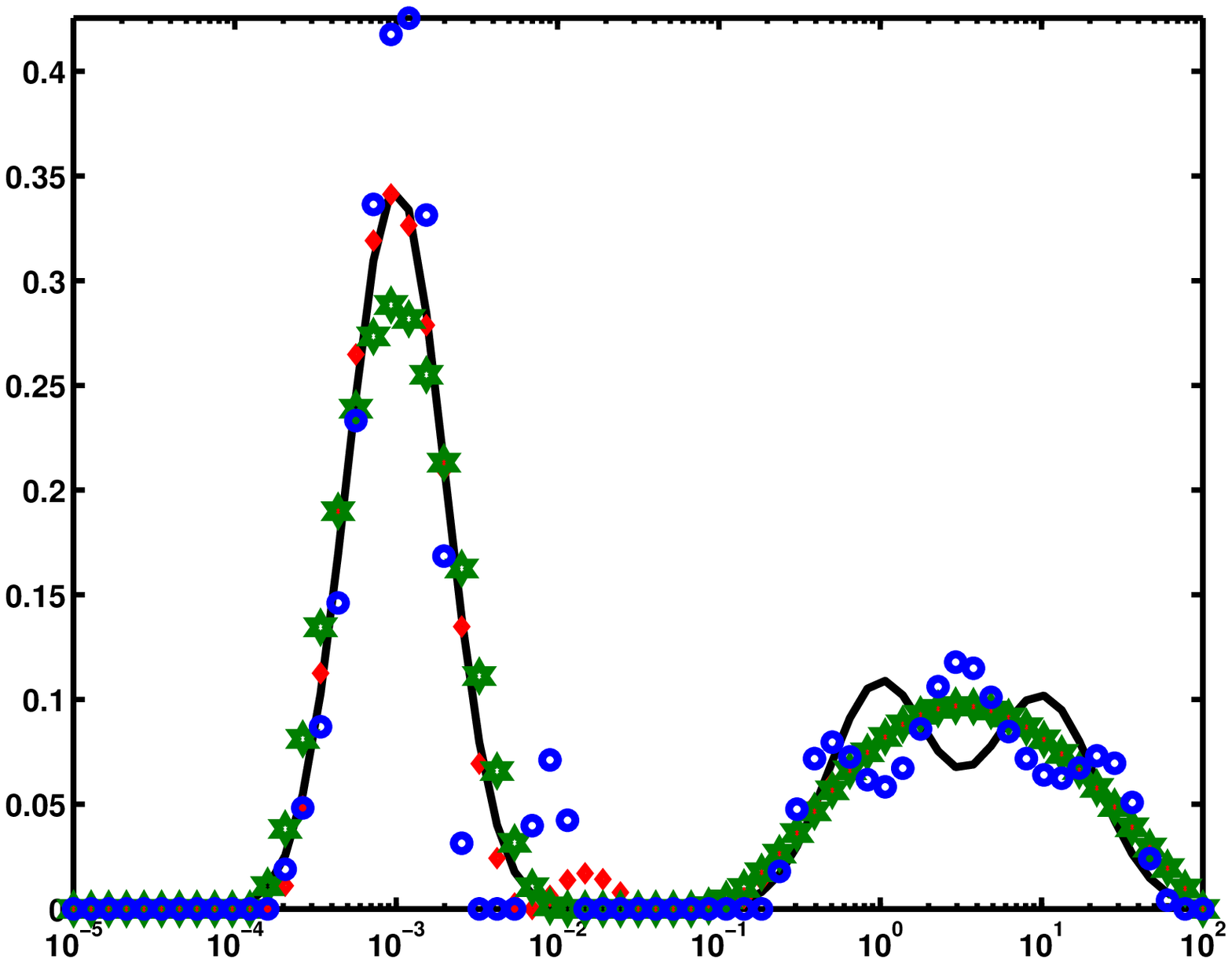}}
\caption{NNLS solutions of LN-C matrix $A_3$. Noise level $1\%$.}
\label{hnfig-lambdachoiceLN6A3HN}
\end{figure}
\clearpage
\subsection{Examples: Noise level $5\%$ $A_4$ NNLS}
 \begin{figure}[!ht]
\centering
\subfigure[$L=I$]{\includegraphics[width=1.7in]{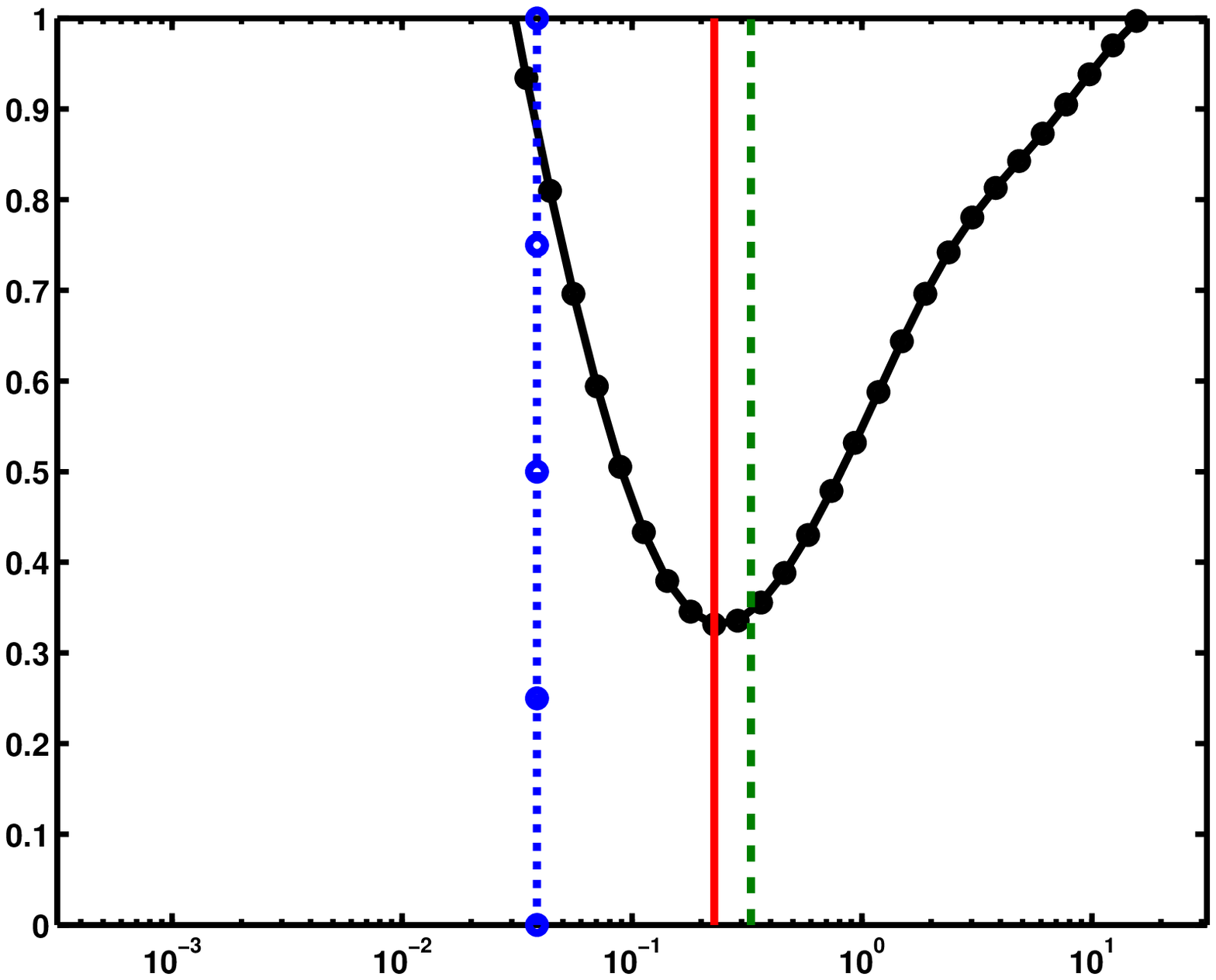}}
\subfigure[$L=L_1$]{\includegraphics[width=1.7in]{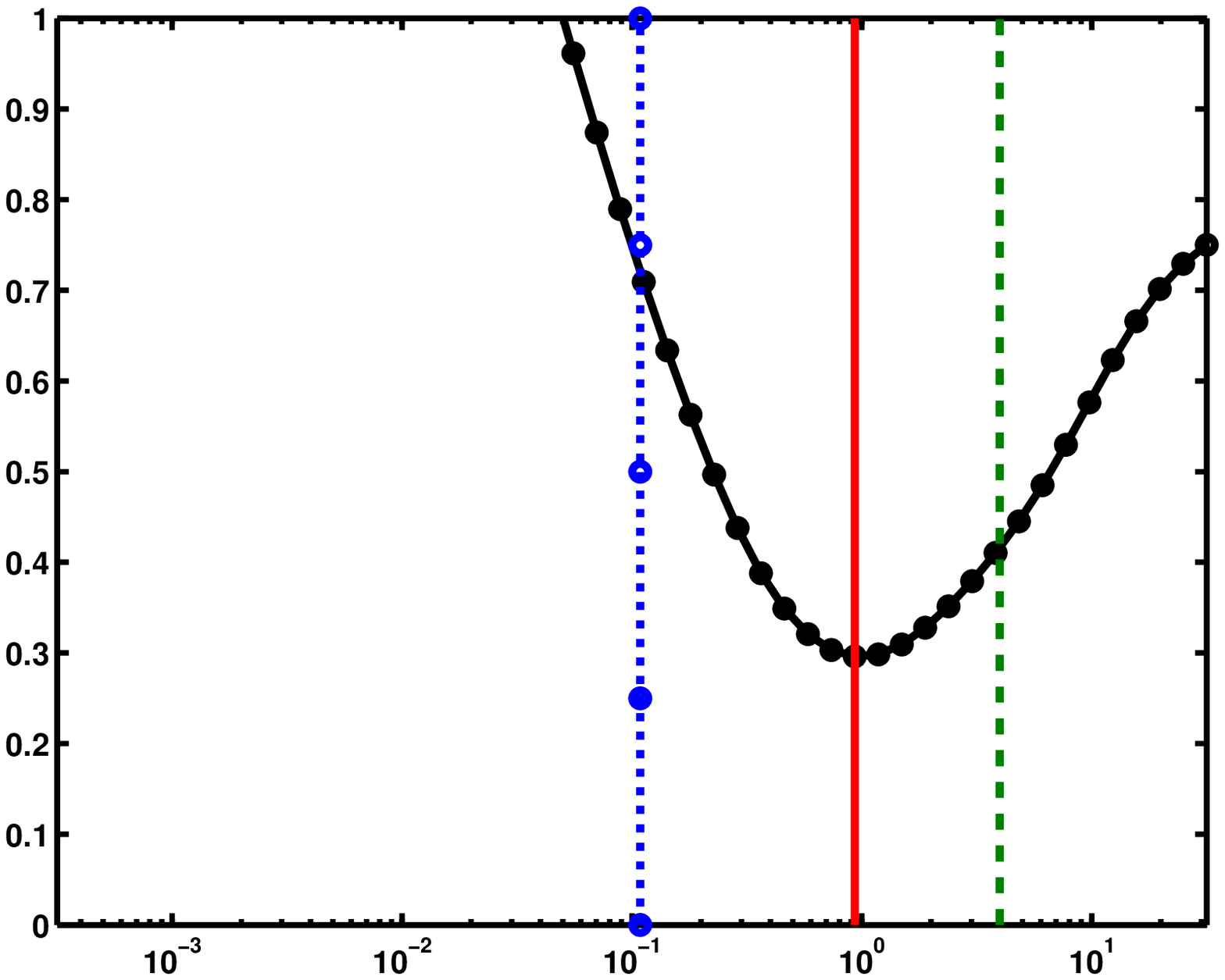}}
\subfigure[$L=L_2$]{\includegraphics[width=1.7in]{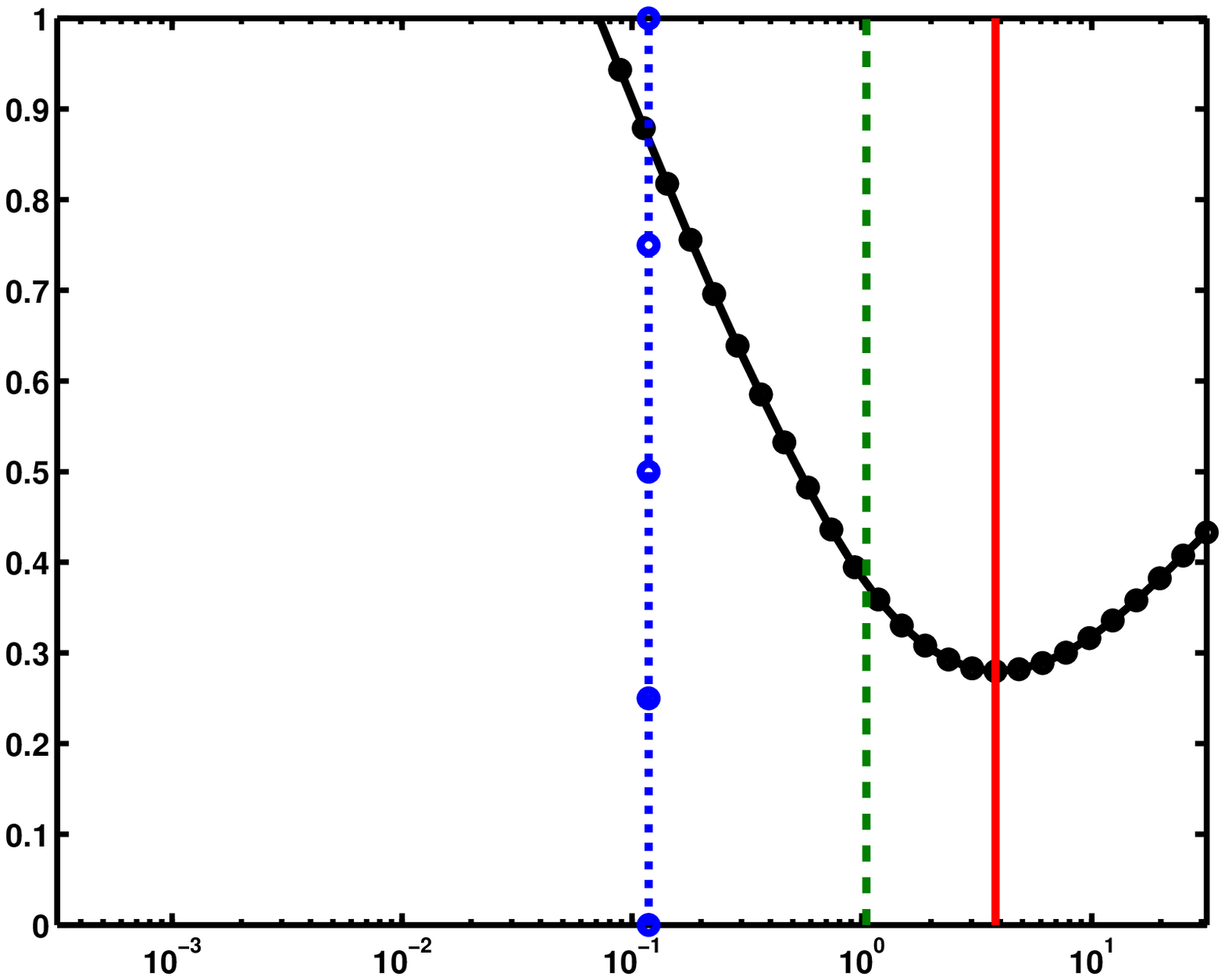}}
\subfigure[$L=I$]{\includegraphics[width=1.7in]{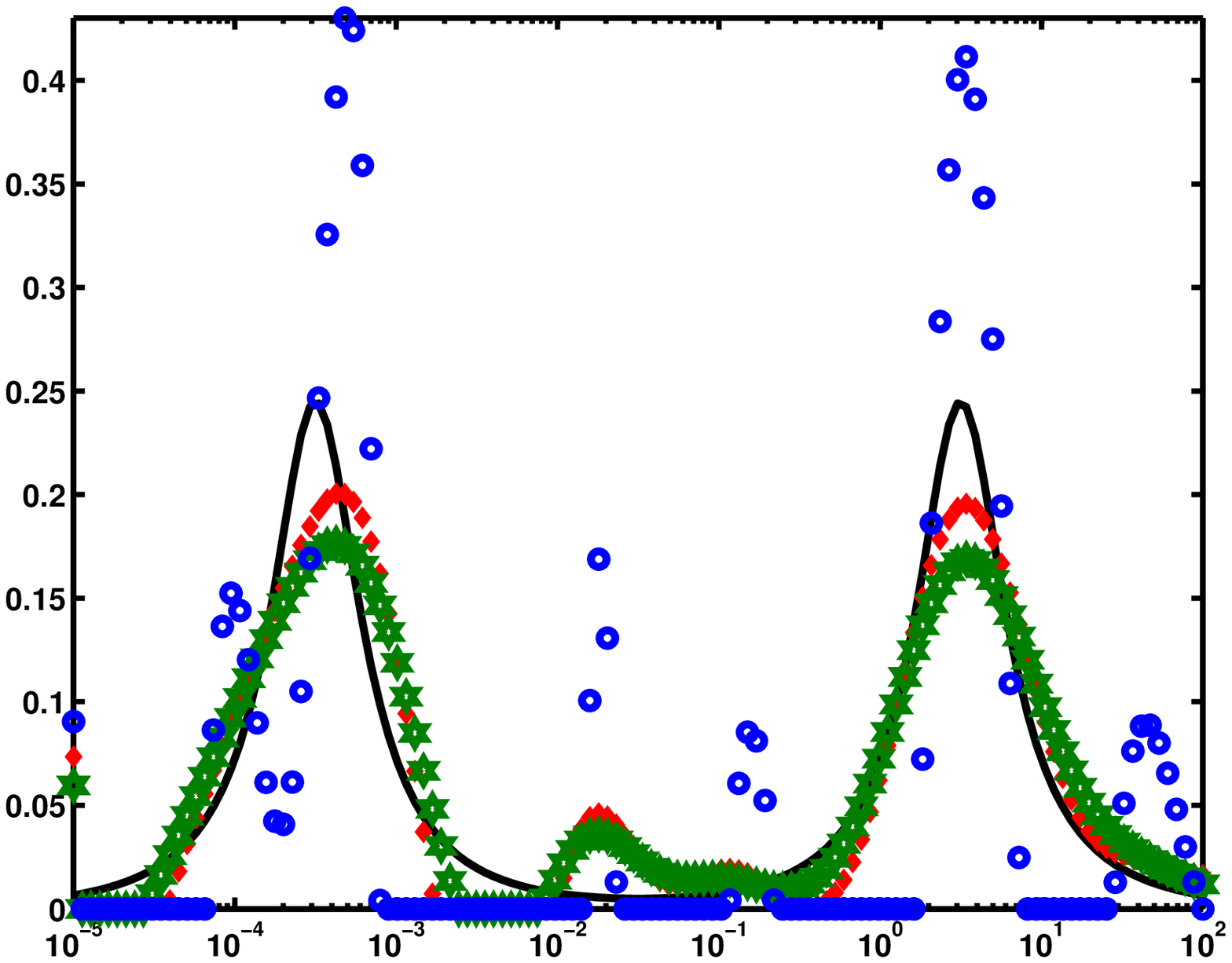}}
\subfigure[$L=L_1$]{\includegraphics[width=1.7in]{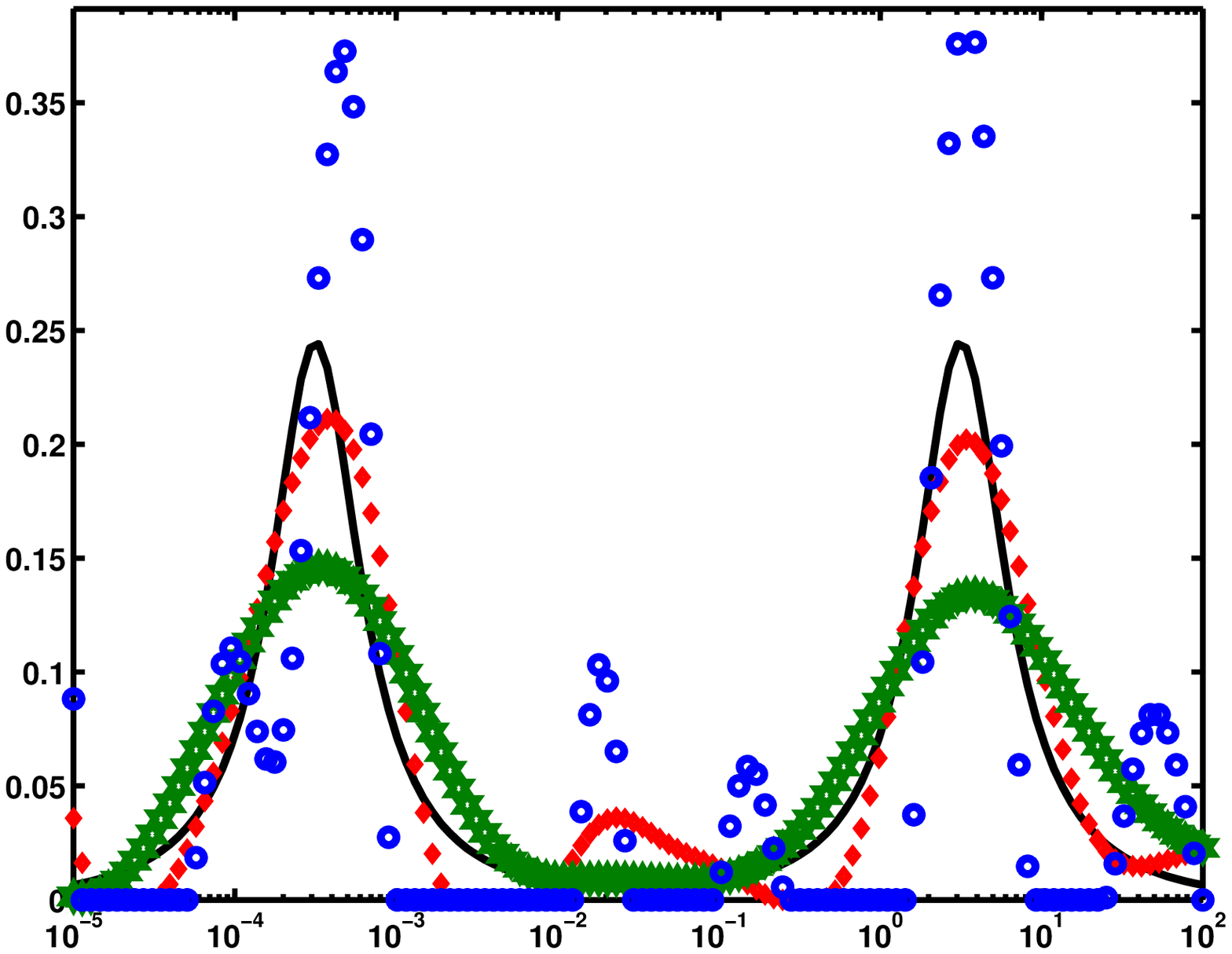}}
\subfigure[$L=L_2$]{\includegraphics[width=1.7in]{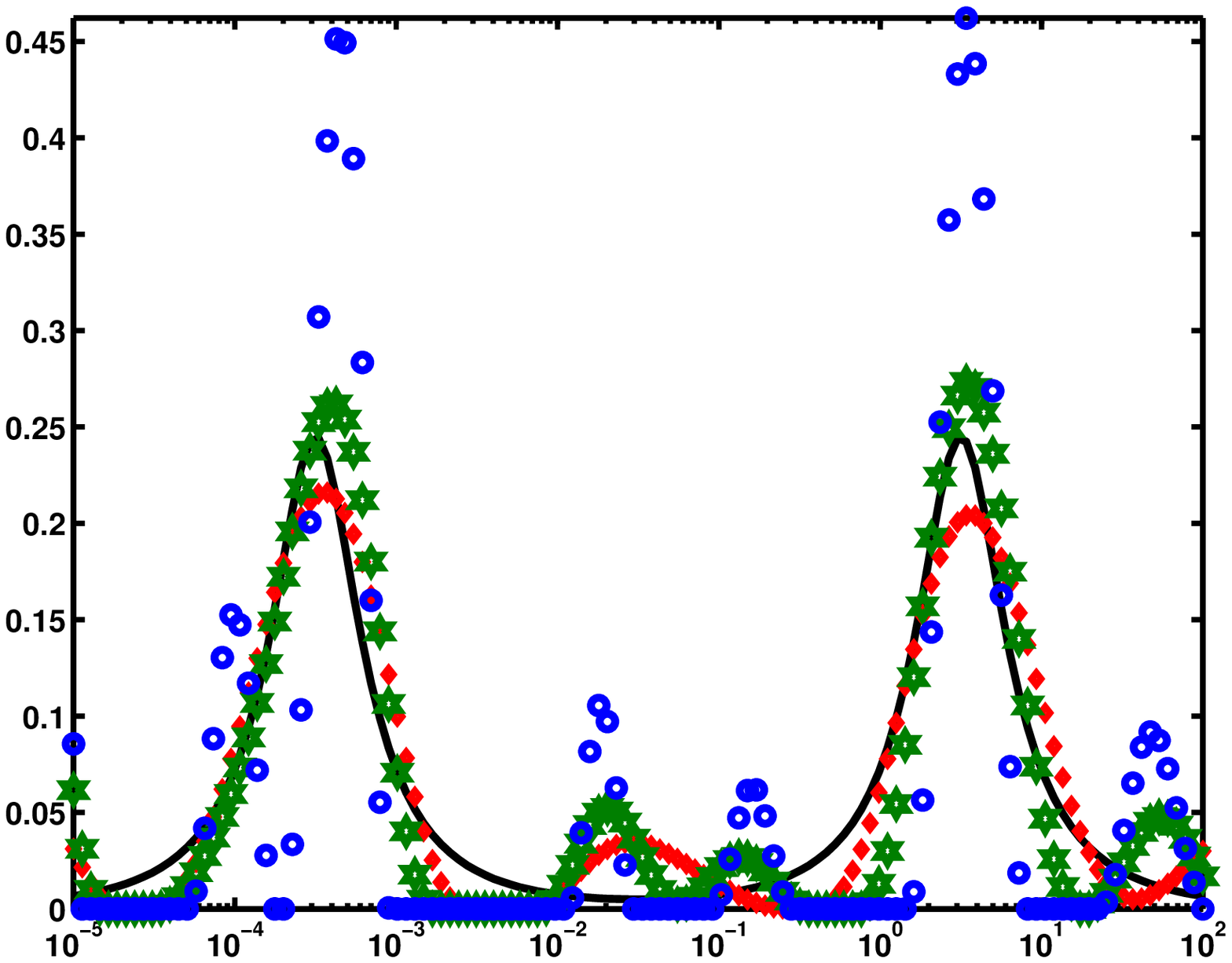}}
\caption{NNLS solutions of RQ-A matrix $A_4$. Noise level $5\%$.}
\label{hnfig-lambdachoiceRQ1A4HN}
\end{figure}
 \begin{figure}[!ht]
\centering
\subfigure[$L=I$]{\includegraphics[width=1.7in]{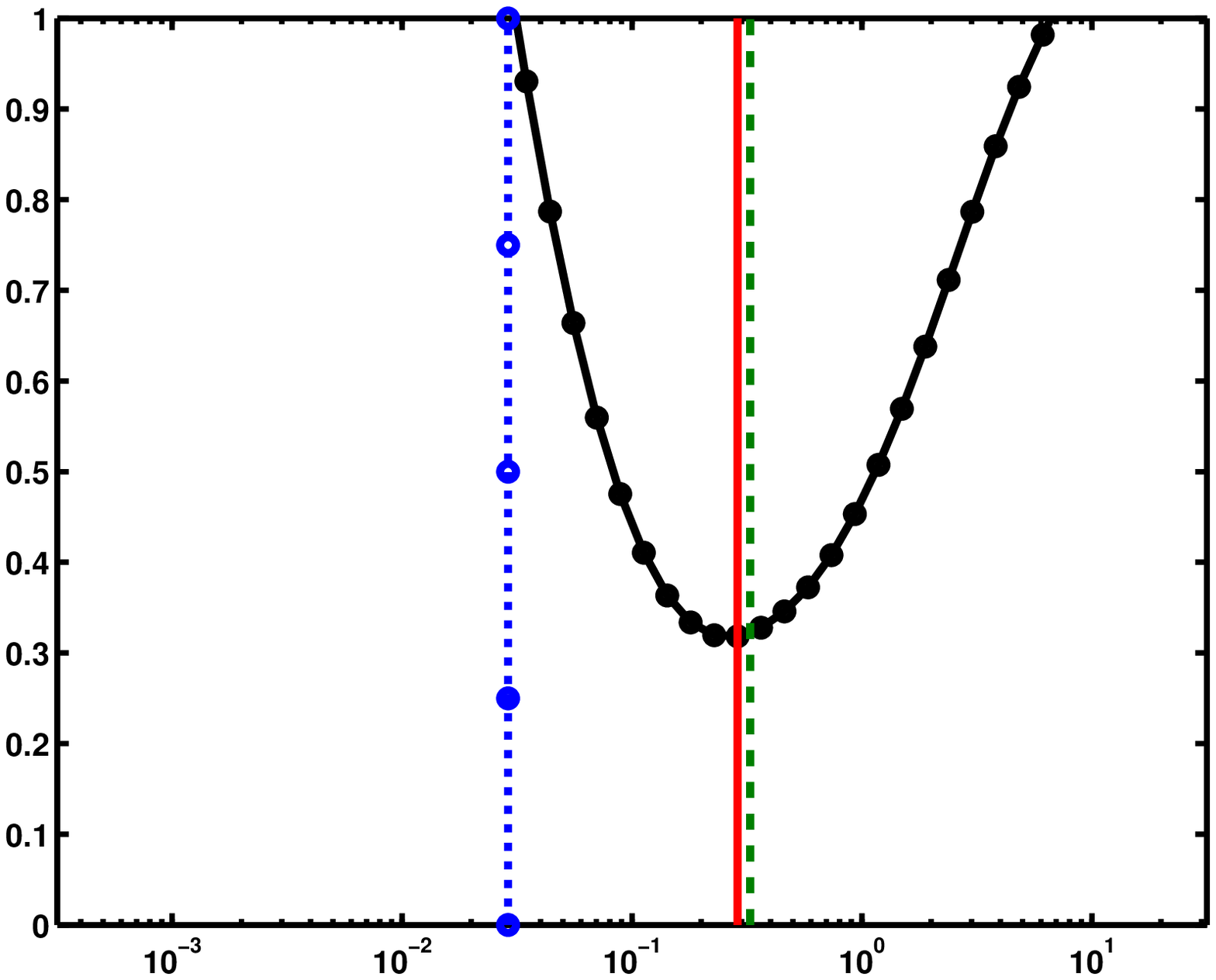}}
\subfigure[$L=L_1$]{\includegraphics[width=1.7in]{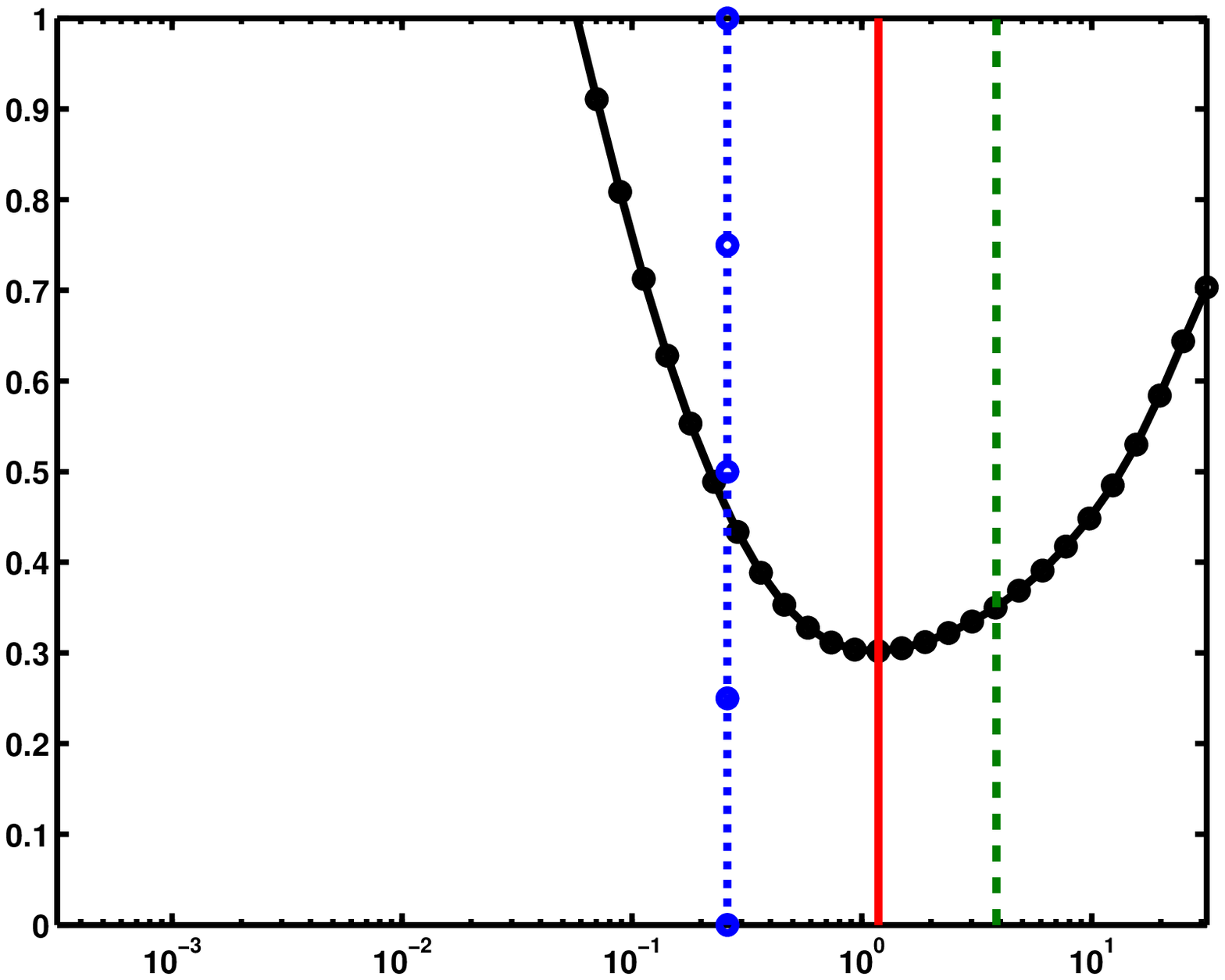}}
\subfigure[$L=L_2$]{\includegraphics[width=1.7in]{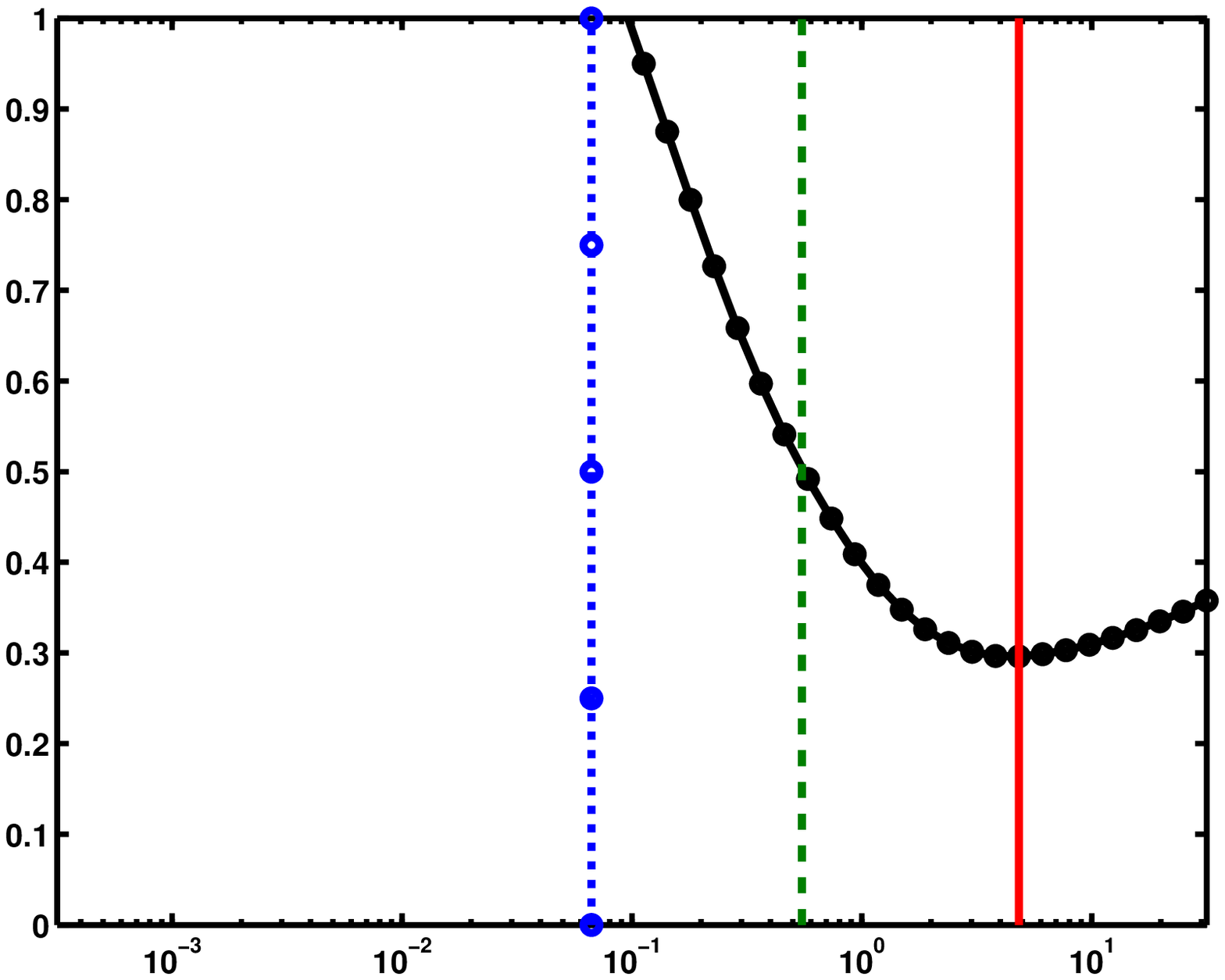}}
\subfigure[$L=I$]{\includegraphics[width=1.7in]{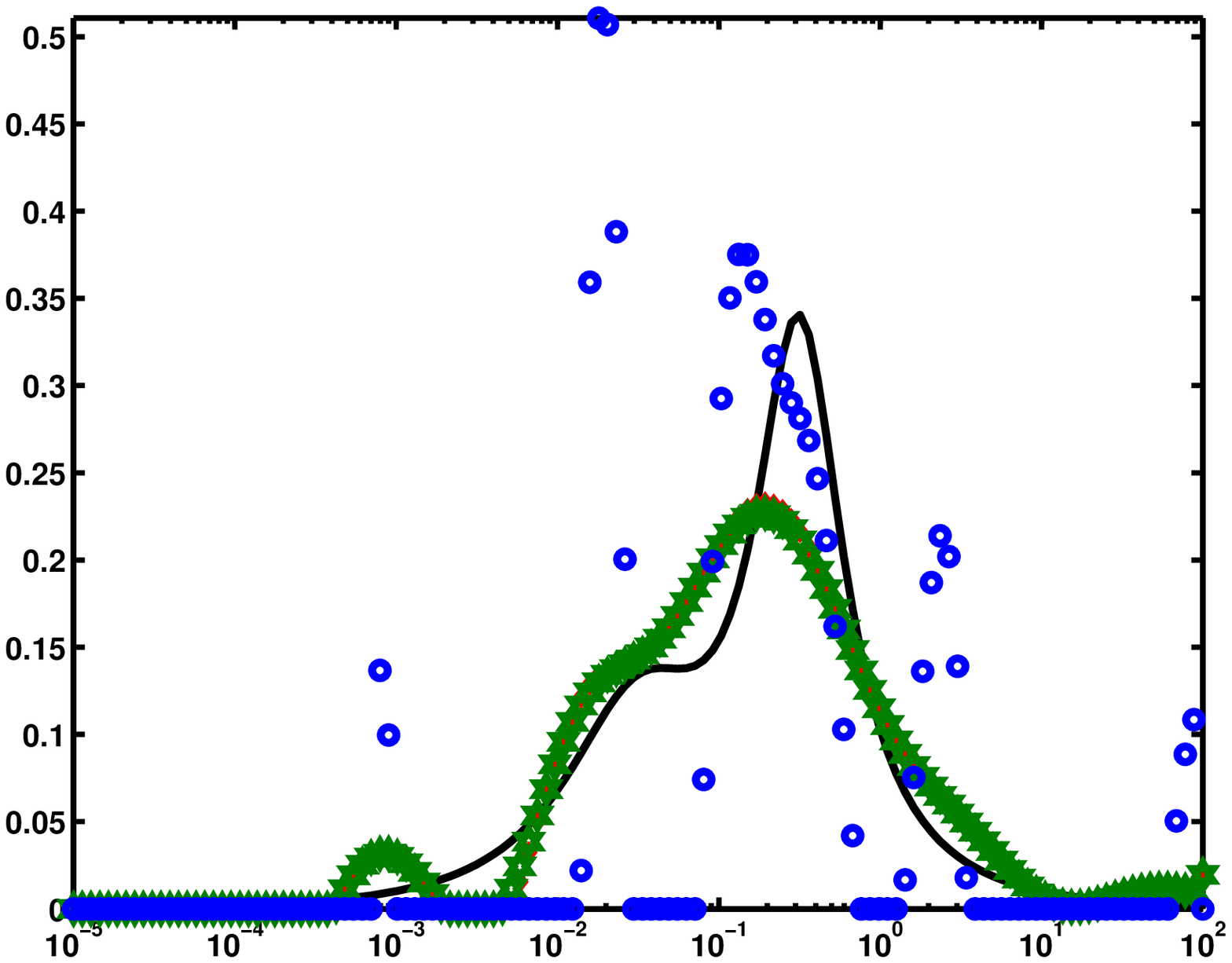}}
\subfigure[$L=L_1$]{\includegraphics[width=1.7in]{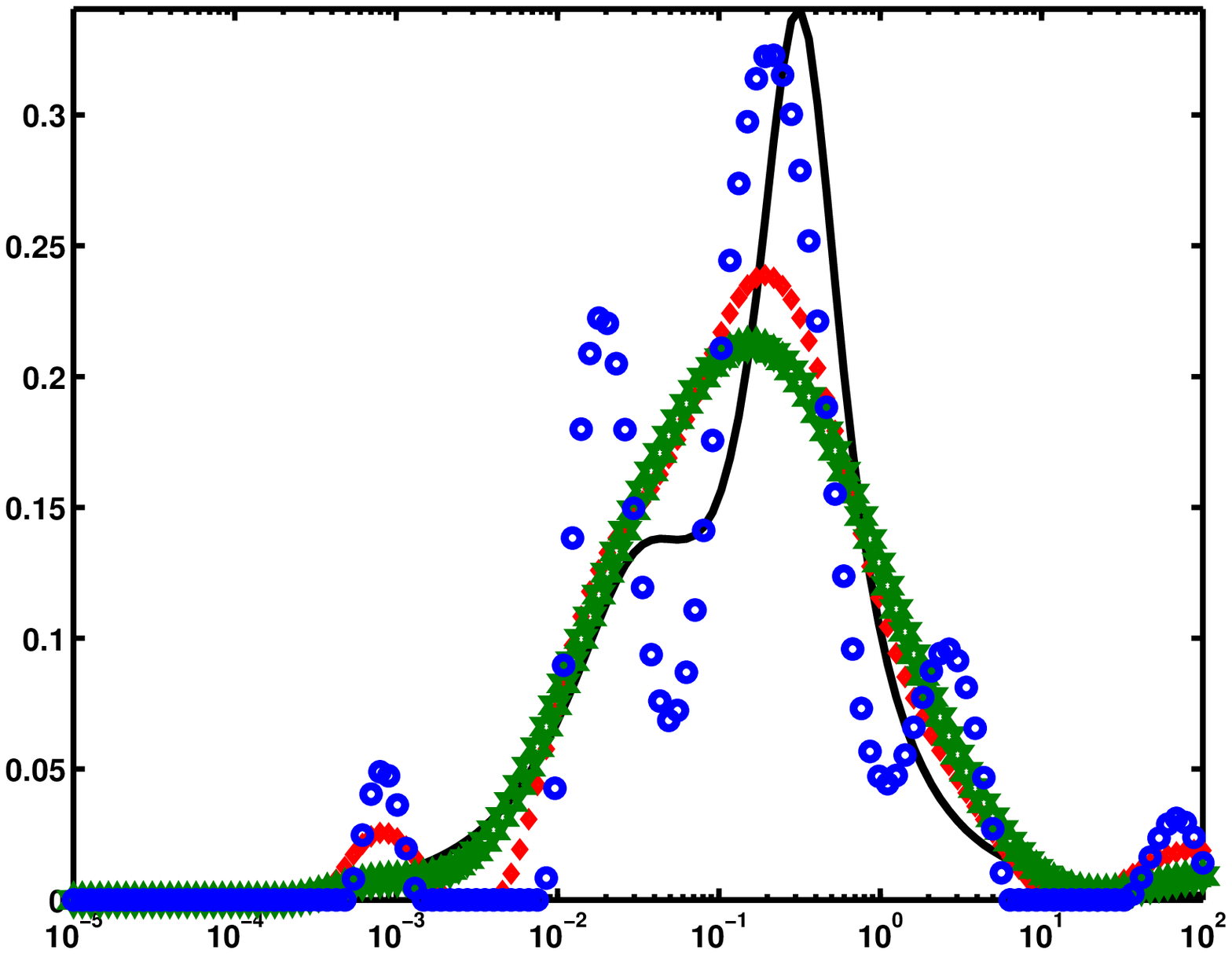}}
\subfigure[$L=L_2$]{\includegraphics[width=1.7in]{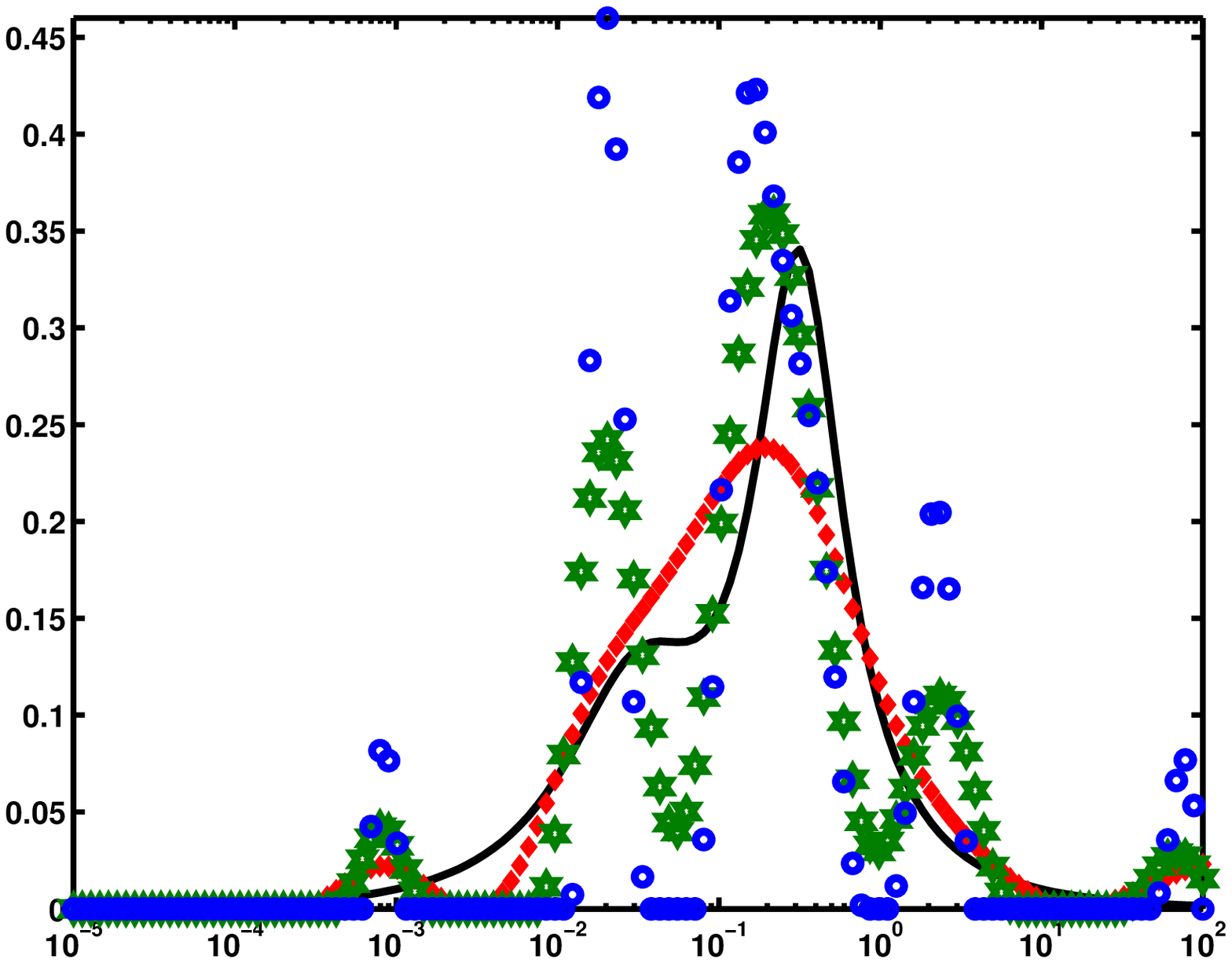}}
\caption{NNLS solutions of RQ-B matrix $A_4$. Noise level $5\%$.}
\label{hnfig-lambdachoiceRQ5A4HN}
\end{figure}
 \begin{figure}[!ht]
\centering
\subfigure[$L=I$]{\includegraphics[width=1.7in]{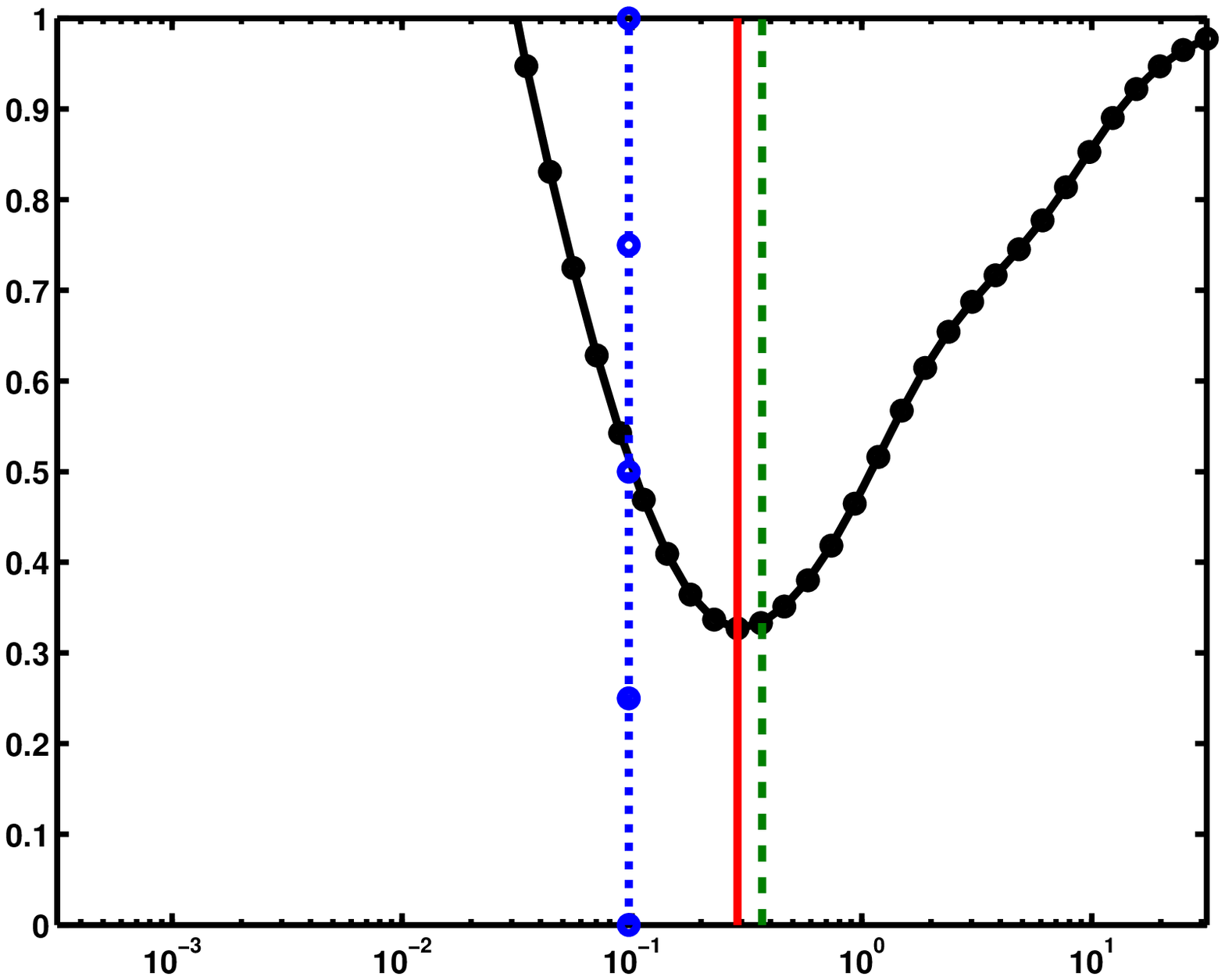}}
\subfigure[$L=L_1$]{\includegraphics[width=1.7in]{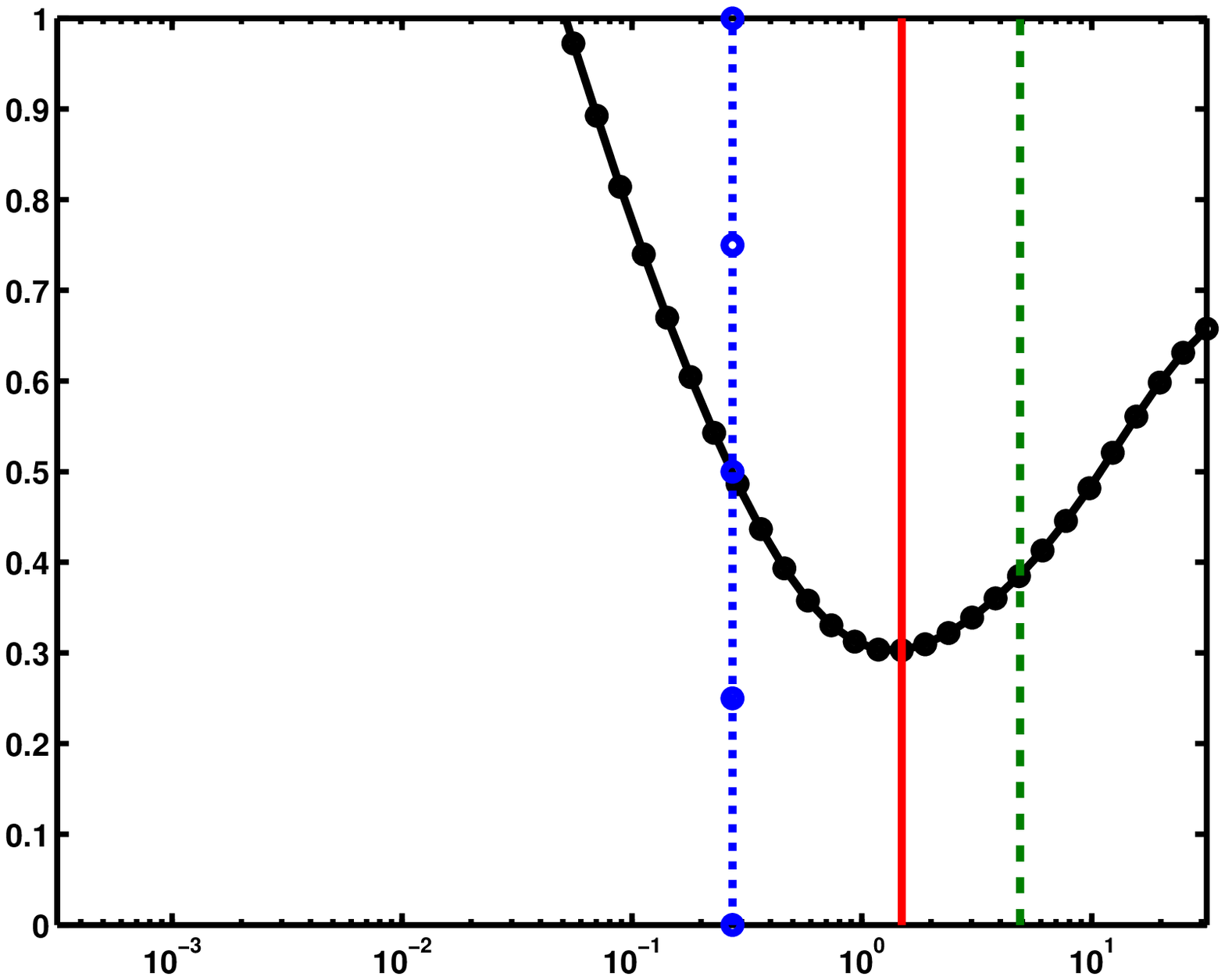}}
\subfigure[$L=L_2$]{\includegraphics[width=1.7in]{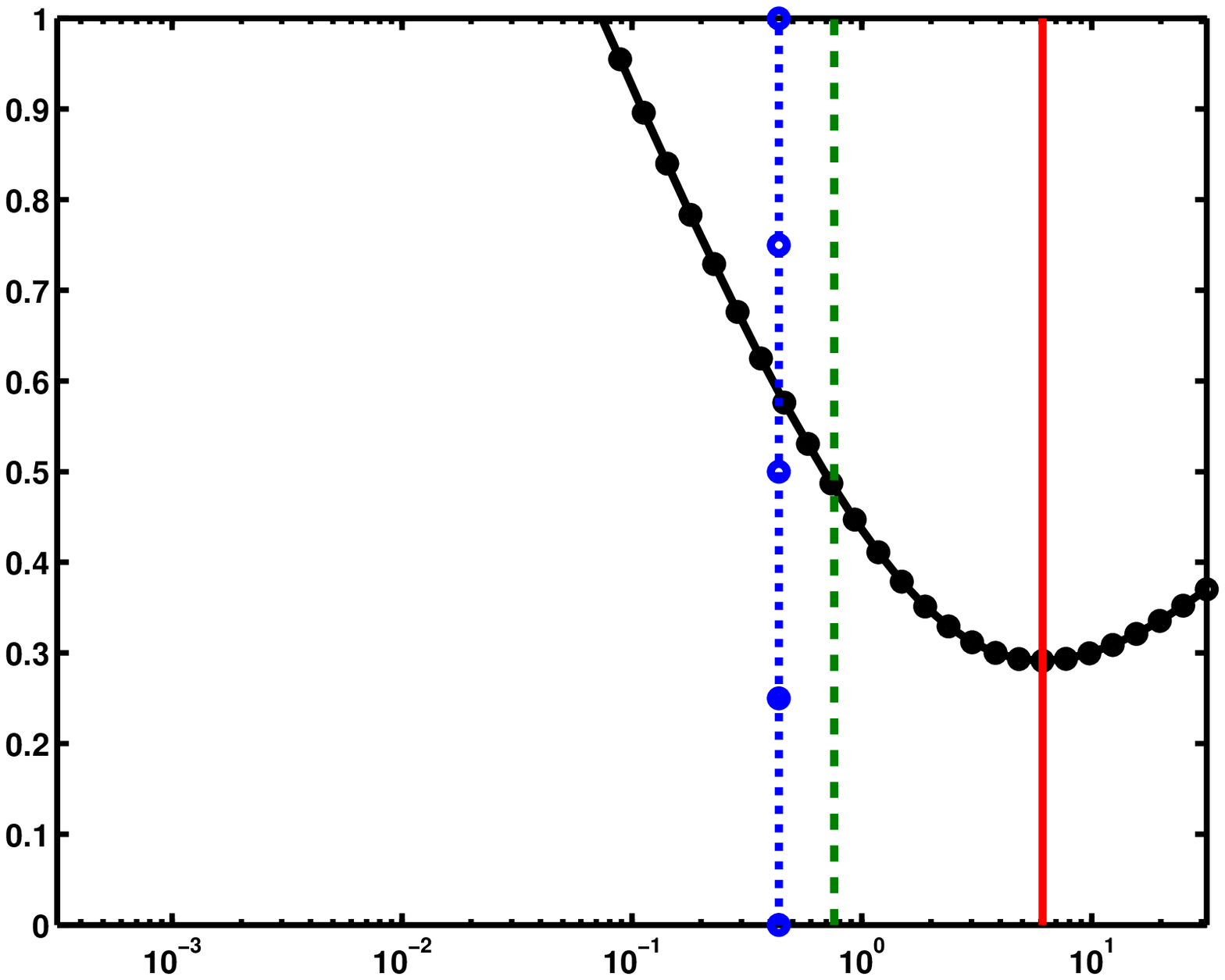}}
\subfigure[$L=I$]{\includegraphics[width=1.7in]{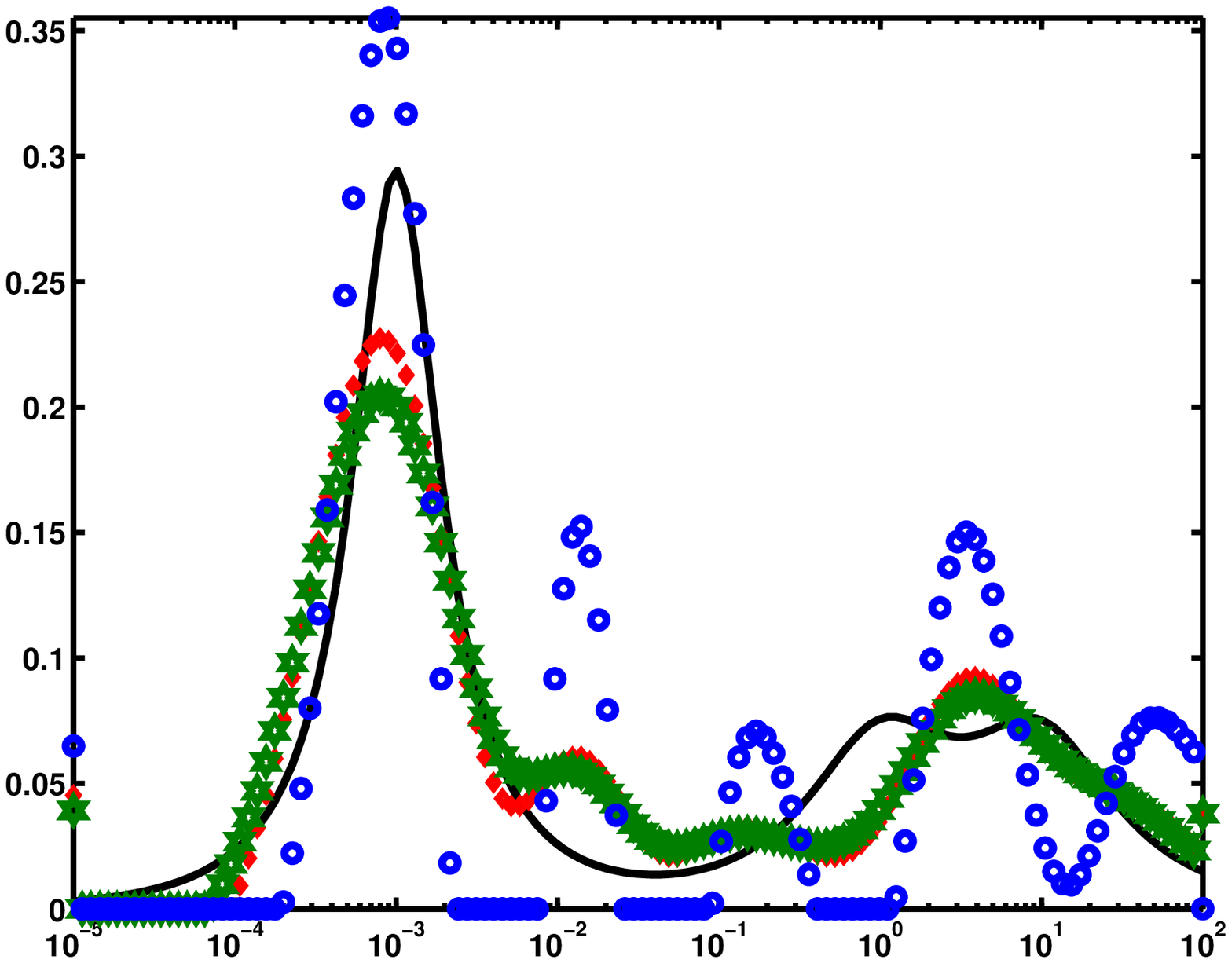}}
\subfigure[$L=L_1$]{\includegraphics[width=1.7in]{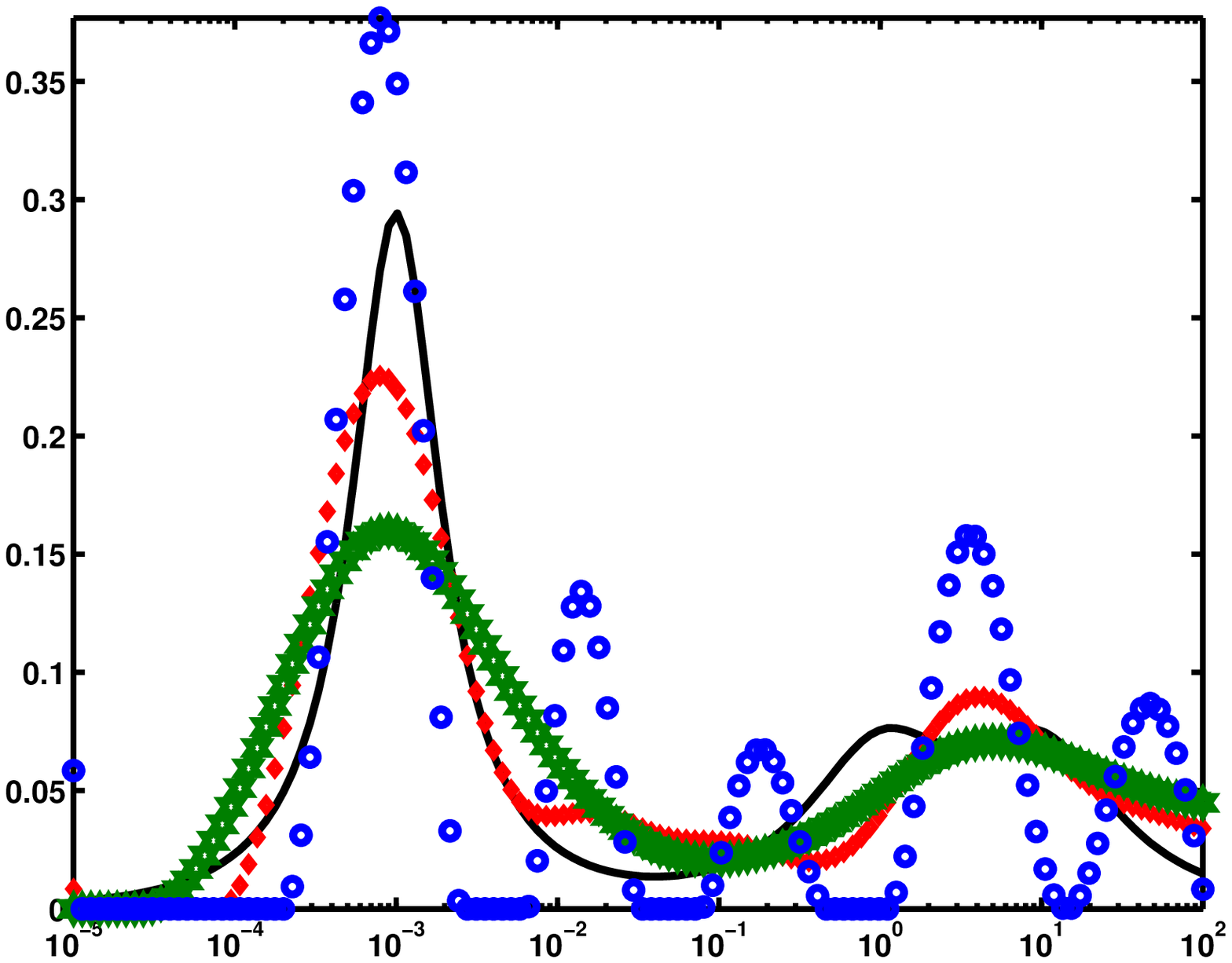}}
\subfigure[$L=L_2$]{\includegraphics[width=1.7in]{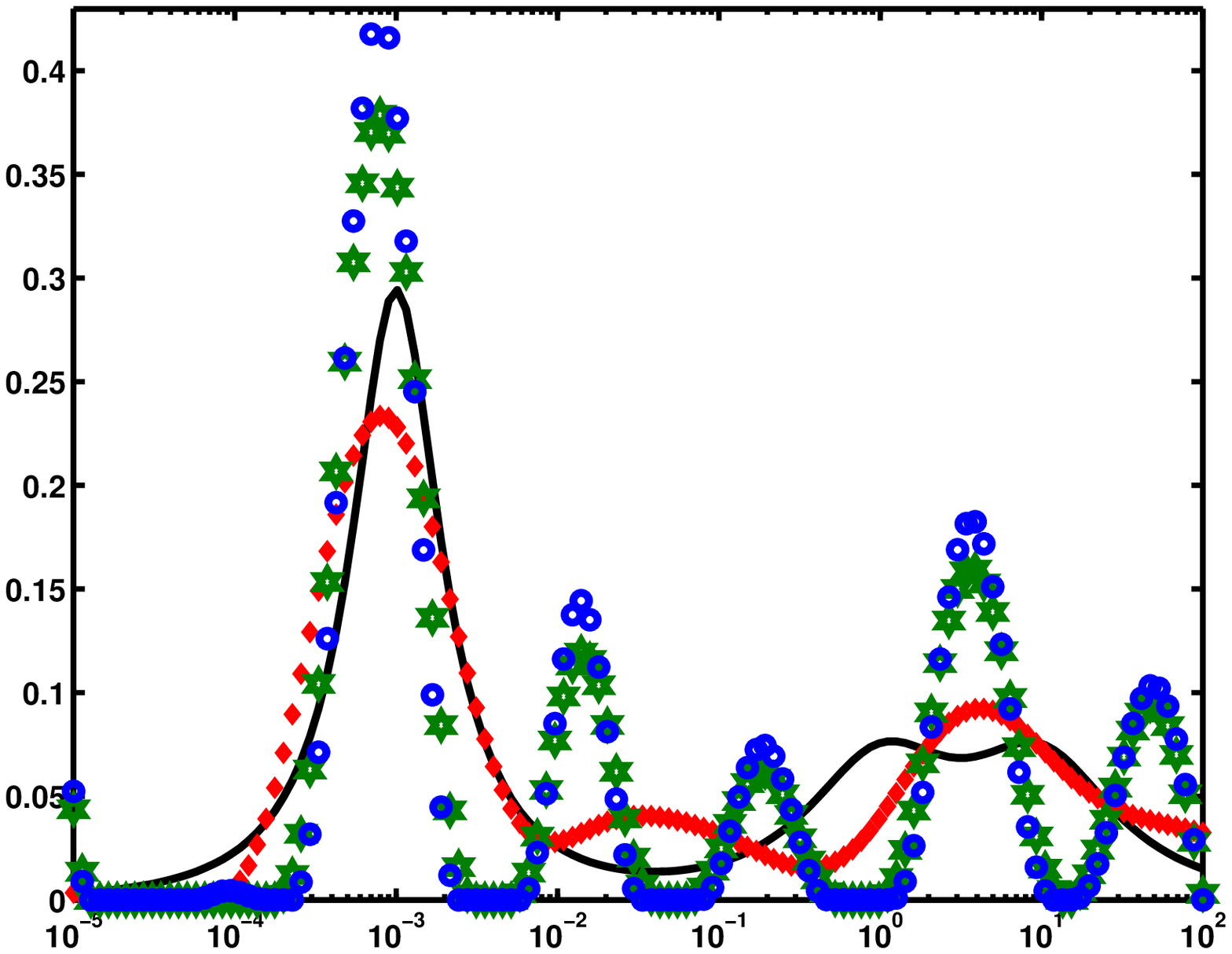}}
\caption{NNLS solutions of RQ-C matrix $A_4$. Noise level $5\%$.}
\label{hnfig-lambdachoiceRQ6A4HN}
\end{figure}
 \begin{figure}[!ht]
\centering
\subfigure[$L=I$]{\includegraphics[width=1.7in]{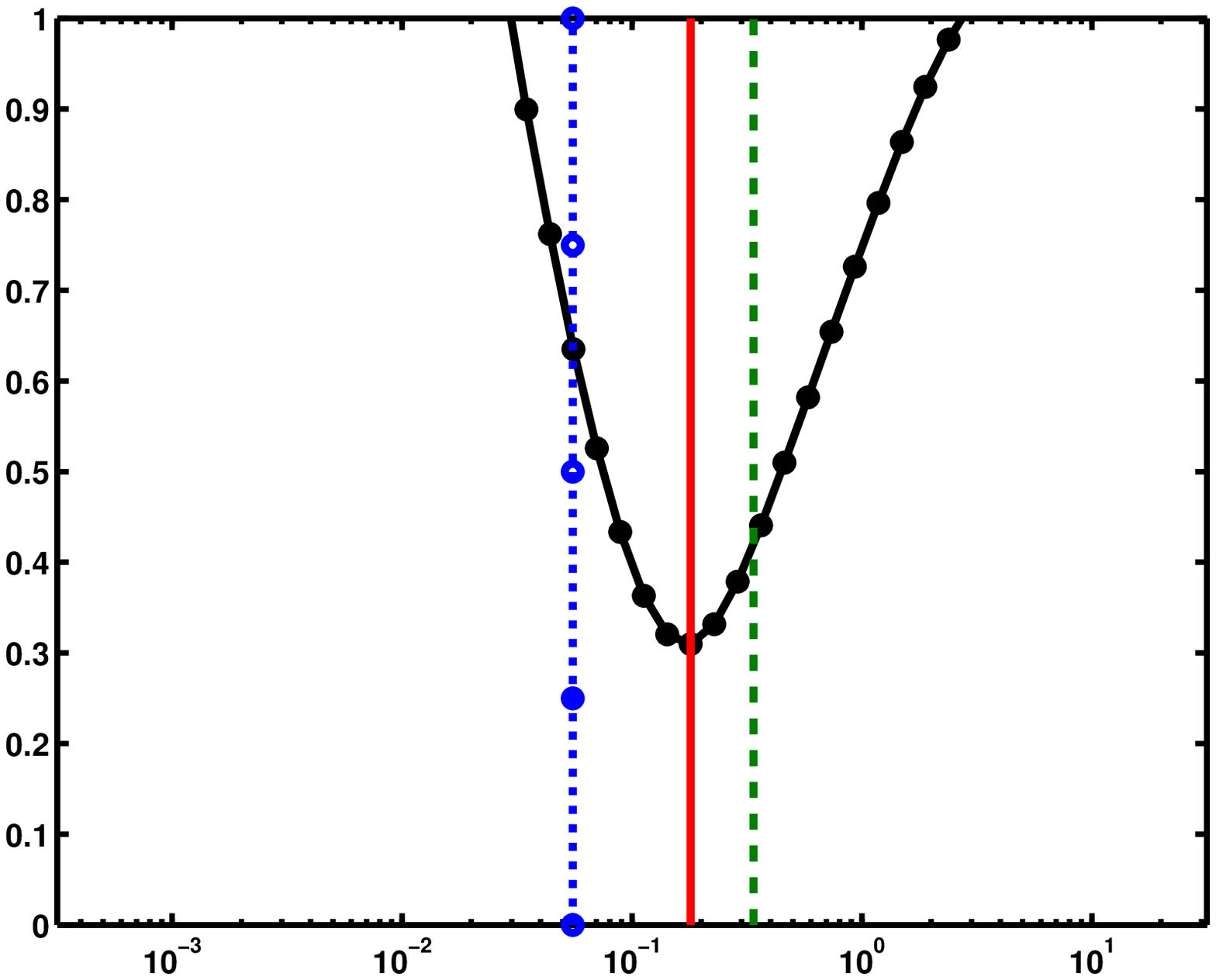}}
\subfigure[$L=L_1$]{\includegraphics[width=1.7in]{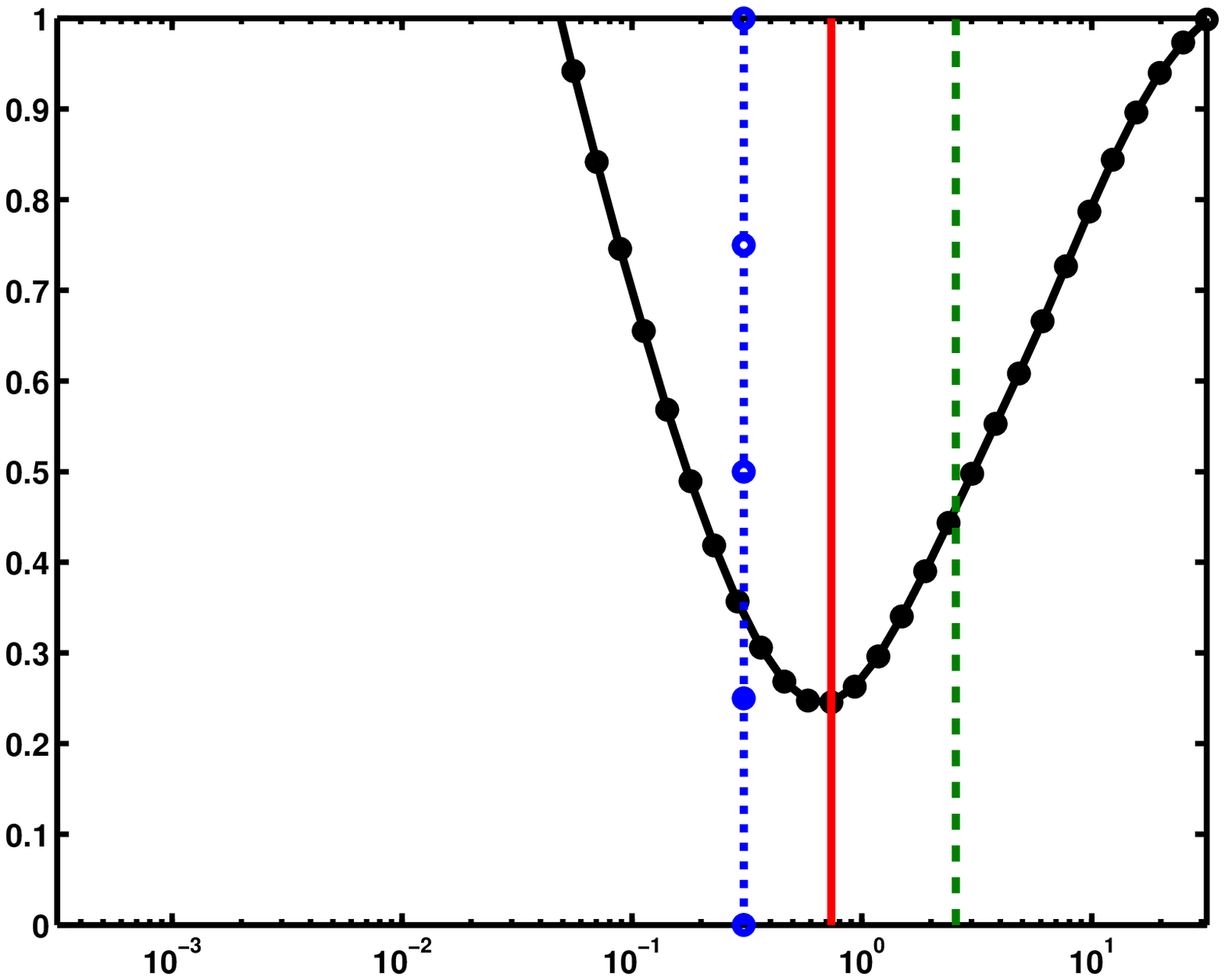}}
\subfigure[$L=L_2$]{\includegraphics[width=1.7in]{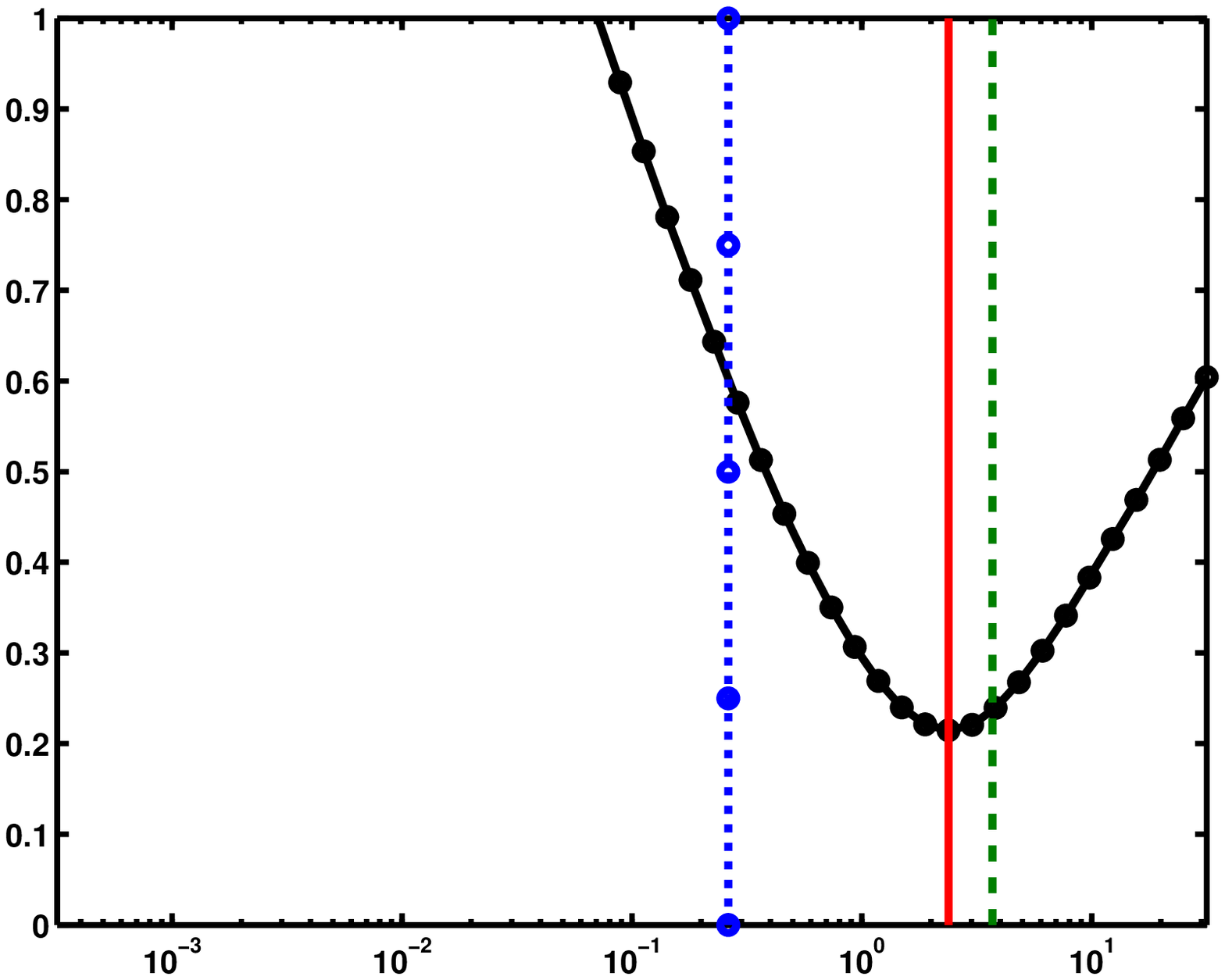}}
\subfigure[$L=I$]{\includegraphics[width=1.7in]{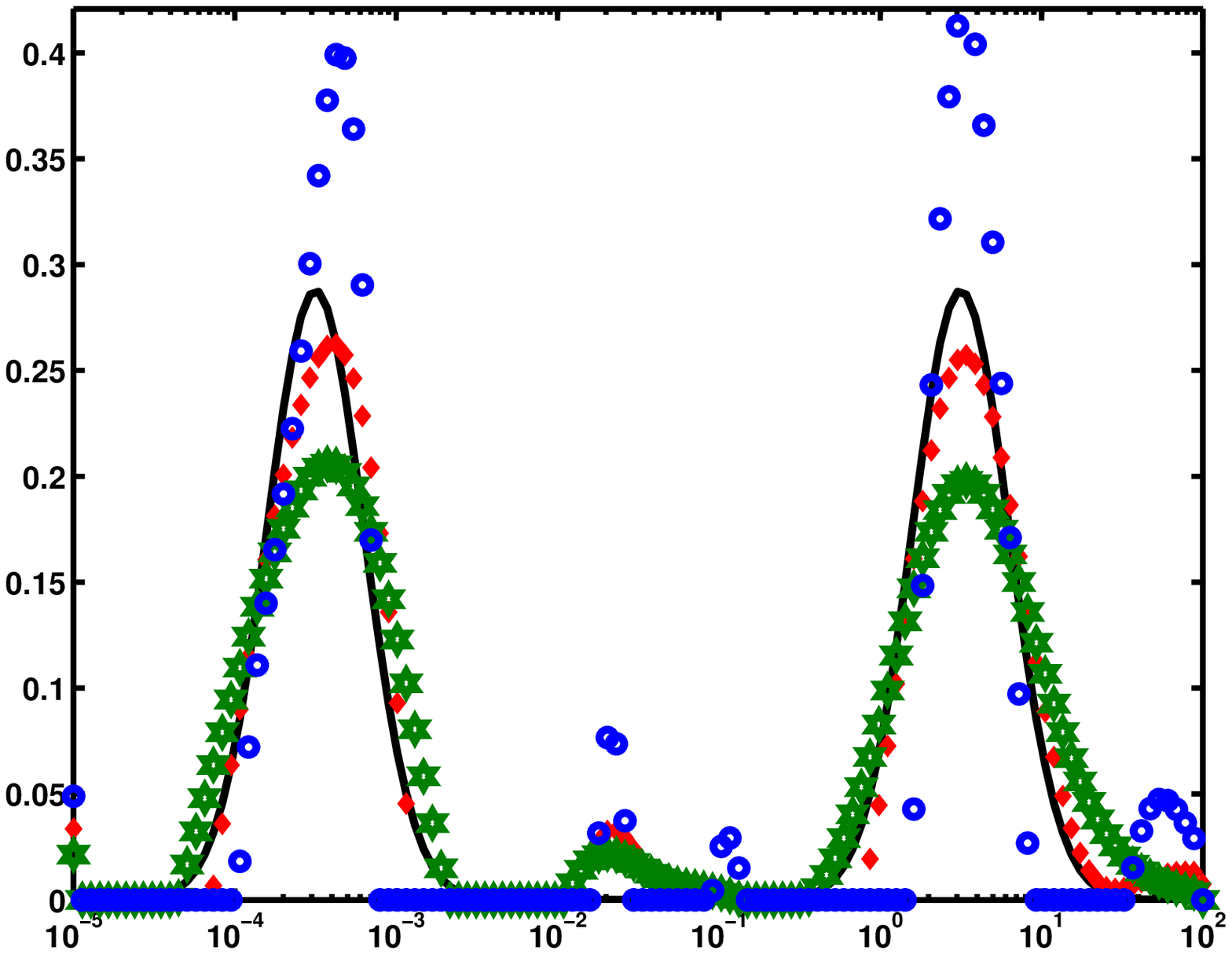}}
\subfigure[$L=L_1$]{\includegraphics[width=1.7in]{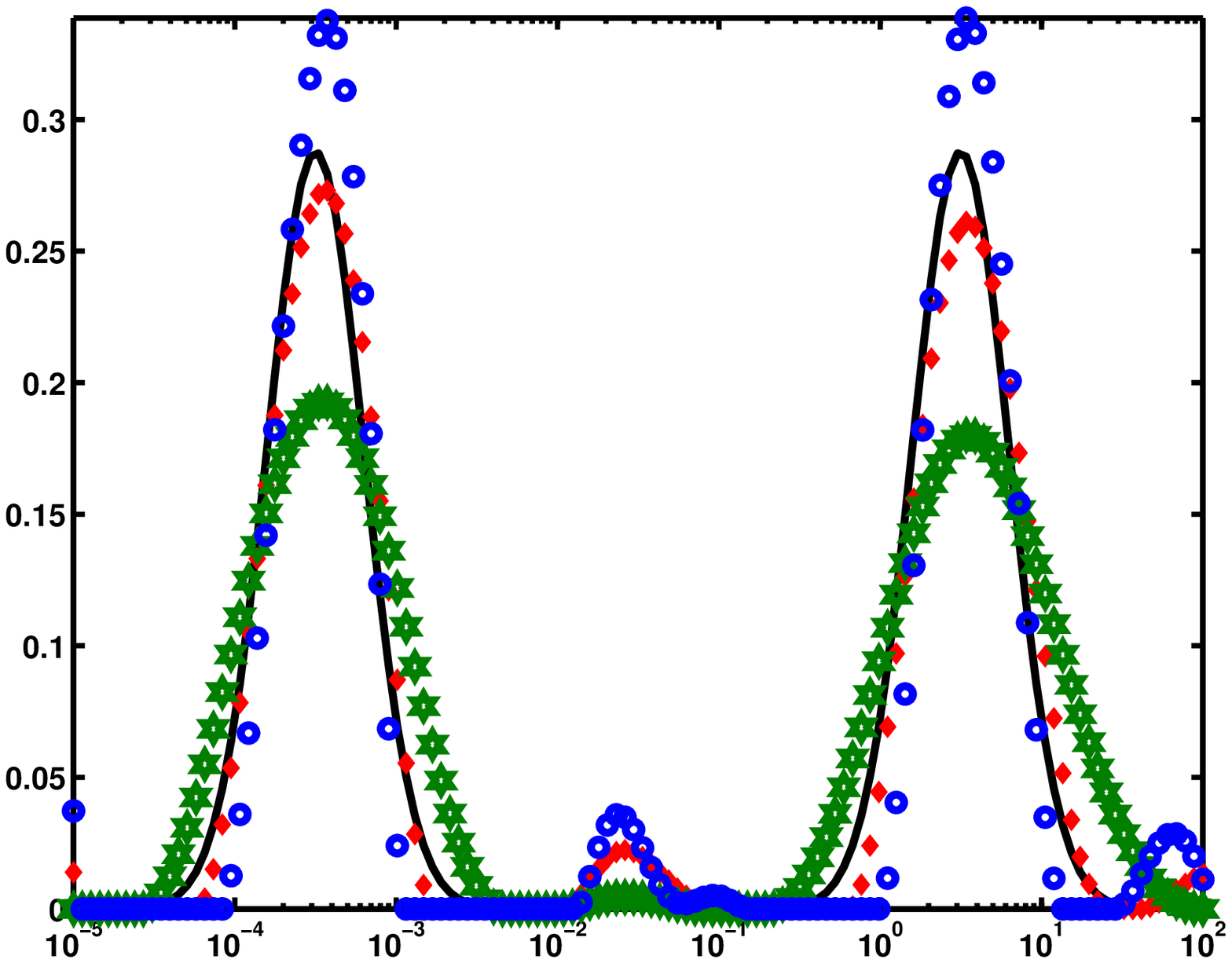}}
\subfigure[$L=L_2$]{\includegraphics[width=1.7in]{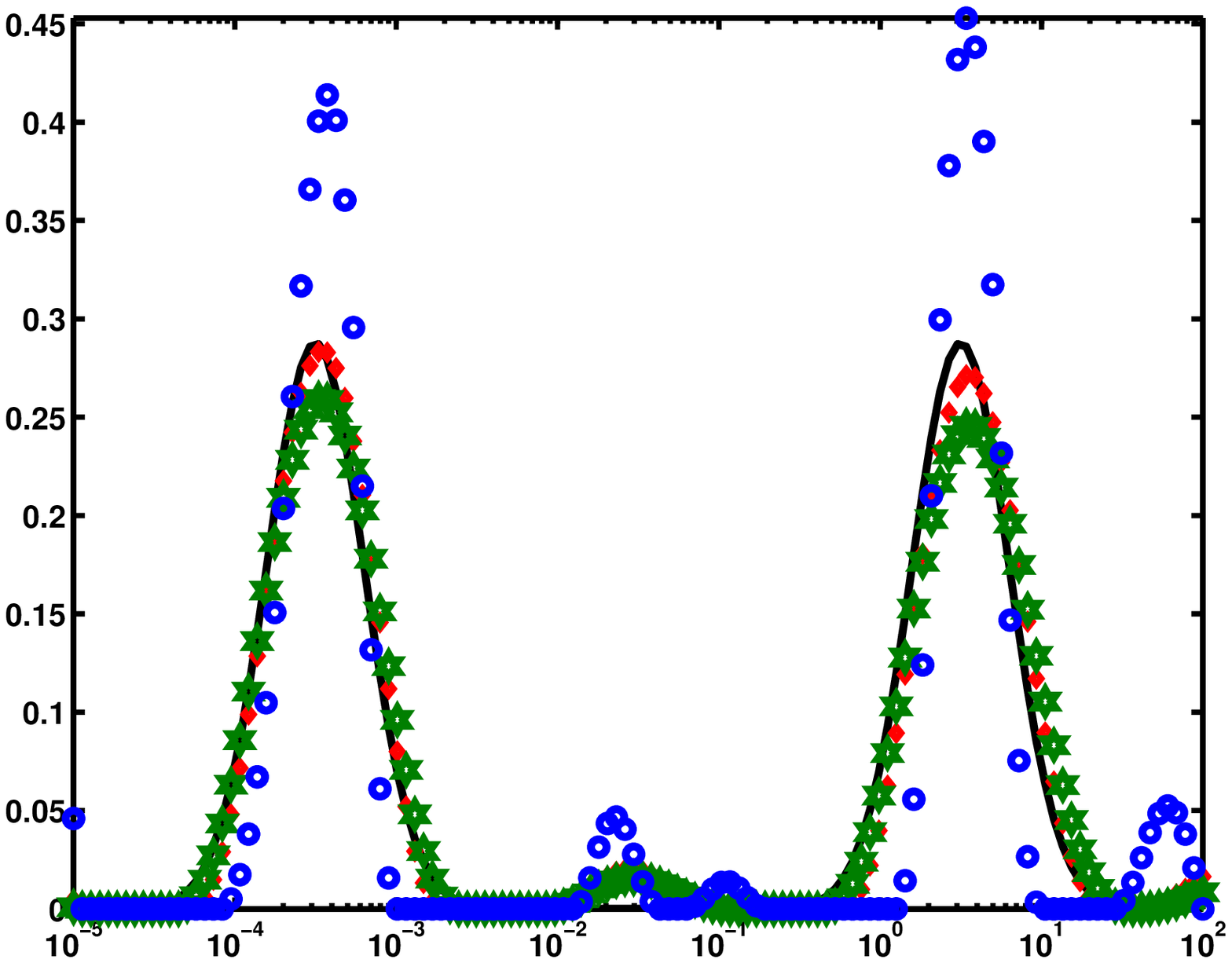}}
\caption{NNLS solutions of LN-A matrix $A_4$. Noise level $5\%$.}
\label{hnfig-lambdachoiceLN2A4HN}
\end{figure}
 \begin{figure}[!ht]
\centering
\subfigure[$L=I$]{\includegraphics[width=1.7in]{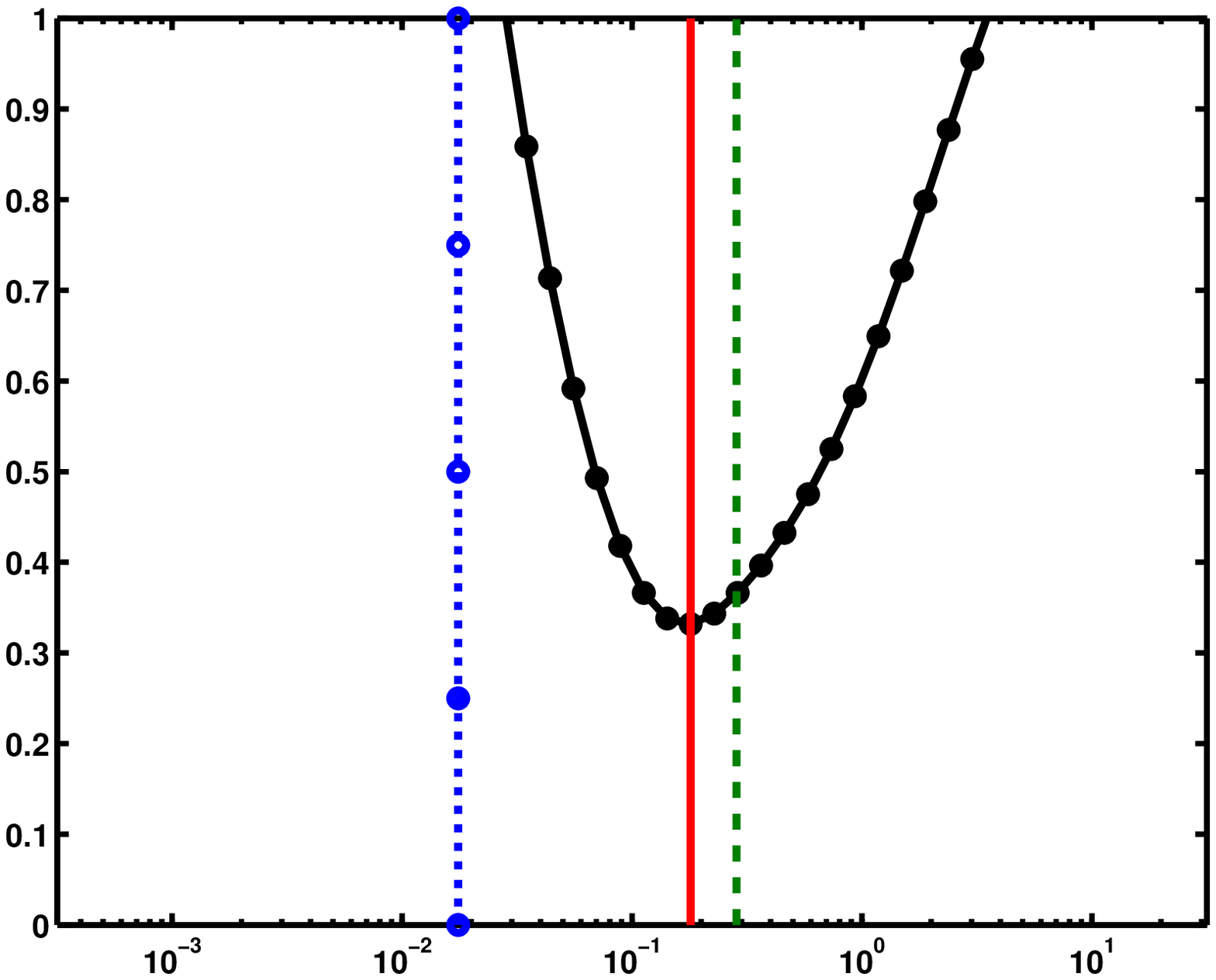}}
\subfigure[$L=L_1$]{\includegraphics[width=1.7in]{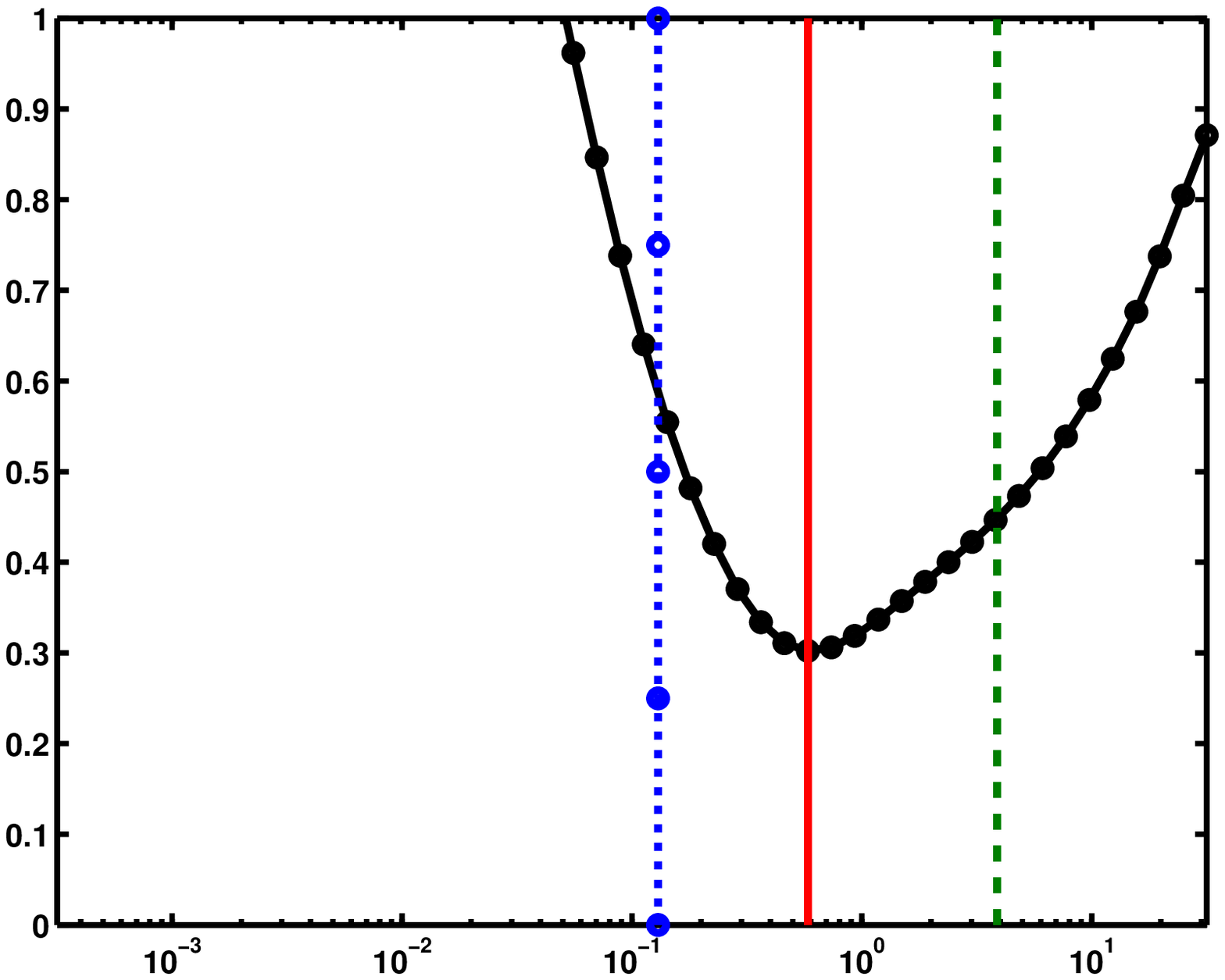}}
\subfigure[$L=L_2$]{\includegraphics[width=1.7in]{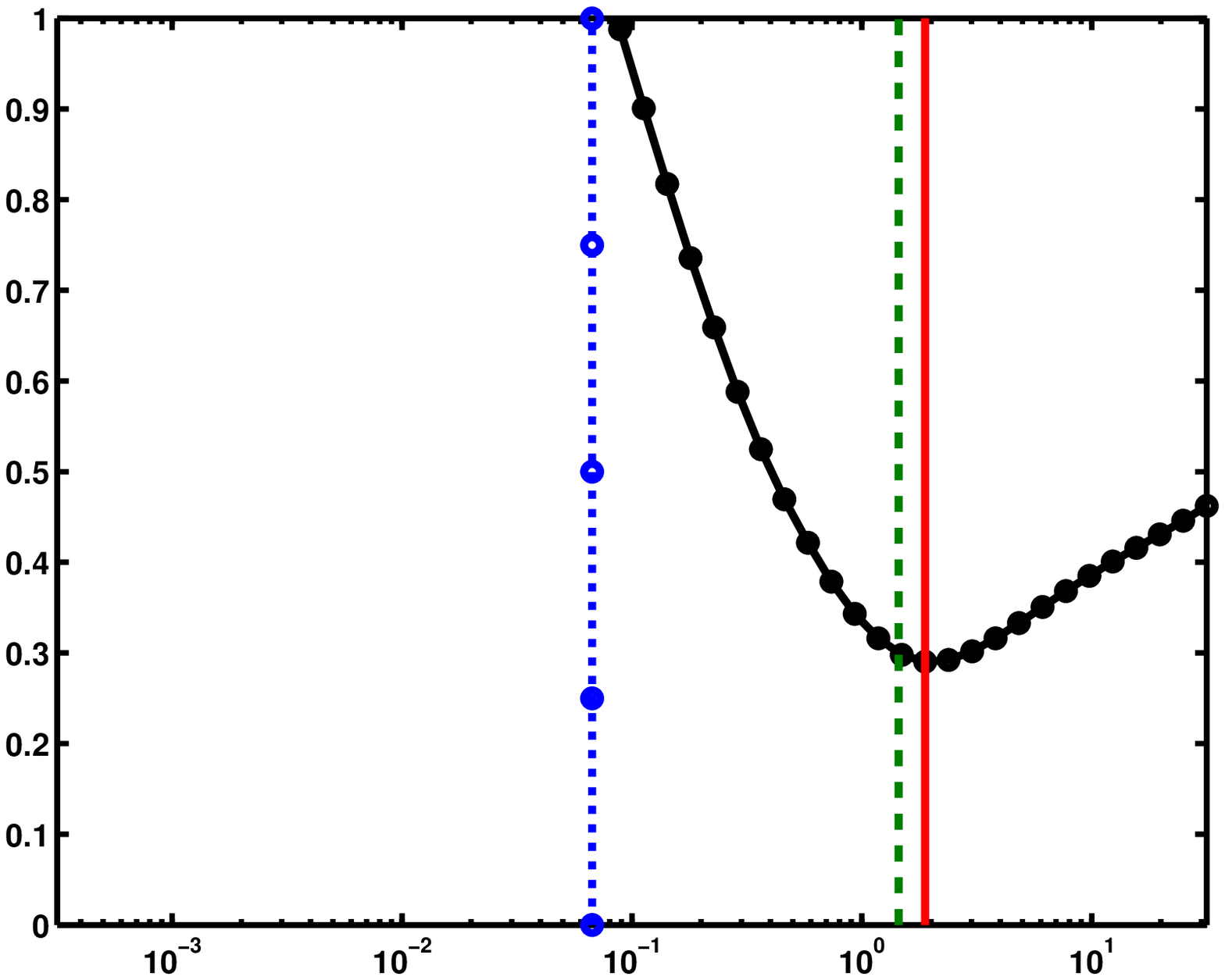}}
\subfigure[$L=I$]{\includegraphics[width=1.7in]{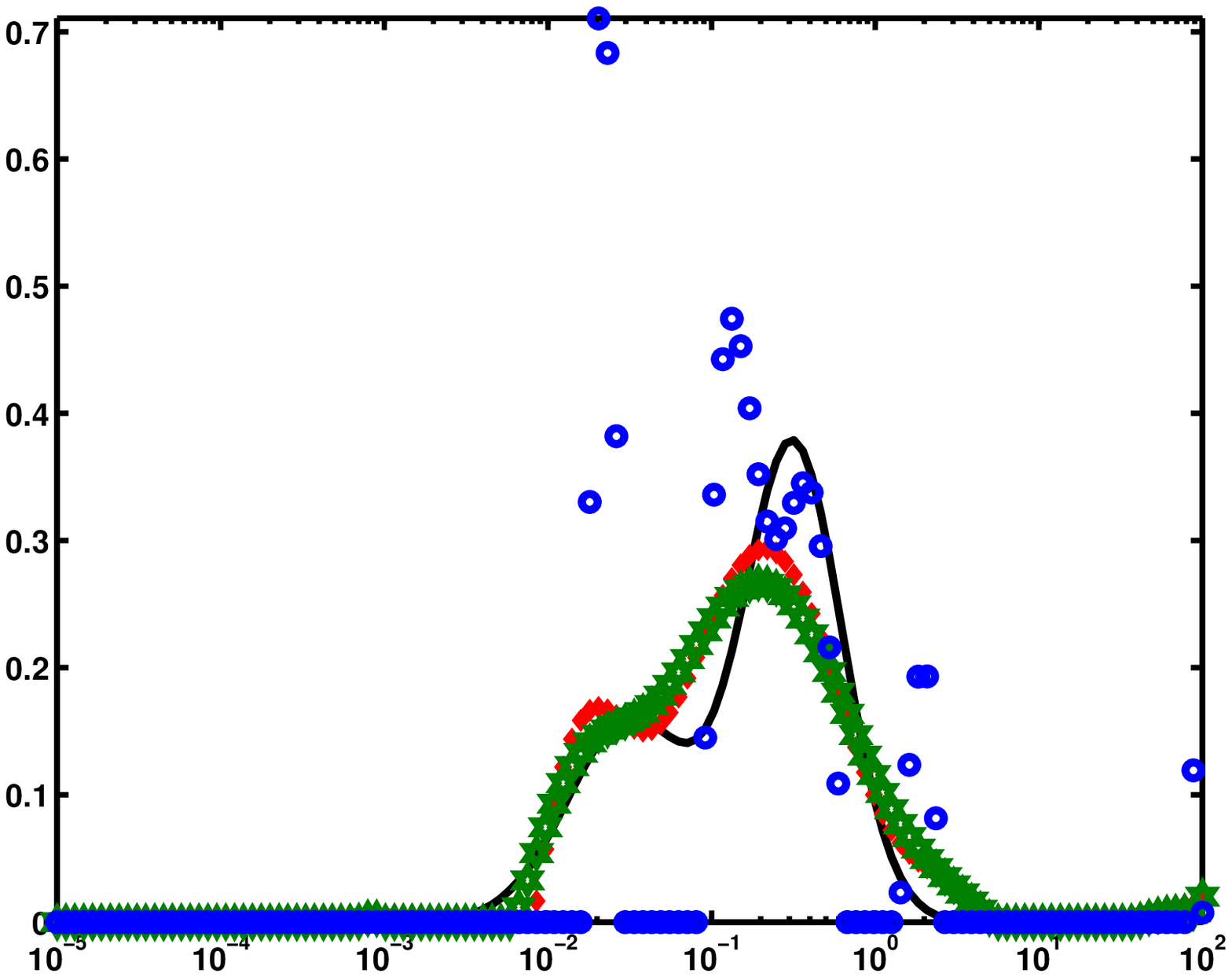}}
\subfigure[$L=L_1$]{\includegraphics[width=1.7in]{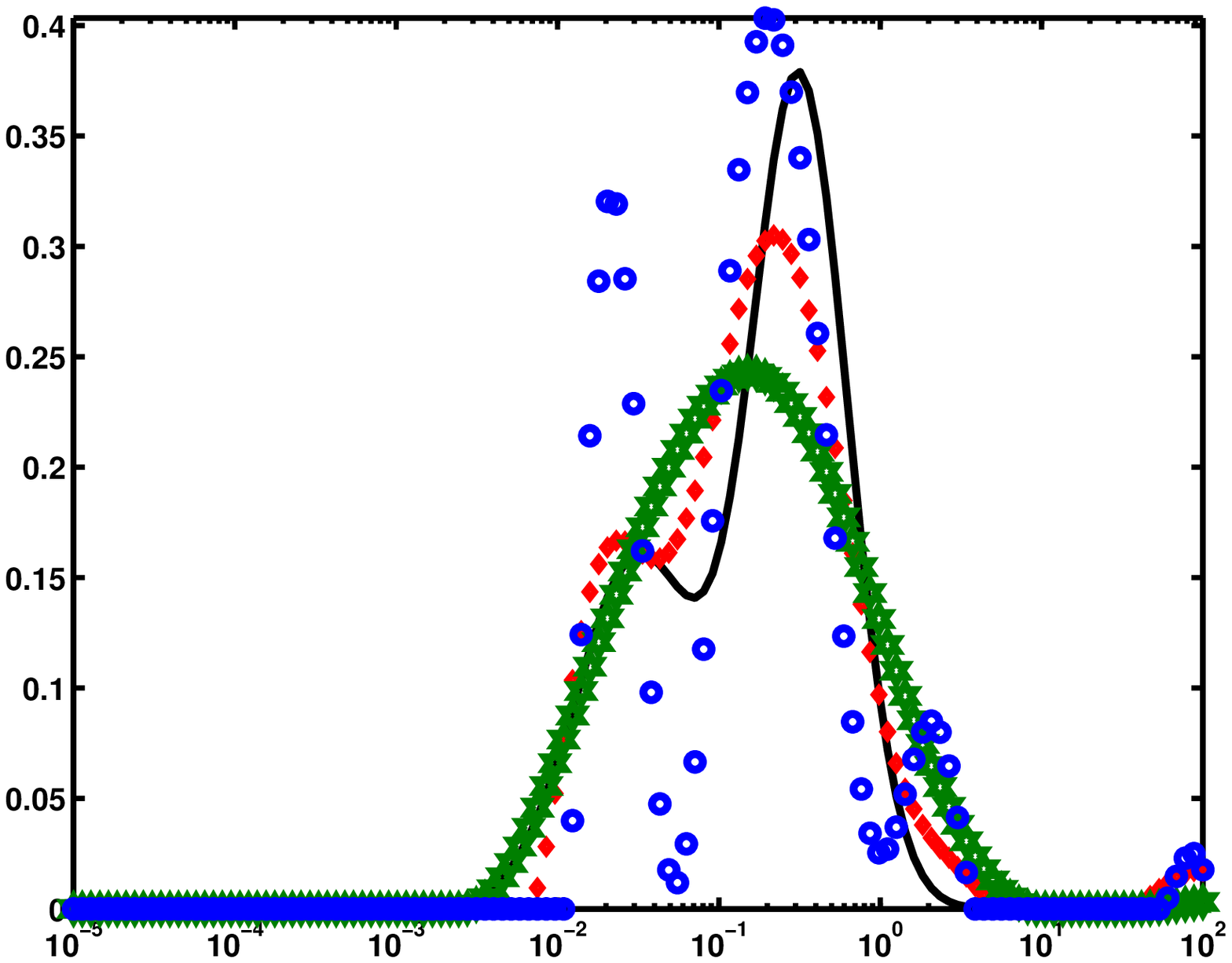}}
\subfigure[$L=L_2$]{\includegraphics[width=1.7in]{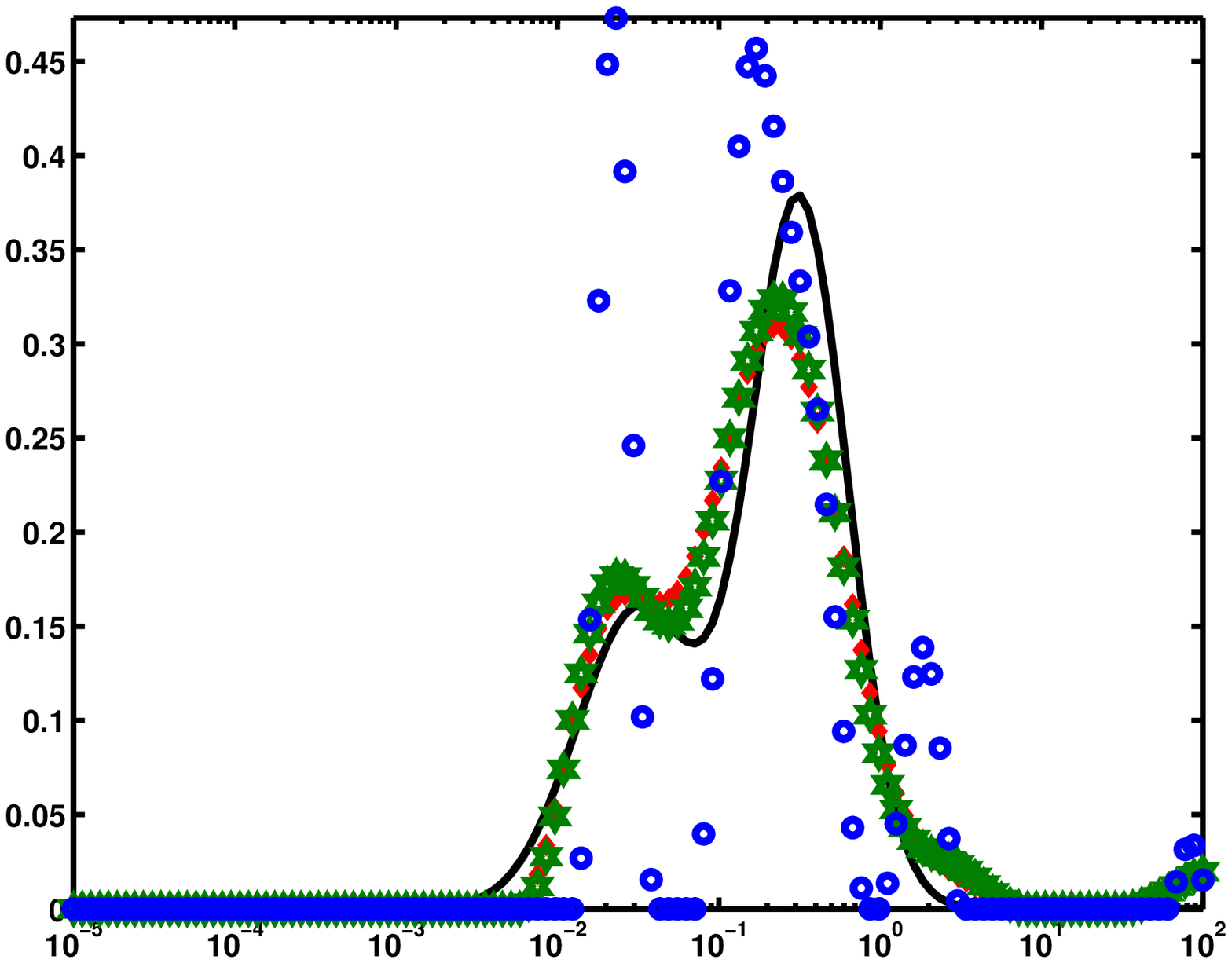}}
\caption{NNLS solutions of LN-B matrix $A_4$. Noise level $5\%$.}
\label{hnfig-lambdachoiceLN5A4HN}
\end{figure}
 \begin{figure}[!ht]
\centering
\subfigure[$L=I$]{\includegraphics[width=1.7in]{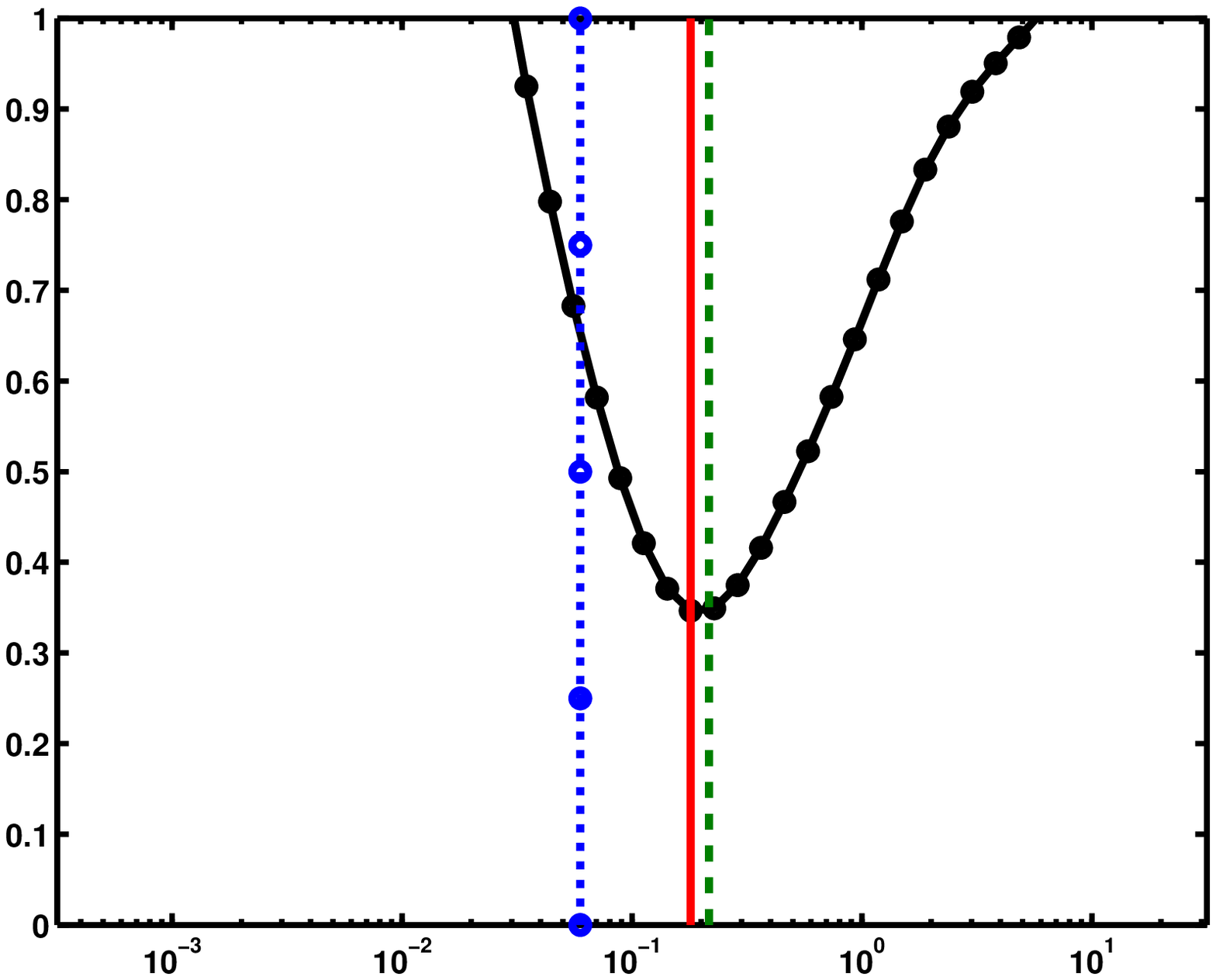}}
\subfigure[$L=L_1$]{\includegraphics[width=1.7in]{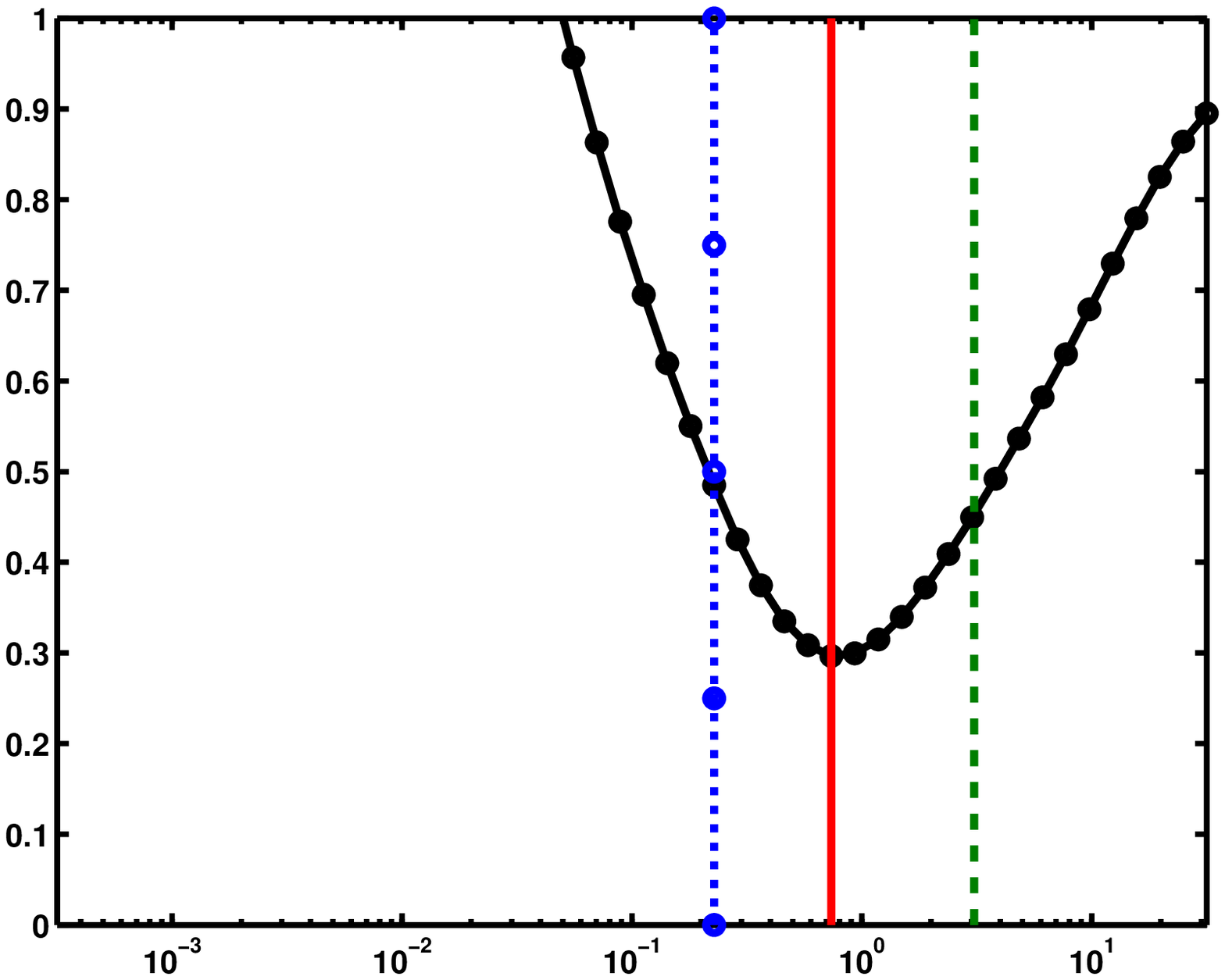}}
\subfigure[$L=L_2$]{\includegraphics[width=1.7in]{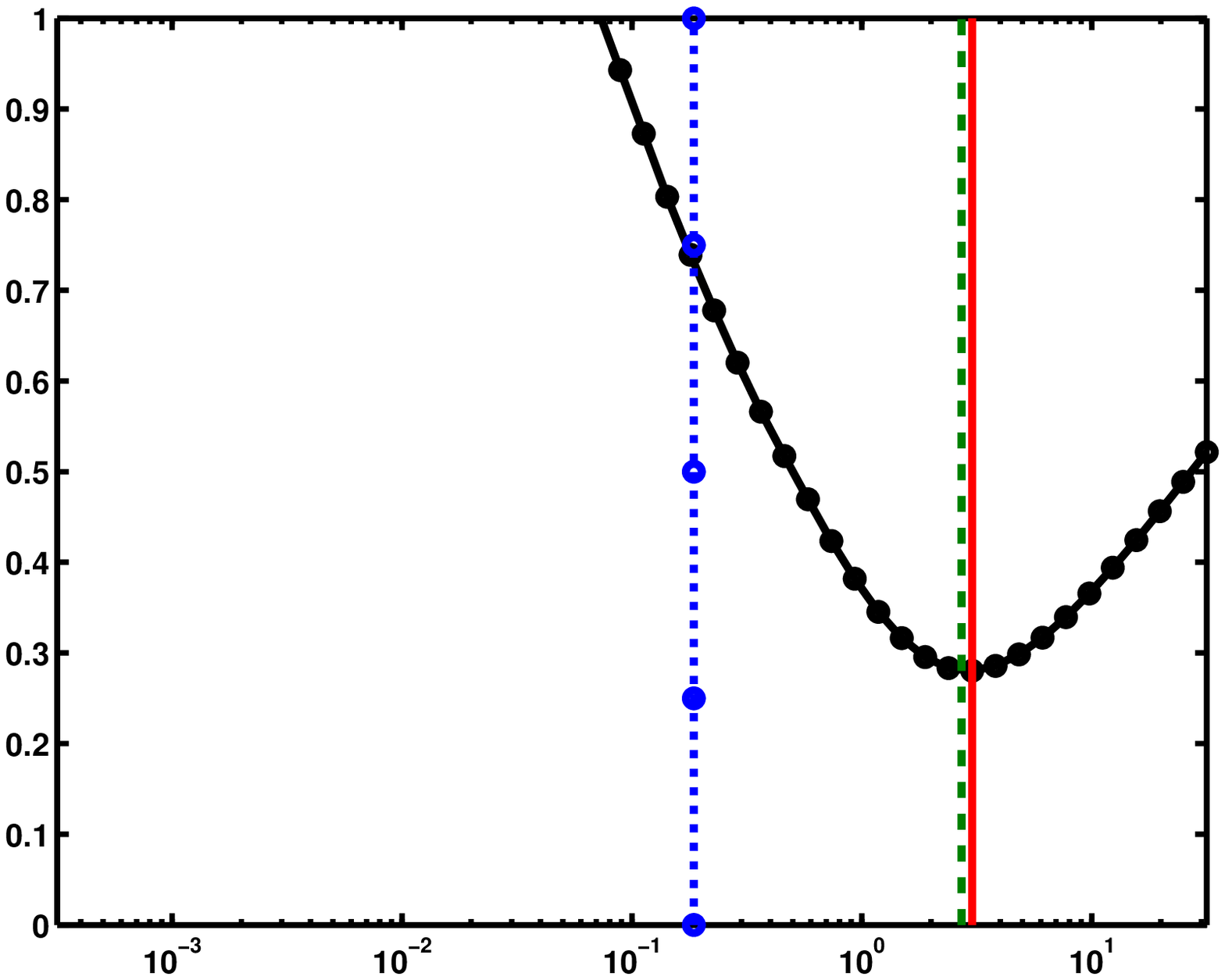}}
\subfigure[$L=I$]{\includegraphics[width=1.7in]{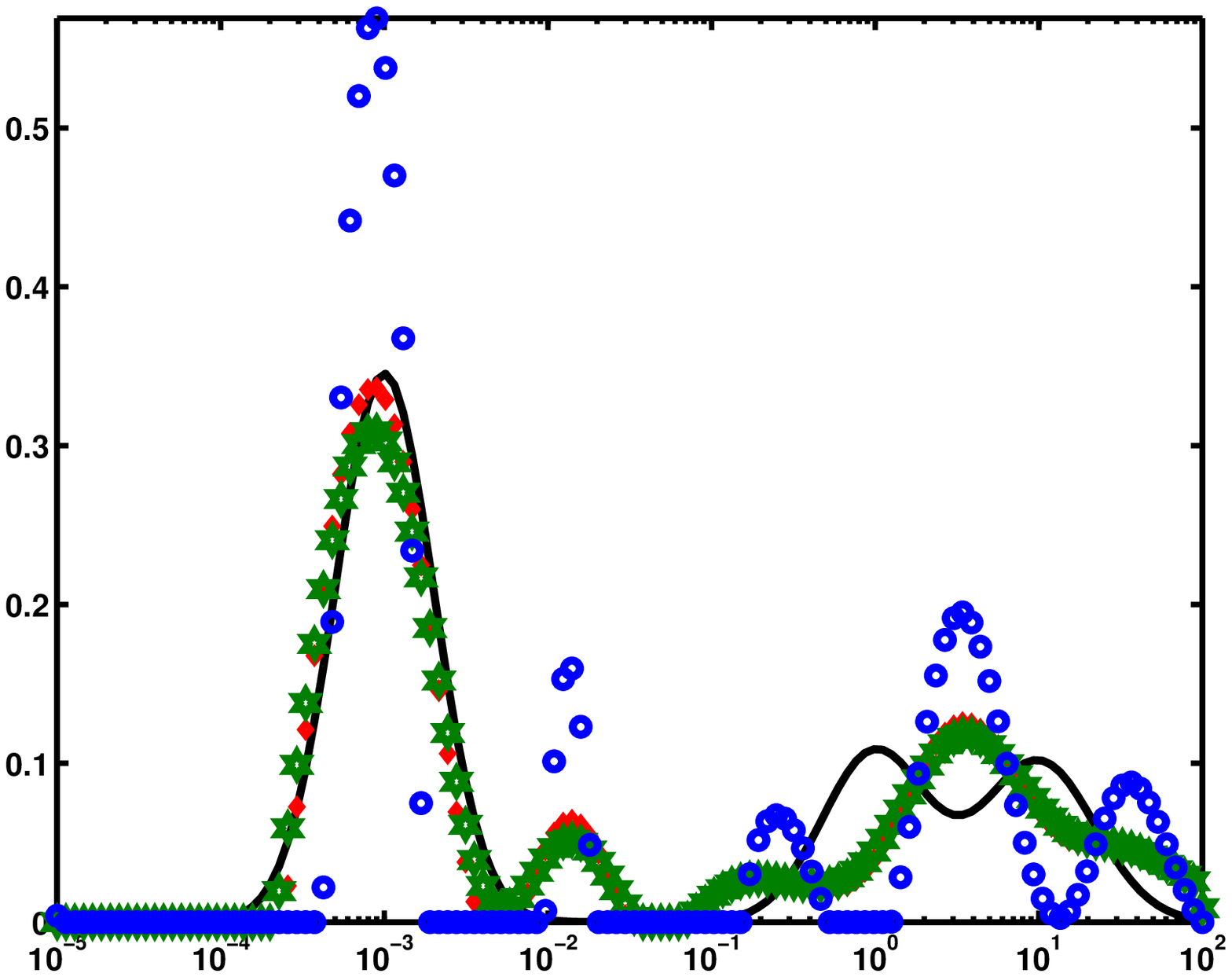}}
\subfigure[$L=L_1$]{\includegraphics[width=1.7in]{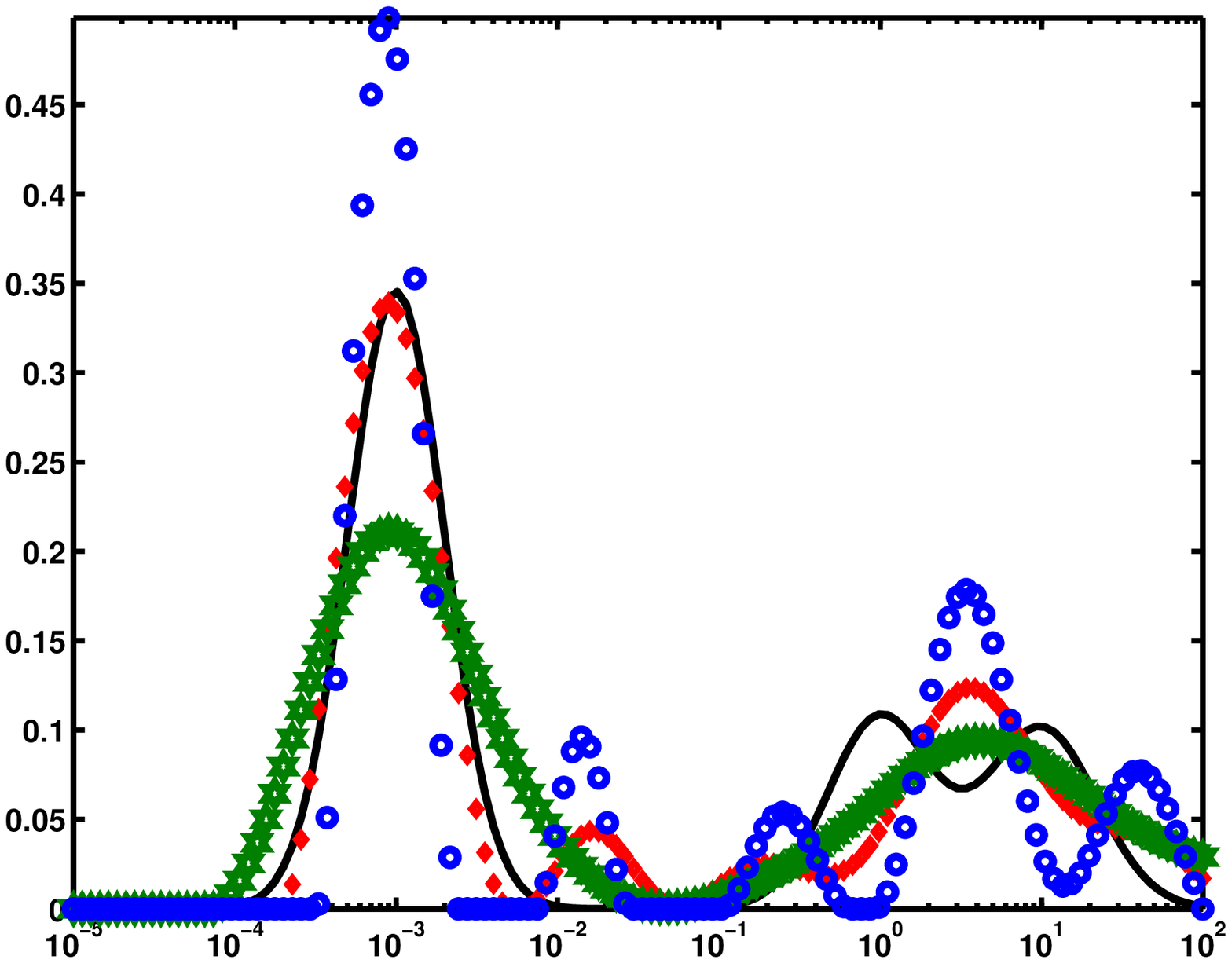}}
\subfigure[$L=L_2$]{\includegraphics[width=1.7in]{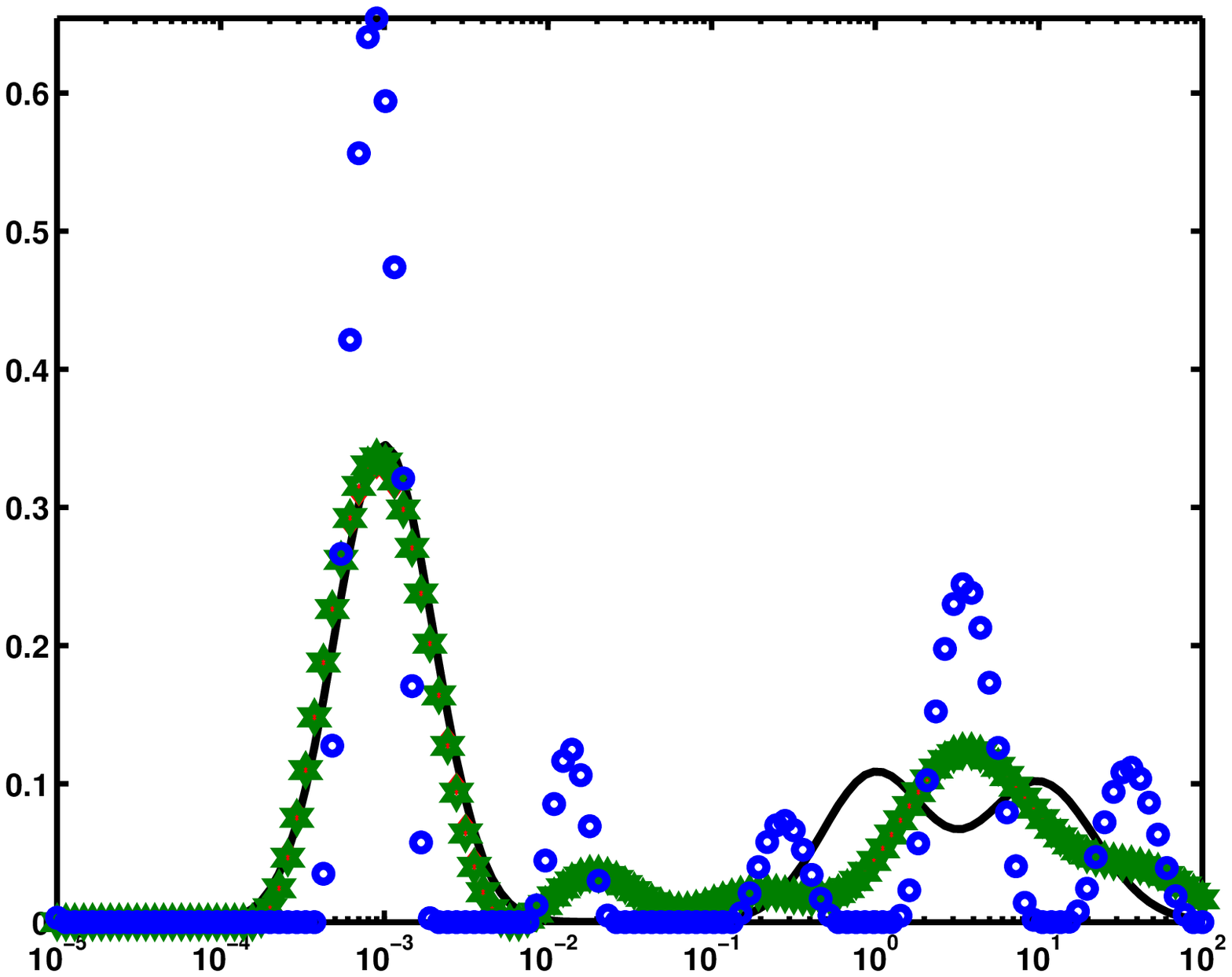}}
\caption{NNLS solutions of LN-C matrix $A_4$. Noise level $5\%$.}
\label{hnfig-lambdachoiceLN6A4HN}
\end{figure}
\clearpage
\subsection{Noise level $5\%$ $A_3$ NNLS}
\begin{figure}[!ht]
\centering
\subfigure[$L=I$]{\includegraphics[width=1.7in]{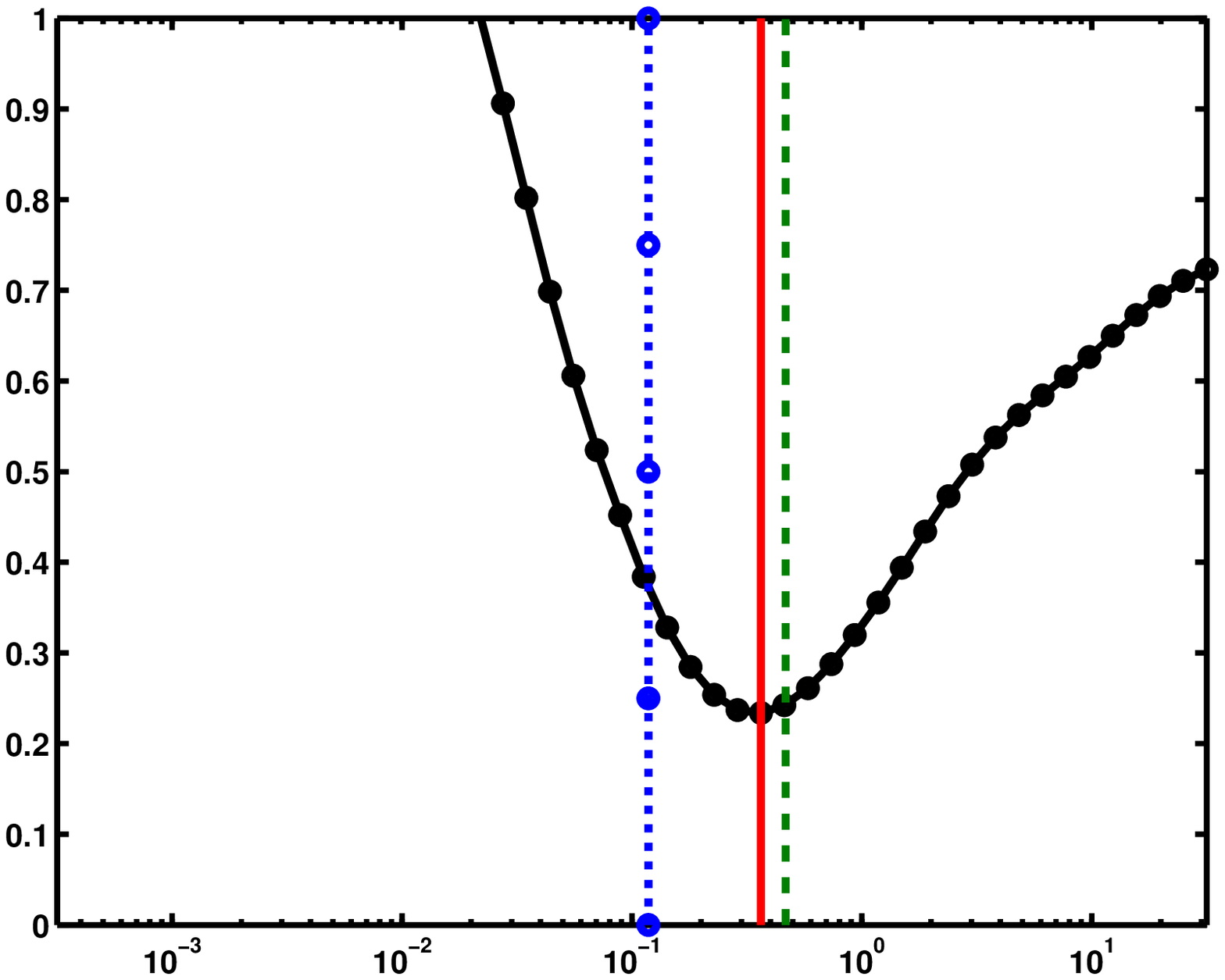}}
\subfigure[$L=L_1$]{\includegraphics[width=1.7in]{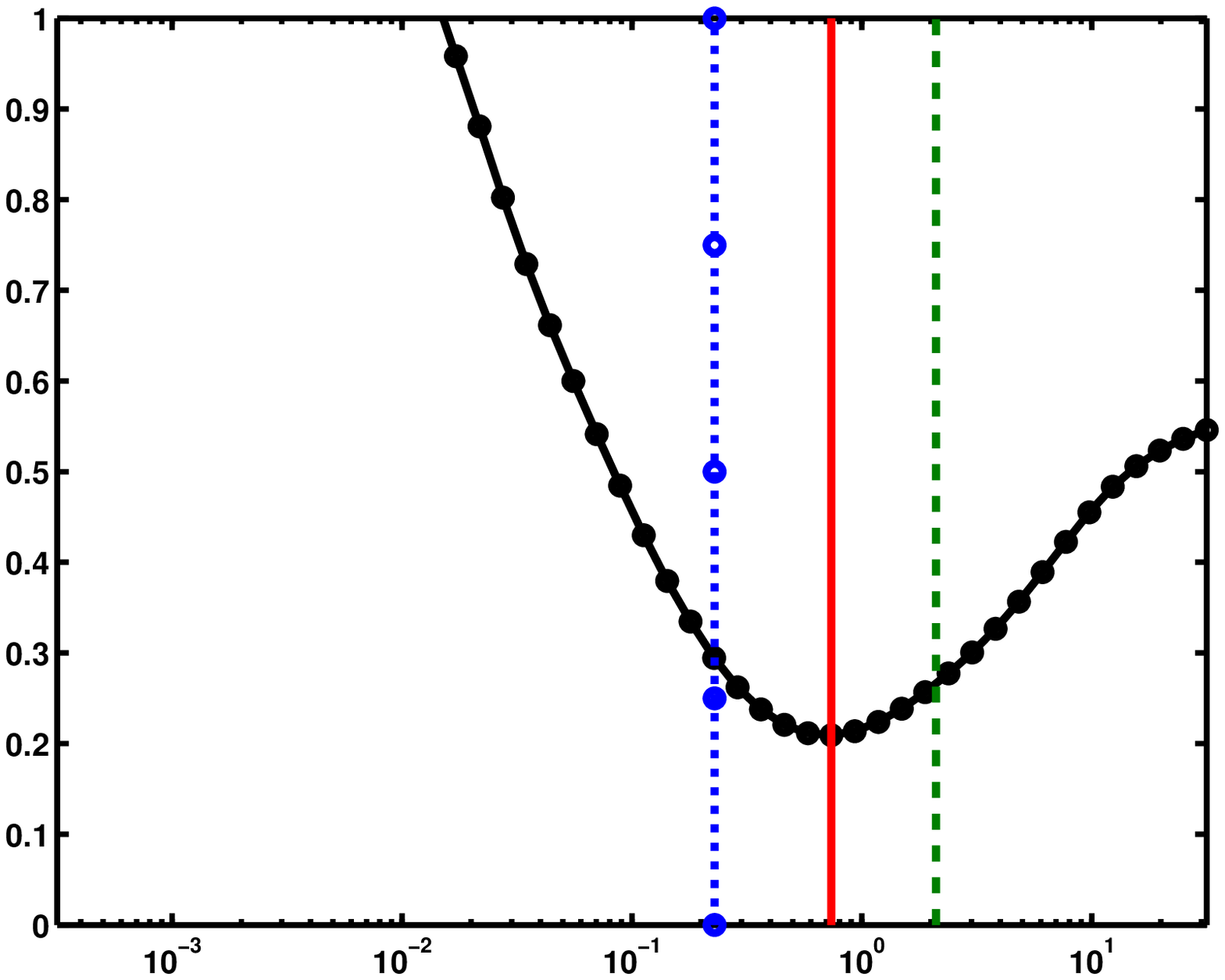}}
\subfigure[$L=L_2$]{\includegraphics[width=1.7in]{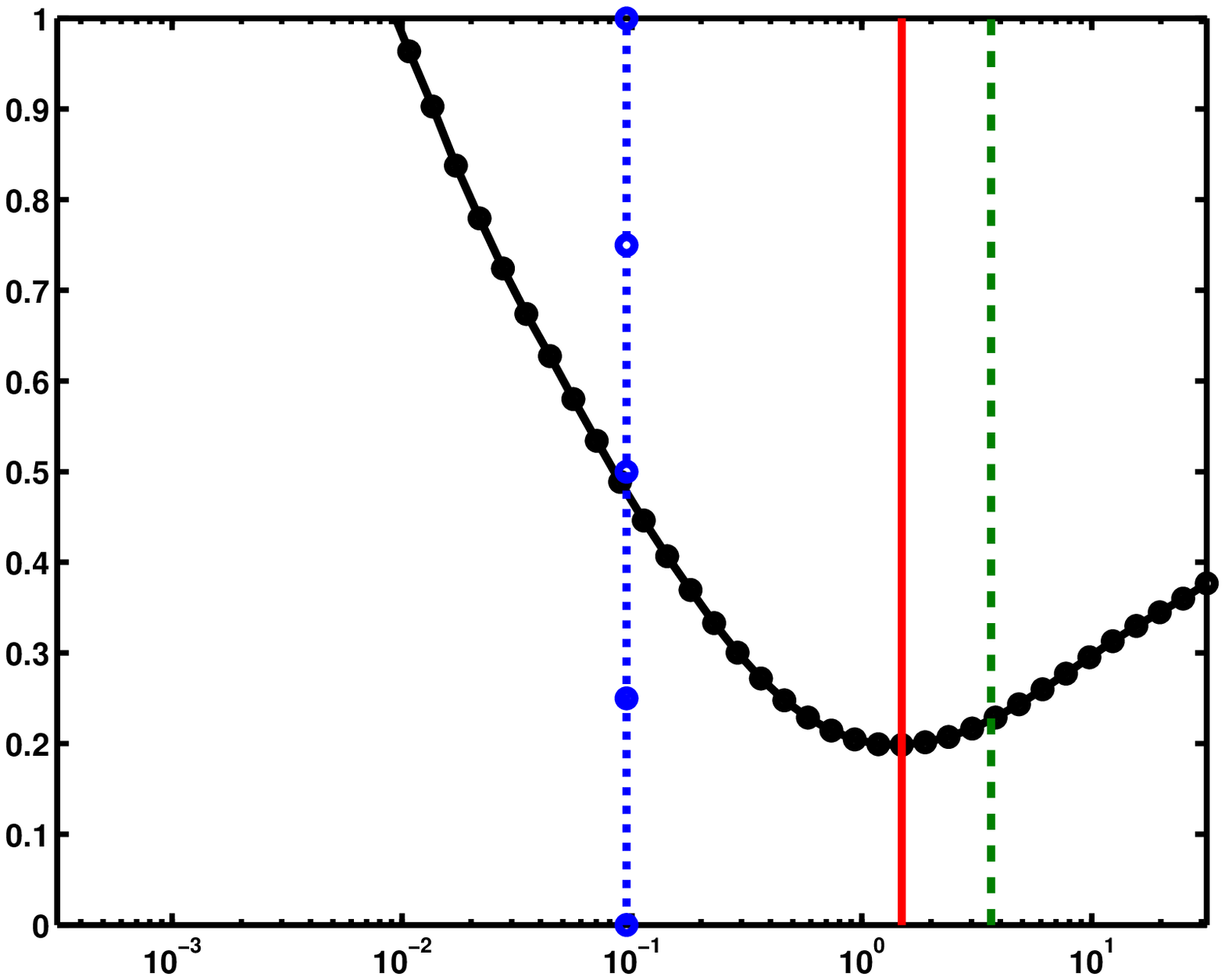}}
\subfigure[$L=I$]{\includegraphics[width=1.7in]{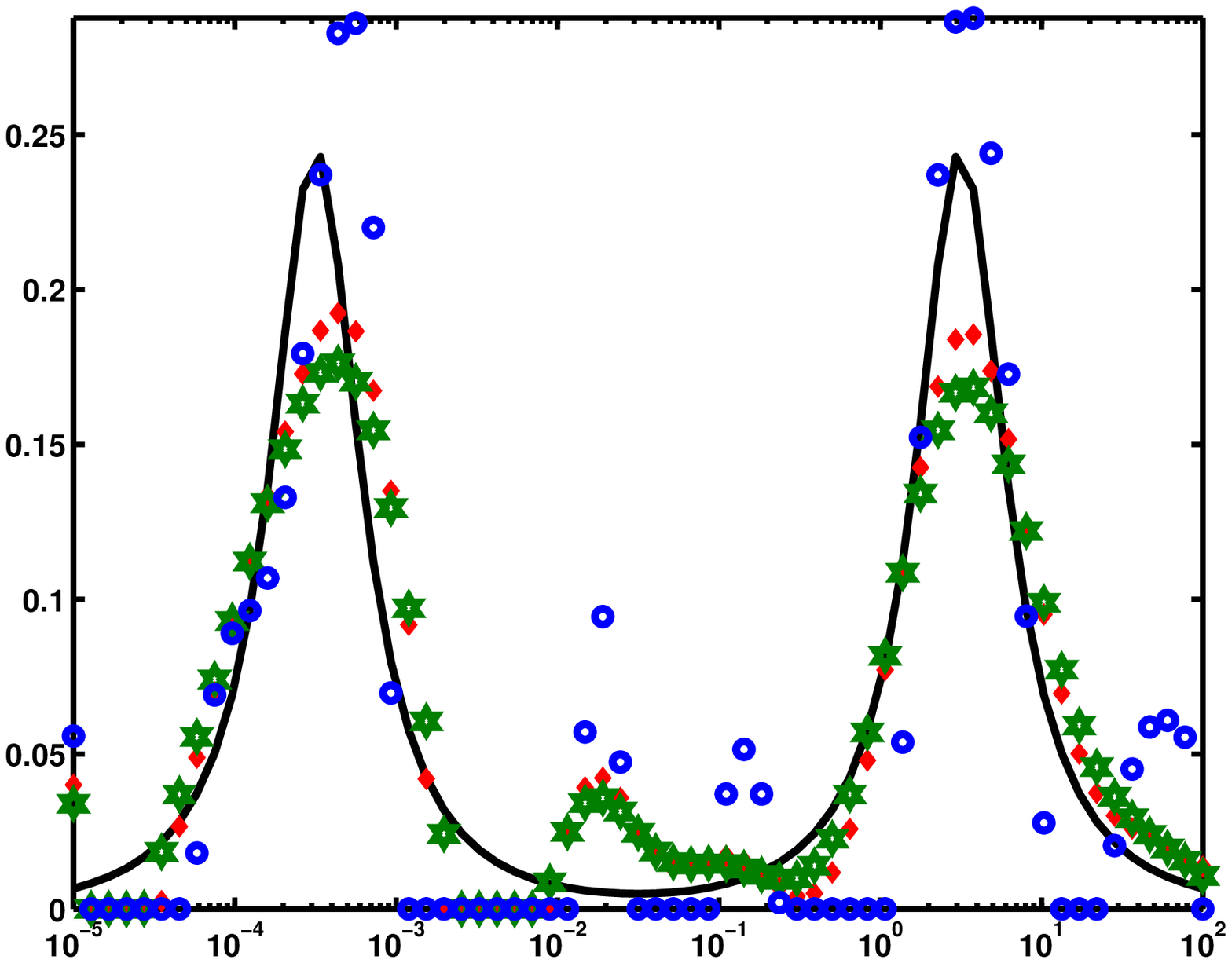}}
\subfigure[$L=L_1$]{\includegraphics[width=1.7in]{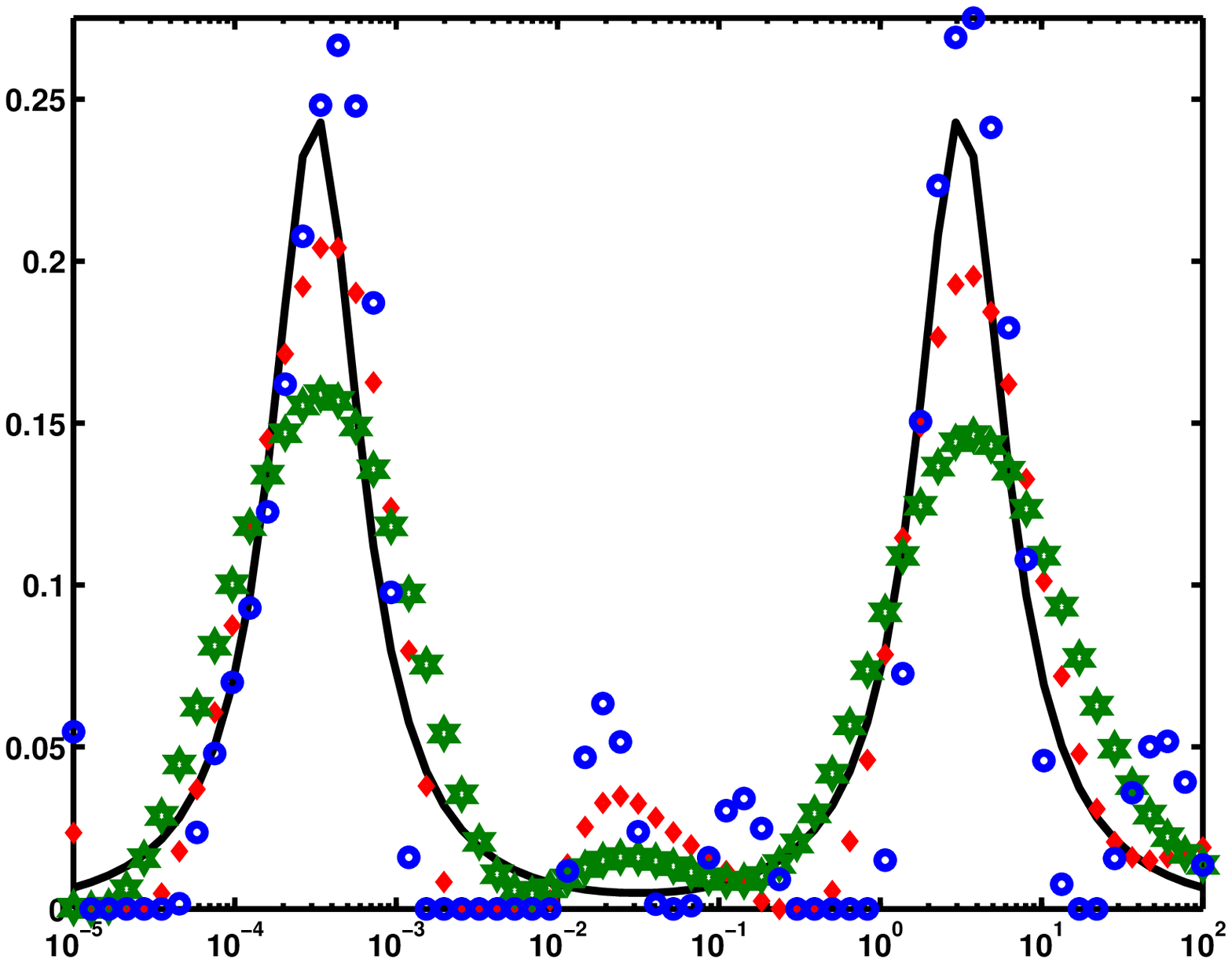}}
\subfigure[$L=L_2$]{\includegraphics[width=1.7in]{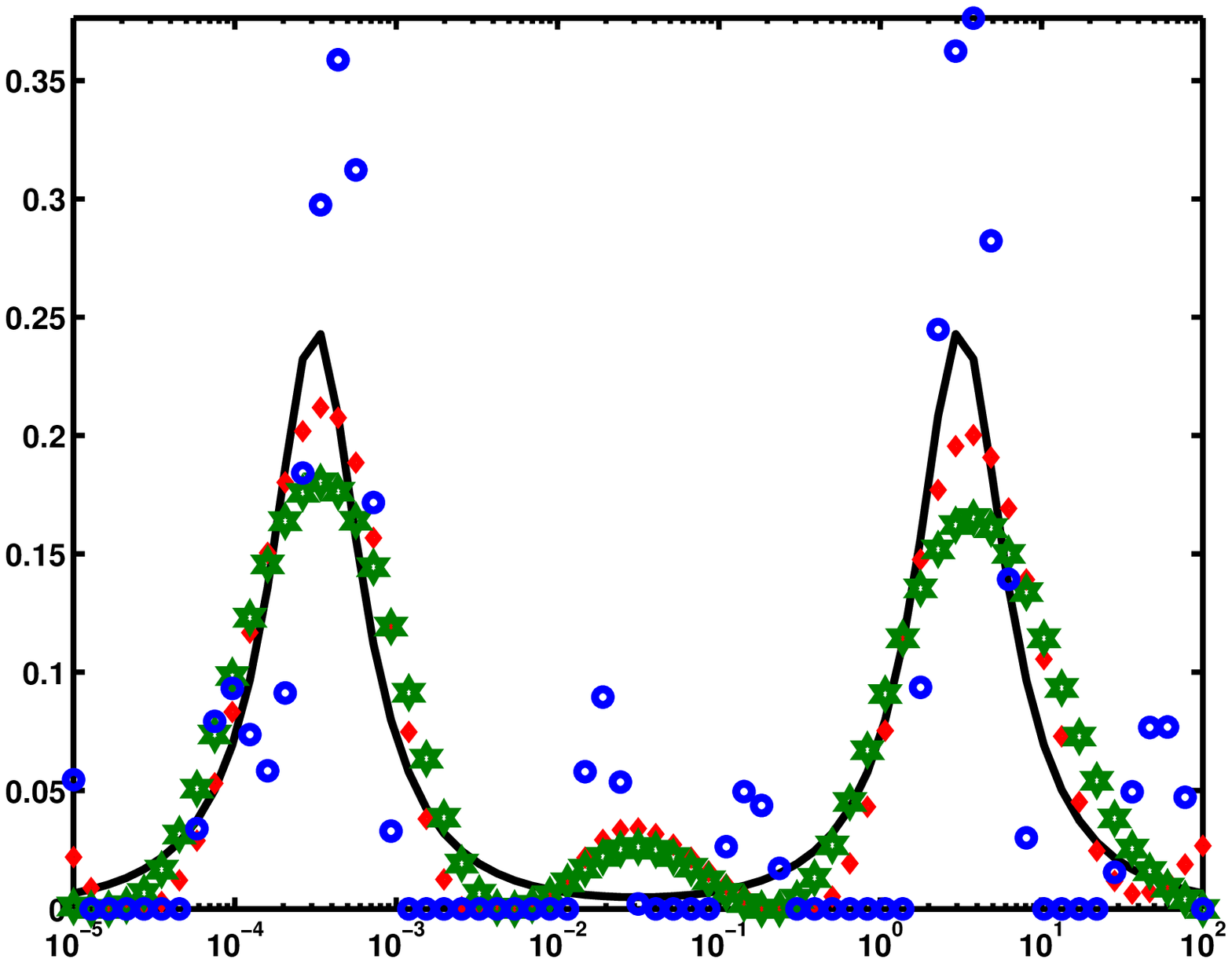}}
\caption{NNLS solutions of RQ-A matrix $A_3$. Noise level $5\%$.}
\label{hnfig-lambdachoiceRQ1A3HN}
\end{figure}

\begin{figure}[!ht]
\centering
\subfigure[$L=I$]{\includegraphics[width=1.7in]{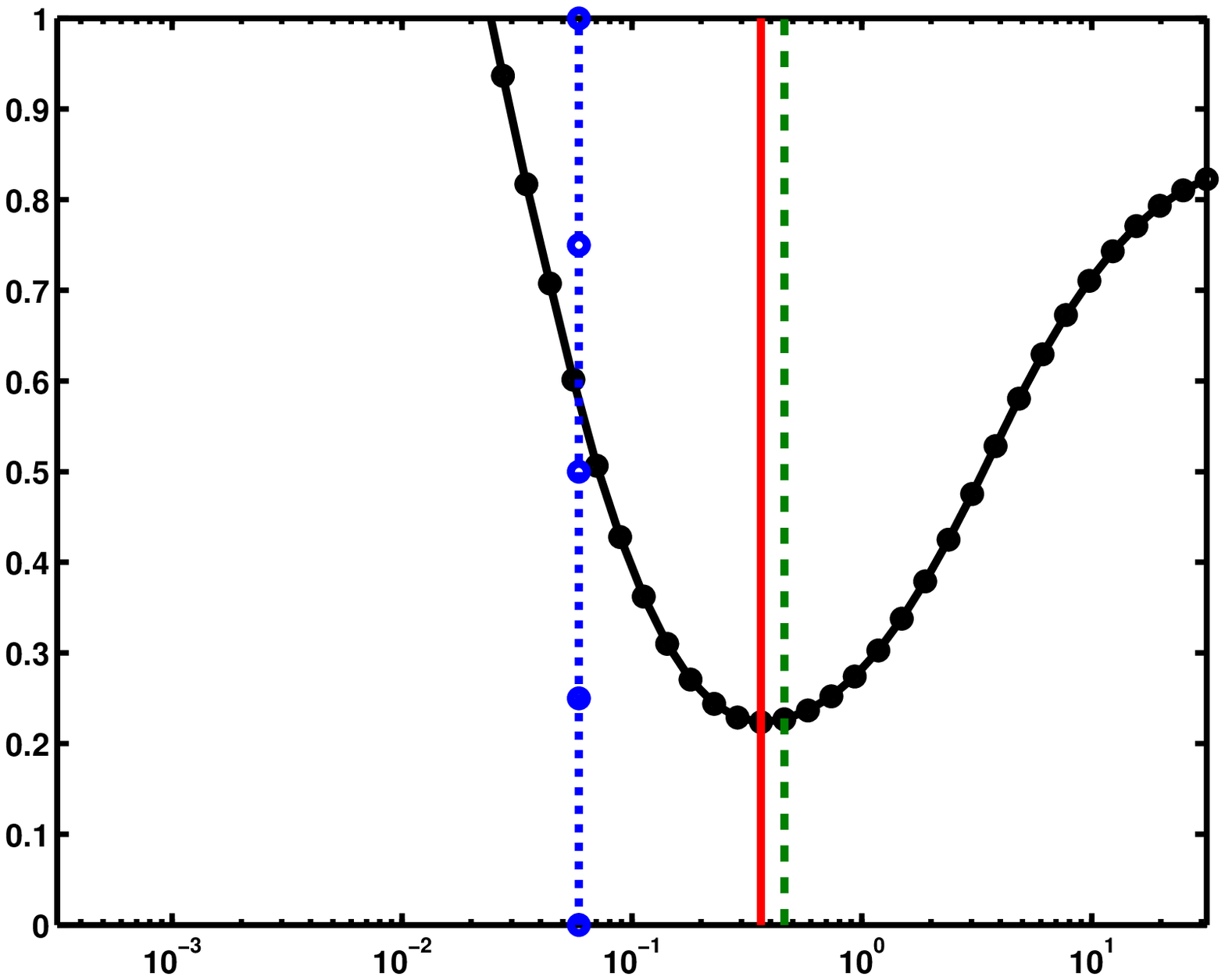}}
\subfigure[$L=L_1$]{\includegraphics[width=1.7in]{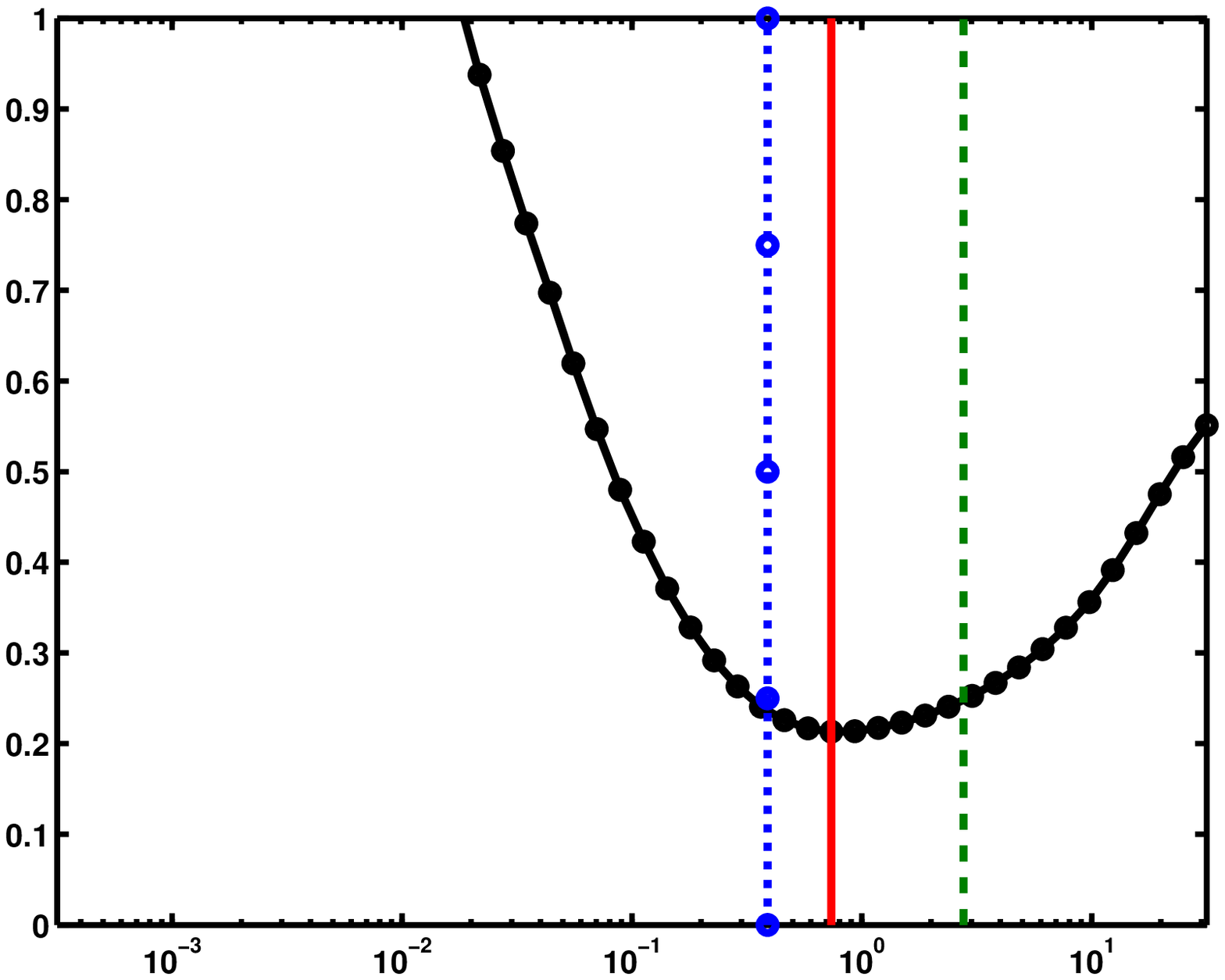}}
\subfigure[$L=L_2$]{\includegraphics[width=1.7in]{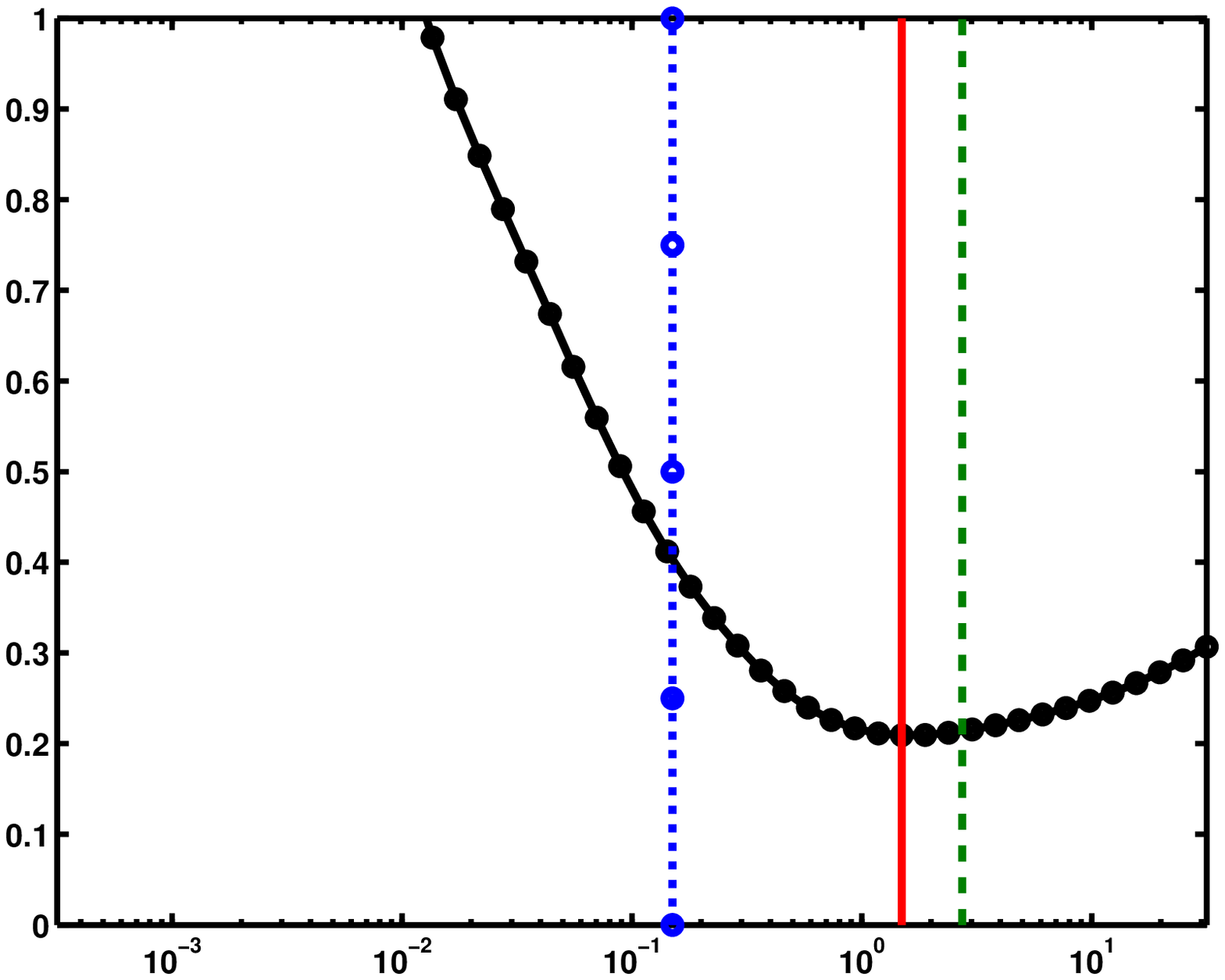}}
\subfigure[$L=I$]{\includegraphics[width=1.7in]{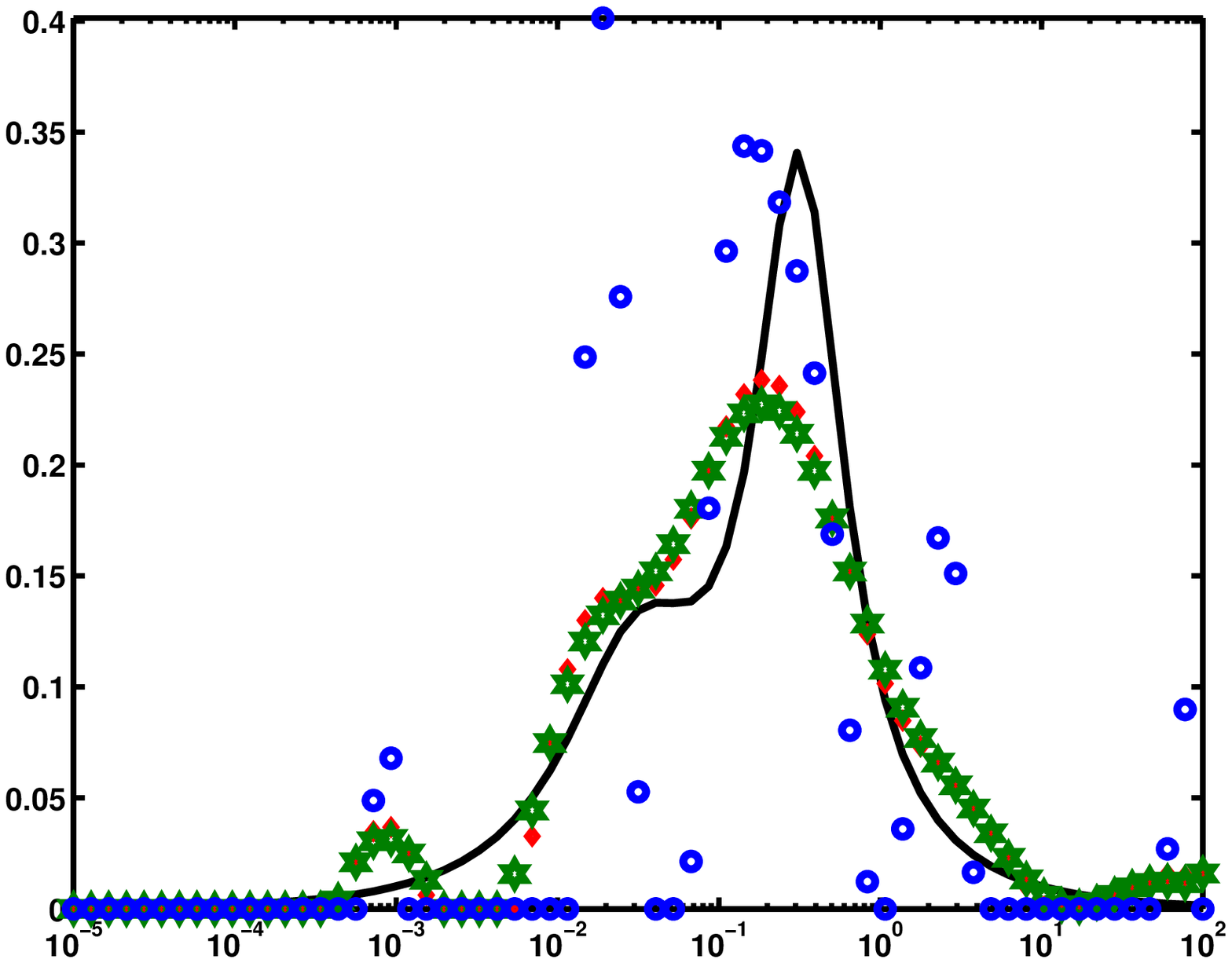}}
\subfigure[$L=L_1$]{\includegraphics[width=1.7in]{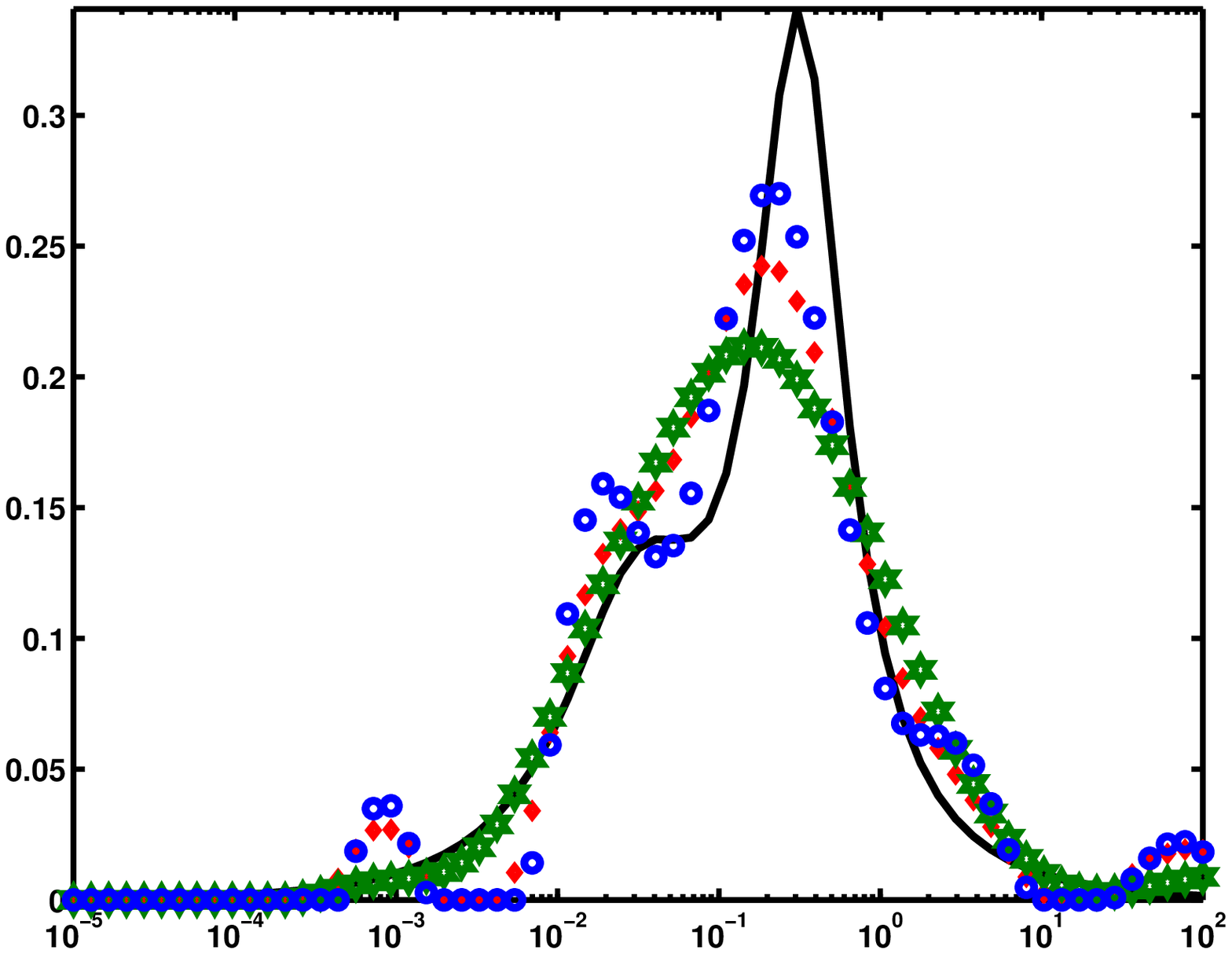}}
\subfigure[$L=L_2$]{\includegraphics[width=1.7in]{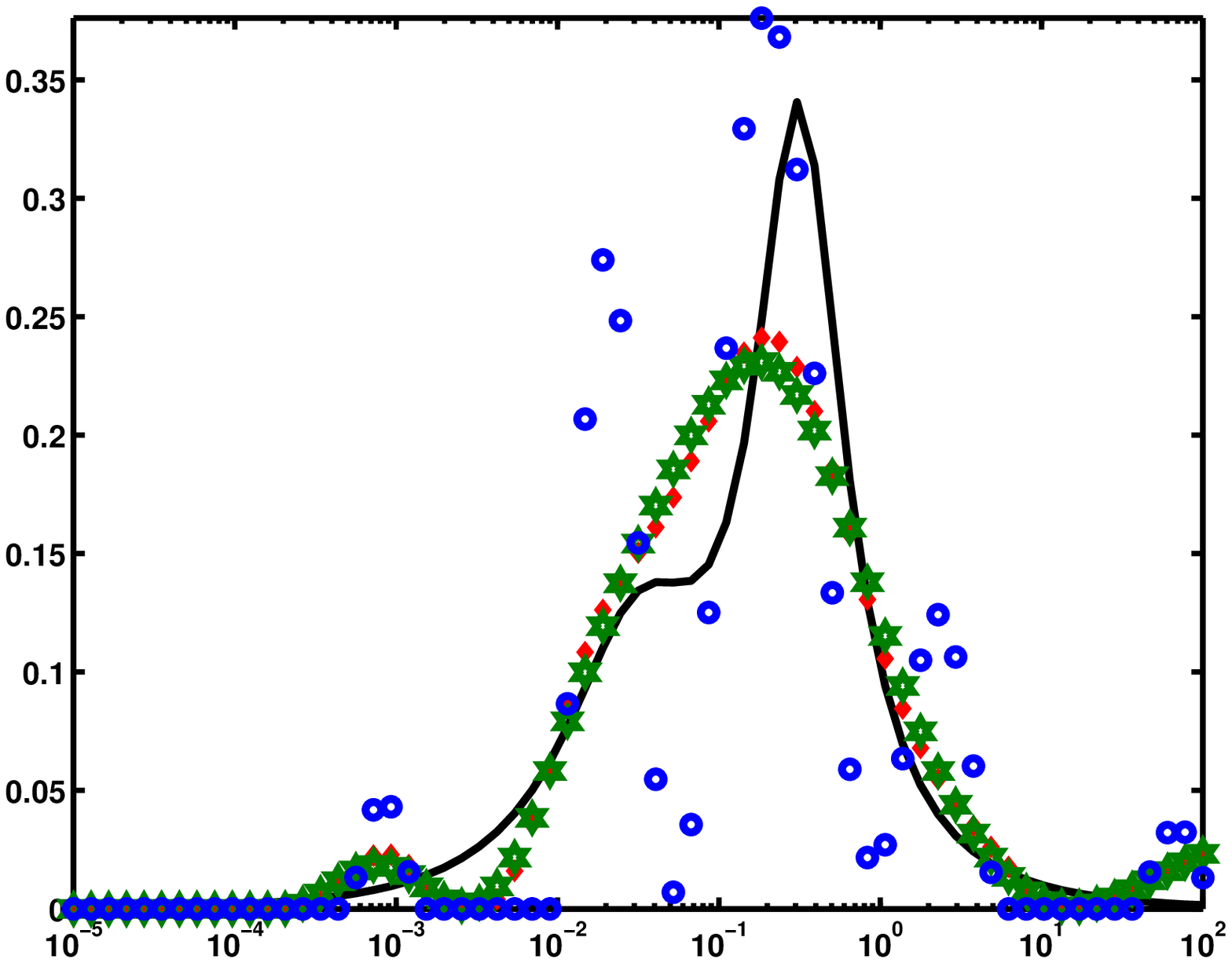}}
\caption{NNLS solutions of RQ-B matrix $A_3$. Noise level $5\%$.}
\label{hnfig-lambdachoiceRQ5A3HN}
\end{figure}

\begin{figure}[!ht]
\centering
\subfigure[$L=I$]{\includegraphics[width=1.7in]{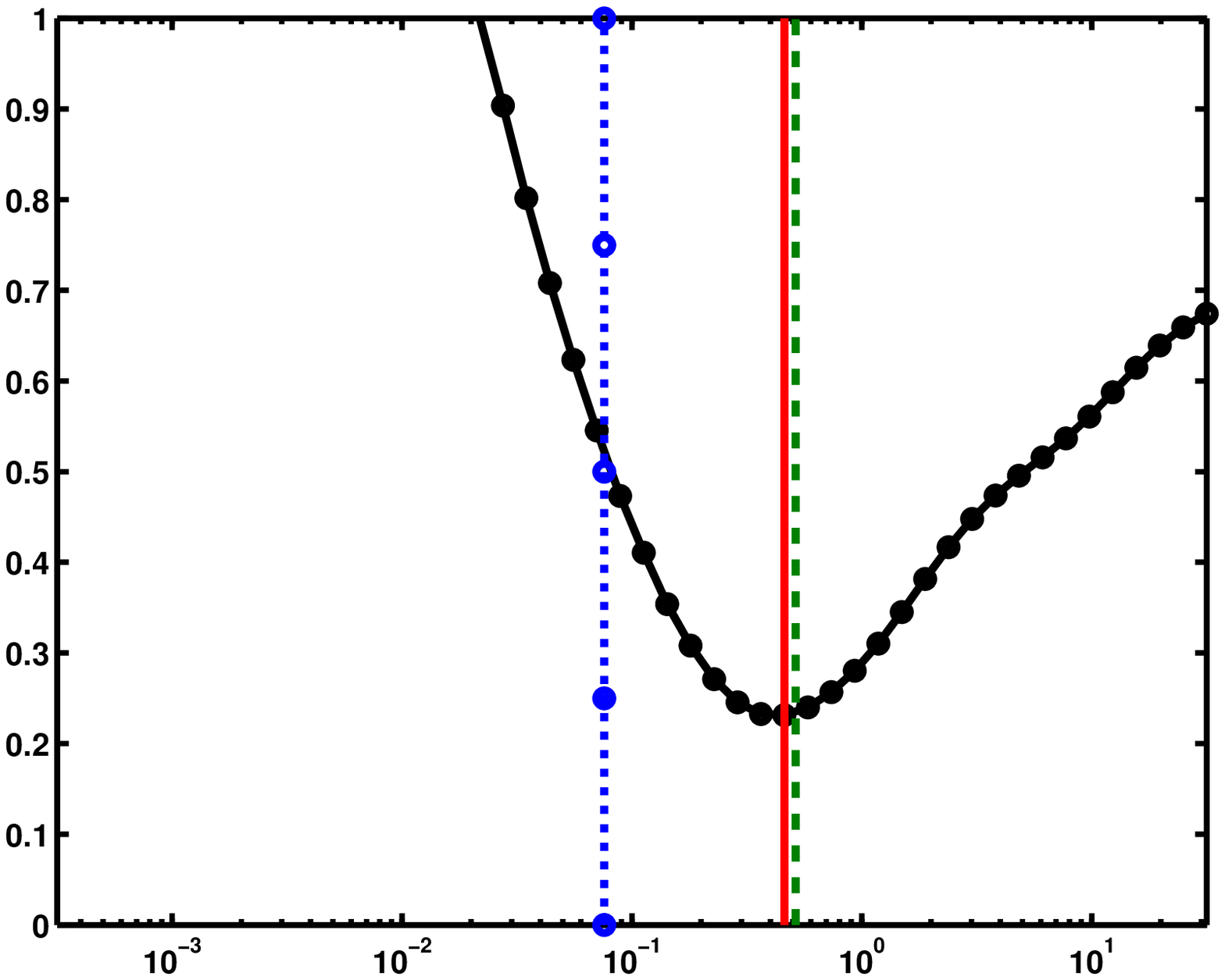}}
\subfigure[$L=L_1$]{\includegraphics[width=1.7in]{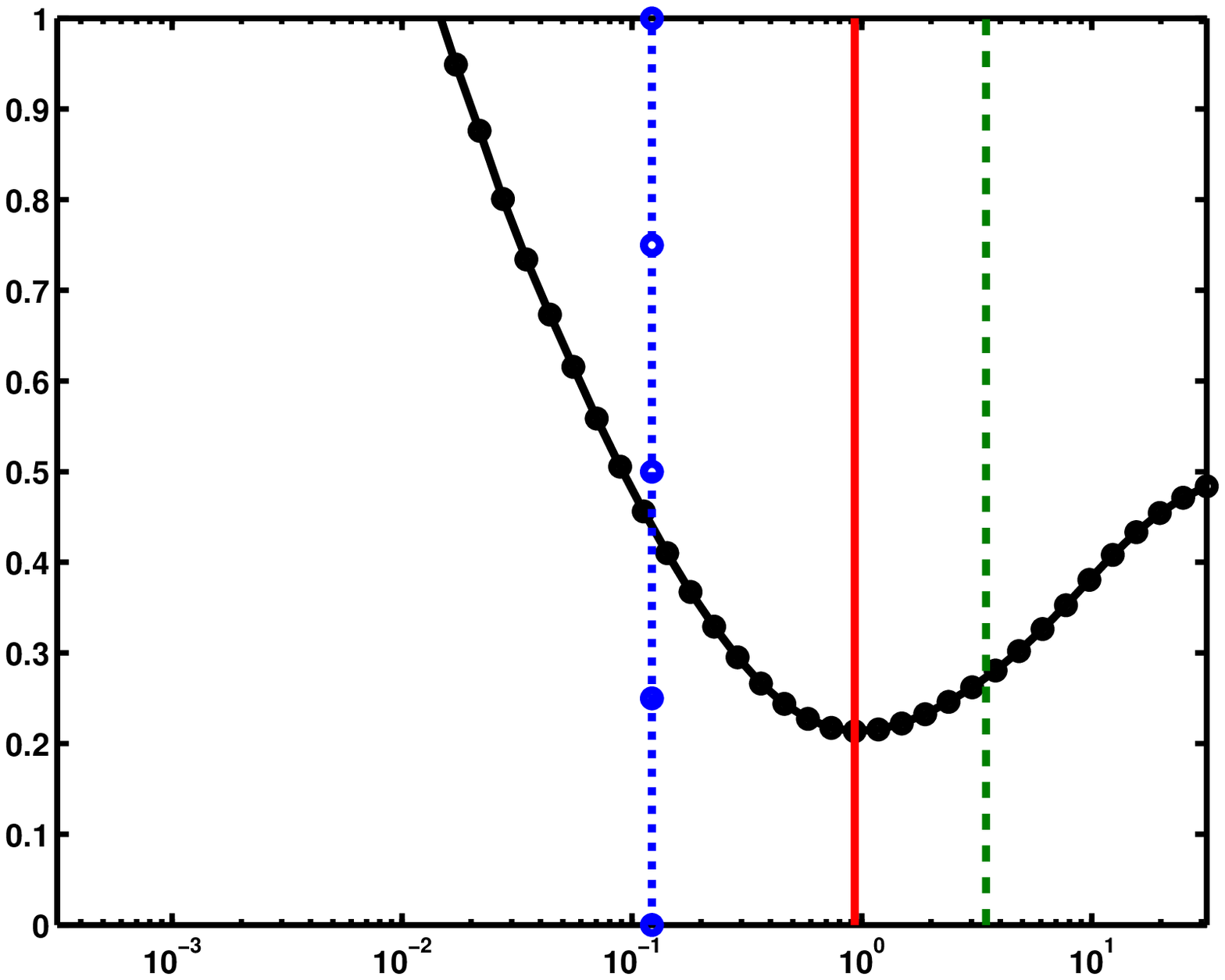}}
\subfigure[$L=L_2$]{\includegraphics[width=1.7in]{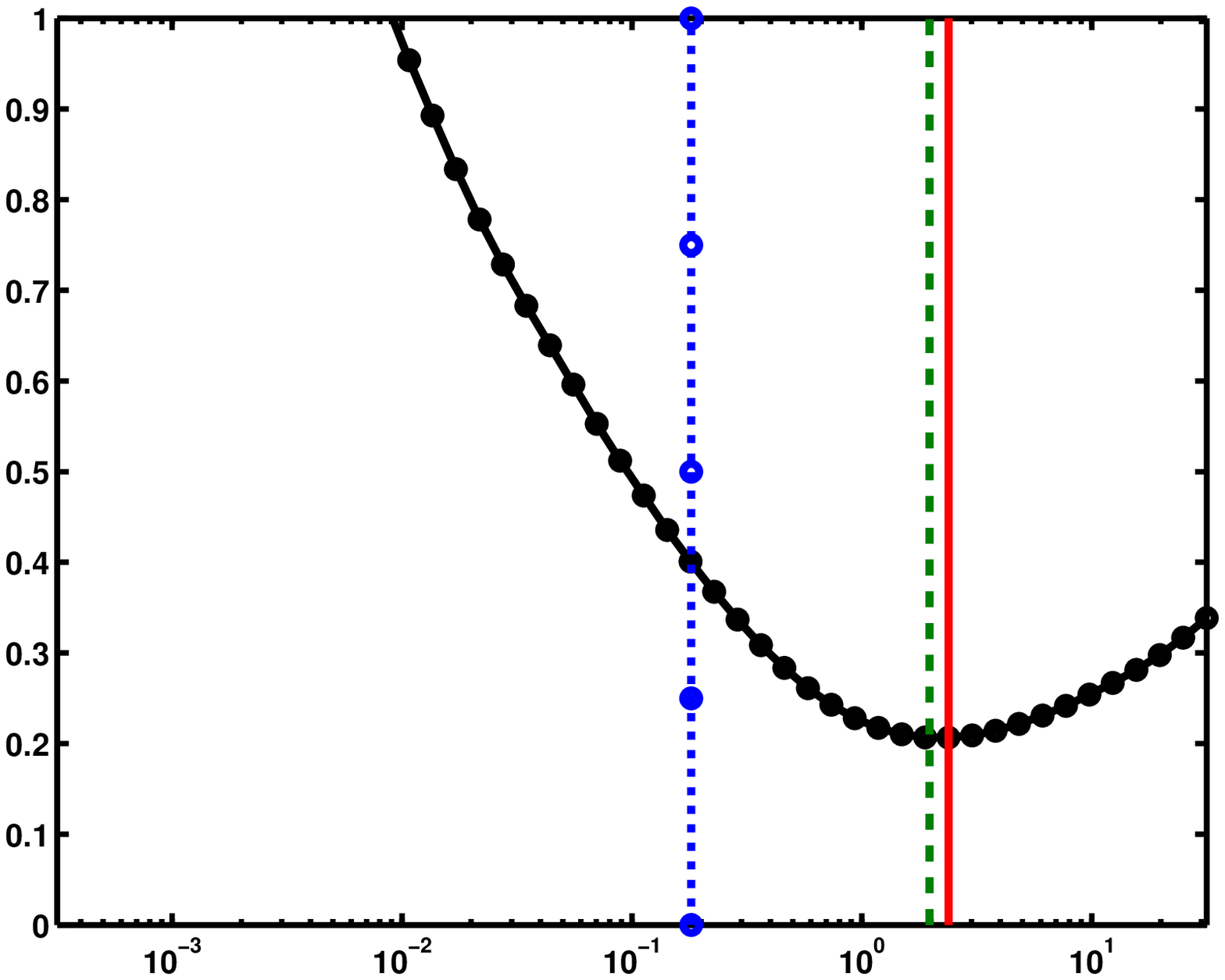}}
\subfigure[$L=I$]{\includegraphics[width=1.7in]{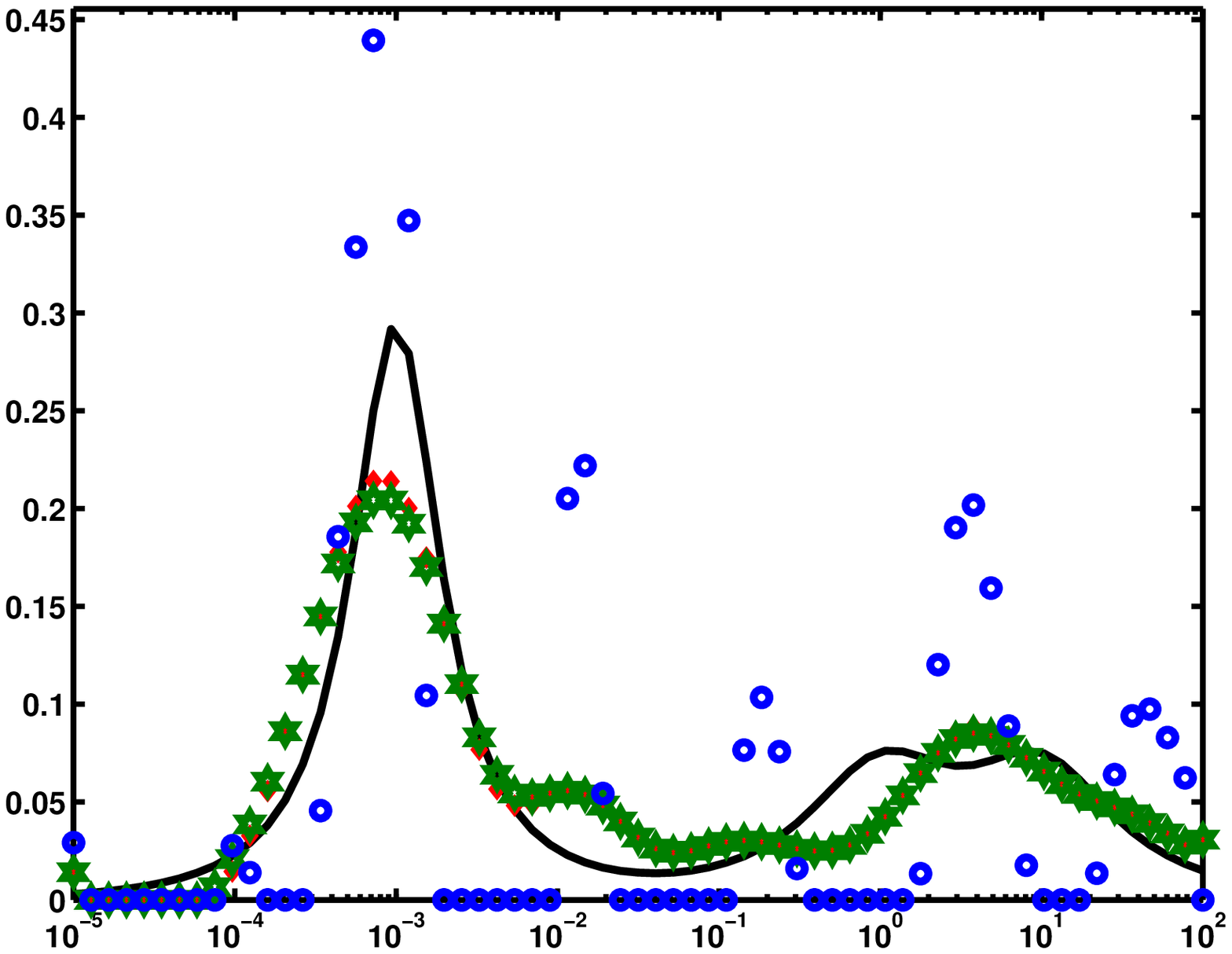}}
\subfigure[$L=L_1$]{\includegraphics[width=1.7in]{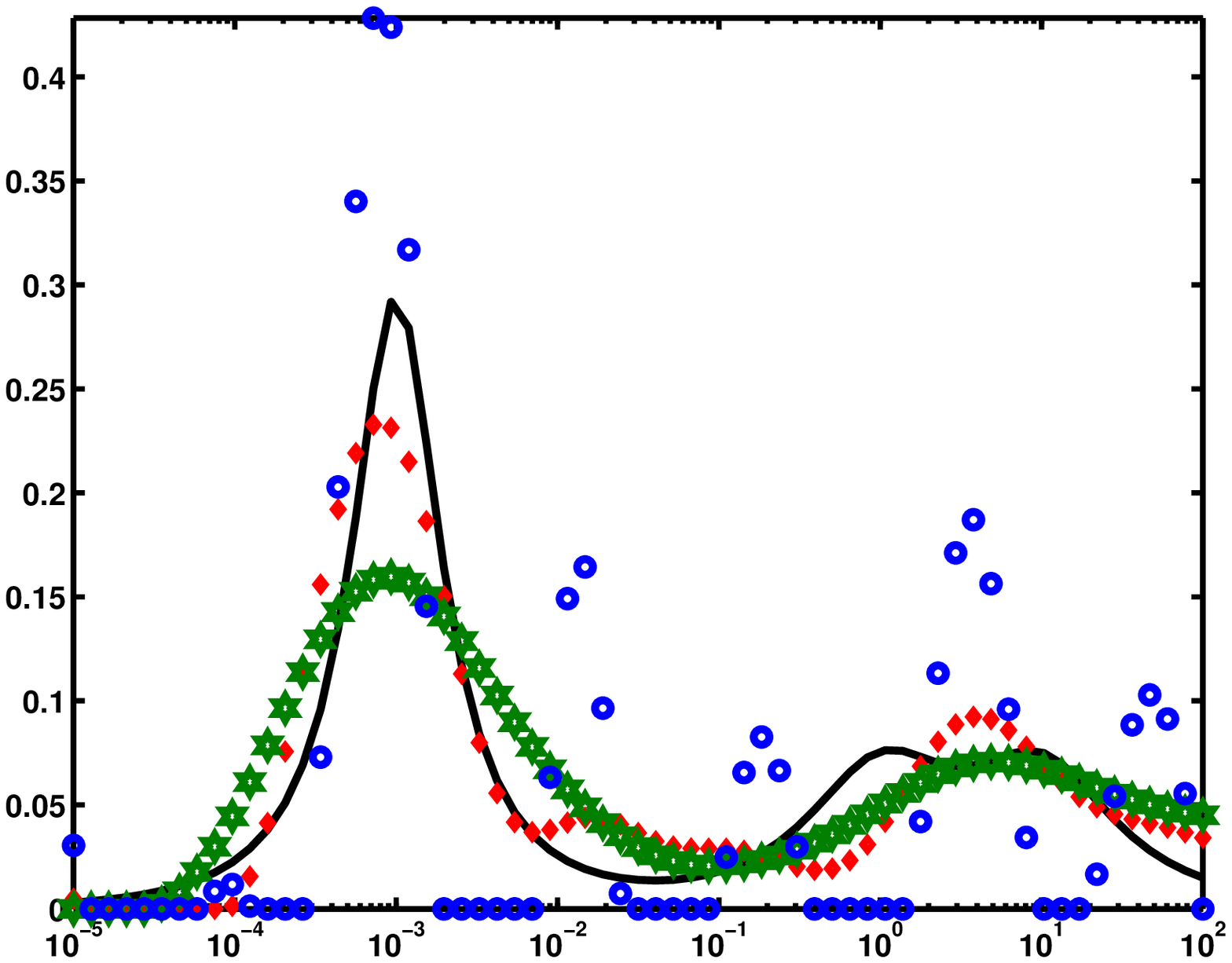}}
\subfigure[$L=L_2$]{\includegraphics[width=1.7in]{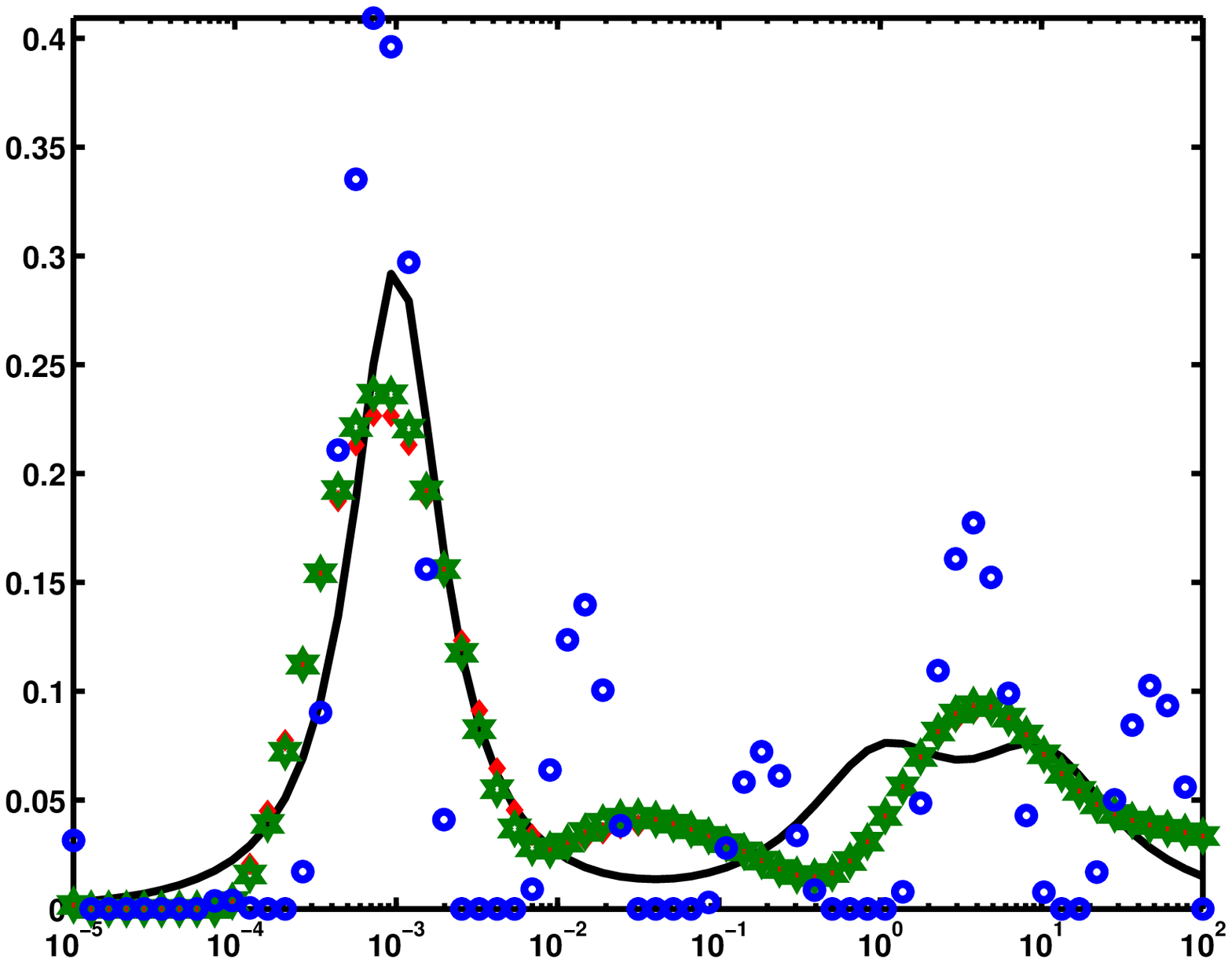}}
\caption{NNLS solutions of RQ-C matrix $A_3$. Noise level $5\%$.}
\label{hnfig-lambdachoiceRQ6A3HN}
\end{figure}

\begin{figure}[!ht]
\centering
\subfigure[$L=I$]{\includegraphics[width=1.7in]{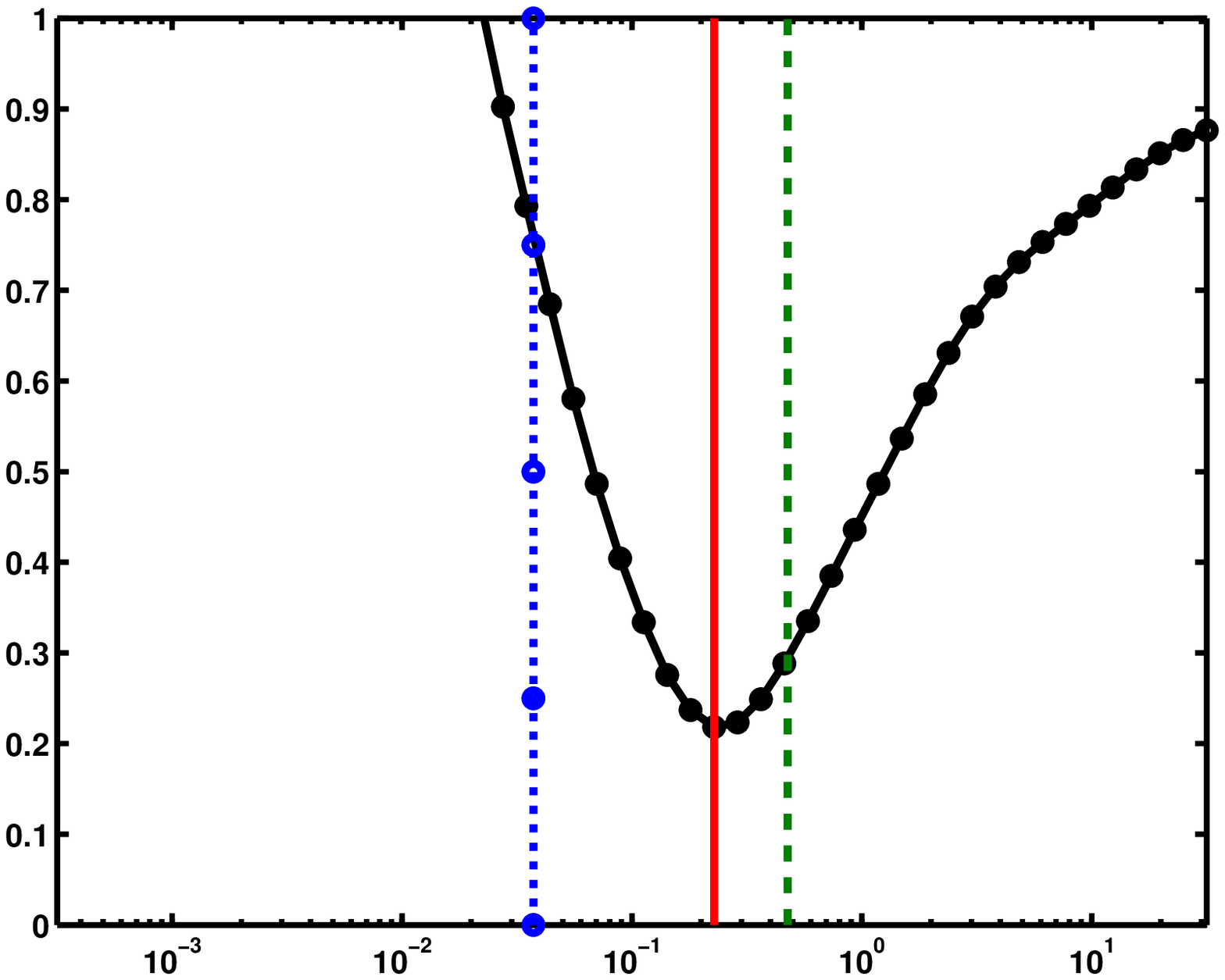}}
\subfigure[$L=L_1$]{\includegraphics[width=1.7in]{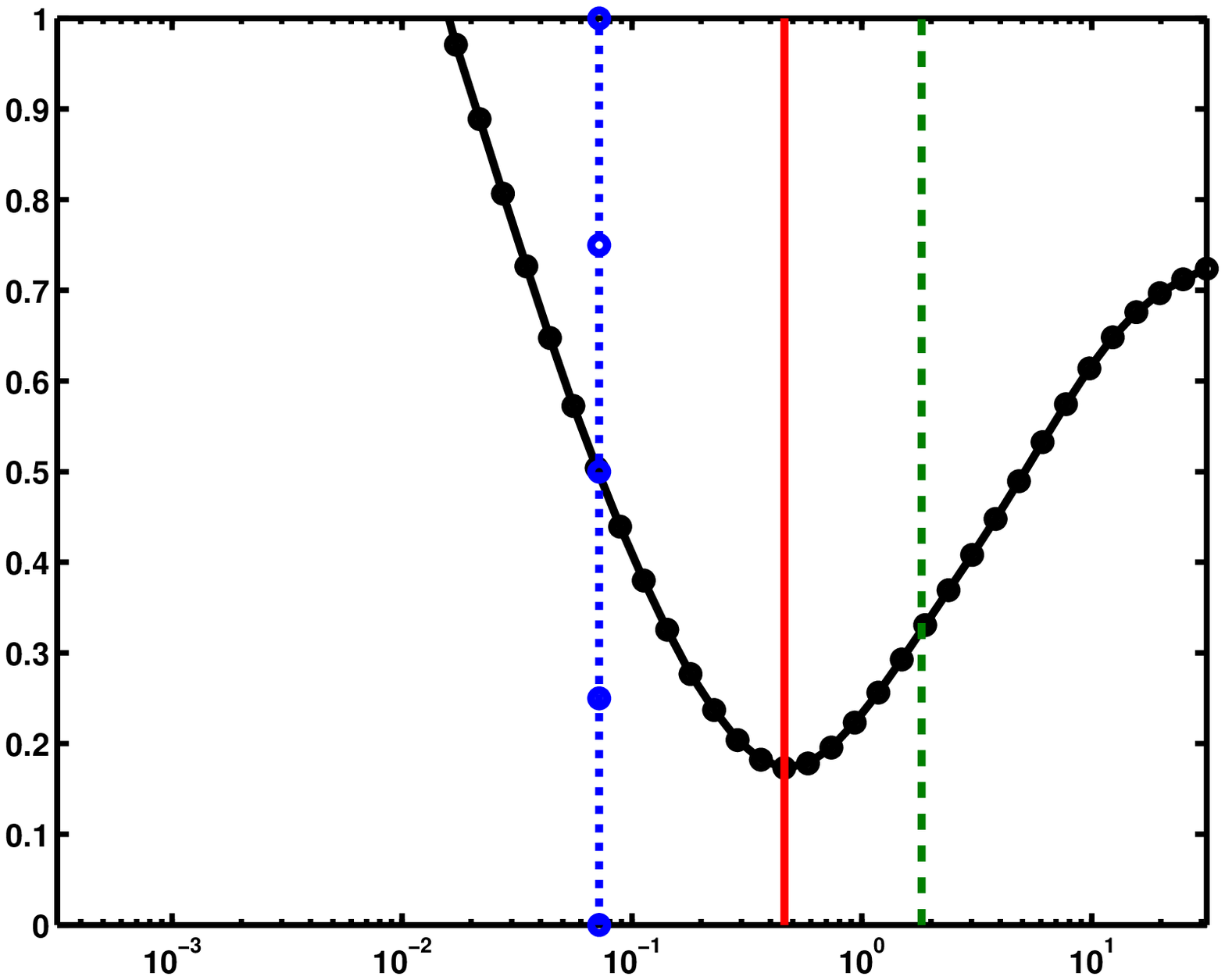}}
\subfigure[$L=L_2$]{\includegraphics[width=1.7in]{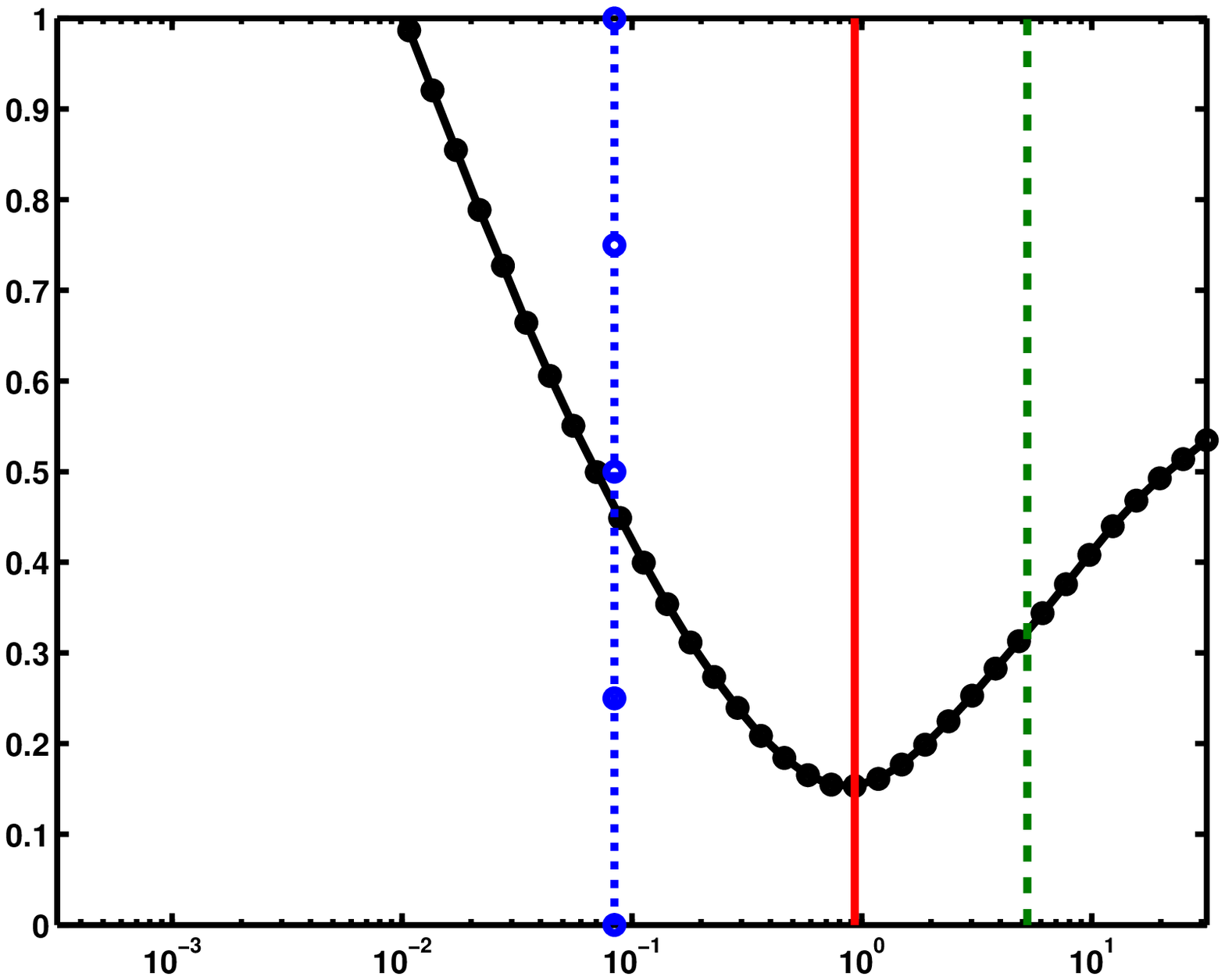}}
\subfigure[$L=I$]{\includegraphics[width=1.7in]{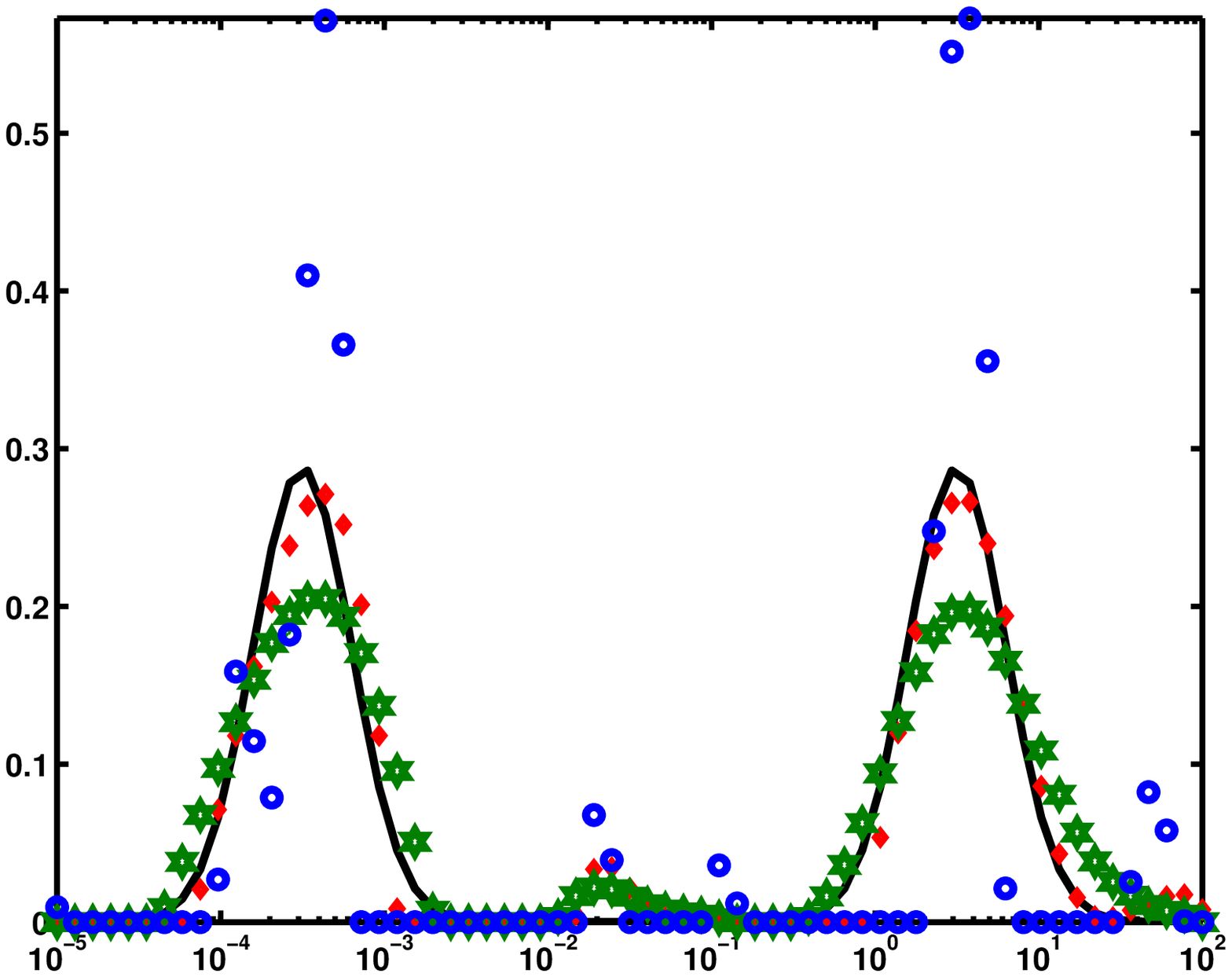}}
\subfigure[$L=L_1$]{\includegraphics[width=1.7in]{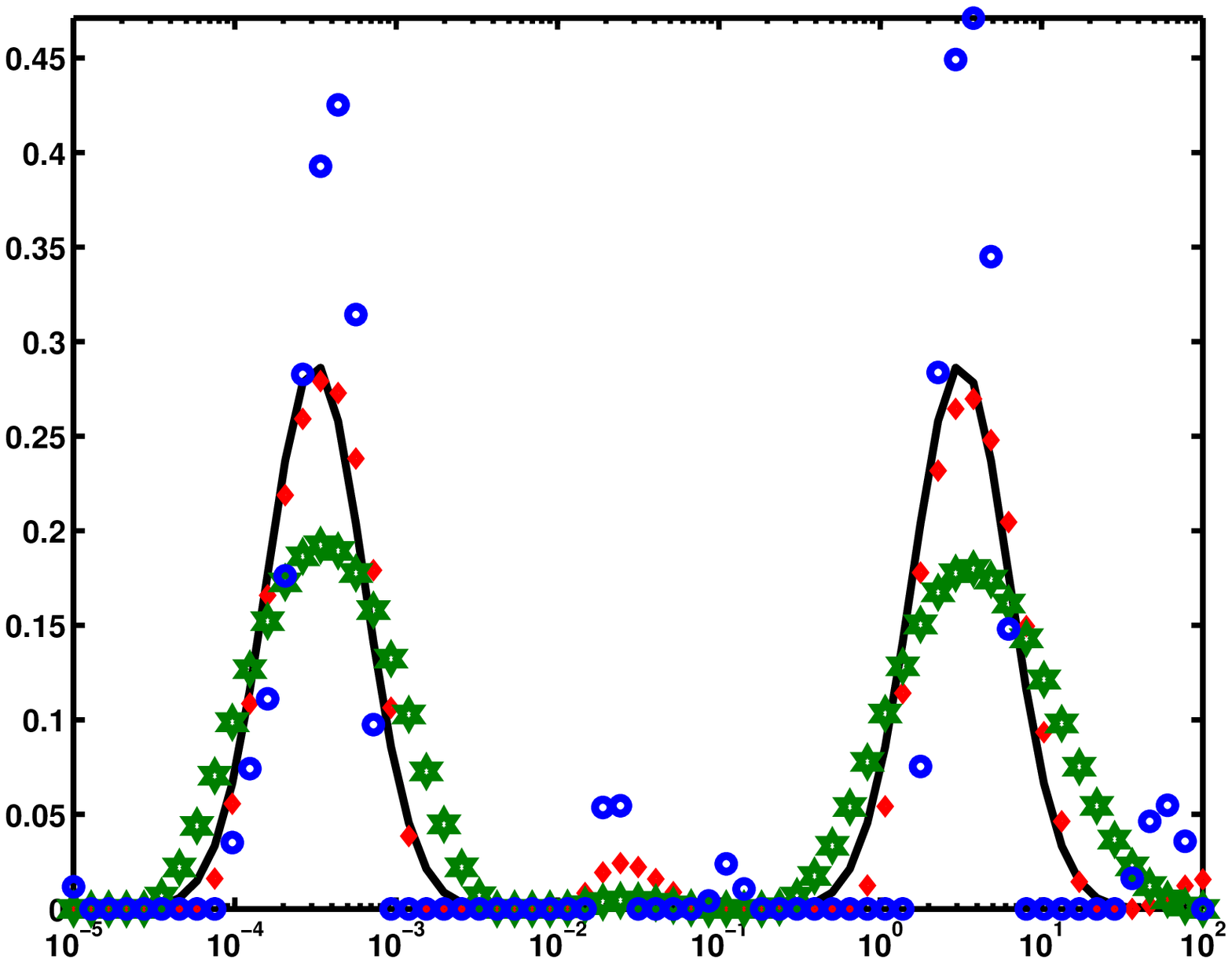}}
\subfigure[$L=L_2$]{\includegraphics[width=1.7in]{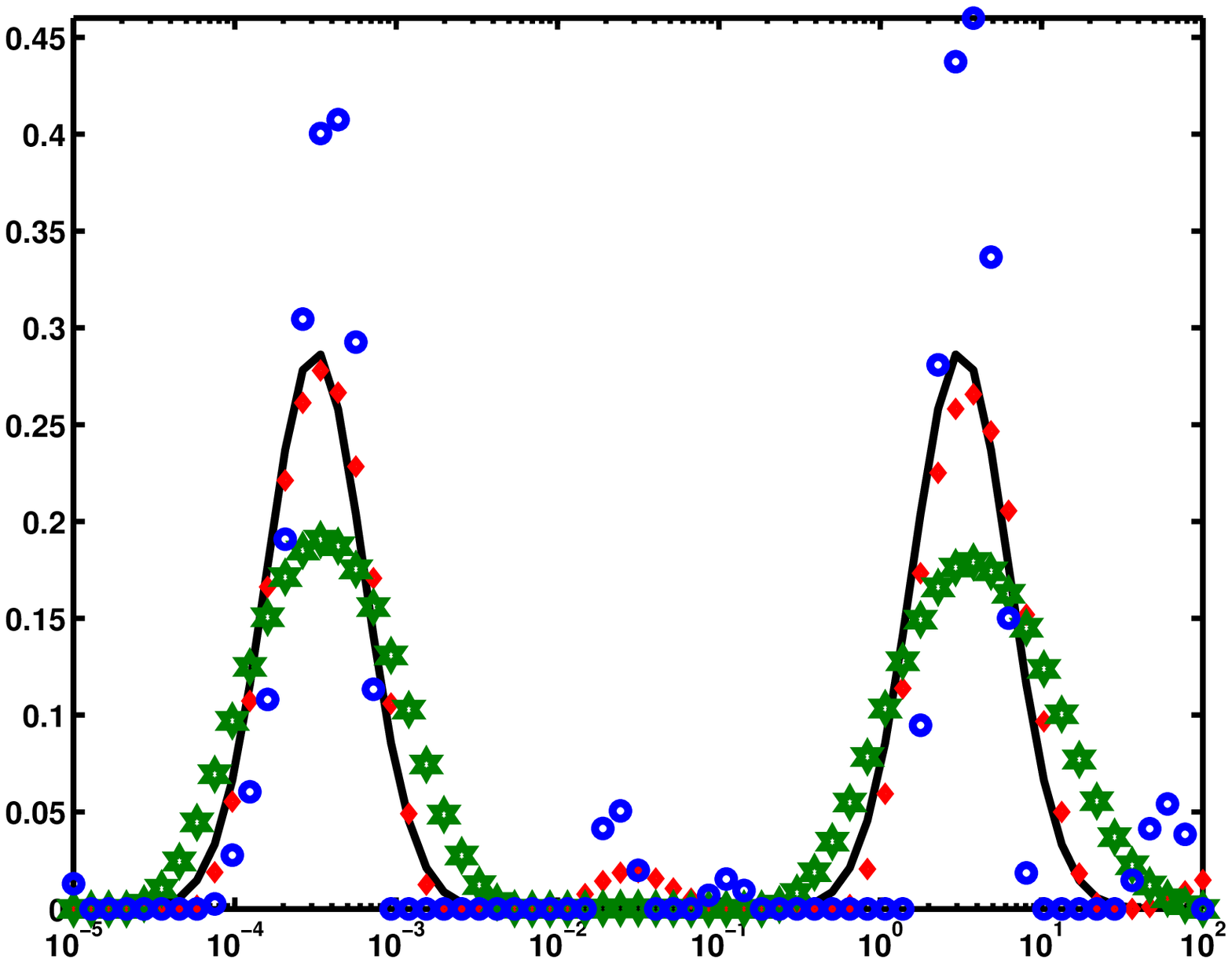}}
\caption{NNLS solutions of LN-A matrix $A_3$. Noise level $5\%$.}
\label{hnfig-lambdachoiceLN2A3HN}
\end{figure}

\begin{figure}[!ht]
\centering
\subfigure[$L=I$]{\includegraphics[width=1.7in]{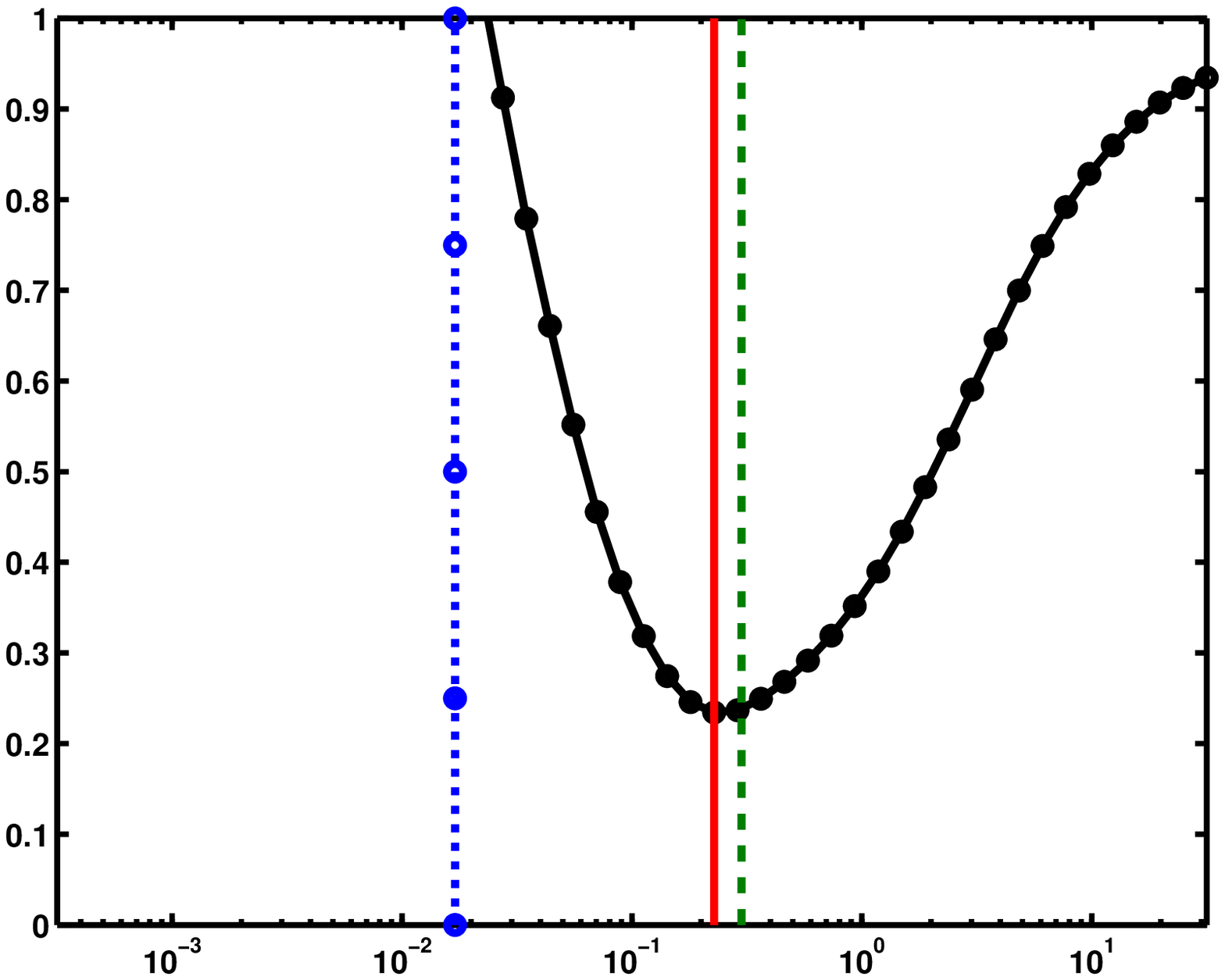}}
\subfigure[$L=L_1$]{\includegraphics[width=1.7in]{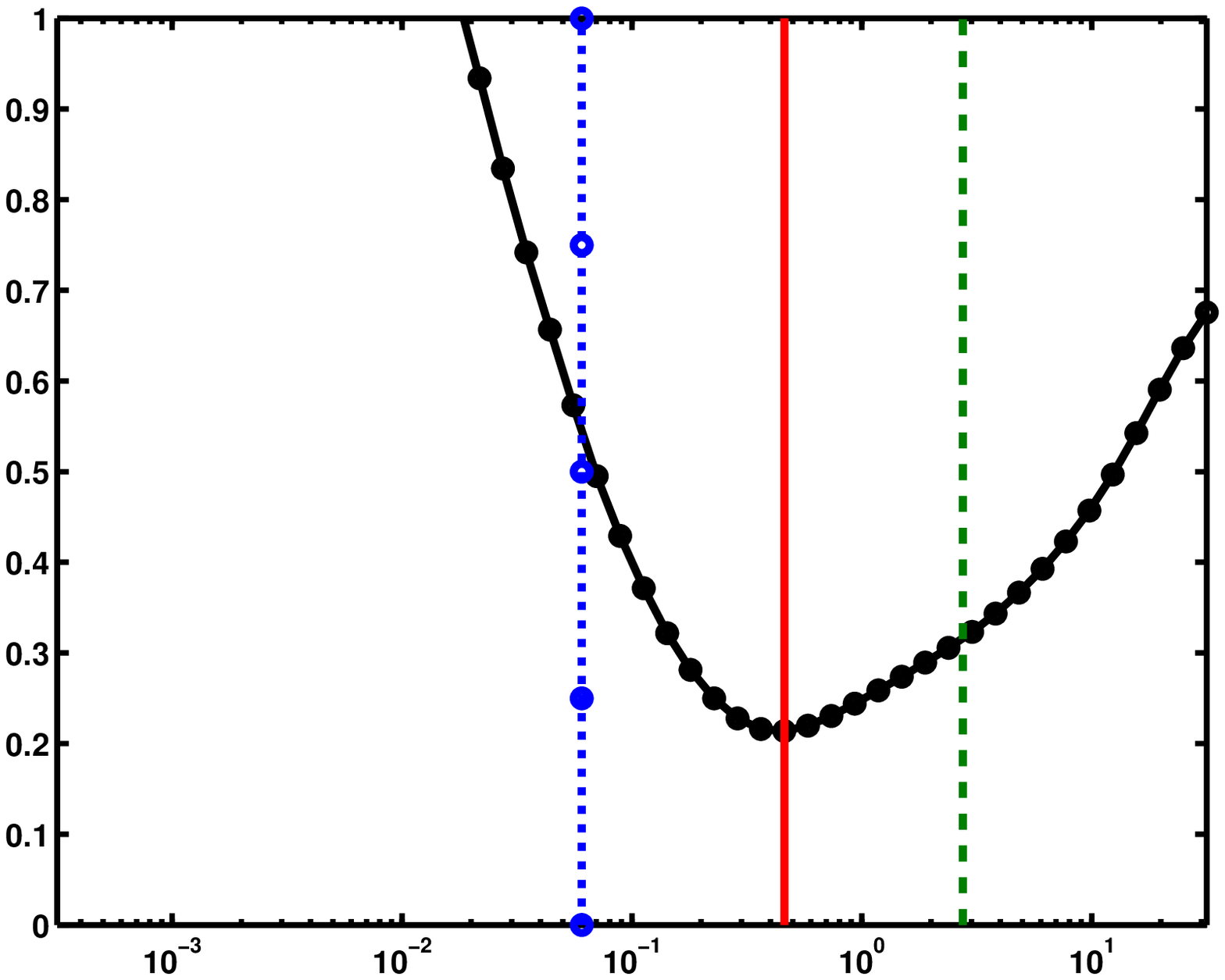}}
\subfigure[$L=L_2$]{\includegraphics[width=1.7in]{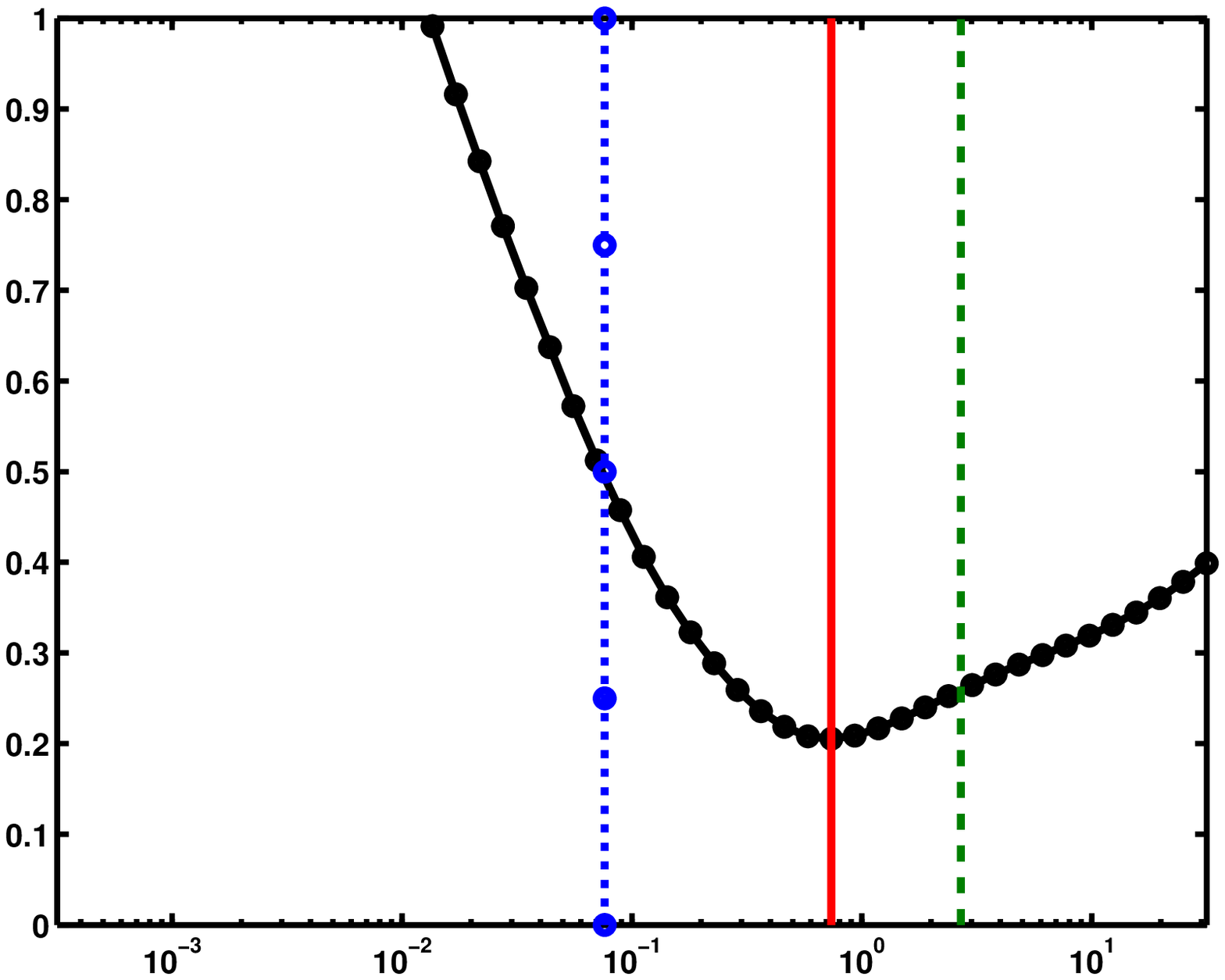}}
\subfigure[$L=I$]{\includegraphics[width=1.7in]{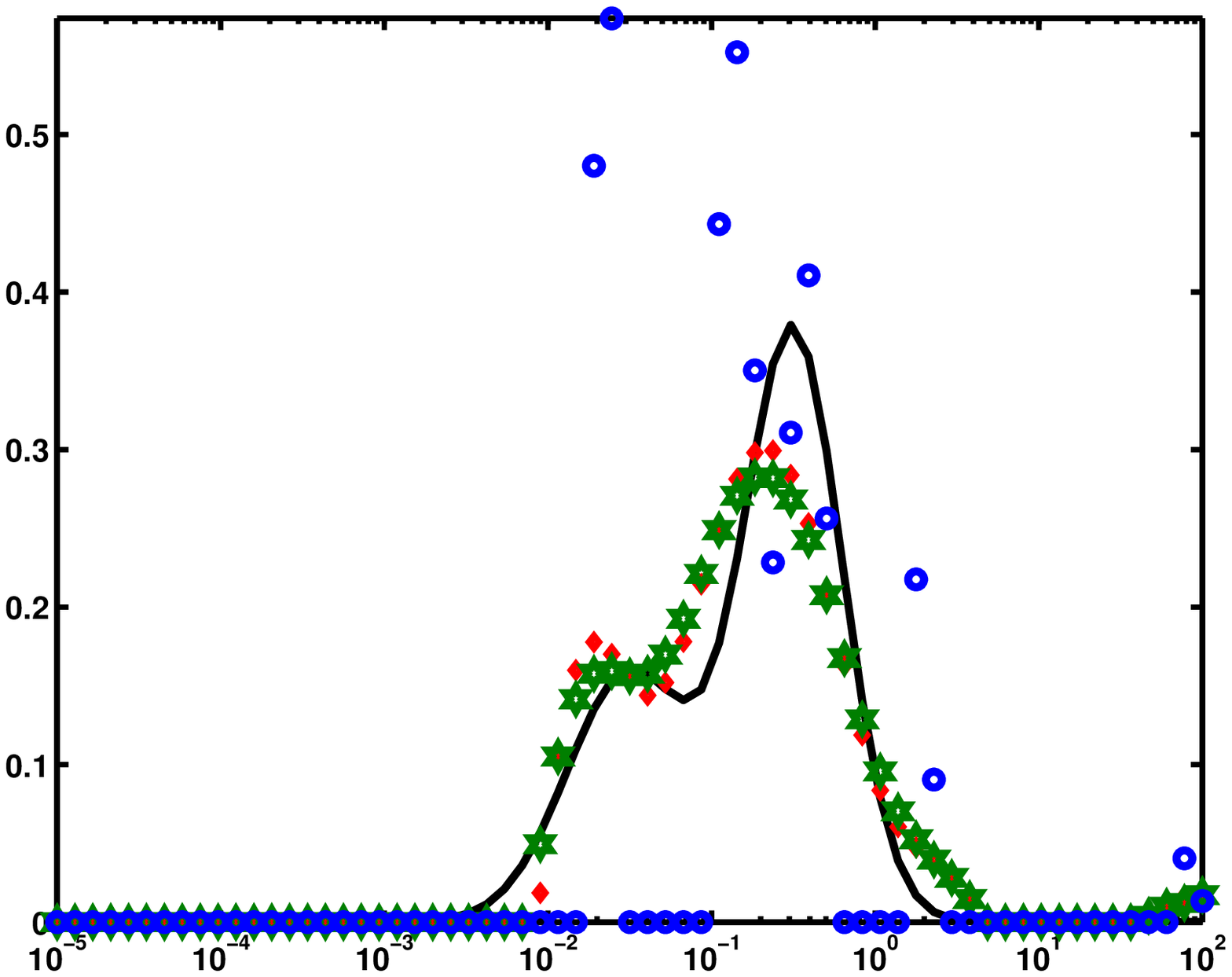}}
\subfigure[$L=L_1$]{\includegraphics[width=1.7in]{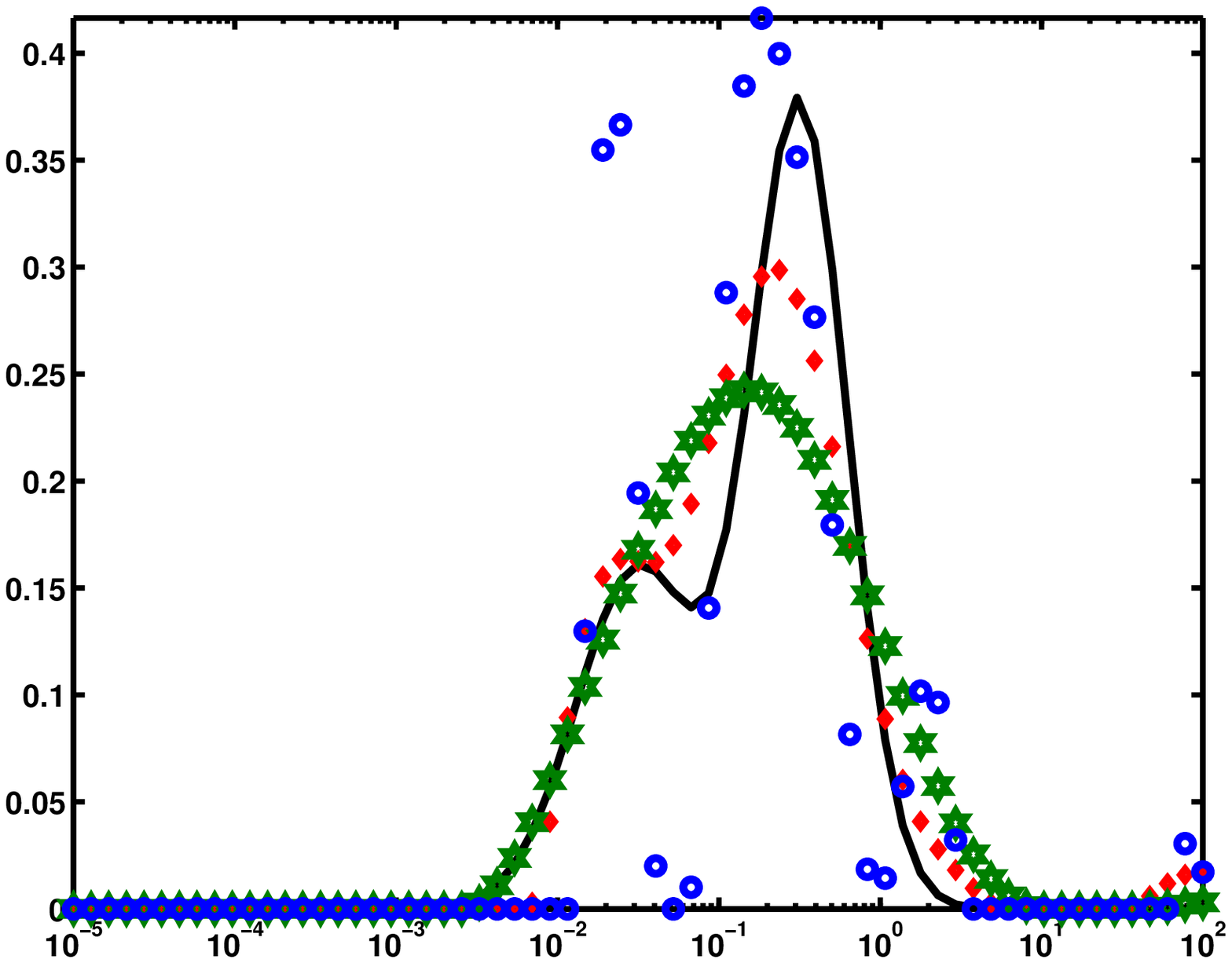}}
\subfigure[$L=L_2$]{\includegraphics[width=1.7in]{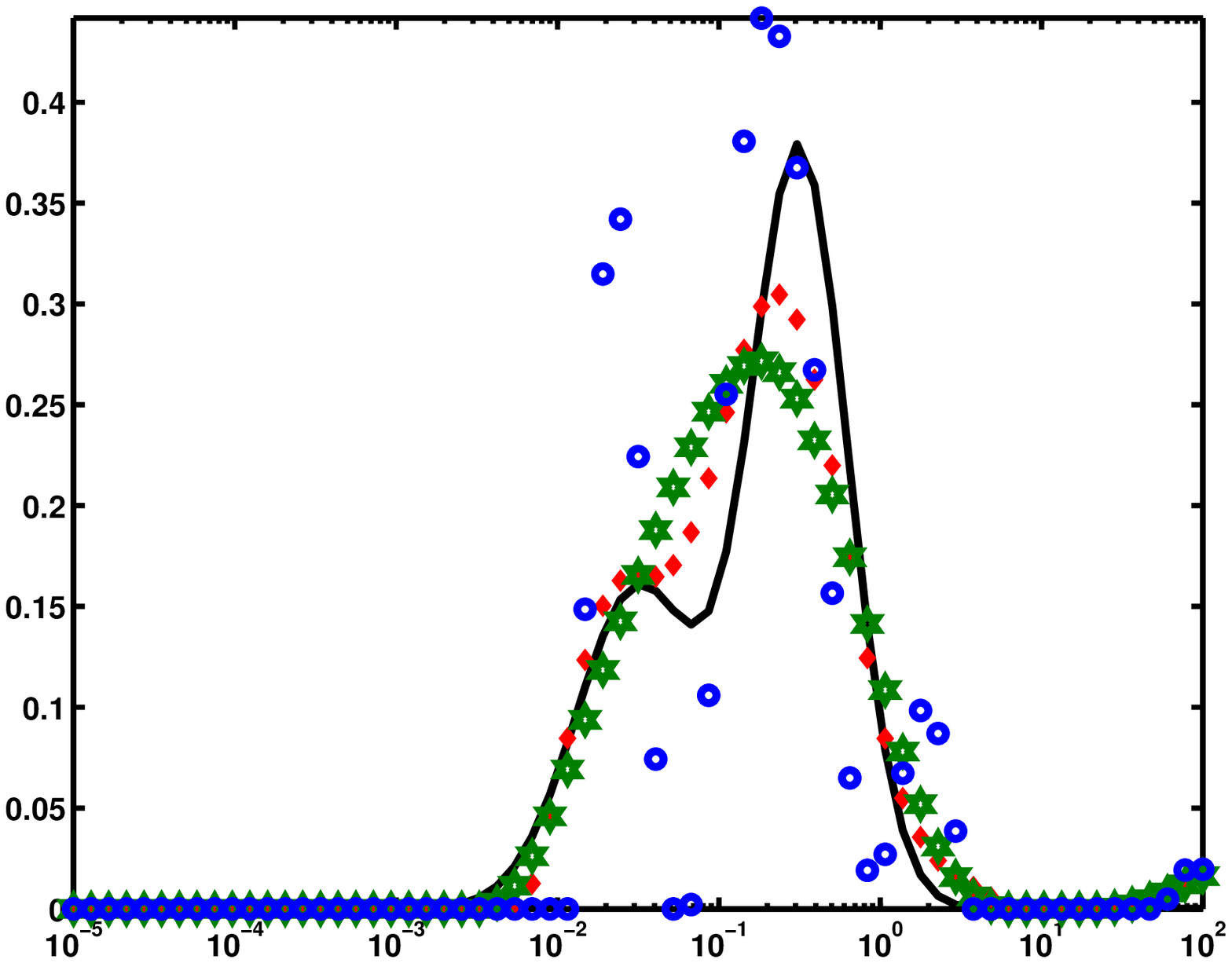}}
\caption{NNLS solutions of LN-B matrix $A_3$. Noise level $5\%$.}
\label{hnfig-lambdachoiceLN5A3HN}
\end{figure}

\begin{figure}[!ht]
\centering
\subfigure[$L=I$]{\includegraphics[width=1.7in]{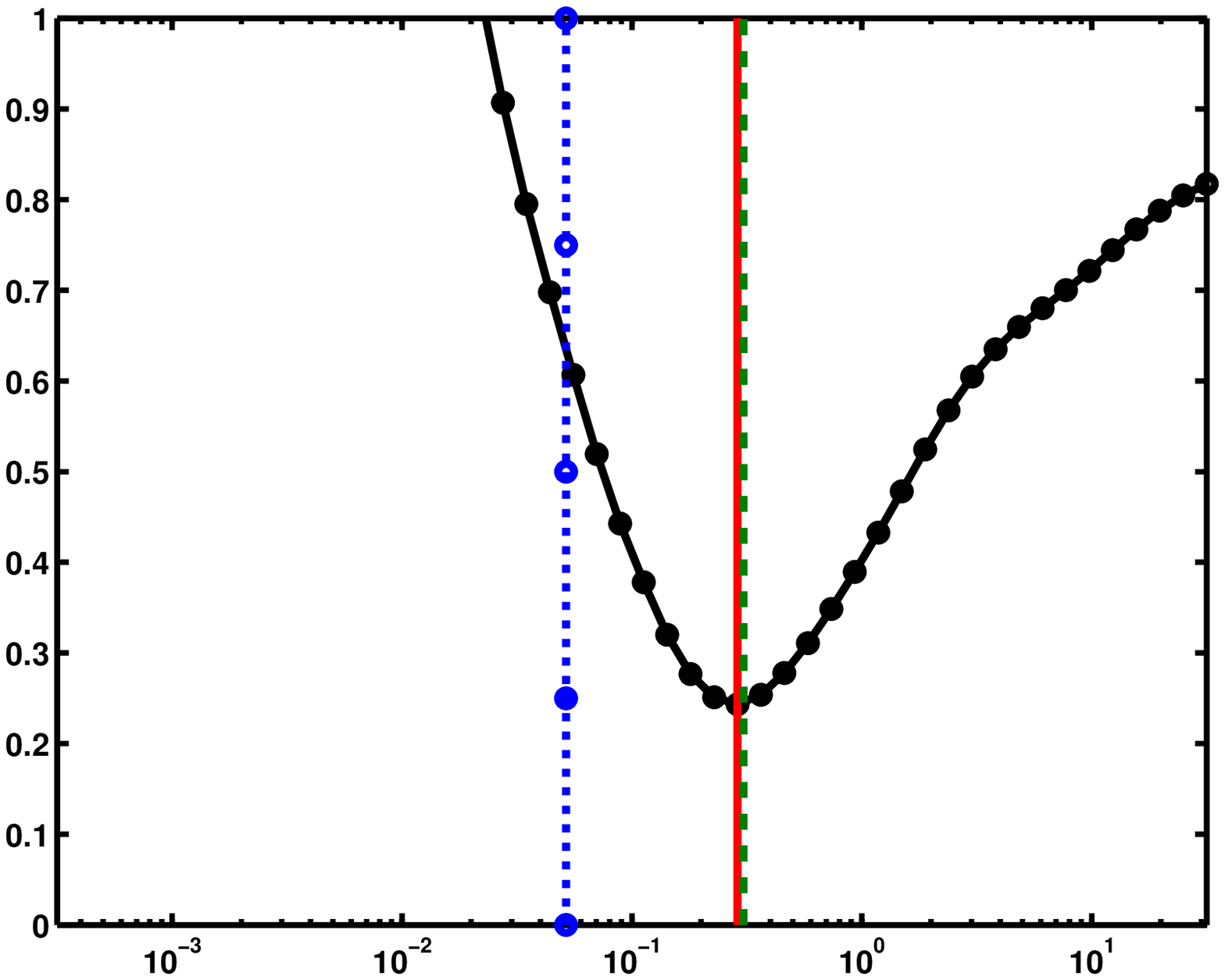}}
\subfigure[$L=L_1$]{\includegraphics[width=1.7in]{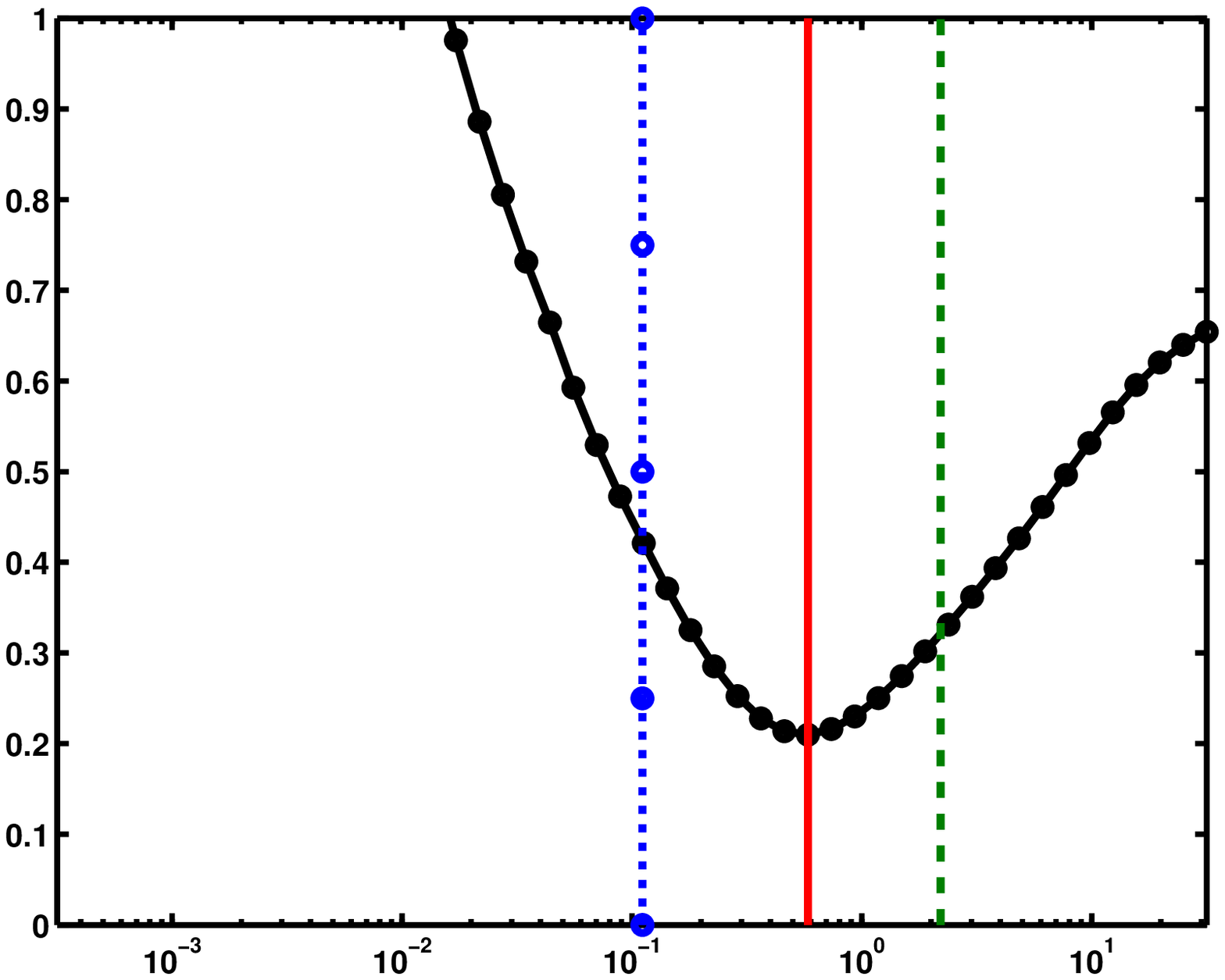}}
\subfigure[$L=L_2$]{\includegraphics[width=1.7in]{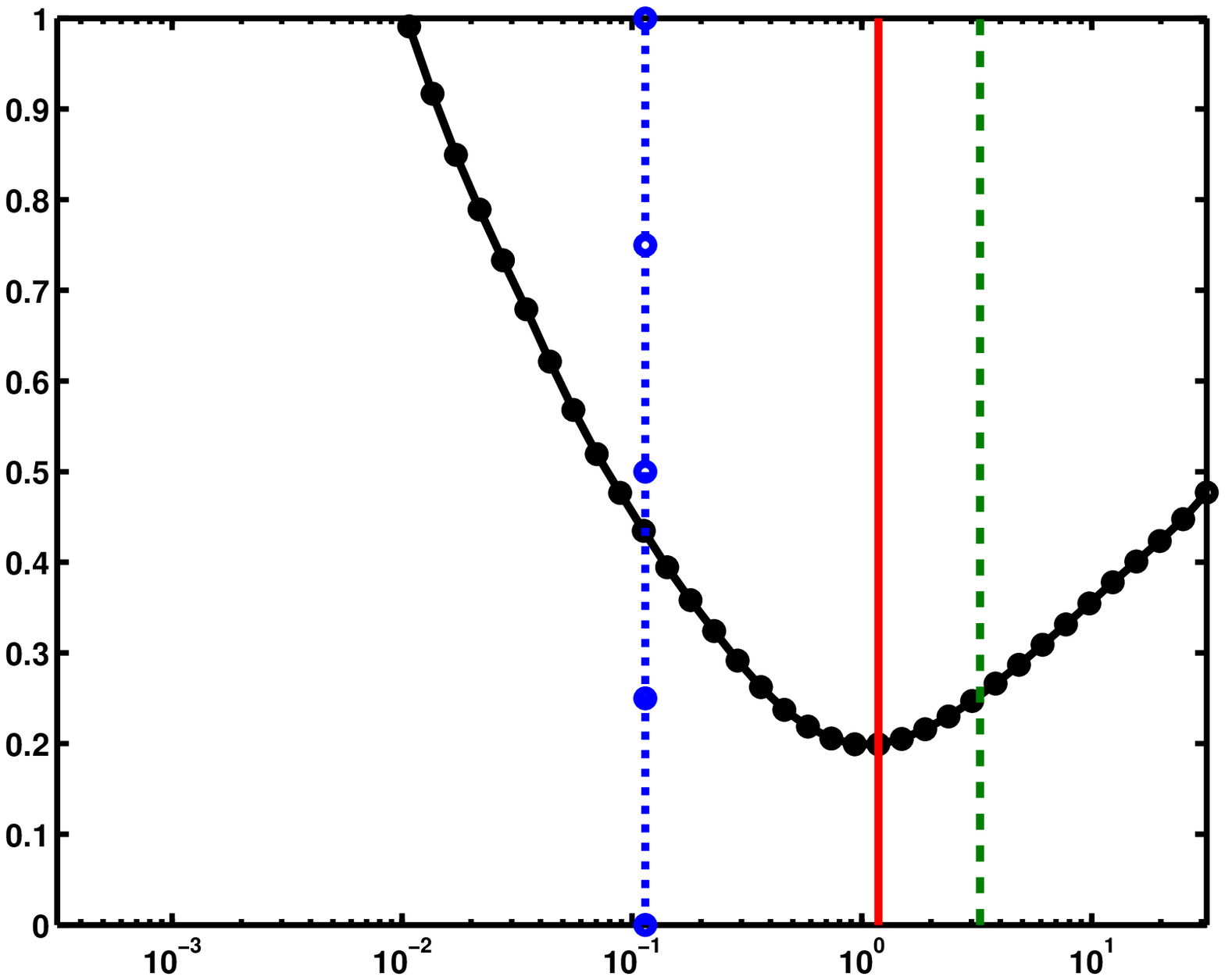}}
\subfigure[$L=I$]{\includegraphics[width=1.7in]{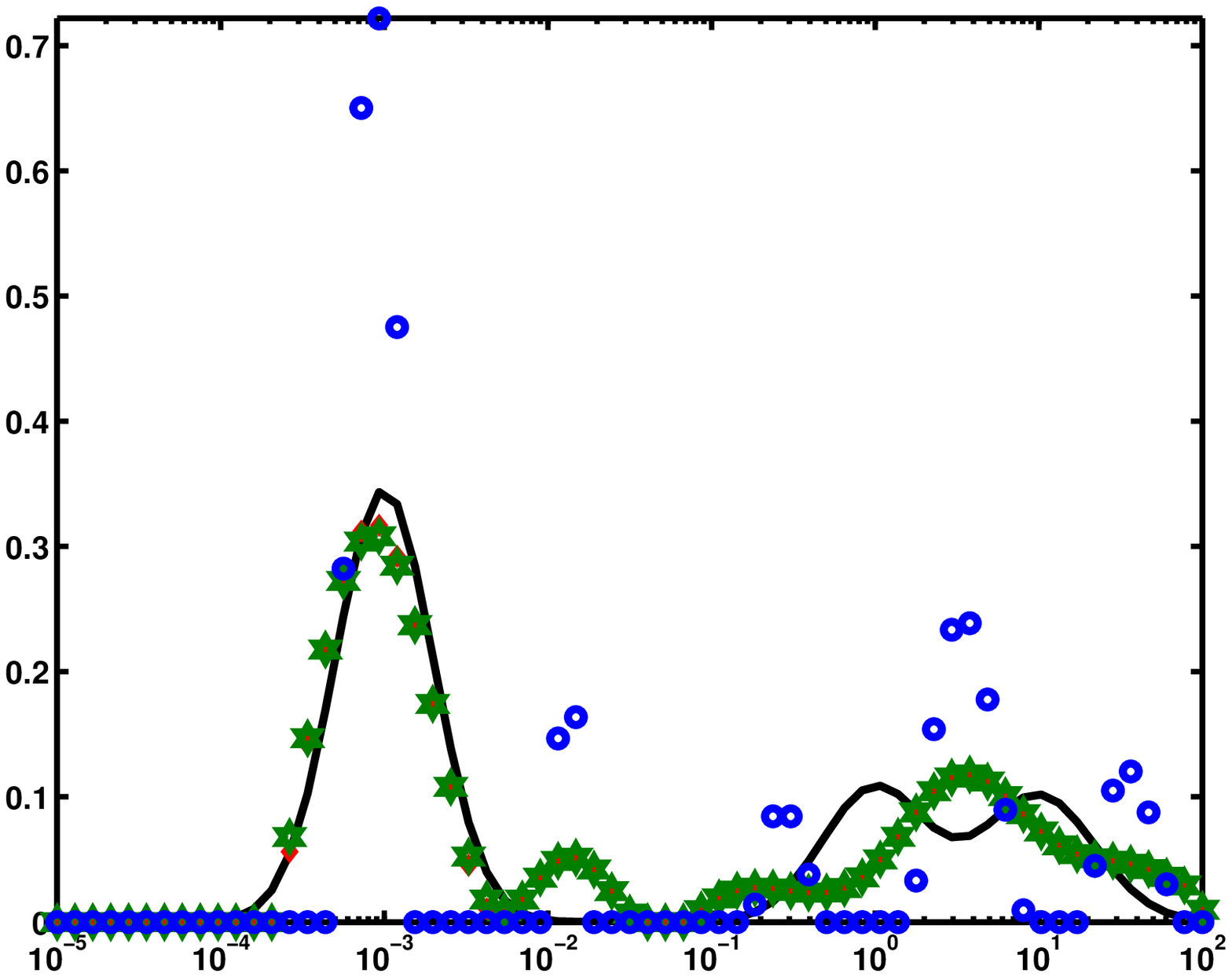}}
\subfigure[$L=L_1$]{\includegraphics[width=1.7in]{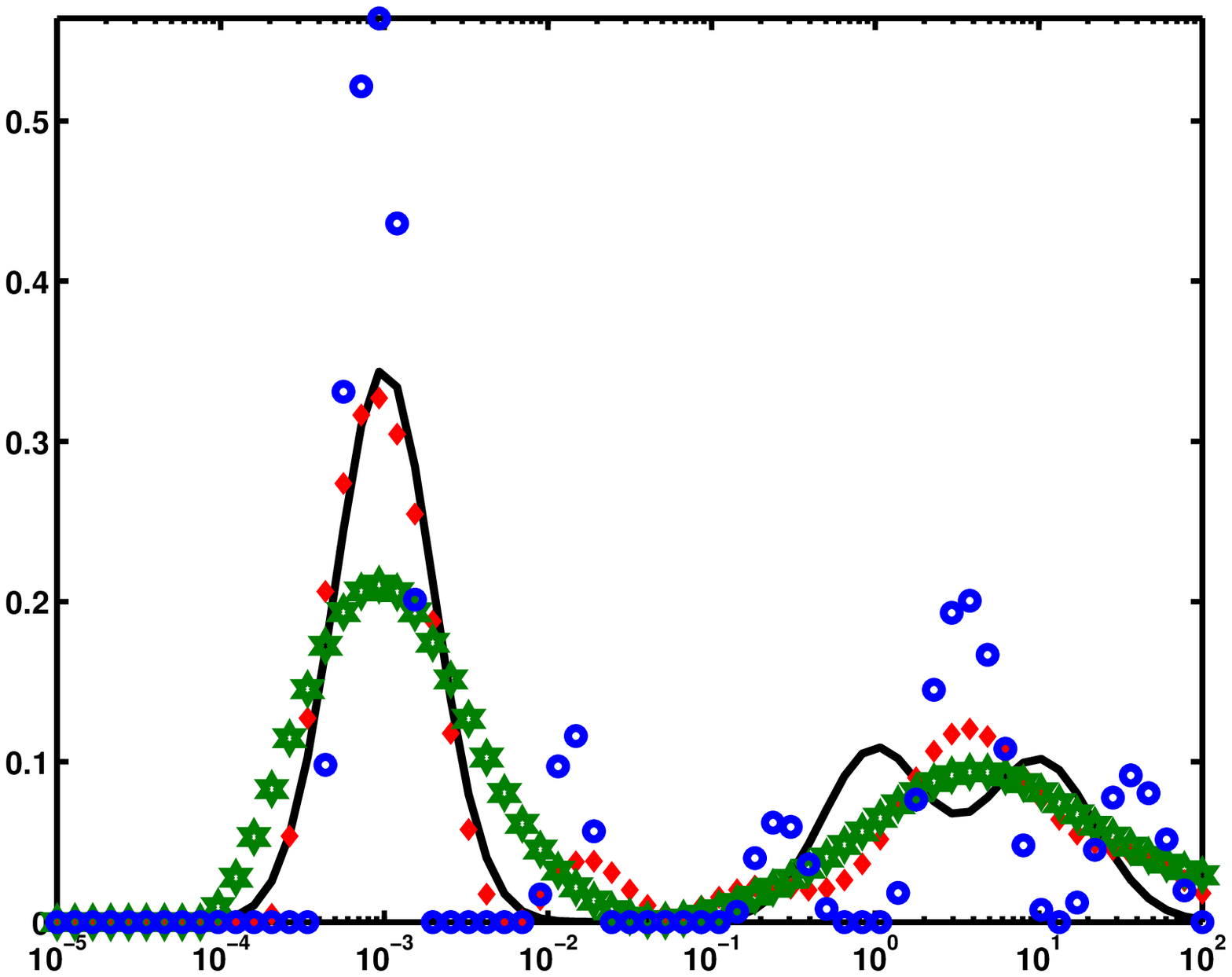}}
\subfigure[$L=L_2$]{\includegraphics[width=1.7in]{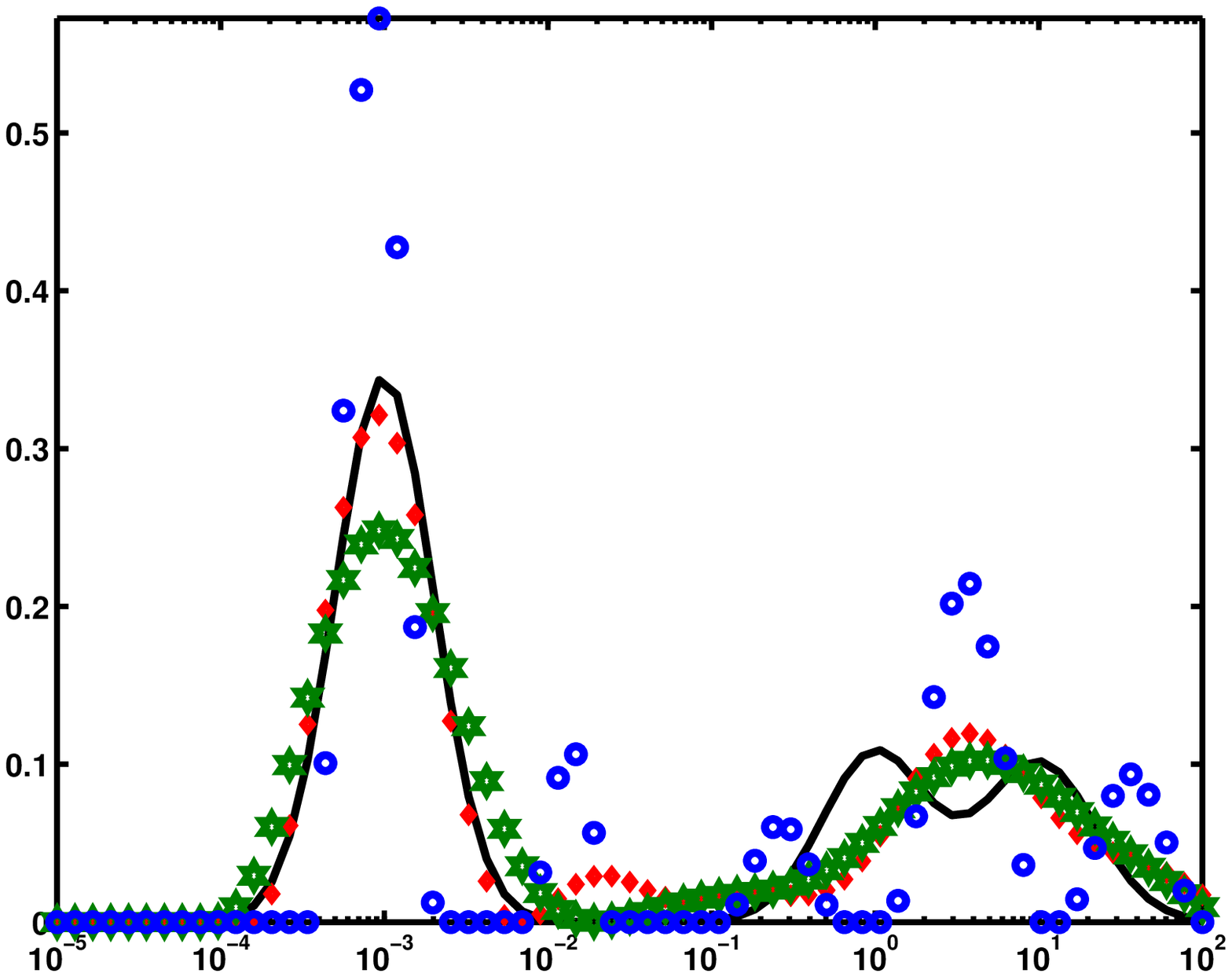}}
\caption{NNLS solutions of LN-C matrix $A_3$. Noise level $5\%$.}
\label{hnfig-lambdachoiceLN6A3HN}
\end{figure}
\clearpage
\subsection{Results using LS $A_3$ Noise level $.1\%$}
 \begin{figure}[!h]
 \centering
\subfigure[$L=I$]{\includegraphics[width=1.7in]{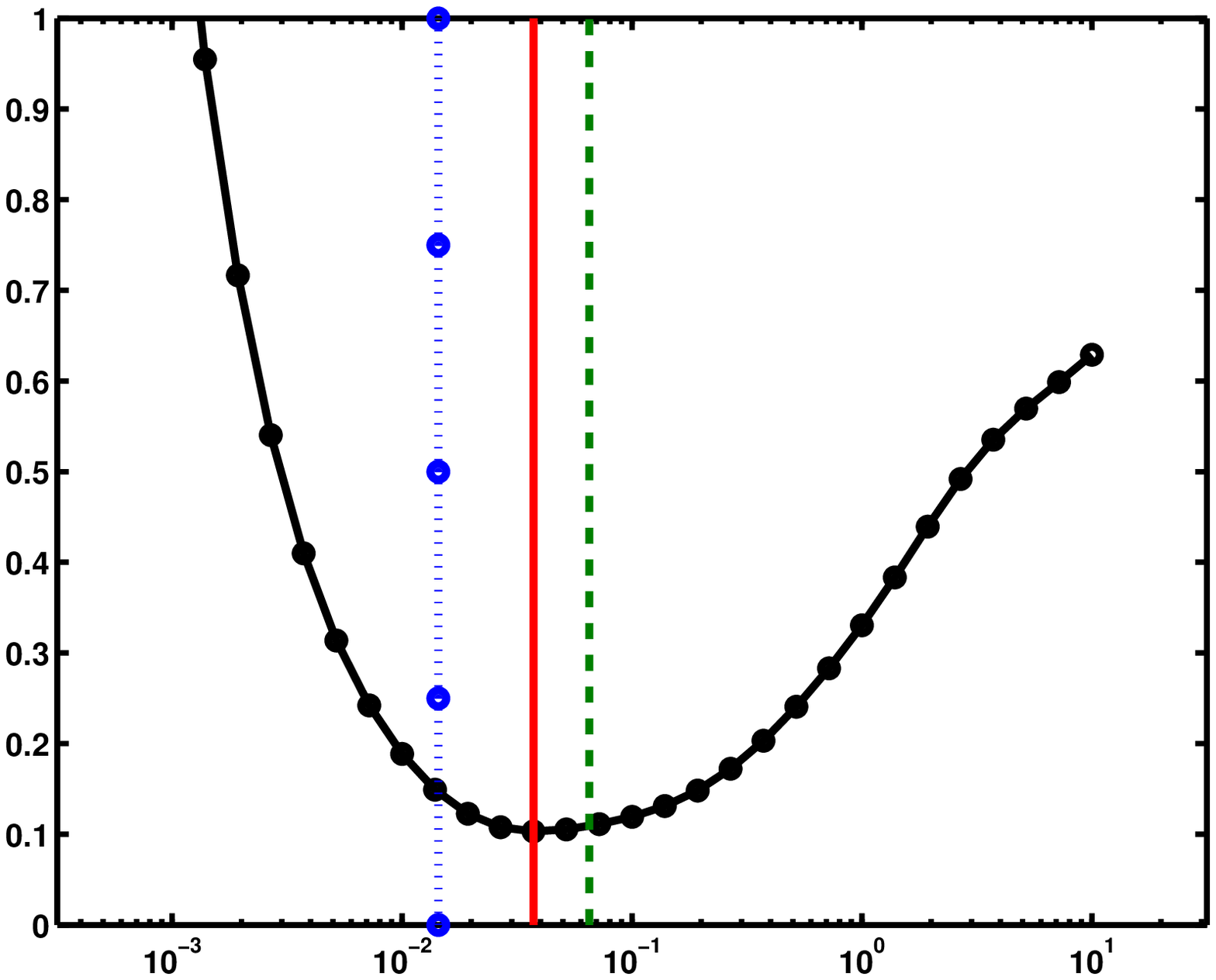}}
\subfigure[$L=L_1$]{\includegraphics[width=1.7in]{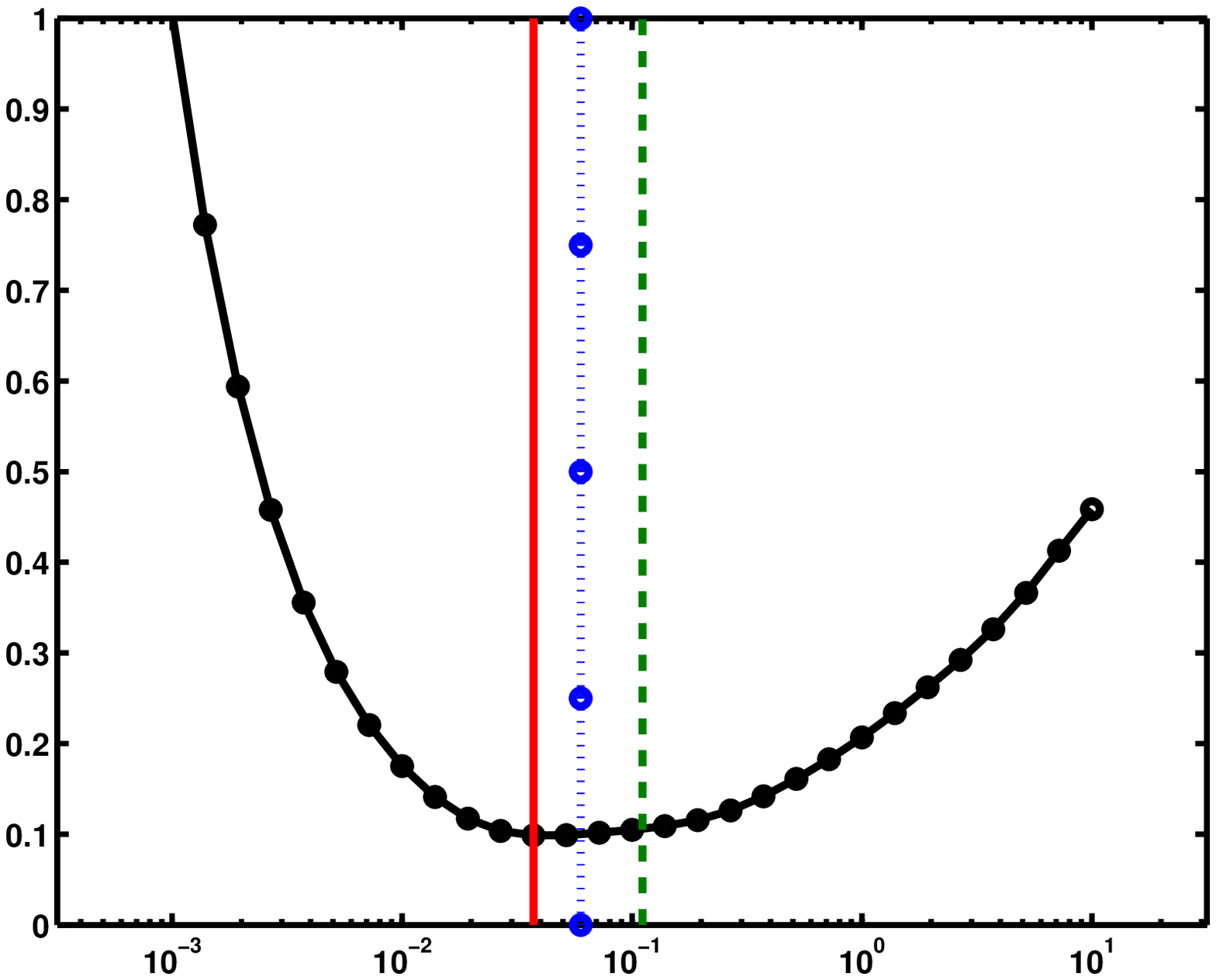}}
\subfigure[$L=L_2$]{\includegraphics[width=1.7in]{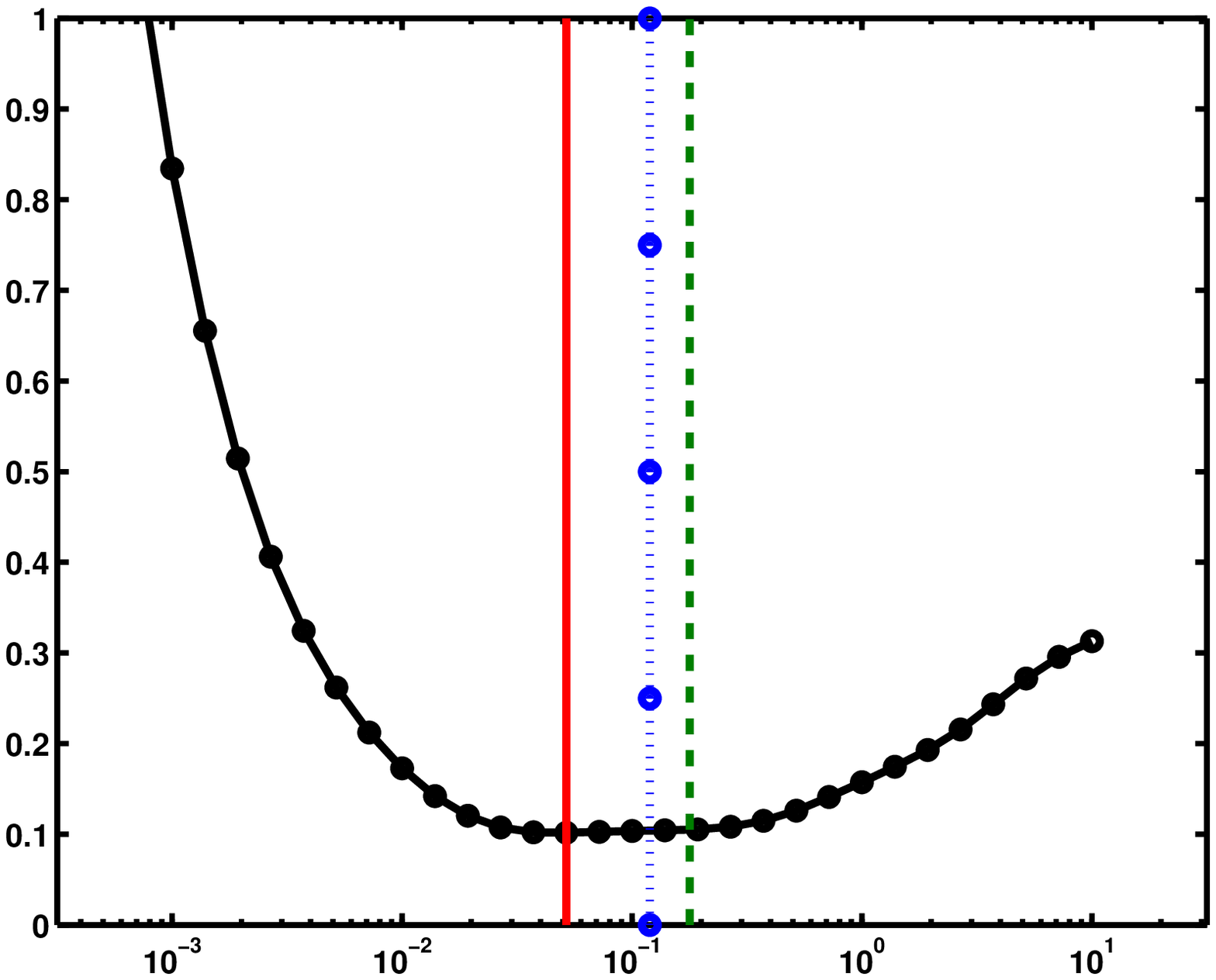}}
\subfigure[$L=I$]{\includegraphics[width=1.7in]{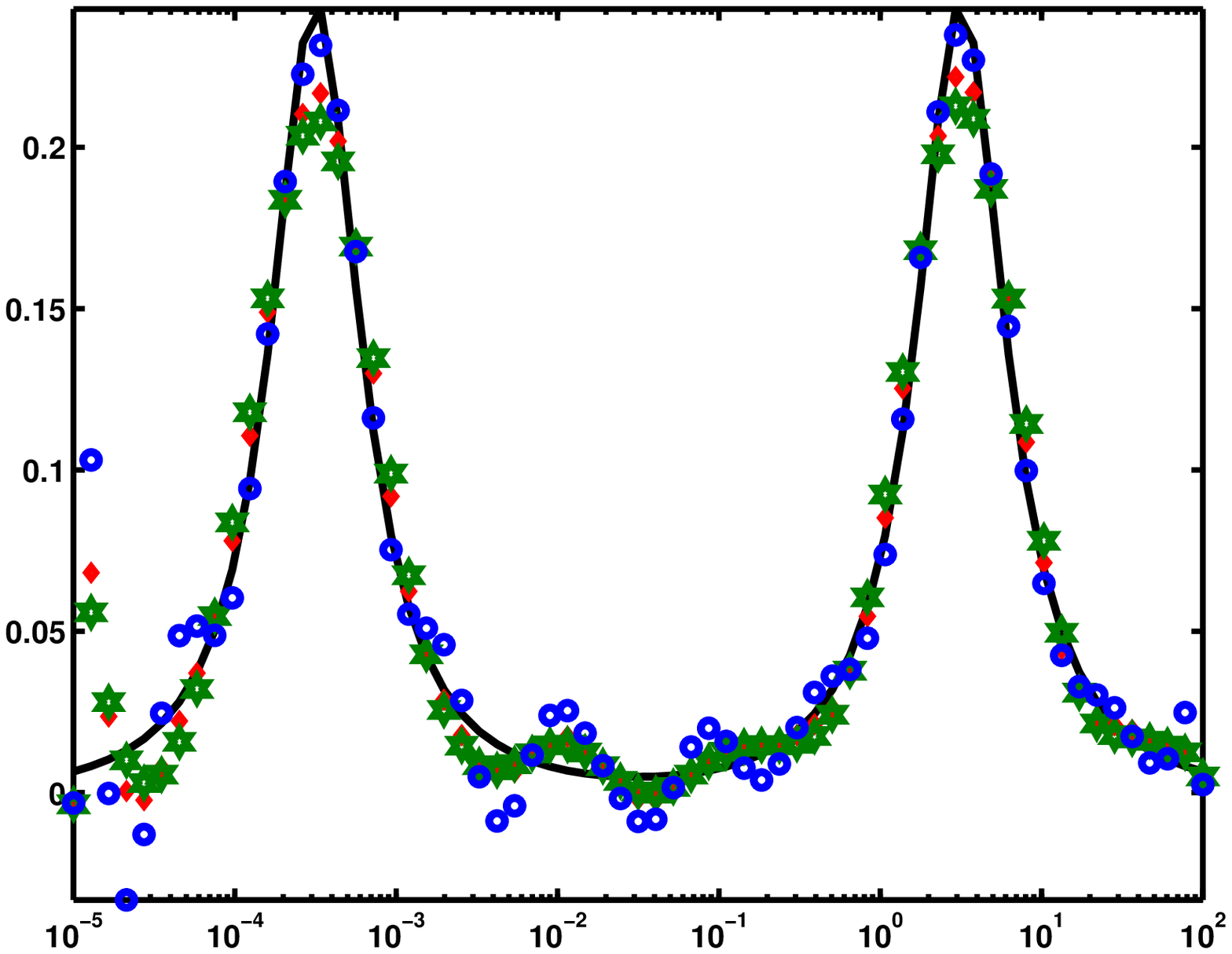}}
\subfigure[$L=L_1$]{\includegraphics[width=1.7in]{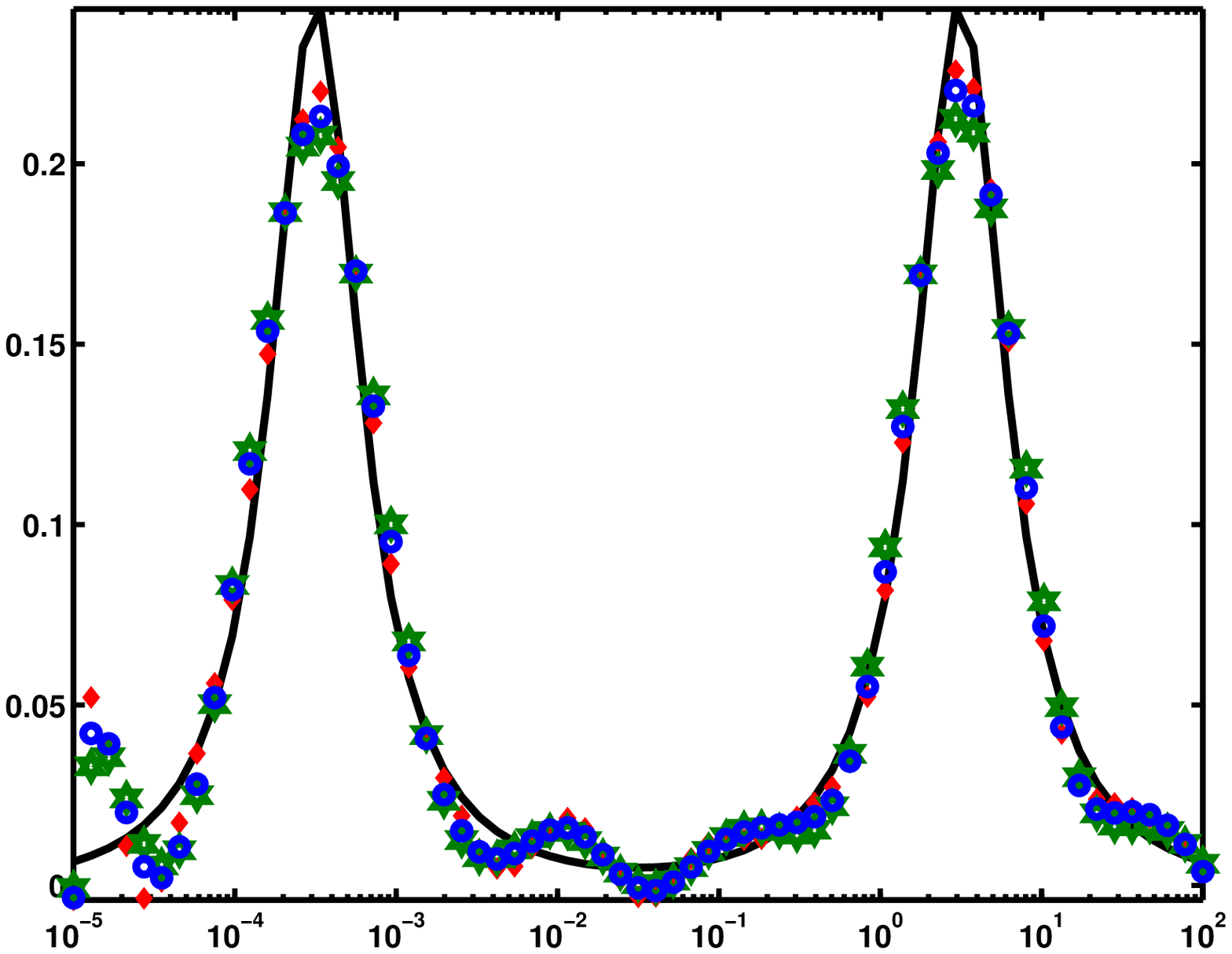}}
\subfigure[$L=L_2$]{\includegraphics[width=1.7in]{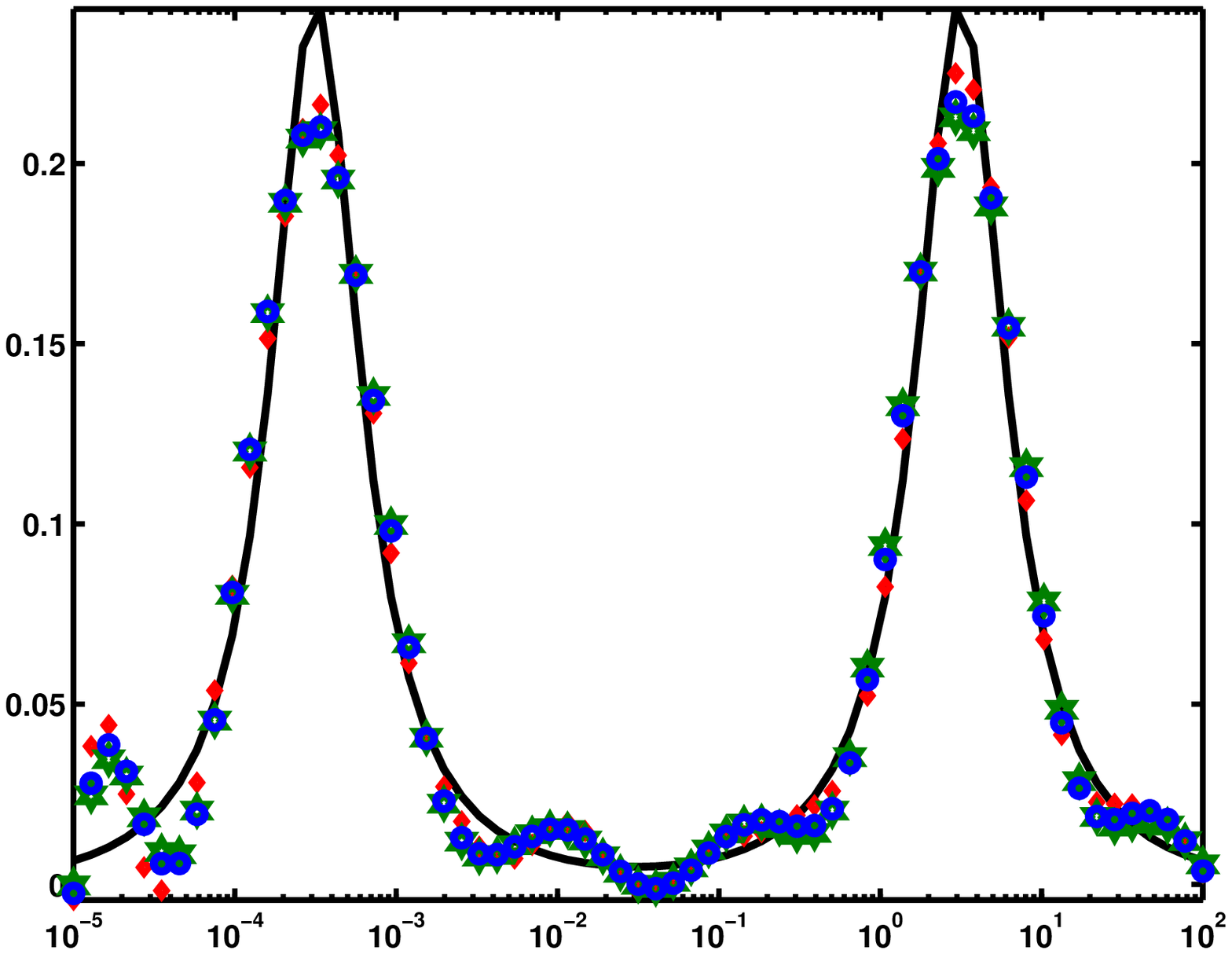}}
\caption{Mean error and example LS solutions.  $.1\%$ noise. RQ-A data set matrix $A_3$}
\label{fig-lambdachoiceRQ1A3LNLS}
\end{figure}

 \begin{figure}[!h]
  \centering
\subfigure[$L=I$]{\includegraphics[width=1.7in]{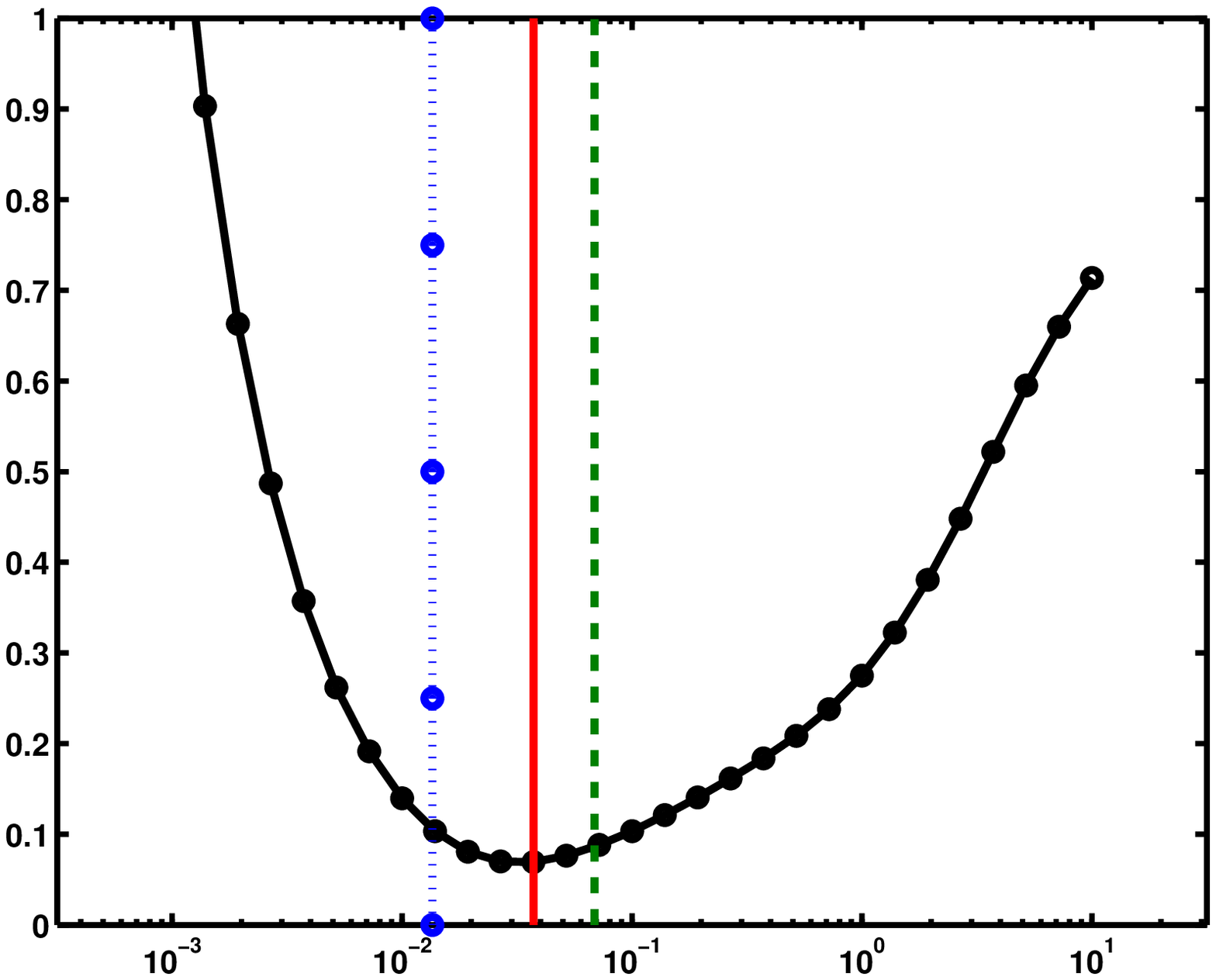}}
\subfigure[$L=L_1$]{\includegraphics[width=1.7in]{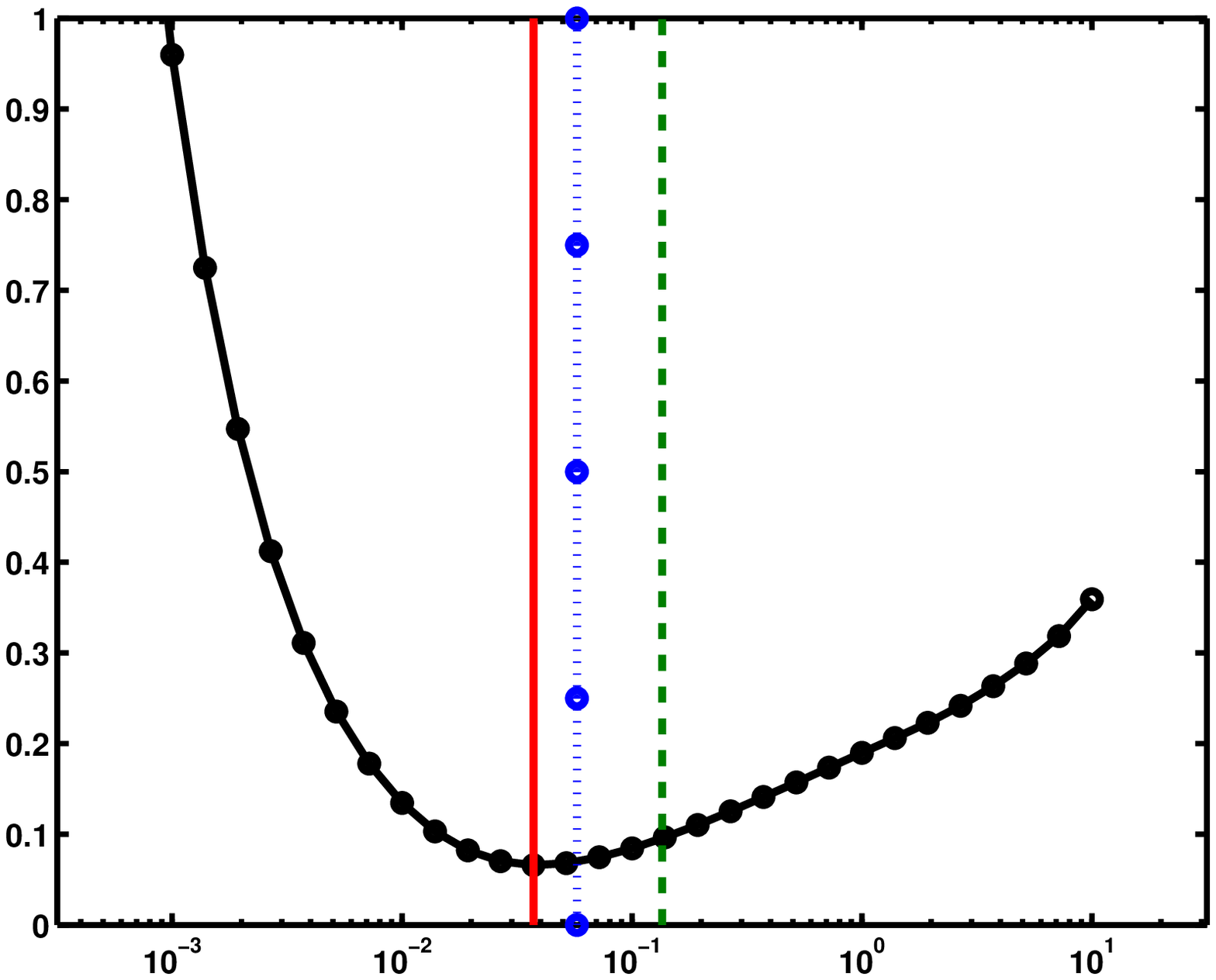}}
\subfigure[$L=L_2$]{\includegraphics[width=1.7in]{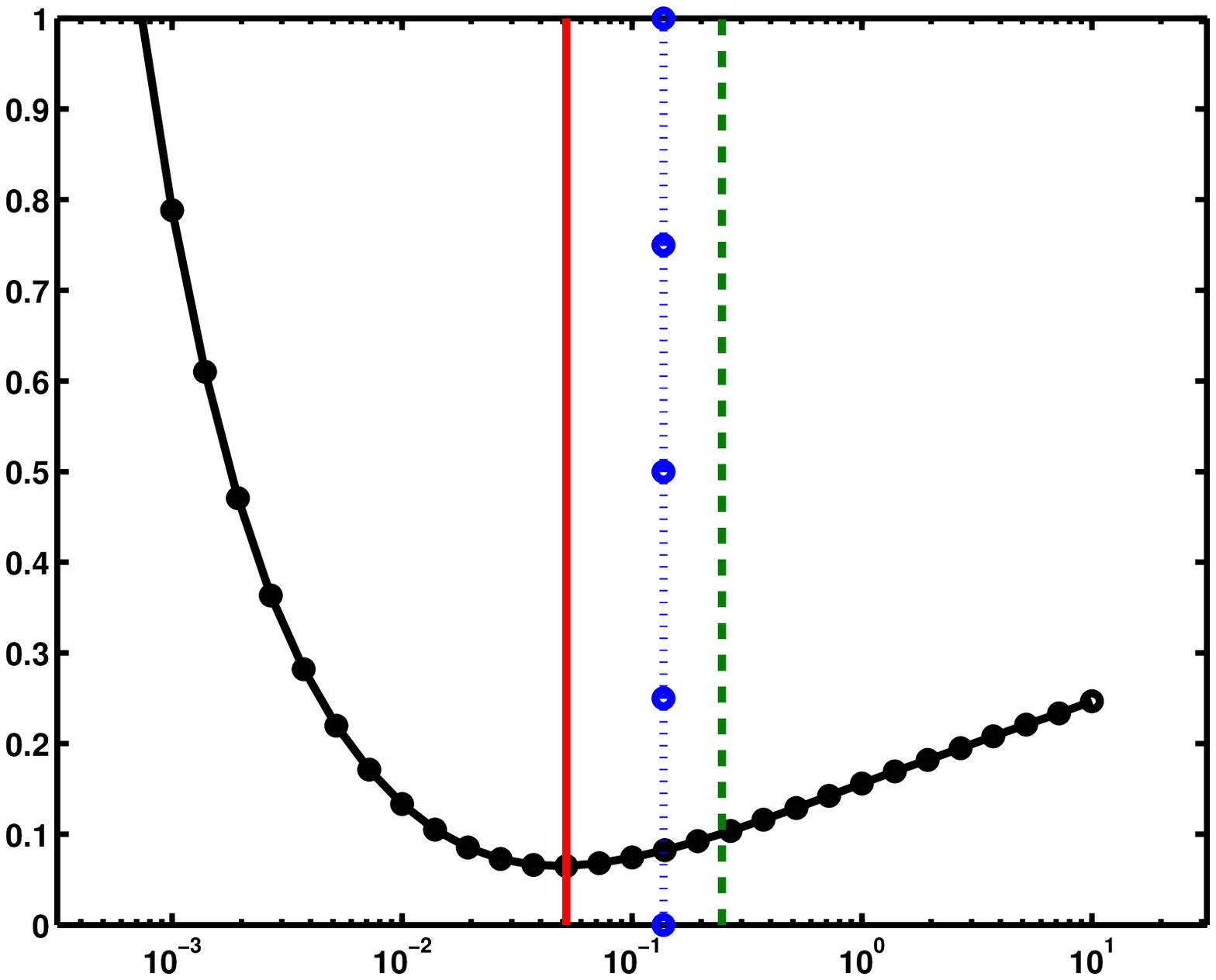}}
\subfigure[$L=I$]{\includegraphics[width=1.7in]{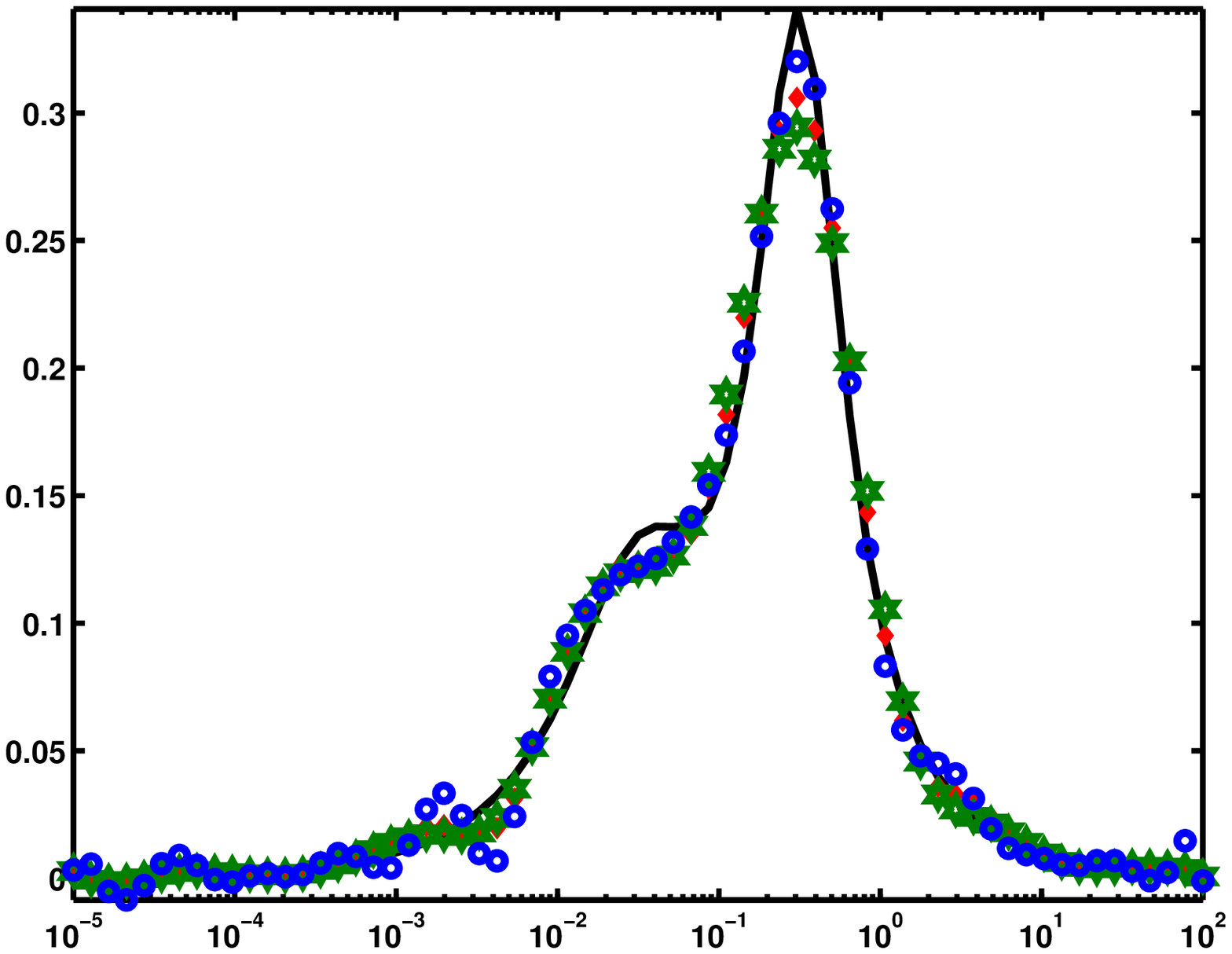}}
\subfigure[$L=L_1$]{\includegraphics[width=1.7in]{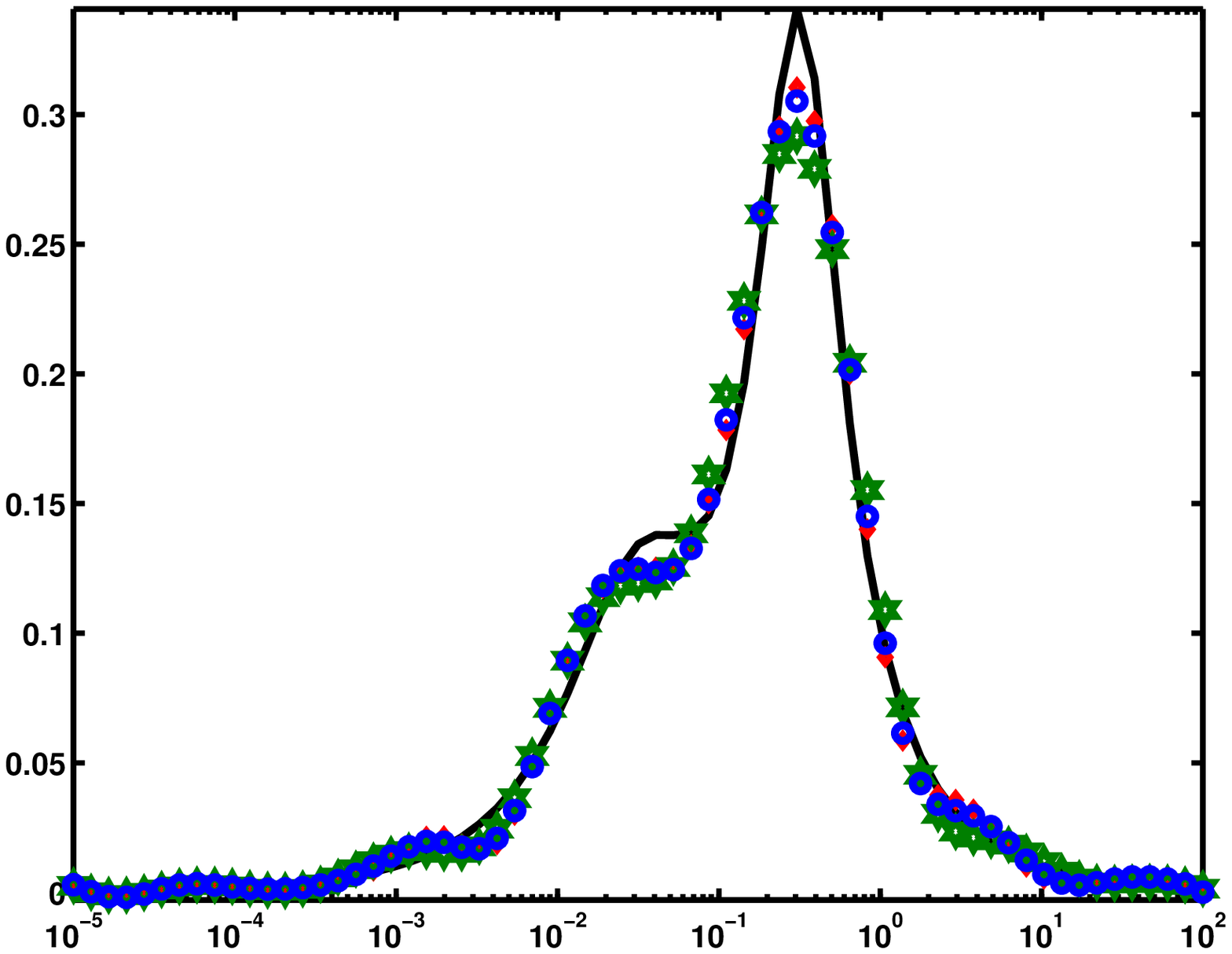}}
\subfigure[$L=L_2$]{\includegraphics[width=1.7in]{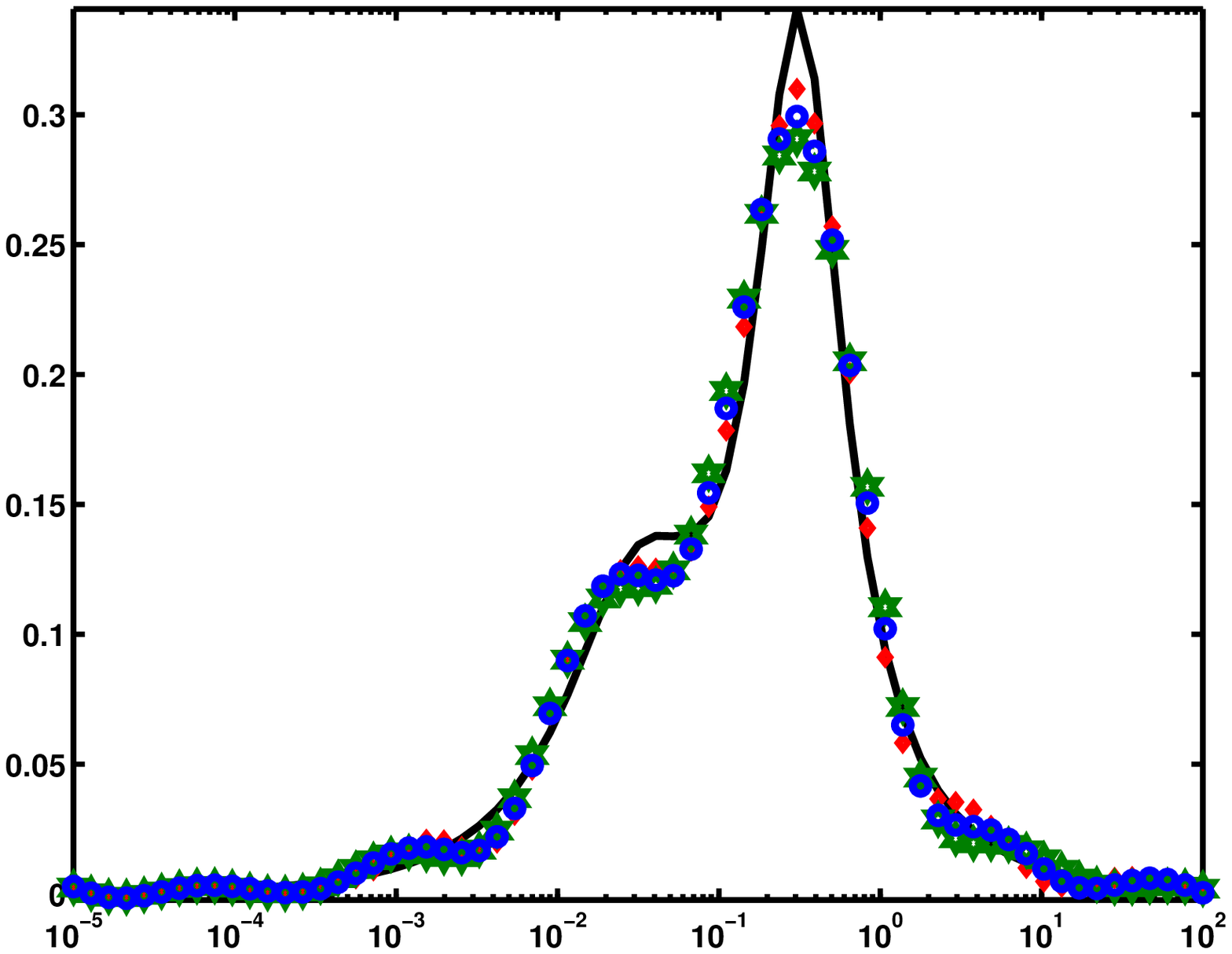}}
\caption{Mean error and example LS solutions. $.1\%$ noise. RQ-B data set matrix $A_3$}
\label{fig-lambdachoiceRQ5A3LNLS}
\end{figure}

 \begin{figure}[!h]
  \centering
\subfigure[$L=I$]{\includegraphics[width=1.7in]{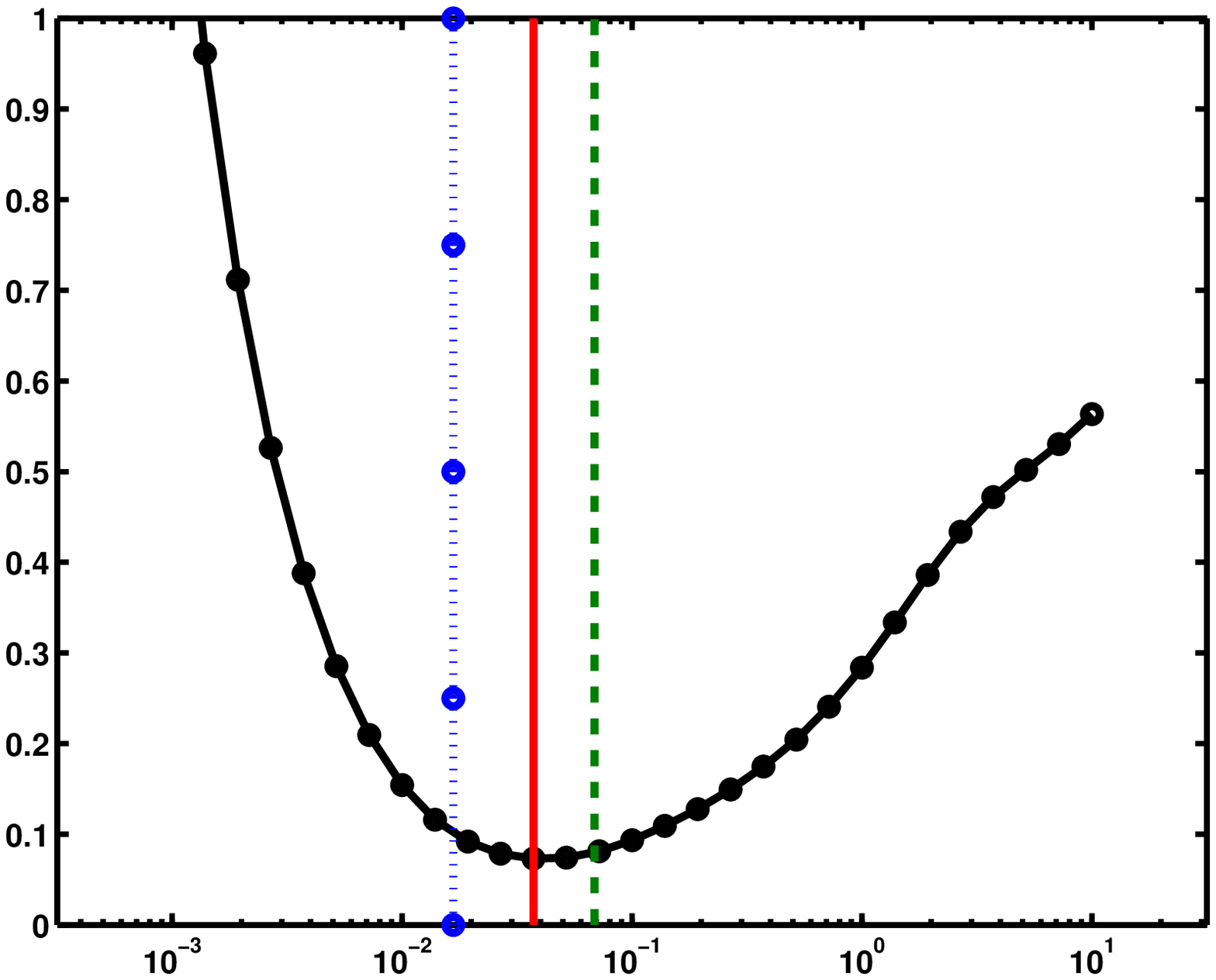}}
\subfigure[$L=L_1$]{\includegraphics[width=1.7in]{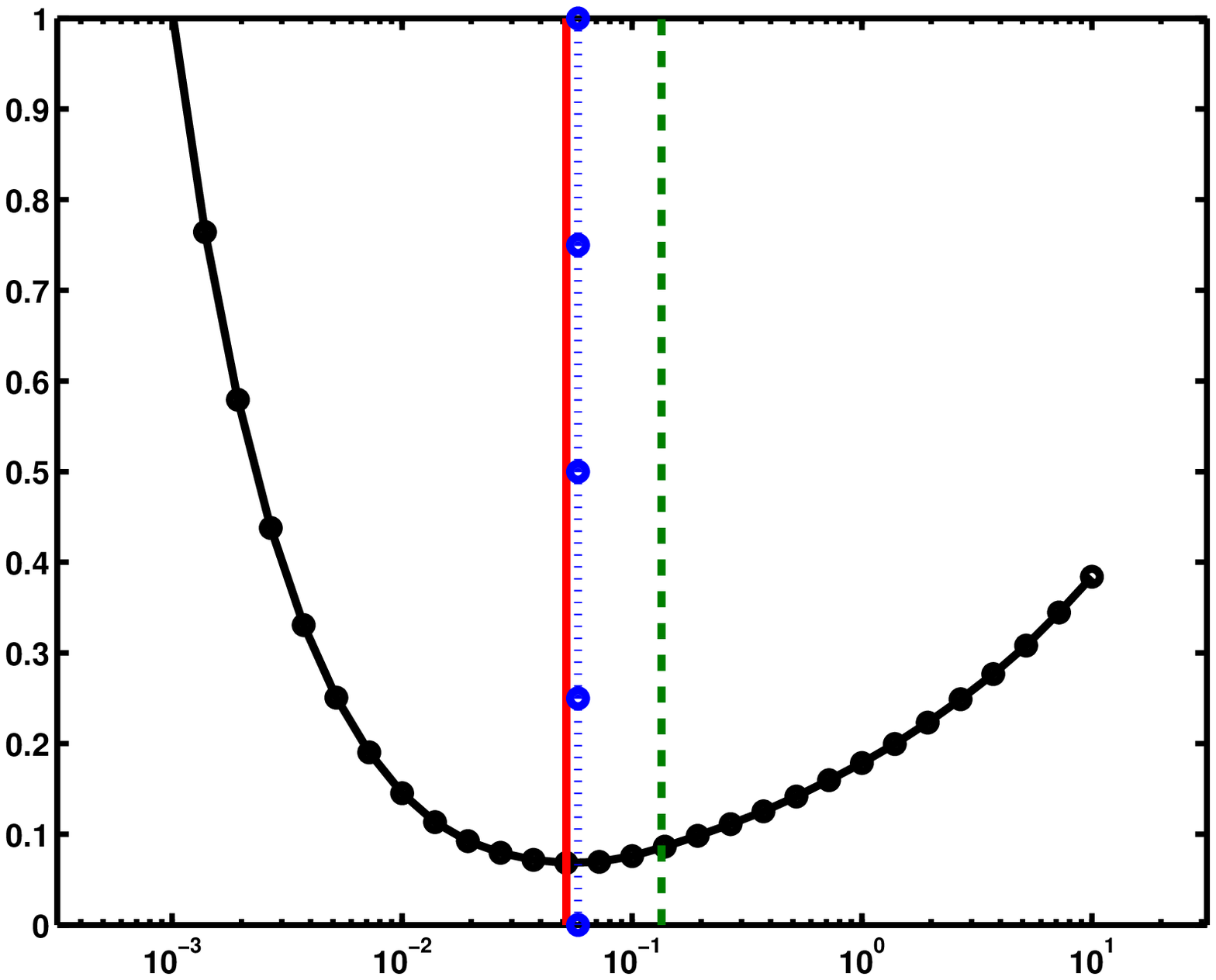}}
\subfigure[$L=L_2$]{\includegraphics[width=1.7in]{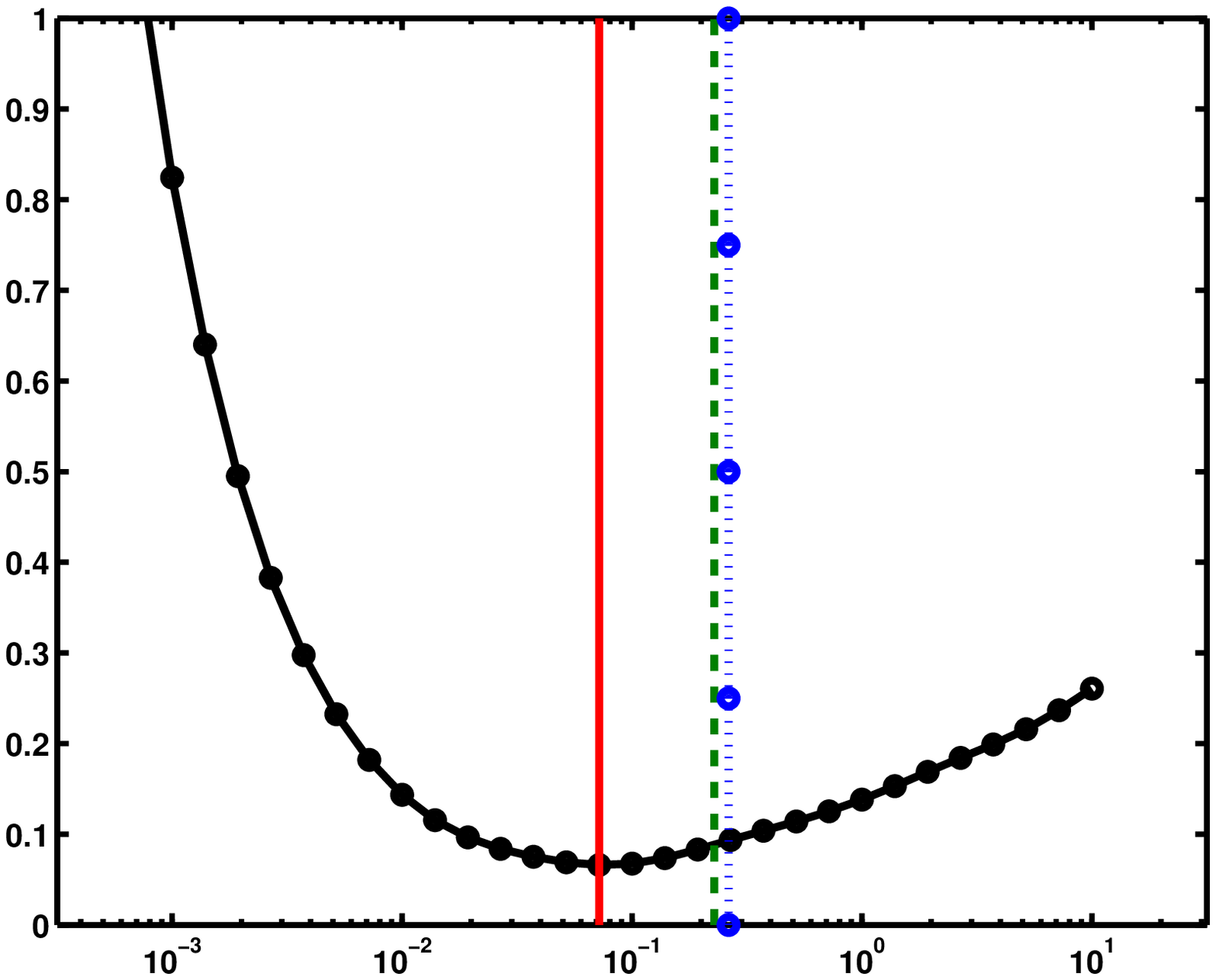}}
\subfigure[$L=I$]{\includegraphics[width=1.7in]{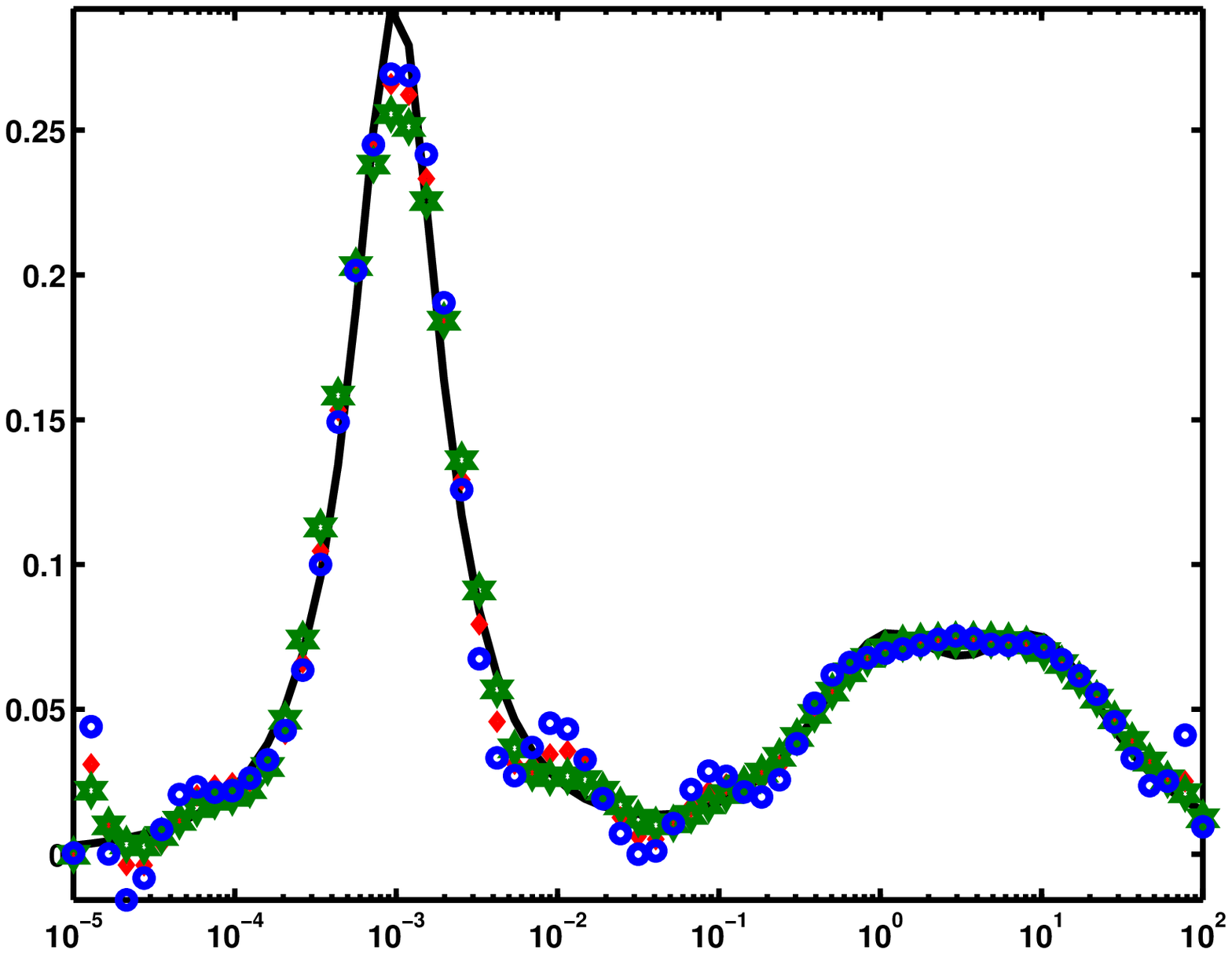}}
\subfigure[$L=L_1$]{\includegraphics[width=1.7in]{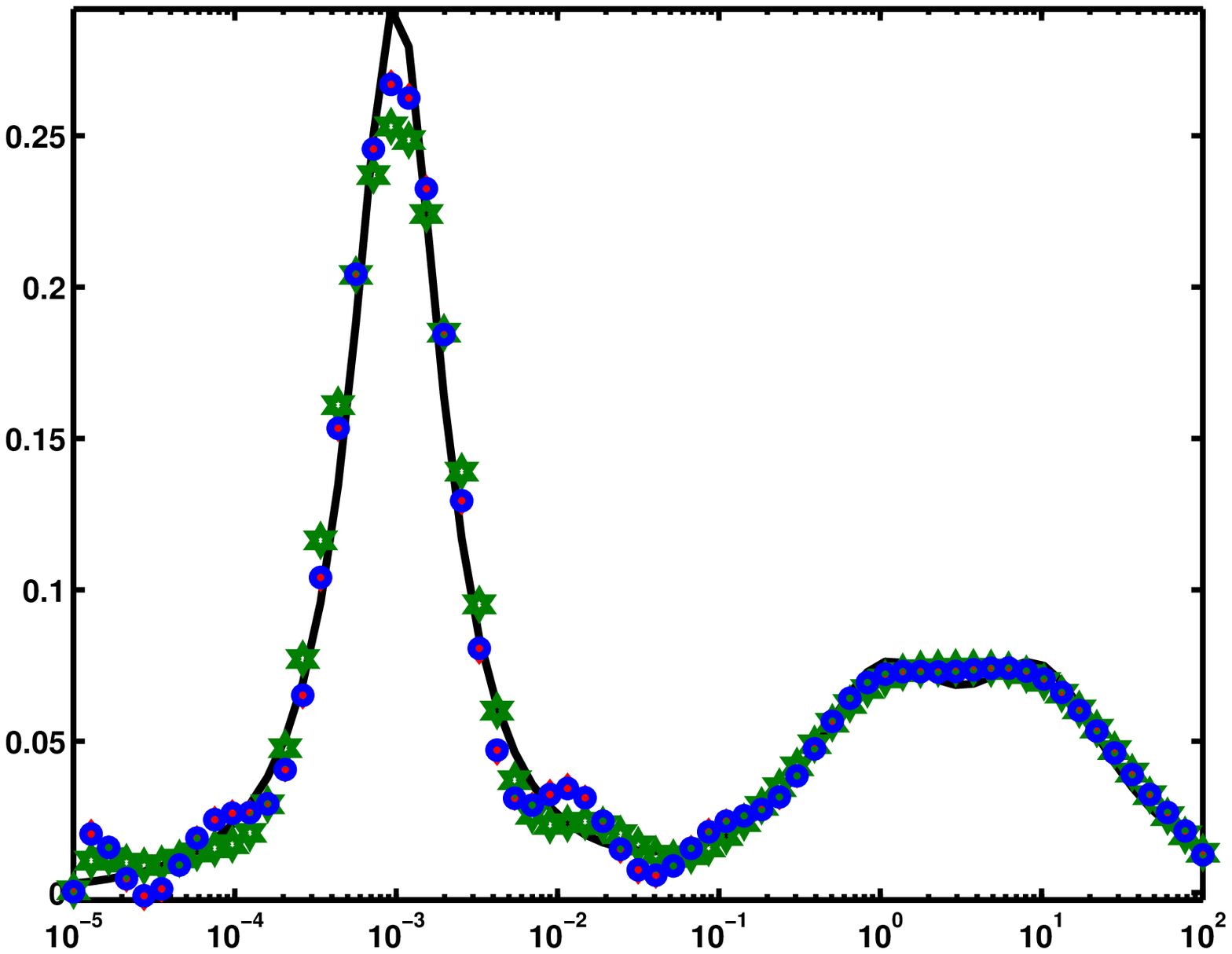}}
\subfigure[$L=L_2$]{\includegraphics[width=1.7in]{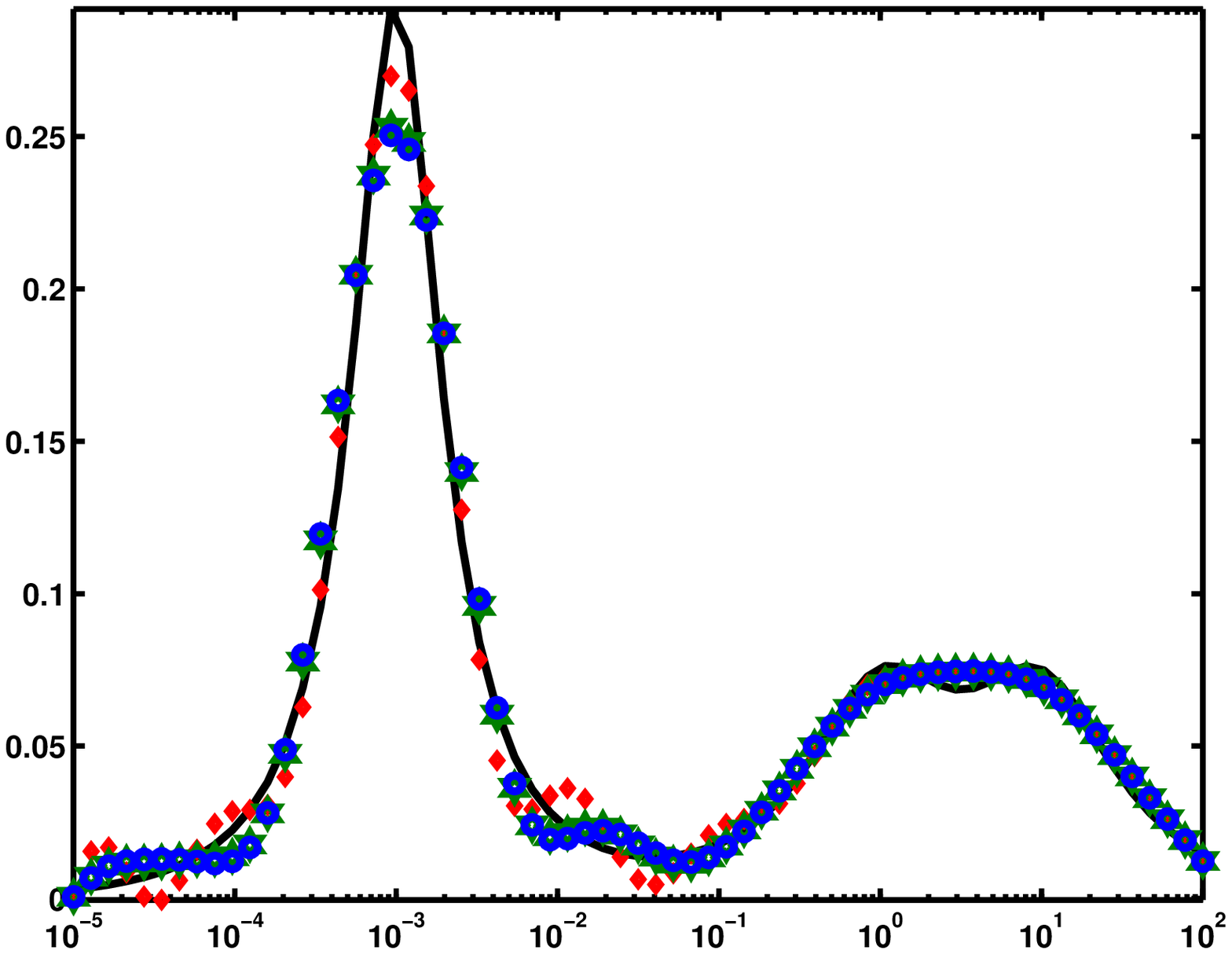}}
\caption{Mean error and example LS solutions.  $.1\%$ noise. RQ-C data set matrix $A_3$}
\label{fig-lambdachoiceRQ6A3LNLS}
\end{figure}

 \begin{figure}[!h]
  \centering
\subfigure[$L=I$]{\includegraphics[width=1.7in]{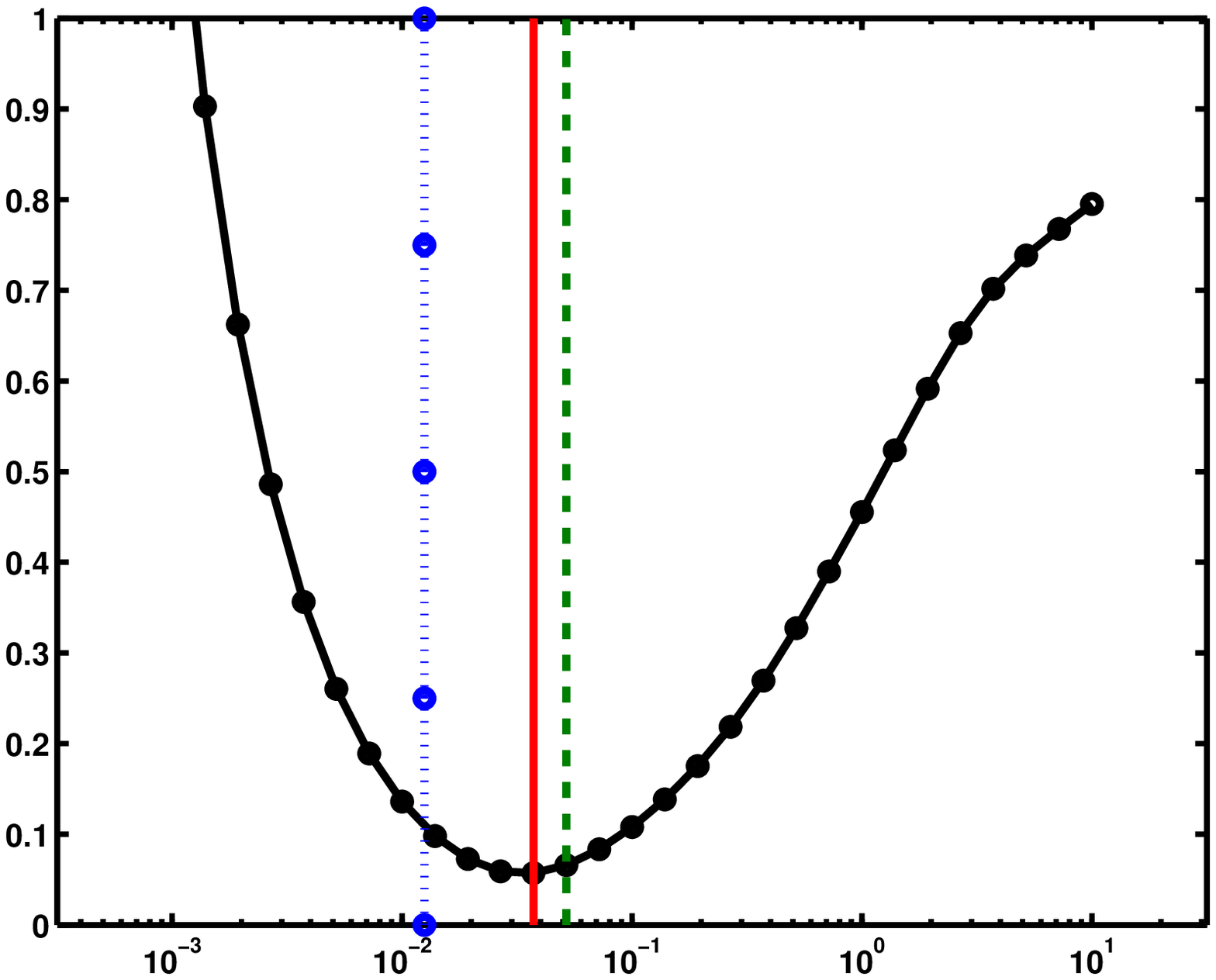}}
\subfigure[$L=L_1$]{\includegraphics[width=1.7in]{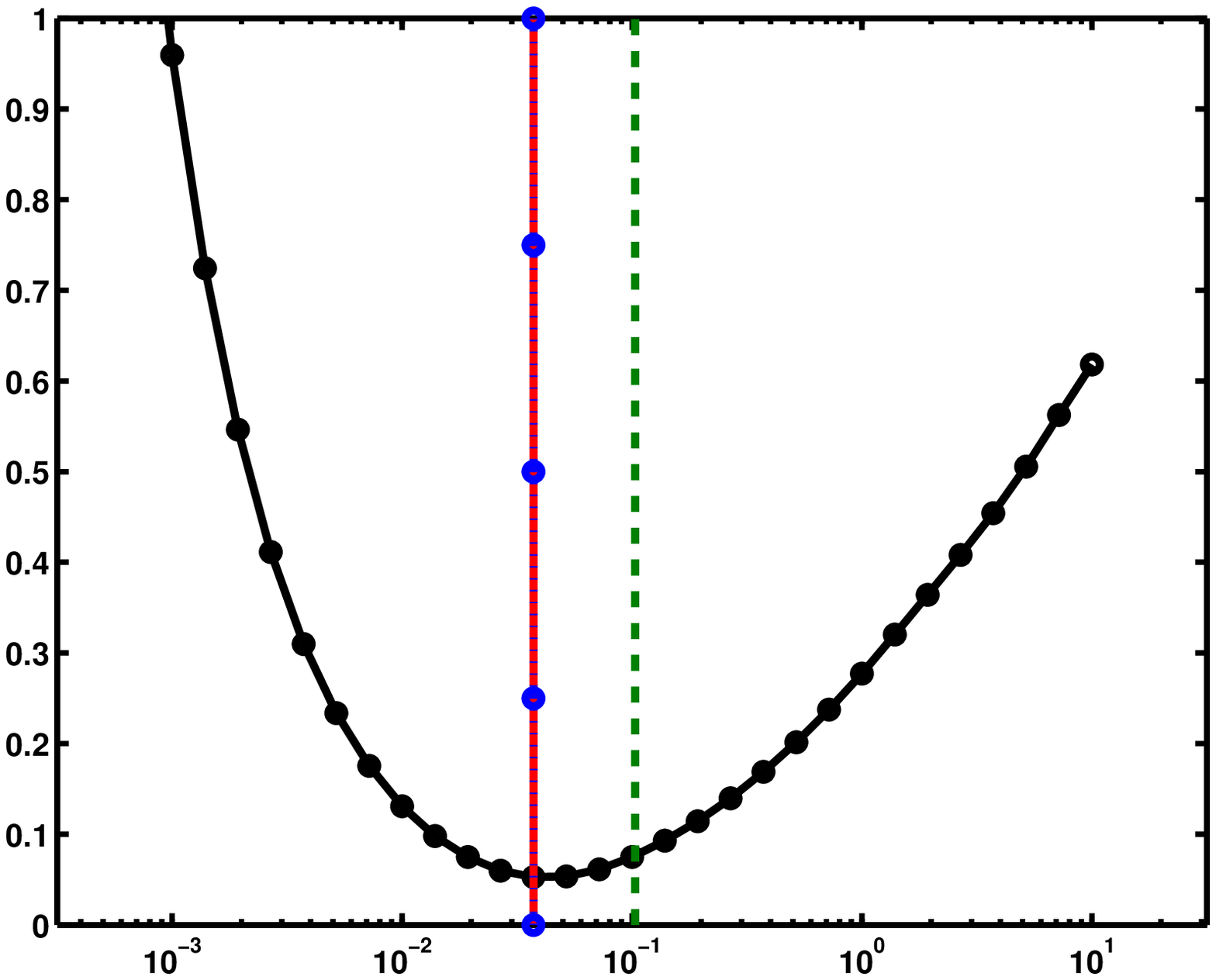}}
\subfigure[$L=L_2$]{\includegraphics[width=1.7in]{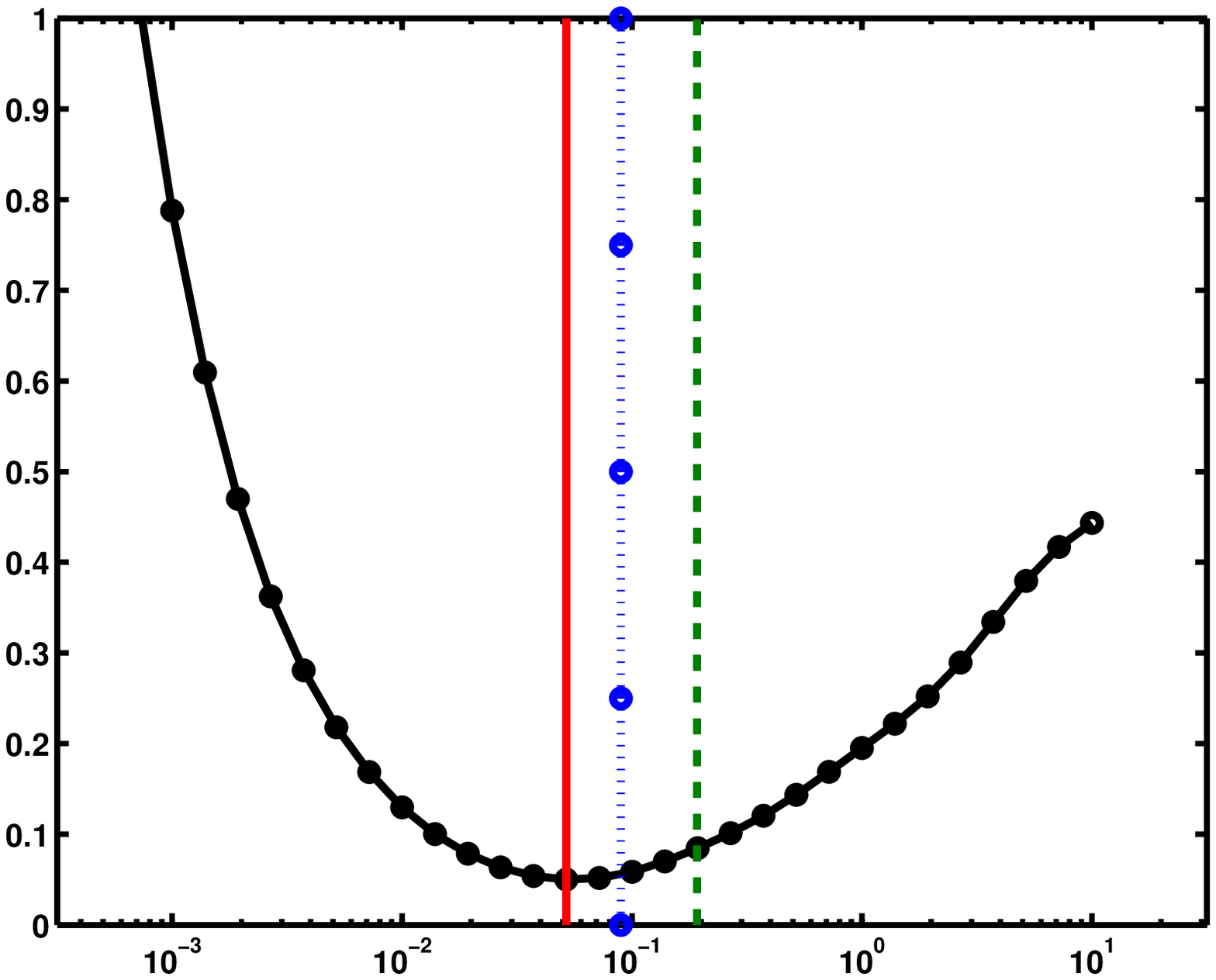}}
\subfigure[$L=I$]{\includegraphics[width=1.7in]{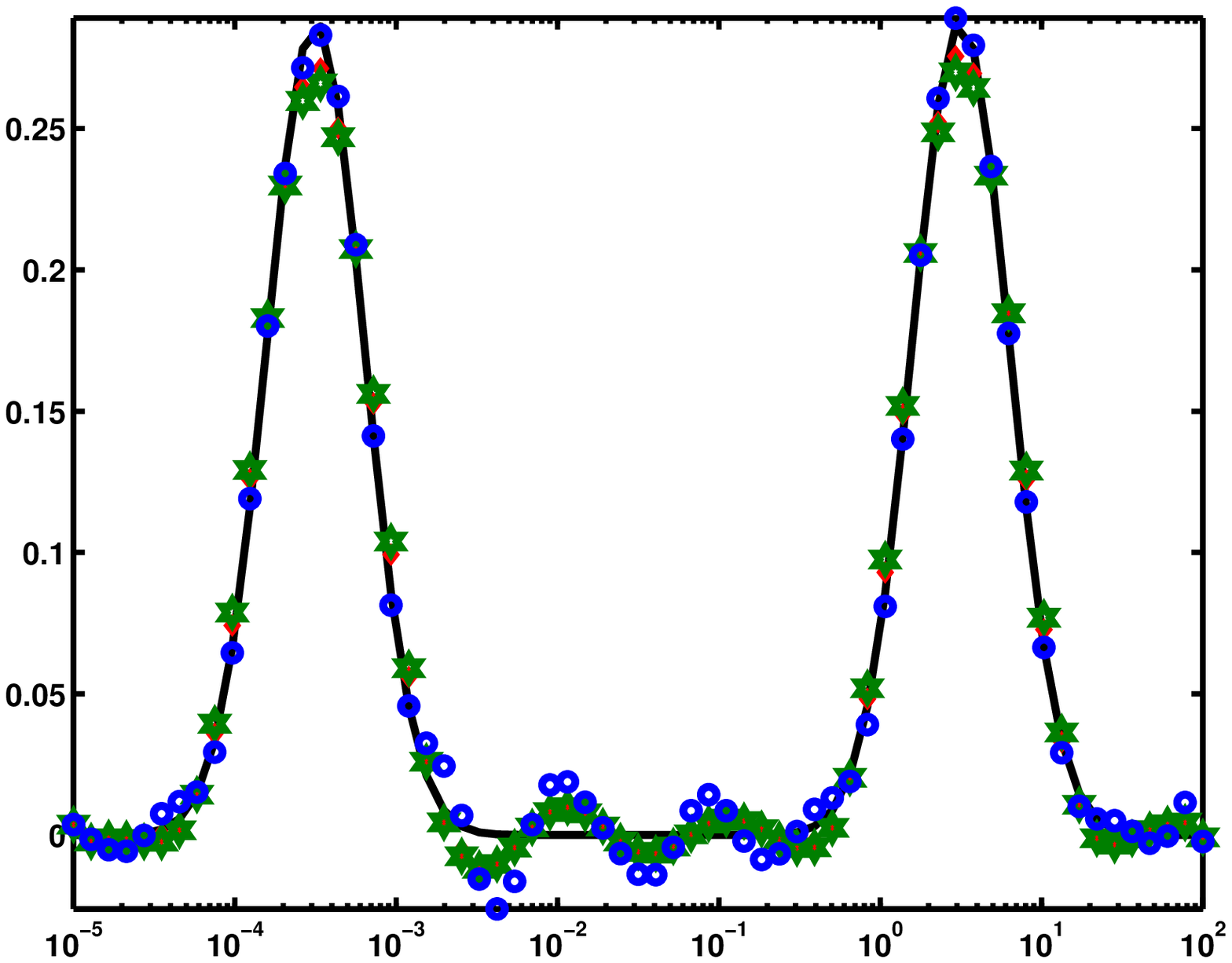}}
\subfigure[$L=L_1$]{\includegraphics[width=1.7in]{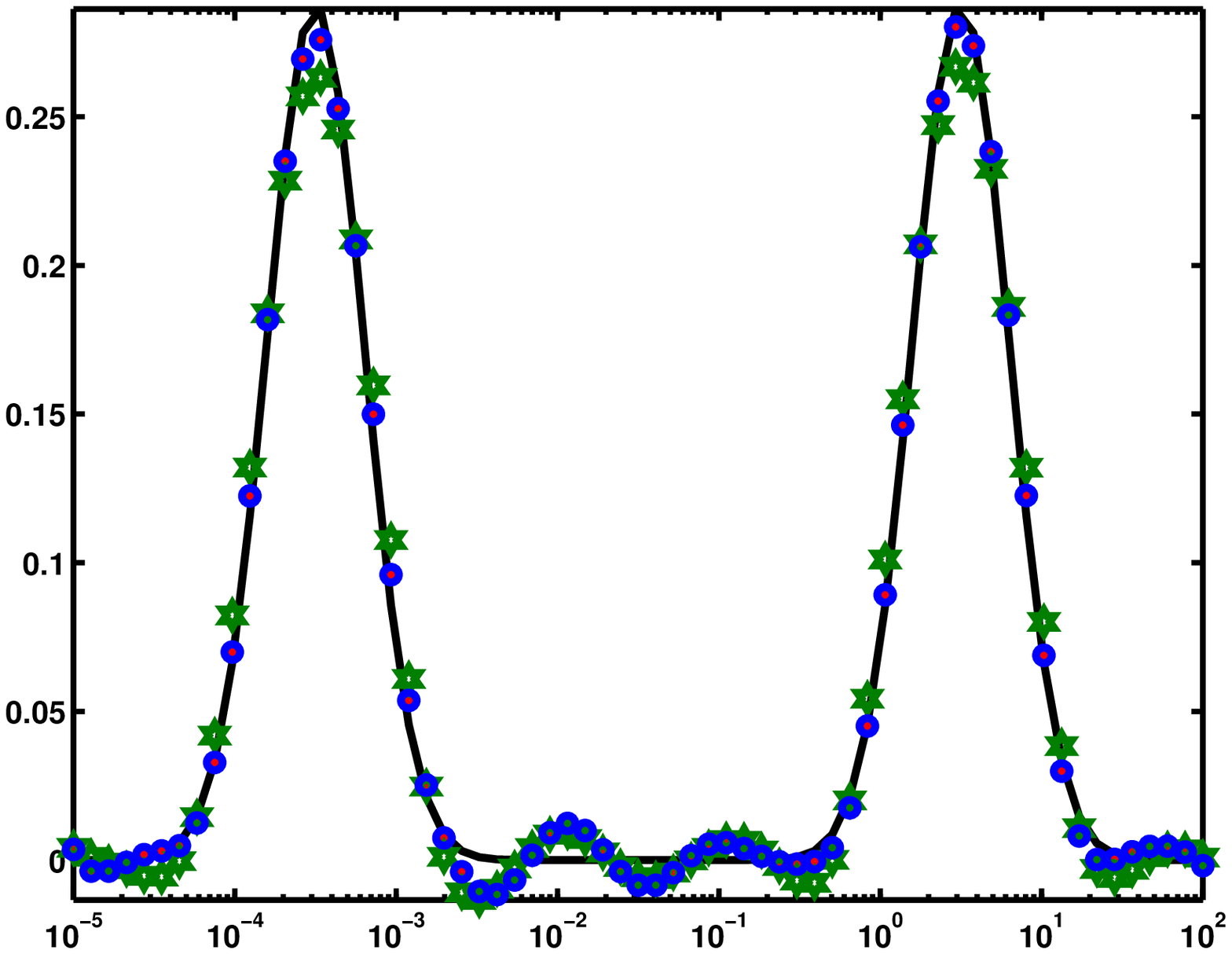}}
\subfigure[$L=L_2$]{\includegraphics[width=1.7in]{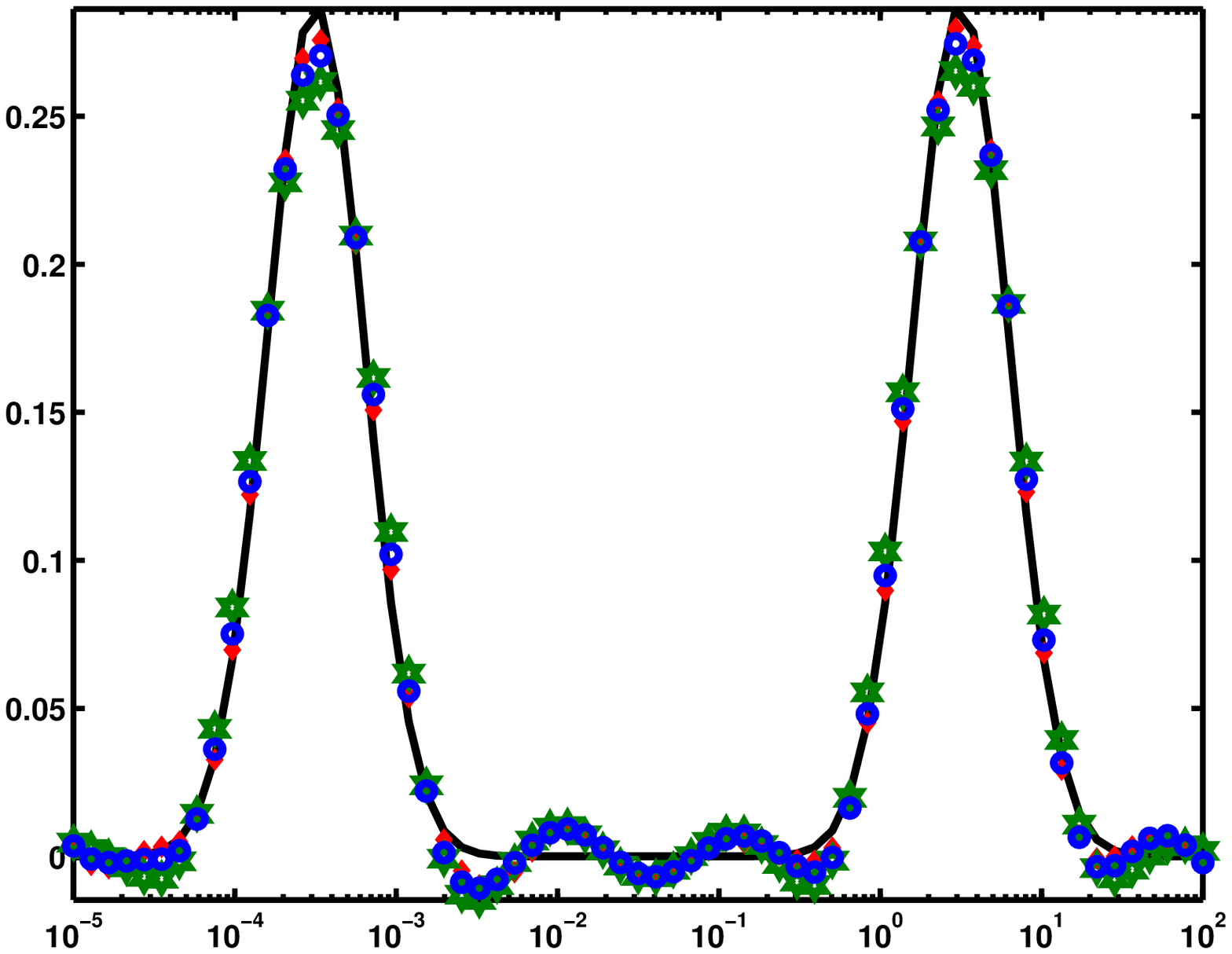}}
\caption{Mean error and example LS solutions.  $.1\%$ noise. LN-A data set matrix $A_3$ }
\label{fig-lambdachoiceLN2A3LNLS}
\end{figure}

 \begin{figure}[!h]
  \centering
\subfigure[$L=I$]{\includegraphics[width=1.7in]{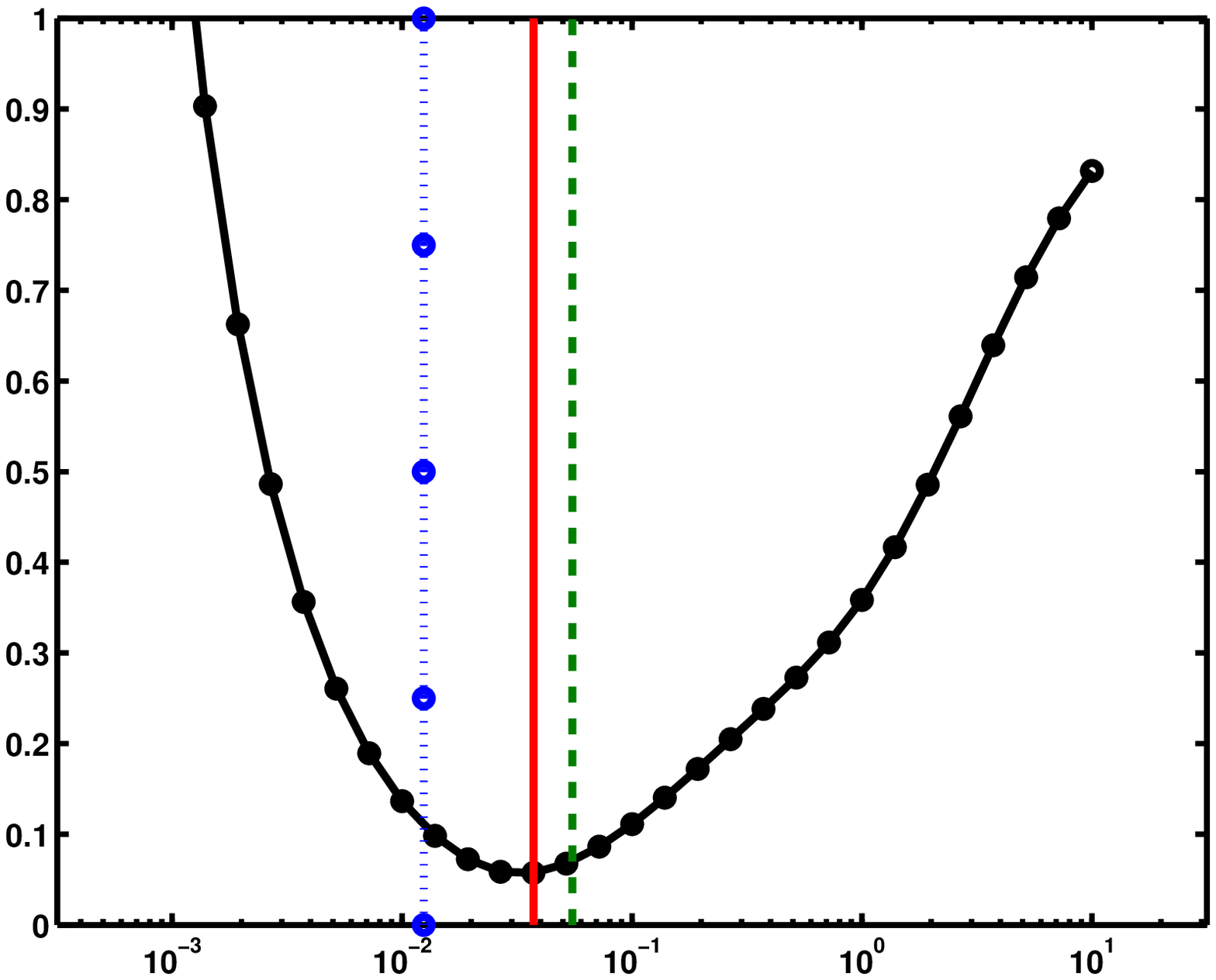}}
\subfigure[$L=L_1$]{\includegraphics[width=1.7in]{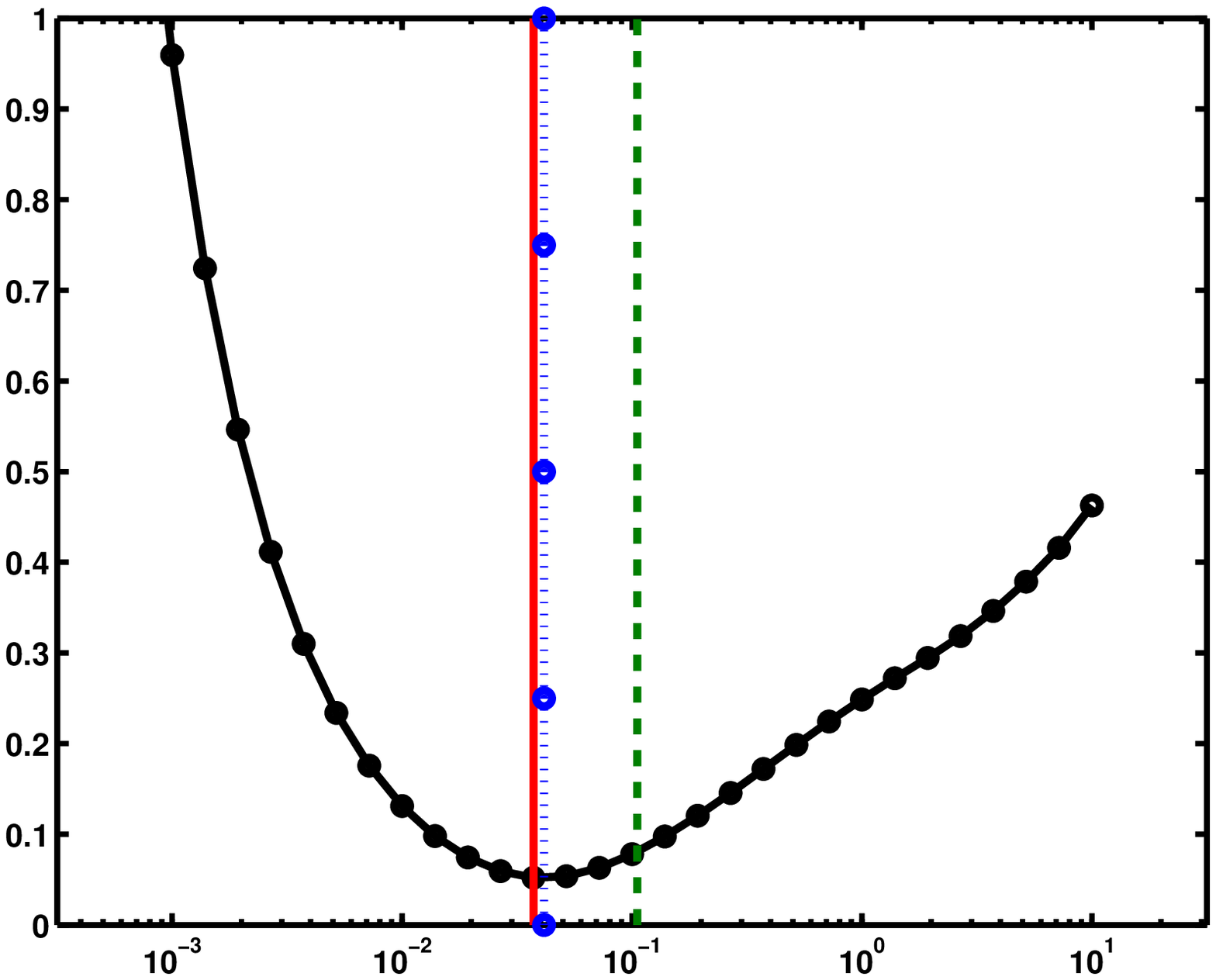}}
\subfigure[$L=L_2$]{\includegraphics[width=1.7in]{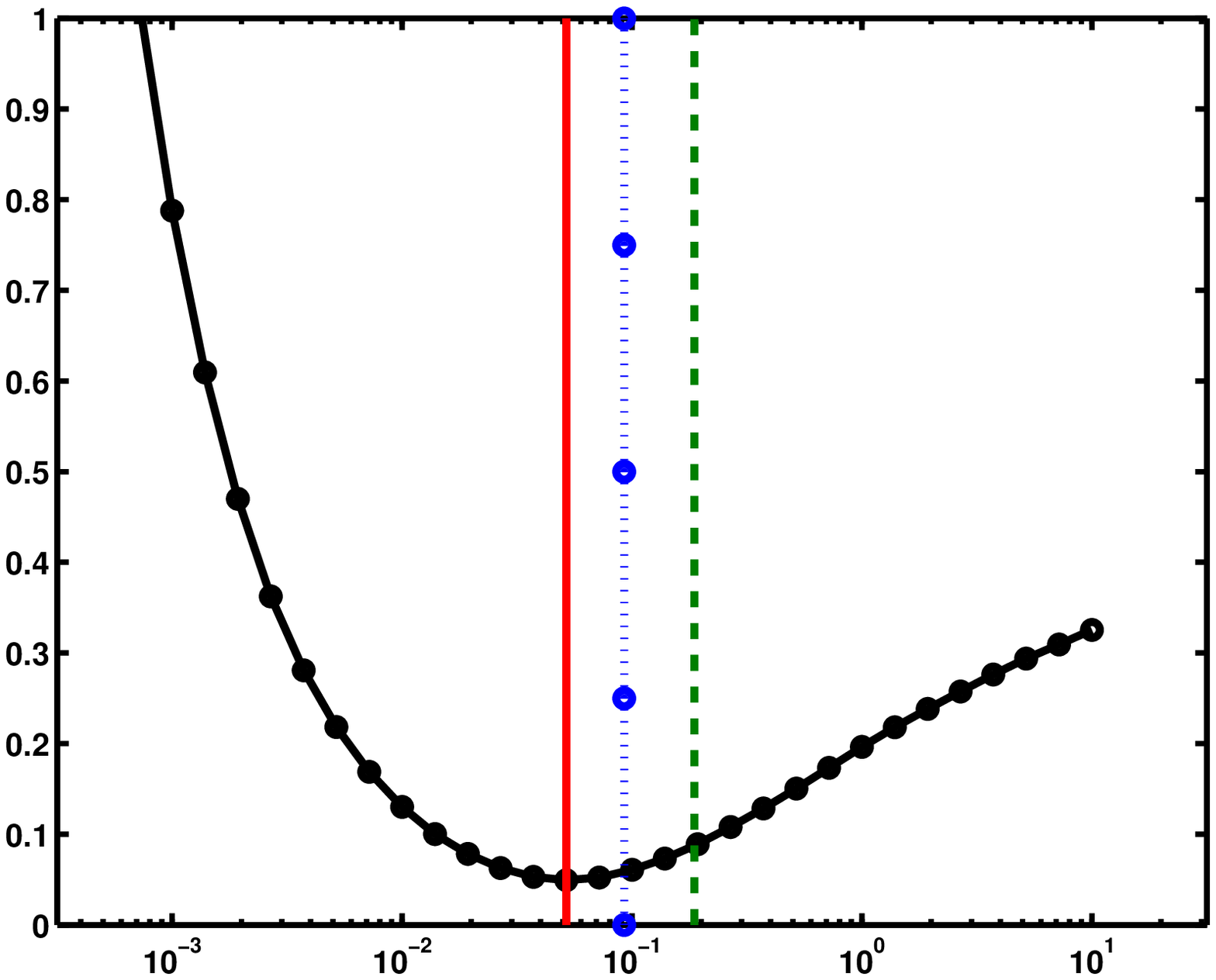}}
\subfigure[$L=I$]{\includegraphics[width=1.7in]{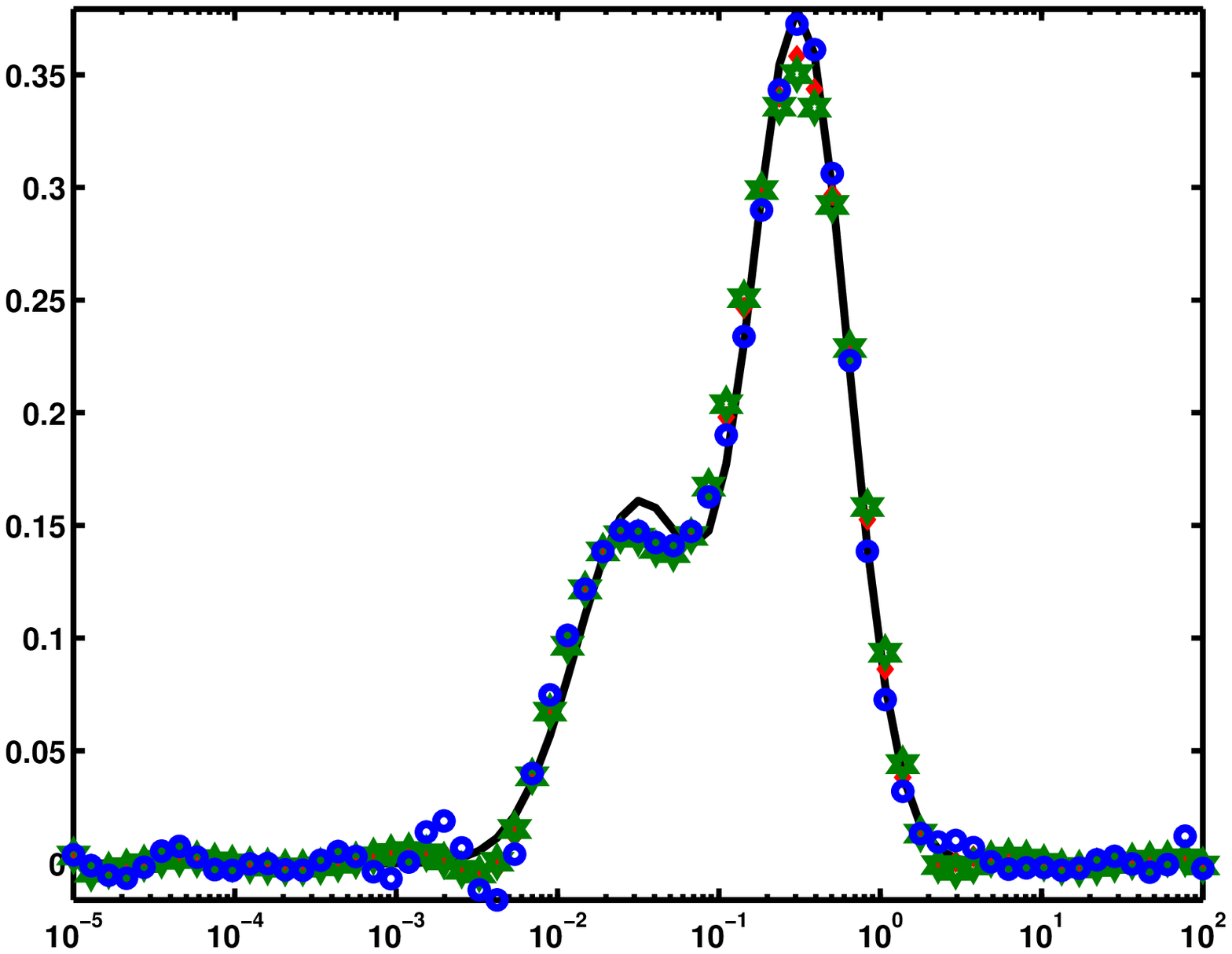}}
\subfigure[$L=L_1$]{\includegraphics[width=1.7in]{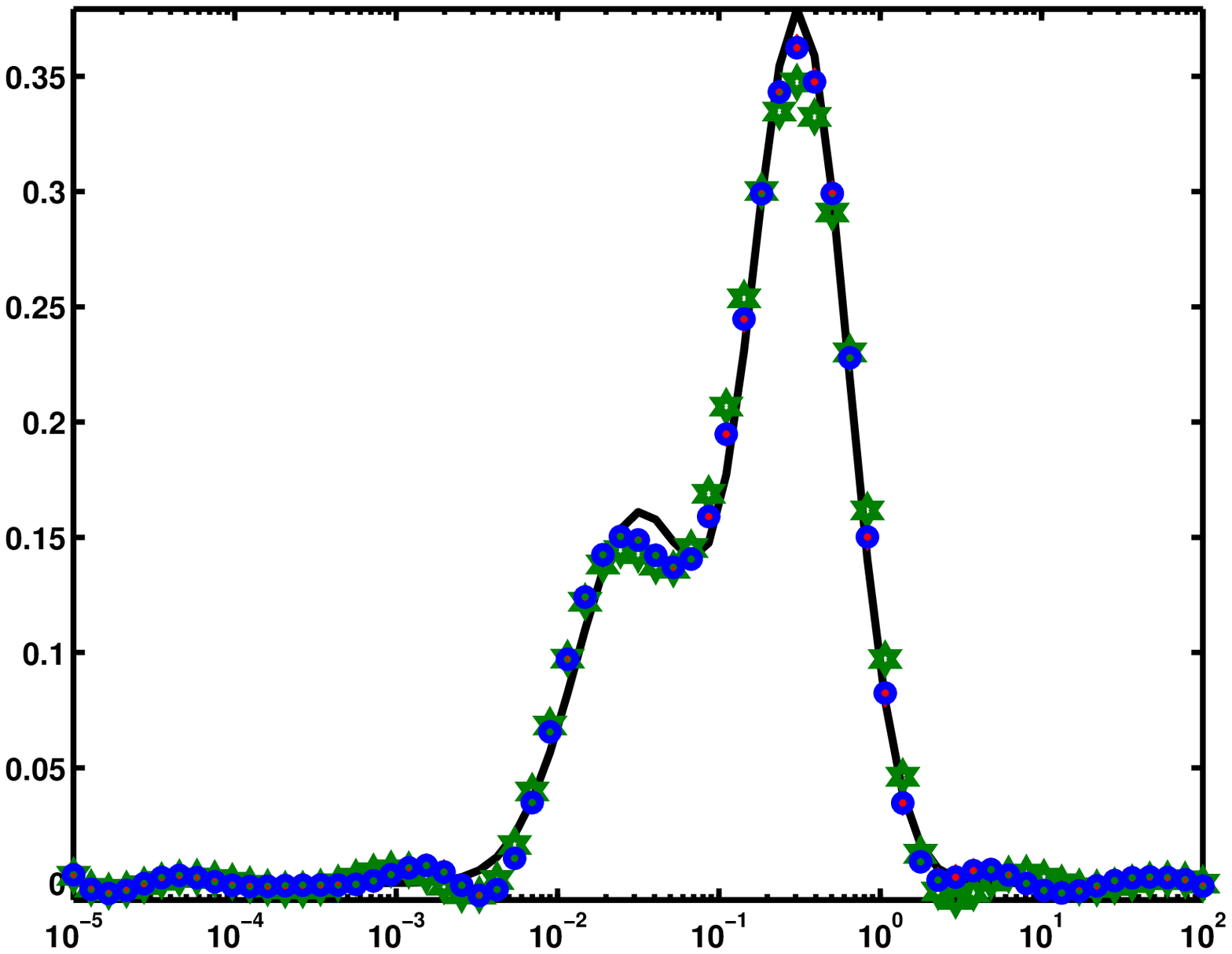}}
\subfigure[$L=L_2$]{\includegraphics[width=1.7in]{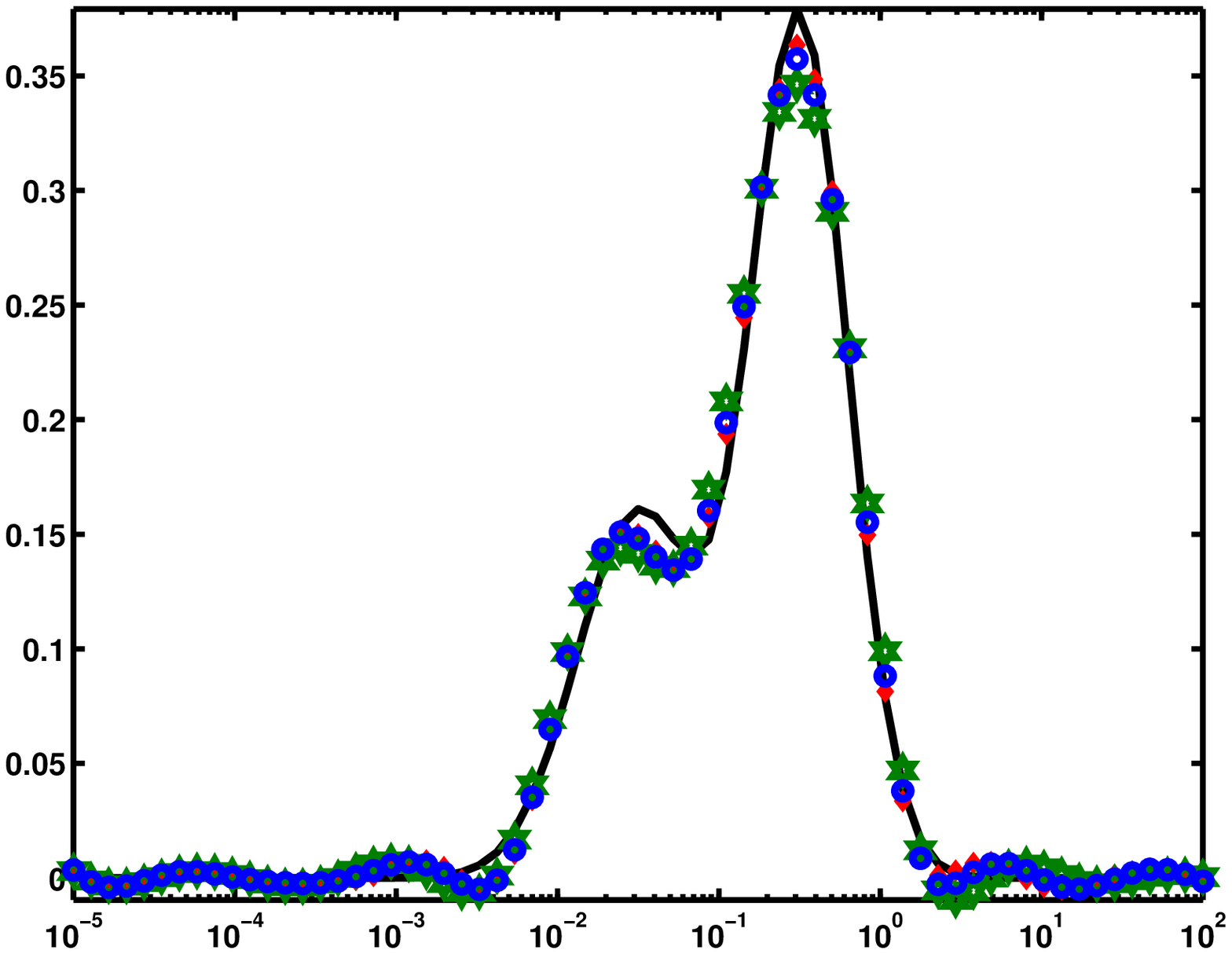}}
\caption{Mean error and example LS solutions.  $.1\%$ noise. LN-B data set matrix $A_3$ }
\label{fig-lambdachoiceLN5A3LNLS}
\end{figure}

 \begin{figure}[!h]
  \centering
\subfigure[$L=I$]{\includegraphics[width=1.7in]{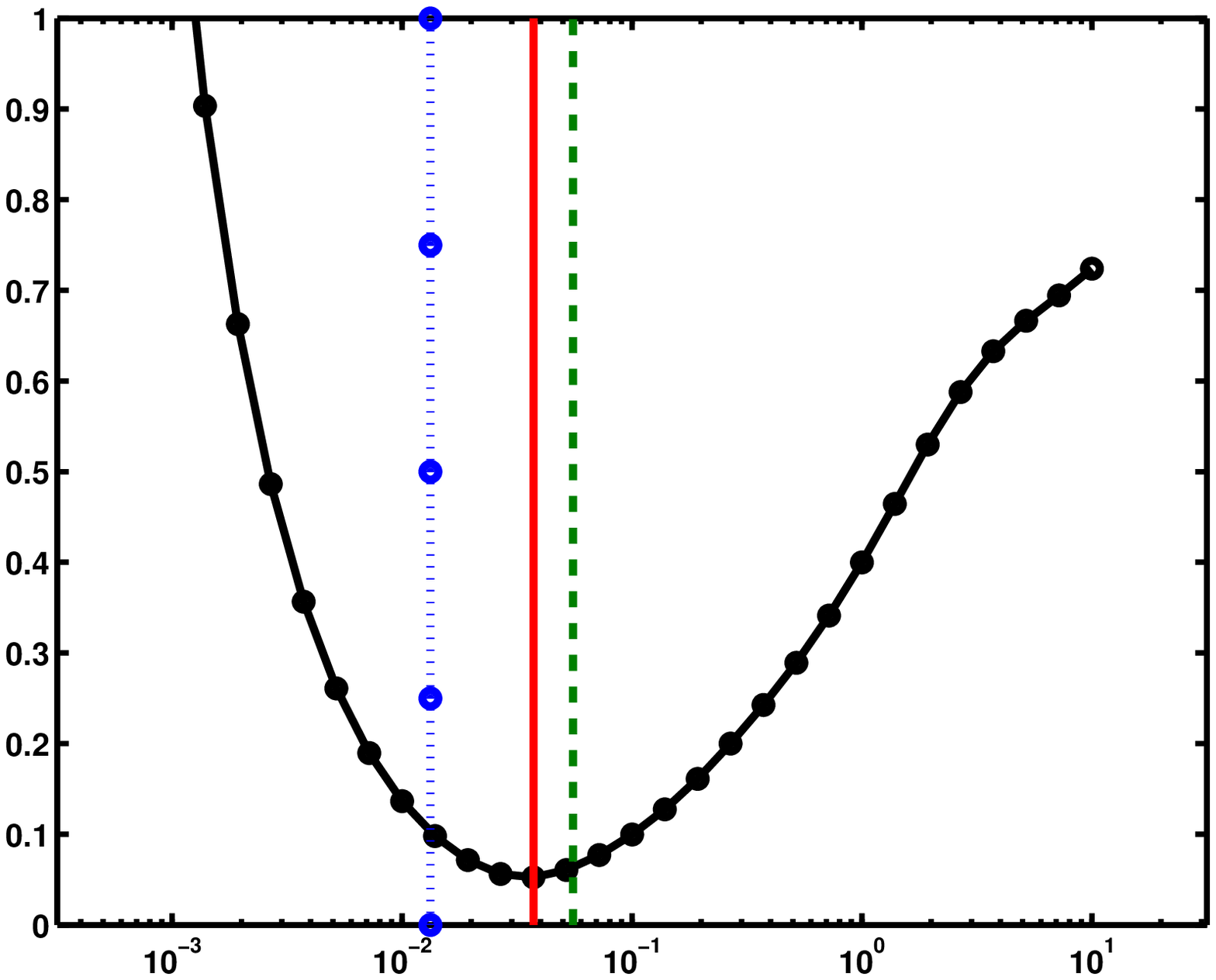}}
\subfigure[$L=L_1$]{\includegraphics[width=1.7in]{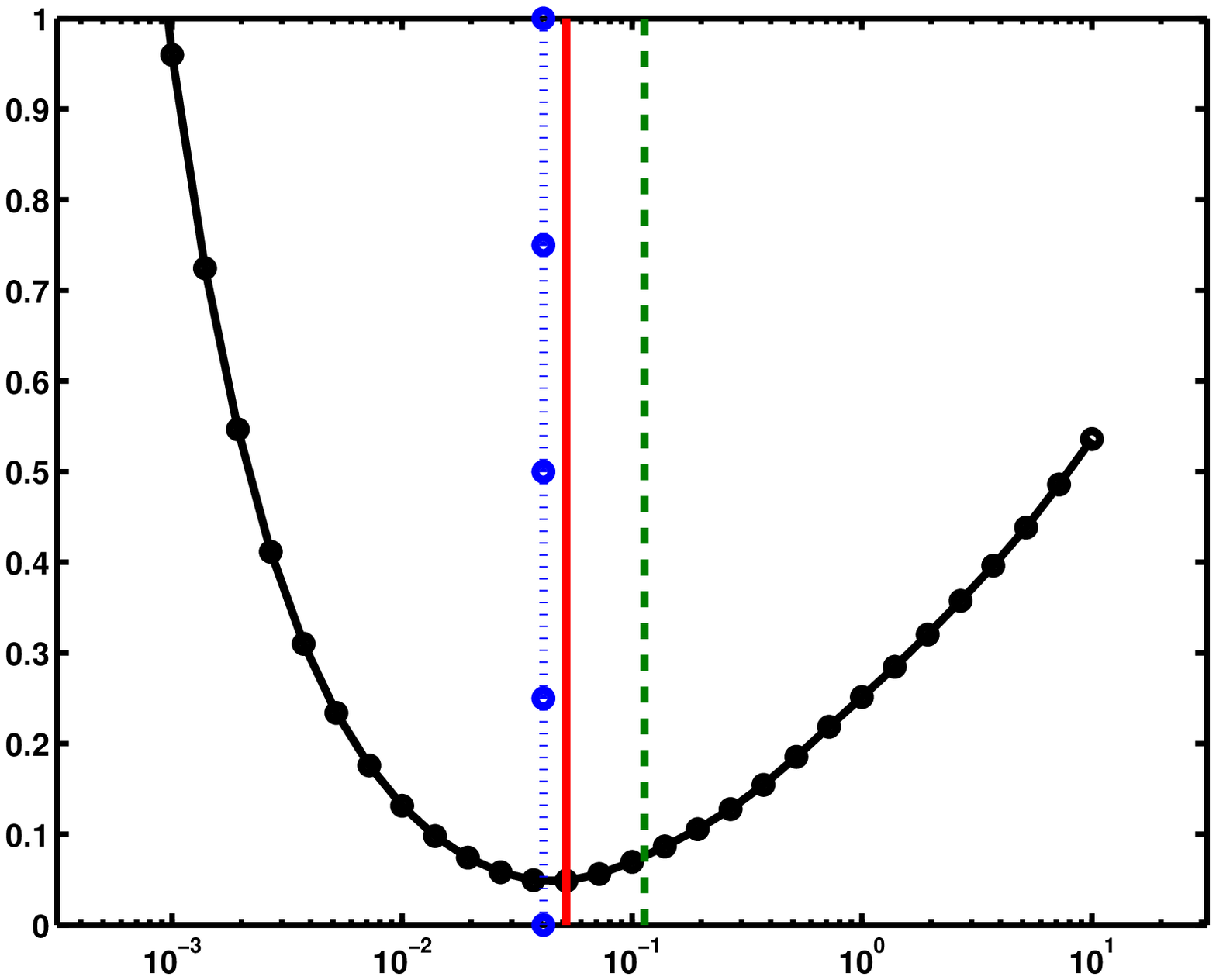}}
\subfigure[$L=L_2$]{\includegraphics[width=1.7in]{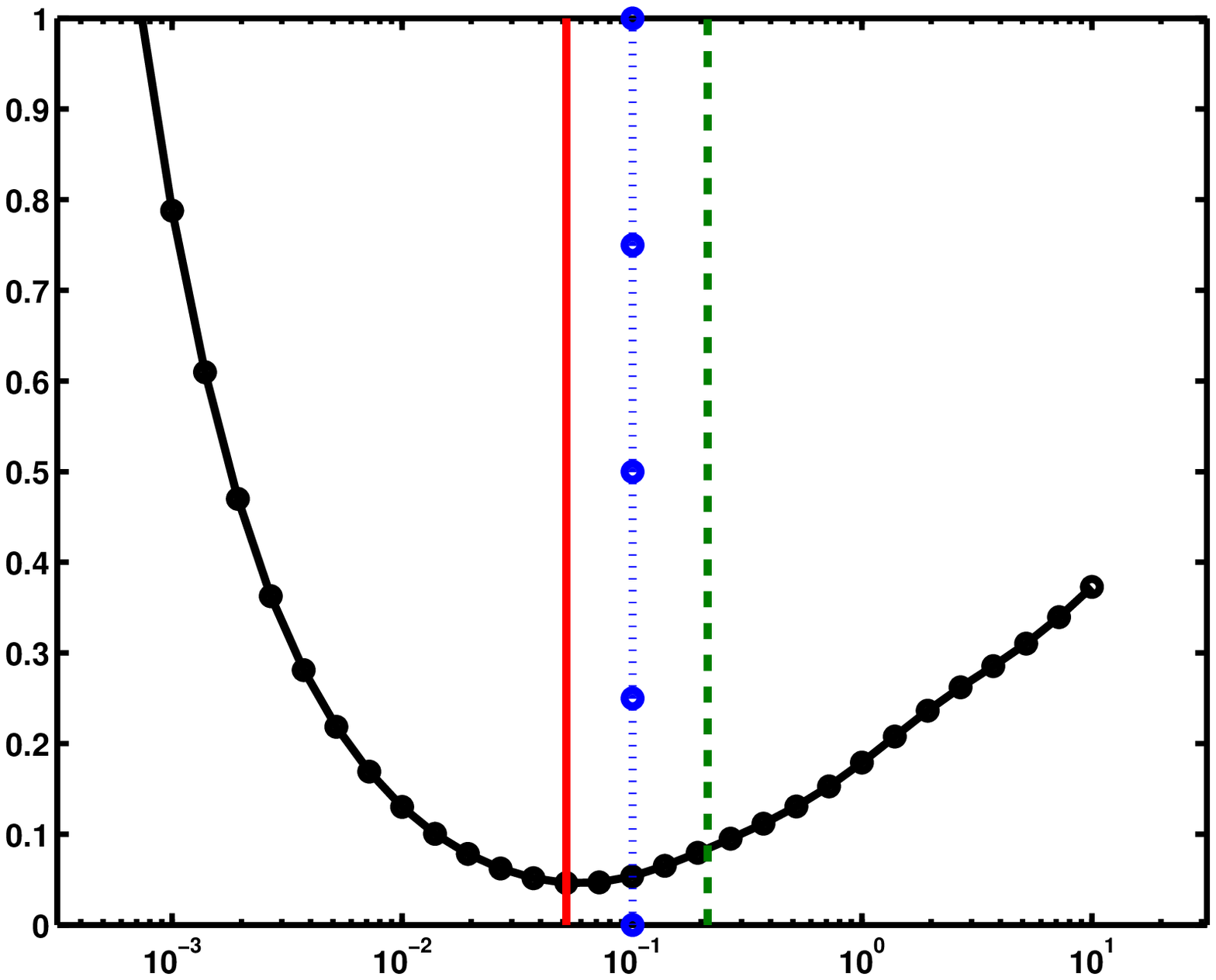}}
\subfigure[$L=I$]{\includegraphics[width=1.7in]{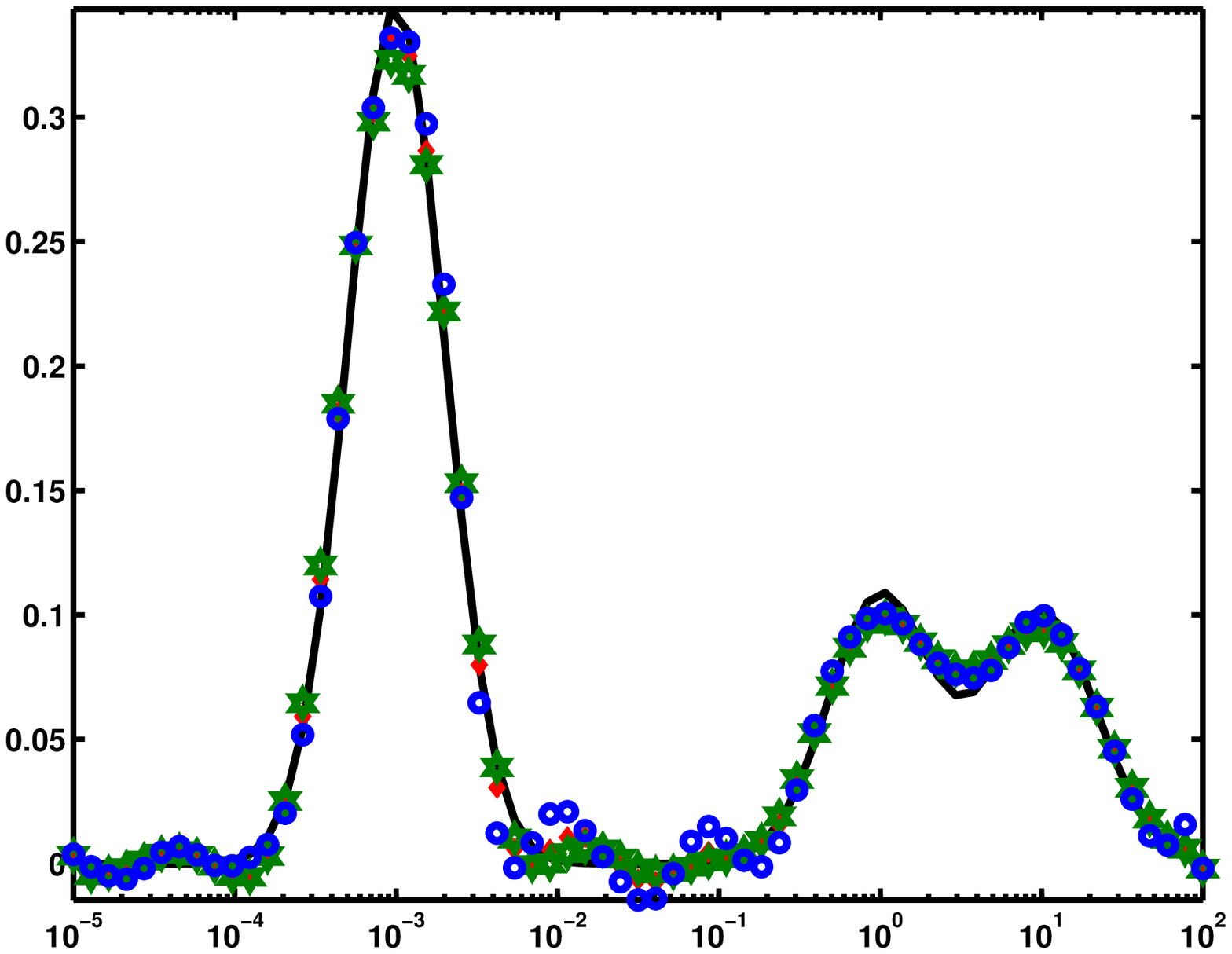}}
\subfigure[$L=L_1$]{\includegraphics[width=1.7in]{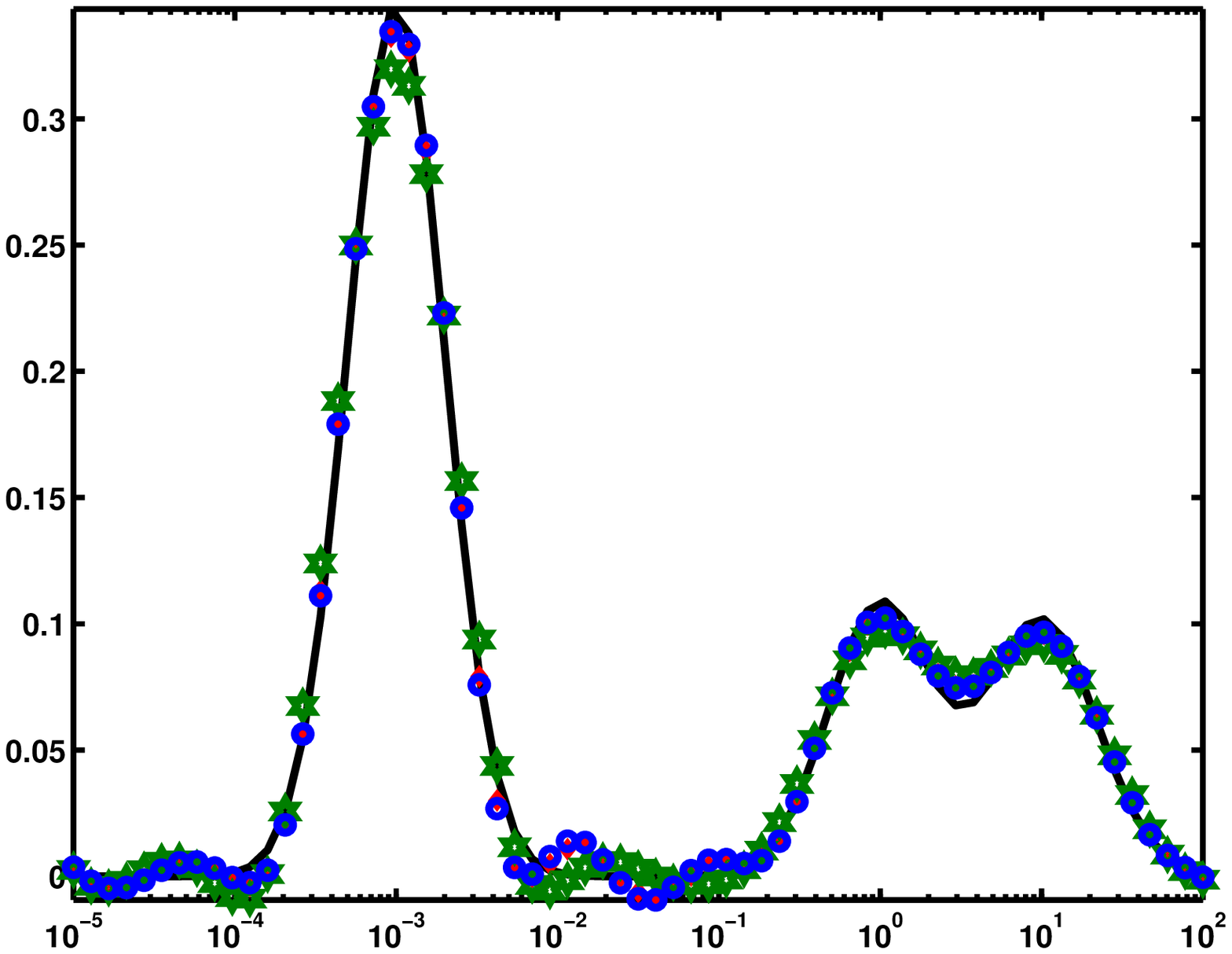}}
\subfigure[$L=L_2$]{\includegraphics[width=1.7in]{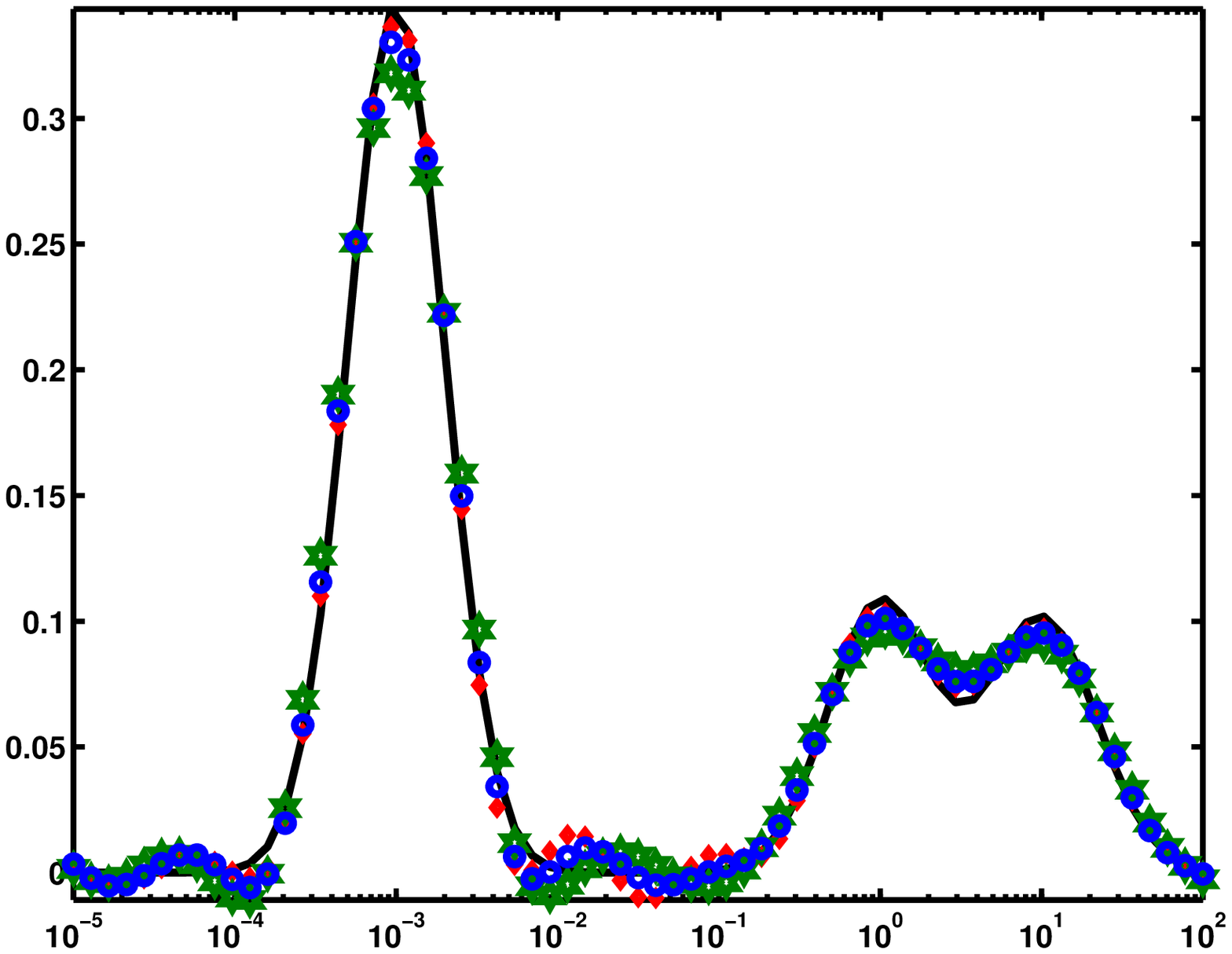}}
\caption{Mean error and example LS solutions.  $.1\%$ noise. LN-C data set matrix $A_3$}
\label{fig-lambdachoiceLN6A3LNLS}
\end{figure}
\clearpage

\subsection{Results using LS $A_3$ Noise level $5\%$}
 \begin{figure}[!h]
 \centering
\subfigure[$L=I$]{\includegraphics[width=1.7in]{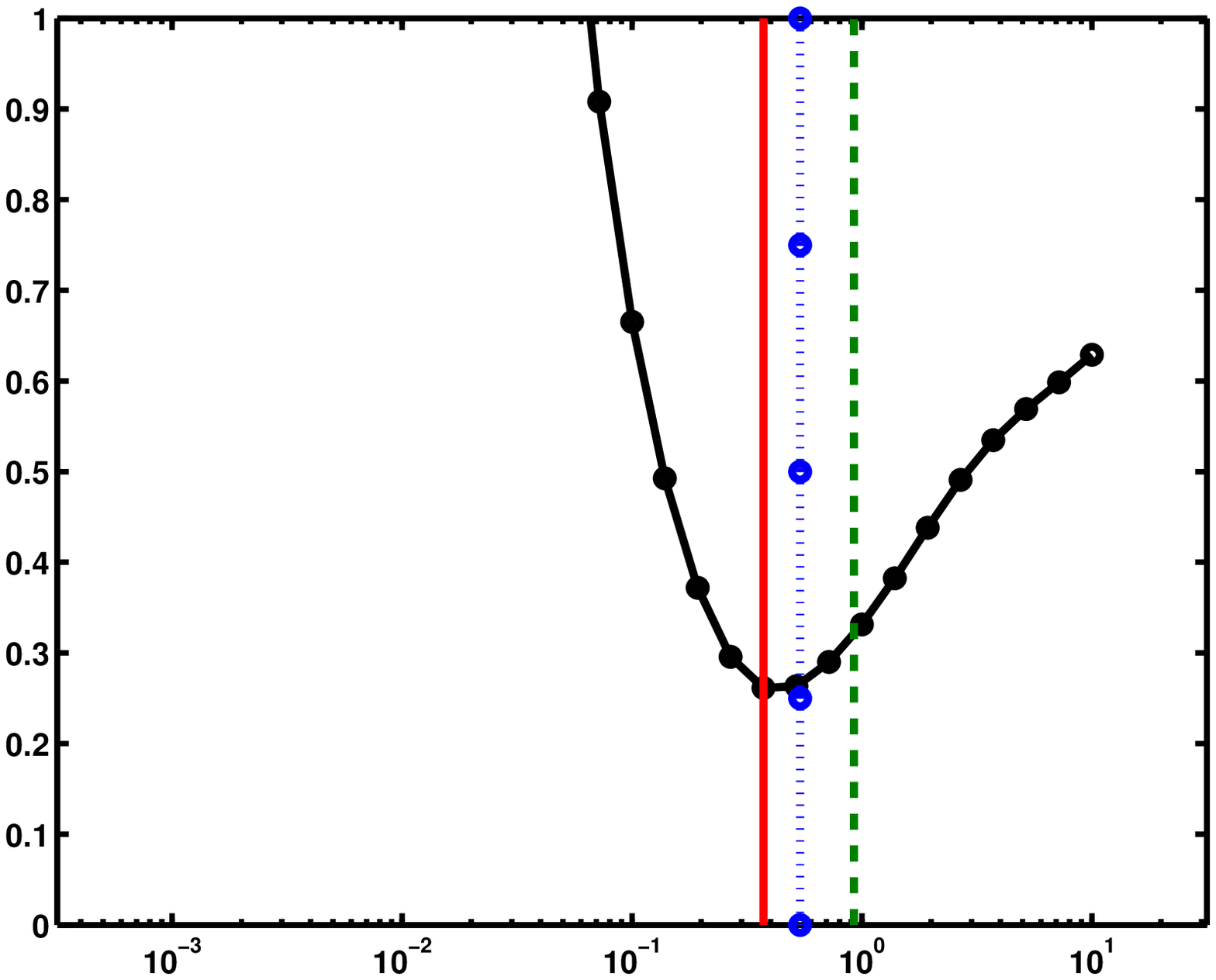}}
\subfigure[$L=L_1$]{\includegraphics[width=1.7in]{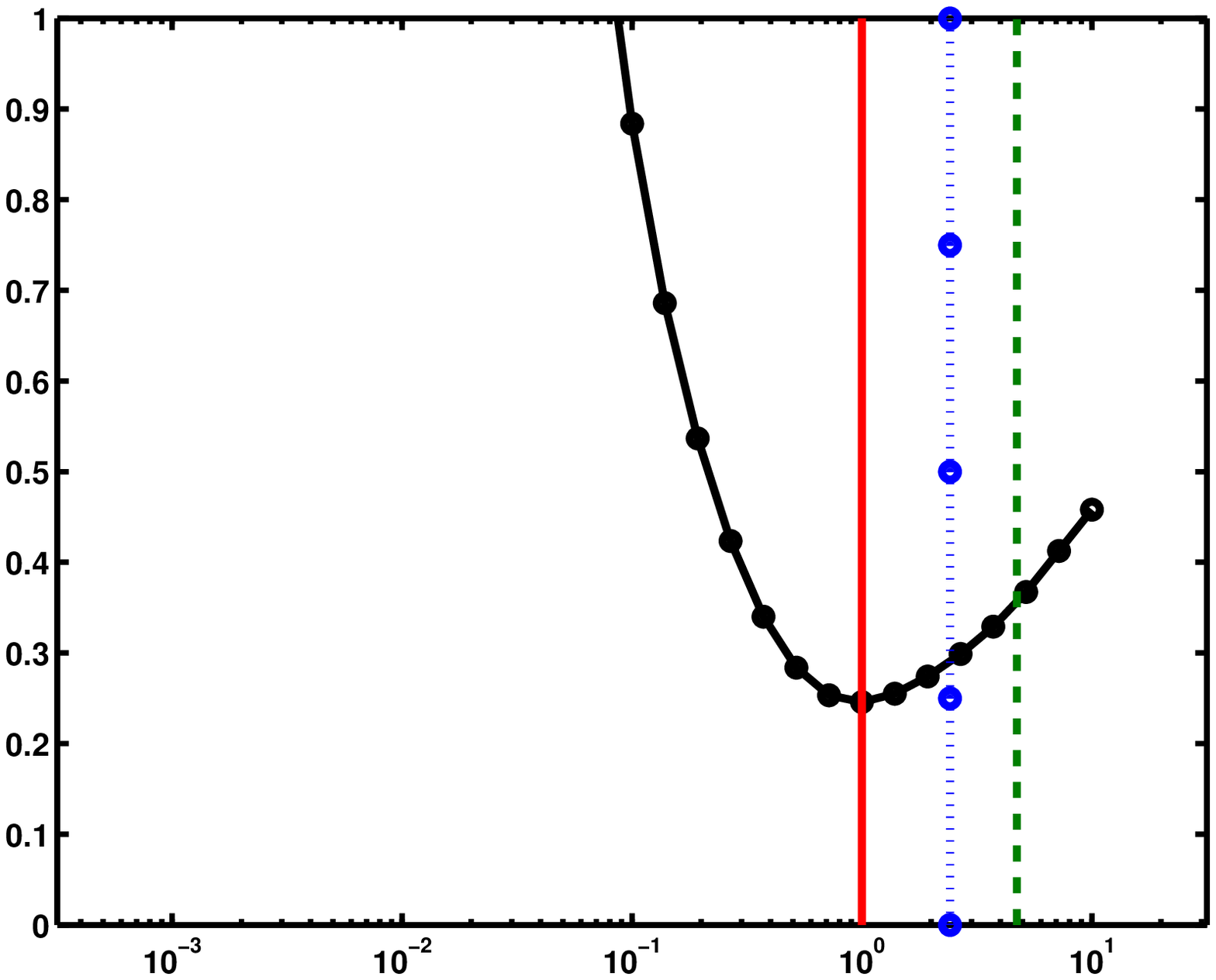}}
\subfigure[$L=L_2$]{\includegraphics[width=1.7in]{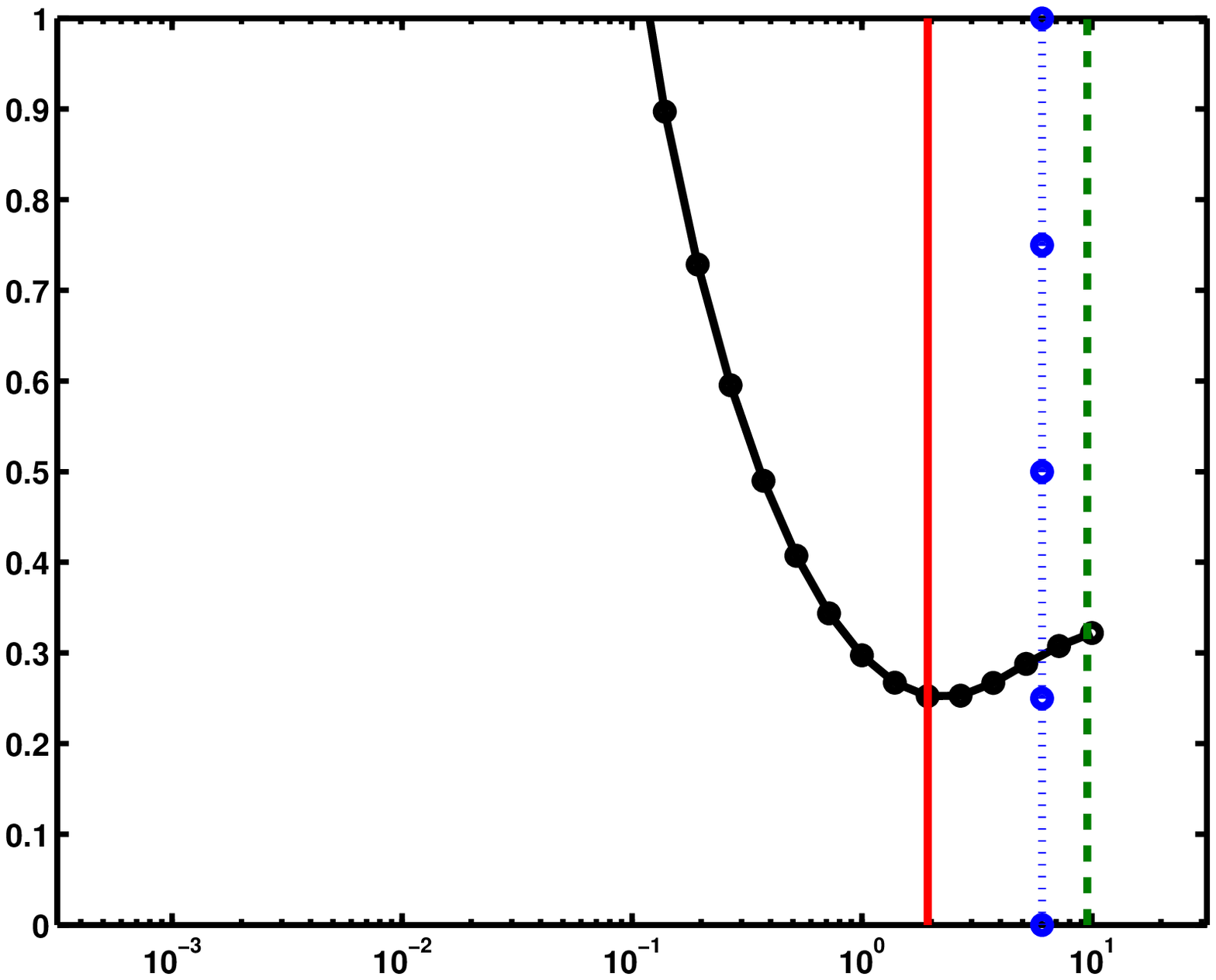}}
\subfigure[$L=I$]{\includegraphics[width=1.7in]{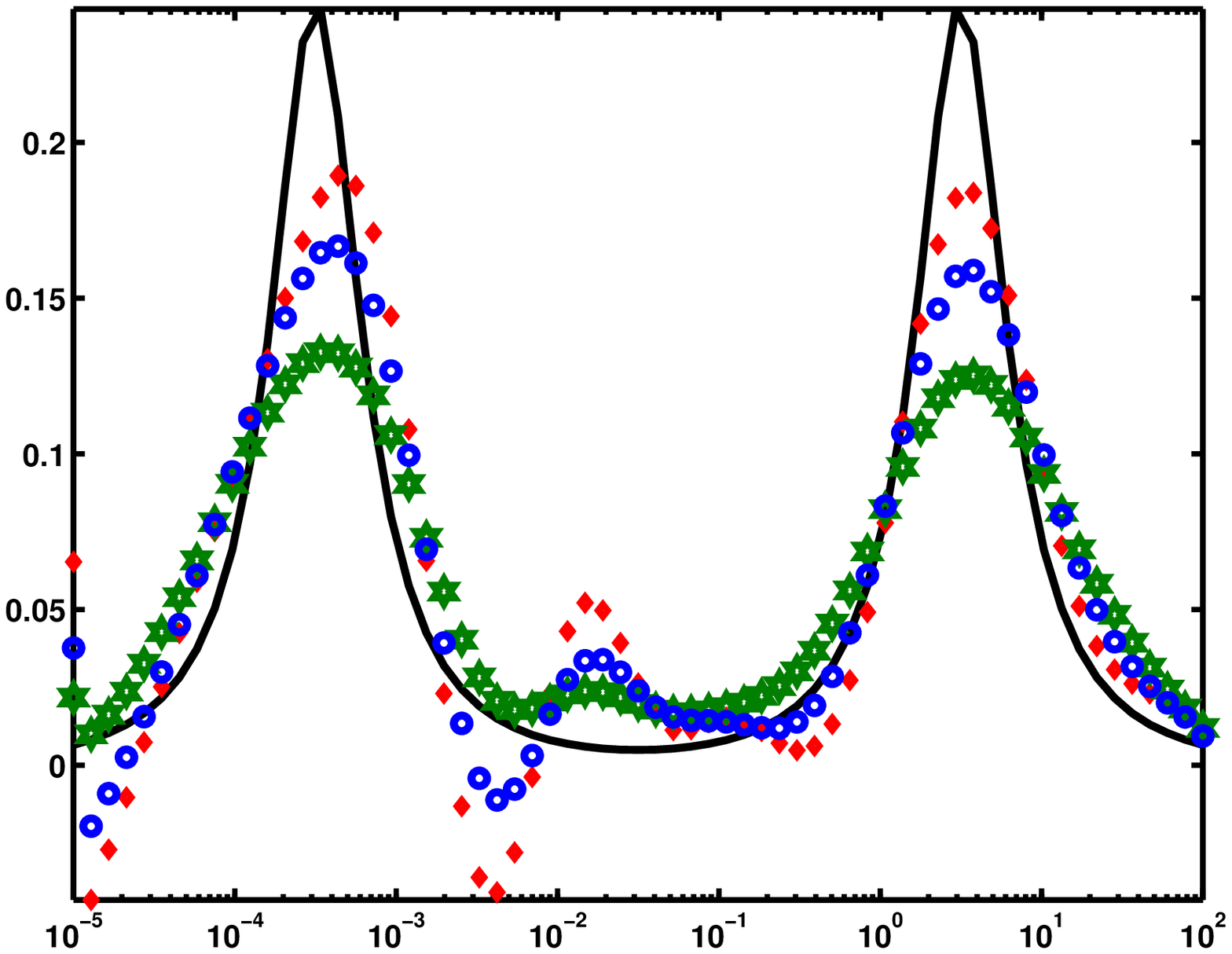}}
\subfigure[$L=L_1$]{\includegraphics[width=1.7in]{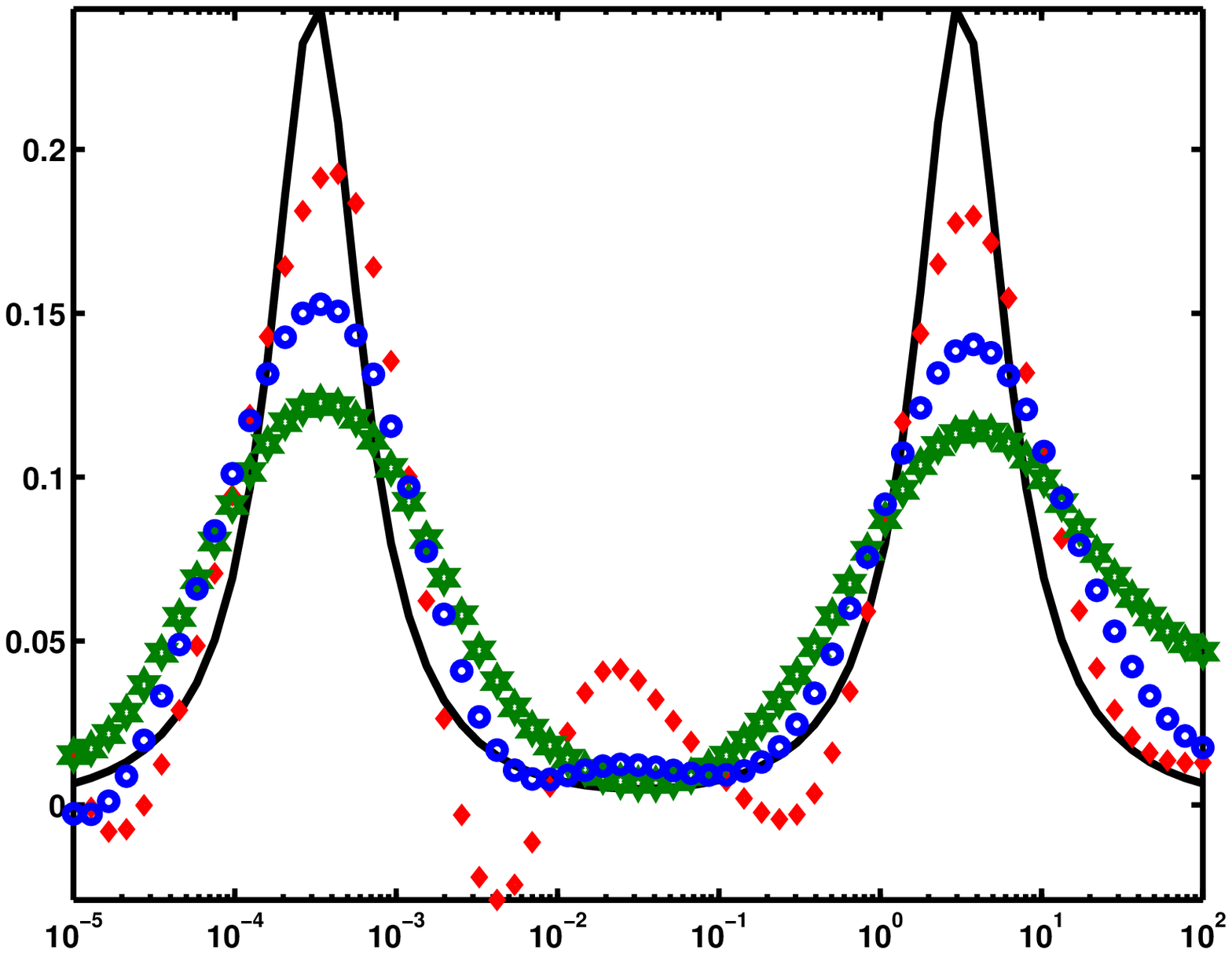}}
\subfigure[$L=L_2$]{\includegraphics[width=1.7in]{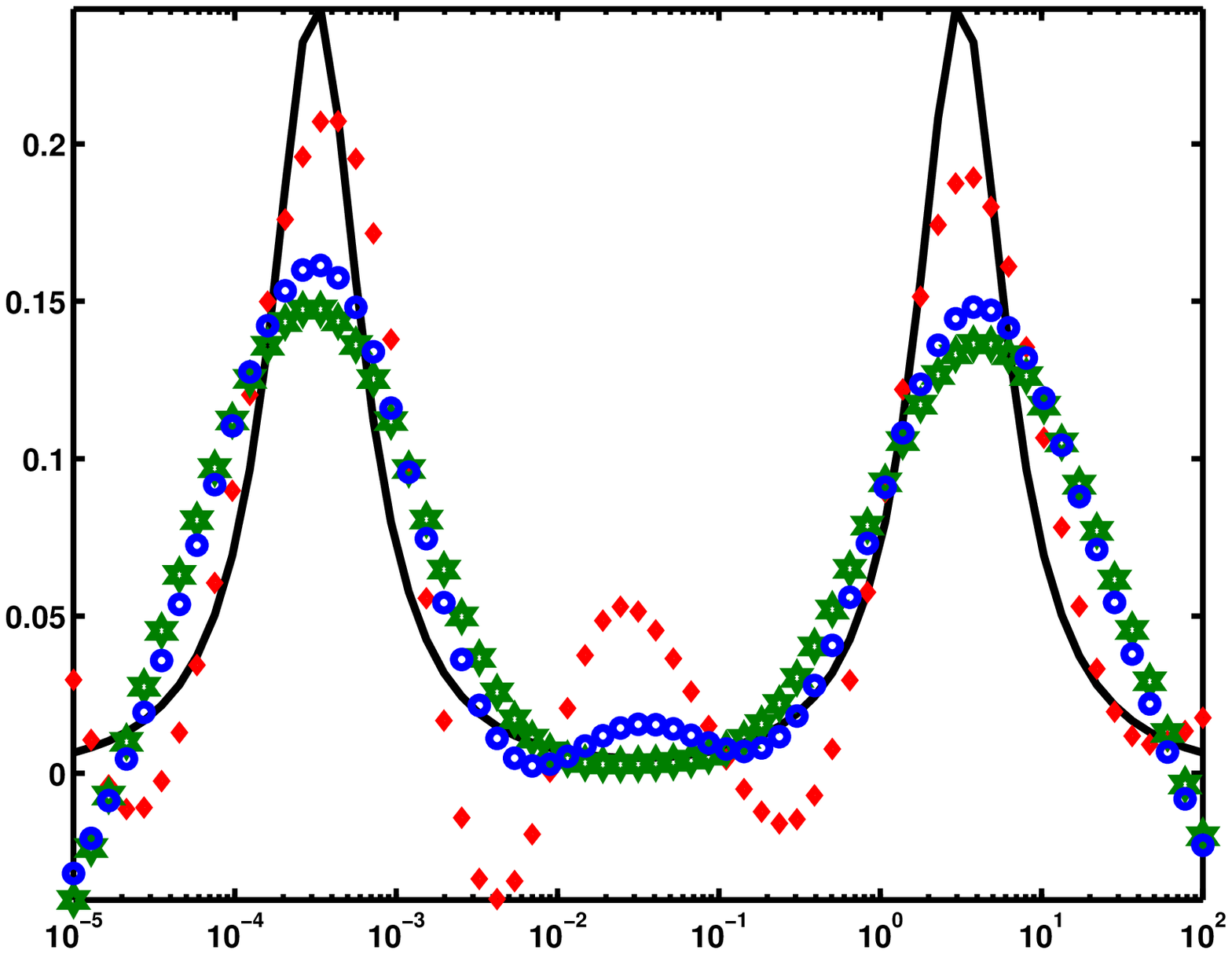}}
\caption{Mean error and example LS solutions.  $5\%$ noise. RQ-A data set matrix $A_3$}
\label{fig-lambdachoiceRQ1A3HNLS}
\end{figure}

 \begin{figure}[!h]
  \centering
\subfigure[$L=I$]{\includegraphics[width=1.7in]{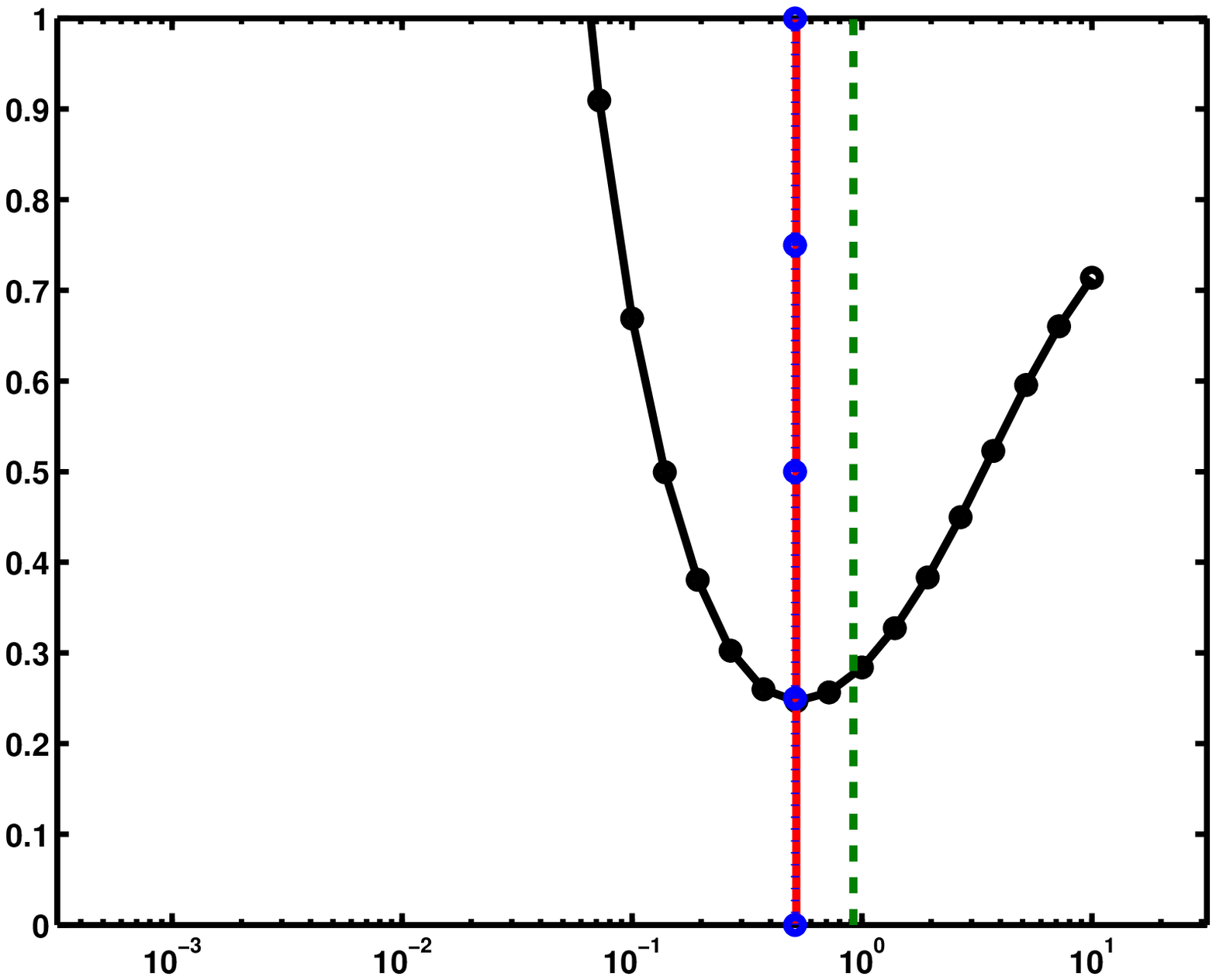}}
\subfigure[$L=L_1$]{\includegraphics[width=1.7in]{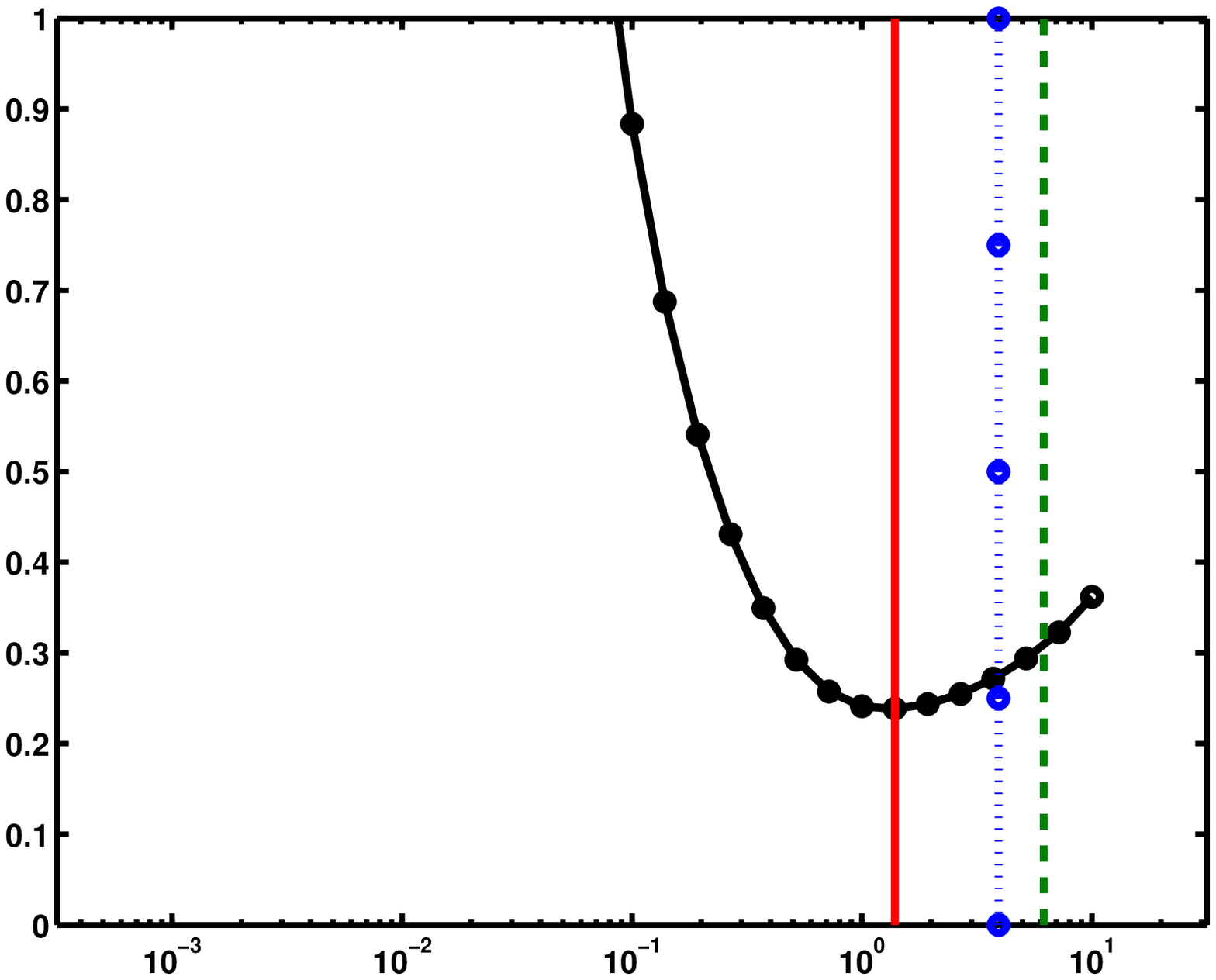}}
\subfigure[$L=L_2$]{\includegraphics[width=1.7in]{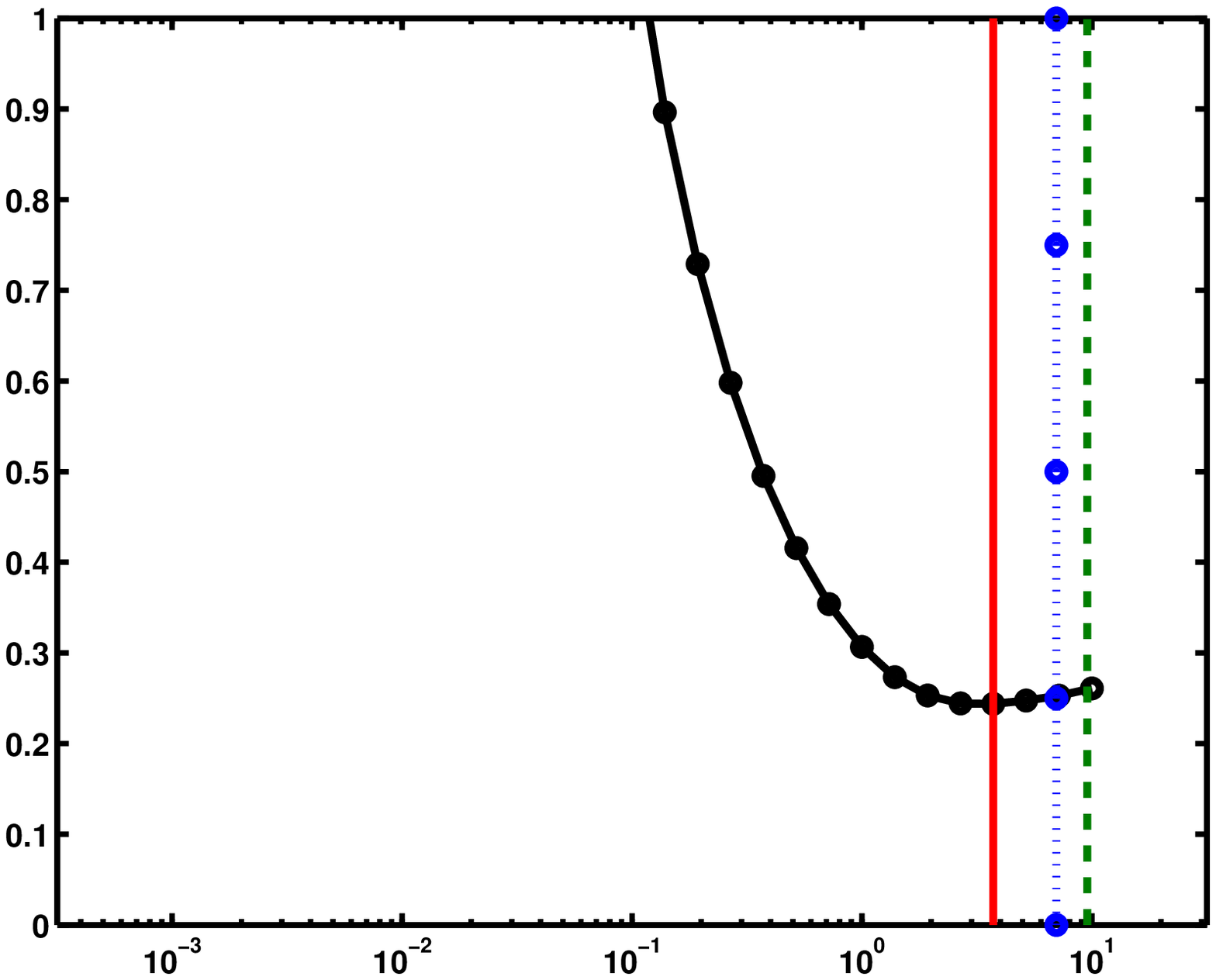}}
\subfigure[$L=I$]{\includegraphics[width=1.7in]{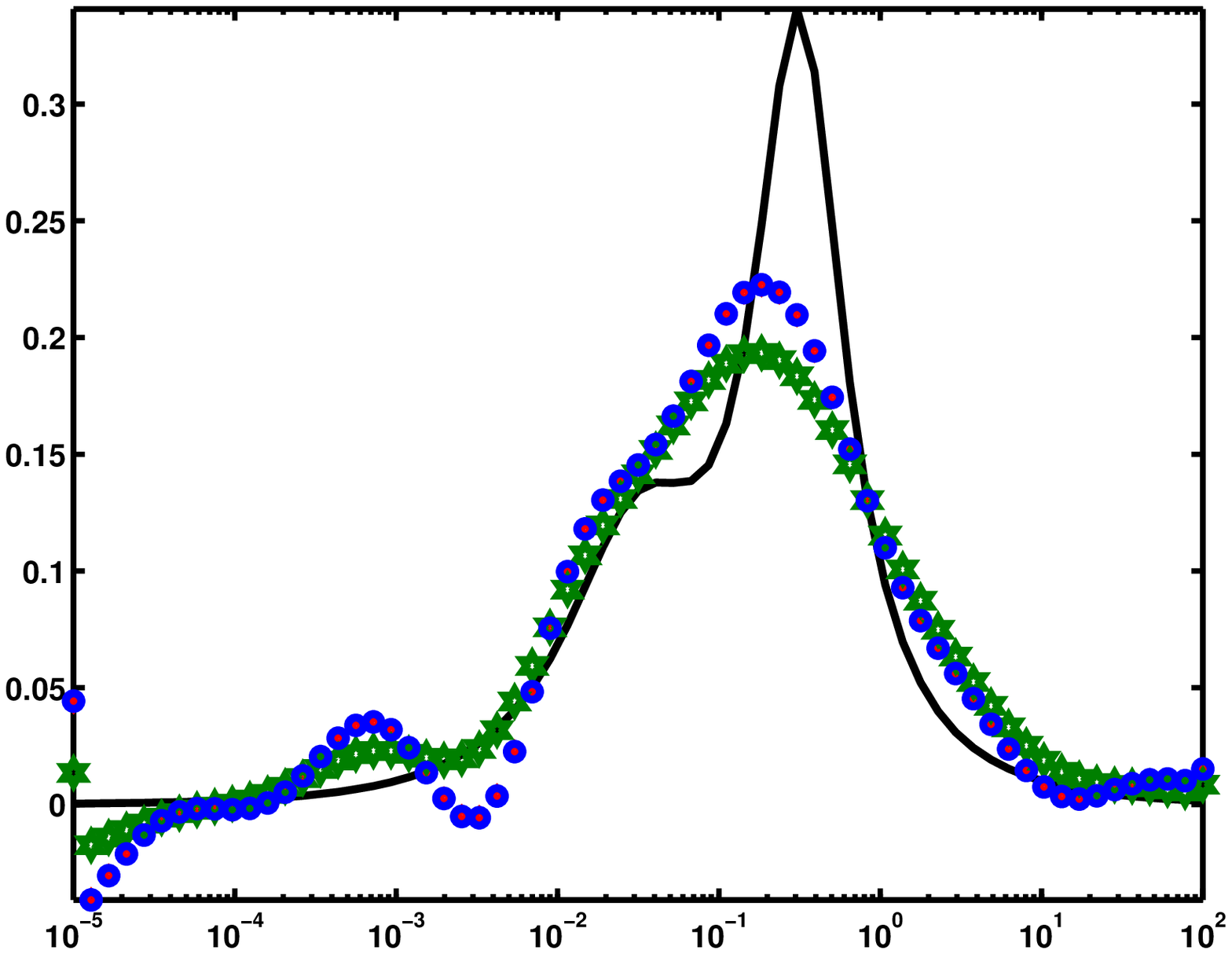}}
\subfigure[$L=L_1$]{\includegraphics[width=1.7in]{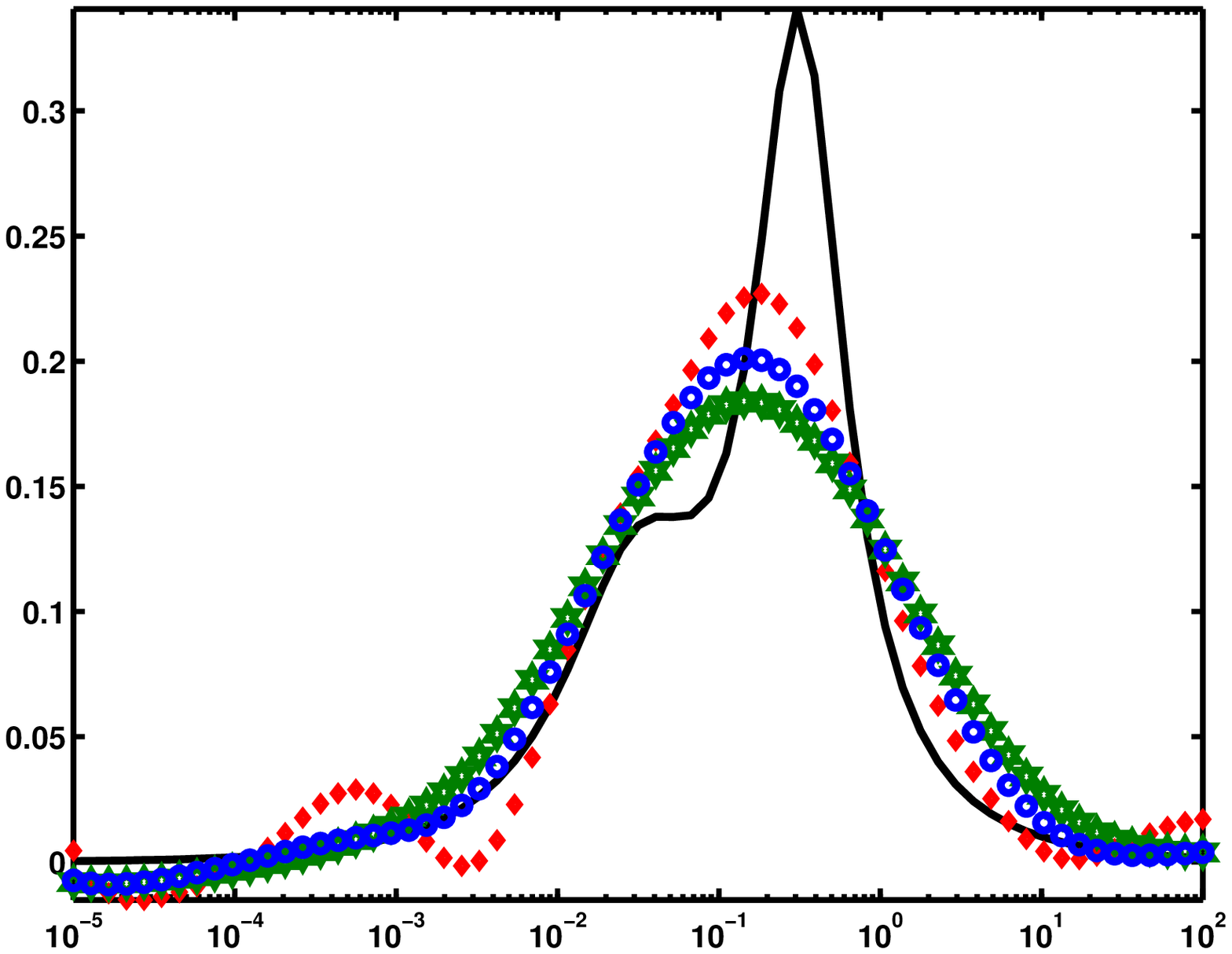}}
\subfigure[$L=L_2$]{\includegraphics[width=1.7in]{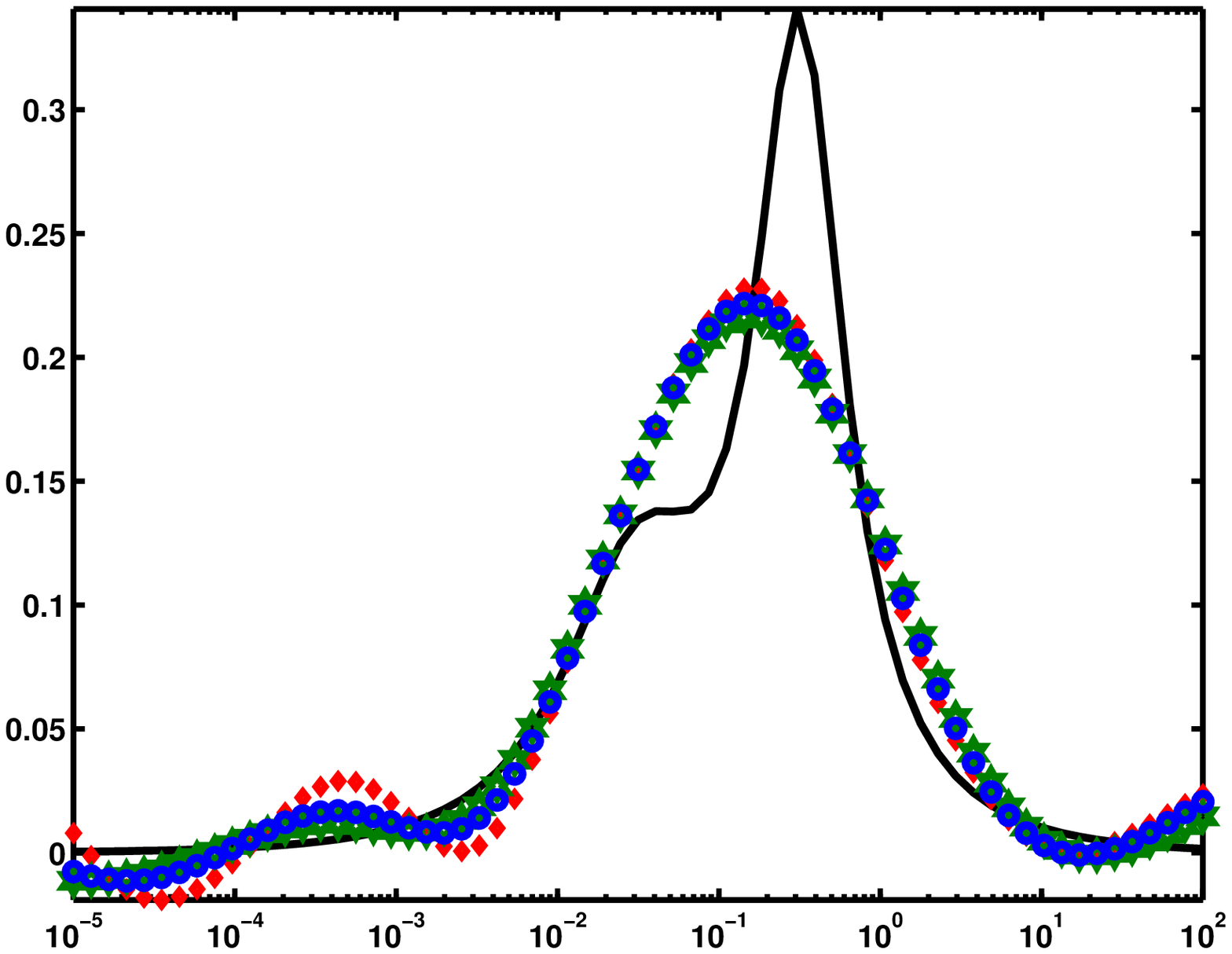}}
\caption{Mean error and example LS solutions. $5\%$ noise. RQ-B data set matrix $A_3$}
\label{fig-lambdachoiceRQ5A3HNLS}
\end{figure}

 \begin{figure}[!h]
  \centering
\subfigure[$L=I$]{\includegraphics[width=1.7in]{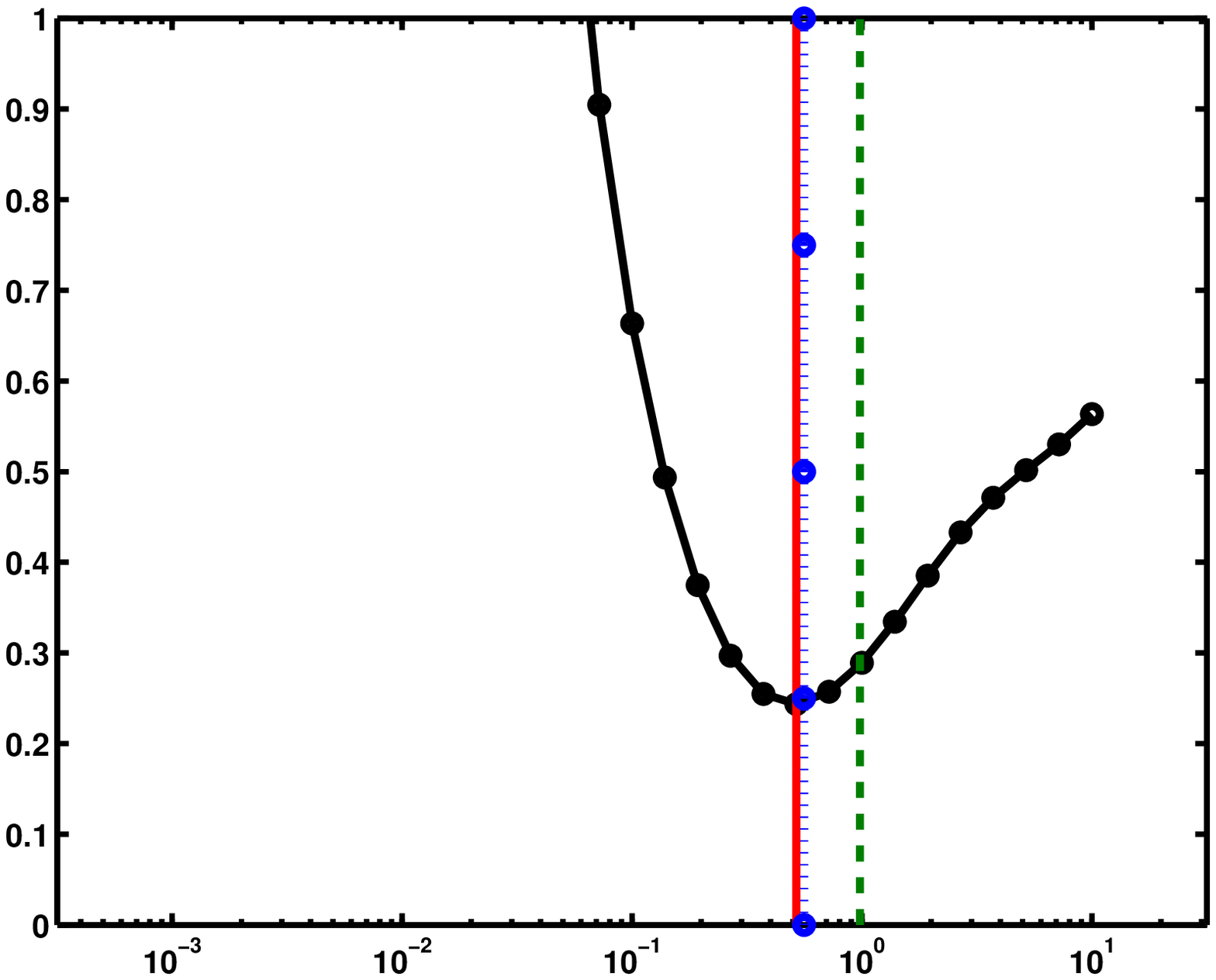}}
\subfigure[$L=L_1$]{\includegraphics[width=1.7in]{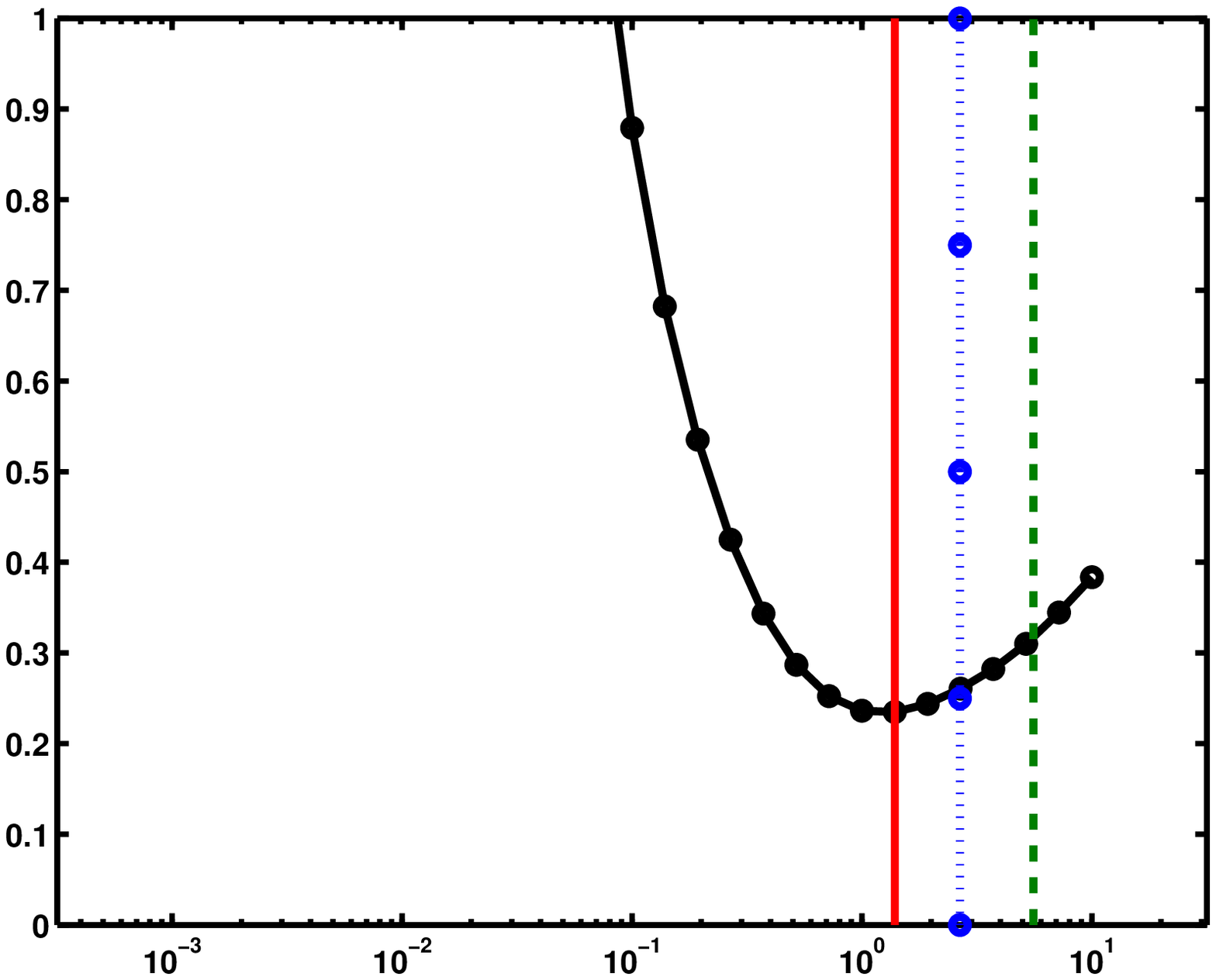}}
\subfigure[$L=L_2$]{\includegraphics[width=1.7in]{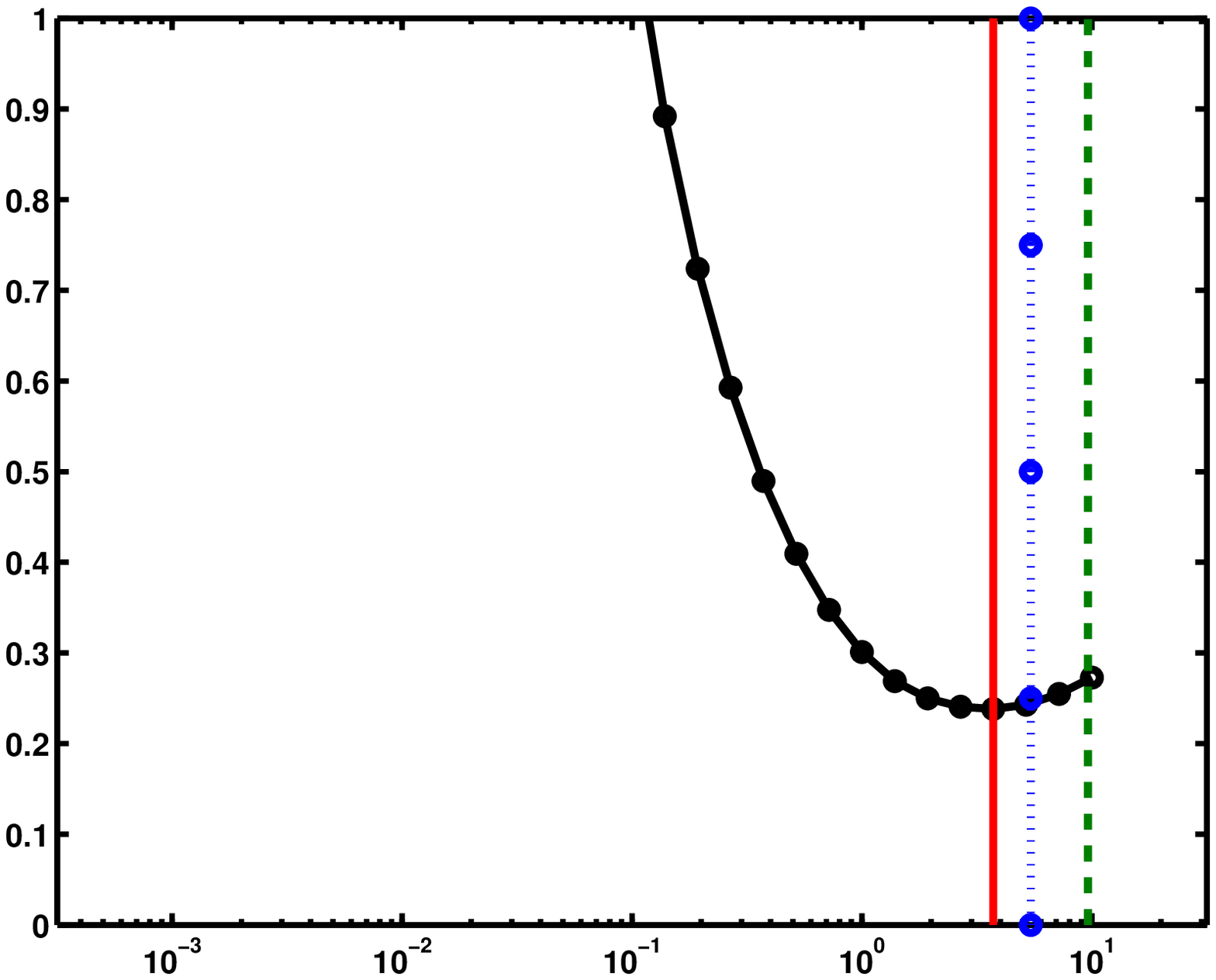}}
\subfigure[$L=I$]{\includegraphics[width=1.7in]{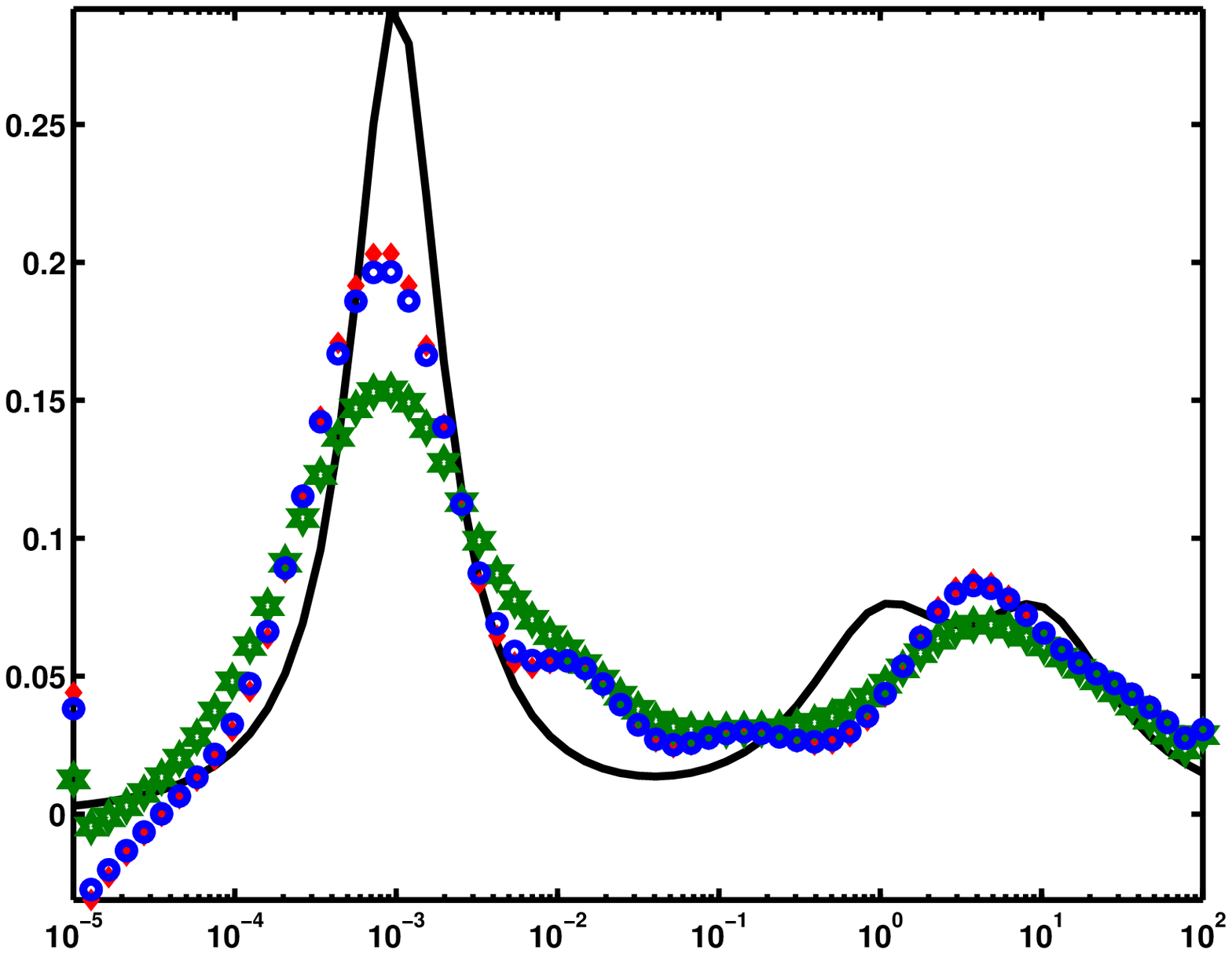}}
\subfigure[$L=L_1$]{\includegraphics[width=1.7in]{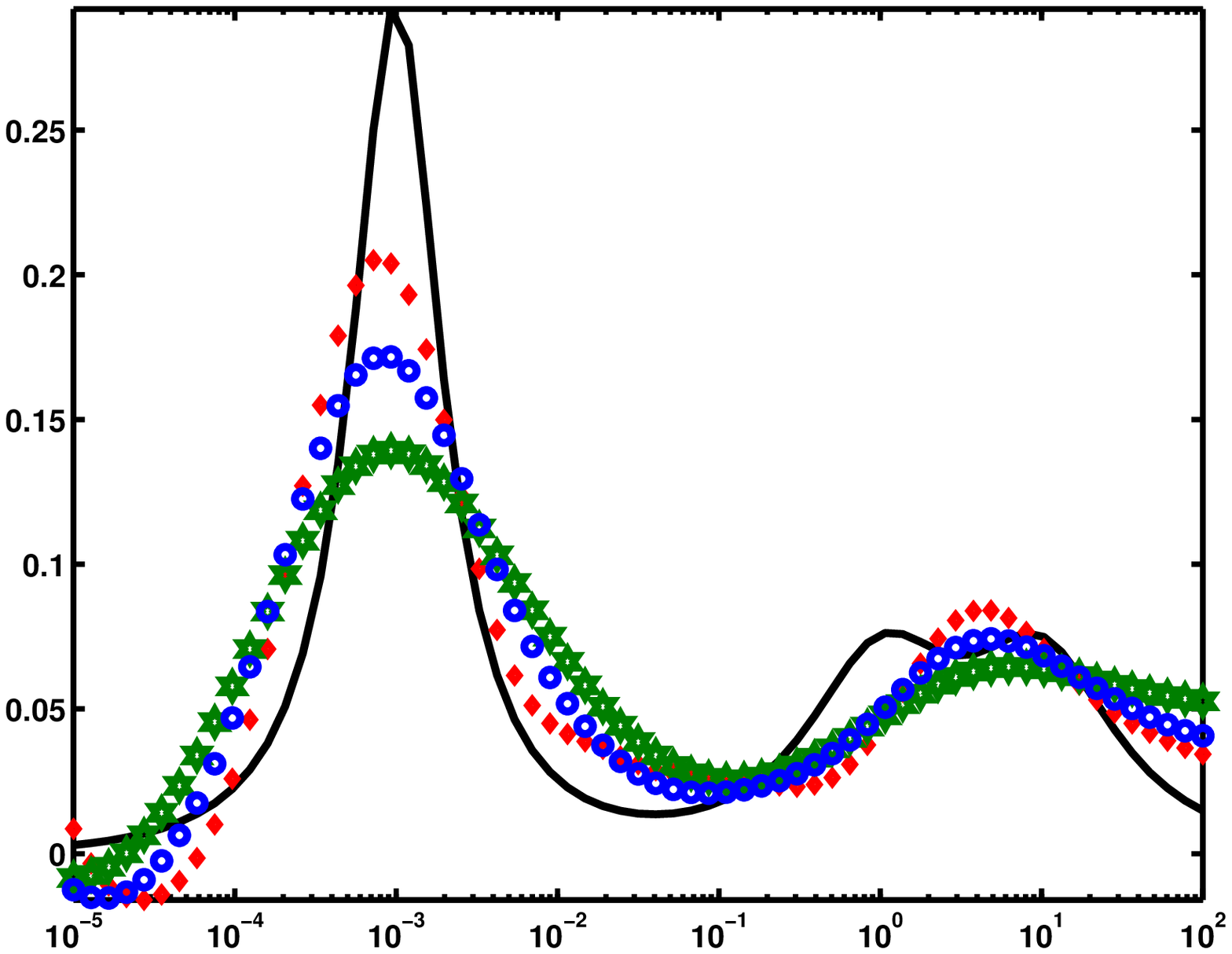}}
\subfigure[$L=L_2$]{\includegraphics[width=1.7in]{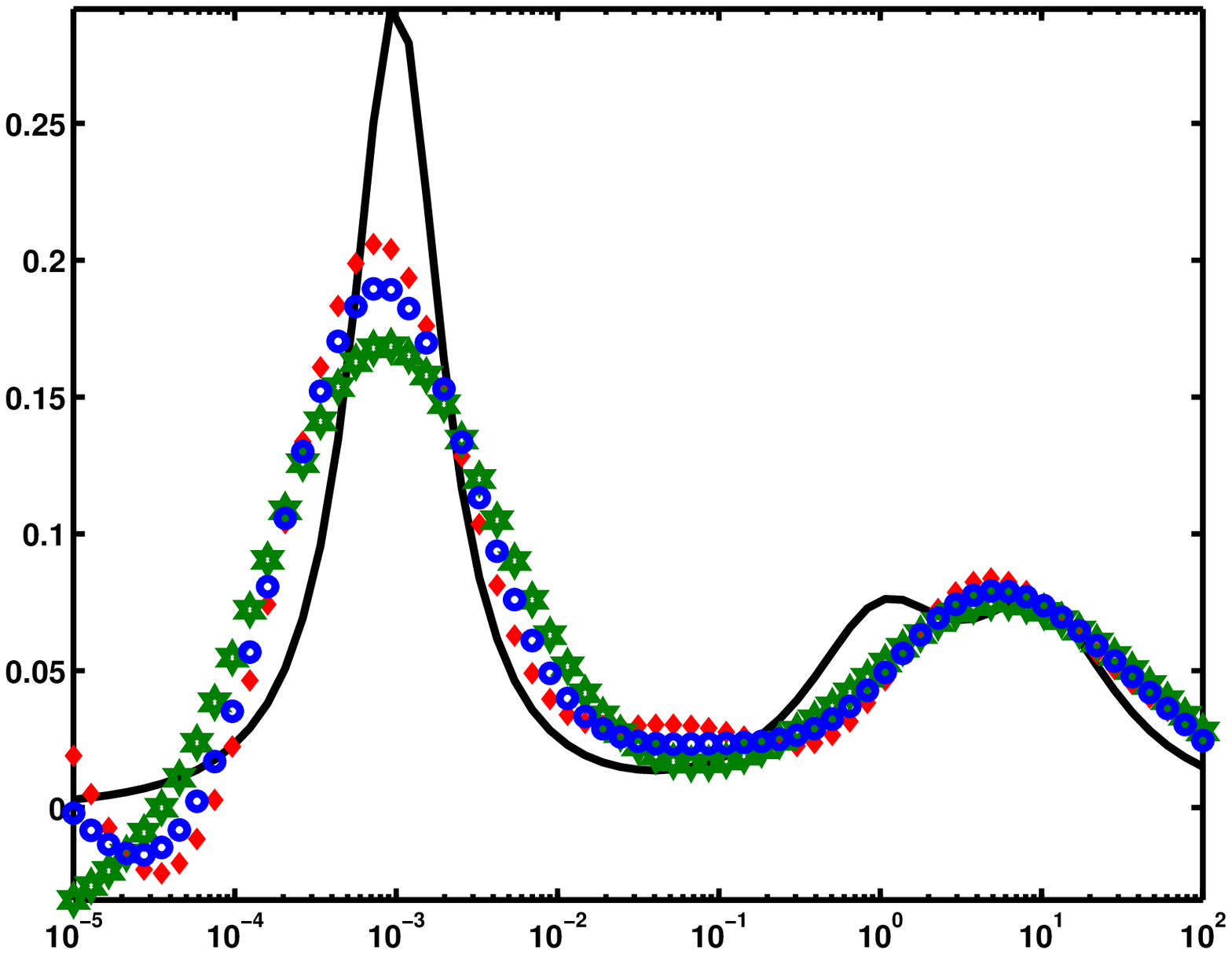}}
\caption{Mean error and example LS solutions.  $5\%$ noise. RQ-C data set matrix $A_3$}
\label{fig-lambdachoiceRQ6A3HNLS}
\end{figure}

 \begin{figure}[!h]
  \centering
\subfigure[$L=I$]{\includegraphics[width=1.7in]{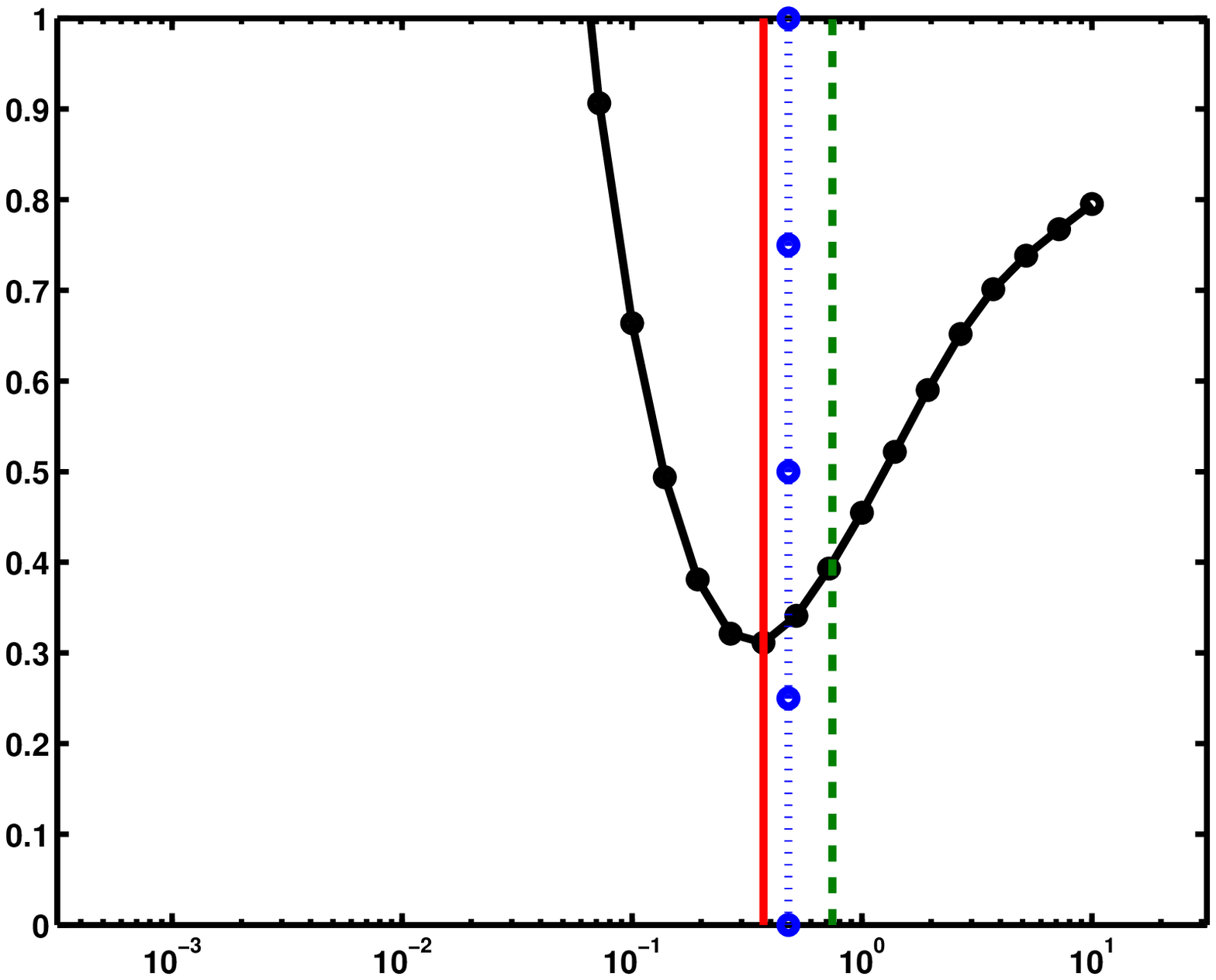}}
\subfigure[$L=L_1$]{\includegraphics[width=1.7in]{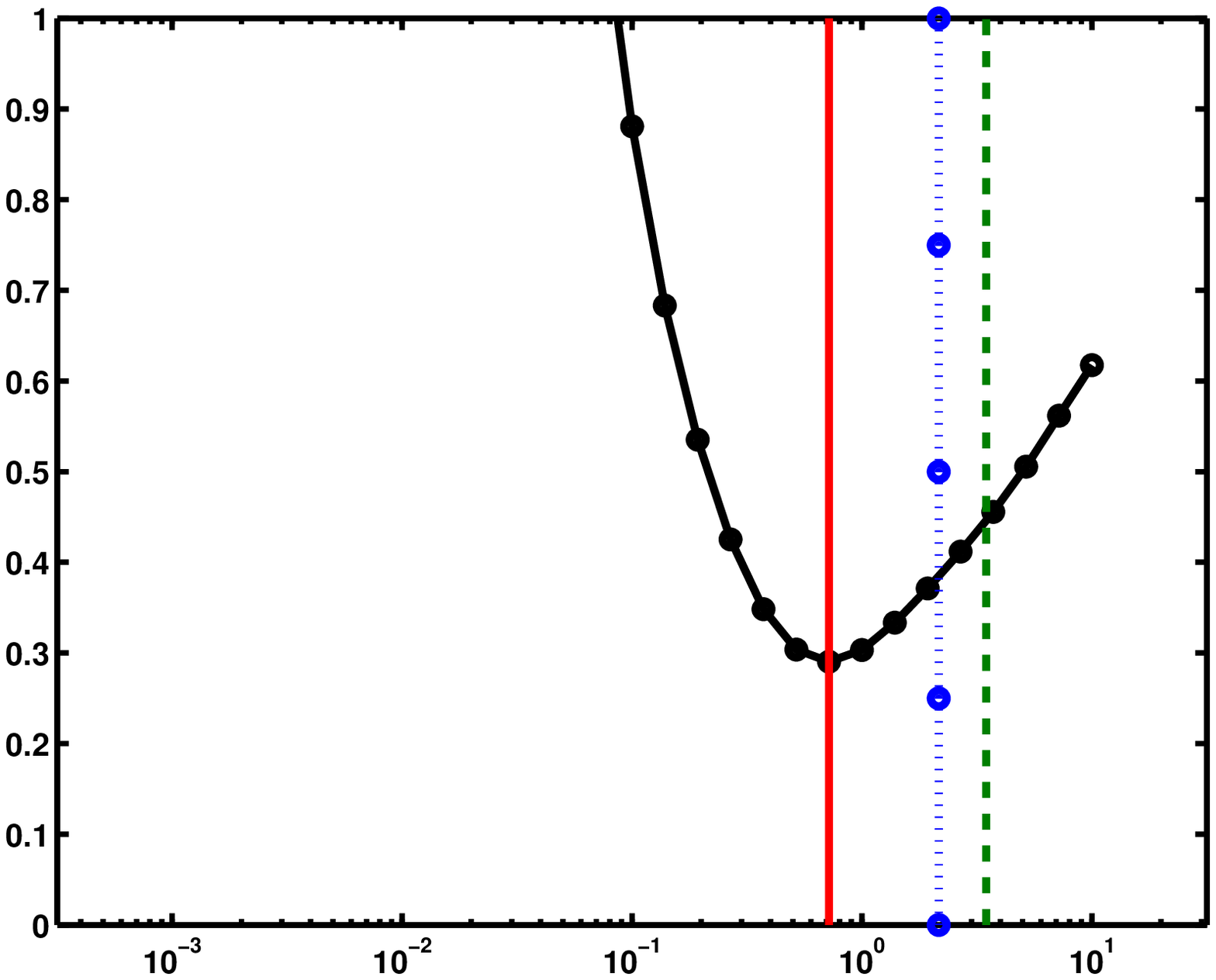}}
\subfigure[$L=L_2$]{\includegraphics[width=1.7in]{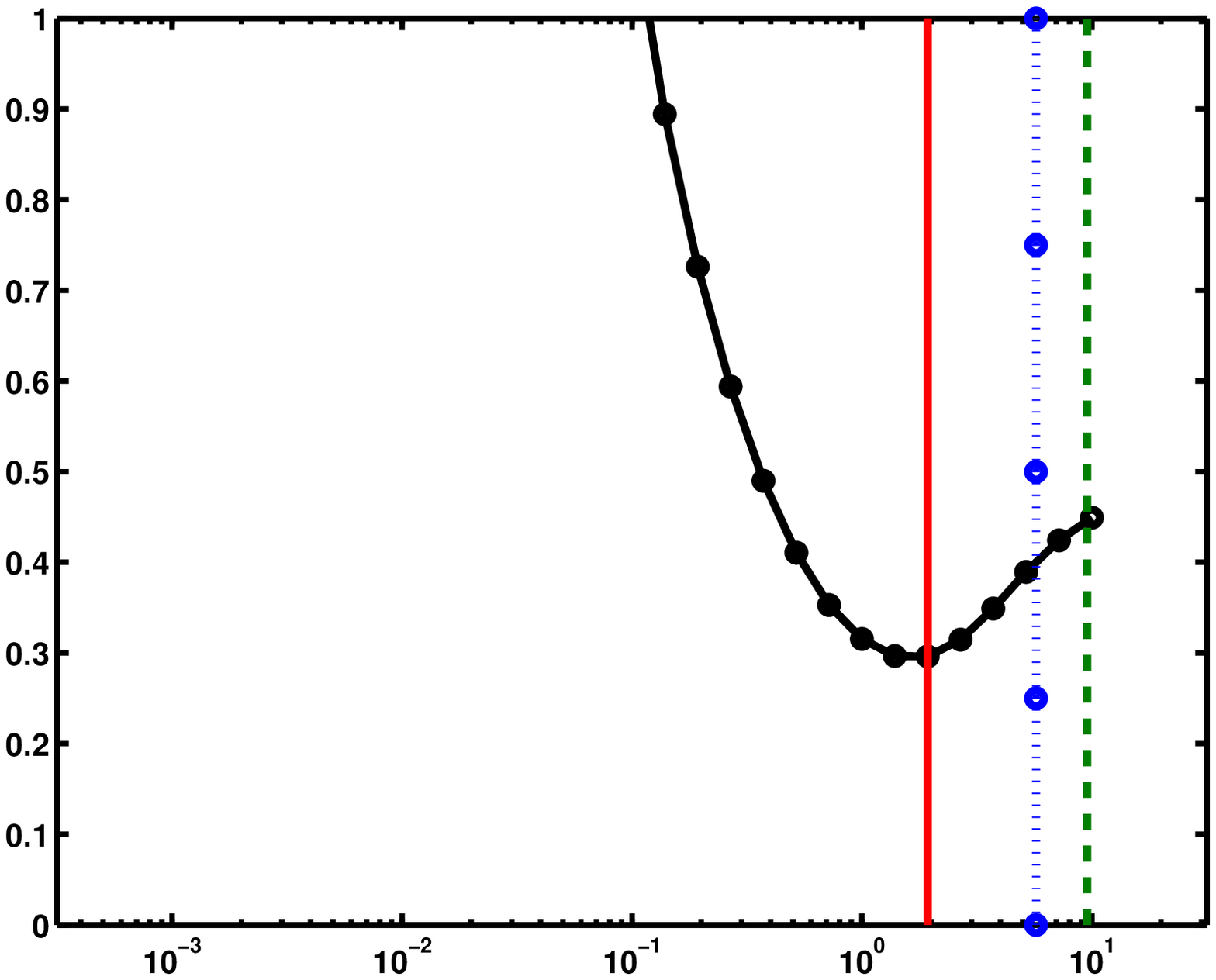}}
\subfigure[$L=I$]{\includegraphics[width=1.7in]{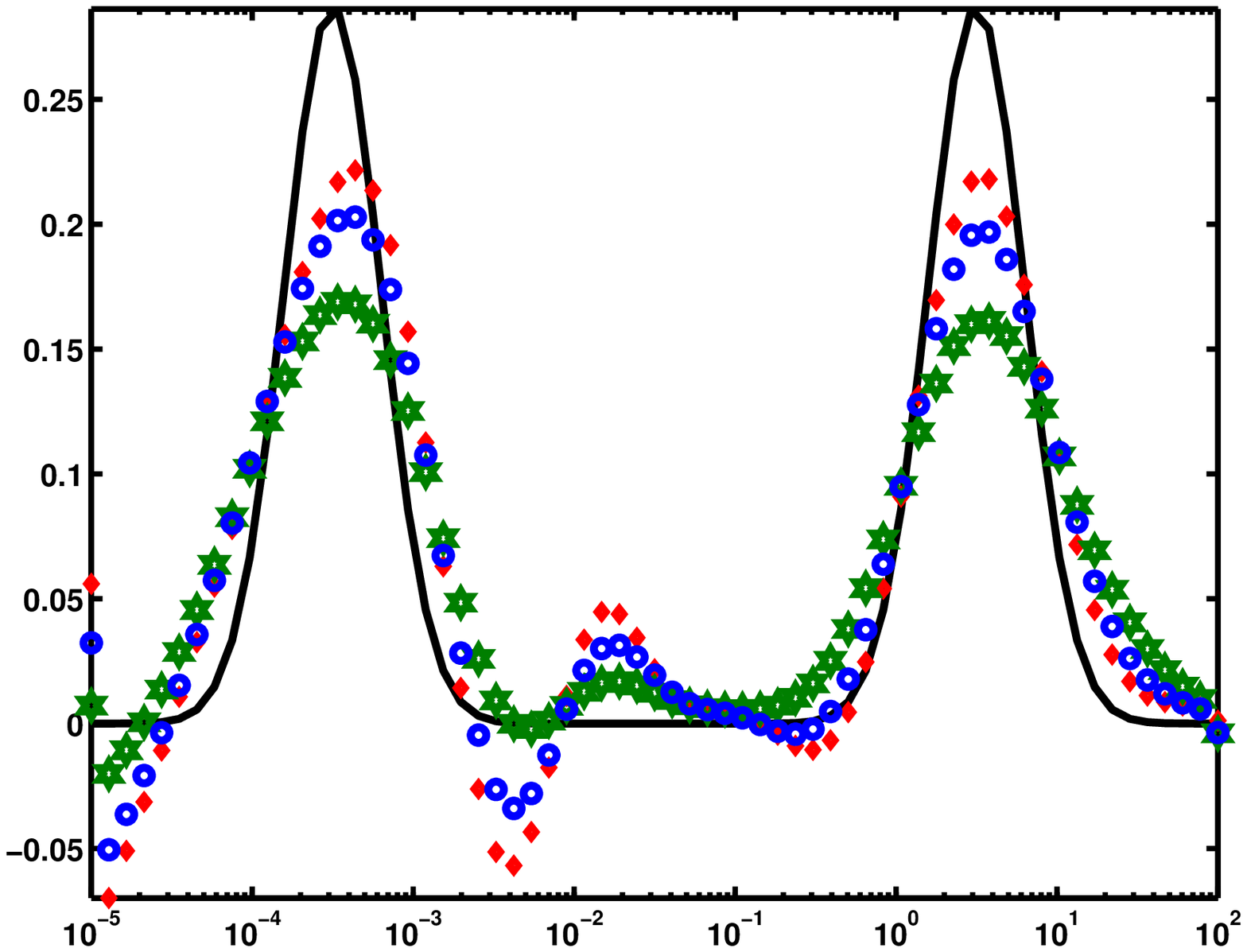}}
\subfigure[$L=L_1$]{\includegraphics[width=1.7in]{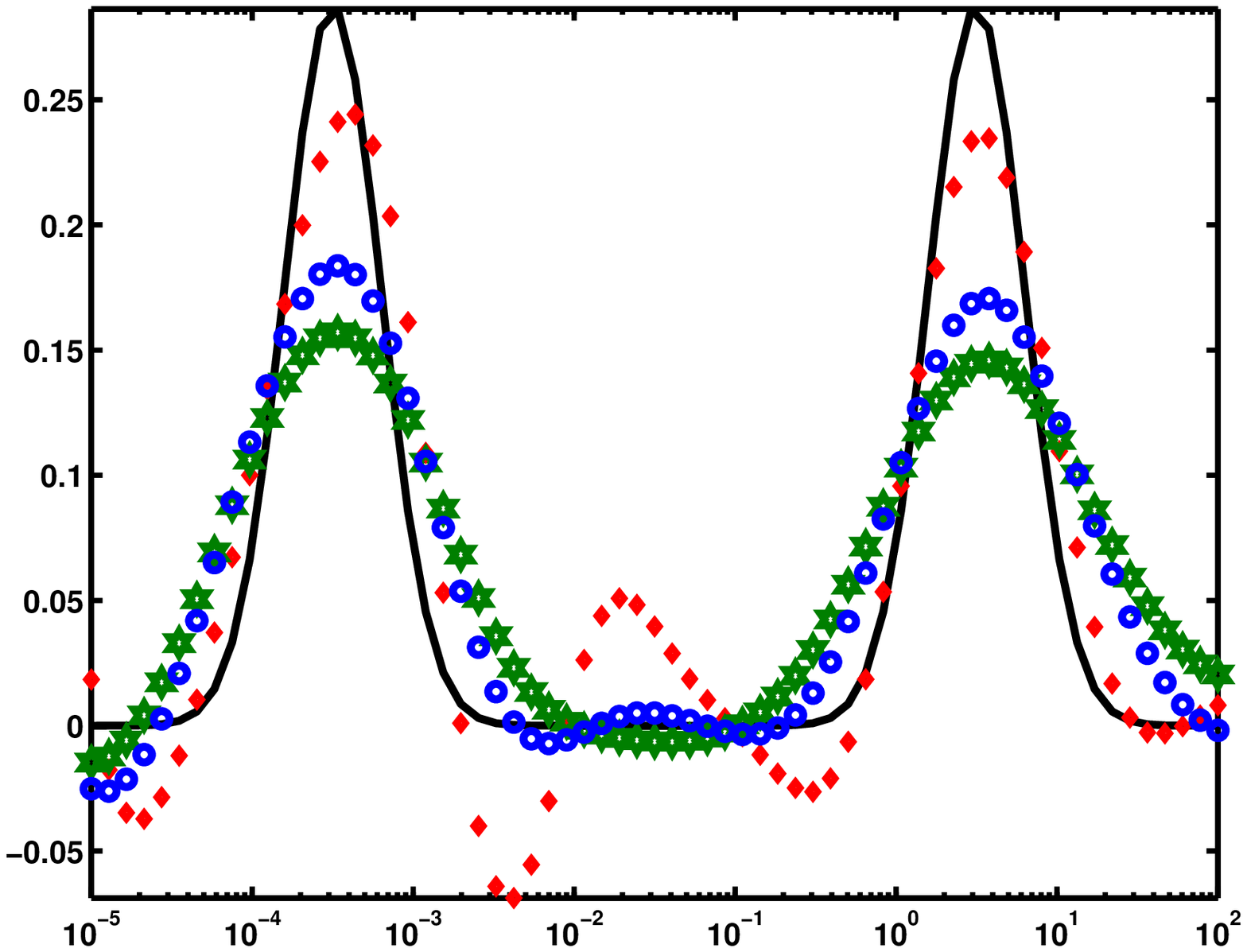}}
\subfigure[$L=L_2$]{\includegraphics[width=1.7in]{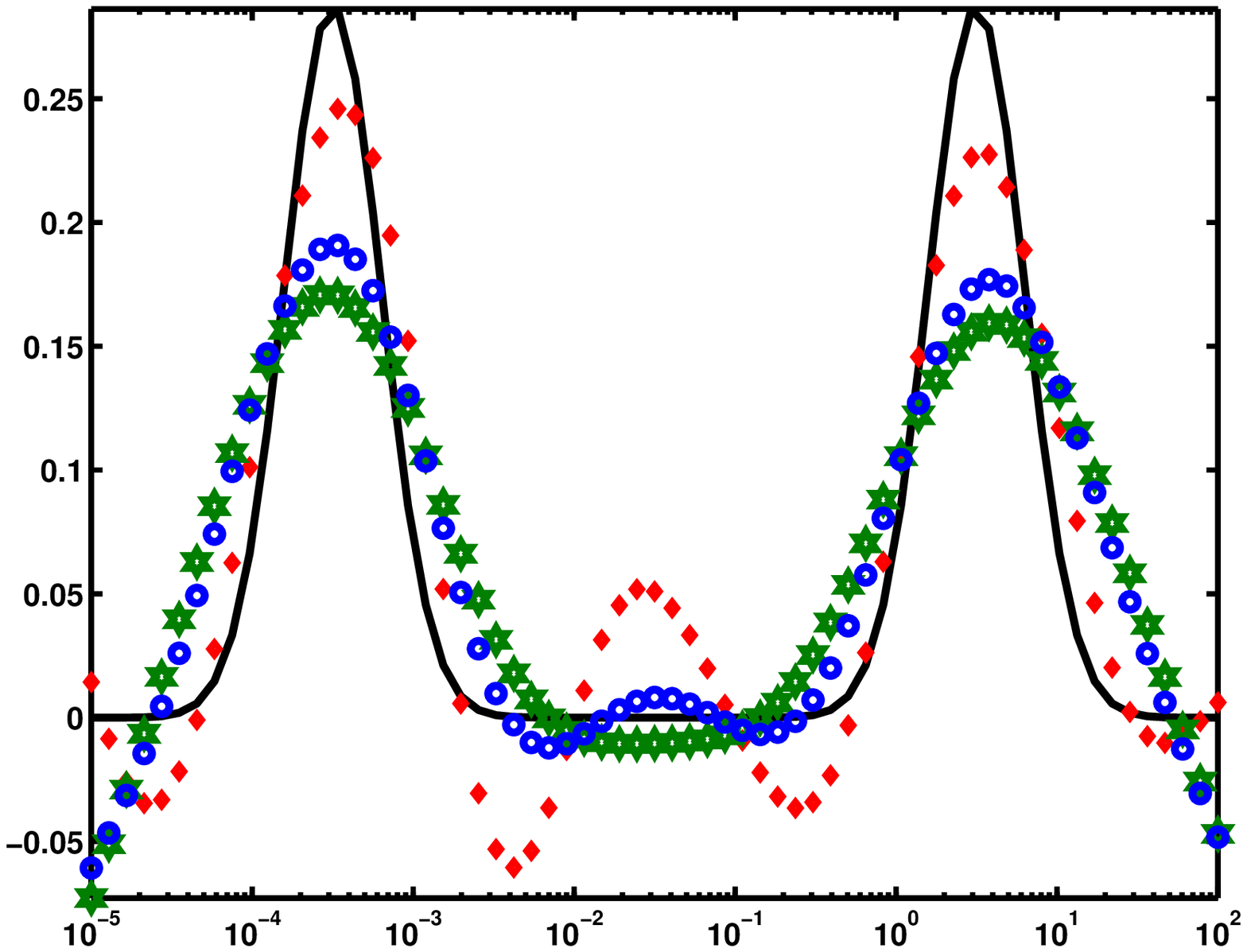}}
\caption{Mean error and example LS solutions.  $5\%$ noise. LN-A data set matrix $A_3$ }
\label{fig-lambdachoiceLN2A3HNLS}
\end{figure}

 \begin{figure}[!h]
  \centering
\subfigure[$L=I$]{\includegraphics[width=1.7in]{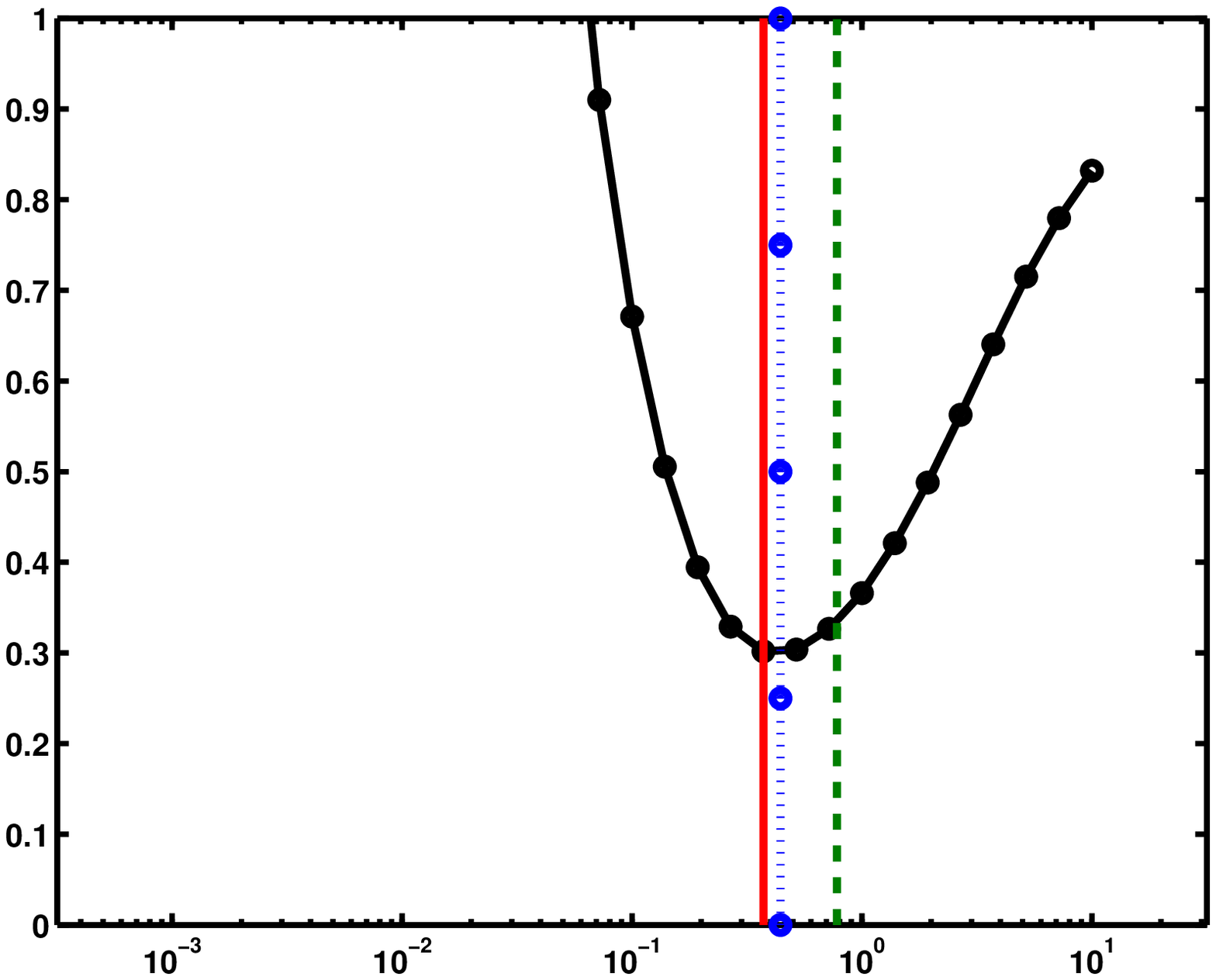}}
\subfigure[$L=L_1$]{\includegraphics[width=1.7in]{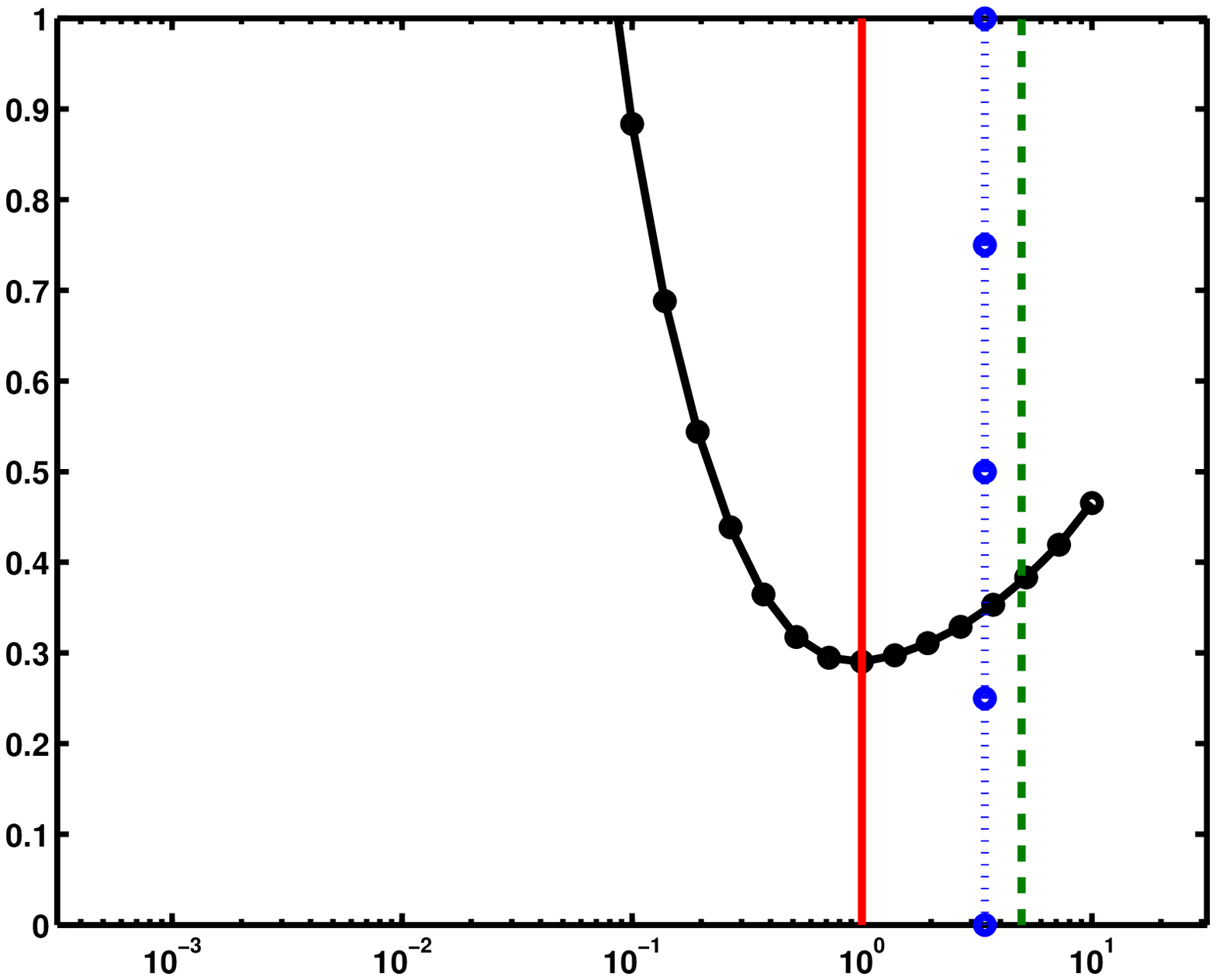}}
\subfigure[$L=L_2$]{\includegraphics[width=1.7in]{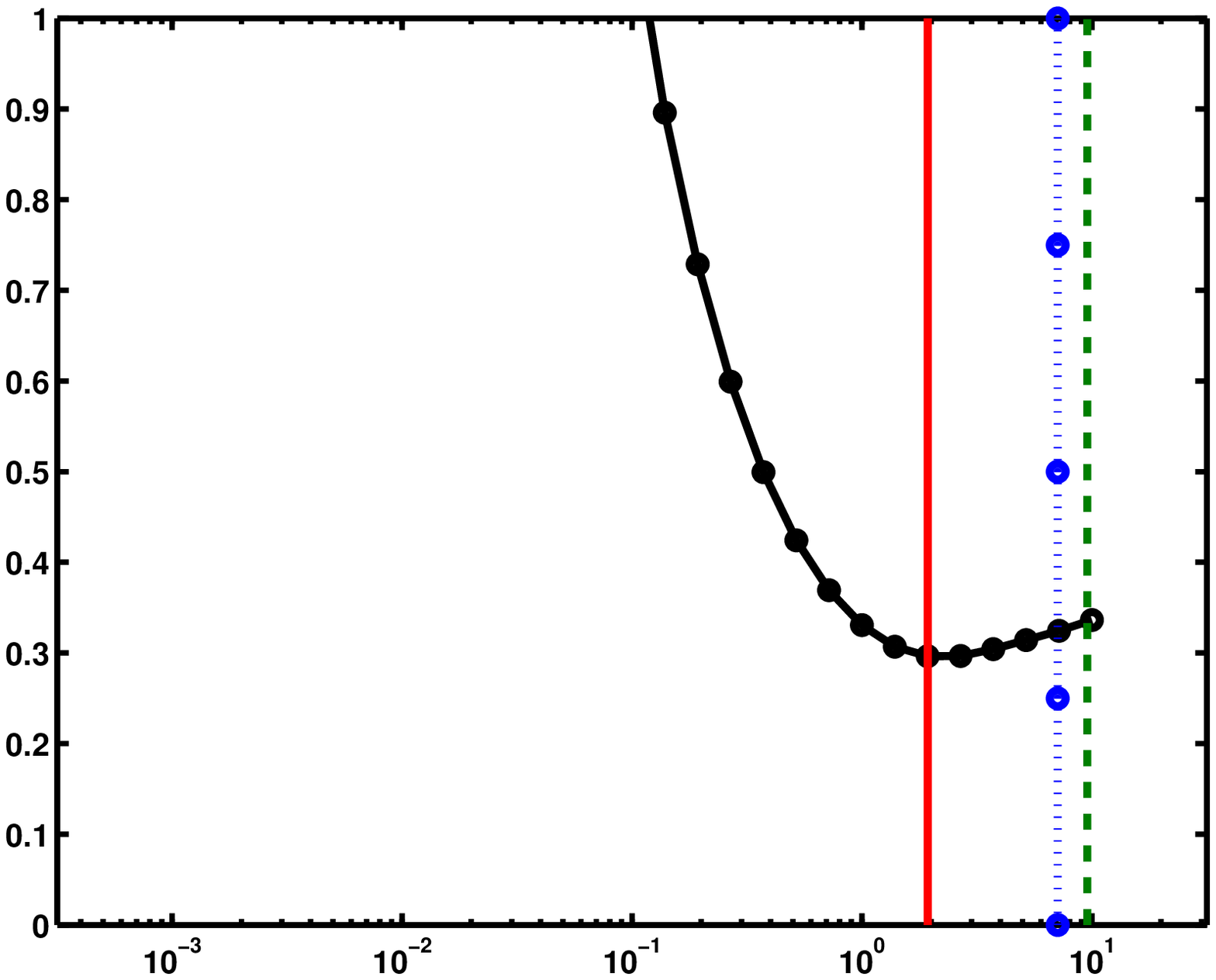}}
\subfigure[$L=I$]{\includegraphics[width=1.7in]{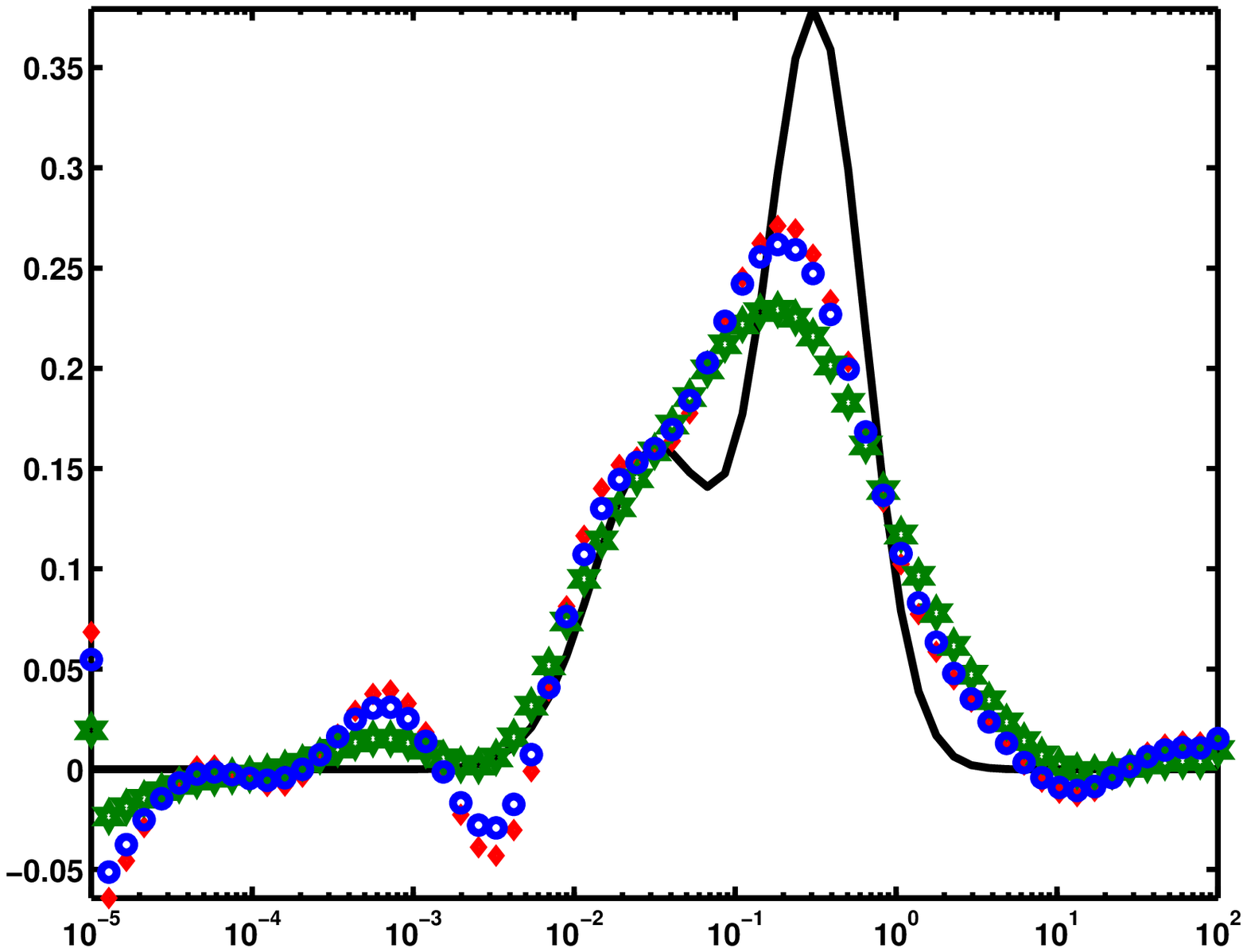}}
\subfigure[$L=L_1$]{\includegraphics[width=1.7in]{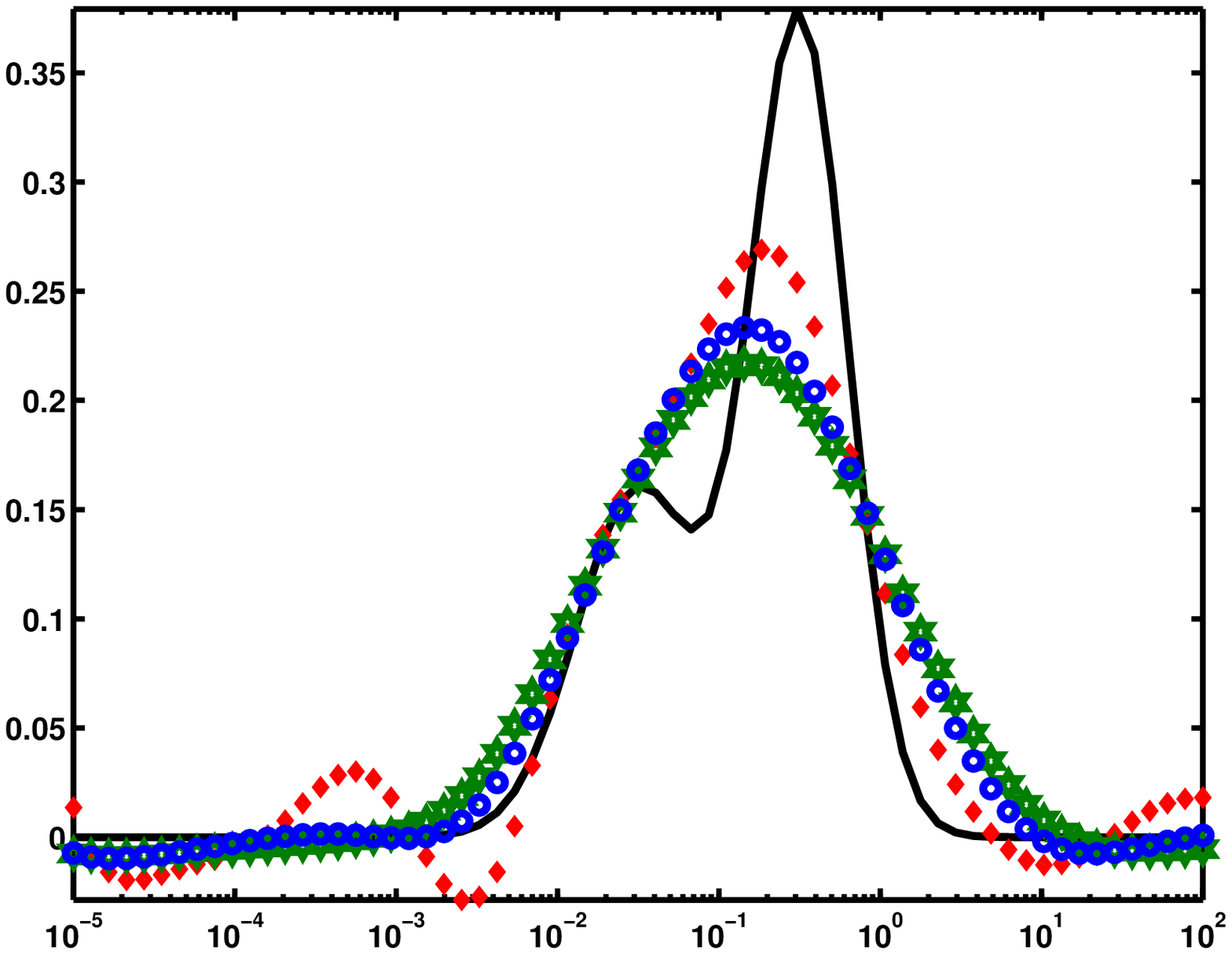}}
\subfigure[$L=L_2$]{\includegraphics[width=1.7in]{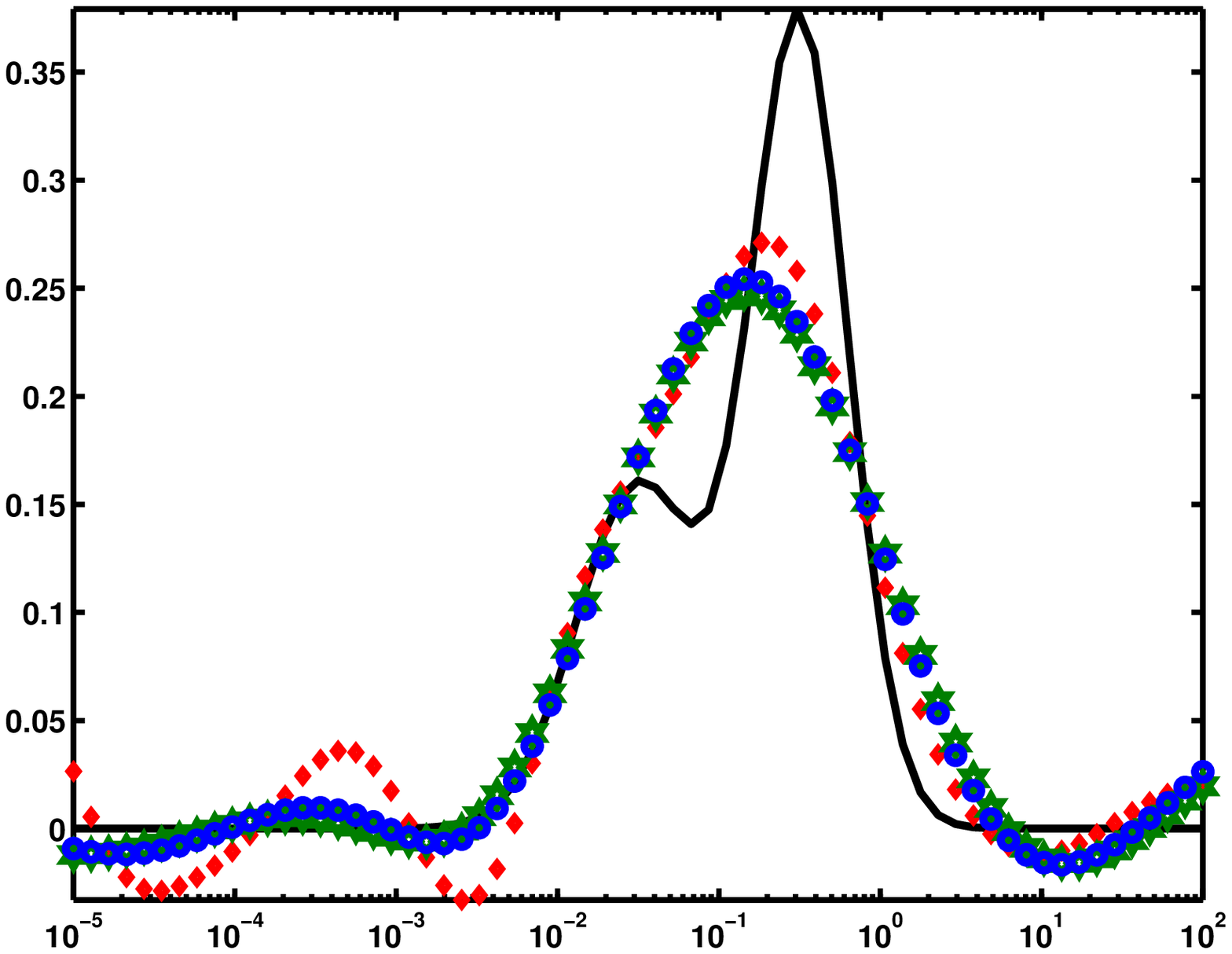}}
\caption{Mean error and example LS solutions.  $5\%$ noise. LN-B data set matrix $A_3$ }
\label{fig-lambdachoiceLN5A3HNLS}
\end{figure}

 \begin{figure}[!h]
  \centering
\subfigure[$L=I$]{\includegraphics[width=1.7in]{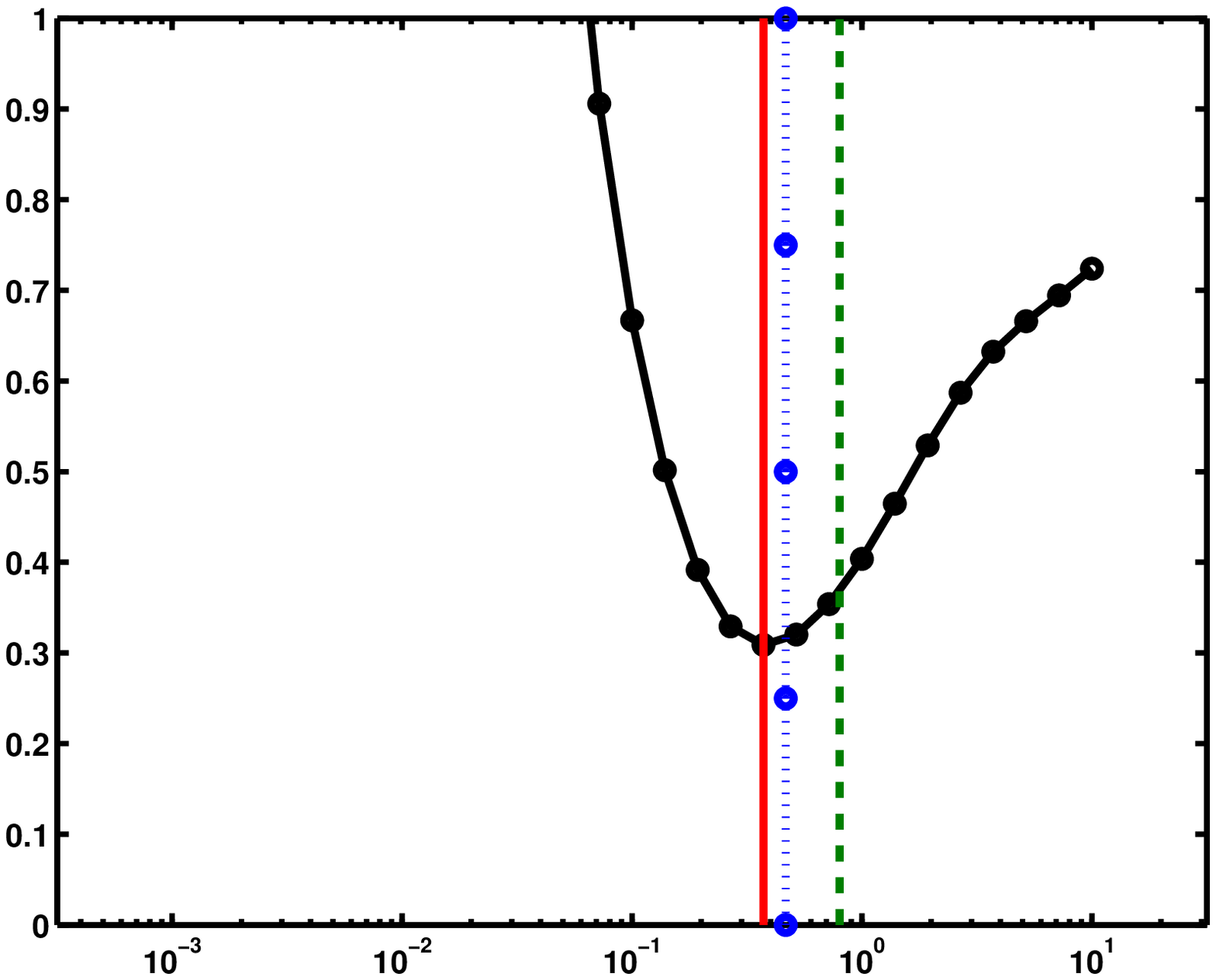}}
\subfigure[$L=L_1$]{\includegraphics[width=1.7in]{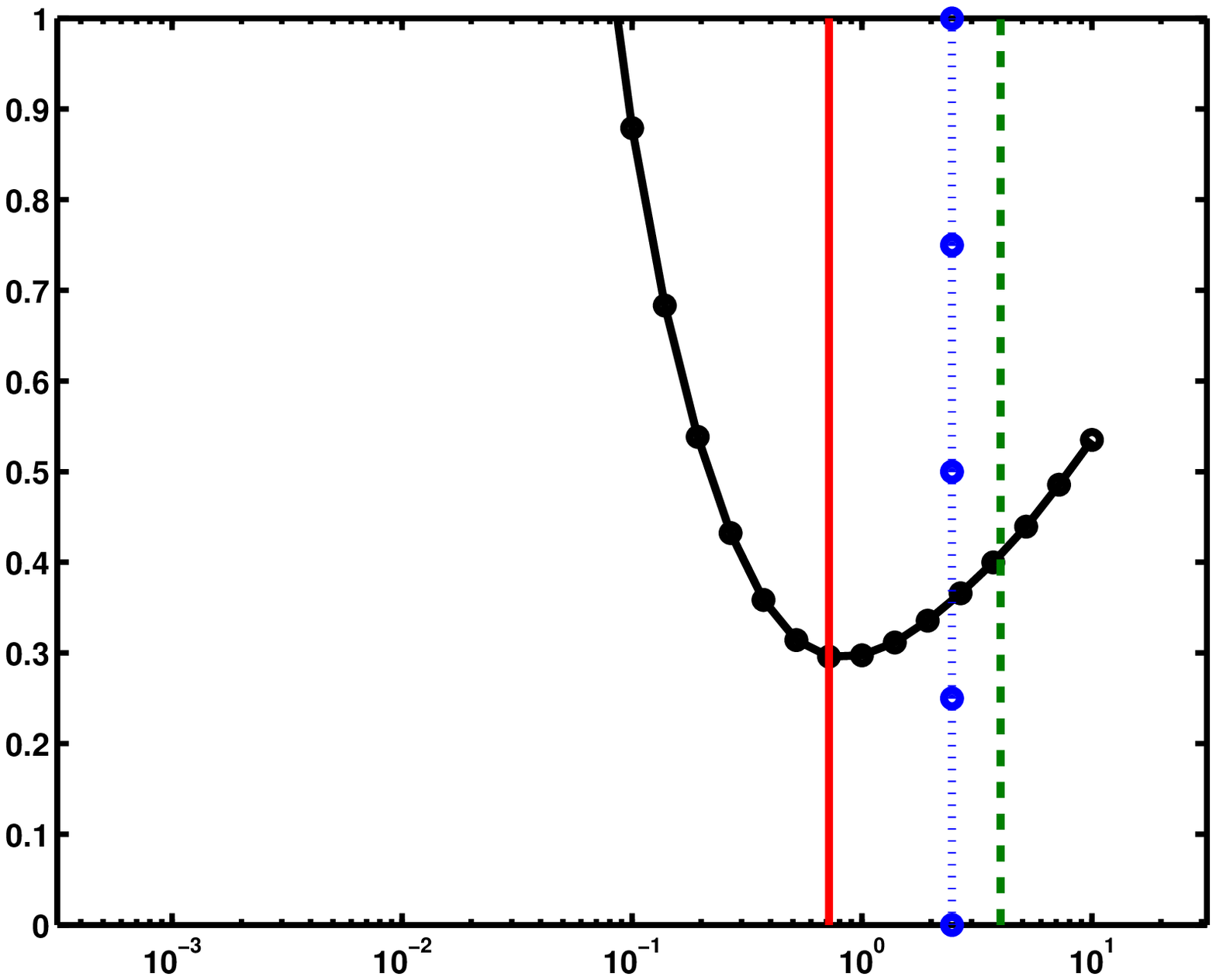}}
\subfigure[$L=L_2$]{\includegraphics[width=1.7in]{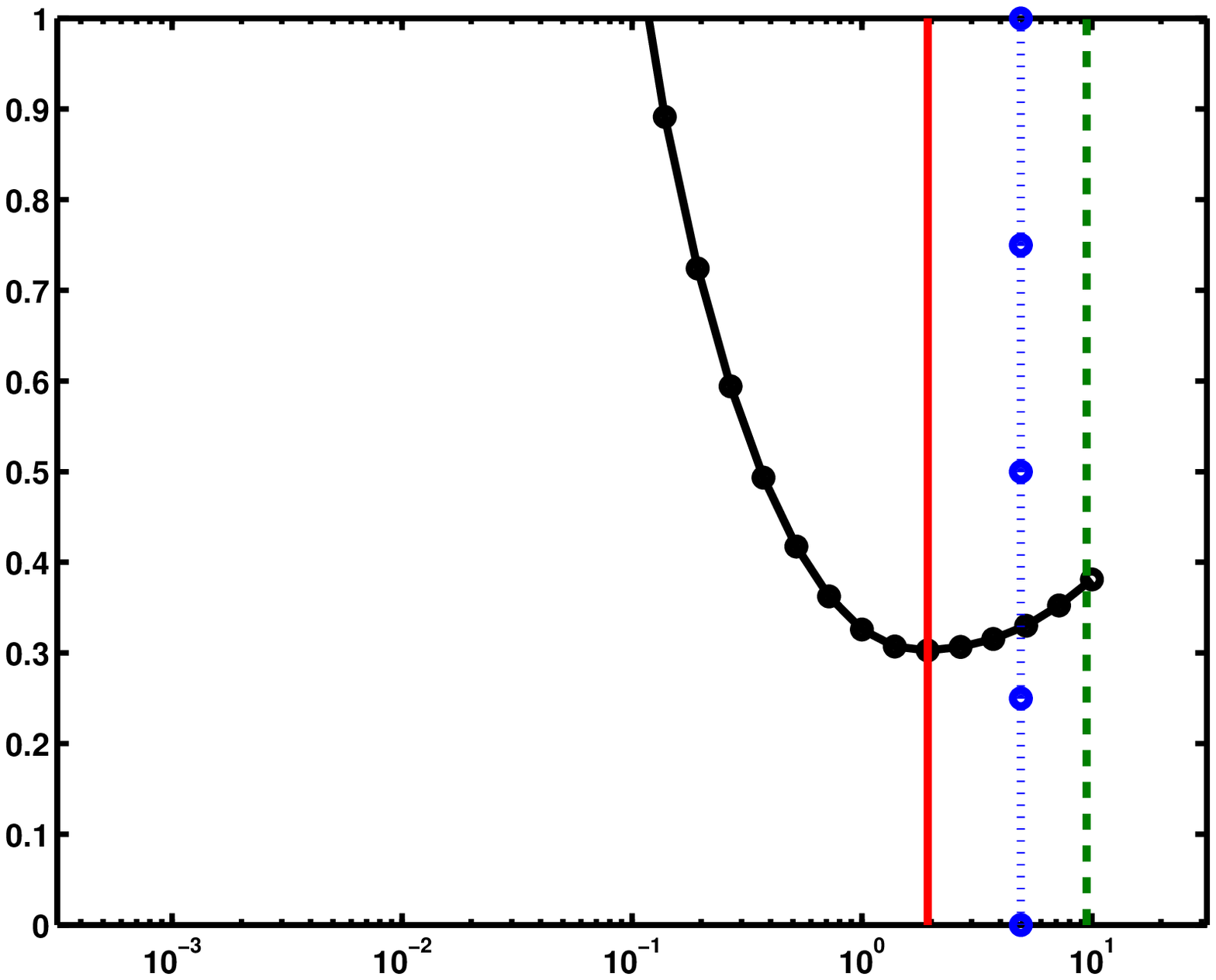}}
\subfigure[$L=I$]{\includegraphics[width=1.7in]{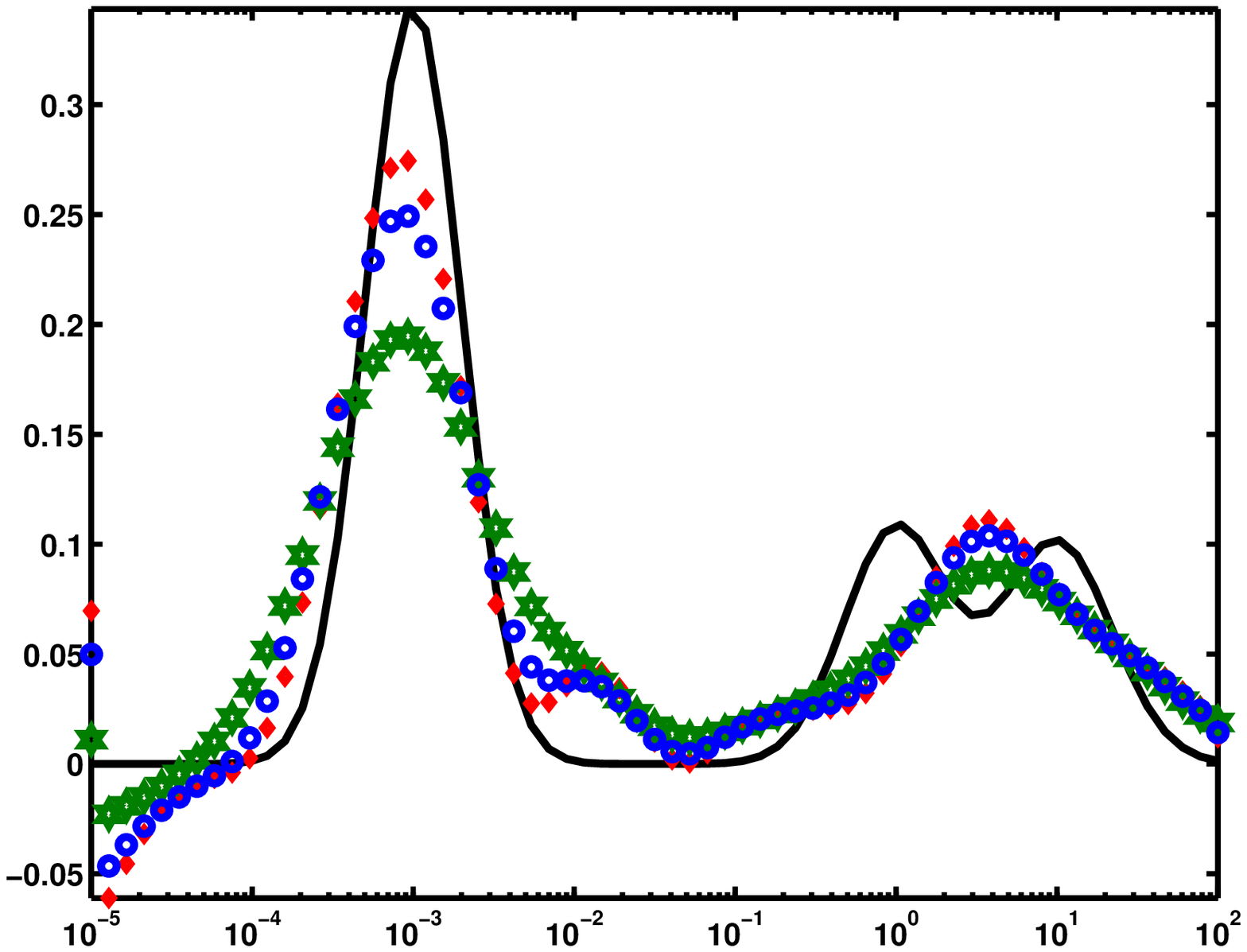}}
\subfigure[$L=L_1$]{\includegraphics[width=1.7in]{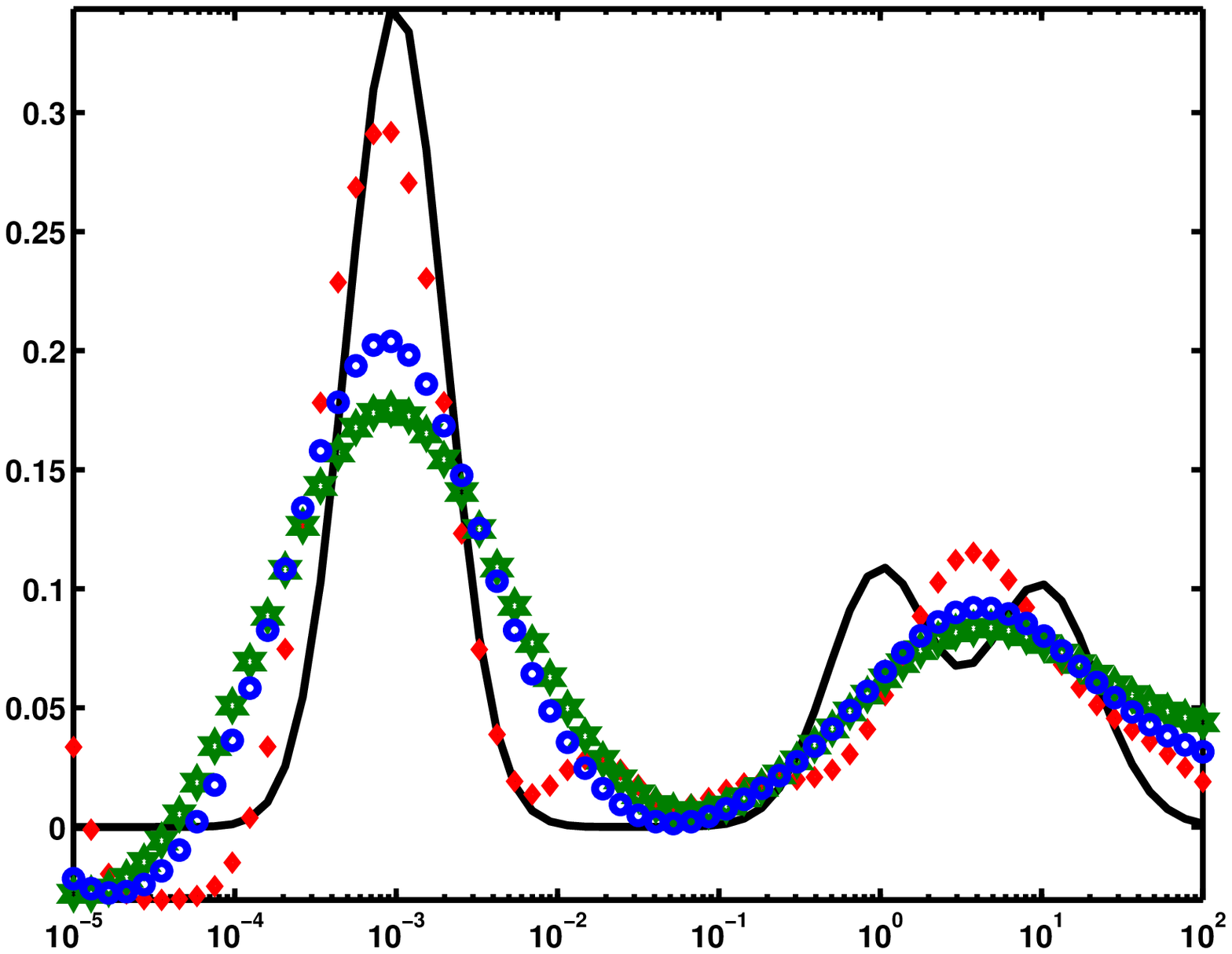}}
\subfigure[$L=L_2$]{\includegraphics[width=1.7in]{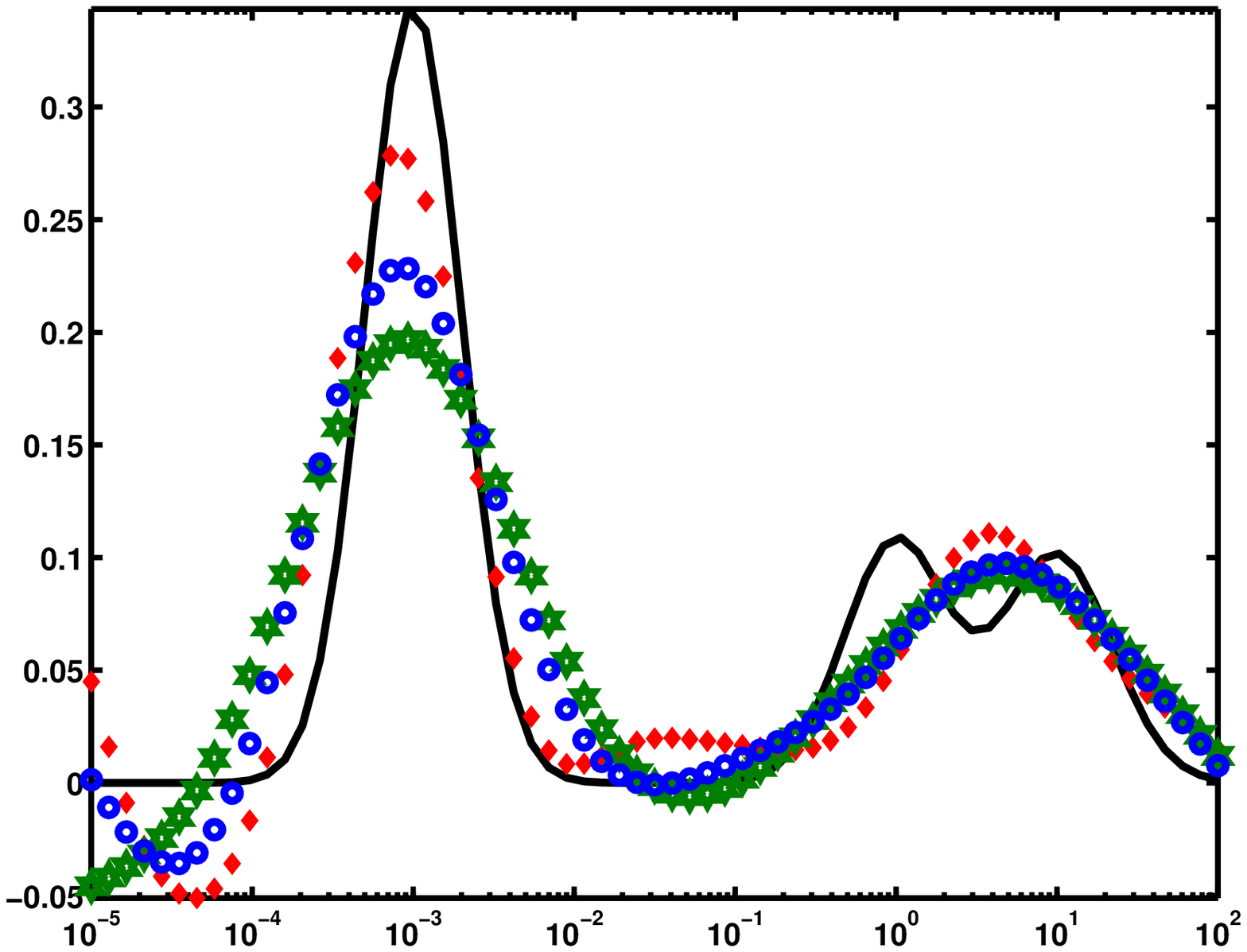}}
\caption{Mean error and example LS solutions.  $5\%$ noise. LN-C data set matrix $A_3$}
\label{fig-lambdachoiceLN6A3HNLS}
\end{figure}
\clearpage

\subsection{Results using LS $A_4$ Noise level $5\%$}
 \begin{figure}[!h]
 \centering
\subfigure[$L=I$]{\includegraphics[width=1.7in]{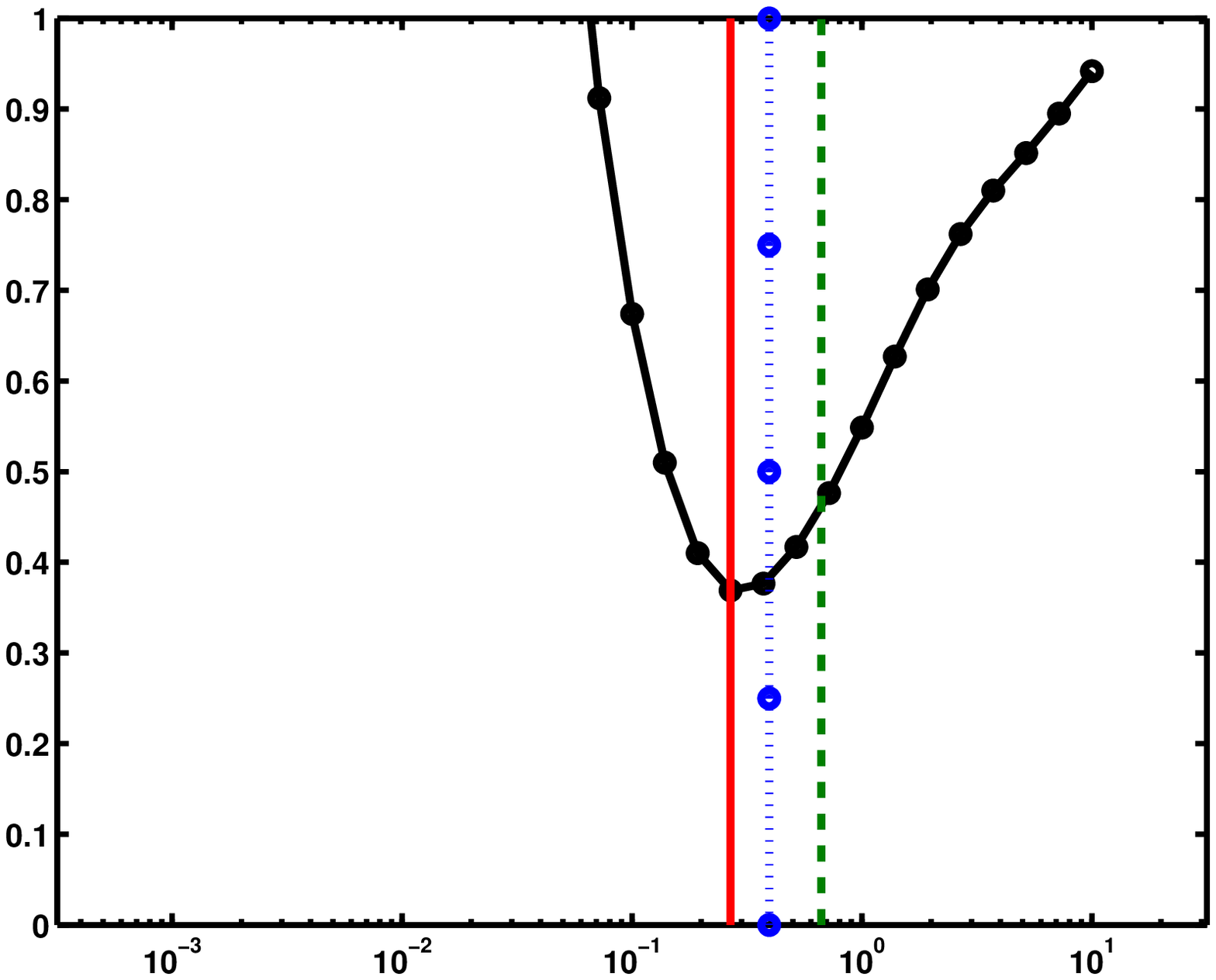}}
\subfigure[$L=L_1$]{\includegraphics[width=1.7in]{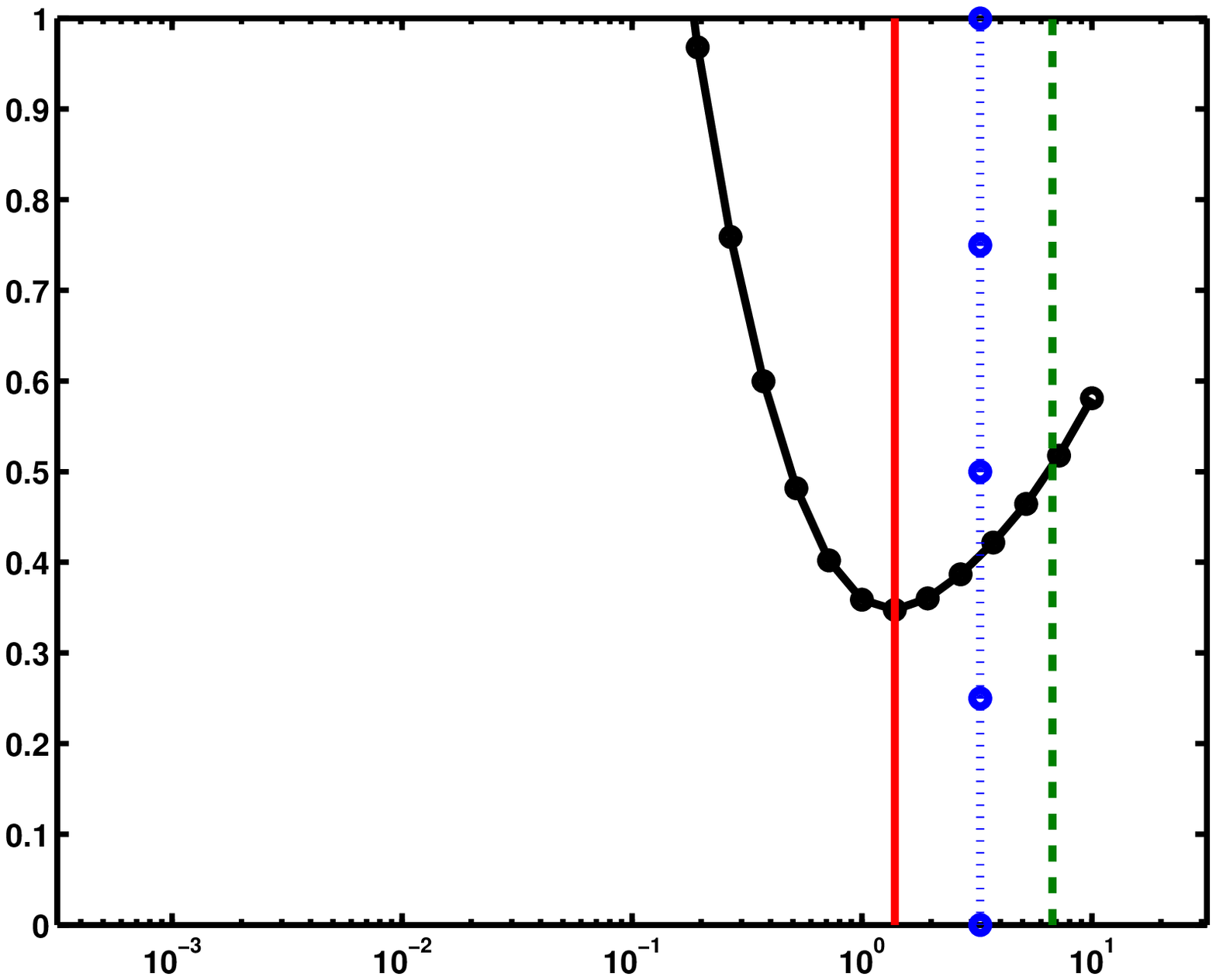}}
\subfigure[$L=L_2$]{\includegraphics[width=1.7in]{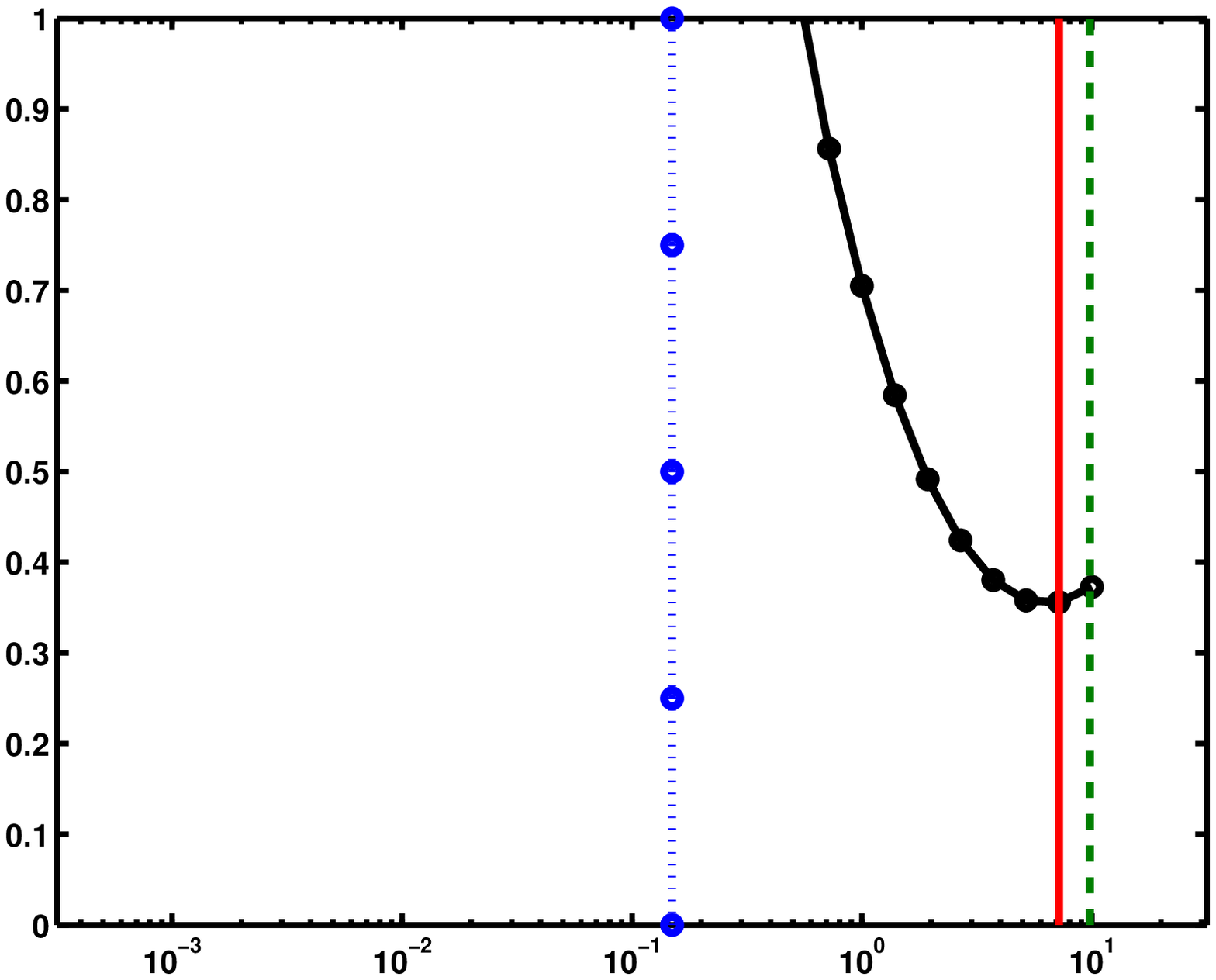}}
\subfigure[$L=I$]{\includegraphics[width=1.7in]{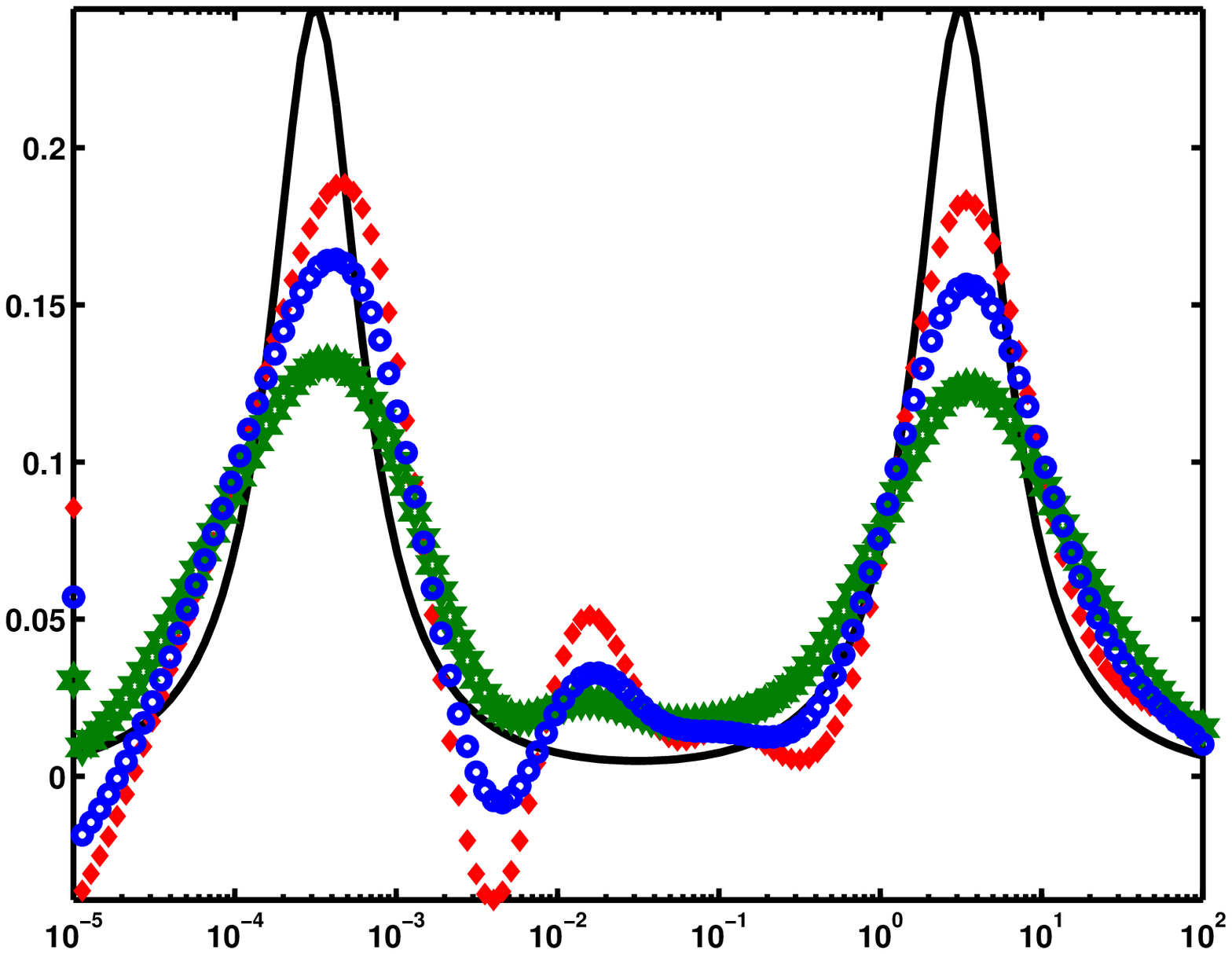}}
\subfigure[$L=L_1$]{\includegraphics[width=1.7in]{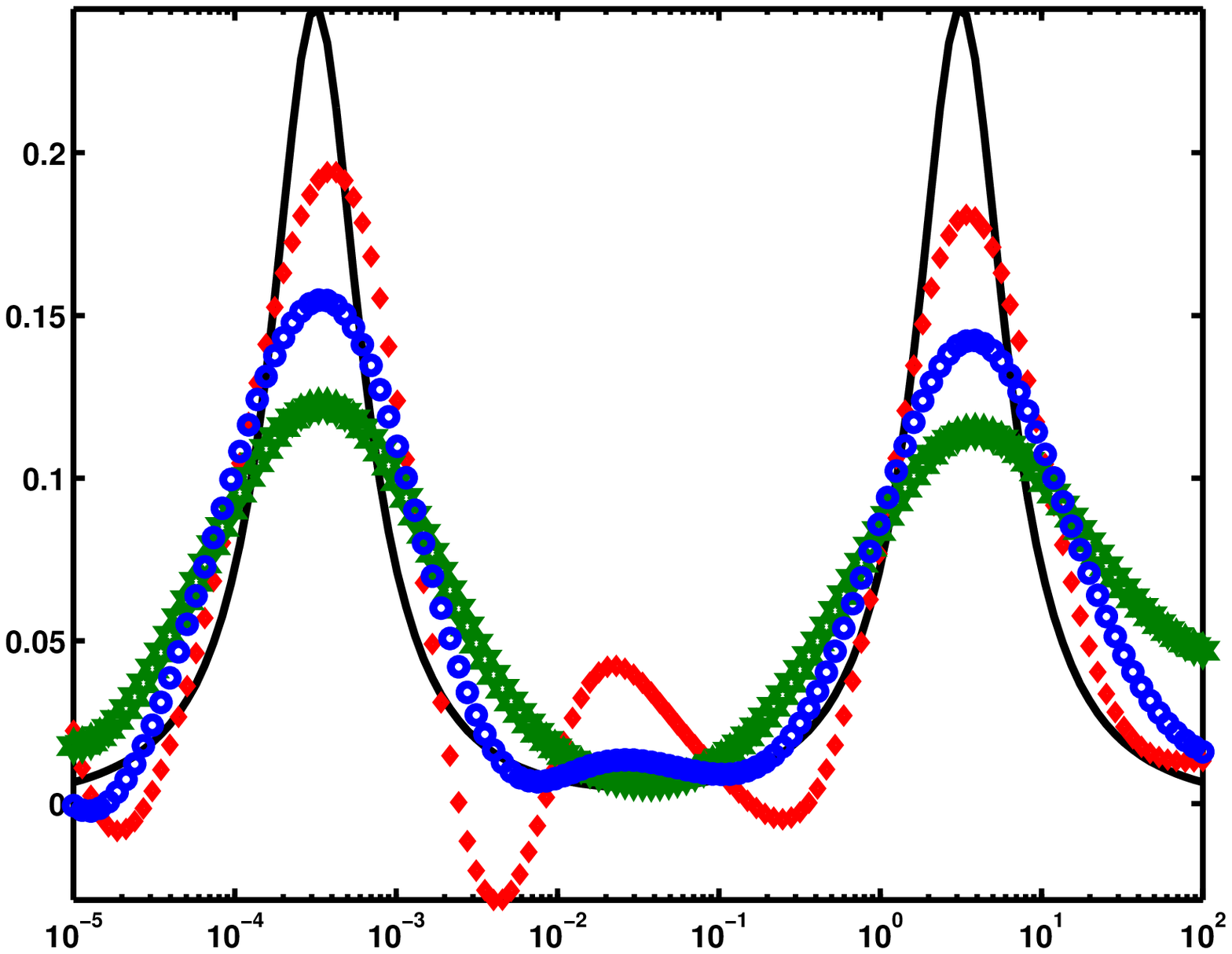}}
\subfigure[$L=L_2$]{\includegraphics[width=1.7in]{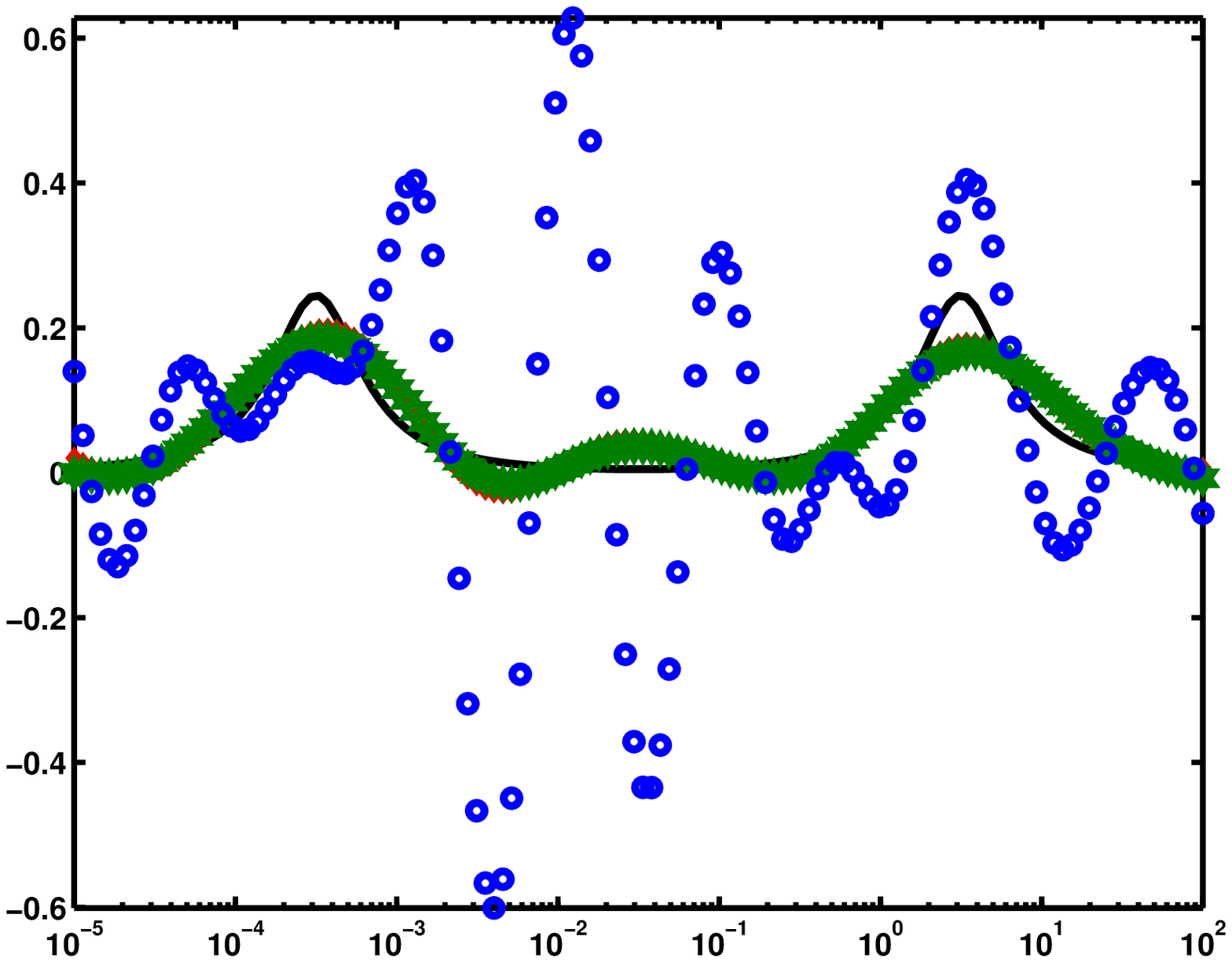}}
\caption{Mean error and example LS solutions.  $5\%$ noise. RQ-A data set matrix $A_4$}
\label{fig-lambdachoiceRQ1A4HNLS}
\end{figure}

 \begin{figure}[!h]
  \centering
\subfigure[$L=I$]{\includegraphics[width=1.7in]{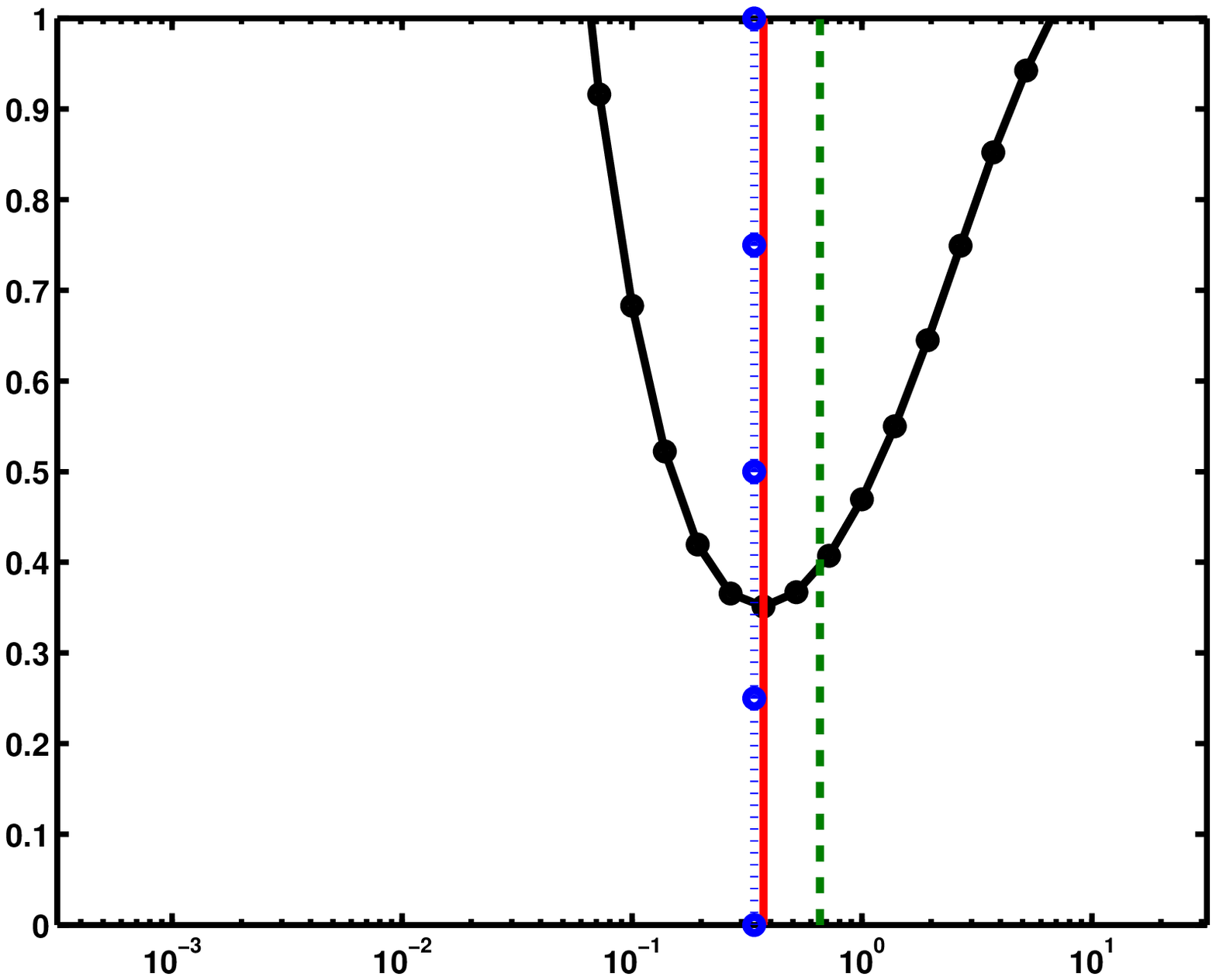}}
\subfigure[$L=L_1$]{\includegraphics[width=1.7in]{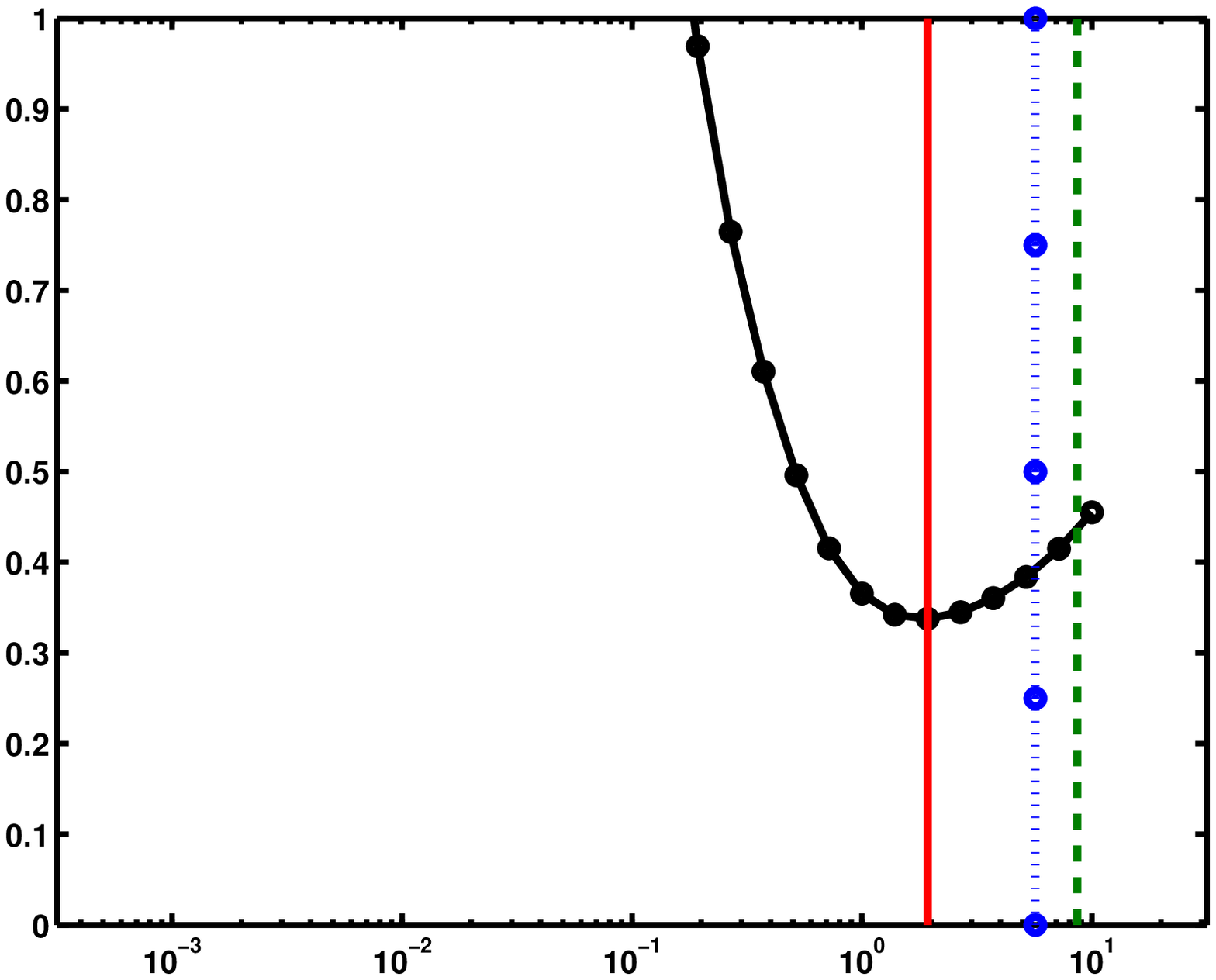}}
\subfigure[$L=L_2$]{\includegraphics[width=1.7in]{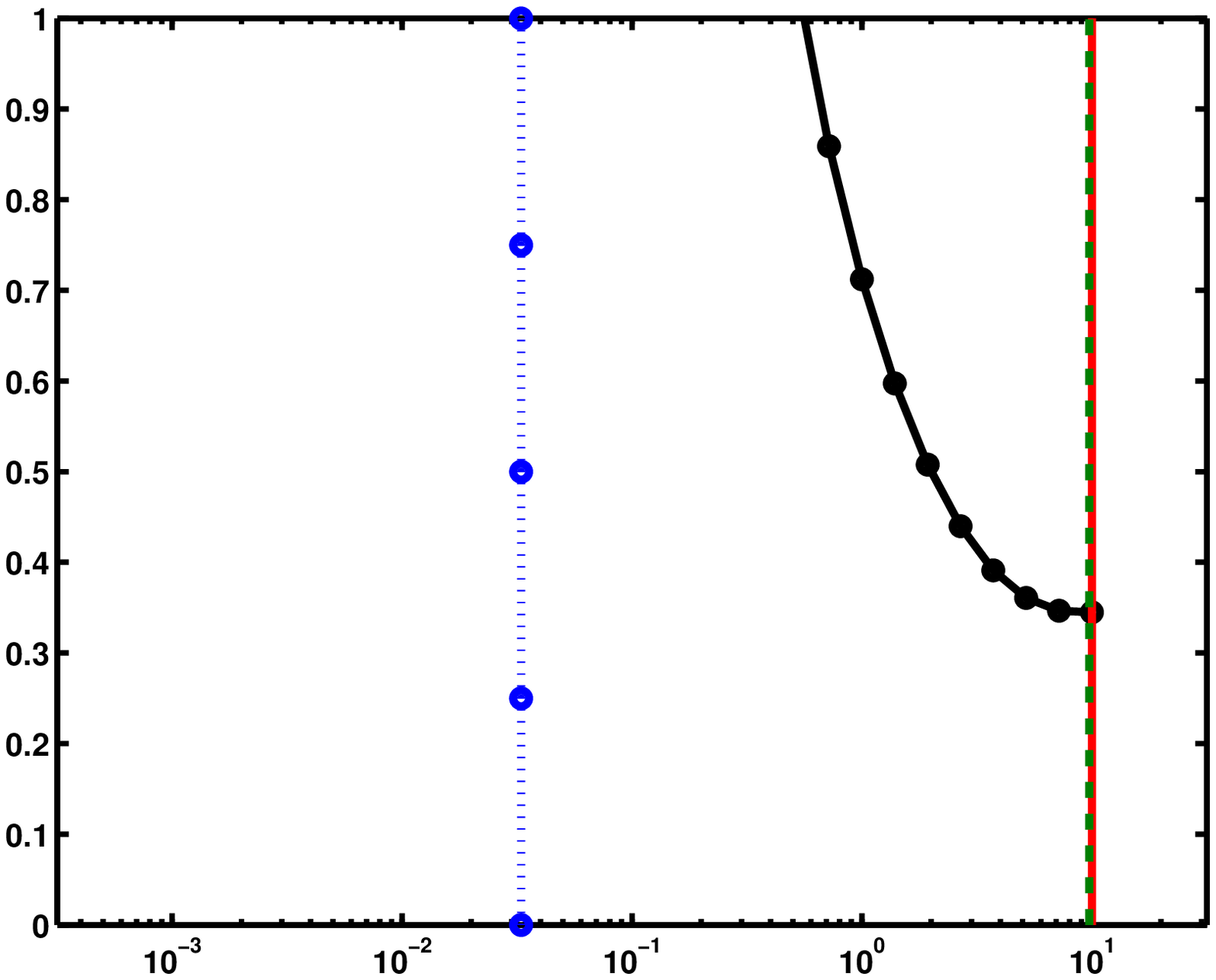}}
\subfigure[$L=I$]{\includegraphics[width=1.7in]{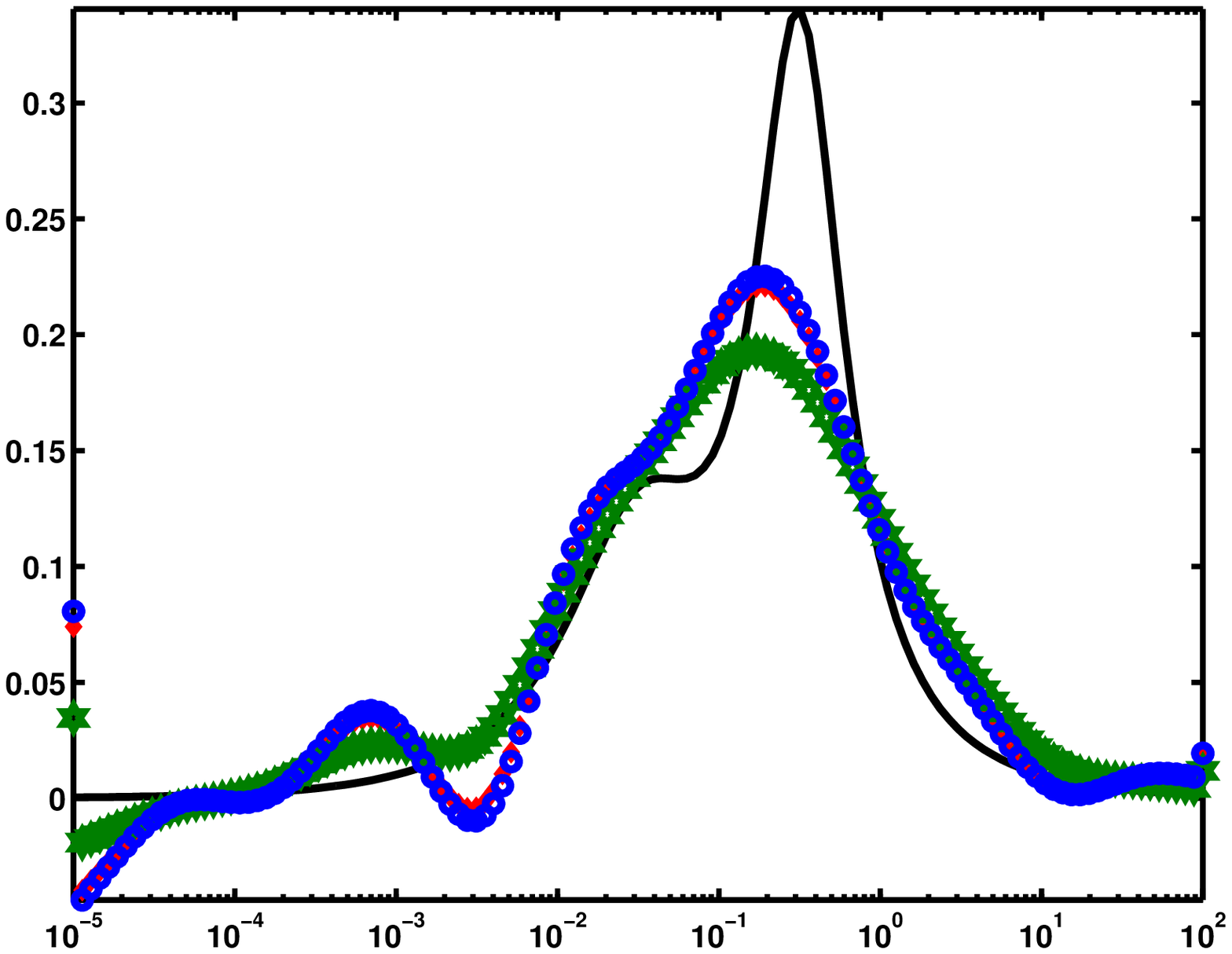}}
\subfigure[$L=L_1$]{\includegraphics[width=1.7in]{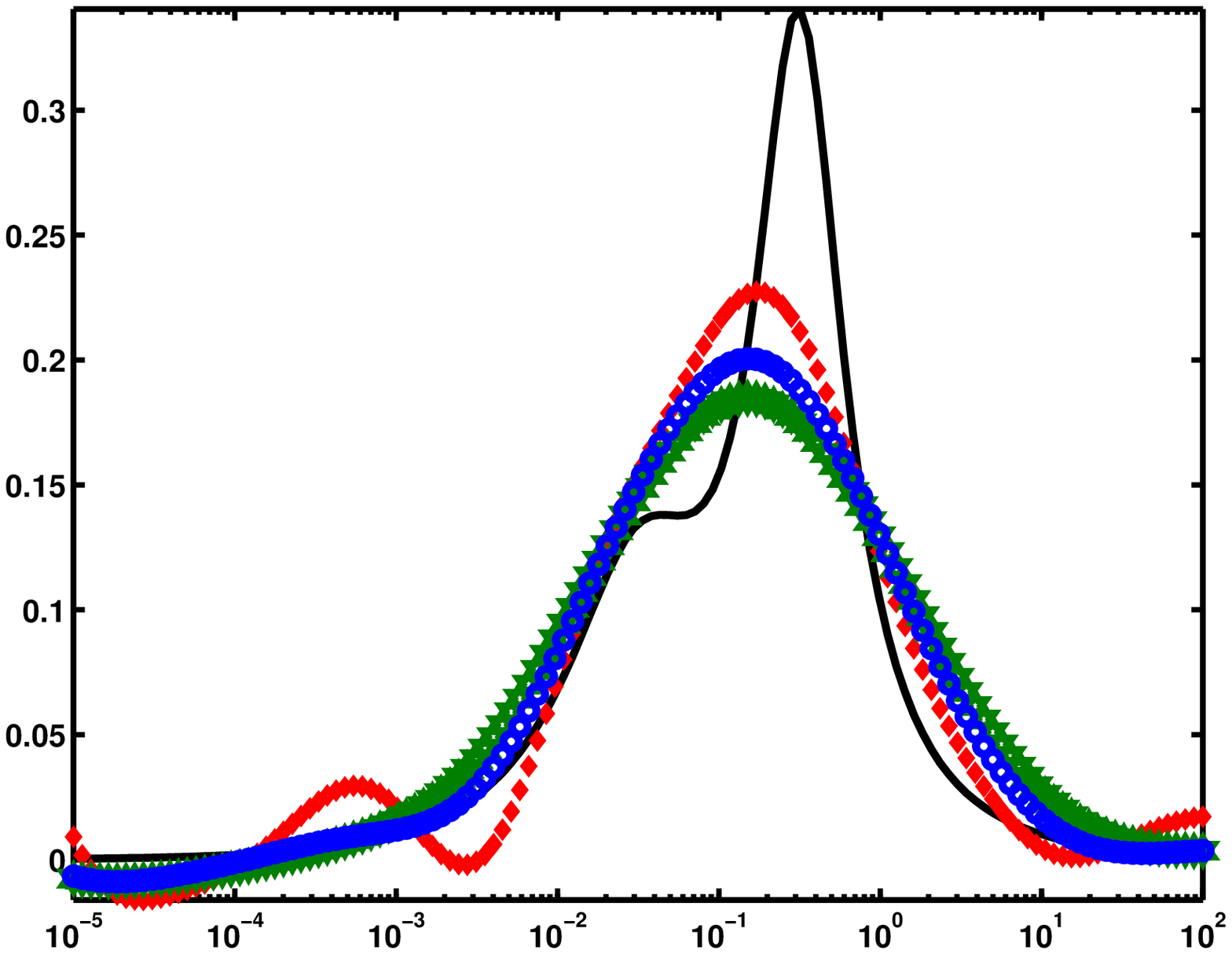}}
\subfigure[$L=L_2$]{\includegraphics[width=1.7in]{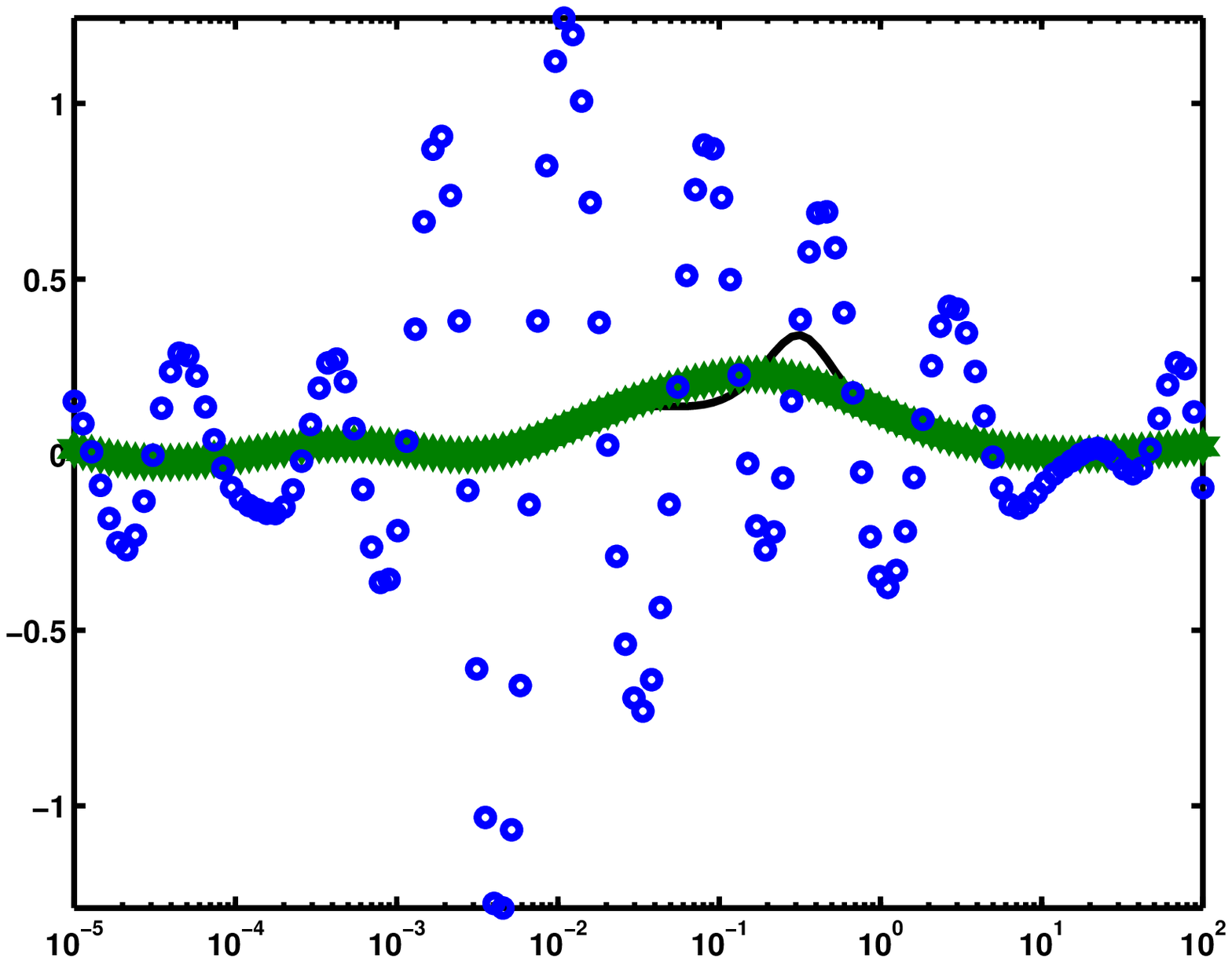}}
\caption{Mean error and example LS solutions. $5\%$ noise. RQ-B data set matrix $A_4$}
\label{fig-lambdachoiceRQ5A4HNLS}
\end{figure}

 \begin{figure}[!h]
  \centering
\subfigure[$L=I$]{\includegraphics[width=1.7in]{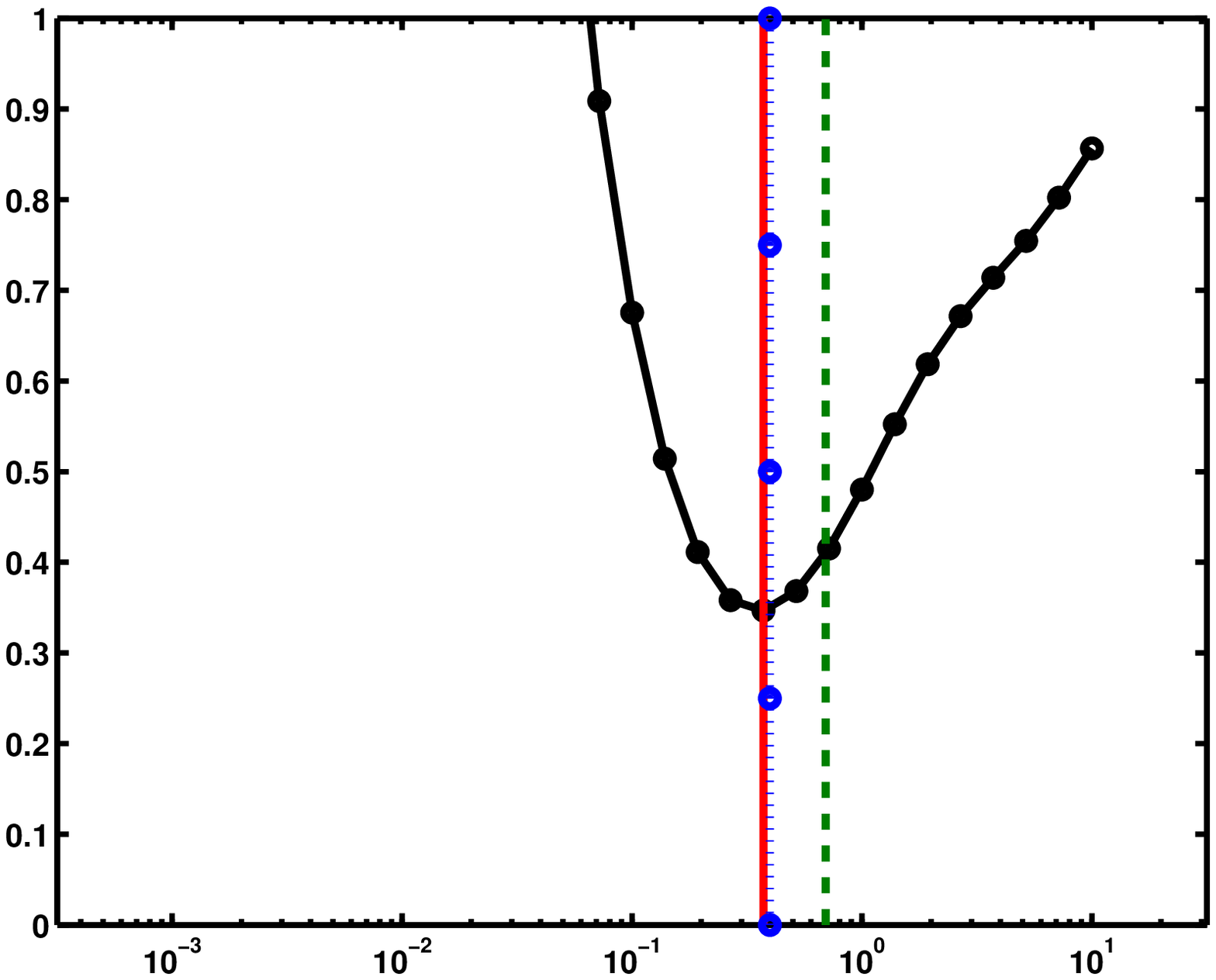}}
\subfigure[$L=L_1$]{\includegraphics[width=1.7in]{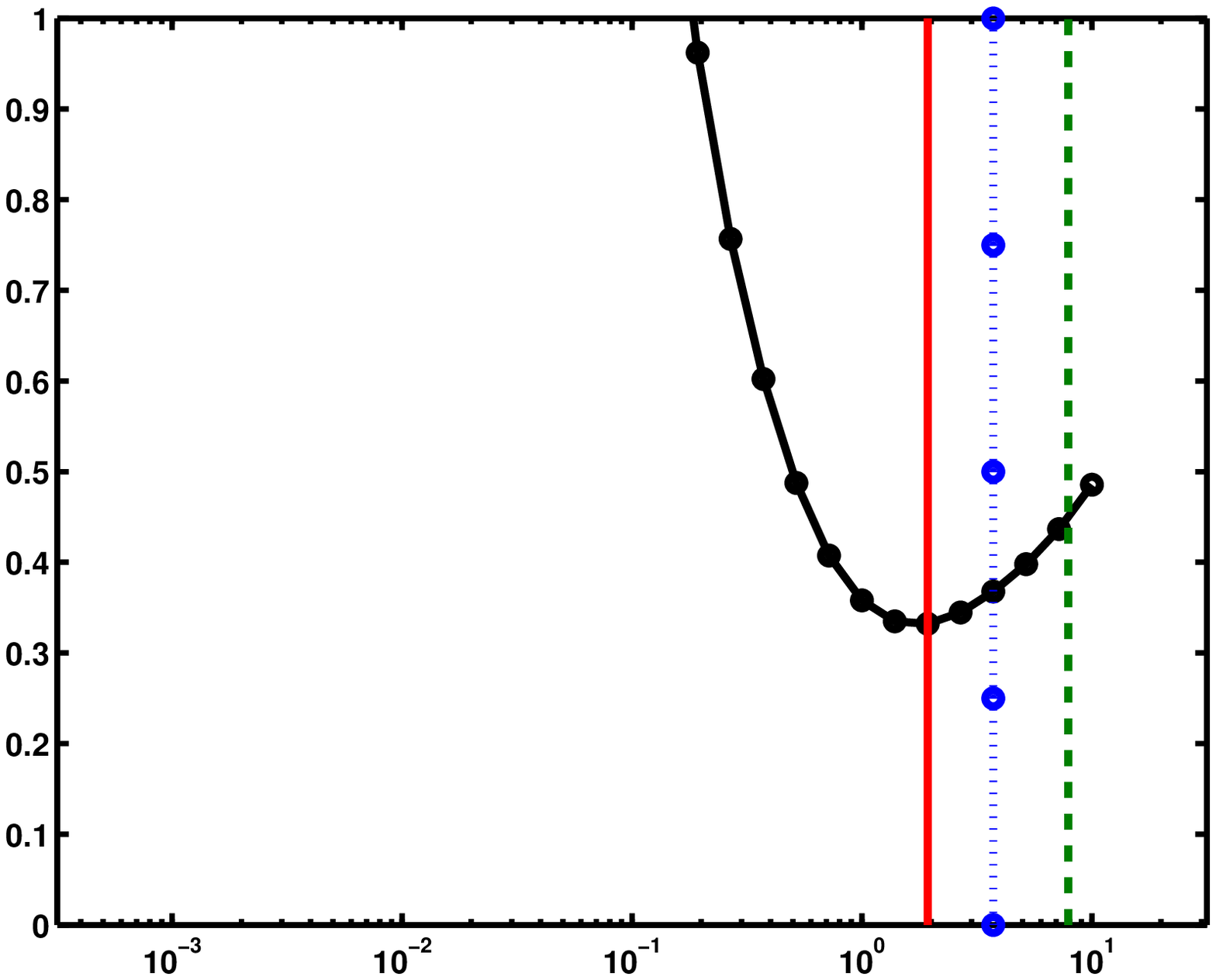}}
\subfigure[$L=L_2$]{\includegraphics[width=1.7in]{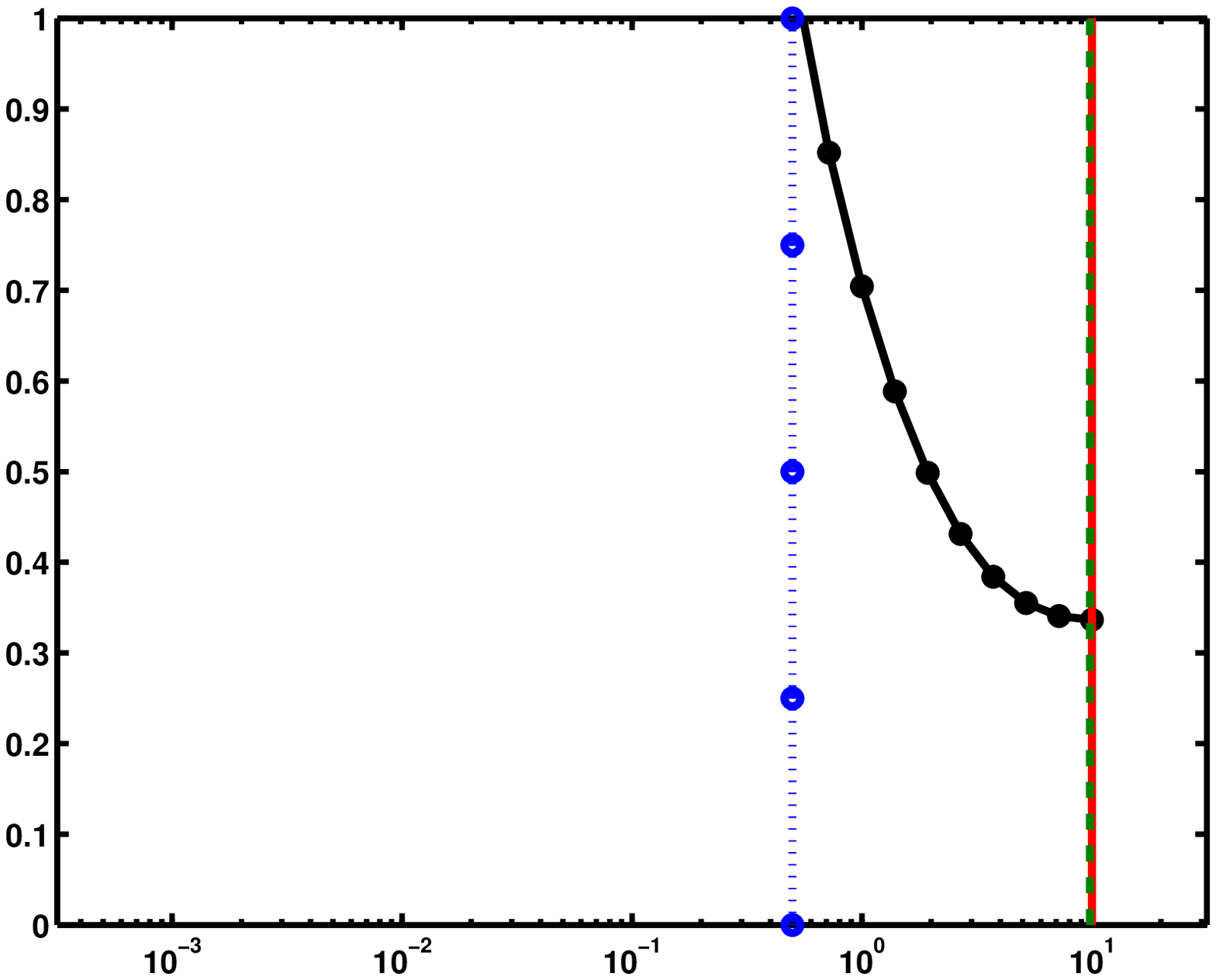}}
\subfigure[$L=I$]{\includegraphics[width=1.7in]{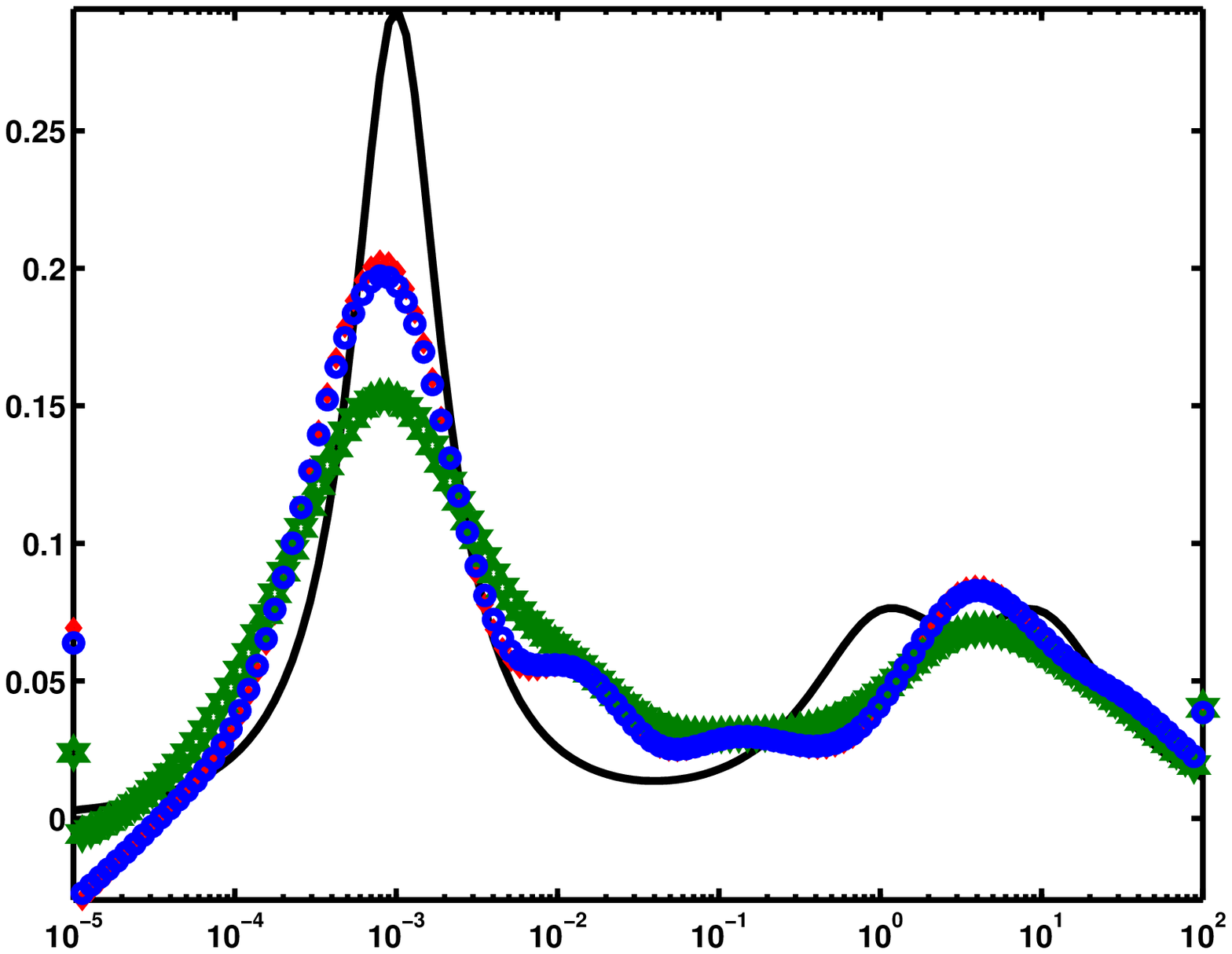}}
\subfigure[$L=L_1$]{\includegraphics[width=1.7in]{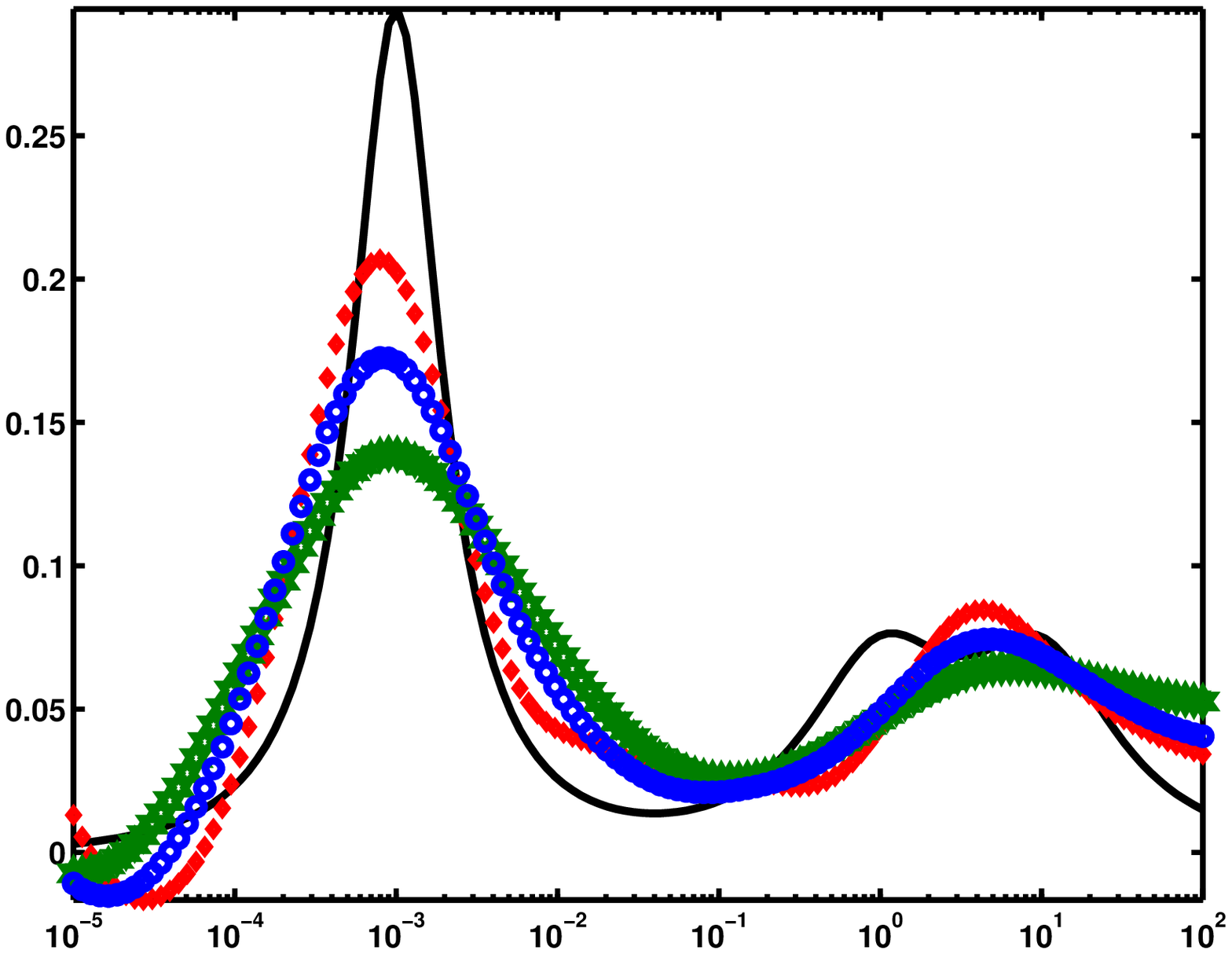}}
\subfigure[$L=L_2$]{\includegraphics[width=1.7in]{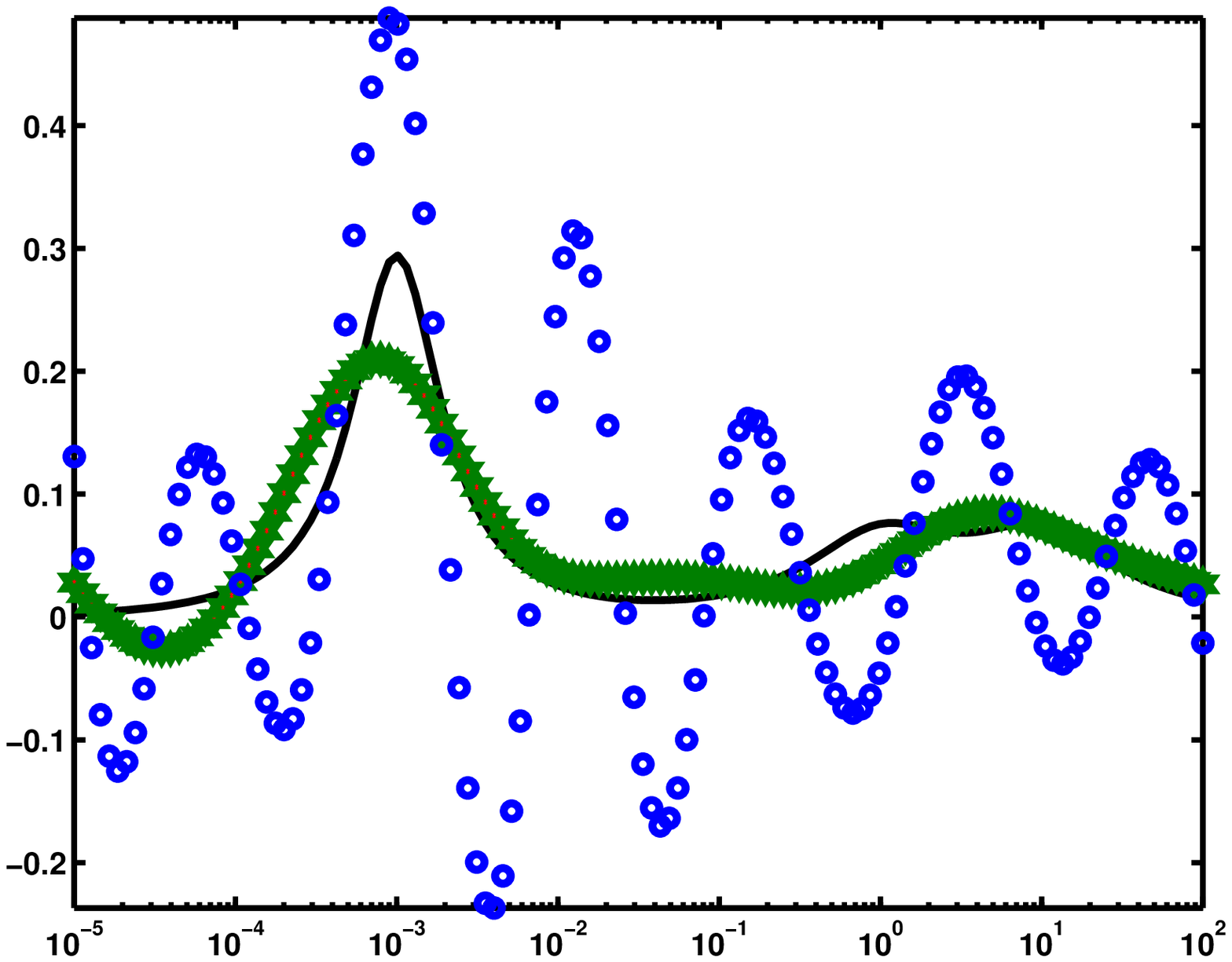}}
\caption{Mean error and example LS solutions.  $5\%$ noise. RQ-C data set matrix $A_4$}
\label{fig-lambdachoiceRQ6A4HNLS}
\end{figure}

 \begin{figure}[!h]
  \centering
\subfigure[$L=I$]{\includegraphics[width=1.7in]{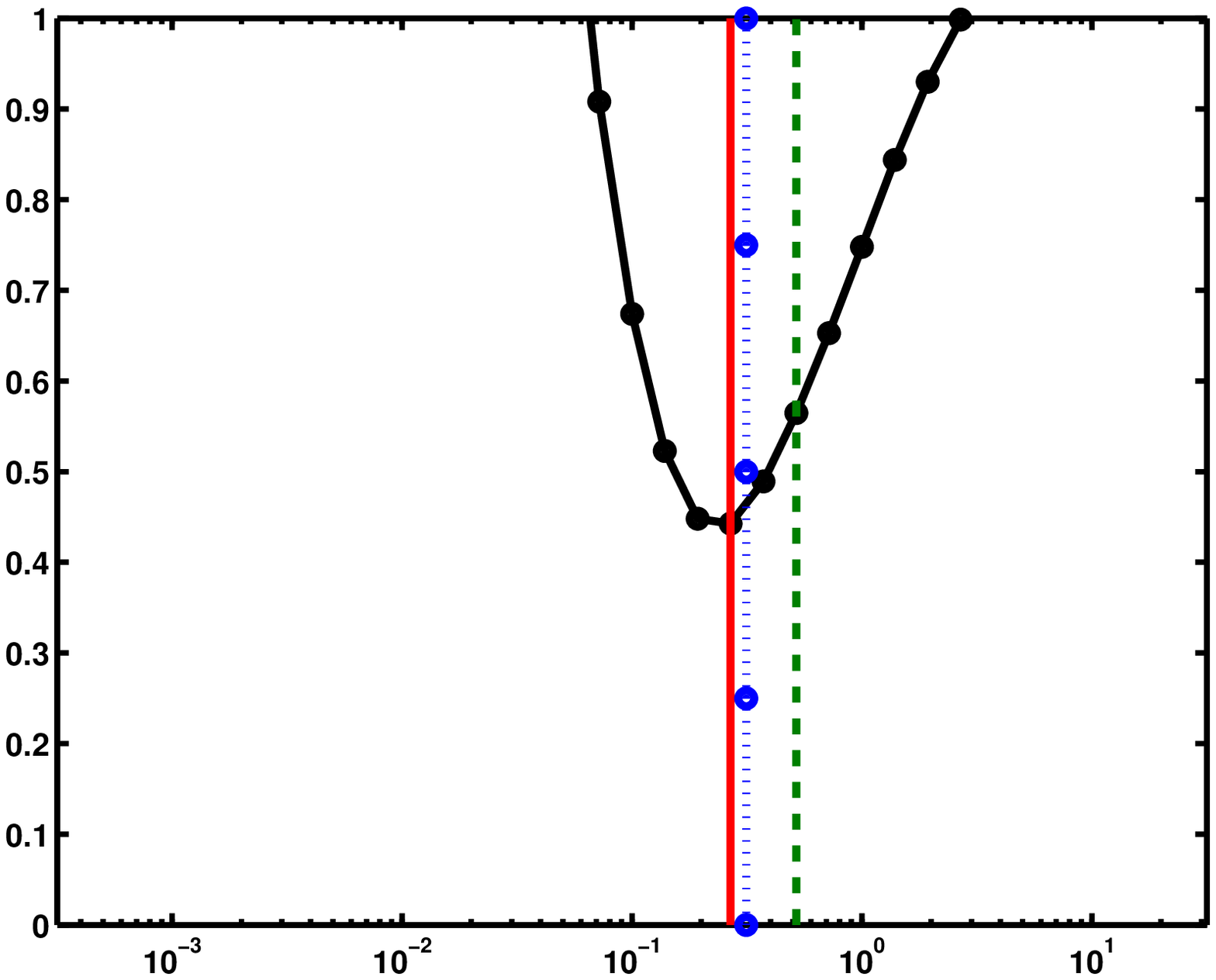}}
\subfigure[$L=L_1$]{\includegraphics[width=1.7in]{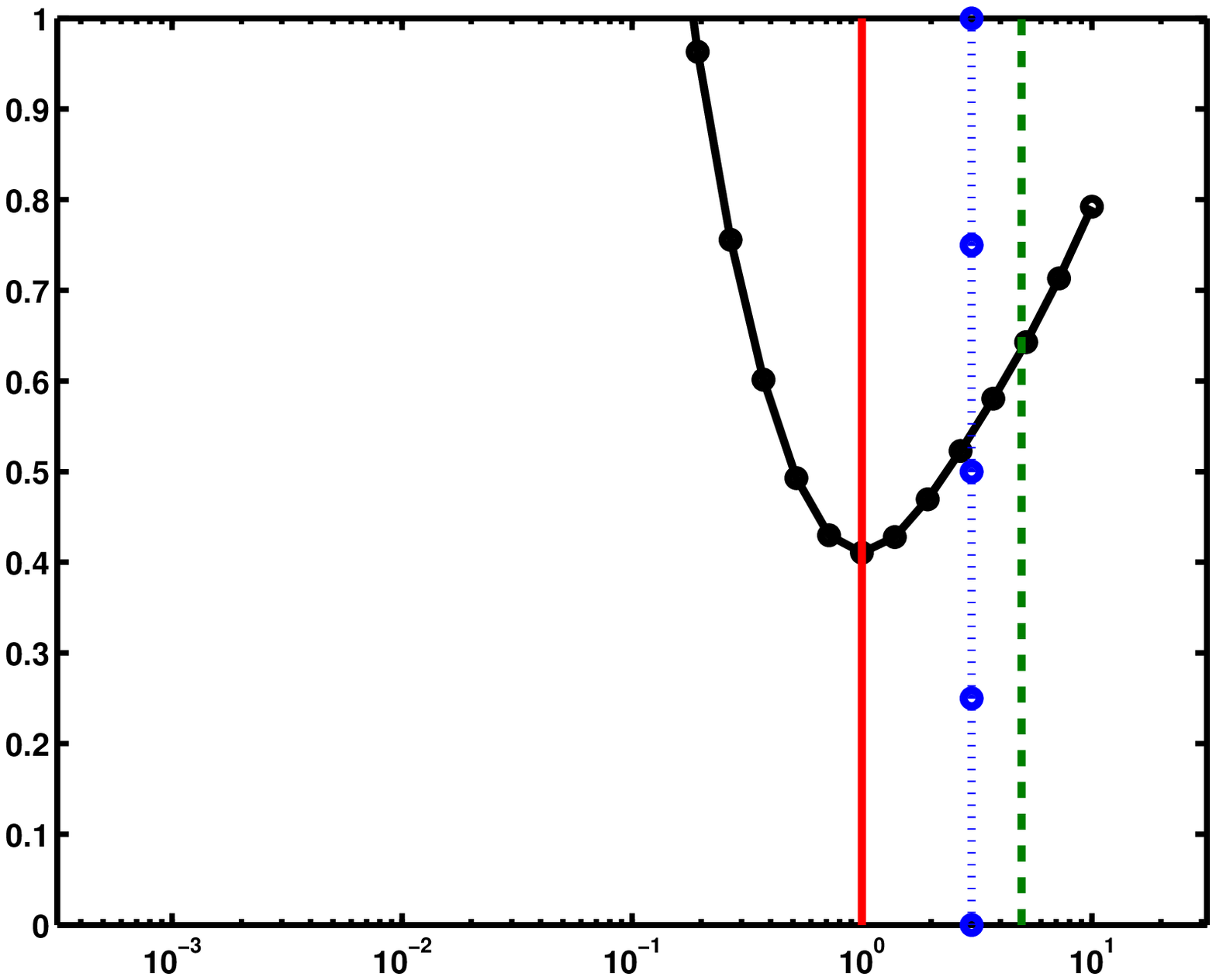}}
\subfigure[$L=L_2$]{\includegraphics[width=1.7in]{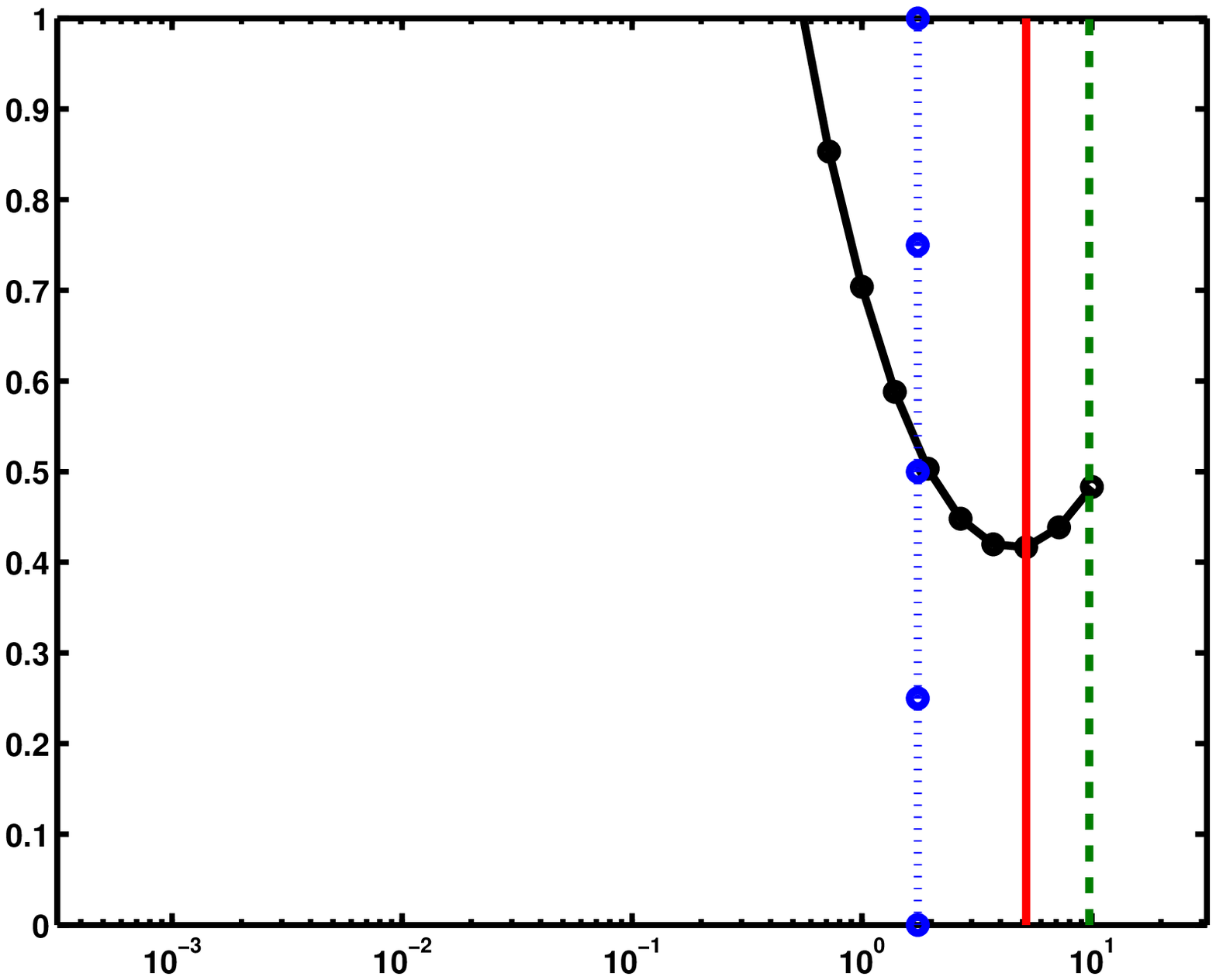}}
\subfigure[$L=I$]{\includegraphics[width=1.7in]{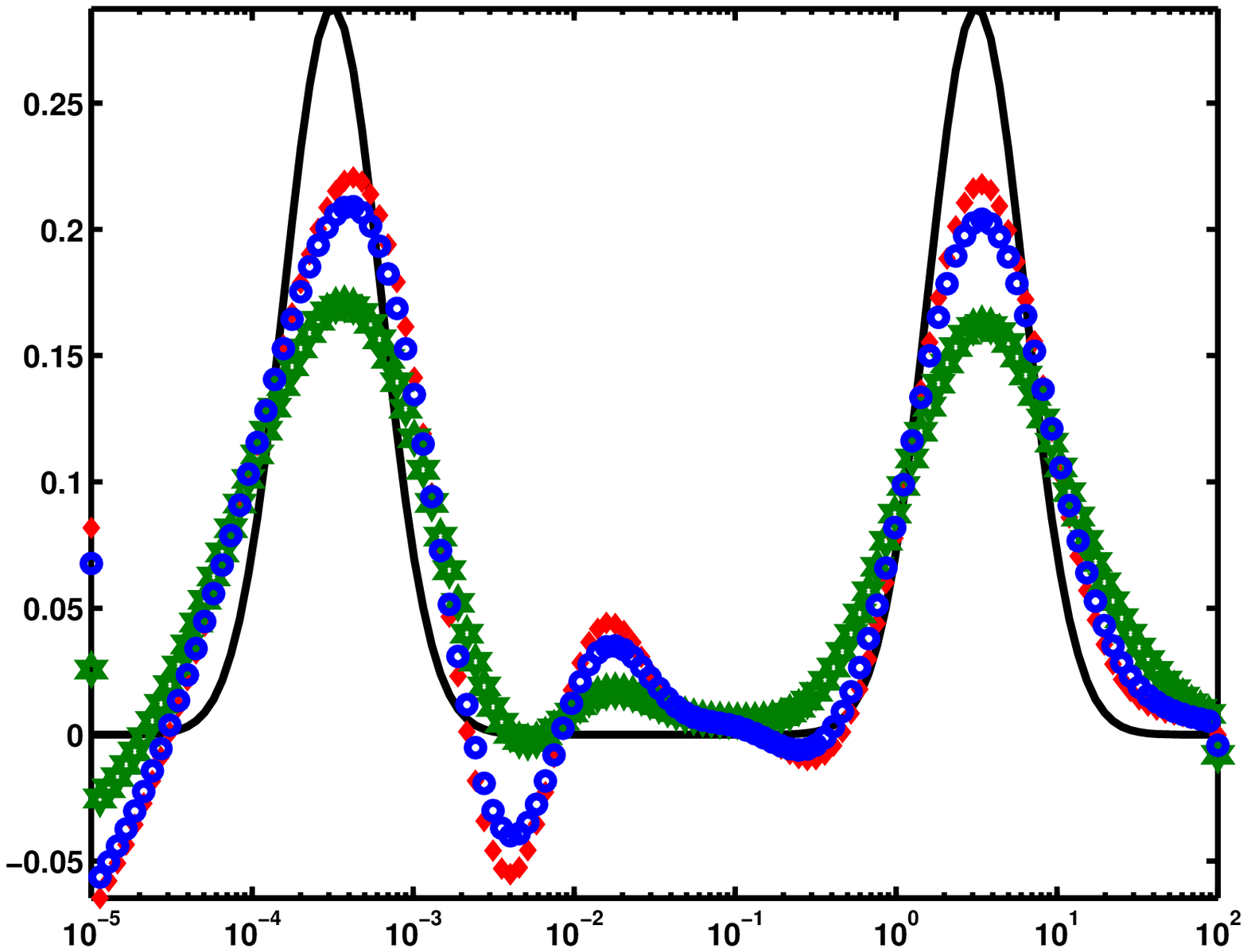}}
\subfigure[$L=L_1$]{\includegraphics[width=1.7in]{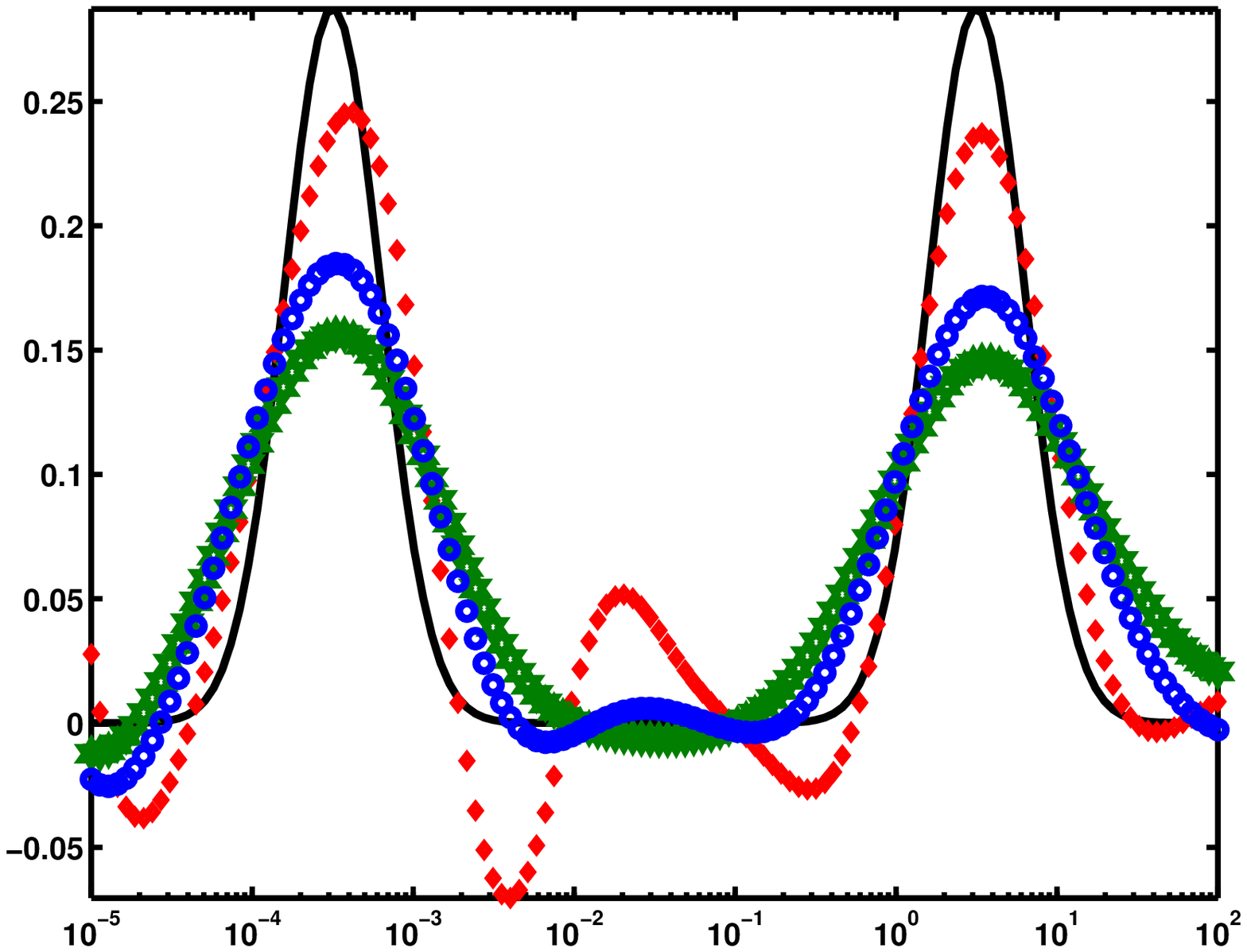}}
\subfigure[$L=L_2$]{\includegraphics[width=1.7in]{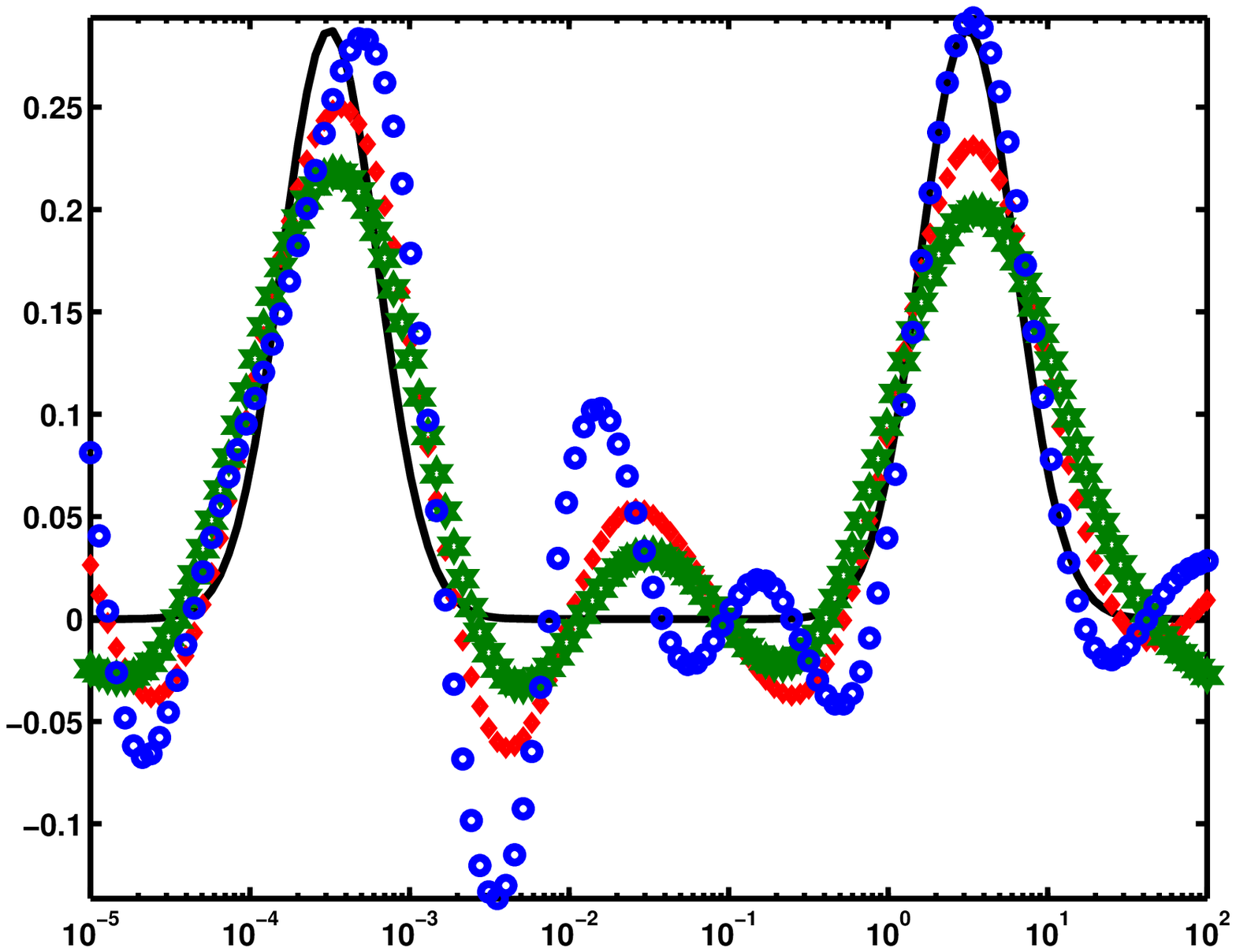}}
\caption{Mean error and example LS solutions.  $5\%$ noise. LN-A data set matrix $A_4$ }
\label{fig-lambdachoiceLN2A4HNLS}
\end{figure}

 \begin{figure}[!h]
  \centering
\subfigure[$L=I$]{\includegraphics[width=1.7in]{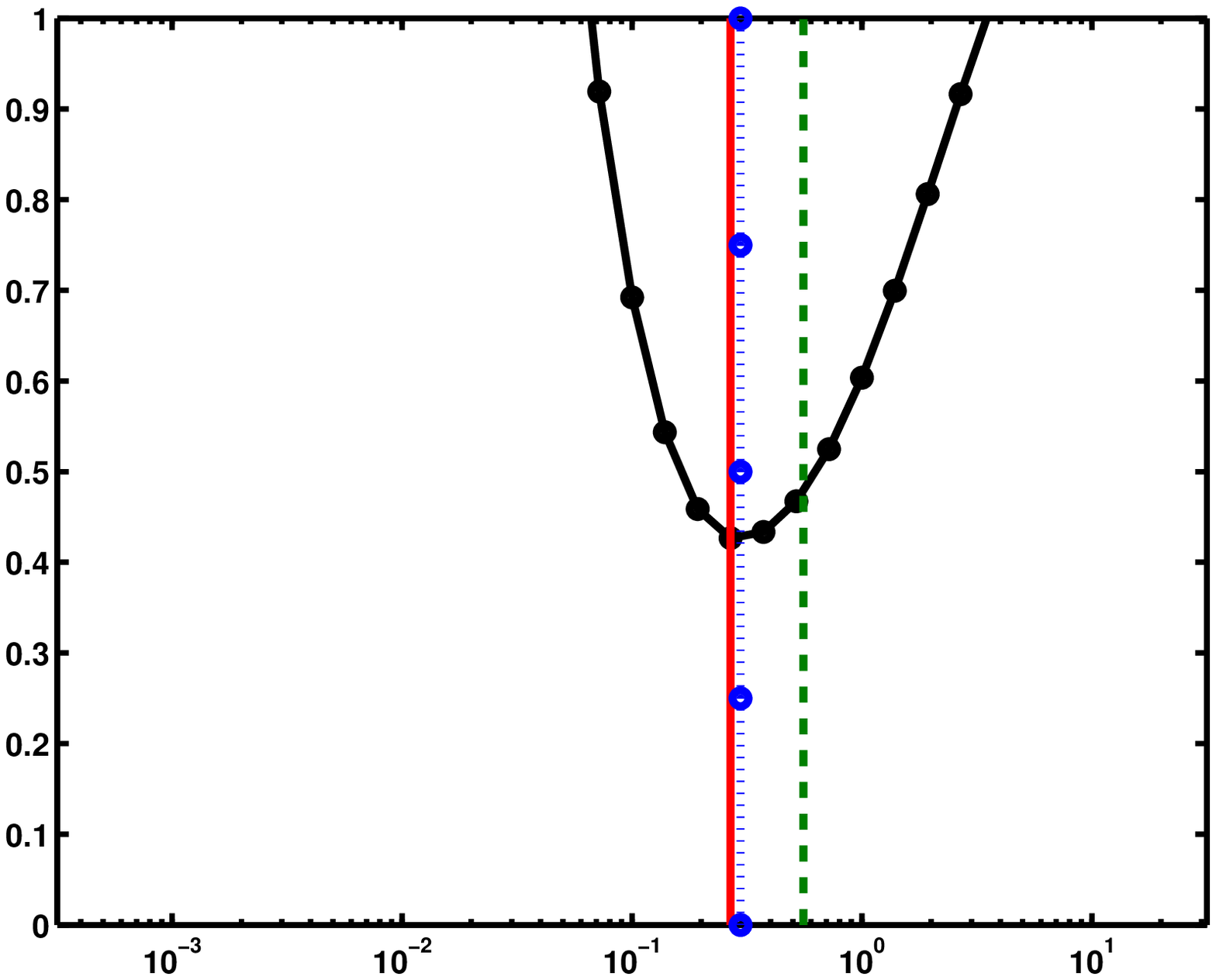}}
\subfigure[$L=L_1$]{\includegraphics[width=1.7in]{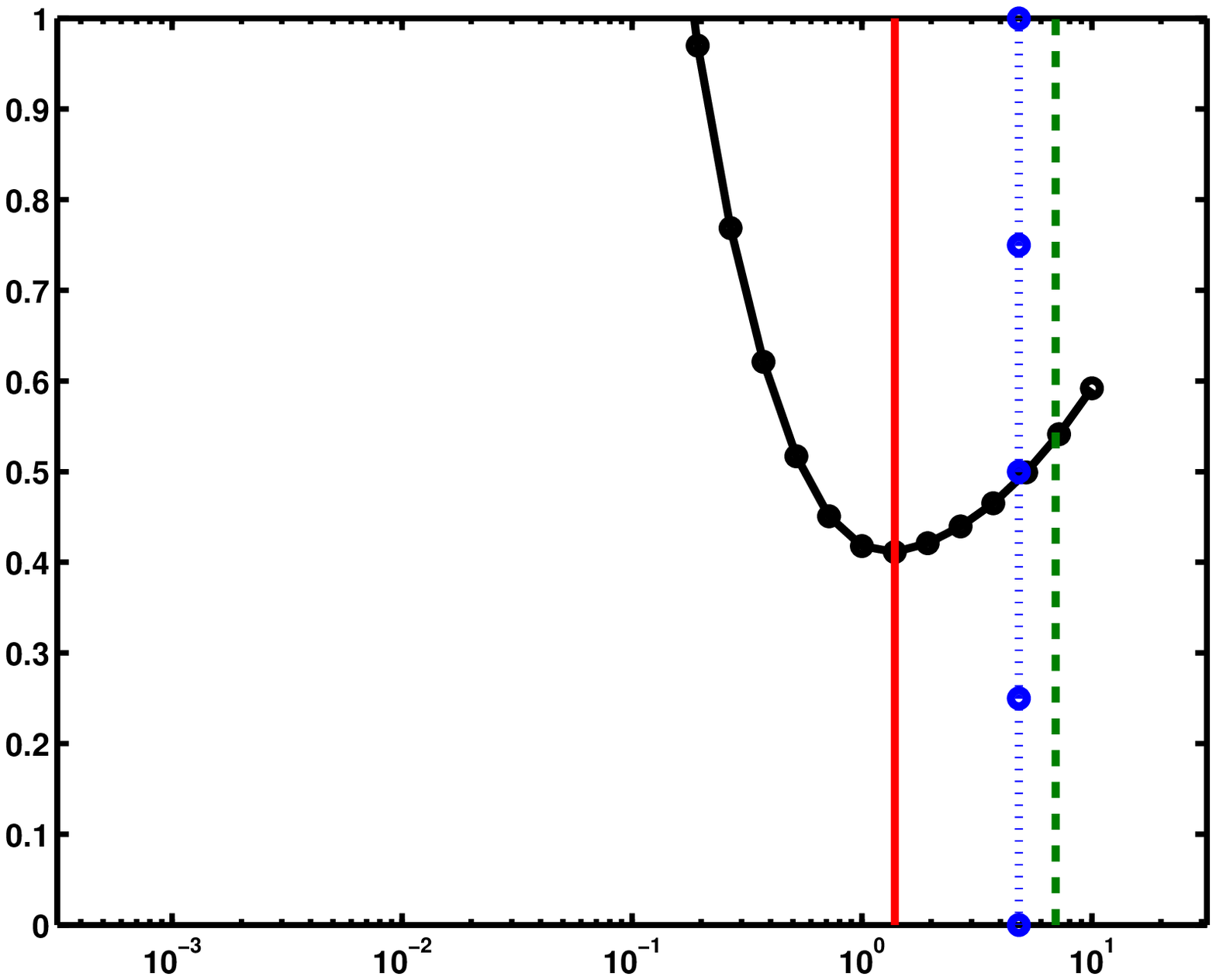}}
\subfigure[$L=L_2$]{\includegraphics[width=1.7in]{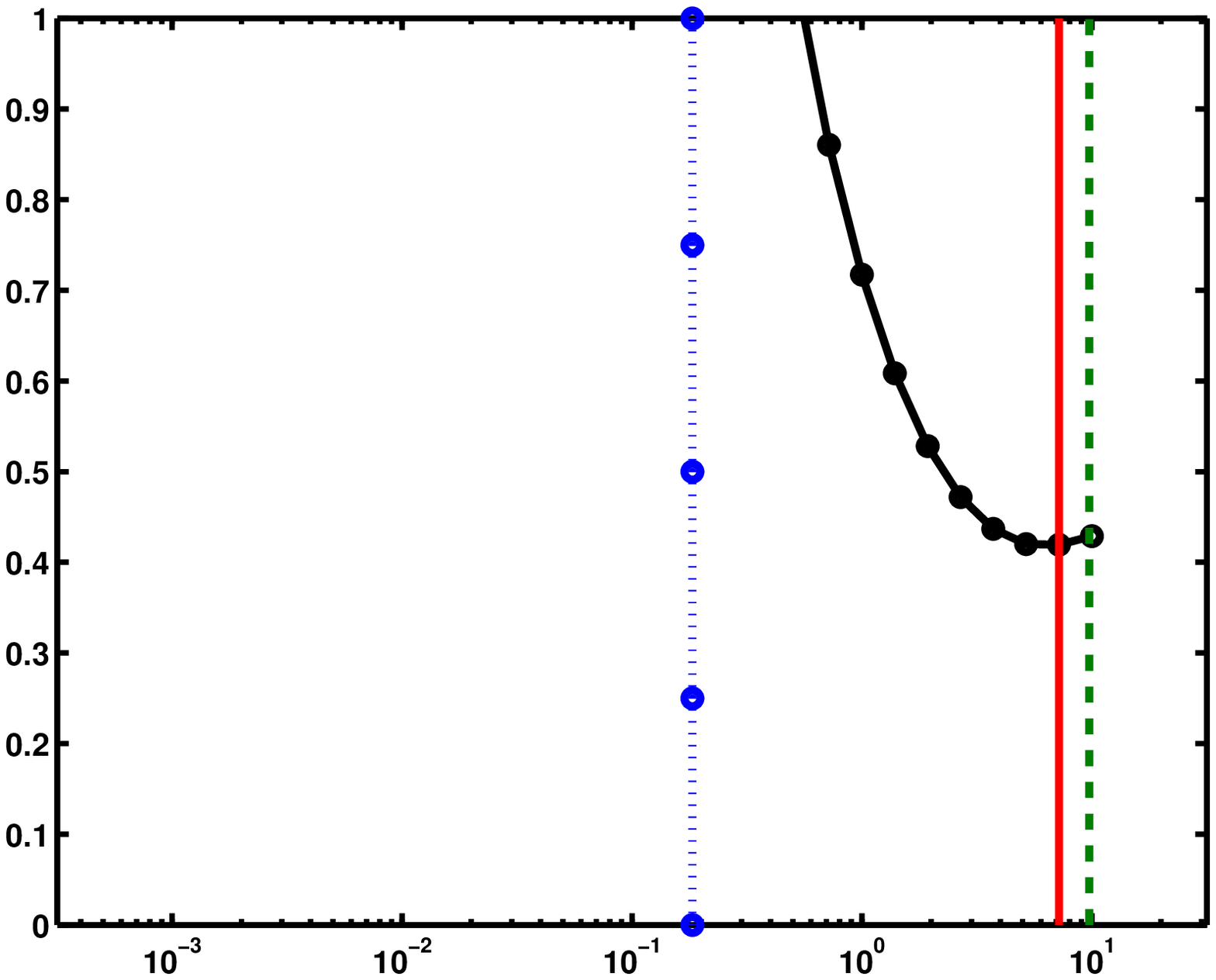}}
\subfigure[$L=I$]{\includegraphics[width=1.7in]{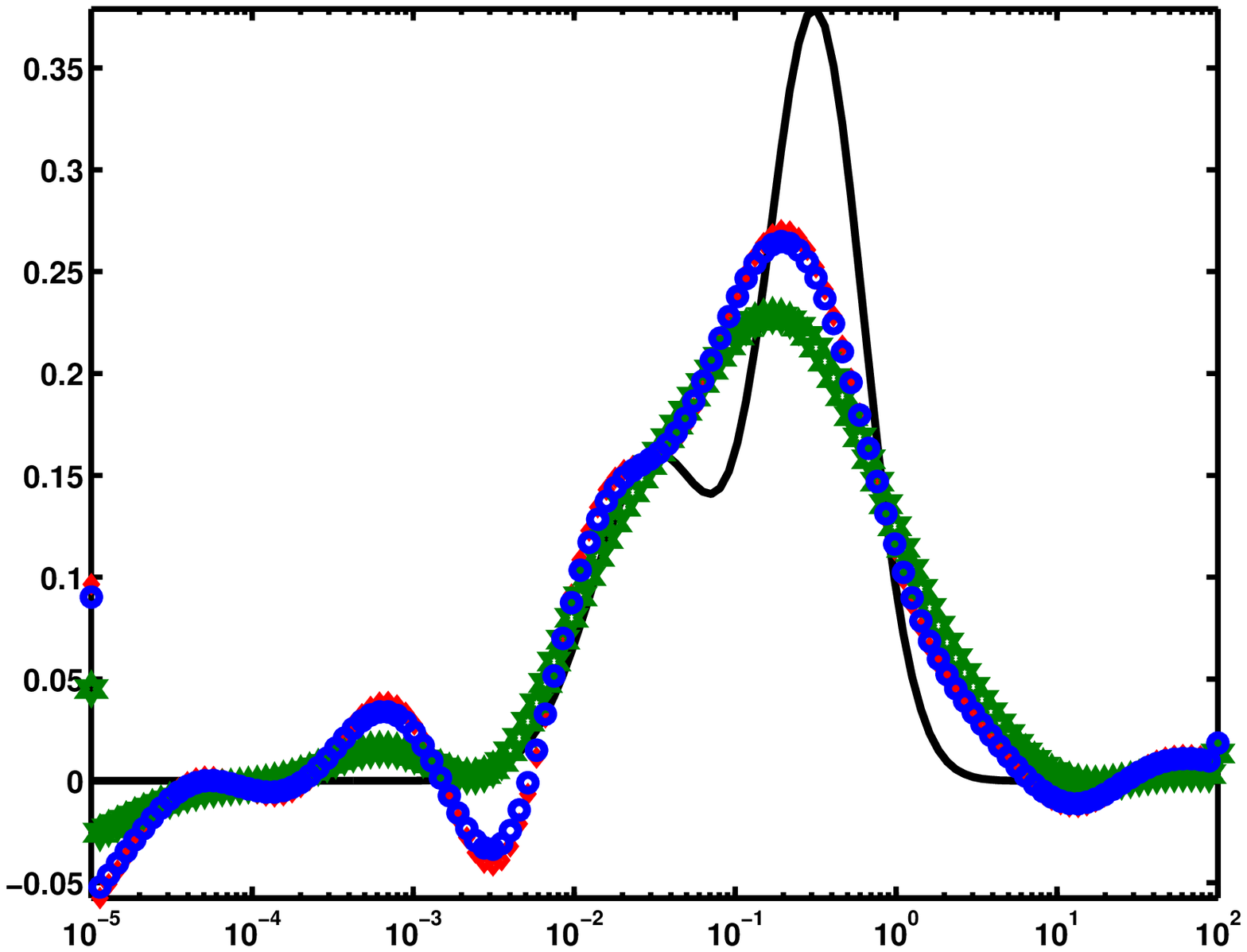}}
\subfigure[$L=L_1$]{\includegraphics[width=1.7in]{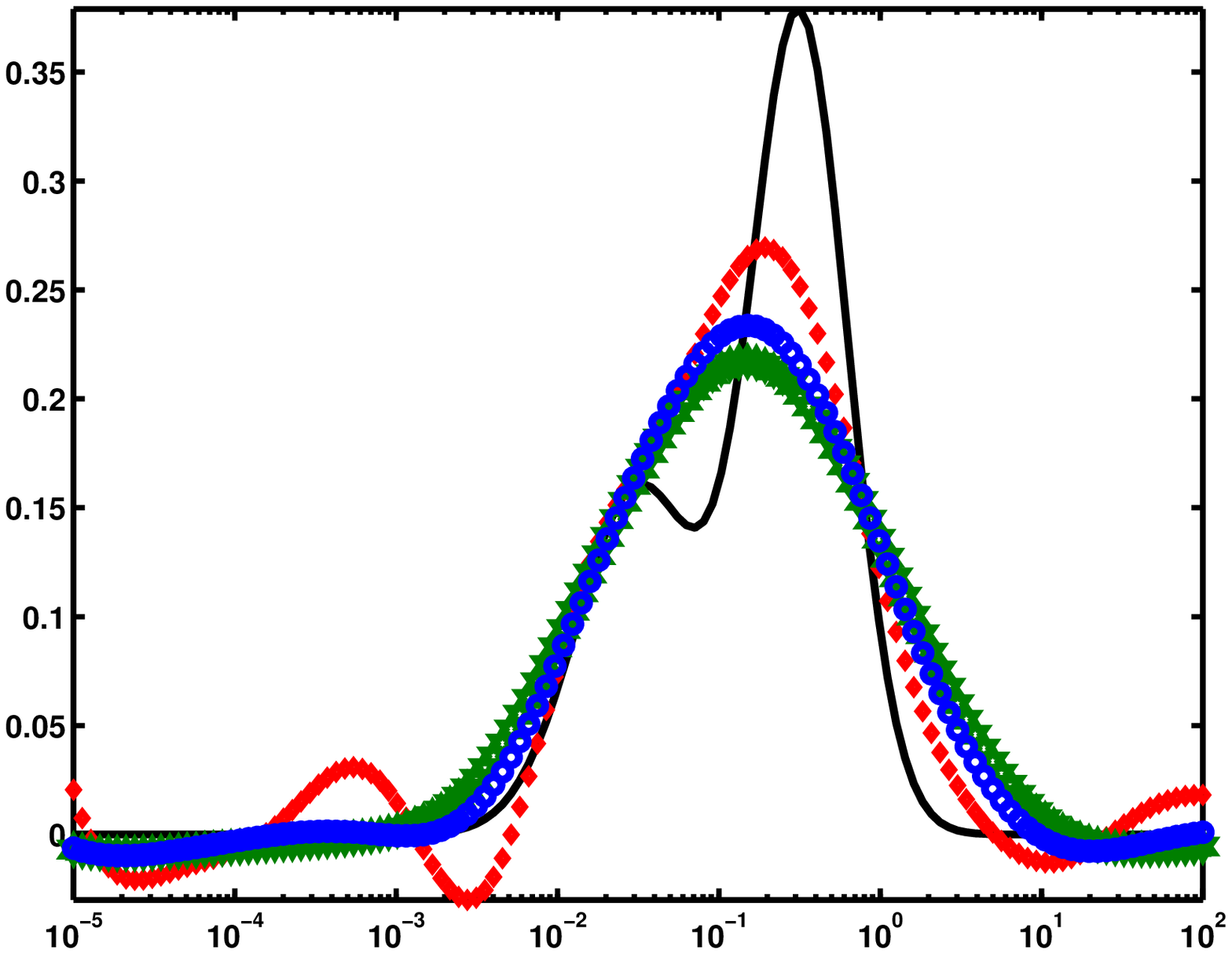}}
\subfigure[$L=L_2$]{\includegraphics[width=1.7in]{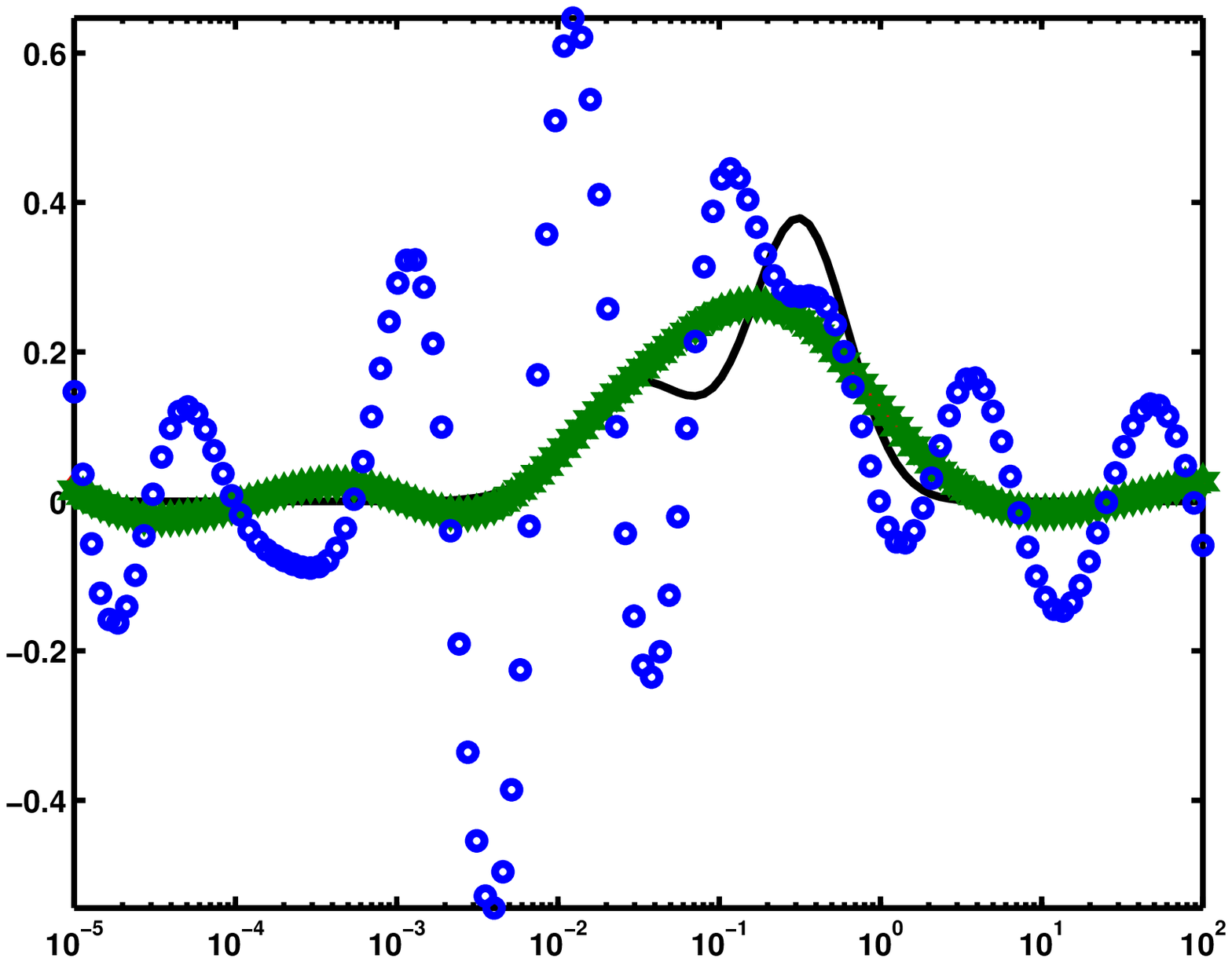}}
\caption{Mean error and example LS solutions.  $5\%$ noise. LN-B data set matrix $A_4$ }
\label{fig-lambdachoiceLN5A4HNLS}
\end{figure}

 \begin{figure}[!h]
  \centering
\subfigure[$L=I$]{\includegraphics[width=1.7in]{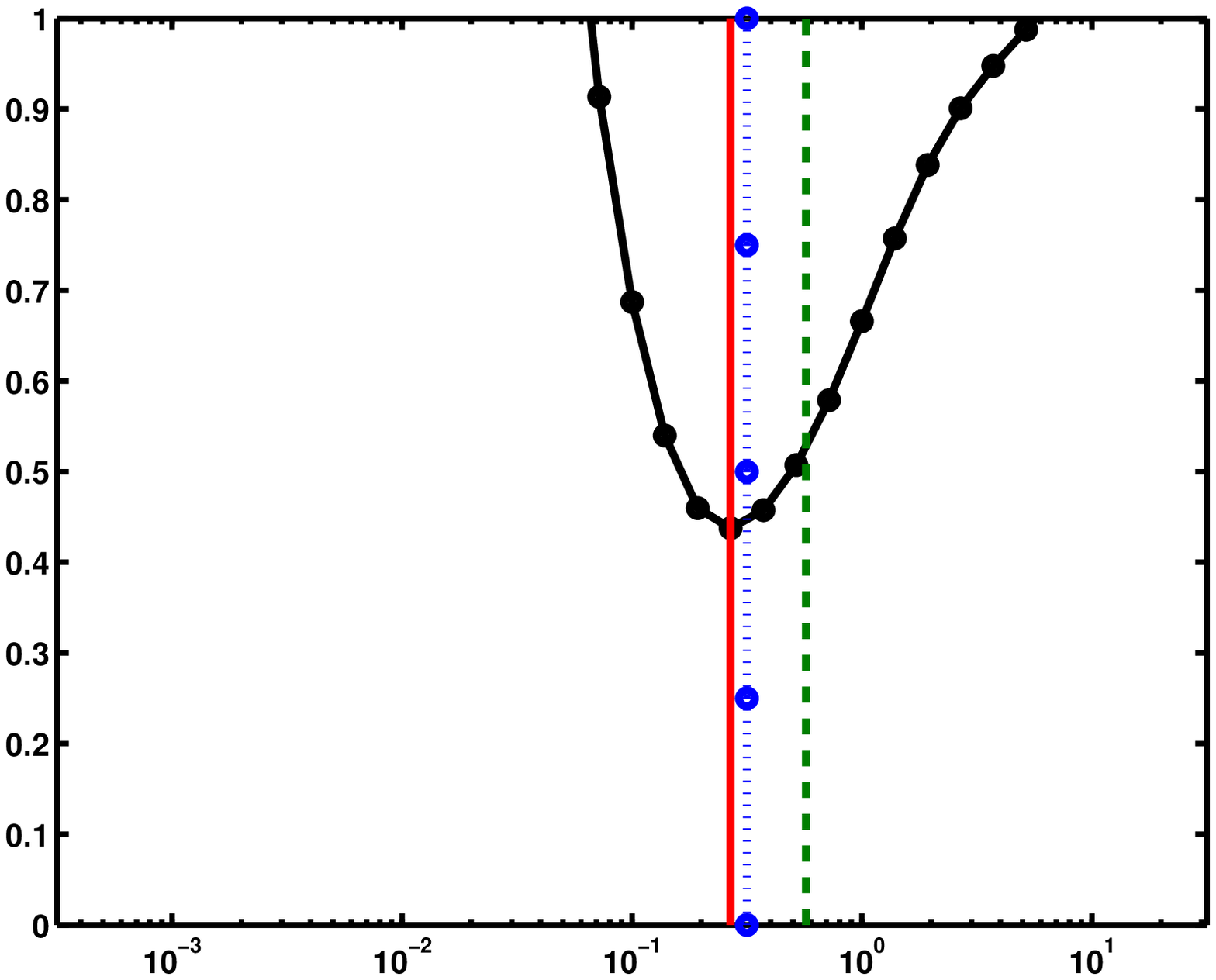}}
\subfigure[$L=L_1$]{\includegraphics[width=1.7in]{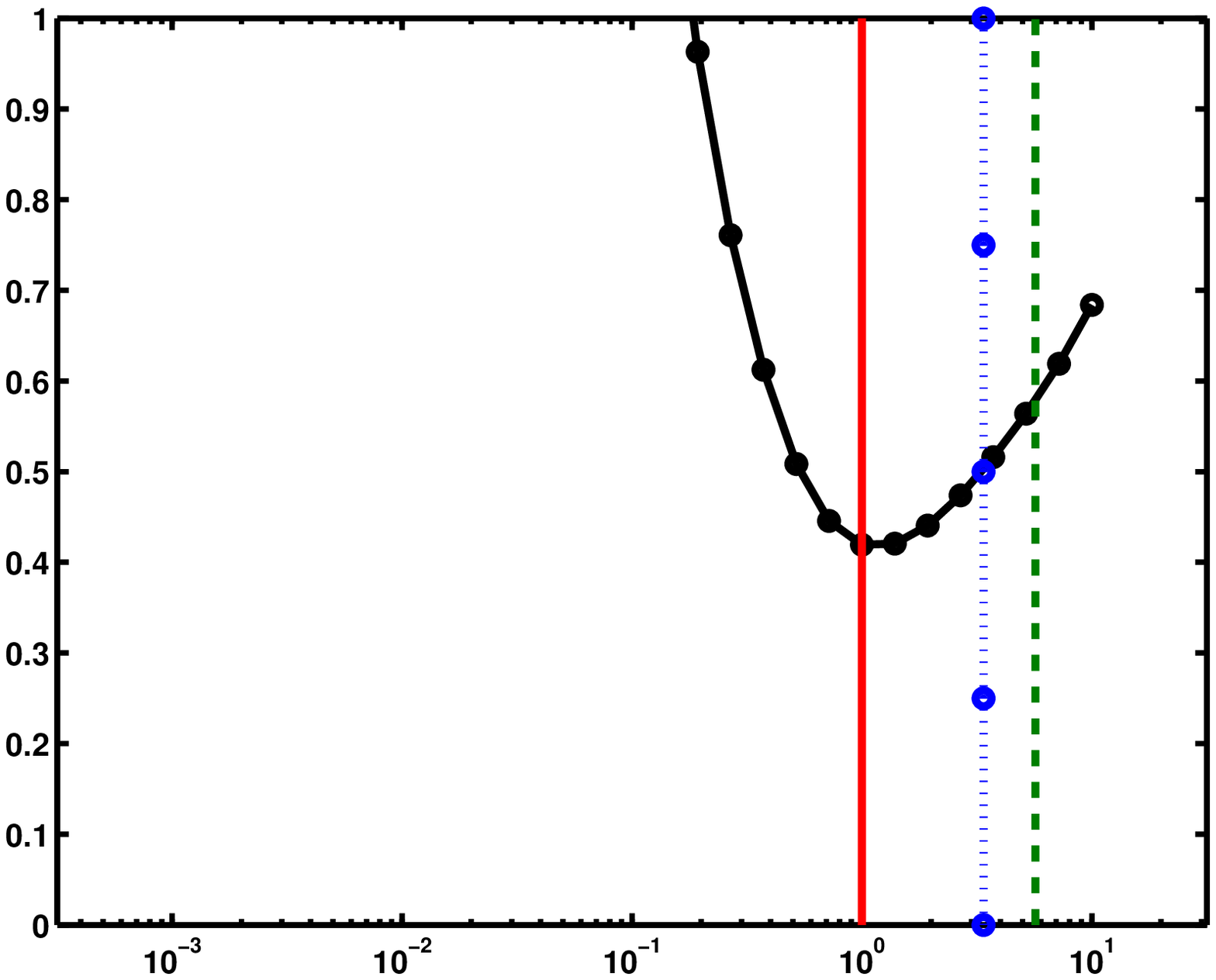}}
\subfigure[$L=L_2$]{\includegraphics[width=1.7in]{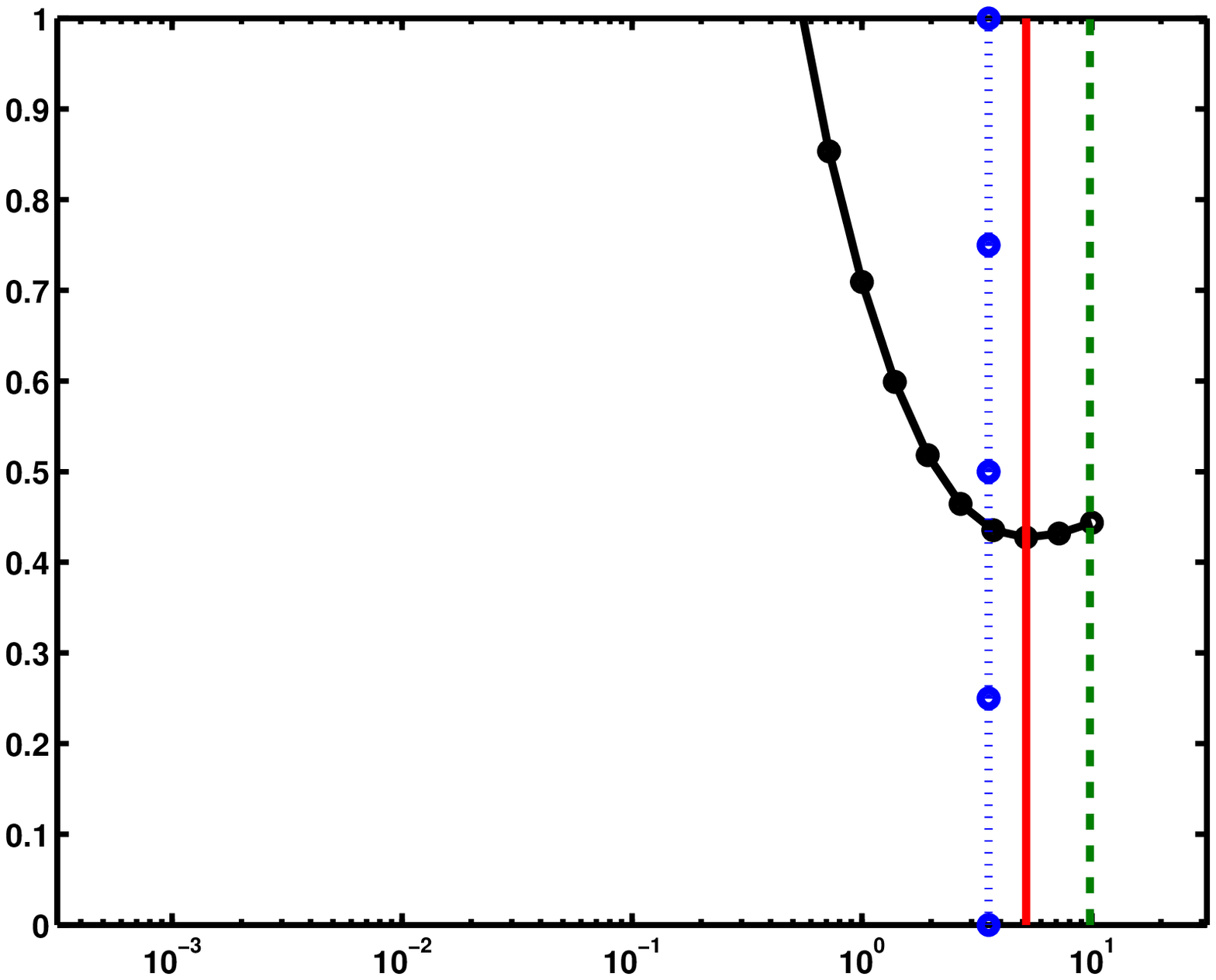}}
\subfigure[$L=I$]{\includegraphics[width=1.7in]{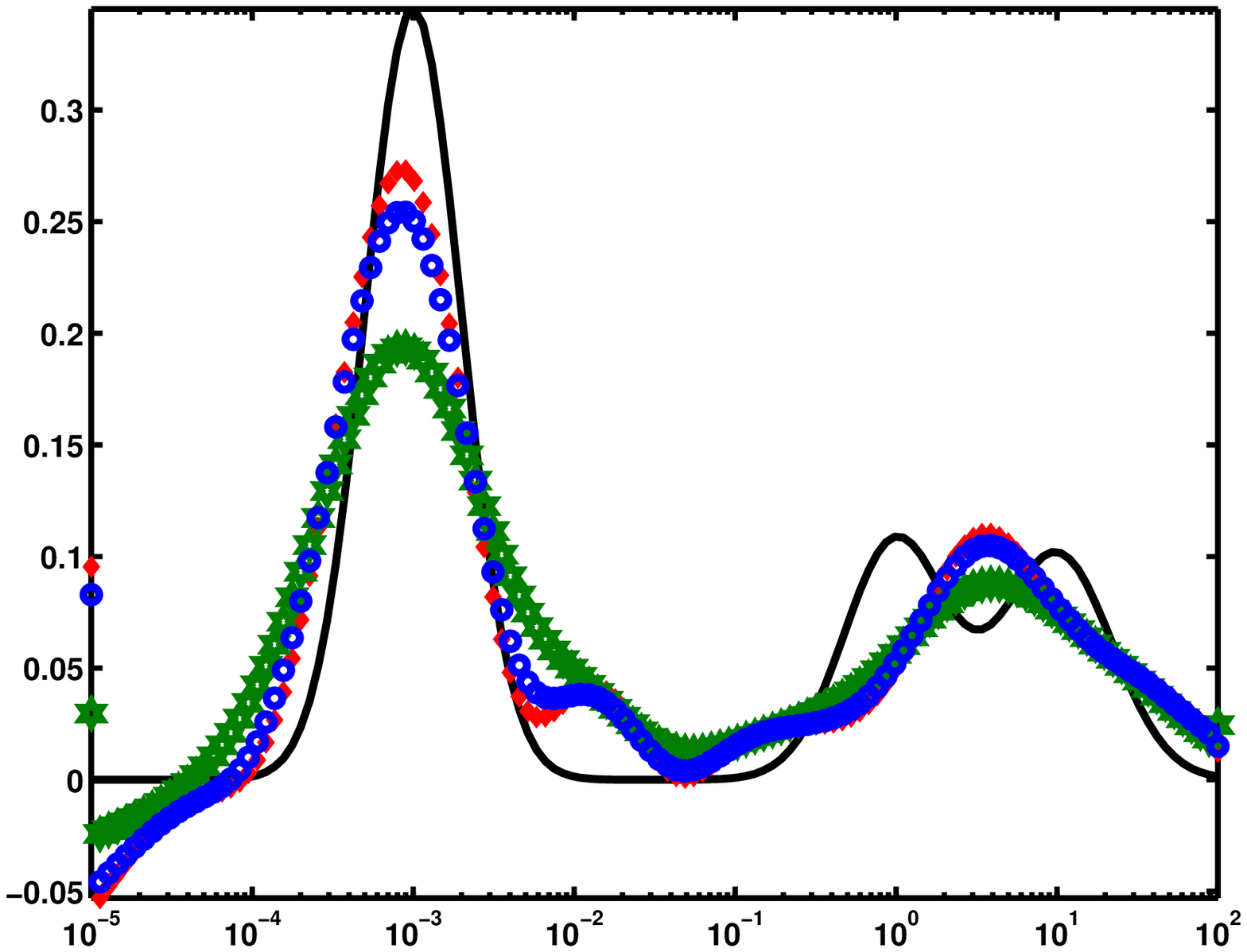}}
\subfigure[$L=L_1$]{\includegraphics[width=1.7in]{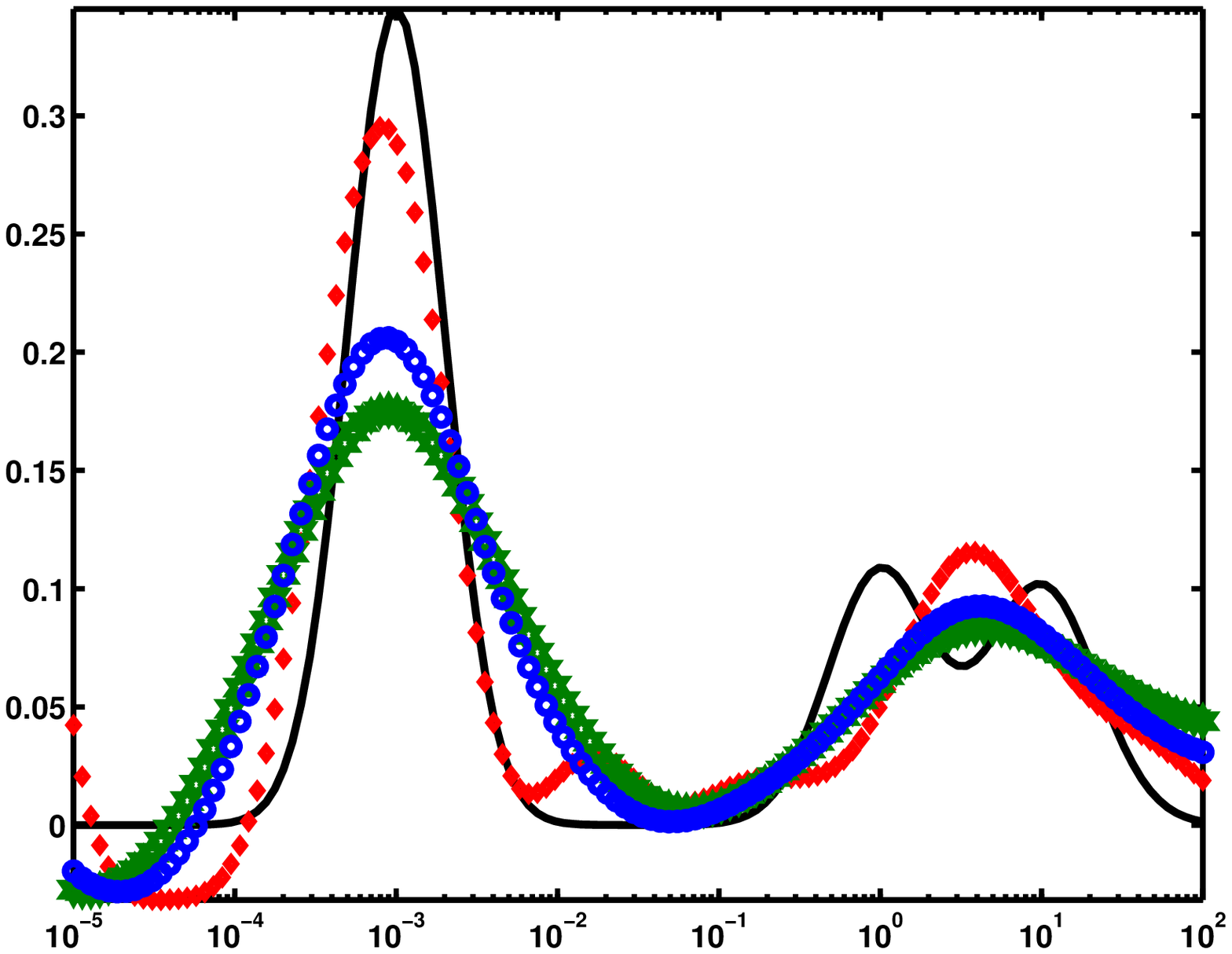}}
\subfigure[$L=L_2$]{\includegraphics[width=1.7in]{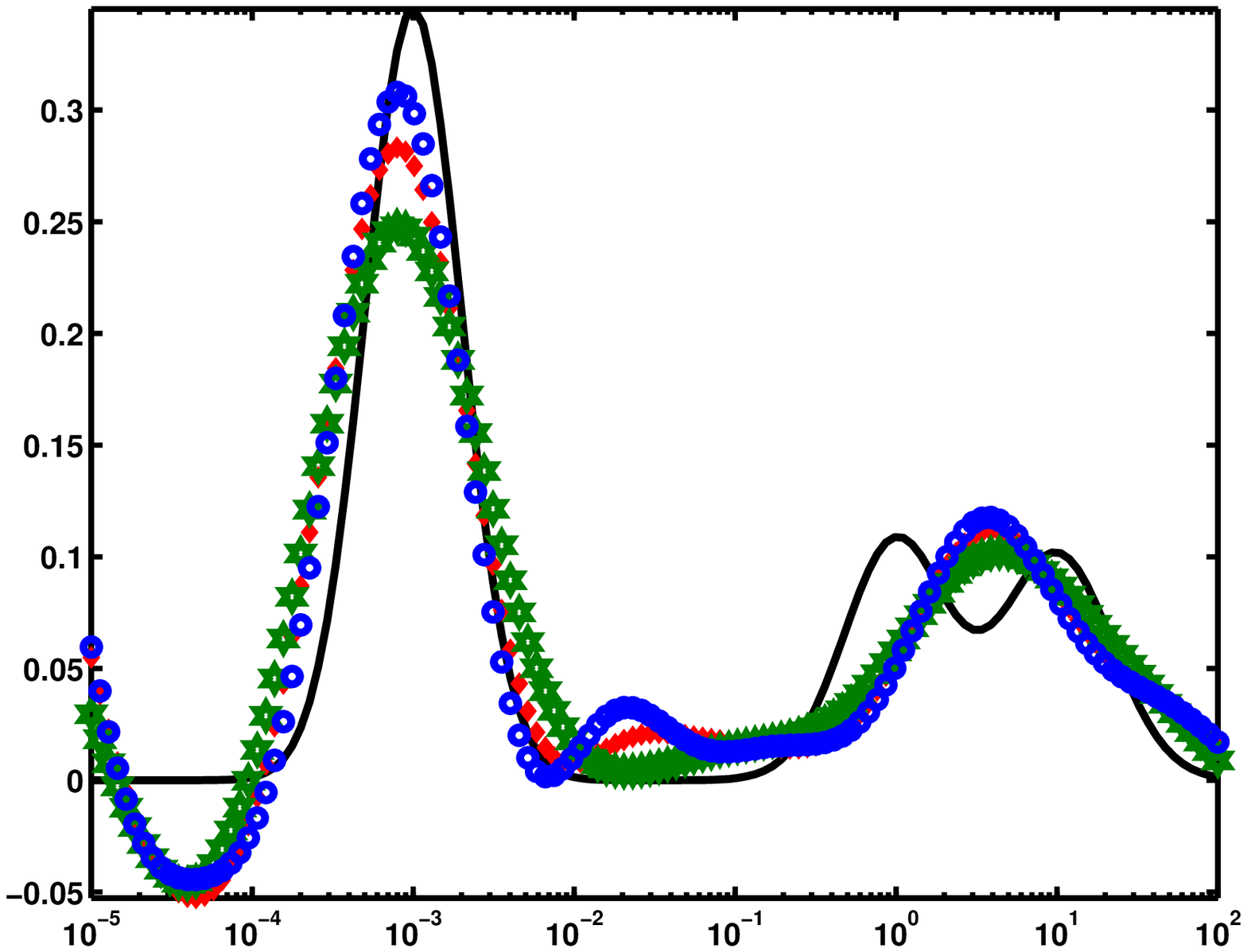}}
\caption{Mean error and example LS solutions.  $5\%$ noise. LN-C data set matrix $A_4$}
\label{fig-lambdachoiceLN6A4HNLS}
\end{figure}

\section{Acknowledgements}
Authors Hansen, Hogue and Sander were supported by  NSF CSUMS grant DMS 0703587: ``CSUMS:
Undergraduate Research Experiences for Computational Math Sciences Majors at ASU". Renaut was supported by NSF MCTP grant DMS 1148771: ``MCTP: Mathematics Mentoring Partnership Between Arizona State University and the Maricopa County Community College District",  NSF grant DMS 121655: 
``Novel Numerical Approximation Techniques for Non-Standard Sampling Regimes",  and AFOSR grant 025717 ``Development and Analysis of Non-Classical Numerical Approximation Methods".

\end{document}